\documentclass[a4paper,12pt]{amsart}
\usepackage[frenchb]{babel}
\usepackage[latin1]{inputenc}
\usepackage{pdfsync}
\usepackage[T1]{fontenc}   
\usepackage{amsmath}
\usepackage{amssymb}
\usepackage{amsxtra}
\usepackage{amscd}
\usepackage{amsthm}
\usepackage{amsfonts}
\usepackage{eucal}
\usepackage[all]{xy}
\usepackage{graphicx}
\usepackage{comment}
\usepackage{epsfig}
\usepackage{psfrag}
\usepackage{mathrsfs}
\usepackage{amscd}
\usepackage{rotating}
\usepackage{lscape}
\usepackage{amsbsy}
\usepackage{verbatim}
\usepackage{moreverb}
\usepackage{url}

\newcommand{\nc}{\newcommand}
\nc{\renc}{\renewcommand}

\nc\restr[2]{{ 
  \left.\kern-\nulldelimiterspace    #1  
  \vphantom{\big|}  
  \right|_{#2}  
  }}

\newtheorem{thm}{Théorème}[section]
\newtheorem{prop}[thm]{Proposition}
\newtheorem{lem}[thm]{Lemme}
\newtheorem{sous-lem}[thm]{Sous-lemme}
\newtheorem{cor}[thm]{Corollaire}
\newtheorem{conj}[thm]{Conjecture}

\newtheorem{def-prop}[thm]{Définition-Proposition}
 \theoremstyle{definition}
\newtheorem{defi}[thm]{Définition}
\newtheorem{rem}[thm]{Remarque}

\newtheorem{construction}[thm]{Construction}
\newtheorem{notation}[thm]{Notation}
\newtheorem{observation}[thm]{Observation}
\numberwithin{equation}{section}

\renc{\sec}{\section}
\nc{\ssec}{\subsection}
\nc{\sssec}{\subsubsection}

\nc{\thmref}[1]{théorème~\ref{#1}}
\nc{\secref}[1]{paragraphe~\ref{#1}}
\nc{\lemref}[1]{lemme~\ref{#1}}
\nc{\defiref}[1]{définition~\ref{#1}}
\nc{\propref}[1]{proposition~\ref{#1}}
\nc{\corref}[1]{corollaire~\ref{#1}}
\nc{\constructionref}[1]{construction~\ref{#1}}
\nc{\conjref}[1]{conjecture~\ref{#1}}
\nc{\remref}[1]{remarque~\ref{#1}}
\nc{\questref}[1]{question~\ref{#1}}
\nc\Omegasour{\hbox{$\buildrel\smile\over{\vrule height 6pt depth 0pt width 0pt \smash \Omega}$}}

\nc{\on}{\operatorname}

\nc\wt{\widetilde}
\nc\wh{\widehat}
\nc\ol{\ov}
\nc{\oc}[1]{{\overset{\circ}{#1}}}
\nc{\ov}[1]{{\overline{#1}}}
\nc{\isor}[1]{{\xrightarrow[\raisebox{0.25 em}{\smash{\ensuremath{\sim}}}]{#1}}}
 \nc{\isol}[1]{{\xleftarrow[\raisebox{0.25 em}{\smash{\ensuremath{\sim}}}]{#1}}}
\nc{\modmod}{/ \! \! /}

\nc{\mc}{\mathcal}
\nc{\mf}{\mathfrak}
\nc{\mr}{\mathrm}
\nc{\mb}{\mathbb}
\nc{\mbf}{\mathbf}

\nc{\R}{{\mathbb R}}
\nc{\Z}{{\mathbb Z}}
\nc{\N}{{\mathbb N}}
\nc{\C}{{\mathbb C}}
\nc{\Q}{{\mathbb Q}}

\nc{\Fq}{{\mathbb F}_q}
\nc{\Fl}{{\mathbb F}_\ell}
\nc{\Fqbar}{\ol{{\mathbb F}_q}}
\nc{\Flbar}{\ol{{\mathbb F}_\ell}}
\nc{\Zl}{{\mathbb Z}_\ell}
\nc{\Zlbar}{\ol{{\mathbb Z}_\ell}}
\nc{\Ql}{{\mathbb Q}_\ell}
\nc{\Qlbar}{\ol{{\mathbb Q}_\ell}}
\nc{\hl}{\overset{\leftarrow}h{}}
\nc{\hr}{\overset{\rightarrow}h{}}
\nc{\Gr}{{\on{Gr}}}
\nc{\Hecke}{\on{Hecke}}
 \nc{\Hom}{\on{Hom}}
 \nc{\Coker}{\on{Coker}}
 \nc{\Ker}{\on{Ker}}
 \nc{\Lie}{\on{Lie}}
\nc{\Loc}{\on{Loc}}
\nc{\Pic}{\on{Pic}}
\nc{\Bun}{\on{Bun}}
\nc{\IC}{\on{IC}}
\nc{\Aut}{\on{Aut}}
\nc{\Perv}{\on{Perv}}
\nc{\pos}{{\on{pos}}}
\nc{\Sym}{\on{Sym}}

\nc{\ta} {{}^\tau}
\nc {\tu}[1]{{}^{\tau^{#1}}\!}

\nc{\Id}{\on{Id}}
\nc{\Fil}{\on{Fil}}
\nc{\pr}{\on{pr}}
\nc{\Res}{\on{Res}}
\nc{\cusp}{\on{cusp}}
\nc{\Frob}{\on{Frob}}
\nc{\diag}{\Delta}
\nc{\gr}{\on{gr}}
\nc{\Inj}{\on{Inj}}
\nc{\Bl}{\on{Bl}}
\nc{\dem}{\noindent {\bf Démonstration. }}
\nc{\cqfd}{{\ }\hfill $\square$ \vskip 1mm}
\nc{\s}[1]{\langle #1 \rangle}
\nc{\Cht}{\on{Cht}}
\nc{\isom}{\overset {\thicksim}{\to}}
\nc{\sm}{\smallsetminus}

\begin{document}

\title[$G$-chtoucas et paramétrisation de Langlands]{
Chtoucas pour les groupes réductifs et paramétrisation de  Langlands globale }

\author{Vincent Lafforgue}
\thanks{L'auteur fait partie de l'ANR-13-BS01-0001-01}
\address{Vincent Lafforgue: CNRS et Institut Fourier, UMR 5582, Université Grenoble Alpes, 
 100 rue des Maths, 38610 Gières, France.}

\date{\today}
\maketitle

      \section*{Introduction}
          
           Soit  $\Fq$ un corps fini et $\ell$ un nombre premier ne divisant pas $q$. 
  Soit  $X$ une courbe projective lisse géométriquement irréductible  sur $\Fq$ et  $F$ son corps de fonctions. 
      Soit  $G$ un groupe réductif connexe sur $F$. 
      
    On montre dans cet article le sens ``automorphe vers Galois'' de la correspondance de Langlands  globale pour $G$ \cite{langlands67}. En fait  on construit une décomposition {\it canonique} de l'espace des formes automorphes cuspidales pour $G$ à valeurs dans $\Qlbar$ (ou plus précisément, lorsque $G$ n'est pas déployé,  d'une somme, indexée par 
    $\ker^{1}(F,G)$, des espaces de formes automorphes cuspidales  pour des formes intérieures de $G$). Cette décomposition canonique est indexée par les  paramètres de Langlands globaux $\ell$-adiques.
    
    On n'obtient  pas de résultat nouveau dans   le cas où $G=GL_{r}$ puisque tout était déjà connu par Drinfeld \cite{drinfeld78,Dr1,drinfeld-proof-peterson,drinfeld-compact} pour $r=2$  et Laurent Lafforgue \cite{laurent-inventiones} pour $r$ arbitraire (voir le chapitre \ref{GL-previous-works} pour le cas de $GL_{r}$).

   Cet article est totalement indépendant de la formule des traces d'Arthur-Selberg. Il  utilise les deux   ingrédients suivants : 
        \begin{itemize}
     \item les champs classifiants de chtoucas, introduits par   Drinfeld pour $GL_{r}$  
     \cite{drinfeld78,Dr1} et généralisés à tous les groupes  réductifs par Varshavsky
     \cite{var} 
   \item l'équivalence de  Satake géométrique de   Lusztig, Drinfeld, Ginzburg, et  Mirkovic--Vilonen 
   \cite{lusztig-satake,ginzburg,hitchin,mv}. 
     \end{itemize}
     
    Pour énoncer le théorème principal on suppose que $G$ est déployé.   
     On note  $\wh G$ le groupe dual de  Langlands de  $G$,  considéré comme un groupe déployé sur  $\Ql$. 
        Ses racines et ses poids sont les coracines et les copoids de $G$, et vice-versa (voir \cite{borel-corvallis}  pour plus de  détails).

Pour des raisons topologiques on doit travailler avec des extensions finies de $\Ql$ au lieu de  $\Qlbar$.  Soit $E$ une extension finie de $\Ql$ contenant une racine carrée de $q$ et $\mc O_{E}$ son anneau d'entiers.

        Soit  $v$ une place de $X$. On note  $\mc O_{v}$ l'anneau local complété en  $v$ et  $F_{v}$ son corps de fractions. 
   On a l'isomorphisme de   Satake    $[V]\mapsto h_{V,v}$ de l'anneau des représentations (de dimension finie) de  $\wh G$ à coefficients dans $E$ vers l'algèbre de  Hecke  
      $C_{c}(G(\mc O_{v})\backslash G(F_{v})/G(\mc O_{v}),E)$ (voir  
      \cite{satake,cartier-satake,gross}). En fait les $h_{V,v}$  pour $V$ irréductible forment une base sur $\mc O_{E}$ de $C_{c}(G(\mc O_{v})\backslash G(F_{v})/G(\mc O_{v}),\mc O_{E})$. 
      On note   $\mb A=\prod_{v\in |X|} ' F_{v}$ l'anneau des adèles  de $F$ et $\mb O=
      \prod_{v\in |X|} \mc O_{v}$. 
             Soit   $N$ un sous-schéma fini de  $X$.   
   On note  
   $\mc O_{N}$ l'anneau des fonctions sur $N$, et 
   \begin{gather}\label{def-K-N}K_{N}=\on{Ker}(G(\mb O)\to G(\mc O_{N}))\end{gather} le sous-groupe compact ouvert de $G(\mb A)$ associé au niveau  $N$. On fixe un réseau   $\Xi\subset Z(F)\backslash Z(\mb A)$ (où $Z$ est le centre de $G$).  
   Une fonction $f\in C_{c} (G(F)\backslash G(\mb A)/K_N \Xi,E)$ 
   est dite cuspidale si pour tout 
    parabolique $P\subsetneq G$, de  Levi $M$ et de radical unipotent $U$,    le terme constant $f_{P}: g\mapsto \int_{U(F)\backslash U(\mb A)}f(ug)$ est  nul comme fonction sur $U(\mb A)M(F)\backslash G(\mb A)/K_{N}\Xi$.  
 On rappelle que le 
 $E$-espace vectoriel  
 $C_{c}^{\rm{cusp}}(G(F)\backslash G(\mb A)/K_N \Xi,E)$ des formes automorphes cuspidales est   de dimension finie. Il est muni d'une structure de module sur l'algèbre de Hecke $C_{c}(K_{N}\backslash G(\mb A)/K_{N},E)$ : on convient  que la fonction caractéristique de $K_{N}$ est une unité et agit par l'identité et pour $f\in   C_{c}(K_{N}\backslash G(\mb A)/K_{N},E)$ on note $$T(f)\in \on{End}(C_{c}^{\rm{cusp}}(G(F)\backslash G(\mb A)/K_N \Xi,E))$$ l'opérateur de Hecke correspondant. 
    
 \subsection{Enoncé du théorème principal} 
   On construira  les opérateurs suivants, dits ``d'excursion''. 
    Soient $I$ un ensemble fini, $f$ une fonction sur 
$\wh G \backslash (\wh G)^{I}/\wh G    $ 
(quotient grossier de $(\wh G)^{I}$ par les translations à gauche et à droite par $\wh G $ diagonal), et 
$(\gamma_{i})_{i\in I}\in (\on{Gal}(\ov F/F))^{I}$. 
On construira  l'opérateur  d'excursion $$S_{I,f,(\gamma_{i})_{i\in I}}\in \mr{End}_{C_{c}(K_{N}\backslash G(\mb A)/K_{N},E)}( C_{c}^{\mr{cusp}}(G(F)\backslash G(\mb A)/K_{N}\Xi,E)). 
    $$ On montrera que ces opérateurs  engendrent une sous-algèbre commutative $\mc B$. 
    
    On ne sait pas si $\mc B$ est réduite mais 
par décomposition spectrale  on obtient néanmoins  une décomposition canonique 
 \begin{gather}\label{intro1-nu-dec-canonique}
 C_{c}^{\mr{cusp}}(G(F)\backslash G(\mb A)/K_{N}\Xi,\Qlbar)=\bigoplus_{\nu}
 \mf H_{\nu} \end{gather}
    où la somme directe dans le membre de droite est indexée par des  caractères $\nu$ de $\mc B$, et où 
   $\mf H_{\nu}$ est l'espace propre généralisé (ou ``espace caractéristique'') associé à $\nu$. On montrera ensuite qu'à tout caractère $\nu$ de $\mc B$ correspond un {\it unique} paramètre  de Langlands $\sigma$ (au sens du théorème suivant), caractérisé par \eqref{relation-fonda} ci-dessous. 
    En posant $ \mf H_{\sigma} =\mf H_{\nu}$, on en déduira le théorème suivant.

 \begin{thm}  \label{intro-thm-ppal}  (\thmref{dec-param-cor-thm})  
  On possède  une décomposition canonique de  
   $C_{c}(K_{N}\backslash G(\mb A)/K_{N},\Qlbar)$-modules
 \begin{gather}\label{intro1-dec-canonique}
 C_{c}^{\mr{cusp}}(G(F)\backslash G(\mb A)/K_{N}\Xi,\Qlbar)=\bigoplus_{\sigma}
 \mf H_{\sigma},\end{gather}
 où la somme directe dans le membre de droite est indexée par des paramètres de Langlands globaux, c'est-à-dire des classes de  $\wh G(\Qlbar)$-conjugaison de  morphismes 
       $\sigma:\on{Gal}(\ov F/F)\to \wh G(\Qlbar)$ 
       définis sur une extension finie de  $\Ql$, continus,  semi-simples et non ramifiés en dehors de  $N$.
       
     Cette décomposition  est caractérisée par la propriété suivante : 
       $ \mf H_{\sigma}$ est égal à l'espace propre généralisé $\mf H_{\nu}$ 
        associé au caractère $\nu$ de $\mc B$ défini par 
       \begin{gather}\label{relation-fonda}\nu(S_{I,f,(\gamma_{i})_{i\in I}})=f((\sigma(\gamma_{i}))_{i\in I}). \end{gather}
       
      Elle est compatible avec l'isomorphisme de  Satake en toute place  $v$  de $X\sm N$, c'est-à-dire que pour toute représentation irréductible $V$ de $\wh G$, 
       $T(h_{V,v})$ agit sur   $\mf H_{\sigma}$
       par multiplication par le scalaire  $\chi_{V}(\sigma(\Frob_{v}))$, où $\chi_{V}$ est le caractère de $V$ et $\Frob_{v}$ est un relèvement arbitraire d'un élément de  Frobenius   en $v$. 
     Elle   est aussi compatible avec la limite sur  $N$.   
    \end{thm}

La compatibilité avec l'isomorphisme de  Satake en les places     de $X\sm N$ montre que ce théorème réalise la ``correspondance'' de Langlands globale dans le sens ``automorphe vers Galois''.

Dans le chapitre  \ref{para-non-deploye}  on 
traitera le cas des groupes réductifs  non nécessairement déployés et on  prouvera le  \thmref{dec-param-cor-thm-non-deploye} qui est similaire  au théorème ci-dessus, à part que le membre de gauche est remplacé par une somme directe indexée par      $ \ker^ {1}(F, G) $ d'espaces de formes automorphes cuspidales pour des formes intérieures de  $ G $, et les   paramètres de Langlands $\sigma$ sont définis à l'aide du $L$-groupe. 

  Les  paragraphes  \ref{para-non-deploye-statement}
et \ref{complements-non-split}, qui  peuvent être lus en complément de cette introduction, contiennent les énoncés dans le cas non déployé, ainsi que les deux conjectures suivantes. 
La   \conjref{conj-arthur} affirme que   les paramètres de Langlands  intervenant dans la  décomposition \eqref{intro1-dec-canonique} (ou la décomposition analogue dans le cas non déployé) proviennent de paramètres d'Arthur elliptiques. 
      D'autre part la  conjecture \ref{conj-algebre-B-Qbar} affirme que cette décomposition   est définie sur  $\overline \Q$ et est indépendante de $\ell$.

Le chapitre \ref{mod-ell} montre que la décomposition \eqref{intro1-dec-canonique} existe aussi à coefficients dans  $\ov {\mathbb F_{\ell}}$. 
Le chapitre  \ref{para-meta} indique comment adapter ces méthodes au cas métaplectique.

 Dans le  cas où $G=GL_{r}$ les théorèmes inverses \cite{inverse-thm} et la formule du produit de Laumon \cite{laumon-produit} fournissent, par récurrence,  le sens ``Galois  vers automorphe'' (voir \cite{laurent-inventiones} ainsi que le chapitre \ref{GL-previous-works}). 
  Pour $G$ quelconque,   les conjectures  de Langlands consistent  plutôt en 
\begin{itemize}
\item une paramétrisation, obtenue dans le théorème ci-dessus, 
\item des formules de multiplicités d'Arthur pour les  $\mf H_{\sigma}$, que nous ne savons pas calculer avec les méthodes de cet article. \end{itemize}

Dans un article   avec Alain Genestier  \cite{genestier-lafforgue}, 
nous montrons  la paramétrisation  de Langlands locale et la compatibilité local-global.

Dans les  paragraphes  \ref{defi-chtou-intro} à \ref{subsection-intro-decomp}
 nous esquissons la preuve du  \thmref{intro-thm-ppal}. 
     
    \subsection{
    Chtoucas de Drinfeld pour les  groupes réductifs, d'après Varshavsky 
   } \label{defi-chtou-intro}
       
       La preuve repose  sur le fait que les champs de chtoucas, qui jouent un rôle analogue aux variétés de Shimura sur les corps de nombres, existent dans une  généralité beaucoup plus grande. En effet, alors que les variétés de Shimura  sont définies sur un ouvert du spectre d'un anneau d'entiers d'un  corps de nombres et sont associées à {\it un copoids  minuscule} du groupe dual, on possède
pour tout ensemble fini $I$, pour tout  niveau $N$ et pour toute représentation irréductible $W$  de $(\wh G)^{I}$  un champ de chtoucas $\on{Cht}_{N,I,W}$ qui est défini  sur $(X\sm N)^{I}$.

    Les chtoucas ont été introduits par  Drinfeld  \cite{drinfeld78,Dr1} pour  $GL_{r}$ (et $I=\{1,2\}, W=\mr{St}\otimes \mr{St^{*}}$)    et  généralisés  aux  groupes réductifs (et aux copoids arbitraires) par  Varshavsky dans  
    \cite{var} (entre-temps le cas des algèbres à  division a été étudié par  Laumon-Rapoport-Stuhler, Laurent Lafforgue, Ngô Bao Châu et  Eike Lau, 
    et des copoids arbitraires introduits simultanément par Ngô Bao Châu et  Eike Lau, 
    voir les  réferences dans le paragraphe \ref{intro-previous-works}).    
    
    Soit   $E$  une extension finie de 
    $\Ql$ contenant une racine carrée de $q$. On note $\mc O_{E}$  son anneau d'entiers.

             Soit $I$  un ensemble fini 
 et  $W$ une    représentation $E$-linéaire irréductible de $(\wh G)^{I}$. On écrit   $W=\boxtimes_{i\in I}W_{i}$ où $W_{i}$ est une représentation irréductible de  $\wh G$.  
  Le champ $\Cht_{N,I,W}^{(I)}$ classifiant les  $G$-chtoucas avec structure de  niveau $N$, a été étudié dans  \cite{var}.  
  Contrairement à \cite{var} nous imposons dans la définition suivante  qu'il soit réduit 
  (bien sûr cela ne change rien pour  la cohomologie étale). 
  
  \noindent{\bf Notation. }  Pour tout schéma  $S$ sur $\Fq$  et pour tout 
  $G$-torseur $\mc G$ sur  $ X\times S$ on note   $\ta \mc G=(\Id_{X}\times \Frob_{S})^{*}(\mc G)$. 
\begin{defi} On définit $\Cht_{N,I,W}^{(I)}$ comme le champ de Deligne-Mumford {\it réduit} dont 
 les points   sur un schéma  $S$ sur $\Fq$ classifient   
  \begin{itemize}
  \item des points  $(x_{i})_{i\in I}: S\to (X\sm N)^{I}$,  
\item un    $G$-torseur $\mc G$ sur  $ X\times S$, 
\item  un isomorphisme 
$$\phi :\restr{\mc G }{(X\times S)\sm(\bigcup_{i\in I }\Gamma_{x_i})}\isom \restr{\ta \mc G}{(X\times S)\sm(\bigcup_{i\in I }\Gamma_{x_i})}$$ 
 où   $\Gamma_{x_i}$ désigne le graphe de $x_{i}$, tel que la position relative en   $x_{i}$ soit bornée par le  copoids dominant de  $G$ correspondant au poids dominant $\omega_{i}$ de   $W_{i}$, 
 \item une trivialisation de   $(\mc G,\phi)$  sur  $N\times S$. 
 \end{itemize}
 \end{defi}
 
      Pour que la condition bornant les positions relatives soit définie sans ambiguïté on demande que $\Cht_{N,I,W}^{(I)}$  soit réduit et égal à l'adhérence de Zariski de son intersection avec l'ouvert où les $x_{i}$ sont deux à deux distincts. 
 
 Cette définition sera généralisée dans la \defiref{defi-Cht-I1-Ik} ci-dessous.

 On note  $\Cht_{ I,W}^{(I)}$ lorsque $N$ est vide et on remarque que  $\Cht_{N,I,W}^{(I)}$ est un $G(\mc O_{N})$-torseur 
 sur $\restr{\Cht_{ I,W}^{(I)}}{(X\sm N)^{I}}$.

 \begin{rem}  Les lecteurs connaissant le programme de Langlands géométrique noteront que $\Cht_{N,I,W}^{(I)}$ est l'intersection d'un champ de Hecke (considéré comme une correspondance entre $\Bun_{G,N}$ et lui-même) avec le graphe du morphisme de Frobenius de $\Bun_{G,N}$. 
 \end{rem}

Les  $x_{i}$ seront appelés les pattes du chtouca. On notera    $$\mf p_{N,I,W} ^{(I)}: \Cht_{N,I,W} ^{(I)}\to (X\sm N)^{I}$$ le morphisme correspondant. 
 
 Pour tout copoids dominant  $\mu$ de $G^{\mr{ad}}$ on note  
$\Cht_{N,I,W}^{(I),\leq\mu}$ l'ouvert de  $\Cht_{N,I,W}^{(I)}$ défini par la condition que le polygone de   Harder-Narasimhan   de $\mc G$
(ou plutôt, pour être précis, du  $G^{\mr{ad}}$-torseur associé) est $\leq \mu$.   
On fixe un réseau   $\Xi\subset Z(F)\backslash Z(\mb A)$. 
Alors $\Xi$ s'envoie dans $\Bun_{Z,N}(\Fq)$ qui agit sur $\Cht_{N,I,W}^{(I)}$ par torsion, et préserve les ouverts $\Cht_{N,I,W}^{(I),\leq\mu}$. 
On montrera que   $\Cht_{N,I,W}^{(I),\leq\mu}/\Xi$ est un champ de  Deligne-Mumford  de type fini. On note $\on{IC}_{\Cht_{N,I,W}^{(I),\leq \mu}/\Xi}$ le faisceau d'intersection  de  $\Cht_{N,I,W}^{(I),\leq \mu}/\Xi$ à coefficients dans  $E$, normalisé relativement à    $(X\sm N)^{I}$. 
On note  $$\mf p_{N,I,W} ^{(I),\leq\mu}: \Cht_{N,I,W}^{(I),\leq\mu}/\Xi\to (X\sm N)^{I}$$ le morphisme déduit de $\mf p_{N,I,W} ^{(I)}$ par restriction à 
   $\Cht_{N,I,W}^{(I),\leq\mu}$ et quotient par $\Xi$.

La définition suivante sera rendue plus canonique dans la \defiref{defi-HNIW-can}. 

\begin{defi} \label{defi-HNIW-naive}
On pose 
$$\mc H_{N,I,W}^{0,\leq\mu,E}=R^{0}\big(\mf p_{N,I,W} ^{(I),\leq\mu}\big)_{!}
 \Big(\on{IC}_{\Cht_{N,I,W}^{(I),\leq \mu}/\Xi}\Big).$$
  \end{defi}

  \noindent Dans le membre de droite  le  faisceau d'intersection est à coefficients dans $E$ et la cohomologie est prise au sens de \cite{laumon-moret-bailly,laszlo-olsson}.  En fait la cohomologie étale des schémas suffirait  
 (en effet, dès que le degré de $N$ est suffisamment grand en fonction de $\mu$, $\Cht_{N,I,W}^{(I),\leq \mu}/\Xi$ est un schéma de type fini).  
     
 Bien sûr $\mc H_{N,I,W}^{0,\leq\mu,E}$ dépend de $\Xi$ mais on omet $\Xi$ de la notation pour raccourcir un peu.

 Quand  $I$ est vide   et $W=\mbf 1$,  on a  
\begin{gather}\label{I-vide-intro} \restr{\varinjlim_{\mu}\mc H_{N,\emptyset,\mbf 1}^{0,\leq\mu,E}}{\Fqbar}=C_{c}(G(F)\backslash G(\mb A)/K_N \Xi,E)\end{gather}
  car  $\Cht_{N,\emptyset,\mbf 1}$ est le champ discret 
  $\Bun_{G,N}(\Fq)$, considéré comme un champ constant  sur $\Fq$, et  par ailleurs 
 $\Bun_{G,N}(\Fq)=G(F)\backslash G(\mb A)/K_N$. 
  On utilise ici  l'hypothèse que   $G$ est déployé, en général $\Bun_{G,N}(\Fq)$ est une réunion finie de quotients adéliques pour des formes intérieures de $G$, comme on le verra dans le  chapitre  \ref{para-non-deploye} (pous plus de détails on renvoie le lecteur aux  remarques \ref{quotient-adelique-deploye} et \ref{rem-une-classe-ker1}).   
   
 \begin{rem} \label{rem-cht-W1} Plus généralement pour tout $I$ et  $W=\mbf  1$, le champ  $\Cht_{N,I,{\mbf  1} }^{(I)}/\Xi$ est simplement le champ constant 
  $G(F)\backslash G(\mb A)/K_{N}\Xi$ sur $(X\sm N)^{I}$.  \end{rem}

On considère $\varinjlim_{\mu}\mc H_{N,I,W}^{0,\leq\mu,E}$ comme un système inductif de    $E$-faisceaux constructibles sur  $(X\sm N)^{I}$.
On va  introduire maintenant les actions des morphismes de  Frobenius partiels et des  opérateurs de Hecke sur ce système inductif (on notera que ces actions augmentent  $\mu$).  
Pour toute partie  $J\subset I$ on note   $$
 \Frob_{J}:(X\sm N)^{I}\to (X\sm N)^{I}$$  le  morphisme qui à  
 $(x_{i})_{i\in I}$ associe   $(x'_{i})_{i\in I}$ avec   $$x'_{i}=\Frob(x_{i})\text{ \  si \  }
  i\in J\text{   \ et  \  }x'_{i}=x_{i} \text{ \ sinon.}$$  
 Alors on possède 
 \begin{itemize}
 \item pour $\kappa$ assez grand en fonction de $W$,    pour tout   $i\in I$ et pour tout $\mu$, un  morphisme  
 \begin{gather}\label{intro-action-Frob-partiel}F_{\{i\}}:\Frob_{\{i\}}^{*}(\mc H_{N,I,W}^{0,\leq\mu,E})\to 
 \mc H_{N,I,W}^{0,\leq\mu+\kappa,E}\end{gather} 
 de faisceaux constructibles sur   $  (X\sm N)^{I}$, de sorte que les 
 $F_{\{i\}}$ commutent entre eux et que leur produit pour $i\in I$ est l'action naturelle du morphisme de Frobenius total de $(X\sm N)^{I}$ sur le faisceau $\mc H_{N,I,W}^{0,\leq\mu,E}$, 
  \item pour tout   $f\in C_{c}(K_{N}\backslash G(\mb A)/K_{N},E)$,  pour $\kappa$ assez grand en fonction de $W$ et de $f$, et pour tout $\mu$, un  morphisme 
 \begin{gather}\label{defi-Tf}T(f):\restr{\mc H_{N,I,W}^{0,\leq\mu,E}}{(X\sm \mf P)^{I}}\to 
 \restr{\mc H_{N,I,W}^{0,\leq\mu+\kappa,E}}{(X\sm \mf P)^{I}}\end{gather} de faisceaux constructibles sur   $  (X\sm \mf P)^{I}$ où $\mf P$ est un ensemble fini de places contenant  $|N|$ et en dehors duquel   $f$ est triviale.  
 \end{itemize}
 
 Les morphismes $T(f)$ sont appelés des ``opérateurs de Hecke'' bien que ce soient des  morphismes de faisceaux. 
Ils sont obtenus grâce à la construction, assez évidente, de correspondances de Hecke entre les champs de chtoucas. 
  On verra après la \propref{prop-coal-frob-cas-part-intro}   que   $T(f)$ peut être étendu naturellement en un  morphisme de faisceaux sur   $(X\sm N)^{I}$,  mais cela n'est pas trivial. 
  Bien sûr lorsque $I=\emptyset$ et $W=\mbf 1$, les morphismes $T(f)$ sont les opérateurs de Hecke habituels sur \eqref{I-vide-intro}. 
 
  Pour construire les actions \eqref{intro-action-Frob-partiel} des morphismes de Frobenius partiels, on a besoin d'une petite généralisation des champs 
 $\Cht_{N,I,W} ^{(I)}$ où l'on demande une factorisation de  $\phi$ en une suite de plusieurs modifications. Soit $(I_{1},...,I_{k})$ une partition (ordonnée)  de $I$. 
 Comme précédemment $W=\boxtimes _{i\in I}W_{i}$ est une représentation irréductible de $(\wh G)^{I}$.

 \begin{defi}\label{defi-Cht-I1-Ik}
 On définit     $\Cht_{N,I,W} ^{(I_{1},...,I_{k})}$ comme le champ de Deligne-Mumford {\it réduit }      dont les points sur 
 un schéma   $S$ sur  $\Fq$  classifient  les données    \begin{gather}\label{intro-donnee-chtouca}\big( (x_i)_{i\in I}, (\mc G_{0}, \psi_{0}) \xrightarrow{\phi_{1}}  (\mc G_{1}, \psi_{1}) \xrightarrow{\phi_{2}}
\cdots\xrightarrow{\phi_{k-1}}  (\mc G_{k-1}, \psi_{k-1}) \xrightarrow{ \phi_{k}}    (\ta{\mc G_{0}}, \ta \psi_{0})
\big)
\end{gather}
avec 
 \begin{itemize}
\item $x_i\in (X\sm N)(S)$ pour $i\in I$, 
\item pour $i\in \{0,...,k-1\}$, $(\mc G_{i}, \psi_{i})\in \Bun_{G,N}(S)$ (c'est-à-dire que  $\mc G_{i}$ est un   $G$-torseur sur  $X\times S$ et 
$\psi_{i} : \restr{\mc G_{i}}{N\times S} 
   \isom 
   \restr{G}{N\times S}$ est une trivialisation au-dessus de  $N\times S$) et on note  
   $(\mc G_{k}, \psi_{k})=(\ta{\mc G_{0}}, \ta \psi_{0})$
 \item  
pour   $j\in\{1,...,k\}$
 $$\phi_{j}:\restr{\mc G_{j-1}}{(X\times S)\sm(\bigcup_{i\in I_{j}}\Gamma_{x_i})}\isom \restr{\mc G_{j}}{(X\times S)\sm(\bigcup_{i\in I_{j}}\Gamma_{x_i})}$$ est un   isomorphisme  tel que la position relative de  $\mc G_{j-1}$ par rapport à  $\mc G_{j}$ en  $x_{i}$ (pour  $i\in I_{j}$) soit bornée par le  copoids dominant de  $G$ correspondant au  poids dominant de  $W_{i}$, 
 \item les $\phi_{j}$, qui induisent des isomorphismes sur $N\times S$, respectent les structures de niveau, c'est-à-dire que 
$\psi_{j}\circ \restr{\phi_{j}}{N\times S}=\psi_{j-1}$ pour tout  $j\in\{1,...,k\}$. 
 \end{itemize}
 \end{defi}
  On note $\Cht_{N,I}^{(I_{1},...,I_{k})}$ l'ind-champ obtenu en oubliant la condition sur les positions relatives. 
 
De plus   on note  
 $\Cht_{N,I,W}^{(I_{1},...,I_{k}),\leq\mu}$   l'ouvert de  $\Cht_{N,I,W}^{(I_{1},...,I_{k})}$ défini par la condition que  le polygone de Harder-Narasimhan  de  $\mc G_{0}$ est  $\leq\mu$. On note   $$\mf p_{N,I,W} ^{(I_{1},...,I_{k})}: \Cht_{N,I,W} ^{(I_{1},...,I_{k})}\to (X\sm N)^{I}$$ le morphisme   qui à un chtouca associe la famille de ses pattes. 
  
    \noindent{\bf Exemple.} Lorsque  $G=GL_r$,  $I=\{1,2\}$ et $W=\mr{St}\boxtimes \mr{St}^{*}$, les champs 
 $\Cht_{N,I,W}^{(\{1\}, \{2\})}$, {\it resp.} 
  $\Cht_{N,I,W}^{(\{2\}, \{1\})}$ 
 sont les champs de chtoucas à gauche, {\it resp.} à droite introduits par  Drinfeld \cite{Dr1}  (et utilisés dans  \cite{laurent-inventiones}), 
 et $x_{1}$ et  $x_{2}$ sont le zéro et le pôle.  
  
On  construit maintenant un     morphisme  lisse  
\eqref{lisse-chtouca-grass-intro} de $\Cht_{N,I,W}^{(I_{1},...,I_{k})}$ vers le  quotient d'une  strate  fermée  d'une  grassmannienne  affine  de Beilinson-Drinfeld  par un schéma   en groupes lisse.   
Les lecteurs familiers avec les variétés de Shimura peuvent considérer ce morphisme comme un ``modèle local'' à condition de noter
\begin{itemize}
\item que l'on est dans une situation de bonne réduction puisque les $x_{i}$ appartiennent à $X\sm N$, 
\item et que pourtant ce modèle local n'est pas lisse (sauf si tous les $I_{j}$ sont des singletons et tous  les copoids sont minuscules). 
\end{itemize}

\begin{defi}
La grassmannienne  affine  de Beilinson-Drinfeld est l'ind-schéma  $\mr{Gr}_{I }^{(I_{1},...,I_{k})}$ sur $X^{I}$ dont les $S$-points classifient la donnée de  
\begin{gather}\label{formule-rem-grassm-intro}\big((x_{i})_{i\in I}, \mc G_{0} \xrightarrow{\phi_{1}}  
\mc G_{1}\xrightarrow{\phi_{2}}
\cdots\xrightarrow{\phi_{k-1}}  \mc G_{k-1} \xrightarrow{ \phi_{k}}   \mc G_{k}\isor{\theta} G_{X\times S} \big)  \end{gather}
où les $\mc G_{i}$ sont des $G$-torseurs sur $X\times S$, $\phi_{i}$ est un isomorphisme sur $(X\times S)\sm(\bigcup_{i\in I_{j}}\Gamma_{x_i})$ et $\theta$ est une trivialisation de $\mc G_{k}$. 
 La strate  fermée   $\mr{Gr}_{I,W}^{(I_{1},...,I_{k})}$ est le  sous-schéma fermé réduit de 
  $\mr{Gr}_{I }^{(I_{1},...,I_{k})}$ défini par la condition que 
  la position relative de  $\mc G_{j-1}$ par rapport à  $\mc G_{j}$ en  $x_{i}$ (pour  $i\in I_{j}$) est  bornée par le  copoids dominant de  $G$ correspondant au  poids dominant $\omega_{i}$ de  $W_{i}$. Plus précisément 
  au-dessus de l'ouvert $\mc U$ de $X^{I}$ où les $x_{i}$ sont deux à deux distincts, $\mr{Gr}_{I }^{(I_{1},...,I_{k})}$ est un produit de grassmanniennes affines usuelles et 
  \begin{itemize}
  \item on définit 
  la restriction de $\mr{Gr}_{I,W}^{(I_{1},...,I_{k})}$ au-dessus de $\mc U$  comme le produit des strates fermées habituelles (notées $\ov {\mr{Gr}_{\omega_{i}}}$ dans  \cite{mv,brav-gaitsgory} où $\omega_{i}$ 
  est le plus haut poids de $W_{i}$), 
  \item puis 
    on   définit  $\mr{Gr}_{I,W}^{(I_{1},...,I_{k})}$ comme l'adhérence de Zariski (dans $\mr{Gr}_{I }^{(I_{1},...,I_{k})}$) de 
  sa restriction au-dessus de $\mc U$. 
  \end{itemize}
  \end{defi}
  
   D'après Beauville-Laszlo \cite{BL} (voir aussi le premier chapitre  
   pour des références complémentaires  dans \cite{hitchin}), 
   $\mr{Gr}_{I }^{(I_{1},...,I_{k})}$ peut aussi être défini comme l'ind-schéma dont les $S$-points classifient 
       \begin{gather}\label{formule-rem-grassm-intro-loc}\big((x_{i})_{i\in I}, \mc G_{0} \xrightarrow{\phi_{1}}  
\mc G_{1}\xrightarrow{\phi_{2}}
\cdots\xrightarrow{\phi_{k-1}}  \mc G_{k-1} \xrightarrow{ \phi_{k}}   \mc G_{k}\isor{\theta} G_{\Gamma_{\sum \infty x_i}} \big)  \end{gather}
 où les $\mc G_{i}$ sont des $G$-torseurs  sur le voisinage formel 
  $\Gamma_{\sum \infty x_i}$ de la réunion des graphes des $x_{i}$ dans
 $  X\times S$, $\phi_{i}$ est un isomorphisme sur $\Gamma_{\sum \infty x_i}\sm(\bigcup_{i\in I_{j}}\Gamma_{x_i})$ et $\theta$ est une trivialisation de $\mc G_{k}$. 
La restriction à la Weil $G_{\sum \infty x_i}$ de $G$ de $\Gamma_{\sum \infty x_i}$ à $S$ agit donc sur  $\mr{Gr}_{I }^{(I_{1},...,I_{k})}$  et $\mr{Gr}_{I,W}^{(I_{1},...,I_{k})}$ par changement de la trivialisation $\theta$. 

On a un morphisme naturel 
    \begin{gather}\label{lisse-chtouca-grass-intro0}
    \Cht_{N,I,W}^{(I_{1},...,I_{k})}\to \mr{Gr}_{I,W}^{(I_{1},...,I_{k})}/
  G_{\sum \infty  x_i}\end{gather} 
  qui associe à un chtouca \eqref{intro-donnee-chtouca} le 
  $G_{\sum \infty x_{i}}$-torseur $\restr{\mc G_{k}}{\Gamma_{\sum \infty x_{i}}}$
  et, pour toute trivialisation $\theta$ de celui-ci, le  point de $\mr{Gr}_{I,W}^{(I_{1},...,I_{k})}$ égal à  \eqref{formule-rem-grassm-intro-loc}.

 \begin{rem} \label{position-relative} La meilleure fa\c con d'énoncer la condition sur les positions relatives dans la \defiref{defi-Cht-I1-Ik} est de {\it définir}   $\Cht_{N,I,W}^{(I_{1},...,I_{k})}$ comme l'image inverse de 
  $\mr{Gr}_{I,W}^{(I_{1},...,I_{k})}/
  G_{\sum \infty x_i}$ par le morphisme   
  $\Cht_{N,I}^{(I_{1},...,I_{k})}\to \mr{Gr}_{I}^{(I_{1},...,I_{k})}/
  G_{\sum \infty x_i}$ construit  comme \eqref{lisse-chtouca-grass-intro0}. 
    \end{rem}

Pour $(n_{i})_{i\in I}\in \N^{I}$ on note $\Gamma_{\sum n_{i} x_i}$ le sous-schéma fermé de $X\times S$ associé au diviseur de Cartier $\sum n_{i} x_i$ qui est  effectif et relatif  sur $S$. On note $G_{\sum n_{i} x_i}$ le schéma en groupes lisse sur $S$ obtenu par restriction à la Weil de $G$ de $\Gamma_{\sum n_{i} x_i}$ à $S$. 
 Alors si les entiers $n_{i}$ sont assez grands en fonction de $W$, l'action de $G_{\sum \infty x_i}$ sur $\mr{Gr}_{I,W}^{(I_{1},...,I_{k})}$ se factorise par 
  $G_{\sum n_{i} x_i}$. 
 Le  morphisme \eqref{lisse-chtouca-grass-intro0} fournit donc un morphisme 
    \begin{gather}\label{lisse-chtouca-grass-intro}
    \Cht_{N,I,W}^{(I_{1},...,I_{k})}\to \mr{Gr}_{I,W}^{(I_{1},...,I_{k})}/
  G_{\sum n_{i} x_i}\end{gather} 
  (qui associe à un chtouca \eqref{intro-donnee-chtouca} le 
  $G_{\sum n_{i}x_{i}}$-torseur $\restr{\mc G_{k}}{\Gamma_{\sum n_{i}x_{i}}}$
  et, pour toute trivialisation $\lambda$ de celui-ci, le  point de $\mr{Gr}_{I,W}^{(I_{1},...,I_{k})}$ égal à  \eqref{formule-rem-grassm-intro-loc} pour toute trivialisation $\theta$ de $\restr{\mc G_{k}}{\Gamma_{\sum \infty x_i}}$ prolongeant  $\lambda$ de $\Gamma_{\sum n_{i} x_i}$ à 
  $\Gamma_{\sum \infty   x_i}$).

 On montrera  dans la \propref{lissite-Cht-Grass}  que le morphisme 
  \eqref{lisse-chtouca-grass-intro} 
    est  lisse de  dimension $\dim G_{\sum n_{i}x_{i}}=(\sum_{i\in I} n_{i})\dim G$ (l'idée est la suivante : il suffit de le montrer dans le  cas où $N$ est vide  
et alors cela résulte  facilement 
 du fait que le morphisme de Frobenius de $\Bun_{G}$   a une dérivée nulle).  
    
 On en déduit  que le  morphisme d'oubli des modifications intermédiaires 
 \begin{gather}
 \label{oubli-Cht}
 \Cht_{N,I,W}^{(I_{1},...,I_{k})}\to \Cht_{N,I,W}^{(I)} \\ \nonumber 
 \text{   qui envoie \eqref{intro-donnee-chtouca} sur  } 
\big( (x_i)_{i\in I}, (\mc G_{0}, \psi_{0}) \xrightarrow{\phi_{k} \cdots \phi_{1}}      (\ta{\mc G_{0}}, \ta \psi_{0})
\big), \end{gather}
 est petit. En effet il est connu que le morphisme analogue 
\begin{gather}\label{mor-Gr-oubli}\mr{Gr}_{I,W}^{(I_{1},...,I_{k})}\to \mr{Gr}_{I,W}^{(I)}  \text{    qui envoie \eqref{formule-rem-grassm-intro-loc} sur }
\big( (x_i)_{i\in I},  \mc G_{0}  \xrightarrow{\phi_{k} \cdots \phi_{1}}       \mc G_{k}\isor{\theta} G_{\Gamma_{\sum \infty x_i}} 
\big)  \end{gather}   
est petit, et d'ailleurs cela joue un rôle essentiel dans   \cite{mv}. De plus l'image inverse de $\Cht_{N,I,W}^{(I),\leq\mu}$ 
par \eqref{oubli-Cht} est exactement $\Cht_{N,I,W}^{(I_{1},...,I_{k}),\leq\mu}$ 
puisque les troncatures ont été définies à l'aide du polygône de Harder-Narasimhan de $\mc G_{0}$. 
 On a donc  $$\mc H_{N,I,W}^{0,\leq\mu,E}=
   R^{0}(\mf p_{N,I,W} ^{(I_{1},...,I_{k}),\leq\mu})_{!}\Big(\on{IC}_{\Cht_{N,I,W}^{(I_{1},...,I_{k}),\leq\mu}/\Xi}\Big)$$ pour {\it toute} partition $(I_{1},...,I_{k})$ de $I$ (alors que la \defiref{defi-HNIW-naive} utilisait la partition grossière $(I)$).  

Le  morphisme de  Frobenius partiel 
$$\on {Fr}_{I_{1}} ^{(I_{1},...,I_{k})}: \Cht_{N,I,W} ^{(I_{1},...,I_{k})}\to \Cht_{N,I,W} ^{(I_{2},...,I_{k},I_{1})},$$ 
 défini par 
\begin{gather}\on {Fr}_{I_{1}} ^{(I_{1},...,I_{k})}\big( (x_i)_{i\in I}, (\mc G_{0}, \psi_{0}) \xrightarrow{\phi_{1}}  (\mc G_{1}, \psi_{1}) \xrightarrow{\phi_{2}}
\cdots\xrightarrow{\phi_{k-1}}  (\mc G_{k-1}, \psi_{k-1}) \xrightarrow{ \phi_{k}}    (\ta{\mc G_{0}}, \ta \psi_{0})
\big)\nonumber \\
= \nonumber
\big( \Frob_{I_{1}}\big((x_i)_{i\in I}\big), (\mc G_{1}, \psi_{1}) \xrightarrow{\phi_{2}}  (\mc G_{2}, \psi_{2}) \xrightarrow{\phi_{3}}
\cdots \xrightarrow{ \phi_{k}}    (\ta{\mc G_{0}}, \ta \psi_{0}) \xrightarrow{\ta  \phi_{1}   } (\ta{\mc G_{1}}, \ta \psi_{1})
\big)  
\end{gather}
est au-dessus du morphisme  $\Frob_{I_{1}}:(X\sm N)^{I}\to (X\sm N)^{I}$. 
Comme   $\on {Fr}_{I_{1}} ^{(I_{1},...,I_{k})}$
est un  homéomorphisme local totalement radiciel, on a un isomorphisme canonique   $$\Big(\on {Fr}_{I_{1}} ^{(I_{1},...,I_{k})}\Big)^{*}\Big(\on{IC}_{\Cht_{N,I,W} ^{(I_{2},...,I_{k},I_{1})}}\Big)=\on{IC}_{ \Cht_{N,I,W} ^{(I_{1},...,I_{k})}} .   $$
L'isomorphisme de changement de base propre fournit alors un morphisme $$F_{I_{1}}:\Frob_{I_{1}}^{*}(\mc H_{N,I,W}^{0,\leq\mu,E})\to 
 \mc H_{N,I,W}^{0,\leq\mu+\kappa,E}$$ pour $\kappa$ assez grand
 (en effet, dans les notations ci-dessus, si le polygône de Harder-Narasimhan de $\mc G_{1}$ est $\leq \mu$, comme la modification entre $\mc G_{0}$ et $\mc G_{1}$ est bornée en fonction de $W$, le polygône de Harder-Narasimhan de $\mc G_{0}$ est $\leq \mu+\kappa$ où $\kappa$ dépend de $W$). En prenant  n'importe quelle partition $(I_{1},...,I_{k})$  telle que  $I_{1}=\{i\}$ on obtient $F_{\{i\}}$ dans  \eqref{intro-action-Frob-partiel}.

 Pour l'instant on a défini $\mc H_{N,I,W}^{0,\leq\mu,E}$ pour les classes d'isomorphismes de représentations irréductibles $W$ de $(\wh G)^{I}$. On va raffiner  cette construction en celle  {\it canonique}  d'un {\it  foncteur} $E$-linéaire
 \begin{gather}\label{fonc-W-HIW-intro}W\mapsto \mc H_{N,I,W}^{0,\leq\mu,E} \end{gather} 
 de la catégorie des représentations $E$-linéaires de dimension finie de $(\wh G)^{I}$ vers la catégorie des $E$-faisceaux constructibles sur $(X\sm N)^{I}$. En particulier pour tout    morphisme $u:W\to W'$ de  représentations $E$-linéaires de $(\wh G)^{I}$  on notera 
 $$\mc H(u): \mc H_{N,I,W}^{0,\leq\mu,E} \to \mc H_{N,I,W'}^{0,\leq\mu,E}$$  
le morphisme de faisceaux constructibles  associé. 

 Le foncteur \eqref{fonc-W-HIW-intro} sera compatible à la coalescence des pattes, au sens suivant. Dans tout cet article nous appelons coalescence la situation où des pattes fusionnent entre elles. 
 Nous aurions pu employer  le mot fusion plutôt que coalescence
 mais nous avons préféré utiliser le mot coalescence pour les pattes
 (qui ne sont que des points sur la courbe) en gardant le mot fusion 
 pour le produit de fusion (qui intervient dans l'équivalence de Satake géométrique et concerne les faisceaux pervers sur les grassmanniennes affines de Beilinson-Drinfeld). 
    Soit  $\zeta: I \to J$ une application.  On note $W^{\zeta}$ la représentation de  $\wh G^{J}$ qui est la composée de la  représentation $W$ avec le morphisme diagonal
 $$ \wh G^{J}\to \wh G^{I}, (g_{j})_{j\in J}\mapsto (g_{\zeta(i)})_{i\in I}.$$
   On note 
   \begin{gather}\label{morph-giad-X-intro}\Delta_{\zeta}: X^{J}\to X^{I},(x_{j})_{j\in J}\mapsto (x_{\zeta(i)})_{i\in I}\end{gather} le   morphisme diagonal. 
On va construire, d'après  \cite{var} et \cite{brav-var}   un isomorphisme de faisceaux constructibles sur $(X\sm N)^{J}$, dit de  {\it coalescence}:  \begin{gather}\label{intro-isom-coalescence}\chi_{\zeta}: \Delta_{\zeta}^{*}(\mc H_{N,I,W}^{0,\leq\mu,E})\isom 
 \mc H^{0,\leq\mu,E}_{N,J,W^{\zeta}}.\end{gather} 
 Cet isomorphisme sera {\it canonique}, au sens où ce sera un isomorphisme de foncteurs en $W$, compatible à la composition de $\zeta$. 
    
     On explique  maintenant  la construction du foncteur \eqref{fonc-W-HIW-intro} et de l'isomorphisme de coalescence \eqref{intro-isom-coalescence}. 
 
 Lorsque $W$ n'est pas irréductible on 
 note $\mr{Gr}_{I,W }^{(I_{1},...,I_{k})}$  la réunion des $\mr{Gr}_{I,V }^{(I_{1},...,I_{k})}\subset \mr{Gr}_{I }^{(I_{1},...,I_{k})}$ pour $V$ constituant irréductible de $W$. On fait de même avec 
  $\Cht_{N,I,W}^{(I_{1},...,I_{k})}$.

   On rappelle maintenant l'équivalence de Satake géométrique, due à Lusztig, Drinfeld, Ginzburg et Mirkovic-Vilonen. 
   Pour des références on cite  \cite{lusztig-satake,ginzburg,hitchin,mv,ga-iwahori,ga-de-jong,richarz,zhu}. On va utiliser  ici     la  forme  expliquée par Gaitsgory dans \cite{ga-de-jong}.  Habituellement l'équivalence de Satake géométrique s'exprime de la manière suivante. Pour tout corps $k$ algébriquement clos de caractéristique première à $\ell$, 
    la catégorie des faisceaux pervers $G(k[[z]])$-équivariants sur la grassmannienne affine $G(k((z)))/G(k[[z]])$ est munie d'une structure tensorielle par le produit de fusion (ou de convolution) et d'un foncteur fibre   donné par la cohomologie totale. Comme cela est expliqué dans \cite{mv, hitchin} et rappelé plus en détails ci-dessous, on modifie un peu la règle des signes dans la contrainte de commutativité pour que le foncteur fibre soit tensoriel à valeurs dans la catégorie des espaces vectoriels (et non pas des super-espaces vectoriels). Pour être plus canonique il faut introduire en plus  une torsion à la Tate, c'est-à-dire tensoriser par $E(\frac{i}{2})$ la partie de   degré cohomologique $i$ pour tout $i\in \Z$. 
   Cette catégorie  est donc équivalente à la catégorie des représentations de dimension finie 
   du groupe algébrique des automorphismes du foncteur fibre, qui s'avère être isomorphe à  
   $\wh G$ (muni d'un épinglage canonique). 
   De plus ces faisceaux pervers sont naturellement équivariants par le groupe des automorphismes de $k[[z]]$. Ceci permet de remplacer 
   $\on{Spf}(k[[z]])$ par un disque formel arbitraire, et en particulier un disque  formel variant sur une courbe. 
   
   Ici on   utilise  seulement un sens de l'équivalence, à savoir le foncteur de la catégorie des représentations de dimension finie de $\wh G$ vers la catégorie des faisceaux pervers 
   $G(k[[z]])$-équivariants  sur la grassmiannenne affine. En revanche  on l'énonce en utilisant la grassmannienne affine de Beilinson-Drinfeld. 
   Dans les notations du théorème ci-dessous, 
    \eqref{fonct-W-thm-Satake-geom} est un foncteur en $W$, et cela est plus fort que si on avait seulement  la propriété d) du théorème, qui   détermine de fa\c con non canonique  
    la classe d'isomorphisme de 
      $\mc S_{I,W,E }^{(I_{1},...,I_{k})}$ pour chaque classe de représentation irréductible $W$ de $(\wh G)^{I}$. Comme on l'a rappelé précédemment, la construction de ce foncteur repose sur le foncteur fibre  donné par la cohomologie totale. 
  Par exemple, dans les notations du théorème ci-dessous, $W$ est {\it égal} à la cohomologie totale  de $\mc S_{I,W,E }^{(I_{1},...,I_{k})}$ dans les fibres de $\mr{Gr}_{I }^{(I_{1},...,I_{k})}$ au-dessus de $X^{I}$ (avec les avertissements précédents concernant la règle des signes et les torsions à la Tate).
  On note  cependant que l'on n'utilise pas cette propriété en tant que telle, seulement à travers le fait qu'elle fournit le foncteur   \eqref{fonct-W-thm-Satake-geom}, et donc la canonicité de nos constructions.

 \begin{thm}  \label{satake-geom-thm} (un sens de l'équivalence de Satake géométrique \cite{hitchin, mv, ga-de-jong}, voir le \thmref{thm-geom-satake} ci-dessous pour les détails). On possède pour tout ensemble fini $I$ et pour toute partition $(I_{1},...,I_{k})$ un   foncteur   $E$-linéaire   
   \begin{gather}\label{fonct-W-thm-Satake-geom}W\mapsto \mc S_{I,W,E }^{(I_{1},...,I_{k})}\end{gather}
de la catégorie des  représentations $E$-linéaires de dimension finie de $(\wh G)^{I}$ vers la  catégorie des    $E$-faisceaux pervers  
$G_{\sum \infty x_{i}}$-équivariants sur $\mr{Gr}_{I }^{(I_{1},...,I_{k})}$ (pour la normalisation perverse  relative à $X^{I}$). De plus $\mc S_{I,W,E}^{(I_{1},...,I_{k})}$ est supporté par 
$\mr{Gr}_{I,W,E}^{(I_{1},...,I_{k})}$ et universellement localement acyclique relativement à $X^{I}$. Ces foncteurs vérifient les propriétés suivantes.  

\begin{itemize}
\item [] a) Compatibilité aux morphismes d'oubli (des modifications intermédiaires)  : $\mc S_{I,W,E}^{(I)}$ est canoniquement isomorphe à l'image directe de 
$\mc S_{I,W,E}^{(I_{1},...,I_{k})}$ par le morphisme d'oubli  
$\mr{Gr}_{I}^{(I_{1},...,I_{k})}\to \mr{Gr}_{I}^{(I)}$ (défini dans \eqref{mor-Gr-oubli}).   
\item [] b)  Compatibilité à la convolution : si $W=\boxtimes_{j\in \{1,...,k\}} W_{j}$ où $W_{j}$ est une  représentation de $(\wh G)^{I_{j}}$, 
  $\mc S_{I,W,E}^{(I_{1},...,I_{k})}$  est canoniquement isomorphe à l'image inverse de $\boxtimes _{j\in \{1,...,k\}} \mc S_{I_{j},W_{j},E}^{(I_{j})}$ par le morphisme 
  \begin{align*}
  \mr{Gr}_{I}^{(I_{1},...,I_{k})}/G_{\sum_{i\in I} \infty x_{i}} &\to 
  \prod_{j=1}^{k}\Big( \mr{Gr}_{I_{j}}^{(I_{j})}/G_{\sum_{i\in I_{j}} \infty x_{i}}\Big) \\
  (\mc G_{0}\to \mc G_{1} \to  \cdots \to  \mc G_{k})& 
  \mapsto \Big(\big(\restr{\mc G_{j-1}}{\Gamma_{\sum_{i\in I_{j}} \infty x_{i}}}\to \restr{\mc G_{j}}{\Gamma_{\sum_{i\in I_{j}} \infty x_{i}}} \big) \Big)_{j=1,...,k}
  \end{align*}
  où les $\mc G_{i}$ sont des $G$-torseurs sur 
  $\Gamma_{\sum_{i\in I} \infty x_{i}}$. 
   \item [] c)  Compatibilité à la fusion: soient  $I,J$ des ensembles finis et  $\zeta: I\to J$ une application. Soit $(J_{1},...,J_{k})$ une partition de $J$. Son image inverse 
   $(\zeta^{-1}(J_{1}), ..., \zeta^{-1}(J_{k}))$ est une partition de $I$. 
On note  $$\Delta_{\zeta} : X^{J}\to X^{I},  \ \ (x_{j})_{j\in J}\mapsto (x_{\zeta(i)})_{i\in I}$$ le   morphisme diagonal associé à 
  $\zeta$. On note encore 
  $\Delta_{\zeta}$    l'inclusion  
 $$\mr{Gr}_{J}^{(J_{1},...,J_{k})} =\mr{Gr}_{I}^{(\zeta^{-1}(J_{1}), ..., \zeta^{-1}(J_{k}))}\times _{X^{I}}X^{J}\hookrightarrow \mr{Gr}_{I}^{(\zeta^{-1}(J_{1}), ..., \zeta^{-1}(J_{k}))}.$$
   Soit  $W$ une   représentation  $E$-linéaire de dimension finie   de  $\wh G^{I}$. On note  $W^{\zeta}$ la   représentation de $\wh G^{J}$ qui est la composée de la   représentation $W$ avec le    morphisme diagonal $$\wh G^{J}\to \wh G^{I},  \ \ (g_{j})_{j\in J}\mapsto (g_{\zeta(i)})_{i\in I}  .$$   On a alors un  isomorphisme canonique 
 \begin{gather}\label{coalescence-Gr-intro} \Delta_{\zeta}^{*}\Big( \mc S_{I,W,E}^{(\zeta^{-1}(J_{1}), ..., \zeta^{-1}(J_{k}))} \Big)\simeq 
 \mc S_{J,W^{\zeta},E}^{(J_{1},...,J_{k})}\end{gather}  
 qui est fonctoriel en  $W$ et compatible avec la composition pour $\zeta$. 
  \item [] d) 
 Lorsque  
$W$ est irréductible, le faisceau pervers  $\mc S_{I,W,E}^{(I_{1},...,I_{k})}$ sur 
 $\mr{Gr}_{I,W}^{(I_{1},...,I_{k})}$ est isomorphe au  faisceau d'intersection   (avec la  normalisation perverse relative à  $X^{I}$). 
 \end{itemize}  
\end{thm}

 Dans le théorème précédent  $\mc S_{I,W,E}^{(I_{1},...,I_{k})}$ est supporté par 
$\mr{Gr}_{I,W}^{(I_{1},...,I_{k})}$ et on peut donc le considérer comme un faisceau pervers 
(à un décalage près)   sur  $\mr{Gr}_{I,W}^{(I_{1},...,I_{k})}/G_{\sum n_{i}x_{i}}$  (avec les entiers $n_{i}$   assez grands).   
  
  Comme on le dira plus en détails dans la  \remref{contrainte-commutativite}
 la contrainte de commutativité est définie par la convention  {\it modifiée},  introduite dans la discussion qui précède la   proposition 6.3 de \cite{mv}. La contrainte de commutativité modifiée 
 est celle que l'on obtiendrait naturellement si tous les faisceaux d'intersection 
   $\mc S_{I,W,E}^{(I_{1},...,I_{k})}$ pour $W$ irréductible étaient normalisés, du point de vue du degré cohomologique modulo $2$, 
   pour être en degré cohomologique pair au point générique de $\mr{Gr}_{I,W}^{(I_{1},...,I_{k})}/G_{\sum n_{i}x_{i}}$, et alors le lemme 3.9 de \cite{mv} 
   affirme que   leur cohomologie totale serait supportée en degrés cohomologiques pairs. Donc avec cette  contrainte de commutativité modifiée 
  le foncteur fibre donné par la cohomologie totale sur  la catégorie tensorielle des  faisceaux pervers  $G(\mc O)$-équivariants   sur la   grassmannienne affine  est un foncteur tensoriel à valeurs dans la catégorie tensorielle des espaces vectoriels (et non pas des super-espaces vectoriels),  d'où  l'équivalence avec  la catégorie tensorielle des représentations de  dimension finie de $\wh G$ (défini comme le groupe d'automorphismes du foncteur fibre).  
   
   \begin{rem}\label{rem-Satake-Gad}
    Dans le théorème précédent $Z_{\sum_{i\in I} \infty x_{i}}\subset 
G_{\sum_{i\in I} \infty x_{i}}$ agit trivialement sur 
 $ \mr{Gr}_{I}^{(I_{1},...,I_{k})}$ et donc (par d)) sur tous les faisceaux $\mc S_{I,W }^{(I_{1},...,I_{k})}$. On note  $G^{\mr{ad}}=G/Z$. On peut donc considérer $\mc S_{I,W}^{(I_{1},...,I_{k})}$  comme un 
 faisceau pervers  (à un décalage près)
$G^{\mr{ad}}_{\sum \infty x_{i}}$-équivariant sur $\mr{Gr}_{I }^{(I_{1},...,I_{k})}$
 ou si on préfère comme  
 un faisceau pervers 
(à un décalage près)   sur  $\mr{Gr}_{I,W}^{(I_{1},...,I_{k})}/G^{\mr{ad}}_{\sum n_{i}x_{i}}$  (avec les entiers $n_{i}$   assez grands).   
    \end{rem}

Voici la construction du  foncteur \eqref{fonc-W-HIW-intro}. 
Le morphisme \eqref{lisse-chtouca-grass-intro} ne se factorise pas par le quotient par $\Xi$ (comme me l'a fait remarquer un rapporteur anonyme), mais c'est  le cas de sa composée avec le morphisme d'oubli $ \mr{Gr}_{I,W}^{(I_{1},...,I_{k})}/
  G_{\sum n_{i} x_i}\to  \mr{Gr}_{I,W}^{(I_{1},...,I_{k})}/
  G^{\mr{ad}}_{\sum n_{i} x_i}$. Autrement dit on possède un morphisme 
  \begin{gather}\label{lisse-chtouca-grass-intro-ad}
    \Cht_{N,I,W}^{(I_{1},...,I_{k})}/\Xi\to \mr{Gr}_{I,W}^{(I_{1},...,I_{k})}/
  G^{\mr{ad}}_{\sum n_{i} x_i} \end{gather} 
et d'après la \remref{rem-Satake-Gad}, $\mc S_{I,W}^{(I_{1},...,I_{k})}$  est un faisceau pervers (à un décalage près) sur l'espace d'arrivée.

\begin{defi}\label{defi-HNIW-can}
On 
 définit le  faisceau pervers (avec la normalisation relative à $(X\sm N)^{I}$) 
 $\mc F_{N,I,W,\Xi,E}^{(I_{1},...,I_{k})}$ sur 
 $\Cht_{N,I,W}^{(I_{1},...,I_{k}) }/\Xi$ 
 comme l'image inverse de 
 $\mc S_{I,W,E}^{(I_{1},...,I_{k})}$ par le morphisme 
   \eqref{lisse-chtouca-grass-intro-ad}.    
   On définit alors le foncteur \eqref{fonc-W-HIW-intro} en posant   
  \begin{gather}\label{defi-can-H}\mc H_{N,I,W}^{0,\leq\mu,E}=
   R^{0}(\mf p_{N,I,W} ^{(I_{1},...,I_{k}),\leq\mu})_{!}\Big(\restr{\mc F_{N,I,W,\Xi,E}^{(I_{1},...,I_{k})}}{\Cht_{N,I,W}^{(I_{1},...,I_{k}), \leq\mu}/\Xi}\Big)\end{gather} pour toute partition $(I_{1},...,I_{k})$ de $I$. 
   \end{defi}
   Grâce au a) du théorème précédent la définition  \eqref{defi-can-H} ne dépend pas du choix de la partition  $(I_{1},...,I_{k})$. 
 
 Lorsque $W$ est irréductible, la lissité du morphisme \eqref{lisse-chtouca-grass-intro}, donc celle du morphisme \eqref{lisse-chtouca-grass-intro-ad},  et le calcul de sa dimension impliquent  que 
   $\mc F_{N,I,W,\Xi,E}^{(I_{1},...,I_{k})}$ est isomorphe au faisceau  d'intersection de 
   $\Cht_{N,I,W}^{(I_{1},...,I_{k})}/\Xi$ (avec la normalisation perverse relative à $(X\sm N)^{I}$).   Donc la définition précédente est cohérente avec 
   la définition \ref{defi-HNIW-naive} (et la raffine en la rendant plus canonique).

 L'action des opérateurs de Hecke et des morphismes  de Frobenius partiels peut se reformuler à l'aide de la définition \eqref{defi-can-H}. 
 La canonicité de la définition \eqref{defi-can-H} est surtout cruciale pour la construction des isomorphismes de coalescence 
 \eqref{intro-isom-coalescence}, que nous allons expliquer maintenant, 
 car les espaces de départ et d'arrivée ne sont pas les mêmes. 
     
 \begin{defi}  L'isomorphisme {\it canonique} \eqref{intro-isom-coalescence} est défini (grâce au théorème de changement de base propre) par l'isomorphisme {\it canonique} 
   entre 
   $\mc F_{N,J,W^{\zeta},\Xi,E}^{(J_{1},...,J_{k})}$   et l'image inverse de 
 $\mc F_{N,I,W,\Xi,E}^{(\zeta^{-1}(J_{1}), ..., \zeta^{-1}(J_{k}))}$ par (le quotient par $\Xi$ de)  l'inclusion 
   $$\Cht_{N,J }^{(J_{1},...,J_{k})}
     =\Cht_{N,I }^{(\zeta^{-1}(J_{1}), ..., \zeta^{-1}(J_{k}))}\times _{(X\sm N)^{I}}(X\sm N)^{J}\hookrightarrow \Cht_{N,I }^{(\zeta^{-1}(J_{1}), ..., \zeta^{-1}(J_{k}))} 
     $$
     qui provient de l'isomorphisme \eqref{coalescence-Gr-intro} dans le  c) du \thmref{satake-geom-thm}. 
     \end{defi}
     La définition précédente est  indépendante du choix de la partition $(J_{1},...,J_{k})$. 
    
    \begin{rem}
    \label{rem-coalescence-frob-compat} Le fait d'avoir pris une partition $(J_{1},...,J_{k})$ arbitraire permet de   montrer la compatibilité entre l'isomorphisme de coalescence \eqref{intro-isom-coalescence} et les morphismes de Frobenius partiels, à savoir
    que pour tout $j\in J$, $\Delta_{\zeta}^{*}(F_{\zeta^{-1}(\{j\})})$ et 
    $F_{\{j\}}$ se correspondent par l'isomorphisme $\chi_{\zeta}$ 
     de \eqref{intro-isom-coalescence}. \end{rem}

            \subsection{Morphismes de création et d'annihilation}
            \label{subsection-crea-annihil-intro}

    Dans ce paragraphe notre but est d'utiliser les   isomorphismes de coalescence  \eqref{intro-isom-coalescence} pour construire des morphismes de création et d'annihilation, puis d'exprimer les   opérateurs de Hecke en les  places de  $X\sm N$ comme la composée   
 \begin{itemize}
 \item d'un morphisme de création,
  \item de  l'action d'un  morphisme de Frobenius partiel,
  \item d'un  morphisme d'annihilation,
  \end{itemize}
et d'utiliser cela pour étendre les  opérateurs de Hecke \eqref{defi-Tf} en des morphismes de faisceaux sur   $(X\sm N)^{I}$ tout entier  et pour obtenir les relations d'Eichler-Shimura.

  Soient   $I$ et $J$ des ensembles finis. 
  On va définir maintenant les morphismes de création et d'annihilation, dont l'idée est la suivante. Les pattes indexées par $I$ restent inchangées et on crée (ou on annihile) les pattes indexées par $J$ en un même point de la courbe (indexé par un ensemble à un élément, que l'on note $\{0\}$).

  On a des applications évidentes 
   $$\zeta_{J} :J \to \{0\},  \ 
  \zeta_{J}^{I}=(\Id_{I},\zeta_{J}): I\cup J\to I\cup \{0\} \text{ \   et   \ } \zeta_{\emptyset}^{I}=(\Id_{I},\zeta_{\emptyset}): I \to I \cup \{0\}.$$ 
  Soient $W$ et $U$ des  représentations  $E$-linéaires de dimension finie de 
$(\wh G)^{I}$ et   $(\wh G)^{J}$ respectivement. 
On rappelle que   $U^{\zeta_{J}}$ est 
la représentation de $\wh G$ obtenue en restreignant 
  $U$ à la diagonale   $\wh G\subset (\wh G)^{J}$. 
Soient  $x\in U$ et $\xi\in U^{*}$ invariants sous  l'action diagonale de $\wh G$.  
Alors $W\boxtimes U$ est une représentation de  $(\wh G)^{I\cup J}$ et 
 $W\boxtimes U^{\zeta_{J}}$ et  $W\boxtimes \mbf 1$ sont des  représentations de  $(\wh G)^{I\cup \{0\}}$ reliées par les   morphismes 
 $$W\boxtimes \mbf 1\xrightarrow{\Id_{W}\boxtimes x}W\boxtimes U^{\zeta_{J}}\text{ \  et \ }  W\boxtimes U^{\zeta_{J}}\xrightarrow{\Id_{W}\boxtimes \xi}W\boxtimes \mbf 1.$$ 
 
 On note  $\Delta:X\to X^{J}$ le   morphisme diagonal et on désigne par $E_{X\sm N}$   le faisceau constant sur  $X\sm N$. 
         \begin{defi}
 On définit  le  morphisme de création
  $
\mc C_{x}^{\sharp}$ comme la  composée   \begin{gather*}
  \mc H _{N,I,W}^{0,\leq\mu,E} \boxtimes E_{(X\sm N)} 
  \isor{\chi_{\zeta_{\emptyset}^{I}}}
 \mc H _{N,I\cup\{0\},W\boxtimes \mbf 1}^{0,\leq\mu,E}  \\
 \xrightarrow{\mc H(\Id_{W}\boxtimes x) }
 \mc H _{N,I\cup\{0\},W\boxtimes U^{\zeta_{J}}}^{0,\leq\mu,E} 
 \isor{ \chi_{\zeta_{J}^{I}}^{-1} }
 \restr{ \mc H _{N,I\cup J,W\boxtimes U}^{0,\leq\mu,E}}{(X\sm N)^{I}\times \Delta(X\sm N)}  
   \end{gather*}
   où $\chi_{\zeta_{\emptyset}^{I}}$ et $\chi_{\zeta_{J}^{I}}$ sont  les   isomorphismes  de coalescence de \eqref{intro-isom-coalescence}. 
   De même on définit le  morphisme d'annihilation $
\mc C_{\xi}^{\flat}$  comme la  composée 
  \begin{gather*}
  \restr{ \mc H _{N,I\cup J,W\boxtimes U}^{0,\leq\mu,E}}{(X\sm N)^{I}\times \Delta(X\sm N)} 
  \isor{ \chi_{\zeta_{J}^{I}}}
 \mc H _{N,I\cup\{0\},W\boxtimes U^{\zeta_{J}}}^{0,\leq\mu,E} \\
 \xrightarrow{ \mc H(\Id_{W}\boxtimes \xi) }
  \mc H _{N,I\cup\{0\},W\boxtimes \mbf 1}^{0,\leq\mu,E} 
 \isor{\chi_{\zeta_{\emptyset}^{I}}^{-1} }
   \mc H _{N,I,W}^{0,\leq\mu,E} \boxtimes E_{(X\sm N)} . 
    \end{gather*}
    \end{defi}
    
   Tous les  morphismes ci-dessus sont des  morphismes de faisceaux sur  
    $(X\sm N)^{I} \times (X\sm N)$.

 Nous allons maintenant utiliser ces morphismes avec $J=\{1,2\}$. 
     Soit  $v$  une  place dans   $|X|\sm |N|$.    On considère  $v$ également  comme un  sous-schéma de  $X$ et on note   $E_{v}$ le faisceau constant sur   $v$.    Soit $V$ une représentation irréductible de $\wh G$.    On note   $\mbf 1\xrightarrow{\delta_{V}} V\otimes V^{*}$ et $ V\otimes V^{*}\xrightarrow{\on{ev}_{V}} \mbf 1$ les  morphismes naturels. 
        
      Pour  $\kappa$ assez grand  (en fonction de $\deg(v),V$), on définit   $S_{V,v}$ comme la  composée 
  \begin{gather}\label{def-SVv-intro1}
 \mc H _{N,I,W}^{0,\leq\mu,E} \boxtimes E_{v} \\ \label{def-SVv-intro2}
  \xrightarrow{ \restr{\mc C_{
    \delta_{V}}^{\sharp}}{(X\sm N)^{I}\times v}}
 \restr{ \mc H _{N,I\cup\{1,2\},W\boxtimes V\boxtimes V^{*}}^{0,\leq\mu,E}}{(X\sm N)^{I}\times \Delta(v)} \\ \label{def-SVv-intro3}
 \xrightarrow{ \restr{(F_{\{1\}})^{\deg(v)} }{(X\sm N)^{I}\times \Delta(v)}}
 \restr{ \mc H _{N,I\cup\{1,2\},W\boxtimes V\boxtimes V^{*}}^{0,\leq\mu+\kappa,E}}{(X\sm N)^{I}\times \Delta(v)} \\ \label{def-SVv-intro4}
  \xrightarrow{ \restr{\mc C_{
    \on{ev_{V}}}^{\flat }}{(X\sm N)^{I}\times v}}
     \mc H _{N,I,W}^{0,\leq\mu+\kappa,E}  \boxtimes E_{v} . 
 \end{gather}
Autrement dit on crée deux nouvelles pattes en  $v$ à l'aide de  
 $ \delta_{V}:\mbf 1\to V\otimes V^{*}$, on applique le  morphisme de Frobenius partiel (à la puissance  $\deg(v)$) à la première, puis on les annihile à l'aide de  $ \on{ev_{V}}: V\otimes V^{*}\to\mbf 1$. 
  
   En tant que  morphisme de faisceaux constructibles sur 
     $(X\sm N)^{I}\times v$,  $S_{V,v}$ commute avec l'action naturelle du morphisme de   Frobenius partiel sur   $E_{v}$ dans  \eqref{def-SVv-intro1} et  \eqref{def-SVv-intro4}, puisque 
     \begin{itemize}
     \item les  morphismes de création et  d'annihilation  entrelacent  cette action avec l'action  de $F_{\{1,2\}}$ sur  \eqref{def-SVv-intro2} et \eqref{def-SVv-intro3}, par la \remref{rem-coalescence-frob-compat}, 
          \item  $F_{\{1\}}$ et donc   $F_{\{1\}}^{\deg(v)}$ commutent avec  $F_{\{1,2\}}=F_{\{1\}}F_{\{2\}}$. 
     \end{itemize}
  \begin{defi} Par abus  on note encore     
        $$S_{V,v}: \mc H _{N,I,W}^{0,\leq\mu,E}\to 
         \mc H _{N,I,W}^{0,\leq\mu+\kappa,E}$$  le  morphisme 
    de faisceaux sur   $(X\sm N)^{I}$ obtenu par descente relativement à $\Z/\deg(v)\Z$ (en prenant les invariants par  l'action naturelle du morphisme de   Frobenius partiel sur   $E_{v}$ dans  \eqref{def-SVv-intro1} et  \eqref{def-SVv-intro4}). 
     \end{defi}
  
 \begin{prop}\label{prop-coal-frob-cas-part-intro} 
  La  restriction de  $S_{V,v}$  à  $(X\sm (N \cup v))^{I}$
est égale, en tant que   morphisme de faisceaux sur  $(X\sm (N \cup v))^{I}$,  à l'opérateur de Hecke $$T(h_{V,v}):  \restr{\mc H _{N,I,W}^{0,\leq\mu,E}}{(X\sm (N \cup v))^{I}}\to 
       \restr{  \mc H _{N,I,W}^{0,\leq\mu+\kappa,E}}{(X\sm (N \cup v))^{I}}
       . $$
       \end{prop}
\noindent
  Il suffit de le démontrer lorsque   $V$  et $W$ sont irréductibles. La preuve est de nature  géométrique.  
  Nous l'esquissons ici dans un cas simple, où elle se réduit à l'intersection de deux sous-champs lisses  dans un champ de Deligne-Mumford lisse et où cette intersection s'avère être transverse. La preuve est plus compliquée en général à cause des singularités. On renvoie à la preuve de la \propref{prop-coal-frob-cas-part}  pour le cas général (mais une solution alternative consisterait à se ramener à une situation d'intersection transverse lisse à l'aide des résolutions de Bott-Samelson). 
  
\noindent {\bf Démonstration lorsque  $V$ est minuscule  et  $\deg(v)=1$.} 
On rappelle qu'une représentation irréductible de $\wh G$ est dite minuscule si tous ses poids sont conjugués par le groupe de Weyl. Cela équivaut au fait que l'orbite correspondante dans la grassmannienne affine est fermée  (et implique donc que la strate fermée correspondante  est lisse). On note $d$ la dimension de l'orbite associée à $V$. 

 Grâce à l'hypothèse que  $\deg(v)=1$ on peut supprimer   
$\boxtimes E_{v}$ partout. 
On considère le   champ de Deligne-Mumford  
$$\mc Z^{(\{1\},\{2\}, I)}=\restr{\Cht_{N,I \cup \{1,2\},W \boxtimes V\boxtimes V^{*}}^{(\{1\},\{2\},I)}}
 {(X\sm (N\cup v))^{I}\times \Delta(v)}.$$
On va construire deux sous-champs fermés   $\mc Y_{1}$ et $\mc Y_{2}$ dans $\mc Z^{(\{1\},\{2\}, I)}$,    munis de   morphismes 
 $\alpha_{1}$ et $\alpha_{2}$ vers  $$ \mc Z^{(I)}=\restr{\Cht_{N,I ,W }^{(I)}}
 {(X\sm (N\cup v))^{I}}$$ de sorte que 
 \begin{itemize}
 \item {\bf A)} la  restriction à  $(X\sm (N\cup v))^{I}$ de la  composée 
\eqref{def-SVv-intro1}$\to$\eqref{def-SVv-intro2}$\to$\eqref{def-SVv-intro3}
du  morphisme de création
 et de l'action du morphisme de Frobenius partiel  est réalisée par une correspondance cohomologique   supportée par la correspondance $\mc Y_{2}$ de  $\mc Z^{(I)}$ vers    $\mc Z^{(\{1\},\{2\}, I)}$, et dont la restriction aux ouverts de lissité est déterminée par son support, avec un coefficient correctif de $q^{-d/2}$, 
  \item {\bf B)} la  restriction à $(X\sm (N\cup v))^{I}$ du morphisme d'annihilation    \mbox{\eqref{def-SVv-intro3}$\to$\eqref{def-SVv-intro4}}  est réalisée par une  correspondance cohomologique   supportée par    la   correspondance $\mc Y_{1}$ de  $\mc Z^{(\{1\},\{2\}, I)}$ vers $\mc Z^{(I)}$, et dont la restriction aux ouverts de lissité est déterminée par son support.   
 \end{itemize}
 Le coefficient correctif de $q^{-d/2}$ est évidemment dû à l'action du morphisme de Frobenius partiel. On justifiera dans la  \remref{rem-justif-norm-crea-annihil}   qu'il n'y a pas de signe à ajouter.

Donc $S_{V,v}$ sera réalisée par une correspondance cohomologique supportée par   le produit 
    $\mc Y_{1}\times_{\mc Z^{(\{1\},\{2\}, I)}}\mc Y_{2}$ de ces correspondances. On verra 
    \begin{itemize}
    \item que le  produit $\mc Y_{1}\times_{\mc Z^{(\{1\},\{2\}, I)}}\mc Y_{2}$ n'est autre que la   correspondance  de Hecke $\Gamma^{(I)}$ de  
$ \mc Z^{(I)}$ dans lui-même  (qui est une correspondance finie  étale)
\item que $S_{V,v}$, qui est donc une correspondance cohomologique supportée par $\Gamma^{(I)}$ est en fait égale à la correspondance cohomologique   supportée par $\Gamma^{(I)}$ avec un coefficient correctif de $q^{-d/2}$  
 (cette correspondance  réalise 
  $T(h_{V,v})$ puisque $V$ est minuscule). \end{itemize}
   Grâce à l'hypothèse que $V$ est minuscule il suffira de faire  le  calcul  sur les ouverts de lissité, et le calcul  sera évident car on verra que sur les ouverts de lissité l'intersection $\mc Y_{1}\times_{\mc Z^{(\{1\},\{2\}, I)}}\mc Y_{2}$ 
   est une intersection transverse de deux sous-champs lisses.

 On construit maintenant tous ces objets. 
 La correspondance de Hecke  $\Gamma^{(I)}$ est le  champ 
classifiant la donnée de   $(x_i)_{i\in I}$ et d'un diagramme 
\begin{gather} \label{intro-diag-Gamma}
 \xymatrix{
 (\mc G', \psi') \ar[r]^-{\phi'} & 
 (\ta{\mc G'}, \ta \psi')   \\
  (\mc G, \psi)  \ar[r]^-{\phi}  \ar[u]_-{\kappa}  &
 (\ta{\mc G}, \ta \psi) \ar[u]_-{\ta \kappa}
 } \end{gather}
 tel que   
 \begin{itemize}
 \item la ligne inférieure $\big( (x_i)_{i\in I}, (\mc G, \psi) \xrightarrow{\phi}   (\ta{\mc G}, \ta \psi)
\big)$
et la ligne supérieure  $\big( (x_i)_{i\in I}, (\mc G', \psi') \xrightarrow{\phi'}   (\ta{\mc G'}, \ta \psi')
\big)
$ appartiennent à  $\mc Z^{(I)}$, 
\item  $\kappa:\restr{\mc G}{(X\sm v)\times S}\isom \restr{\mc G'}{(X\sm v)\times S}$  est un  isomorphisme 
tel que la   position relative de $\mc G$ par rapport à  $\mc G'$ en   $v$ est {\it égale} au poids dominant de $V$ (on rappelle que  $V$ est minuscule),
\item la restriction de $\kappa$ à $N\times S$, qui est un isomorphisme,  entrelace les structures de niveau $\psi$ et $\psi'$.   
\end{itemize}
De plus les deux   projections 
 $\Gamma^{(I)}\to \mc Z^{(I)}$ sont les  morphismes qui conservent les lignes inférieures et supérieures de   \eqref{intro-diag-Gamma}. 
  
 Comme les  pattes indexées par   $I$ varient dans   $X\sm (N\cup v)$ et restent   disjointes  des  pattes $1$ et   $2$ fixées en  $v$, on 
 peut changer la partition $(\{1\},\{2\}, I)$ en $(\{1\},I,\{2\})$ et on a donc
    \begin{gather*}\mc Z^{(\{1\},\{2\}, I)}
=
 \restr{\Cht_{N,I \cup \{1,2\},W \boxtimes V\boxtimes V^{*}}^{(\{1\},I,\{2\})}}
 {(X\sm (N\cup v))^{I}\times \Delta(v)}.\end{gather*}
 Autrement dit le champ  $\mc Z^{(\{1\},\{2\}, I)}$ classifie la donnée de  $(x_i)_{i\in I}$ et d'un diagramme 
 \begin{gather} \label{intro-diag-W}
 \xymatrix{
 & (\mc G_{1}, \psi_{1}) \ar[r]^-{\phi'_{2}} \ar[d]^-{\phi_{2}}& 
 (\mc G'_{2}, \psi'_{2}) \ar[d]^-{\phi'_{3}} & 
 (\ta{\mc G_{1}}, \ta \psi_{1}) \\
 (\mc G_{0}, \psi_{0}) \ar[ru]^-{\phi_{1}} &
 (\mc G_{2}, \psi_{2}) \ar[r]^-{\phi_{3}}    &
 (\ta{\mc G_{0}}, \ta \psi_{0})\ar[ru]^-{\ta \phi_{1}} &
 } \end{gather}
 avec 
\begin{gather*}\big( (x_i)_{i\in I}, (\mc G_{0}, \psi_{0}) \xrightarrow{\phi_{1}}   (\mc G_{1}, \psi_{1}) \xrightarrow{\phi_{2}} (\mc G_{2}, \psi_{2}) \xrightarrow{\phi_{3}}    (\ta{\mc G_{0}}, \ta \psi_{0})
\big)\\
\in \restr{\Cht_{N,I \cup \{1,2\},W \boxtimes V\boxtimes V^{*}}^{(\{1\},\{2\},I)}}
 {(X\sm (N\cup v))^{I}\times \Delta(v)}\end{gather*} et  \begin{gather*}\big( (x_i)_{i\in I}, (\mc G_{0}, \psi_{0}) \xrightarrow{\phi_{1}}   (\mc G_{1}, \psi_{1}) \xrightarrow{\phi'_{2}} (\mc G'_{2}, \psi'_{2}) \xrightarrow{\phi'_{3}}    (\ta{\mc G_{0}}, \ta \psi_{0})
\big)\\
\in \restr{\Cht_{N,I \cup \{1,2\},W \boxtimes V\boxtimes V^{*}}^{(\{1\},I,\{2\})}}
 {(X\sm (N\cup v))^{I}\times \Delta(v)}. \end{gather*}  
 Les flèches  obliques, verticales et horizontales du  diagramme  \eqref{intro-diag-W} sont respectivement les   modifications associées à la  patte  $1$, à  la  patte  $2$ et aux pattes indexées par   $I$. La flèche $\ta \phi_{1}$ à droite du diagramme \eqref{intro-diag-W} est déterminée par $\phi_{1}$, mais on l'a dessinée car elle servira pour définir $\mc Y_{2}$ ci-dessous.

 On note  $\mc Y_{1} \overset{\iota_{1}}{\hookrightarrow} 
 \mc Z^{(\{1\},\{2\}, I)}$ le sous-champ fermé défini par la condition que dans le diagramme  \eqref{intro-diag-W},  $\phi_{2}\phi_{1}: 
 \restr{\mc G_{0}}{(X-v)\times S}\to \restr{\mc G_{2}}{(X-v)\times S}$ 
 s'étend en un  isomorphisme sur  $X\times S$. 
On a un   morphisme $$\alpha_{1}: \mc Y_{1}\to \mc Z^{(I)}
 $$
qui envoie    \begin{gather*}
 \xymatrix{
 & (\mc G_{1}, \psi_{1}) \ar[r]^-{\phi'_{2}} \ar[d]^-{\phi_{2}}& 
 (\mc G'_{2}, \psi'_{2}) \ar[d]^-{\phi'_{3}} & 
 (\ta{\mc G_{1}}, \ta \psi_{1}) \\
 (\mc G_{0}, \psi_{0}) \ar[ru]^-{\phi_{1}} \ar[r]^-{\sim}&
 (\mc G_{2}, \psi_{2}) \ar[r]^-{\phi_{3}}    &
 (\ta{\mc G_{0}}, \ta \psi_{0})\ar[ru]^-{\ta \phi_{1}} &
 } \end{gather*}
sur la ligne du bas, c'est-à-dire   
 \begin{gather}\label{intro-ligne-bas}\big( (x_i)_{i\in I}, (\mc G_{0}, \psi_{0}) \xrightarrow{\phi_{3}(\phi_{2}\phi_{1})}    (\ta{\mc G_{0}}, \ta \psi_{0})
\big). \end{gather}

L'assertion  B)  ci-dessus vient d'un énoncé similaire concernant les faisceaux de Mirkovic-Vilonen. En effet 
\begin{itemize}\item par le a) du \thmref{satake-geom-thm} l'image directe de 
 $\mc S_{\{1,2\},V\boxtimes V^{*} ,E}^{(\{1\}, \{2\})}$ (qui est le faisceau constant $E$ décalé) par le morphisme d'oubli
 (de la modification intermédiaire) 
 $\mr{Gr}_{\{1,2\},V\boxtimes V^{*} }^{(\{1\}, \{2\})}\to 
 \mr{Gr}_{\{1,2\},V\boxtimes V^{*} }^{(\{1,2\})}$ est égale à 
 $\mc S_{\{1,2\},V\boxtimes V^{*} ,E}^{(\{1,2\})}$, 
 \item  par le c) du \thmref{satake-geom-thm} la restriction de  $\mc S_{\{1,2\},V\boxtimes V^{*} ,E}^{(\{1,2\})}$ au-dessus de la diagonale 
 (et donc en particulier au-dessus de $\Delta(v)$) est égale à 
  $\mc S_{\{0\},V\otimes V^{*} ,E}^{(\{0\})}$ que l'on envoie dans 
  le faisceau gratte-ciel $\mc S_{\{0\},\mbf 1,E}^{(\{0\})}$ par $\on{ev}_{V}:V\otimes V^{*}\to \mbf 1$
\end{itemize}
et par le théorème de changement de base propre cela donne lieu à une correspondance cohomologique entre $\restr{\mr{Gr}_{\{1,2\},V\boxtimes V^{*} }^{(\{1\}, \{2\})}}{\Delta(v)}$ et le point, et on vérifie que celle-ci est la correspondance cohomologique  évidente supportée par le sous-schéma fermé lisse de $\restr{\mr{Gr}_{\{1,2\},V\boxtimes V^{*} }^{(\{1\}, \{2\})}}{\Delta(v)}$ formé des 
$ (\mc G_{0} \xrightarrow{\phi_{1}}  
\mc G_{1}\xrightarrow{\phi_{2}}
  \mc G_{2}  \isom G ) $ tels que $\phi_{2}\phi_{1}$ soit un isomorphisme. 
  
  On note  
$\mc Y_{2} \overset{\iota_{2}}{\hookrightarrow} \mc Z^{(\{1\},\{2\}, I)}$ 
le sous-champ fermé  défini par la condition que  dans  le diagramme  \eqref{intro-diag-W},   $\ta \phi_{1}\phi_{3}': 
 \restr{\mc G'_{2}}{(X-v)\times S}\to \restr{\ta{\mc G_{1}}}{(X-v)\times S}$ s'étend un  isomorphisme sur $X\times S$. 
 On a un   morphisme $$\alpha_{2}: \mc Y_{2}\to \mc Z^{(I)}
 $$ 
qui envoie  \begin{gather}\label{intro-diag-élément-Y2}
 \xymatrix{
 & (\mc G_{1}, \psi_{1}) \ar[r]^-{\phi'_{2}} \ar[d]^-{\phi_{2}}& 
 (\mc G'_{2}, \psi'_{2}) \ar[d]^-{\phi'_{3}} \ar[r]^-{\sim}& 
 (\ta{\mc G_{1}}, \ta \psi_{1}) \\
 (\mc G_{0}, \psi_{0}) \ar[ru]^-{\phi_{1}} &
 (\mc G_{2}, \psi_{2}) \ar[r]^-{\phi_{3}}    &
 (\ta{\mc G_{0}}, \ta \psi_{0})\ar[ru]^-{\ta \phi_{1}} &
 } \end{gather}
sur la ligne du haut, c'est-à-dire   
 \begin{gather}\label{intro-ligne-haut}\big( (x_i)_{i\in I}, (\mc G_{1}, \psi_{1}) \xrightarrow{(\ta \phi_{1}\phi'_{3})\phi'_{2}}    (\ta{\mc G_{1}}, \ta \psi_{1})
\big). \end{gather}
La justification de l'assertion A)  ci-dessus se fait
\begin{itemize}
\item par un argument similaire à celui donné pour justifier B)
mais concernant cette fois-ci  $\delta_{V}:\mbf 1\to V\otimes V^{*}$ et le champ
$\restr{\Cht_{N,I \cup \{1,2\},W \boxtimes V\boxtimes V^{*}}^{(I,\{2\},\{1\})}}
 {(X\sm (N\cup v))^{I}\times \Delta(v)}$
\item par le fait que la restriction à  $ (X\sm (N\cup v))^{I}\times \Delta(v)$ du morphisme de Frobenius partiel
$$\on {Fr}_{\{1\}} ^{(\{1\},I,\{2\})}:  \Cht_{N,I \cup \{1,2\},W \boxtimes V\boxtimes V^{*}}^{(\{1\},I,\{2\})} \to 
 \Cht_{N,I \cup \{1,2\},W \boxtimes V\boxtimes V^{*}}^{(I,\{2\},\{1\})}$$
envoie \eqref{intro-diag-W} sur 
\begin{gather}  \nonumber 
 \xymatrix{
   (\mc G_{1}, \psi_{1}) \ar[r]^-{\phi'_{2}}  & 
 (\mc G'_{2}, \psi'_{2}) \ar[d]^-{\phi'_{3}} & 
 (\ta{\mc G_{1}}, \ta \psi_{1}) \\
      &
 (\ta{\mc G_{0}}, \ta \psi_{0})\ar[ru]^-{\ta \phi_{1}} &
 } \end{gather}
\end{itemize}

Le fait qu'on n'ait pas besoin d'introduire de signe  
dans A) et B) ci-dessus sera justifié dans la \remref{rem-justif-norm-crea-annihil}, que le lecteur peut lire dès à présent s'il le souhaite. 

D'autre part on a un  isomorphisme canonique 
\begin{gather}\label{intro-Gamma-produit-fibre}
\Gamma^{(I)}\simeq \mc Y_{1} \times_{\mc Z^{(\{1\},\{2\}, I)}}\mc Y_{2}. 
\end{gather}
En effet un  point  de  $\mc Z^{(\{1\},\{2\}, I)}$ appartenant à   $\mc Y_{1} $ et  à $\mc Y_{2}$ est donné par un  diagramme 
  $$ \xymatrix{
 & (\mc G_{1}, \psi_{1}) \ar[r]^-{\phi'_{2}} \ar[d]^-{\phi_{2}}& 
 (\mc G'_{2}, \psi'_{2}) \ar[d]^-{\phi'_{3}} \ar[r]^-{\sim}& 
 (\ta{\mc G_{1}}, \ta \psi_{1}) \\
 (\mc G_{0}, \psi_{0}) \ar[ru]^-{\phi_{1}} \ar[r]^-{\sim}&
 (\mc G_{2}, \psi_{2}) \ar[r]^-{\phi_{3}}    &
 (\ta{\mc G_{0}}, \ta \psi_{0})\ar[ru]^-{\ta \phi_{1}} &
 } $$
Il équivaut donc  à la donnée d'un point de  $ \Gamma^{(I)}$, car  
 en contractant les deux   isomorphismes du   diagramme précédent on obtient  le diagramme 
  \begin{gather*}
 \begin{CD} 
  (\mc G_{1}, \psi_{1}) 
  @>(\ta \phi_{1}\phi'_{3}) \phi'_{2}>> 
  (\ta{\mc G_{1}}, \ta \psi_{1})  \\
 @AA\phi_{1}A 
 @AA\ta \phi_{1}A \\
(\mc G_{0}, \psi_{0})
 @>\phi_{3}(\phi_{2}\phi_{1})>>  
 (\ta{\mc G_{0}}, \ta \psi_{0}) 
 \end{CD}
 \end{gather*}
que l'on identifie au diagramme \eqref{intro-diag-Gamma}. 

On a des morphismes naturels des champs $\mc Z^{(\{1\},\{2\}, I)}, \mc Z^{(I)} , \mc Y_{1}, \mc Y_{2}$  et $\Gamma^{(I)}$ vers $ \mr{Gr}_{I,W}^{(I)}/
  G_{\sum n_{i} x_i}$. Comme $V$ est minuscule 
  ces morphismes sont lisses. Donc les 
   ouverts de lissité 
  ${}^{\circ}{\mc Z}^{(\{1\},\{2\}, I)}, {}^{\circ}{\mc Z}^{(I)},  {}^{\circ}{\mc Y}_{1}, {}^{\circ}{\mc Y}_{2}, {}^{\circ}{\Gamma}^{(I)}$
sont les images inverses de $ {}^{\circ}{\mr{Gr}}_{I,W}^{(I)}/
  G_{\sum n_{i} x_i}$ où $ {}^{\circ}{\mr{Gr}}_{I,W}^{(I)}$ désigne l'ouvert de lissité de 
$ \mr{Gr}_{I,W}^{(I)}$. 

Un calcul d'espaces tangents  montre que 
 ${}^{\circ}\mc Y_{1}$ et ${}^{\circ}\mc Y_{2}$ sont des sous-champs lisses dans le 
champ de Deligne-Mumford  lisse ${}^{\circ}\mc Z^{(\{1\},\{2\}, I)}$ et s'y s'intersectent transversalement, et de plus il résulte de    \eqref{intro-Gamma-produit-fibre} que leur  intersection est 
${}^{\circ}\Gamma^{(I)}$.   
On a donc l'égalité de  correspondances cohomologiques
entre $S_{V,v}$ et $T(h_{V,v})$ sur ${}^{\circ}\Gamma^{(I)}$ mais comme 
$\Gamma^{(I)}$ est une correspondance   étale  entre 
$\mc Z^{( I)}$ et lui-même l'égalité a lieu partout
(en effet un morphisme du faisceau pervers $\on{IC}_{\Gamma^{(I)}}$ dans lui-même est déterminé par sa restriction à ${}^{\circ}{\Gamma}^{(I)}$). 
\cqfd

    Une  conséquence de la   \propref{prop-coal-frob-cas-part-intro} est que l'on possède,  pour  tout    $f\in C_{c}(K_{N}\backslash G(\mb A)/K_{N},E)$ et  $\kappa$ assez grand,  
    une extension naturelle 
  du  morphisme $T(f)$ (introduit dans  \eqref{defi-Tf})  en un  morphisme  
 $T(f):\mc H_{N,I,W}^{0,\leq\mu,E}\to 
 \mc H_{N,I,W}^{0,\leq\mu+\kappa,E}$ de faisceaux constructibles sur  $  (X\sm N)^{I}$ tout entier, de fa\c con compatible avec la  composition des   opérateurs de Hecke. En effet, en notant $K_{N}=\prod K_{N,v}$,  il suffit de le montrer pour toute place $v$ et pour   $f\in C_{c}(K_{N,v}\backslash G(F_{v})/K_{N,v},E)$.
  Il n'y a rien à faire  si $v\in N$. Si $v\not\in N$ 
  il suffit de traiter le cas où $f=h_{V,v}$ et alors 
  l'extension est donnée par $S_{V,v}$ grâce à la  \propref{prop-coal-frob-cas-part-intro}. Pour plus de détails, on renvoie au \corref{cor-hecke-etendus-compo}.  
 
 Pour les variétés de Shimura sur les corps de nombres de tels prolongements ont été définis dans de nombreux cas, de fa\c con modulaire, par adhérence de  Zariski ou à l'aide de  cycles proches (voir \cite{deligne-bki-mod,faltings-chai, genestier-tilouine}). 

     Comme $S_{V,v}$ est le prolongement de 
   $T(h_{V,v})$, la proposition suivante exprime exactement la 
  relation d'Eichler-Shimura. 
 Cette relation affirme que pour toute place finie et pour tout $i\in I$ le morphisme de Frobenius partiel en $v$ est annulé par un polynôme en les opérateurs de Hecke en $v$ (à coefficients dans $\mc O_{E}$).

  On utilise encore $\{0\}$ pour noter un ensemble à un élement (indexant la patte à laquelle s'applique la relation  d'Eichler-Shimura). 
 
 \begin{prop}\label{Eichler-Shimura-intro} (\propref{prop-eichler-shimura}) 
 Soient $I,W$ comme ci-dessus et $V$  une   représentation irréductible de  $\wh G$. Alors  $$F_{\{0\}}^{\deg(v)}: \varinjlim _{\mu}\restr{\mc H _{N, I\cup\{0\}, W\boxtimes V}^{0,\leq\mu,E}}{(X\sm N)^{I}\times v}\to \varinjlim _{\mu}\restr{\mc H _{N, I\cup\{0\}, W\boxtimes V}^{0,\leq\mu,E}}{(X\sm N)^{I}\times v}$$ est annulé par un  polynôme de degré 
 $\dim(V)$ dont les  coefficients sont les 
  restrictions à $(X\sm N)^{I}\times v$ des morphismes $S_{\Lambda^{i}V,v}$. 
Plus précisément on a  
  \begin{gather*}  \sum_{i=0}^{\dim V} (-1)^{i} (F_{\{0\}}^{\deg(v)})^{i}\circ \restr{S_{\Lambda^{\dim V-i}V,v}}{(X\sm N)^{I}\times v}=0  . \end{gather*} 
    \end{prop}
    On rappelle que  $S_{\Lambda^{i}V,v}$ étend  l'opérateur de Hecke 
 $T(h_{\Lambda^{i}V,,v})$ $$\text{de \ \ }  (X\sm (N\cup v))^{I\cup \{0\}}\text{ \ \  à \ \ }(X\sm N)^{I\cup \{0\}}$$ et on remarque que cette  extension est absolument nécessaire pour prendre la restriction à  $(X\sm N)^{I}\times v$. Grâce à la définition des morphismes   $S_{\Lambda^{i}V,v}$ par  \eqref{def-SVv-intro1}-\eqref{def-SVv-intro4}, 
 la preuve de la \propref{Eichler-Shimura-intro} est un simple calcul d'algèbre tensorielle  (inspiré d'une démonstration  du théorème de  Hamilton-Cayley, et fondée uniquement sur le fait que $\Lambda^{\dim V+1}V=0$). 
 On renvoie au chapitre \ref{para-Relations d'Eichler-Shimura} pour cette preuve.

 \subsection{Propriétés des faisceaux de cohomologie des champs de chtoucas}
              
       On rappelle  que si  $\ov x$ est un point géométrique d'un schéma $Y$
       son localisé strict  (ou   hensélisé strict) $Y_{(\ov x)}$ est défini comme 
       la limite projective des  voisinage étales  $\ov x$-pointés de  $x$.
       Si $\ov y$ est un autre point géométrique  on 
       appelle flèche de spécialisation  $\on{\mf{sp}}:\ov x\to \ov y$   un  morphisme
     $Y_{(\ov x)}\to Y_{(\ov y)}$, ou de fa\c con équivalente un    morphisme  
 $\ov x\to Y_{(\ov y)}$. Pour tout faisceau   $\mc F$ sur $Y$ (pour la topologie étale), $\on{\mf{sp}}$ induit 
       un homomorphisme de spécialisation $\on{\mf{sp}}^{*}: \mc F_{\ov y}\to \mc F_{\ov x}$ (voir le paragraphe 7 de \cite{grothendieck-sga4-2-VIII}). 
          
  On fixe une clôture algébrique $\ov F$ de  $F$ et on  note 
      $\ov \eta=\on{Spec}(\ov F)$ le point géométrique correspondant au-dessus du   point  générique
 $\eta$    de $X$.  
 
  On note $\Delta:X\to X^{I}$ le   morphisme diagonal.  
 On fixe un  point géométrique $\ov{\eta^{I}}$ au-dessus du point  générique $\eta^{I}$ de $X^{I}$ et une flèche de spécialisation   
 $\on{\mf{sp}}: \ov{\eta^{I}}\to \Delta(\ov \eta)$. Le rôle de  $\on{\mf{sp}}$  est de rendre le foncteur fibre en  $\ov{\eta^{I}}$ plus canonique, et en particulier compatible avec la  coalescence des pattes (cette dernière affirmation est claire lorsque   $\on{\mf{sp}}^{*}$ est un isomorphisme et en pratique nous serons dans cette situation). Donc $\ov{\eta^{I}}$ et $\on{\mf{sp}}$  vont ensemble et ci-dessous les énoncés faisant intervenir $\ov{\eta^{I}}$ dépendent du choix de $\on{\mf{sp}}$.

 Un résultat  fondamental de Drinfeld (théorème 2.1 de \cite{drinfeld78} et proposition 6.1 de \cite{drinfeld-compact}) est rappelé  dans le  
 lemme suivant (voir le chapitre \ref{para-sous-faisceaux-Frob-partiels} pour d'autres références, notamment \cite{eike-lau}).
  On notera toujours les  $\mc O_{E}$-modules et  les $\mc O_{E}$-faisceaux   par  des lettres gothiques. 
   
 \begin{lem} \label{lem-Dr-intro} (Drinfeld)  Si 
 $\mf E$ est un $\mc O_{E}$-faisceau   lisse   constructible 
 sur un ouvert dense     de  
 $(X\sm N)^{I}$, muni d'une action des  morphismes de  Frobenius partiels, c'est-à-dire  d'isomorphismes
 $$F_{\{i\}}:\restr{\Frob_{\{i\}}^{*}(\mf E)}{\eta^{I}}\to \restr{\mf E}{\eta^{I}}$$
 commutant entre eux et dont la composée est  l'isomorphisme naturel 
    $\Frob^{*}(\mf E)\isom \mf E$, 
alors il s'étend en un faisceau lisse sur  $U^{I}$, où  $U$ est  un ouvert dense assez petit  de  $X\sm N$, et la fibre  $\restr{\mf E}{\Delta(\ov \eta)}$ est munie d'une action de   $\pi_{1}(U,\ov\eta)^{I}$. De plus, si on fixe  $U$,  le  foncteur
 $\mf E\mapsto \restr{\mf E}{\Delta(\ov \eta)}$ 
 fournit  une équivalence 
 \begin{itemize}
 \item de la  catégorie des   $\mc O_{E}$-faisceaux
lisses  constructibles  sur   $U^{I}$ munis d'une action des morphismes de  Frobenius partiels 
\item
 vers la  catégorie des 
 représentations continues de   $\pi_{1}(U,\ov\eta)^{I}$ sur des  $\mc O_{E}$-modules de   type fini, 
 \end{itemize}
 de fa\c con  compatible avec la  coalescence (c'est-à-dire avec l'image inverse par le morphisme $U^{J}\to U^{I}$ associé à toute application $I\to J$).  
 \end{lem}
       
       \begin{rem} \label{rem-lem-Dr} Dans la situation du lemme précédent, $\on{\mf{sp}}^{*}: 
       \restr{\mf E}{\Delta(\ov \eta)}\to \restr{\mf E}{\ov{\eta^{I}}}$ est un isomorphisme, donc $\restr{\mf E}{\ov{\eta^{I}}}$ est muni lui aussi d'une action de $\pi_{1}(U,\ov\eta)^{I}$. 
       \end{rem}
       
       Soit $I$ un ensemble fini et  $W=\boxtimes_{i\in I}W_{i}$ une   représentation irréductible de $(\wh G)^{I}$.  
       
 On ne peut pas appliquer directement le lemme précédent, car l'action des  morphismes de  Frobenius partiels augmente $\mu$, et d'autre part la limite  inductive  $\varinjlim_{\mu}\mc H_{N,I,W}^{0,\leq\mu,E}$ n'est pas constructible  (car ses fibres sont de  dimension infinie). 
 Mais on pourra l'appliquer à la partie ``Hecke-finie'', au sens suivant.

 \begin{defi}
Soit  $\ov x$ un  point géométrique de  $(X\sm N)^{I}$. 
Un élément de  $\varinjlim _{\mu}\restr{\mc H _{N, I, W}^{0,\leq\mu,E}}{\ov{x}}$ est dit   Hecke-fini s'il appartient à un  sous-$\mc O_{E}$-module  de type fini de  $\varinjlim _{\mu}\restr{\mc H _{N, I, W}^{0,\leq\mu,E}}{\ov{x}}$ qui est 
stable par  $T(f)$ pour tout  $f\in C_{c}(K_{N}\backslash G(\mb A)/K_{N},\mc O_{E})$. 
\end{defi}
On note  $\Big( \varinjlim _{\mu}\restr{\mc H _{N, I, W}^{0,\leq\mu,E}}{\ov{x}}\Big)^{\mr{Hf}}$ l'ensemble de tous les  éléments Hecke-finis. 
C'est un sous-$E$-espace vectoriel de   $ \varinjlim _{\mu}\restr{\mc H _{N, I, W}^{0,\leq\mu,E}}{\ov{x}}$ et il est stable 
 par  $\pi_{1}(x,\ov{x})$ et  $C_{c}(K_{N}\backslash G(\mb A)/K_{N},E)$. 

\begin{rem}
La définition ci-dessus sera appliquée avec  $\ov x$ égal à 
$\Delta(\ov{\eta})$ ou $\ov{\eta^{I}}$. Dans ce cas l'action des opérateurs de Hecke $T(f)$ sur $\varinjlim _{\mu}\restr{\mc H _{N, I, W}^{0,\leq\mu,E}}{\ov{x}}$ est évidente (et ne nécessite pas leur extension en des morphismes de faisceaux sur $(X\sm N)^{I}$  réalisée après la \propref{Eichler-Shimura-intro}). 
\end{rem}
 
 On possède l'homomorphisme de spécialisation
 \begin{gather}\label{sp*-sans-Hf-intro}\on{\mf{sp}}^{*}: 
 \varinjlim _{\mu}\restr{\mc H _{N, I, W}^{0,\leq\mu,E}}{\Delta(\ov{\eta})} \to
  \varinjlim _{\mu}\restr{\mc H _{N, I, W}^{0,\leq\mu,E}}{\ov{\eta^{I}}}\end{gather} 
 où les deux membres sont considérés comme des $E$-espaces vectoriels
 (limites inductives de $E$-espaces vectoriels de dimension finie).

 \begin{lem}\label{lem-Hf-union-stab} (\propref{prop-action-Hf})
 L'espace   $\Big( \varinjlim _{\mu}\restr{\mc H _{N, I, W}^{0,\leq\mu,E}}{\ov{\eta^{I}}}\Big)^{\mr{Hf}}$  est la réunion  de sous-$\mc O_{E}$-modules $\mf M=\restr{\mf G}{\ov{\eta^{I}}}$ où 
   $\mf G$ est un sous-$\mc O_{E}$-faisceau constructible   de 
  $\varinjlim _{\mu}\restr{\mc H _{N, I, W}^{0,\leq\mu,E}}{\eta^{I}}$ 
     stable sous l'action des morphismes de  Frobenius partiels. 
 \end{lem}
  \noindent {\bf Démonstration.} On renvoie à la démonstration de   la \propref{prop-action-Hf} pour plus de détails. Il suffit de traiter le cas où $W=\boxtimes_{i\in I}W_{i}$ est irréductible. 
Pour toute famille  $(v_{i})_{i\in I}$ de points fermés de  $X\sm N$, on note  ${\times_{i\in I} v_{i}}$ leur produit, qui est une réunion finie de  points fermés de  $(X\sm N)^{I}$. 
  Soit $\check{\mf M}$ un  sous-$\mc O_{E}$-module de type fini de 
  $\varinjlim _{\mu}\restr{\mc H _{N, I, W}^{0,\leq\mu,E}}{\ov{\eta^{I}}}$ stable par  
$\pi_{1}(\eta^{I},\ov{\eta^{I}})$ et   $C_{c}(K_{N}\backslash G(\mb A)/K_{N},\mc O_{E})$. 
On va construire $\mf M\supset \check{\mf M}$ vérifiant les propriétés de l'énoncé. 
Comme  $\check{\mf M}$ est de type fini, il existe $\mu_0$ tel que  $\check{\mf M}$ soit inclus dans  l'image de  $\restr{\mc H _{N, I, W}^{0,\leq\mu_0,E}}{\ov{\eta^{I}}}$ dans  $\varinjlim _{\mu}\restr{\mc H _{N, I, W}^{0,\leq\mu,E}}{\ov{\eta^{I}}}$. Quitte à augmenter  $\mu_0$, on peut supposer que   $\check{\mf M}$ est un sous-$\mc O_{E}$-module  de  $\restr{\mc H _{N, I, W}^{0,\leq\mu_0,E}}{\ov{\eta^{I}}}$. 
Soit $\Omega_0$ un ouvert dense de  $X^{I}$ sur lequel  
$ \mc H _{N, I, W}^{0,\leq\mu_0,E}$ est lisse. Il existe un  unique 
sous-$\mc O_{E}$-faisceau lisse $\check{\mf G}\subset  \restr{\mc H _{N, I, W}^{0,\leq\mu_0,E}}{\Omega_0}$ sur  $\Omega_0$ tel que  $\restr{\check{\mf G}}{\ov{\eta^{I}}}=\check{\mf M}$. On choisit  
 $(v_{i})_{i\in I}$ tel que  ${\times_{i\in I} v_{i}}$ soit inclus dans  $\Omega_0$. Pour tout  $i$, la  relation d'Eichler-Shimura  
 (\propref{Eichler-Shimura-intro}) 
 implique alors que   \begin{gather}\label{intro-ES-inclusion}(F_{\{i\}}^{\deg(v_{i})})^{\dim W_{i}}(\restr{\check{\mf G}}{\times_{i\in I} v_{i}})
   \subset 
   \sum_{\alpha=0}^{\dim W_{i}-1}  (F_{\{i\}}^{\deg(v_{i})})^{\alpha}(S_{\Lambda^{\dim W_{i}-\alpha}W_{i},v_{i}}(\restr{\check{\mf G}}{\times_{i\in I} v_{i}}))
  \end{gather} dans  $
  \varinjlim_{\mu} \restr{\mc H _{N, I, W}^{0,\leq\mu,E}}{\times_{i\in I} v_{i}}$. 
Grâce à la lissité de 
$(\Frob_{\{i\}}^{\deg(v_{i})\dim W_{i}})^{*}(\check{\mf G})$
en ${\times_{i\in I} v_{i}}$, l'inclusion   \eqref{intro-ES-inclusion} se propage  en  $\eta^{I}$, c'est-à-dire  que 
 \begin{gather}\label{incl-F-Frob-Hecke-compl}
 F_{\{i\}}^{\deg(v_{i})\dim W_{i}}
 (\restr{(\Frob_{\{i\}}^{\deg(v_{i})\dim W_{i}})^{*}(\check{\mf G}}{\eta^{I}}) )
 \\  \nonumber  \subset 
   \sum_{\alpha=0}^{\dim W_{i}-1}  F_{\{i\}}^{\deg(v_{i})\alpha}(\Frob_{\{i\}}^{\deg(v_{i})\alpha})^{*}(\restr{S_{\Lambda^{\dim W_{i}-\alpha}W_{i},v_i}(\check{\mf G})}{\eta^{I}})
  \end{gather} dans  $
  \varinjlim_{\mu} \restr{\mc H _{N, I, W}^{0,\leq\mu,E}}{\eta^{I}}$. 
 Or   $\restr{\check{\mf G}}{\eta^{I}}$ est stable par  $S_{\Lambda^{\dim W_{i}-\alpha}W_{i},v_i}=T(h_{\Lambda^{\dim W_{i}-\alpha}W_{i},v_i})$ 
 puisque  
 $$h_{\Lambda^{\dim W_{i}-\alpha}W_{i},v_i}\in C_{c}(G(\mc O_{v_{i}})\backslash G(F_{v_{i}})/G(\mc O_{v_{i}}),\mc O_{E})\subset C_{c}(K_{N}\backslash G(\mb A)/K_{N},\mc O_{E}). $$
   Par conséquent \eqref{incl-F-Frob-Hecke-compl} se simplifie en 
 \begin{gather*}
 F_{\{i\}}^{\deg(v_{i})\dim W_{i}}
 (\restr{(\Frob_{\{i\}}^{\deg(v_{i})\dim W_{i}})^{*}(\check{\mf G}}{\eta^{I}}) )
   \subset 
   \sum_{\alpha=0}^{\dim W_{i}-1}  F_{\{i\}}^{\deg(v_{i})\alpha}(\Frob_{\{i\}}^{\deg(v_{i})\alpha})^{*}(\restr{\check{\mf G}}{\eta^{I}})
  \end{gather*} dans  $
  \varinjlim_{\mu} \restr{\mc H _{N, I, W}^{0,\leq\mu,E}}{\eta^{I}}$. 
  On en déduit que   $$\mf G=\sum_{(n_{i})_{i\in I}\in \prod _{i\in I}\{0,...,\deg(v_{i})\dim(W_{i})-1\}}\restr{\prod _{i\in I}F_{\{i\}}^{n_{i}}\Big(\prod _{i\in I}\Frob_{\{i\}}^{n_{i}}\Big)^{*}(\check{\mf G})}{\eta^{I}}$$
 est  un  sous-$\mc O_{E}$-faisceau constructible  de 
  $\varinjlim _{\mu}\restr{\mc H _{N, I, W}^{0,\leq\mu,E}}{\eta^{I}}$ 
qui  est stable sous  l'action des morphismes de  Frobenius partiels. Comme 
$\Big( \varinjlim _{\mu}\restr{\mc H _{N, I, W}^{0,\leq\mu,E}}{\ov{\eta^{I}}}\Big)^{\mr{Hf}}$  est la réunion  de sous-$\mc O_{E}$-modules $\check{\mf M}$ comme au début de la démonstration et que $ \mf M =\restr{ \mf G }{\ov{\eta^{I}}}$ contient $\check{\mf M}$ on obtient l'énoncé du lemme. 
 \cqfd

\begin{prop}\label{cor-action-Hf-intro}  
L'espace  $\Big( \varinjlim _{\mu}\restr{\mc H _{N, I, W}^{0,\leq\mu,E}}{\ov{\eta^{I}}}\Big)^{\mr{Hf}}$ est muni d'une action naturelle  de  $\pi_{1}(\eta,\ov{\eta})^{I}$. Plus précisément c'est une réunion de sous-$E$-espaces vectoriels de dimension finie   munis d'une action continue de 
$\pi_{1}(\eta,\ov{\eta})^{I}$. 
 \end{prop}
\dem Pour tout 
   sous-$\mc O_{E}$-faisceau constructible $\mf G$   de 
  $\varinjlim _{\mu}\restr{\mc H _{N, I, W}^{0,\leq\mu,E}}{\eta^{I}}$ 
     stable sous l'action des morphismes de  Frobenius partiels, 
      le lemme de  Drinfeld \lemref{lem-Dr-intro}
 fournit  (grâce à $\mf{sp}$ et à la \remref{rem-lem-Dr}) une action continue de   $\pi_{1}(\eta,\ov{\eta})^{I}$ sur 
$\mf M=\restr{\mf G}{\ov{\eta^{I}}}$. D'après le \lemref{lem-Hf-union-stab}, 
  $\Big( \varinjlim _{\mu}\restr{\mc H _{N, I, W}^{0,\leq\mu,E}}{\ov{\eta^{I}}}\Big)^{\mr{Hf}}$ est la réunion de  tels  $\mf M$. \cqfd

\begin{rem}\label{rem-WeilF} L'action de $\pi_{1}(\eta,\ov{\eta})^{I}$
sur $\Big( \varinjlim _{\mu}\restr{\mc H _{N, I, W}^{0,\leq\mu,E}}{\ov{\eta^{I}}}\Big)^{\mr{Hf}}$  est    déterminée de manière unique par les actions de 
$\pi_{1}(\eta^{I},\ov{\eta^{I}})$ et   des morphismes de Frobenius partiels. Cela résulte du \lemref{lem-Dr-intro}
mais voici une autre fa\c con de le voir  (pour plus de détails on renvoie au chapitre \ref{para-sous-faisceaux-Frob-partiels}). 
 Suivant \cite{drinfeld78}, on va définir 
   un groupe    $\on{FWeil}(\eta^{I},\ov{\eta^{I}})$ 
       \begin{itemize}
     \item qui   est une  extension de $\Z^{I}$ par   $\on{Ker}(\pi_{1}(\eta^{I},\ov{\eta^{I}})\to \wh \Z)$, 
     \item et qui, lorsque  $I$ est un   singleton, 
 s'identifie au groupe de   Weil   usuel 
   $\on{Weil} (\eta,\ov{\eta})=\pi_{1}(\eta,\ov{\eta})\times_{\wh \Z}\Z$.      
       \end{itemize}
       
   On note  $F^{I}$ le corps des fonctions de $X^{I}$, $(F^{I})^{\mr{perf}}$ son 
  perfectisé et $\ov{F^{I}}$ la clôture algébrique de $F^{I}$ telle que 
  $\ov{\eta^{I}}=\on{Spec}(\ov{F^{I}})$. 
   On définit alors 
          \begin{gather*}\on{FWeil}(\eta^{I},\ov{\eta^{I}})=
       \big\{\varepsilon \in \on{Aut}_{\ov\Fq}((\ov{F^{I}})), \exists (n_{i})_{i\in I}\in \Z^{I}, \restr{\varepsilon}{(F^{I})^{\mr{perf}}}=\prod_{i\in I}(\Frob_{\{i\}})^{n_{i}}\big\} 
  .\end{gather*}
Le choix de 
  $\mf{sp}$ fournit 
  une inclusion   $\ov{F}\otimes_{\ov\Fq} \cdots 
 \otimes_{\ov\Fq}\ov{F}\subset \ov{F^{I}} 
  $.  Par restriction des automorphismes, on  en déduit   un 
   morphisme  surjectif 
    \begin{gather}\label{morph-Weil-I-intro}\on{FWeil}(\eta^{I},\ov{\eta^{I}})\to  \on{Weil} (\eta,\ov{\eta}) ^{I}\end{gather} (dépendant du choix de  $\mf{sp}$). L'énoncé de la \propref{cor-action-Hf-intro}  se reformule  alors en disant que l'action naturelle de $\on{FWeil}(\eta^{I},\ov{\eta^{I}})$ 
    sur $\Big( \varinjlim _{\mu}\restr{\mc H _{N, I, W}^{0,\leq\mu,E}}{\ov{\eta^{I}}}\Big)^{\mr{Hf}}$
    se factorise par le morphisme  \eqref{morph-Weil-I-intro}, et même à travers 
    $\pi_{1}(\eta,\ov{\eta})^{I}$.  
      \end{rem}

 \begin{prop}\label{prop2admise}   (\corref{bijectivite-Hecke-fini}). La restriction de l'homomorphisme $\on{\mf{sp}}^{*}$ de \eqref{sp*-sans-Hf-intro} aux parties Hecke-finies est  un isomorphisme 
   \begin{gather}\label{isom-sp*-avec-Hf-intro} 
  \Big( \varinjlim _{\mu}\restr{\mc H _{N, I, W}^{0,\leq\mu,E}}{\Delta(\ov{\eta})}\Big)^{\mr{Hf}} \isor{\on{\mf{sp}}^{*}}
 \Big(  \varinjlim _{\mu}\restr{\mc H _{N, I, W}^{0,\leq\mu,E}}{\ov{\eta^{I}}}\Big)^{\mr{Hf}}. \end{gather} 
\end{prop}
  
  Pour la preuve de la \propref{prop2admise}   le lecteur peut lire dès à présent, s'il le souhaite, les énoncés et les démonstrations  des propositions  \ref{surjectivite-Hecke-fini} et \ref{injectivite-sp}   et du \corref{bijectivite-Hecke-fini}.

  Les propositions  \ref{cor-action-Hf-intro}  et \ref{prop2admise} permettent de définir maintenant les $E$-espaces vectoriels $H_{I,W}$  (on omet $N$ dans la notation $ H_{I,W}$ pour limiter la taille des diagrammes dans le paragraphe suivant). Dans la définition suivante on utilise le membre de gauche 
  de \eqref{isom-sp*-avec-Hf-intro} 
  car il est plus canonique et ne dépend pas du choix de $\ov{\eta^{I}}$ et $\on{\mf{sp}}$. 
   
   \begin{defi}
    On définit $H_{I,W}$ comme le {\it membre de gauche} de \eqref{isom-sp*-avec-Hf-intro}. 
  \end{defi}
 
 L'action de $\on{Gal}(\ov F/F)^{I}=\pi_{1}(\eta,\ov\eta)^{I}$ sur 
 $H_{I,W}$ 
 ne dépend pas du choix de $\ov{\eta^{I}}$ et $\on{\mf{sp}}^{*}$. 
 En effet on peut reformuler ce qui précède en disant que d'après la \propref{cor-action-Hf-intro} 
      on peut trouver  
 \begin{itemize}
 \item une réunion croissante 
 (indexée par $\lambda\in \N$) 
 de sous-$\mc O_{E}$-faisceaux constructibles 
 $\mf F_{\lambda}\subset \varinjlim _{\mu}\restr{\mc H _{N, I, W}^{0,\leq\mu,E}}{\eta^{I}}$ stables par les morphismes de Frobenius partiels
 (auxquels s'applique donc le lemme de Drinfeld) 
 \item une suite décroissante d'ouverts  denses $U_{\lambda}\subset X\sm N$ tels que 
 $\mf F_{\lambda}$ se prolonge en un faisceau lisse sur $(U_{\lambda})^{I}$
 \end{itemize}
 de sorte que $\bigcup _{\lambda\in \N}  \restr{\mf F_{\lambda}}{\ov{\eta^{I}}} 
 =
  \Big(  \varinjlim _{\mu}\restr{\mc H _{N, I, W}^{0,\leq\mu,E}}{\ov{\eta^{I}}}\Big)^{\mr{Hf}}$. 
   Alors la \propref{prop2admise}   
   implique   que  le morphisme naturel 
  \begin{gather}\label{mor-HIW-Flambda-intro}
  H_{I,W}= \Big( \varinjlim _{\mu}\restr{\mc H _{N, I, W}^{0,\leq\mu,E}}{\Delta(\ov{\eta})}\Big)^{\mr{Hf}} \to \bigcup _{\lambda\in \N}  \restr{\mf F_{\lambda}}{\Delta(\ov{\eta})} \end{gather}
(qui vient de la lissité de $\mf F_{\lambda}$ sur $(U_{\lambda})^{I}\ni 
\Delta(\ov{\eta})$)
  est un isomorphisme. Or l'action de $\on{Gal}(\ov F/F)^{I}$   sur le membre de droite de \eqref{mor-HIW-Flambda-intro}, qui est donnée par le lemme de Drinfeld, ne dépend pas  du choix de $\ov{\eta^{I}}$ et $\on{\mf{sp}}$, et donc l'action  de $\on{Gal}(\ov F/F)^{I}$ sur le membre de gauche n'en dépend pas non plus. 
  
 \begin{rem}
Dans cet article nous montrons seulement que  $H_{I,W}$ est une limite inductive de $E$-espaces vectoriels de dimension finie munis de représentations continues de $\on{Gal}(\ov F/F)^{I}$. En fait Cong Xue a montré dans \cite{these-cong} que 
 $H_{I,W}$ est de dimension finie. 
 Les arguments de notre article nous permettent de montrer notre résultat principal en nous passant de ce résultat, qui est de démonstration délicate, et est seulement disponible dans le cas des groupes déployés pour le moment. 
  \end{rem}
 
   Pour toute application  $\zeta:I\to J$, l'isomorphisme de coalescence 
   \eqref{intro-isom-coalescence} respecte trivialement les  parties Hecke-finies et induit donc  un isomorphisme 
    \begin{gather} \label{isom-chi-zeta-intro}
    H_{I,W}= \Big( \varinjlim _{\mu}\restr{\mc H _{N, I, W}^{0,\leq\mu,E}}{\Delta(\ov{\eta})}\Big)^{\mr{Hf}}
    \isor{\chi_{\zeta}} 
    \Big( \varinjlim _{\mu}\restr{\mc H _{N, J,W^{\zeta}}^{0,\leq\mu,E}}{\Delta(\ov{\eta})}\Big)^{\mr{Hf}}=H_{J,W^{\zeta}}
    \end{gather}
    où $\Delta$ désigne le morphisme diagonal $X\to X^{I}$ ou $X\to X^{J}$. 
    
     On  note 
    \begin{gather}\label{def-chi-zeta-intro} \chi_{\zeta}:    H_{I,W} \isom H_{J,W^{\zeta}}\end{gather}  l'isomorphisme \eqref{isom-chi-zeta-intro} ci-dessus. Il 
     est $\on{Gal}(\ov F/F)^{J}$-équivariant,  où $\on{Gal}(\ov F/F)^{J}$ agit sur le membre de gauche par le   morphisme diagonal  
\begin{gather}\nonumber 
\on{Gal}(\ov F/F)^{J}\to \on{Gal}(\ov F/F)^{I},  \ (\gamma_{j})_{j\in J}\mapsto (\gamma_{\zeta(i)})_{i\in I}. 
\end{gather}
En effet, si $\Delta_{\zeta}:X^{J}\to X^{I}$ est le morphisme diagonal 
    \eqref{morph-giad-X-intro} et si la suite  $(\mf F_{\lambda})_{\lambda\in \N}$ est comme ci-dessus relativement à $I$ et $W$, alors la suite  $(\Delta_{\zeta}^{*}(\mf F_{\lambda}))_{\lambda\in \N}$ vérifie les mêmes propriétés relativement à $J$ et $W^{\zeta}$, donc 
       $$\chi_{\zeta}: H_{I,W}=\bigcup _{\lambda\in \N}  \restr{\mf F_{\lambda}}{\Delta(\ov{\eta})} =\bigcup _{\lambda\in \N}  \restr{\Delta_{\zeta}^{*}(\mf F_{\lambda})}{\Delta(\ov{\eta})}=H_{J,W^{\zeta}}
       $$ est $\on{Gal}(\ov F/F)^{J}$-équivariant. 
    
     \begin{prop}\label{prop-a-b-c}  (proposition  \ref{prop-a-b-texte}) 
            Les $H_{I,W}$ vérifient   les 
  propriétés   suivantes : 
        \begin{itemize}
    \item[] {\bf a)} pour tout ensemble fini   $I$,      $$W\mapsto  
    H_{I,W},  \ \ u\mapsto \mc H(u)$$  est un foncteur  $E$-linéaire  de la  catégorie des  représentations $E$-linéaires de dimension finie de  $(\wh G)^{I}$ vers la  catégorie des  limites inductives de représentations $E$-linéaires continues de dimension finie de     $\on{Gal}(\ov F/F)^{I}$,    
              \item[] {\bf b)} pour toute application   $\zeta: I\to J$, 
 on possède  un isomorphisme 
      \begin{gather}\nonumber
                 \chi_{\zeta}: H_{I,W}\isom 
 H_{J,W^{\zeta}},\end{gather} 
 qui est 
 \begin{itemize}
 \item  fonctoriel en   $W$, où  $W$ est une  représentation de $(\wh G)^{I}$ et  $W^{\zeta}$ désigne la   représentation de $(\wh G)^{J}$ sur  $W$ obtenue en composant avec le  morphisme  diagonal $$  (\wh G)^{J}\to (\wh G)^{I}, (g_{j})_{j\in J}\mapsto (g_{\zeta(i)})_{i\in I} $$ 
 \item $\on{Gal}(\ov F/F)^{J}$-équivariant, où $\on{Gal}(\ov F/F)^{J}$ agit sur le membre de gauche par le   morphisme diagonal  
 \begin{gather}\nonumber 
 \on{Gal}(\ov F/F)^{J}\to \on{Gal}(\ov F/F)^{I},  \ (\gamma_{j})_{j\in J}\mapsto (\gamma_{\zeta(i)})_{i\in I}, 
\end{gather}
 \item   et compatible avec la  composition, c'est-à-dire   que pour 
 $I\xrightarrow{\zeta} J\xrightarrow{\eta} K$ on a 
 $\chi_{\eta\circ \zeta}=\chi_{\eta}\circ\chi_{\zeta}$,
  \end{itemize}
     \item[] {\bf c)} pour $I=\emptyset$ et  $W=\mbf  1$, on a un isomorphisme     \begin{gather}\nonumber
       H_{\emptyset,\mbf  1}=C_{c}^{\mr{cusp}}(G(F)\backslash G(\mb A)/K_{N}\Xi,E). \end{gather}
    \end{itemize}
    
  Par ailleurs les $H_{I,W}$ sont des modules sur 
  $C_{c}(K_{N}\backslash G(\mb A)/K_{N},E)$, de fa\c con compatible  avec 
les propriétés a), b), c) ci-dessus. 
  \end{prop}

\noindent {\bf Démonstration de la \propref{prop-a-b-c}.  }    Les propriétés a) et b) ont déjà été expliquées. 
   On remarque qu'en appliquant   b) à l'application évidente  $\zeta_{\emptyset}: \emptyset \to \{0\}$, on obtient un  isomorphisme \begin{gather}
    \label{isom-chi-singleton-empty-intro}\chi_{\zeta_{\emptyset}}: H_{\emptyset,\mbf 1}\isom 
 H_{\{0\},\mbf 1} \end{gather} 
que l'on connaissait  déjà par la \remref{rem-cht-W1} (appliquée à $I=\{0\}$).  
 
La propriété c) résulte du fait que la partie Hecke-finie de 
\eqref{I-vide-intro} est formée exactement des formes automorphes cuspidales, c'est-à-dire que 
$$\Big(C_{c} (G(F)\backslash G(\mb A)/K_N \Xi,E)\Big)^{\mr{Hf}}=C_{c}^{\rm{cusp}}(G(F)\backslash G(\mb A)/K_N \Xi,E).$$

\noindent {\bf Preuve de $\supset$. } Toute fonction cuspidale est Hecke-finie  
 car  le $\mc O_{E}$-module $ C_{c}^{\rm{cusp}}(G(F)\backslash G(\mb A)/K_N \Xi,\mc O_{E}) $
  est   de type fini et  stable par tous les opérateurs $T(f)$ pour $f\in C_{c}(K_{N}\backslash G(\mb A)/K_{N},\mc O_{E})$. 
 
 \noindent {\bf Preuve de $\subset$. }
 On raisonne par l'absurde et on suppose qu'une  fonction Hecke-finie $f$  n'est pas  cuspidale. Il existe alors un parabolique $P\subsetneq G$, de Levi $M$ et radical unipotent $U$,  tel  que le terme constant $f_{P}: g\mapsto \int_{U(F)\backslash U(\mb A)}f(ug)$ soit non nul.   Soit $v$ une  place de $X\sm N$. 
       Comme l'anneau des représentations (de dimension finie) de   $\wh M$ est un module de type fini sur  l'anneau des représentations de $\wh G$, $f_{P}$ est également Hecke-finie relativement aux opérateurs de  Hecke pour $M$ en $v$. 
       Ceux-ci comprennent comme cas particuliers les translations par les éléments de $Z_{M}(F_{v})$. On possède une application degré (relativement à 
       $M/Z$), de $U(\mb A)M(F)\backslash G(\mb A)/K_{N}\Xi$ à valeurs dans un $\Z$-module libre de type fini, sur lequel $Z_{M}(F_{v})$ agit par des translations non triviales. Or le support de $f_{P}$ est inclus dans le translaté d'un cône dans 
       ce $\Z$-module libre de type fini, et cela contredit le fait que $f_{P}$ appartient à 
       un espace vectoriel de dimension finie stable par  $Z_{M}(F_{v})$.
     On renvoie à la \propref{prop-cusp-hecke-finies} pour une preuve plus détaillée de c).   \cqfd

         Dans le prochain paragraphe,  nous expliquerons l'idée de la preuve du \thmref{intro-thm-ppal} à l'aide des propriétés   a), b), c) de la \propref{prop-a-b-c}. Cette esquisse sera complétée par la preuve de la compatibilité avec  l'isomorphisme de Satake dans le paragraphe  
      \ref{subsection-intro-decomp}. 
      
    \subsection{Idée de la preuve  du  \thmref{intro-thm-ppal}  partir de la  \propref{prop-a-b-c}} 
    \label{intro-idee-heurist}
    
      L'idée  se résume ainsi :  grâce à \eqref{isom-chi-singleton-empty-intro} et au   c) de la \propref{prop-a-b-c}, on a 
    $$H_{\{0\},\mbf 1}=C_{c}^{\mr{cusp}}(G(F)\backslash G(\mb A)/K_{N}\Xi,E).$$
   Pour obtenir la décomposition \eqref{intro1-dec-canonique}  il est donc  équivalent de construire (quitte à augmenter $E$) 
   une  décomposition canonique 
     \begin{gather}\label{intro3-dec-canonique}
  H_{\{0\},\mbf 1}=\bigoplus_{\sigma}
 \mf H_{\sigma}.\end{gather}
Cette dernière sera obtenue par  décomposition spectrale  d'une famille commutative d'endomorphismes de $H_{\{0\},\mbf 1}$, appelés opérateurs d'excursion, que nous allons construire et étudier à l'aide 
  des  propriétés a) et  b) de la \propref{prop-a-b-c}.

     Soit  $I$   un ensemble fini et  $W$ une représentation  $E$-linéaire de dimension finie de  
    $(\wh G)^{I}$. 
    On note $\zeta_{I}:I\to \{0\}$ l'application évidente, si bien que 
    $ W^{\zeta_{I}}$ est simplement  $W$ muni de l'action diagonale  de $\wh G$. Soit 
$x: \mbf 1\to W^{\zeta_{I}}$ et  $\xi :  W^{\zeta_{I}}\to  \mbf 1$
des morphismes de  représentations de  $\wh G$ (autrement dit  $x\in W$ et  $\xi\in W^{*}$  sont invariants sous  l'action diagonale de  $\wh G$). Soit   $(\gamma_{i})_{i\in I}\in \on{Gal}(\ov F/F)^{I}$. 
 
 \begin{defi}On définit l'opérateur  
  \begin{gather*}S_{I,W,x,\xi,(\gamma_{i})_{i\in I}}\in 
  \on{End}(H_{\{0\},\mbf 1})  \end{gather*}
    comme la composée 
  \begin{gather}\label{excursion-def-intro}
  H_{\{0\},\mbf  1}\xrightarrow{\mc H(x)}
 H_{\{0\},W^{\zeta_{I}}}\isor{\chi_{\zeta_{I}}^{-1}} 
  H_{I,W}
  \xrightarrow{(\gamma_{i})_{i\in I}}
  H_{I,W} \isor{\chi_{\zeta_{I}}} H_{\{0\},W^{\zeta_{I}}}  
  \xrightarrow{\mc H(\xi)} 
  H_{\{0\},\mbf  1}.
    \end{gather}\end{defi}
    
    Cet opérateur sera appelé ``opérateur d'excursion''. En paraphrasant  \eqref{excursion-def-intro} il est la composée  
    \begin{itemize}
    \item d'un  opérateur de création associé à $x$, dont l'effet est de créer des pattes en le même point (générique) de la courbe, 
    \item d'une action de Galois, 
    qui promène les pattes sur la courbe indépendamment les unes des autres, puis les ramène au même point (générique) de la courbe, 
       \item d'un  opérateur d'annihilation associé à $\xi$, qui annihile les pattes. 
         \end{itemize}

    La remarque suivante  propose une description conjecturale
  des $H_{I,W}$ qui permet de mieux comprendre, de fa\c con heuristique, le sens des opérateurs d'excursion. Bien sûr cette description conjecturale n'intervient nulle part dans les raisonnements.

 \begin{rem}  \label{descr-conj-HIW}   
 On conjecture qu'il existe un ensemble fini $\Sigma$ (dépendant de $N$) 
 de paramètres de Langlands semi-simples (bien déterminés à conjugaison près), et que, quitte à augmenter $E$,  ils soient définis sur $E$ et   qu'on possède pour tout $\sigma\in \Sigma$
     une  représentation  $E$-linéaire   $A_{\sigma}$ 
      du centralisateur   $S_{\sigma}$  de l'image de $\sigma$ dans $\wh G$  (triviale sur $Z_{\wh G}$), de telle sorte que 
  pour tout    $I$ et  $W$  \begin{gather}\label{intro-dec-I-W}
 H_{I,W}\overset{?}{=}\bigoplus_{\sigma\in \Sigma}
\Big(A_{\sigma} \otimes _{E}W_{\sigma^{I}}\Big)^{S_{\sigma}}, \end{gather}
 où 
$W_{\sigma^{I}}$ désigne la   représentation de  
$\on{Gal}(\ov F/F)^{I}$     obtenue en composant la   représentation $W$ avec le  morphisme   
$\sigma^{I}: \on{Gal}(\ov F/F)^{I} \to (\wh G(E))^{I}$. 
De plus $A_{\sigma}$ doit être un module sur 
$C_{c}(K_{N}\backslash G(\mb A)/K_{N},E)$, 
et \eqref{intro-dec-I-W} doit être un isomorphisme de 
$C_{c}(K_{N}\backslash G(\mb A)/K_{N},E)$-modules. 
Dans le cas particulier où $I= \emptyset$  et $W=1$, 
\eqref{intro-dec-I-W} doit être la décomposition  
\eqref{intro3-dec-canonique}  et   on doit avoir 
$ \mf H_{\sigma}=(A_{\sigma})^{S_{\sigma}}$. 
  
   Ces conjectures sont bien connues des experts, par extrapolation des 
conjectures  d'Arthur \cite{arthur-ast} et de  Kottwitz \cite{kottwitz3}    sur les multiplicités dans les espaces de formes automorphes et dans  la cohomologie des variétés de  Shimura, et grâce à l'égalité montrée par Cong Xue \cite{these-cong} entre une variante de $H_{I,W}$ (qui lui est conjecturalement égale) et la ``cohomologie cuspidale''. 
Dans le cas de $GL_{r}$ on s'attend à ce que 
$\Sigma$ soit l'ensemble des représentations irréductibles de dimension $r$ de $\pi_{1}(X\sm N, \ov\eta)$ et que pour tout 
$\sigma\in \Sigma$, 
$S_{\sigma}=\mathbb G_{m}=Z_{\wh G}$  et $A_{\sigma}=(\pi_{\sigma})^{K_{N}}$ où $\pi_{\sigma}$ est la représentation automorphe cuspidale correspondant à $\sigma$ (voir \cite{laurent-inventiones} et  la conjecture 2.35 de \cite{var}). En général si $\sigma$ est associé à un paramètre d'Arthur elliptique  $\psi$ (voir 
 le paragraphe \ref{subsubsection-Arthur} ci-dessous), $A_{\sigma}$ devrait être induit d'une représentation de dimension finie du sous-groupe de  $S_{\sigma}$ engendré par le centralisateur de $\psi$ et par le sous-groupe diagonal  $\mathbb{G}_{m}\subset SL_{2}$ (parce que nous considérons seulement la cohomologie en degré $0$).

On conjecture de plus que   \eqref{intro-dec-I-W} est fonctoriel en   $W$ et que pour toute application $\zeta:I\to J$ il entrelace   $\chi_{\zeta}$ avec  
  \begin{gather}\nonumber \Id:  \bigoplus_{\sigma} \Big(A_{\sigma} \otimes _{E}W_{\sigma^{I}}\Big)^{S_{\sigma}}\to
 \bigoplus_{\sigma} 
 \Big(A_{\sigma} \otimes _{E}(W^{\zeta})_{\sigma^{J}}\Big)^{S_{\sigma}}
   \end{gather}
   (comme $W_{\sigma^{I}}$ et $(W^{\zeta})_{\sigma^{J}}$ sont tous les deux égaux à $W$ comme $E$-espaces vectoriels, le morphisme $\Id$ a bien un sens et est   $\on{Gal}(\ov F/F)^{J}$-équivariant). 
 Sous ces hypothèses,   la composée \eqref{excursion-def-intro} (qui   définit  $S_{I,f,(\gamma_{i})_{i\in I}}$) 
 agit sur   
    $ \mf H_{\sigma}= (A_{\sigma})^{S_{\sigma}} \subset  H_{\{0\},\mbf 1}$ par   la composée  
      \begin{gather} \nonumber (A_{\sigma})^{S_{\sigma}}
 \xrightarrow{\Id_{A_{\sigma}}\otimes x} 
\Big(A_{\sigma} \otimes _{E}W_{\sigma^{I}}\Big)^{S_{\sigma}}
\xrightarrow{\on{Id}_{A_{\sigma}}\otimes (\gamma_{i})_{i\in I}}
\Big(A_{\sigma} \otimes _{E}W_{\sigma^{I}}\Big)^{S_{\sigma}}
\xrightarrow{\Id_{A_{\sigma}}\otimes \xi} 
 (A_{\sigma})^{S_{\sigma}}
    \end{gather}
   c'est-à-dire  par le  produit par le scalaire  $\s{\xi, (\sigma(\gamma_{i}))_{i\in I} \cdot x }=f((\sigma(\gamma_{i}))_{i\in I} )$. 
      
   La conjecture \eqref{intro-dec-I-W} n'est pas démontrée mais suivant une idée de Drinfeld on peut montrer que les propriétés a) et b) de la \propref{prop-a-b-c} 
   suffisent à  impliquer   une décomposition 
 assez proche de \eqref{intro-dec-I-W} 
   (mais plus difficile à énoncer, car rempla\c cant la donnée de $\Sigma$ et des $A_{\sigma}$ par celle d'un ``$\mc O$-module sur le champ des paramètres de Langlands'').    
   
          D'après l'observation \ref{observation-f-W-x-xi} et  la \propref{intro-Xi-n} ci-dessous, la connaissance de   
   $\s{\xi, (\sigma(\gamma_{i}))_{i\in I} \cdot x }$ (pour  $I,W,x,\xi$ et 
   $(\gamma_{i})_{i\in I}$  arbitraires) détermine  $\sigma$ à conjugaison près.  Par conséquent si l'on croit à la conjecture  
    \eqref{intro-dec-I-W}, la décomposition \eqref{intro3-dec-canonique} s'obtient par  diagonalisation simultanée des  opérateurs d'excursion. 
  En fait nous ne savons pas montrer qu'ils sont diagonalisables, et nous obtiendrons  la décomposition \eqref{intro3-dec-canonique} 
  par décomposition spectrale, c'est-à-dire ``trigonalisation simultanée'' des opérateurs d'excursion.

     De plus on devine que les opérateurs d'excursion  vérifient les propriétés 
   énoncées dans le lemme suivant.     
    
Cette remarque   était heuristique et à partir de maintenant on oublie la conjecture \eqref{intro-dec-I-W} 
   (sauf dans la  remarque \ref{rem-intro-heuristique-avant-prop}). 
       \end{rem}

   Le lemme suivant va résulter des propriétés a) et  b) de la \propref{prop-a-b-c}. 
   \begin{lem} 
   \label{lem-intro-ptes-SIW}  (\lemref{lem-ptes-SIW})
   Les  opérateurs d'excursion $S_{I,W,x,\xi,(\gamma_{i})_{i\in I}}$ vérifient 
    les propriétés suivantes : 
     \begin{gather}
     \label{SIW-p0-intro}      S_{I,W,x,{}^{t} u(\xi'),(\gamma_i)_{i\in I}}=S_{I,W',u(x),\xi',(\gamma_i)_{i\in I}} 
 \end{gather}
où $u:W\to W'$ est un   morphisme $(\wh G)^{I}$-équivariant et  $x\in W$ et $\xi'\in (W')^{*}$ sont $\wh G$-invariants,  
  \begin{gather}   \label{SIW-p1-intro}
     S_{J,W^{\zeta},x,\xi,(\gamma_j)_{j\in J}}=S_{I,W,x,\xi,(\gamma_{\zeta(i)})_{i\in I}},
     \\
      \label{SIW-p2-intro}
 S_{I_{1}\cup I_{2},W_{1}\boxtimes W_{2},x_{1}\boxtimes x_{2},\xi_{1}\boxtimes \xi_{2},(\gamma^{1}_i)_{i\in I_{1}}\times (\gamma^{2}_i)_{i\in I_{2}}}= S_{I_{1},W_{1},x_{1},\xi_{1},(\gamma^{1}_i)_{i\in I_{1}}}\circ 
S_{I_{2},W_{2},x_{2},\xi_{2},(\gamma^{2}_i)_{i\in I_{2}}}, 
\\ 
\label{SIW-p3-intro} 
S_{I,W,x,\xi,(\gamma_i(\gamma'_i )^{-1}\gamma''_i)_{i\in I}}=
    S_{I\cup I \cup I,W\boxtimes W^{*}\boxtimes W,\delta_{W} \boxtimes x,
    \xi \boxtimes \on{ev}_{W},
    (\gamma_i)_{i\in I} \times (\gamma'_i)_{i\in I} \times (\gamma''_i)_{i\in I}
    }
             \end{gather}
  où la plupart des  notations sont évidentes, 
    $I_{1}\cup I_{2}$ et  $I\cup I\cup I$ désignent des réunions disjointes, et  $\delta_{W}: \mbf 1\to W\otimes W^{*}$ et $\on{ev}_{W}:W^{*} \otimes W \to \mbf 1$ sont les  morphismes naturels.      
  \end{lem}
    La démonstration est très simple (elle utilise seulement les propriétés a) et b) de la \propref{prop-a-b-c}). Le lecteur peut lire dès à présent, s'il le souhaite, la preuve du \lemref{lem-ptes-SIW}.

      On note 
$\mc B \subset \on{End}_{C_{c}(K_{N}\backslash G(\mb A)/K_{N},E)}( H_{\{0\},\mbf 1}) $ 
la sous-algèbre engendrée par tous les  opérateurs d'excursion $S_{I,W,x,\xi,(\gamma_{i})_{i\in I}}$.
En vertu de   \eqref{SIW-p2-intro},  $\mc B$ est commutative. 
   
   Dans la suite de l'introduction on considère $\wh G$ comme un schéma en groupes défini sur $E$. 
     
 \begin{observation} \label{observation-f-W-x-xi} Les fonctions \begin{gather}\label{intro-def-f}
   f: (g_{i})_{i\in I}\mapsto \s{\xi, (g_{i})_{i\in I}\cdot x}\end{gather} que l'on obtient en faisant varier  $W$, 
 $x$, et  $\xi $ sont exactement les fonctions régulières sur  le quotient grossier     de  $(\wh G)^{I}$ par  translation à gauche et à droite par 
 $\wh G$ diagonal, que l'on notera $\wh G\backslash (\wh G)^{I}/\wh G$. 
 \end{observation}
 
 \begin{lem}\label{f-W,x,xi} (\lemref{f-W,x,xi-texte})
L'opérateur $S_{I,W,x,\xi,(\gamma_{i})_{i\in I}}$ dépend seulement de   $I$, $f$, et $(\gamma_{i})_{i\in I}$, où $f$ est donnée par \eqref{intro-def-f}. 
\end{lem}
\noindent Cela résulte facilement de \eqref{SIW-p0-intro}. 
Le lecteur peut lire dès à présent, s'il le souhaite, la démonstration du \lemref{f-W,x,xi-texte}. 

Le lemme précédent permet d'introduire la notation suivante. 
  
 \begin{defi}
 Pour toute fonction  $f\in  \mc O(\wh G\backslash (\wh G)^{I}/\wh G)$ on pose  
 \begin{gather} S_{I,f,(\gamma_{i})_{i\in I}}=S_{I,W,x,\xi,(\gamma_{i})_{i\in I}}\in \mc B \end{gather}
  où $W,x,\xi$ sont tels que  $f$ satisfasse  
 \eqref{intro-def-f}. 
 \end{defi}

 Cette nouvelle  notation va permettre, dans les assertions (i) à (iv) de la  proposition suivante, de formuler de fa\c con plus synthétique les propriétés     \eqref{SIW-p1-intro}, \eqref{SIW-p2-intro} et \eqref{SIW-p3-intro}. 
 
\begin{rem}
\label{rem-intro-heuristique-avant-prop}
 Si on se place dans l'heuristique de la \remref{descr-conj-HIW}, on peut deviner facilement les propriétés énoncées dans la proposition suivante en remarquant que d'après  la conjecture    \eqref{intro-dec-I-W}, $S_{I,f,(\gamma_{i})_{i\in I}}$ devrait agir sur    $\mf H_{\sigma}$ par  multiplication par le scalaire  $f((\sigma(\gamma_{i}))_{i\in I})$.  
\end{rem}

 \begin{prop} \label{prop-SIf-i-ii-iii} (\propref{art-prop-SIf-i-ii-iii}) 
 Les  opérateurs d'excursion  
   $S_{I,f,(\gamma_{i})_{i\in I}}$ vérifient les propriétés suivantes:   
  \begin{itemize}
  \item [] (i) pour tout   $I$ et 
 $(\gamma_{i})_{i\in I}\in  \on{Gal}(\ov F/F)^{I}$, 
  $$f\mapsto 
  S_{I,f,(\gamma_{i})_{i\in I}}$$ est un  morphisme 
  d'algèbres  commutatives  $\mc O(\wh G\backslash (\wh G)^{I}/\wh G)\to \mc B$, 
  \item [] (ii) pour toute application 
  $\zeta:I\to J$, toute fonction  $f\in \mc O(\wh G\backslash (\wh G)^{I}/\wh G)$ et  tout 
  $(\gamma_{j})_{j\in J}\in  \on{Gal}(\ov F/F)^{J}$, on a 
  $$S_{J,f^{\zeta},(\gamma_{j})_{j\in J}}=S_{I,f,(\gamma_{\zeta(i)})_{i\in I}}$$
   où $f^{\zeta}\in \mc O(\wh G\backslash (\wh G)^{J}/\wh G)$ est définie par    $$f^{\zeta}((g_{j})_{j\in J})=f((g_{\zeta(i)})_{i\in I}),$$
   \item [] (iii) 
  pour tout   $f\in \mc O(\wh G\backslash (\wh G)^{I}/\wh G)$
  et  $(\gamma_{i})_{i\in I},(\gamma'_{i})_{i\in I},(\gamma''_{i})_{i\in I}\in  \on{Gal}(\ov F/F)^{I}$ on a     $$S_{I\cup I\cup I,\wt f,(\gamma_{i})_{i\in I}\times (\gamma'_{i})_{i\in I}\times (\gamma''_{i})_{i\in I}}=
  S_{I,f,(\gamma_{i}(\gamma'_{i})^{-1}\gamma''_{i})_{i\in I}}$$
   où  $I\cup I\cup I$ est une réunion disjointe et 
   $\wt f\in \mc O(\wh G\backslash (\wh G)^{I\cup I\cup I}/\wh G)$ est définie  par  
   $$\wt f((g_{i})_{i\in I}\times (g'_{i})_{i\in I}\times (g''_{i})_{i\in I})=f((g_{i}(g'_{i})^{-1}g''_{i})_{i\in I}).$$
      \item [] (iv) pour tout $I$ et tout $f$, le morphisme 
      \begin{gather}\label{mor-excursion-f}\on{Gal}(\ov F/F)^{I}\to \mc B, \ \ (\gamma_{i})_{i\in I}\mapsto S_{I,f,(\gamma_{i})_{i\in I}}\end{gather} est continu, avec $\mc B$ munie de la topologie $E$-adique. 
                \end{itemize}
         \end{prop}

      \noindent {\bf Esquisse de démonstration.}    
         Cela résulte formellement  de la proposition  \ref{prop-a-b-c}. En effet, on déduit   (ii) de \eqref{SIW-p1-intro}. 
      Pour montrer   (i) on invoque 
\eqref{SIW-p2-intro} avec $I_{1}=I_{2}=I$,  et on applique   \eqref{SIW-p1-intro} à l'application évidente $\zeta :I\cup I\to I$. La propriété  (iii) résulte de   \eqref{SIW-p3-intro}, après avoir remarqué que 
    $$\s{ \xi \boxtimes \on{ev}_{W}, \big((g_{i})_{i\in I} 
  \boxtimes  (g'_{i})_{i\in I} \boxtimes  (g''_{i})_{i\in I} \big)\cdot ( \delta_{W} \boxtimes x)}
  =\s{\xi, (g_{i}(g'_{i})^{-1}g''_{i})_{i\in I} \cdot x}.$$
  Enfin (iv) résulte du fait que $H_{I,W}$ est une limite inductive de représentations continues de dimension finie de $(\on{Gal}(\ov F/F))^{I}$.
        \cqfd

On ne sait pas  si $\mc B$ est réduite. On possède néanmoins  une décomposition spectrale 
  (c'est-à-dire une décomposition en espaces propres généralisés, ou ``espaces caractéristiques'') 
  \begin{gather}\label{dec-nu-intro-intro}
  H_{\{0\},\mbf 1} =\bigoplus_{\nu} \mf H_{\nu}\end{gather}
 où dans le membre de droite la somme directe est indexée  par les  caractères  $\nu$ de  $\mc B$. Quitte à augmenter $E$ on suppose que tous les caractères de $\mc B$ sont définis sur $E$. 
 
  La proposition suivante permet d'obtenir la décomposition \eqref{intro3-dec-canonique}  à partir de \eqref{dec-nu-intro-intro} en associant à chaque caractère $\nu$ un paramètre de Langlands $\sigma$.

 \begin{prop} \label{intro-Xi-n}  Pour tout caractère  $\nu$ de  $\mc B$ il existe un    morphisme 
       $\sigma:\on{Gal}(\ov F/F)\to \wh G(\Qlbar)$
tel que  
\begin{itemize}
\item [] (C1) $\sigma$ prend ses valeurs  dans  $\wh G(E')$, où $E'$ est une extension finie  de $E$, et il est continu, 
\item [] (C2) $\sigma$ est semi-simple, c'est-à-dire que si 
son image est incluse dans un parabolique elle est incluse dans un Levi associé
(comme $\Qlbar$ est de caractéristique $0$ cela équivaut à dire que 
 l'adhérence  de  Zariski de son  image est réductive \cite{bki-serre}),  
\item[] (C3) 
pour tout  $I$ et    $f\in \mc O(\wh G\backslash (\wh G)^{I}/\wh G)$,  on a  
$$\nu(S_{I,f,(\gamma_{i})_{i\in I}})= f\big((\sigma(\gamma_{i}))_{i\in I} \big).
$$
\end{itemize}
De plus   $\sigma$ est unique à conjugaison près par   $\wh G(\Qlbar)$. 
\end{prop}

   \noindent {\bf Démonstration.}   On renvoie à la preuve de la \propref{Xi-n}  pour quelques détails suplémentaires. 
   La preuve utilise uniquement   la  \propref{prop-SIf-i-ii-iii}. 
    Soit $\nu$ un caractère de   $\mc B$.

Pour tout $n\in \N$ 
  on note  
   $(\wh G)^{n}\modmod \wh G$ le  quotient grossier de  $(\wh G)^{n}$ par l'action  de  $\wh G$  par  conjugaison diagonale, c'est-à-dire   
   $$h.(g_{1},...,g_{n})=(hg_{1}h^{-1},...,hg_{n}h^{-1}).$$ 
   Alors le morphisme 
   \begin{gather*} (\wh G)^{n}\to (\wh G)^{\{0,...,n\}}, (g_{1},...,g_{n})\mapsto 
   (1,g_{1},...,g_{n})\end{gather*}
   induit  un  isomorphisme  $$\beta: (\wh G)^{n}\modmod \wh G\isom \wh G\backslash (\wh G)^{\{0,...,n\}}/\wh G ,$$ d'où un   isomorphisme  d'algèbres   $$\mc O( (\wh G)^{n}\modmod \wh G)\isom \mc O(\wh G\backslash (\wh G)^{\{0,...,n\}}/\wh G ), \ \ f\mapsto f \circ \beta^{-1}. $$ 

On introduit  \begin{align} \nonumber \Theta_{n}^{\nu}: \mc O((\wh G)^{n}\modmod \wh G)&\to 
C(\on{Gal}(\ov F/F)^{n}, E)\\ \nonumber 
f& \mapsto [(\gamma_{1},...,\gamma_{n})\mapsto 
\nu(S_{I,f \circ \beta^{-1},(1,\gamma_{1},...,\gamma_{n})})]\end{align}

La condition (C3) que doit vérifier $\sigma$ se reformule de la fa\c con suivante : pour tout $n$ et pour tout $f\in \mc O((\wh G)^{n}\modmod \wh G)$, 
\begin{gather}\label{rel-sigma-nu}\Theta_{n}^{\nu}(f)=[(\gamma_{1},...,\gamma_{n})\mapsto f((\sigma(\gamma_{1}),...,\sigma(\gamma_{n})))   ]. \end{gather}

On déduit immédiatement  de la  \propref{prop-SIf-i-ii-iii} que la suite 
$(\Theta_{n}^{\nu})_{n\in \N^{*}}$ 
vérifie les propriétés suivantes
\begin{itemize} \item 
pour tout  $n$, $\Theta_{n}^{\nu}$ est un  morphisme d'algèbres, 
                                                 \item  la suite  $(\Theta_{n}^{\nu})_{n\in \N^{*}}$ est fonctorielle par rapport à  toutes les applications entre les ensembles  $\{1,...,n\}$, c'est-à-dire  que  pour $m,n\in \N^{*}$, 
               $$\zeta: \{1,...,m\}\to \{1,...,n\}$$ arbitraire,   
               $f\in \mc O((\wh G)^{m}\modmod \wh G)$ et 
               $(\gamma_{1},...,\gamma_{n })\in \on{Gal}(\ov F/F)^{n }$, 
               on a  
             $$\Theta_{n}^{\nu}( f^{\zeta}) ((\gamma_{j})_{j\in \{1,...,n\}})=
             \Theta_{m}^{\nu}(f)((\gamma_{\zeta(i)})_{i\in \{1,...,m\}})$$ 
            où  $f^{\zeta}\in \mc O((\wh G)^{n}\modmod \wh G)$  est définie par 
               $$f^{\zeta}((g_{j})_{j\in \{1,...,n\}})=f((g_{\zeta(i)})_{i\in \{1,...,m\}}), $$
 \item     pour $n\geq 1$, 
  $f\in \mc O((\wh G)^{n}\modmod \wh G)$  
    et   $(\gamma_{1},...,\gamma_{n+1})\in \on{Gal}(\ov F/F)^{n+1}$ on a  
   $$\Theta_{n+1}^{\nu}( \wh f)(\gamma_{1},...,\gamma_{n+1})=
   \Theta_{n}^{\nu}( f)(\gamma_{1},...,\gamma_{n}\gamma_{n+1}) $$
 où $\wh f\in  \mc O((\wh G)^{n+1}\modmod \wh G)$ est définie par 
   $$\wh f(g_{1},...,g_{n+1})=f(g_{1},...,g_{n}g_{n+1}). $$
\end{itemize}
 Pour justifier la dernière propriété, on applique la  propriété 
   (iii) de la  \propref{prop-SIf-i-ii-iii} à  $$I=\{0,...,n \}, 
   (\gamma_{i})_{i\in I}=(1,\gamma_{1},...,\gamma_{n}), (\gamma'_{i})_{i\in I}=(1)_{i\in I}, (\gamma''_{i})_{i\in I}=(1, ...,1,\gamma_{n+1})$$ et on utilise   (ii) pour supprimer tous les   $1$ sauf le premier dans  
   $(\gamma_{i})_{i\in I}\times (\gamma'_{i})_{i\in I}\times (\gamma''_{i})_{i\in I}$.

On va montrer que ces propriétés de 
 la suite 
$(\Theta_{n}^{\nu})_{n\in \N^{*}}$ 
entraînent l'existence et l'unicité de   $\sigma$ vérifiant (C1), (C2)  et (C3) (c'est-à-dire
\eqref{rel-sigma-nu}). 

Pour $G=GL_{r}$  le résultat est déjà connu: la suite $(\Theta_{n}^{\nu})_{n\in \N^{*}}$  
est déterminée par $\Theta_{1}^{\nu}(\on{Tr})$ 
(qui doit être le caractère de $\sigma$) 
et 
$\Lambda^{r+1}\mr{St}=0$ implique la relation de pseudo-caractère 
d'où l'existence de $\sigma$ par     \cite{taylor}.  On renvoie à la \remref{rem-pseudo-car-Taylor}  pour plus de détails. 

En général on utilise des résultats de  Richardson  \cite{richardson}. 
On dit qu'un  $n$-uplet   $(g_{1},...,g_{n})\in
 \wh G(\Qlbar)^{n}$ est semi-simple si 
 tout parabolique le contenant 
 possède un sous-groupe de Levi associé le contenant. Comme $  \Qlbar$ est de caractéristique $0$ cela équivaut à la condition que 
 l'adhérence de   Zariski  
  $\ov{<g_{1},...,g_{n}>}$ du sous-groupe  $<g_{1},...,g_{n}>$ engendré par   $g_{1},...,g_{n}$ est réductive \cite{bki-serre}. 
 D'après le théorème 3.6 de  \cite{richardson} la $\wh G$-orbite (par conjugaison) de $(g_{1},...,g_{n})$ est fermée dans 
 $(\wh G)^{n}$ si et seulement si $(g_{1},...,g_{n})$ est semi-simple. 
 Donc les points sur $\Qlbar$ du quotient  grossier $(\wh G)^{n}\modmod \wh G$ (qui correspondent aux 
    $\wh G$-orbites fermées définies sur    $\Qlbar$ dans   $(\wh G)^{n}$)  sont   en bijection avec  les classes de conjugaison par  $\wh G(\Qlbar)$  de   $n$-uplets  semi-simples $(g_{1},...,g_{n})\in
 \wh G(\Qlbar)^{n}$.

 Pour tout $n$-uplet 
   $(\gamma_{1},...,\gamma_{n})\in \on{Gal}(\ov F/F)^{n}$ on note  
  $\xi_{n} (\gamma_{1},...,\gamma_{n})$ le point défini sur $\Qlbar$ du 
    quotient  grossier $(\wh G)^{n}\modmod \wh G$ associé au caractère 
   $$\mc O((\wh G)^{n}\modmod \wh G)\to \Qlbar,  \ f\mapsto \Theta_{n}^{\nu}( f)(\gamma_{1},...,\gamma_{n}).$$ 
   On note 
 $\xi_{n}^{\mr{ss}}(\gamma_{1},...,\gamma_{n})$ la  classe de conjugaison de   
  $n$-uplets  semi-simples correspondant à  $\xi_{n}(\gamma_{1},...,\gamma_{n})$ par le résultat de  \cite{richardson} rappelé ci-dessus.  

La relation  \eqref{rel-sigma-nu} équivaut à la condition que pour tout $n$ et pour tout 
$(\gamma_{1},...,\gamma_{n})$, 
$ (\sigma(\gamma_{1}),...,\sigma(\gamma_{n}))\in (\wh G(\Qlbar))^{n}$ 
(qui n'est pas en général semi-simple)
est au-dessus de 
$\xi_{n}(\gamma_{1},...,\gamma_{n})$. 

\noindent {\bf Unicité de $\sigma$ (à conjugaison près).}  On choisit $n$ et $(\gamma_{1},...,\gamma_{n})$ tels que $\sigma(\gamma_{1}),...,\sigma(\gamma_{n})$ engendrent un sous-groupe Zariski dense dans  $\ov{ \mr{Im}(\sigma)}$. 
Comme $\sigma$ est supposé semi-simple, $(\sigma(\gamma_{1}),...,\sigma(\gamma_{n}))$ est semi-simple. 
On fixe 
$(g_{1},...,g_{n})$ dans $\xi_{n}^{\mr{ss}}(\gamma_{1},...,\gamma_{n})$. 
Donc $(\sigma(\gamma_{1}),...,\sigma(\gamma_{n}))$ est conjugué à 
$(g_{1},...,g_{n})$ et quitte à conjuguer $\sigma$ on peut supposer qu'il lui est égal. Alors $\sigma$ est déterminé de fa\c con unique car pour tout $\gamma$, 
$\sigma(\gamma)$  appartient à l'adhérence de Zariski du sous-groupe engendré par $(g_{1},...,g_{n})$ et 
$(g_{1},...,g_{n},\sigma(\gamma)) \in 
\xi_{n+1}^{\mr{ss}}(\gamma_{1},...,\gamma_{n},\gamma)$, donc 
la connaissance de $\xi_{n+1} (\gamma_{1},...,\gamma_{n},\gamma)$ détermine uniquement $\sigma(\gamma)$. 

\noindent {\bf Existence de $\sigma$.}  
  Pour tout $n$ et tout $(\gamma_{1},...,\gamma_{n})\in \on{Gal}(\ov F/F)^{n}$ on choisit $(g_{1},...,g_{n})\in \xi_{n}^{\mr{ss}}(\gamma_{1},...,\gamma_{n})$ (bien défini à conjugaison près).  On choisit alors $n$ et 
   $(\gamma_{1},...,\gamma_{n})\in \on{Gal}(\ov F/F)^{n}$ tels que
       \begin{itemize}
   \item    (H1)
  la  dimension de 
  $\ov{<g_{1},...,g_{n}>}$ est la plus grande  possible 
  \item     (H2) le centralisateur   $C(g_{1},...,g_{n})$ de   $<g_{1},...,g_{n}>$ est le plus petit  possible (dimension minimale puis  nombre de composantes connexes minimal).  
  \end{itemize}
  On fixe  $(g_{1},...,g_{n})\in \xi_{n}^{\mr{ss}}(\gamma_{1},...,\gamma_{n})$ pour le reste de la démonstration et on construit une application   $$\sigma: \on{Gal}(\ov F/F)\to \wh G(\Qlbar)$$  en demandant que pour tout  $\gamma\in \on{Gal}(\ov F/F)$, $\sigma(\gamma)$ est l'unique élément $g$ de  $\wh G(\Qlbar)$ tel que  $(g_{1},...,g_{n},g)\in \xi_{n+1}^{\mr{ss}}(\gamma_{1},...,\gamma_{n},\gamma)$.  
  L'existence et l'unicité de $g$ sont justifiées de la fa\c con suivante. 
  \begin{itemize}
  \item  {\bf A) Existence de $g$} : pour $(h_{1},...,h_{n},h)\in \xi_{n+1}^{\mr{ss}}(\gamma_{1},...,\gamma_{n},\gamma)$, $(h_{1},...,h_{n})$ est forcément semi-simple. En effet $(h_{1},...,h_{n})$ est au-dessus de $\xi_{n}(\gamma_{1},...,\gamma_{n})$ et $(g_{1},...,g_{n})\in \xi_{n}^{\mr{ss}}(\gamma_{1},...,\gamma_{n})$
  donc d'après le théorème 5.2 de \cite{richardson}, $\ov{<h_{1},...,h_{n}>}$ admet un sous-groupe de Levi isomorphe à $\ov{<g_{1},...,g_{n}>}$, or 
  $$\dim(\ov{<h_{1},...,h_{n}>})\leq \dim(\ov{<h_{1},...,h_{n},h>})\leq \dim(\ov{<g_{1},...,g_{n}>})$$ où la deuxième inégalité est assurée 
  par (H1). 
  
  Donc quitte à conjuguer $(h_{1},...,h_{n},h)$ on peut supposer que $(h_{1},...,h_{n})=(g_{1},...,g_{n})$ et on prend alors $g=h$.   
  \item {\bf B) Unicité de $g$} : on a $C(g_{1},...,g_{n},g)\subset C(g_{1},...,g_{n})$ et l'égalité a  lieu par (H2), donc $g$ commute avec $C(g_{1},...,g_{n})$ et comme il était bien déterminé modulo conjugaison par $C(g_{1},...,g_{n})$ il est unique. 
  \end{itemize}
  
Puis on montre que l'application  $\sigma$ que l'on vient de construire est un morphisme de groupes. 
En effet soient $\gamma, \gamma'\in \on{Gal}(\ov F/F)$. Le même argument que dans A) ci-dessus montre qu'il existe $g,g'$ tels que 
\begin{gather}\label{semi-simple-g-g'}(g_{1},...,g_{n},g,g') \in \xi_{n+2}^{\mr{ss}}(\gamma_{1},...,\gamma_{n},\gamma,\gamma').\end{gather}
   Grâce aux  propriétés vérifiées par la suite 
      $(\Theta_{n}^{\nu})_{n\in \N^{*}}$ on voit que 
      $\xi_{n+1}(\gamma_{1},...,\gamma_{n},\gamma\gamma')$ est l'image de  $\xi_{n+2}(\gamma_{1},...,\gamma_{n},\gamma,\gamma')$ par le morphisme 
   $$(\wh G)^{n+2}\modmod \wh G \to (\wh G)^{n+1}\modmod \wh G, 
   (h_{1},...,h_{n},h,h')\mapsto (h_{1},...,h_{n},hh').$$
       On en déduit que 
        $   (g_{1},...,g_{n},gg')$ est au-dessus de 
        $  \xi_{n+1} (\gamma_{1},...,\gamma_{n},\gamma\gamma')$. De plus 
        $ (g_{1},...,g_{n},gg')$ est semi-simple par le même argument que dans A), car 
       \begin{gather*}  \dim(\ov{< g_{1},...,g_{n},gg' >}) 
        \leq  \dim(\ov{< g_{1},...,g_{n},g,g' >}) 
         \leq  \dim(\ov{< g_{1},...,g_{n}  >}) \end{gather*}
          (où la dernière inégalité vient de (H1)). Donc 
          $ (g_{1},...,g_{n},gg')$ appartient à  
        $  \xi_{n+1}^{\mr{ss}} (\gamma_{1},...,\gamma_{n},\gamma\gamma')$
          et 
           $gg'=\sigma(\gamma \gamma')$. 
        Les mêmes arguments montrent que 
     \begin{gather}\label{gg'-gamma-gamma'}g=\sigma(\gamma)\text{ \ \   et \ \ }  g'=\sigma(\gamma').\end{gather} 
On a finalement montré que $\sigma(\gamma \gamma') =\sigma(\gamma)\sigma(\gamma')$. 
  
  Donc $\sigma$ est un morphisme de groupes à valeurs dans $\wh G(E')$
  (où $E'$ est une extension finie de $E$ telle que $g_{1},...,g_{n}$ appartiennent à $\wh G(E')$).    
L'argument pour montrer que $\sigma$ est continu est le suivant. 
On sait que  pour toute fonction $f$ sur $(\wh G_{E'})^{n+1}\modmod \wh G_{E'}$, l'application
$$\on{Gal}(\ov F/F)\to E',  \ \ \gamma\mapsto f(g_{1},...,g_{n}, \sigma(\gamma))= 
\Theta_{n+1}^{\nu}( f)(\gamma_{1},...,\gamma_{n},\gamma)$$
est continue. 
Or le morphisme 
 \begin{align*}\mc O((\wh G_{E'})^{n+1}\modmod \wh G_{E'}) &\to  \mc O(\wh G_{E'}\modmod C(g_{1},...,g_{n})) \\ f & \mapsto  [g\mapsto f(g_{1},...,g_{n},g)]
\end{align*} est surjectif
(parce que $(g_{1},...,g_{n})$ est semi-simple, donc son orbite par conjugaison est une sous-variété affine fermée de 
$(\wh G_{E'})^{n}$, 
isomorphe à $\wh G_{E'}/C(g_{1},...,g_{n})$). 
De plus, en notant 
$D(g_{1},...,g_{n})$ le centralisateur de $C(g_{1},...,g_{n})$ (qui contient l'image de $\sigma$), le morphisme de restriction 
$$\mc O(\wh G_{E'}\modmod C(g_{1},...,g_{n})) =
\mc O(\wh G_{E'} )^{C(g_{1},...,g_{n})} \to \mc O(D(g_{1},...,g_{n}))$$ est surjectif  parce que la   restriction $\mc O(\wh G_{E'}) \to \mc O(D(g_{1},...,g_{n}))$  est  évidemment   surjective, 
qu'elle le reste   lorsqu'on prend  les invariants par le groupe réductif   $C(g_{1},...,g_{n})$ (agissant par conjugaison), et que 
$C(g_{1},...,g_{n})$ agit trivialement sur  $\mc O(D(g_{1},...,g_{n}))$.  
Donc pour toute fonction $h\in \mc O(D(g_{1},...,g_{n}))$, l'application 
$$\on{Gal}(\ov F/F)\to E', \ \ \gamma\mapsto h(\sigma(\gamma))$$ est continue, et on a montré que $\sigma$ est continu. 
 
 Il reste à montrer \eqref{rel-sigma-nu}, c'est-à-dire que pour 
   $m\in \N^{*}$, $f\in \mc O((\wh G)^{m}\modmod \wh G)$ et  $(\delta_{1},...,\delta_{m})\in \on{Gal}(\ov F/F)^{m}$, on a  
\begin{gather}\nonumber f(\sigma(\delta_{1}),...,\sigma(\delta_{m}))=
 \big(\Theta^{\nu}_{m}(f)\big)(\delta_{1},...,\delta_{m}).\end{gather}
Par les mêmes  arguments que  pour \eqref{semi-simple-g-g'} et \eqref{gg'-gamma-gamma'} on montre que  
 $$(g_{1},...,g_{n},\sigma(\delta_{1}),...,\sigma(\delta_{m}))\in \xi_{n+m}^{ss}(\gamma_{1},...,\gamma_{n},\delta_{1},...,\delta_{m} ). $$
 Donc 
 $(\sigma(\delta_{1}),...,\sigma(\delta_{m}))$ est au-dessus de  $\xi_{m}(\delta_{1},...,\delta_{m} )$. \cqfd

Donc on a obtenu la décomposition \eqref{intro1-dec-canonique}. Ceci achève la démonstration du \thmref{intro-thm-ppal}, à condition d'admettre les deux résultats suivants, qui seront justifiés  dans le prochain paragraphe:  
\begin{itemize}
\item
tout $\sigma$  apparaissant dans la décomposition \eqref{intro1-dec-canonique}  est non ramifié en dehors de $N$. 
\item 
    la décomposition \eqref{intro1-dec-canonique} est compatible  avec l'isomorphisme de Satake en toutes les places de $X\sm N$. 
    \end{itemize}

    \subsection{Compatibilité avec l'isomorphisme de Satake aux places non ramifiées} 
    \label{subsection-intro-decomp}
    Le but de ce paragraphe est de montrer les deux  résultats admis à la fin du paragraphe précédent. 

\begin{lem}\label{lem-sigma-nonram} 
Tout paramètre  $\sigma$ apparaissant dans \eqref{intro1-dec-canonique} est non ramifié  sur $X\sm N$. 
\end{lem}
\noindent{\bf Démonstration}    (pour un énoncé plus fort  et plus de détails on renvoie à la \propref{prop-harris}). 
         Soit  $v$ une  place de $X\sm N$. On fixe  un plongement $\ov F\subset \ov {F_{v}}$, d'où une  inclusion 
            $ \on{Gal}(\ov {F_{v}}/F_{v})\subset  \on{Gal}(\ov F/F)$. Soit  $I_{v}=\on{Ker}( \on{Gal}(\ov {F_{v}}/F_{v})\to \wh \Z)$ le groupe d'inertie en   $v$. Alors 
pour $I, W,x,\xi$ comme dans \eqref{excursion-def-intro}, l'image de  
la composée $  H_{\{0\},\mbf  1}\xrightarrow{\mc H(x)}
 H_{\{0\},W^{\zeta_{I}}}\isor{\chi_{\zeta_{I}}^{-1}} 
  H_{I,W}$ (qui est le début de \eqref{excursion-def-intro}) est formée d'éléments invariants par 
    $(I_{v})^{I}$, car les opérateurs 
    de création sont des morphismes de faisceaux  sur $\Delta(X\sm N)$ tout entier (et en particulier en $\Delta(v)$).    Donc pour $(\gamma_{i})_{i\in I}\in \on{Gal}(\ov F/F)^{I}$ et 
    $(\delta_{i})_{i\in I}\in (I_{v})^{I}$ on a 
    \begin{gather}\label{rel-S-gammai-deltai}S_{I,W,x,\xi,(\gamma_{i})_{i\in I}}=S_{I,W,x,\xi,(\gamma_{i}\delta_{i})_{i\in I}}.\end{gather} 
     Grâce à la preuve de l'unicité de $\sigma$ incluse dans la démonstration de la \propref{intro-Xi-n}, la relation \eqref{rel-S-gammai-deltai}   implique 
   que pour tout $\sigma$ correspondant à un caractère $\nu$ de $\mc B$, on a   $I_{v}\subset \Ker \sigma$ et donc $\sigma$ est non ramifié en $v$.  
      \cqfd

  Le lemme suivant  montre que les 
    opérateurs de Hecke en les places non ramifiées sont des cas particuliers d'opérateurs d'excursion. 
    
    Soit 
      $v$ une place dans  $X\sm N$.   On fixe un plongement 
       $\ov F\subset \ov F_{v}$. 
    Comme précédemment  $\mbf 1\xrightarrow{\delta_{V}} V\otimes V^{*}$  et   $ V\otimes V^{*}\xrightarrow{\on{ev}_{V}} \mbf 1$ sont les morphismes naturels.

     \begin{lem} \label{S-non-ram-intro}     Pour tout $d\in \N$ et tout  $\gamma\in \on{Gal}(\ov {F_{v}}/F_{v})\subset  \on{Gal}(\ov F/F)$ tel que $\deg(\gamma)=d$, $
   S_{\{1,2\},V \boxtimes V^{*},\delta_{V},\on{ev}_{V},(\gamma,1)}$  dépend seulement de   $d$, 
  et si    $d=1$ il est égal à  $T(h_{V,v})$. 
     \end{lem}
       \noindent{\bf Démonstration.}
 On fixe un  point géométrique  $\ov v$ au-dessus de  $v$ et une  flèche de spécialisation 
    $\on{\mf{sp}}_{v}:\ov \eta\to \ov v$, associés au plongement $\ov F\subset \ov F_{v}$ choisi ci-dessus. On note encore   $\on{\mf{sp}}_{v}$ la  flèche de spécialisation    
    $\Delta(\ov \eta)\to \Delta(\ov v)$ égale à son image   par $\Delta$. Pour que le diagramme suivant tienne dans la page on   pose    $I=\{1,2\}$ et $W=V\boxtimes V^{*}$. 
    Le diagramme  
       $$  \xymatrix{  C_{c}^{\mr{cusp}}(G(F)\backslash G(\mb A)/K_N \Xi,E)\ar[d]_-{\restr{\mc C_{\delta_{V}}^{\sharp}}{\ov v}} \ar[dr]^{ \mc C_{\delta_{V}}^{\sharp}} & &
       \\ \Big( \varinjlim _{\mu} \restr{\mc H _{N, I, W}^{0,\leq\mu,E}}{\Delta(\ov v)}\Big)^{\mr{Hf}} \ar[r]^-{\mf{sp}_{v}^{*}} \ar[d]^{F_{\{1\}}^{\deg(v)d}} 
       & 
       \Big( \varinjlim _{\mu} \restr{\mc H _{N, I, W}^{0,\leq\mu,E}}{\Delta(\ov\eta)}\Big)^{\mr{Hf}} 
               \ar@{=}[r]  & 
       H_{I,W}
       \ar[d]^{(\gamma,1)}
       \\  \Big( \varinjlim _{\mu} \restr{\mc H _{N, I, W}^{0,\leq\mu,E}}{\Delta(\ov v)}\Big)^{\mr{Hf}} \ar[r]^-{\mf{sp}_{v}^{*}}\ar[d]_-{\restr{\mc C_{\on{ev}_{V}}^{\flat}}{\ov v}}   & 
        \Big( \varinjlim _{\mu} \restr{\mc H _{N, I, W}^{0,\leq\mu,E}}{\Delta(\ov\eta)}\Big)^{\mr{Hf}}\ar@{=}[r] \ar[dl]^{ \mc C_{\on{ev}_{V}}^{\flat}} & 
      H_{I,W}
       \\ C_{c}^{\mr{cusp}}(G(F)\backslash G(\mb A)/K_N \Xi,E)  & &
       } $$
   est commutatif (la commutativité du grand rectangle sera justifiée  
   dans le \lemref{S-non-ram-prelim2}). 
   Or    $S_{\{1,2\},V \boxtimes V^{*},\delta_{V},\on{ev}_{V},(\gamma,1)}$  est égal par définition à la composée par le chemin le plus à droite.    Donc      il  est égal à la  composée donnée par la colonne de gauche. Par conséquent  il dépend seulement de   $d$. 
   Lorsque $d=1$ la composée donnée par la colonne de gauche est égale par définition à $S_{V,v}$, et donc à $T(h_{V,v})$  par  la \propref{prop-coal-frob-cas-part-intro}. \cqfd

       \begin{rem} On n'a calculé la composée de la colonne de gauche que pour $d=1$ mais pour d'autres valeurs de $d$ elle n'apporte rien de nouveau 
       car on pourrait montrer qu'elle est égale à une combinaison de $S_{W,v}$ avec $W$ représentation irréductible de $\wh G$. 
       \end{rem}
       
   La proposition suivante affirme  la compatibilité de la décomposition \eqref{intro1-dec-canonique} avec 
l'isomorphisme de Satake  en les places de $X\sm N$.           

\begin{prop}\label{S-non-ram-concl-intro} 
Soit  $\sigma$ apparaissant 
 dans \eqref{intro1-dec-canonique}  et  
$v$ une place    de $X\sm N$. Alors  
$\sigma$ est non ramifié en $v$ et 
 pour toute représentation irréductible $V$ de $\wh G$, 
       $T(h_{V,v})$ agit sur   $\mf H_{\sigma}$
       par multiplication par le scalaire  $\chi_{V}(\sigma(\Frob_{v}))$, où $\chi_{V}$ est le caractère de $V$ et $\Frob_{v}$ est un relèvement arbitraire d'un élément de  Frobenius en $v$. 
 \end{prop}
 \dem 
 Le fait que $\sigma$ est non ramifié en $v$ a déjà été établi dans le \lemref{lem-sigma-nonram}. On reprend les notations du \lemref{S-non-ram-intro}. 
 Comme 
    $  \s{\on{ev}_{V}, (\sigma(\gamma), 1) . \delta_{V} }=\chi_{V}(\sigma(\gamma))$ ce lemme implique que  
        pour tout     $\gamma\in  \on{Gal}(\ov {F_{v}}/F_{v})$  avec  $\deg(\gamma)=1$, 
      et toute   représentation  irréductible $V$ de $\wh G$, 
   $\mf H_{\sigma}$ est inclus dans l'espace propre généralisé (ou espace caractéristique)   de $T(h_{V,v})$ pour la valeur propre 
             $\chi_{V}(\sigma(\gamma))$. Or on sait  que les opérateurs de Hecke aux places non ramifiées sont diagonalisables
      (car ce sont des opérateurs normaux sur l'espace hermitien des formes automorphes cuspidales à coefficients dans $\C$). Donc $T(h_{V,v})$ agit sur $\mf H_{\sigma}$  par homothétie de rapport 
      $\chi_{V}(\sigma(\gamma))$. 
\cqfd

 Cela termine  la preuve  du \thmref{intro-thm-ppal}.

   \subsection{Remarques supplémentaires}
   \label{intro-rem-suppl}
      La décomposition \eqref{intro1-dec-canonique} est certainement plus fine  en général que celle obtenue  par diagonalisation des opérateurs de Hecke en les places  non ramifiées. Même en prenant en compte les classes d'isomorphisme  de représentations de 
 $C_{c}(K_{N}\backslash G(\mb A)/K_{N},\Qlbar)$ on ne récupère pas en général la décomposition \eqref{intro1-dec-canonique}, et bien que les  formules de multiplicités d'Arthur fassent intervenir une somme  sur les  paramètres d'Arthur, une telle  décomposition canonique semble inconnue  en général dans le cas des corps de nombres. 
    En effet d'après des  exemples de \cite{blasius,lapid}, pour certains groupes $ G $ (y compris déployés),   la même représentation de 
 $C_{c}(K_{N}\backslash G(\mb A)/K_{N},\Qlbar)$ peut apparaître dans des espaces $\mf H_{\sigma}$ différents, à cause du phénomène suivant. Il y a des exemples de groupes finis $ \Gamma $ et de morphismes $ \tau, \tau ': \Gamma \to \wh G (\Qlbar) $ tels que $ \tau $ et $ \tau' $ ne soient pas conjugués mais que pour tout  $ \gamma \in \Gamma $, $ \tau (\gamma) $ et $ \tau '(\gamma) $ soient conjugués.
 On s'attend alors à ce qu'il existe un    morphisme surjectif  $ \rho: \on {Gal} (\ov F / F) \to \Gamma $  partout non ramifié et une représentation $ (H_ {\pi}, \pi) $  de $ G (\mb A) $ tels que $ (H_ {\pi}) ^ {K_ {N}} $ apparaisse à la fois dans  $ \mf H_ {\tau \circ \rho} $ et $ \mf H_ {\tau' \circ \rho} $.
 Jusqu'à présent le seul moyen systématique de distinguer ces copies de  $ \pi $ était le programme de Langlands géométrique  (c'est d'ailleurs  ce qui explique que notre approche permette aussi de le faire,  grâce au lien avec le programme de Langlands géométrique qui est expliqué dans le chapitre \ref{subsection-link-langl-geom}). 
 
 Les exemples de Blasius et Lapid sont pour $ G = SL_{r} $, $ r \geq 3$
  (en fait dans ce cas on peut retrouver {\it a posteriori} la décomposition \eqref{intro1-dec-canonique} à l'aide du plongement  $SL_{r}\hookrightarrow GL_{r}$, cf la \remref{rem-lapid}). 
       On renvoie à  \cite{larsen1,larsen2} pour une liste de groupes  $ \wh G $ pour lesquels existent des exemples $ \Gamma, \tau, \tau '$ comme ci-dessus  (et  pour certains de ces groupes, par exemple $E_{8}$, 
 nous ne savons pas comment retrouver la décomposition \eqref{intro1-dec-canonique} autrement que par le programme de Langlands géométrique ou par les méthodes du présent  article, qui ne marchent que sur les corps de fonctions). Plus récemment  S. Wang
    \cite{wang,wang12a,wang12b} a donné des exemples  de $ (\Gamma, \tau, \tau ') $ comme ci-dessus, mais avec  $ \Gamma $ semi-simple connexe.

 \subsection{Lien avec les travaux antérieurs}
  \label{intro-previous-works}

 Le lien avec le programme de Langlands géométrique
  est extrêmement étroit, et fait l'objet du chapitre \ref{subsection-link-langl-geom}.

      Au contraire, les méthodes utilisées dans ce travail sont complètement différentes de celles fondées sur la formule des traces qui ont été développées notamment 
  par Drinfeld 
   \cite{drinfeld78,Dr1,drinfeld-proof-peterson,drinfeld-compact}, Laumon, Rapoport et Stuhler \cite{laumon-rapoport-stuhler}, 
  Laumon \cite{laumon-drinfeld-modular,laumon-cetraro},    Laurent  Lafforgue \cite{laurent-asterisque,laurent-jams,laurent-inventiones,laurent-tata}, 
 Ngô Bao Châu \cite{ngo-jacquet-ye-ulm,ngo-modif-sym},  Eike Lau \cite{eike-lau,eike-lau-duke}, Ngo Dac Tuan \cite{ngo-dac-ast,ngo-dac-09,ngo-dac-11},  Ngô Bao Châu et  Ngo Dac Tuan  \cite{ngo-ngo-elliptique}, Kazhdan et  Varshavsky~\cite{kvar,var-SANT} et 
 Badulescu et Roche \cite{badulescu}. 
  
       Cependant  l'action sur la cohomologie des groupes de permutations des pattes des chtoucas apparaît déjà dans 
les travaux de  Ngô Bao Châu,  Ngo Dac Tuan et Eike Lau que nous venons de citer. Ces actions des groupes de permutations jouent par ailleurs un rôle essentiel  dans le programme de Langlands géométrique, et notamment dans la preuve par Gaitsgory de la conjecture d'annulation \cite{ga-vanishing}. 
D'autre part
 la coalescence des pattes  apparaît dans la thèse de Eike Lau \cite{eike-lau} et elle est aussi  utilisée dans le preprint   \cite{brav-var} de Braverman et Varshavsky  
 (afin de montrer la non nullité de certains des morphismes \eqref{trace-brav-var-piIW}).  
 Les chapitres  
 \ref{para-1-chtoucas} et \ref{para-def-coh}  du présent article contiennent des rappels de l'article \cite{var} de Varshavsky, dont les résultats généraux sur la  structure locale et les  propriétés des champs  
 de $G$-chtoucas sont fondamentaux, et aussi du preprint  très éclairant  \cite{brav-var} de Braverman et  Varshavsky. 
   Dans le cadre du lien entre cet article et le programme de Langlands géométrique on mentionnera dans le  chapitre \ref{subsection-link-langl-geom}  l'analogie entre les résultats de cet article et  le   corollaire  4.5.5 de 
 \cite{dennis-laumon}.

\subsection{Considérations historiques}  L'idée de caractériser les classes de conjugaison de représentations de groupes de monodromie par un nombre fini d'invariants numériques apparaît déjà dans les travaux de  Poincaré. Dans  \cite{poincare} il considère une équation différentielle linéaire sur $\mb P^{1}-\{x_{0},...,x_{n}\}$ dont les coefficients sont des fonctions rationnelles et dont l'espace des solutions est localement de dimension $r$. Poincaré caractérise les classes de conjugaison par $GL_r(\C)$ des représentations  de monodromie par les polynômes caractéristiques d'un nombre fini d'éléments de $\pi_{1}(\mb P^{1}-\{x_{0},...,x_{n}\})$.  
 Il note que les polynômes caractéristiques des éléments de $\pi_{1}(\mb P^{1}-\{x_{0},...,x_{n}\})$ qui sont conjugués à un petit cercle autour d'un  des $x_{i}$ se calculent de fa\c con locale (un peu comme pour les $\Frob_{v}$ dans notre problème) alors que les autres invariants sont globaux.  Comme $\pi_{1}(\mb P^{1}-\{x_{0},...,x_{n}\})$ est un groupe libre à $n$ générateurs, 
 les invariants numériques considérés par Poincaré sont associés à des fonctions dans   
 $\mc O((GL_r)^{n}\modmod GL_r)$. C'est Hilbert qui a montré dans \cite{hilbert} que  de telles algèbres sont  de type fini.

\subsection{Plan de l'article} La première partie est consacrée
 à l'étude des champs de chtoucas, de leur cohomologie et des  morphismes de création et d'annihilation.  Le chapitre \ref{rappels-Hecke-Gr-satake} rappelle l'équivalence de Satake géométrique. 
 Les chapitres \ref{para-1-chtoucas} à \ref{para-def-coh} rappellent la définition et les propriétés élémentaires des champs classifiant les  $G$-chtoucas. 
En particulier les chapitres \ref{para-1-chtoucas} et \ref{para-def-coh}
  reprennent presque mot à mot des parties de l'article \cite{var} de Varshavsky. Le chapitre \ref{para-creation-annih}
 contient la construction des morphismes de création et d'annihilation. 
 Le chapitre \ref{frob-coalescence-cas part} est consacré à la preuve de la  \propref{prop-coal-frob-cas-part-intro}. Dans le chapitre \ref{para-Relations d'Eichler-Shimura} on montre les relations d'Eichler-Shimura.     
   
  La deuxième partie fait  intervenir les actions des groupes de Galois. 
  Dans le chapitre \ref{para-sous-faisceaux-Frob-partiels} on utilise les  relations d'Eichler-Shimura  pour munir 
$\Big( \varinjlim _{\mu}\restr{\mc H _{ N, I, W}^{0,\leq\mu,E}}{\ov{\eta^{I}}}\Big)^{\mr{Hf}}$ d'une action de 
$\pi_{1}(\eta,\ov{\eta})^{I}$. 
 Dans le chapitre 
  \ref{construction-S-proprietes}
   on construit les morphismes 
  d'excursion et  on étudie leurs propriétés. On en déduit dans le chapitre 
  \ref{para-dec-param-Langlands} la décomposition \eqref{intro1-dec-canonique}. Pour ne pas mélanger les difficultés, on a rédigé les chapitres \ref{rappels-Hecke-Gr-satake} à  \ref{para-dec-param-Langlands} en supposant  $G$  déployé. Le chapitre \ref{para-non-deploye} apporte les modifications nécessaires pour traiter le cas des groupes non nécessairement déployés.  Le chapitre \ref{mod-ell} montre que la décomposition \eqref{intro1-dec-canonique} existe aussi à coefficients dans $\ov {\mathbb F_{\ell}}$. 
    Enfin le chapitre \ref{para-meta} indique  les modifications nécessaires pour traiter le cas des groupes   métaplectiques,  le chapitre
\ref{subsection-link-langl-geom} explique le lien avec le programme de Langlands  géométrique et le chapitre \ref{GL-previous-works} concerne le cas de $GL_r$. 
  
  \subsection{Remerciements}

 La première partie  de cet article (jusqu'au chapitre \ref{frob-coalescence-cas part})  est  issue d'un travail en cours en commun avec Jean-Benoît Bost, qui est à l'origine du travail exposé dans cet article mais n'est pas encore achevé. 
 Ce travail en commun a été commencé plus d'un an avant que n'émergent les idées   du présent article, 
 qui n'existerait évidemment pas sans lui. 
     De très longues  discussions avec Alain Genestier ont  été déterminantes et sans elles cet article n'existerait pas non plus.      Je remercie aussi Jean-Benoît Bost et Alain Genestier pour leur aide dans la mise au point  de nombreux arguments.

     Cet article repose doublement sur des travaux de Vladimir Drinfeld, sur les chtoucas et  sur le programme de 
 Langlands géométrique, et je lui exprime ma grande reconnaissance.

    Je remercie  Yakov Varshavsky  pour ses explications et pour m'avoir communiqué le preprint non publié~\cite{brav-var}. 
 Je remercie Vladimir Drinfeld et Dennis Gaitsgory  pour leurs nombreuses explications.  
 Je remercie Laurent Lafforgue pour des discussions et l'idée du titre. 
 Je remercie Sergey Lysenko pour son aide dans le cas métaplectique et de nombreuses discussions. Je remercie  Jochen Heinloth pour son aide concernant les modèles de Bruhat-Tits. Je remercie 
  Gebhard B\"ockle,  Michael Harris, Chandrashekhar Khare et Jack Thorne pour avoir trouvé l'énoncé de la \propref{prop-harris} qui rend la rédaction plus limpide et est de plus utilisé dans  leur preprint \cite{boeckle-harris...}. 
 Je remercie   Ga\"etan Chenevier, Pierre Deligne,  Ed Frenkel,  Nicholas Katz,  Erez Lapid, Gérard Laumon,   Colette Moeglin, Sophie Morel, Ngô Bao Châu, Ngo Dac Tuan, Jean-Pierre Serre,   Jean-Loup Waldspurger,  Cong Xue et Xinwen Zhu pour leurs explications et leurs remarques.   Je remercie les rapporteurs anonymes pour avoir relu  le texte de fa\c con extrêmement soigneuse,  corrigé  des erreurs et apporté des améliorations. 
       
 Je tiens à remercier vivement le CNRS. Le programme de Langlands est  un sujet très différent de ma spécialité d'origine et je n'aurais  pas pu  m'y consacrer sans la grande liberté laissée aux chercheurs pour mener à bien leurs travaux. 
 Je remercie aussi mes collègues du MAPMO et de l'Institut Fourier pour leur soutien. 
   
        \subsection{Notations et conventions}
     Les notations et conventions suivantes  concernent surtout les chapitres \ref{rappels-Hecke-Gr-satake} à \ref{para-dec-param-Langlands}.  
     Elles 
     reprennent celles de \cite{var}. 
   
   1) Soit $G$ un  groupe réductif connexe déployé sur un corps fini $\Fq$
   (les notations pour les groupes non nécessairement déployés seront introduites seulement au chapitre \ref{para-non-deploye}). 
   On note   $G^{\mr{der}}$ le groupe dérivé de $G$, $G^{\mr{sc}}$ le revêtement simplement connexe de $G^{\mr{der}}$, $G^{\mr{ab}}:=G/G^{\mr{der}}$  l'abélianisé de $G$, et $G^{\mr{ad}}$ le groupe adjoint de  $G$.  Soit $B\supset T\supset Z$ un sous-groupe de Borel, un tore maximal  et le centre de $G$, respectivement. On note  $B^{\mr{sc}}\supset T^{\mr{sc}}\supset Z^{\mr{sc}}$ 
les sous-groupes correspondants de  $G^{\mr{sc}}$, et de même pour $G^{\mr{der}}$ et  $G^{\mr{ad}}$. On note $X^{*}(T),X_{*}(T)$ les groupes 
des poids et copoids de $G$, et $X^{*}_{+}(T),X_{*}^{+}(T)$ les sous-ensembles de poids et copoids dominants. Les poids (resp.  copoids) de $G$ sont munis de l'ordre (standard) suivant : $\lambda_{1}\leq \lambda_{2}$ si et seulement si la différence $\lambda_{1}- \lambda_{2}$ est une combinaison à coefficients rationnels positifs  de racines (resp. coracines) simples de $G$. On adopte des notations semblables pour $G^{\mr{sc}}$, $G^{\mr{der}}$ et $G^{\mr{ad}}$. On note $\rho$ la demi-somme des coracines positives de $G$. 
On note $\wh G$ le groupe dual de Langlands de $G$, considéré comme  un groupe réductif connexe déployé défini sur $\Ql$ (et même sur $\Zl$ dans le chapitre \ref{mod-ell}). Pour toute extension 
$E$ de $\Ql$ on note $\wh G_{E}$ le groupe algébrique défini sur $E$ et $\wh G(E)$ les points à valeurs dans $E$. Les représentations irréductibles de $\wh G$ sont en bijection avec $X_{*}^{+}(T)$. 
 
 2) Pour tout poids dominant $\lambda $ de $G$ on note $V_{\lambda}$ le module de Weyl de $G$ de plus haut poids $\lambda$.    
  
 3) Pour tout schéma fini  $N$ sur $\Fq$, on note 
$\mc{O}_N:=\Fq[N]$. On note $G_{N}$ la restriction à la Weil de $G$ de $N$ à $\Fq$ (c'est un schéma en groupes lisse de dimension $\deg(N)$ sur $\Fq$).  
 
  4) Soit $X$ une courbe projective lisse géométriquement irréductible  sur $\Fq$. 
On note $F$ le corps des fonctions rationnelles sur $X$. 
Pour 
tout point fermé $v$ de $X$, on note $\mc O_{v}$  l'anneau complété du faisceau structurel en $v$, $F_{v}$ son corps des fractions et $k(v)$ son corps résiduel, de sorte que $v=\on{Spec}(k(v))$ est un sous-schéma de $X$. 
On note $\mb A=\prod ' F_{v}$ l'anneau des adèles. On note $\eta=\on{Spec}(F)$ le point générique de $X$. Pour tout ensemble fini $I$ on note 
$F^{I}$ le corps des fonctions de $X^{I}$ et $\eta^{I}=\on{Spec}(F^{I})$ le  point générique. 

5) Pour tout  $S$-point $x$ d'un schéma $X$, on note  $\Gamma_x\subset X\times
S$ le graphe de  $x$.

6) Sauf mention explicite du contraire, $E$ désignera une extension finie de $\Ql$ contenant une racine carrée de $q$, et $\mc O_{E}$ son anneau d'entiers. 

7) Pour tout champ $\mc Y$ sur un corps fini $\Fq$,  on définit le faisceau d'intersection $\IC^{E}_{\mc Y}$ comme le prolongement intermédiaire  du $E$-faisceau pervers constant sur un sous-champ ouvert $\mc Y^{0}$ de $\mc Y$ tel que le champ réduit correspondant 
$(\mc Y^{0})_{\mr{red}}$ soit lisse. Ce faisceau d'intersection $\IC^E_{\mc Y}$ est normalisé pour être pur de poids $0$. Cependant on se trouvera souvent dans une situation où $\IC^E_{\mc Y}$ est universellement localement acyclique relativement à un morphisme de  $\mc Y$ vers une base lisse 
(qui sera typiquement $(X\sm N)^{I}$), et où on  normalisera le degré et le poids de $\IC^E_{\mc Y}$ relativement 	à ce morphisme. 

 8) Pour tout champ  $\mc Y$ sur $\Fq$, on note 
   $\Frob_{\mc Y/\Fq}:\mc Y\to \mc Y$ le morphisme de Frobenius absolu sur $\Fq$. On notera souvent 
   $\Frob _{\mc Y}$ ou simplement $\Frob$ au lieu de $\Frob_{\mc Y/\Fq}$. 
   
   9) 
Pour tout champ  $S$ sur $\Fq$ et  pour tout faisceau cohérent ou pour tout  $G$-torseur 
$\mc{F}$ sur  $X\times S$, on écrira $\ta\mc{F}$ au lieu de 
$(\Id_X\times \Frob_{S} )^*(\mc{F})$. On adopte une notation semblable pour les morphismes.

11) On adopte les conventions du chapitre 1.1 de \cite{weil2} concernant les $A$-faisceaux constructibles, pour $A=\mc O_{E},E$ ou $\Qlbar$. Tous les $A$-faisceaux considérés dans cet article seront constructibles. On ne considérera jamais  les limites inductives  autrement que comme des systèmes inductifs abstraits, sauf pour leurs fibres : si  $\varinjlim \mc F_{\mu}$ est un système inductif 
 de $A$-faisceaux constructibles  sur une variété $Y$ et $\ov x$ est un point géométrique de $Y$ au-dessus d'un point $x$, $\restr{\varinjlim \mc F_{\mu}}{\ov x}$ sera considéré comme un $A$-module, muni d'une action de 
 $\pi_{1}(x,\ov{x})$. 
 
 12) Si $k$ est un corps, $\ov k$ désignera une clôture  algébrique de $k$. 
 
  \tableofcontents   
 
\section{Champs de Hecke, grassmannienne affine de Beilinson-Drinfeld, et   équivalence de  Satake  géométrique }
\label{rappels-Hecke-Gr-satake}

Ce chapitre constitue un rappel de l'équivalence de Satake géométrique,  
 due à  Lusztig, Drinfeld, Ginzburg, et Mirkovic--Vilonen   
   \cite{lusztig-satake,ginzburg,hitchin,mv}. 
On renvoie à  \cite{mv,hitchin,ga-de-jong,ga-iwahori,var, brav-var} pour plus de détails. 

 Soit $G$ un  groupe réductif connexe déployé sur un corps fini $\Fq$. 
Soit $X$
une courbe lisse projective et géométriquement irréductible sur
 $\Fq$.

 Soit  $\Bun_{G}$ le champ lisse sur $\Fq$ classifiant les  $G$-torseurs sur  $X$, c'est-à-dire que 
     $$\Bun_{G}(S)=\{G\text{-torseur  sur }X\times S, \text{localement trivial pour la topologie étale}\}. $$
    D'après le théorème 2 de \cite{DS} tout  $G$-torseur sur   $X\times S$ devient  localement trivial pour la   topologie de Zariski de $X\times S$ après  un changement de base étale convenable sur   $S$. 
  Cependant nous n'aurons pas besoin de ce résultat.

 Si $N$ est un sous-schéma fini de  $X$ on note   $\Bun_{G,N}$ le champ classifiant les   $G$-torseurs  sur  $X$ avec structure de niveau $N$, c'est-à-dire que    \begin{gather}\label{defi-Bun-GN} \Bun_{G,N}(S)=\{\mc G\in \Bun_{G}(S), \psi : \restr{\mc G}{N\times S} 
   \isom 
  \restr{G}{N\times S} \}.\end{gather}
  
  Plus généralement on utilisera une   version en famille de $\Bun_{G,N}$ quand  $N$ se déplace sur la courbe. 
  Pour tout schéma $T$ sur $\Fq$ et tout  sous-schéma fermé  $Q\subset X\times T$ qui est fini sur $T$ et localement défini par une équation (c'est-à-dire que $Q$ est un diviseur de Cartier relatif effectif), on note   $\Bun_{G,Q}$ le champ lisse sur  $T$ tel que pour tout schéma  $S$ sur  $T$,  
  \begin{gather}\label{defi-Bun-GQ} \Bun_{G,Q}(S)=\{\mc G\in \Bun_{G}(S), \psi : \restr{\mc G}{Q\times _{(X\times T)} (X\times S)} 
   \isom 
  \restr{G}{Q\times _{(X\times T)} (X\times S)} \}.\end{gather}
On note  $G_{Q}$ le  schéma en groupes lisse sur  $T$ défini  comme la restriction à la Weil  de $G$ de $Q$ à $T$. Il est de dimension relative  
 $\deg(Q)\dim G$ sur $T$, où  $\deg(Q)$ est le  degré de  $Q$ (qui est  une fonction localement constante sur  $T$). 
On remarque que  $\Bun_{G,Q}$ est un $G_{Q}$-torseur sur $\Bun_{G}\times T$. 
 
    Les composantes irréductibles  de  $\Bun_{G}$ ne sont pas  de type fini mais  $\Bun_{G}$ est une réunion d'ouverts $\Bun_{G}^{\leq \mu}$, définis  par troncature par le polygone de Harder-Narasimhan   de $\mc G$, et  dont les composantes irréductibles sont de type fini. Plus précisément, 
 pour tout  $\mu\in X_{*}(T^{\mr{ad}})$, on pose  
   \begin{gather} 
   \label{defi-Harder-Naram-split}
   \Bun_{G}^{\leq \mu}(S) = \{\mc G\in \Bun_{G}(S) |
    \text{ pour tout point géométrique } s\in S, \\ \nonumber
     \text{ toute $B$-structure $\mc B$ sur $\mc G_{s}$ et tout  } \lambda\in X_{+}^{*}(T^{\mr{ad}}), 
   \on{deg} 
  \mc B_{\lambda}
  \leq \s{\mu,\lambda}
   \},  \end{gather}
où $\mc B_{\lambda}$ est le fibré en droites correspondant. 
Il s'agit donc d'une troncature de Harder-Narasimhan pour le $G^{\mr{ad}}$-torseur déduit de $\mc G$. 
D'après le  lemme A.3 de  \cite{var}, $\Bun_{G}^{\leq \mu}$ est ouvert dans   $\Bun_{G}$ et ses composantes  irréductibles (ou connexes) sont de type fini.  Cela nous suffit mais nous mentionnons qu'une  étude fine de la stratification de Harder-Naramsinham-Shatz est réalisée dans   \cite{schieder}. 
Dans le chapitre \ref{para-non-deploye} consacré aux groupes non nécessairement déployés on définira des troncatures à l'aide d'un plongement de $G^{\mr{ad}}$ dans $SL_{r}$ et on pourrait évidemment faire la même chose ici, ce qui aurait l'avantage de n'utiliser les troncatures    de Harder-Naramsimhan que pour $GL_{r}$.

Pour tout  sous-schéma fini $N$ de  $X$, on note   $\Bun_{G,N} ^{\leq \mu}$ l'image inverse de  $\Bun_{G}^{\leq \mu}$ dans  $\Bun_{G,N}$. On omettra souvent la lettre  $G$ dans toutes ces notations. 

\begin{rem} (\cite{hitchin} et lemme 2.2 de \cite{var}) \label{L:conn}
L'ensemble  $\pi_0(Bun_G)$ des composantes connexes de  $Bun_G$ est canoniquement  isomorphe à  $\pi_1(G):=X_*(T)/X_*(T^{\mr{sc}})$  
(le  quotient du réseau des copoids par le réseau des coracines) qui est lui-même   canoniquement  isomorphe au groupe des  caractères de $Z_{\wh G}$, le centre du  groupe dual de Langlands. On note   $[\omega]\in\pi_1(G)$ la classe  de  $\omega\in X_*(T)$.
\end{rem}

La définition suivante généralise légèrement la  définition  2.4 de \cite{var}  (voir aussi la remarque  2.7 c)  de \cite{var}). 
    
\begin{defi}\label{defi-Hecke}
a) Soit $I$   un ensemble fini, $k\in \N$ et $I_{1}, ..., I_{k}$ des parties de  $I$ telles que  $\{I_{1}, ..., I_{k}\}$ forme une  partition de $I$. Soit  $N$ un 
sous-schéma fini  de $X$. On note  $\Hecke_{N,I} ^{(I_{1},...,I_{k})}$ le champ tel que pour tout schéma  $S$ sur $\Fq$, $\Hecke_{N,I}^{(I_{1},...,I_{k})}(S)$ classifie  les  données consistant en 
\begin{itemize}
\item[]
 i) des points $x_i\in (X\sm N)(S)$ pour $i\in I$, 
\item[]
ii) $(\mc G_{0}, \psi_{0}),..., (\mc G_{k}, \psi_{k})\in \Bun_{G,N}(S)$, 
 \item[]
iii) pour $j\in\{1,...,k\}$, un isomorphisme $$\phi_{j}:\restr{\mc G_{j-1}}{(X\times S)\sm(\bigcup_{i\in I_{j}}\Gamma_{x_i})}\isom \restr{\mc G_{j}}{(X\times S)\sm(\bigcup_{i\in I_{j}}\Gamma_{x_i})},$$ préservant les structures de niveau, c'est-à-dire que  $\psi_{j}\circ \restr{\phi_{j}}{N\times S}=\psi_{j-1}$.
\end{itemize}
On omet  $N$ de la notation lorsque $N=\emptyset$.

b)  Pour tout $I$-uplet $\underline{\omega}=(\omega_i)_{i\in I}$ de copoids dominants de $G$, on note  $\Hecke_{N,I,\on{\lesssim} \underline{\omega}}^{(I_{1},...,I_{k})}$ le  sous-champ fermé de  $\Hecke_{N,I}^{(I_{1},...,I_{k})}$ défini par la condition que, pour tout  $j\in\{1,...,k\}$,  la  ``position relative'' de $\mc G_{j-1}$ (ou plutôt $\phi_{j}(\mc G_{j-1})$) par rapport à  $\mc G_{j}$ en les points   $x_i$ (pour $i\in I_{j}$) est bornée par     $\omega_i$ 
 dans le sens faible suivant: 

iii)$_{\underline{\omega}}$ $\phi_{j}((\mc G_{j-1})_{{V_{\lambda}}})\subset (\mc G_{j})_{{V_{\lambda}}}(\sum_{i\in I_{j}}
\s{\lambda,\omega_i}  \Gamma_{x_i})$ pour tout  poids dominant  $\lambda$ de $G$. 
\end{defi}

\begin{rem} D'après le lemme 3.1 de  \cite{var}, $\Hecke_{N,I,\on{\lesssim}\overline{\omega}}^{(I_{1},...,I_{k})}$ est un champ algébrique, localement de type fini sur  $\Fq$.  On va bientôt introduire  un  sous-champ fermé $\Hecke_{N,I,\underline{\omega}}^{(I_{1},...,I_{k})}\subset \Hecke_{N,I,\on{\lesssim} \underline{\omega}}^{(I_{1},...,I_{k})}$,  qui est celui que l'on utilisera. Lorsque  $G^{\mr{der}}$ est simplement connexe il est  égal au sous-champ réduit, en général il serait égal au sous-champ réduit du  sous-champ fermé de $\Hecke_{N,I,\on{\lesssim} \underline{\omega}}^{(I_{1},...,I_{k})}$ donné par une  condition analogue à la condition 
iii$')_{\underline{\omega}}$ de la définition  2.4 de \cite{var} pourvu que l'on remplace, pour tout  point géométrique $s$ de $S$ la condition sur la composante connexe de  $\Bun_{G}$ par une condition sur la composante connexe  de la  grassmannienne affine en tout point de l'ensemble $\{x_{i}(s)\}$. La condition 
iii$')_{\underline{\omega}}$ de la définition  2.4 de \cite{var} (qui figurait aussi  dans les versions 1 à 3 de cet article sur arXiv)  ne donnait donc sans doute pas la bonne définition de $\Hecke_{N,I,\underline{\omega}}^{(I_{1},...,I_{k})}$ (même à réduction près) pour $G^{\mr{der}}$ non simplement connexe. 
\end{rem}

Dans la suite le  point  de 
$\Hecke_{N,I}^{(I_{1},...,I_{k})}(S)$ fourni par les données i), ii), iii) ci-dessus sera noté sous la forme plus concise 
\begin{gather}\label{donnee-Hecke}\big( (x_i)_{i\in I}, (\mc G_{0}, \psi_{0}) \xrightarrow{\phi_{1}}  (\mc G_{1}, \psi_{1}) \xrightarrow{\phi_{2}}
\cdots\xrightarrow{\phi_{k-1}}  (\mc G_{k-1}, \psi_{k-1}) \xrightarrow{ \phi_{k}}    (\mc G_{k}, \psi_{k})
\big).
\end{gather}

\begin{rem} Le champ $\Hecke_{N,I} ^{(I_{1},...,I_{k})}$ dépend de  l'ordre des parties  $I_{1},...,I_{k}$, c'est pourquoi nous écrivons  $(I_{1},...,I_{k})$ et non $\{I_{1},...,I_{k}\}$. 
Pour simplifier nous appellerons  $(I_{1},...,I_{k})$  une partition de  $I$. On remarque que si certaines parties $I_{j}$ sont vides, on peut les supprimer. 
\end{rem}

 Soit 
$k'\in \{1,...,k\}$ et  $(I'_{1},...,I'_{k'})$ une  partition de  $I$ obtenue à partir de    $(I_{1},...,I_{k})$  en réunissant certaines parties dont les indices sont  adjacents. 
Plus précisément on choisit des entiers  $0=j_{0}<j_{1}< \cdots < j_{k'}=k$
et on pose  $I'_{j'}=\bigcup_{j_{j'-1}< j\leq j_{j'}}I_{j}$.

On considère alors le morphisme d'oubli  
\begin{gather}\label{oubli-Hecke}\pi^{(I_{1},...,I_{k})}_{(I'_{1},...,I'_{k'})}: 
\Hecke_{N,I}^{(I_{1},...,I_{k})} \to 
\Hecke_{N,I}^{(I'_{1},...,I'_{k'})}\end{gather} 
qui envoie 
\eqref{donnee-Hecke} sur   $$\big( (x_i)_{i\in I}, (\mc G'_{0}, \psi'_{0}) \xrightarrow{\phi'_{1}}  (\mc G'_{1}, \psi'_{1}) \xrightarrow{\phi'_{2}}
\cdots\xrightarrow{ \phi'_{k'}}    (\mc G'_{k'}, \psi'_{k'})
\big)$$
  avec   $(\mc G'_{j'}, \psi'_{j'})=(\mc G_{j_{j'}}, \psi_{j_{j'}})$ et 
 $\phi'_{j'}=\phi_{j_{j'}}\circ \cdots \circ \phi_{j_{j'-1}+1}$. 
 Autrement dit on oublie certaines étapes intermédiaires entre  $(\mc G_{0}, \psi_{0})$ et 
$(\mc G_{k}, \psi_{k})$.

 Un cas particulier intéressant est le suivant:  
 $$\pi^{(I_{1},...,I_{k})}_{(I)}: 
\Hecke_{N,I}^{(I_{1},...,I_{k})} \to 
\Hecke_{N,I}^{(I)}$$ 
envoie  \eqref{donnee-Hecke} sur  
$$\big( (x_i)_{i\in I}, (\mc G_{0}, \psi_{0}) \xrightarrow{\phi_{k}\circ \cdots \circ \phi_{1}}    (\mc G_{k}, \psi_{k})
\big).$$

\begin{notation} \label{N:proj}
On note  $p_{0}$  le  morphisme  $\Hecke_{N,I}^{(I_{1},...,I_{k})}\to \Bun_{G,N}$ qui envoie  \eqref{donnee-Hecke} sur  
$(\mc{G}_{0},\psi_{0})$. On note 
$\Hecke_{N,I}^{(I_{1},...,I_{k}), \leq\mu}$  l'image inverse de   $Bun^{\leq\mu}_{G,N}$ par $p_{0}$. 
\end{notation}

La raison pour laquelle on utilise dans cet article les  troncatures par les polygones de  
Harder-Narasimhan  de $\mc G_{0}$ est qu'elles sont préservées par les morphismes d'oubli  $\pi^{(I_{1},...,I_{k})}_{(I'_{1},...,I'_{k'})}$.

Voici la définition de  la  {\em grassmannienne affine de Beilinson-Drinfeld } (sur $X^I$).

\begin{defi} \label{defi-gr}
On note  $\mr{Gr}_{I}^{(I_{1},...,I_{k})}$ 
({\it resp.} $\mr{Gr}_{I,\on{\lesssim}\underline{\omega}}^{(I_{1},...,I_{k})}$) 
l'ind-schéma ({\it resp.} le schéma) classifiant les mêmes données  que 
$\Hecke_{I}^{(I_{1},...,I_{k})}$ 
({\it resp.} $\Hecke_{I,\on{\lesssim}\underline{\omega}}^{(I_{1},...,I_{k})}$), plus une  trivialisation de $\mc{G}_{k}$. 
\end{defi}

Quand  $I$ est un singleton, $\mr{Gr}_{I}^{(I)}$ sera aussi noté $\mr{Gr}$. Sa fibre  $\mr{Gr}_{x}$ en un  point géométrique $x$ de $X$ est la grassmannienne affine habituelle, c'est-à-dire  le quotient  fpqc   $G(F_{x})/G(\mc O_{x})$ où $\mc O_{x}$ désigne la complétion en  $x$ de l'anneau local des fonctions sur  $X$, et $F_{x}$ est son corps des fractions. 

\begin{notation}\label{notation-disque-epointe}
Soit $S$ un  schéma et  $(x_{i})_{i\in I}$ une famille de  $S$-points de $X$.  
On note $\Gamma_{\sum \infty x_i}$ le voisinage formel de  $\cup_{i\in I}\Gamma_{ x_{i}}$ dans $X\times S$. 
Un $G$-torseur $\mc G$ sur $\Gamma_{\sum \infty x_i}$ est la même chose qu'une limite projective de $G$-torseurs sur $\Gamma_{\sum n x_i}$ pour $n$ tendant vers l'infini. 

Soit $J$ une partie de $I$ et $\mc G$ et $\mc G'$ deux $G$-torseurs   sur $\Gamma_{\sum \infty x_i}$. Sans chercher à donner un sens au  voisinage épointé $\Gamma_{\sum \infty x_i}\sm(\bigcup_{i\in J}\Gamma_{x_i})$, on 
définit 
un isomorphisme \begin{gather}\label{abus-G-G'-isom}\phi :\restr{\mc G }{\Gamma_{\sum \infty x_i}\sm(\bigcup_{i\in J}\Gamma_{x_i})}\isom \restr{\mc G' }{\Gamma_{\sum \infty x_i}\sm(\bigcup_{i\in J}\Gamma_{x_i})}\end{gather} comme la donnée pour toute représentation de  dimension finie 
$V$ de $G$ d'un morphisme 
$$V_{\mc G}\to V_{\mc G'}(N\sum _{i\in J} x_{i})$$  où l'ordre $N$ du pôle dépend de $V$, de telle sorte que ces morphismes soient fonctoriels en $V$ et compatibles au produit tensoriel et au dual. 
Cette notation \eqref{abus-G-G'-isom} 
est utilisée aussi par exemple dans 
  \cite{ga-iwahori}
et dans le paragraphe 2.4 de \cite{ga-de-jong}.  
\end{notation}

\begin{construction}\label{rem-grassm}
Soit $S$ un  schéma et  $(x_{i})_{i\in I}$ une famille de  $S$-points de $X$.  
Alors un  point  de $\mr{Gr}_{I}^{(I_{1},...,I_{k})}(S)$ au-dessus 
$(x_{i})_{i\in I}$ équivaut à la donnée
\begin{itemize}
\item  de  $G$-torseurs $\mc G_{0},... , \mc G_{k}$ sur $\Gamma_{\sum \infty x_i}$, 
\item 
d'isomorphismes $\phi_{j}:\restr{\mc G_{j-1}}{\Gamma_{\sum \infty x_i}\sm(\bigcup_{i\in I_{j}}\Gamma_{x_i})}\isom \restr{\mc G_{j}}{\Gamma_{\sum \infty x_i}\sm(\bigcup_{i\in I_{j}}\Gamma_{x_i})}$, 
\item d'une trivialisation $\theta :  \mc G_{k}\isom G$ sur  $\Gamma_{\sum \infty x_i}$. 
\end{itemize}
En effet pour tout  $j\in \{0,...,k\}$, on étend     $\mc G_{j}$ en un   $G$-torseur sur 
$X\times S$ en le recollant  sur  $\Gamma_{\sum \infty x_i}\sm \bigcup_{i\in I}\Gamma_{x_{i}}$ avec le   $G$-torseur trivial  sur  $(X\times S)\sm \bigcup_{i\in I}\Gamma_{x_{i}}$,  grâce à   $\theta\circ \phi_{k}\circ \cdots\circ \phi_{j+1}$. 
Ce recollement est justifié par la remarque  2.3.7 et le théorème  2.12.1 de 
 \cite{hitchin} qui  généralise le  lemme de descente  de  Beauville-Laszlo \cite{BL}. Bien que ce  théorème concerne les modules cohérents, on peut  l'appliquer ici car un  $G$-torseur est aussi un  foncteur tensoriel de la catégorie des représentations de  dimension finie de $G$ vers la  catégorie des fibrés vectoriels (voir la preuve  du théorème  2.3.4 de \cite{hitchin}).  
Quand nous ferons référence à cette   construction, nous dirons qu'un tel  $S$-point de 
$\mr{Gr}_{I}^{(I_{1},...,I_{k})}$ est associé à 
\begin{gather}\label{formule-rem-grassm}\big((x_{i})_{i\in I}, \mc G_{0} \xrightarrow{\phi_{1}}  
\mc G_{1}\xrightarrow{\phi_{2}}
\cdots\xrightarrow{\phi_{k-1}}  \mc G_{k-1} \xrightarrow{ \phi_{k}}   \mc G_{k}\isor{\theta} \restr{G}{\Gamma_{\sum \infty x_i}}\big). \end{gather}
\end{construction}

La  grassmannienne affine de Beilinson-Drinfeld  a été  introduite dans  \cite{hitchin} et joue un  rôle fondamental dans l'équivalence de Satake géométrique \cite{mv} (voir aussi \cite{ga-de-jong}). La fibre  de $\mr{Gr}_{I}^{(I)}$ en un  point géométrique $(x_i)_{i\in I}\in X^{I}$ est le  produit de la  grassmannienne affine habituelle en les points de l'ensemble  $\{x_{i},i\in I\}$. Le cas particulier où $\sharp I=2$ est étudié en détail  pour le produit de  fusion dans  \cite{mv}. Plus généralement, la  grassmannienne affine possède la  propriété de factorisation suivante.

\begin{rem}\label{grassmann-disj-prod}
Soit $\zeta: I\to J $ une application surjective et $(I_{1},...,I_{k})$ une  partition de  $I$. Soit  $U_{\zeta}$ l'ouvert de   $X^{I}$ formé des  $(x_{i})_{i\in I}$ tels que $x_{i}\neq x_{j}$ si $\zeta(i)\neq \zeta(j)$. Alors  $\mr{Gr}_{I}^{(I_{1},...,I_{k})}$  et $\prod_{j\in J}\mr{Gr}_{\zeta^{-1}(\{j\})}^{(I_{1}\cap\zeta^{-1}(\{j\}), ..., I_{k}\cap \zeta^{-1}(\{j\}))}$ sont  canoniquement isomorphes au-dessus de $U_{\zeta}$. 
\end{rem}

Pour tout $(n_{i})_{i\in I}\in \N^{I}$ on note  
$\Gamma_{\sum n_{i}x_{i}}$ le  sous-schéma fermé  de $X\times X^{I}$ dont l'idéal est engendré, localement pour la topologie de Zariski, par $\prod_{i\in I}t_{i}^{n_{i}}$, où $t_{i}$ est une équation du graphe $\Gamma_{x_{i}}$. 
On note  $G_{\sum n_{i}x_{i}}$ le  schéma en groupes lisse sur  $X^{I}$ égal à la restriction à la Weil de  $G$ 
du  sous-schéma fermé  $\Gamma_{\sum_{i\in I} n_{i}x_{i}}\subset X\times X^{I}$ à $X^{I}$.
C'est évidemment un  quotient de $G_{\sum \infty x_i}$ (défini  comme la restriction à la Weil de  $G$ de $\Gamma_{\sum \infty x_i}$ à $X^{I}$). 

On a une action évidente de  $G_{\sum \infty x_i}$ sur $\mr{Gr}_{I}^{(I_{1},...,I_{k})}$ par changement de trivialisation de $\mc G_{k}$. Autrement dit, pour tout schéma $S$ et toute section 
$\gamma$ de  $G$ sur  $\Gamma_{\sum \infty x_i}$, l'action de  $\gamma$ 
 envoie  \eqref{formule-rem-grassm}  
sur   
\begin{gather}\nonumber
\big( (x_i)_{i\in I}, \mc G_{0} \xrightarrow{\phi_{1}}  \mc G_{1} \xrightarrow{\phi_{2}}
\cdots\xrightarrow{\phi_{k}}  \mc G_{k}\isor{L_{\gamma}\circ \theta}  \restr{G}{\Gamma_{\sum \infty x_i}}\big)
\end{gather}
(où $L_{\gamma}$ est l'automorphisme du  $G$-torseur trivial donné par la multiplication à gauche par $\gamma$).  

\begin{prop}\label{action-Gnx-Grass}
Si les entiers   $n_{i}$  sont  assez grands  en fonction des copoids  $\omega_{i}$, cette action de $G_{\sum \infty x_i}$ sur $\mr{Gr}_{I,\on{\lesssim}\underline{\omega}}^{(I_{1},...,I_{k})}$ se factorise à travers le quotient $G_{\sum n_{i}x_{i}}$. \cqfd
\end{prop}

\begin{rem}\label{rem-apres-action-Gnx-Grass}
Par conséquent, si les entiers   $n_{i}$ sont comme dans la proposition précédente, pour tout schéma $S$, 
à la donnée de  
\begin{itemize}
\item [] a) $S$-points 
$(x_i)_{i\in I}$,
\item [] b) $\mc G_{0} \xrightarrow{\phi_{1}}  \mc G_{1} \xrightarrow{\phi_{2}}
\cdots\xrightarrow{\phi_{k}}  \mc G_{k}$, où les $ \mc G_{i}$ sont des  $G$-torseurs sur $\Gamma_{\sum \infty x_i}$ et les $\phi_{i}$ satisfont la  condition 
iii)$_{\underline{\omega}}$ de la 
\defiref{defi-Hecke},  
\item [] c) une  trivialisation de $\restr{\mc G_{k}}{\Gamma_{\sum n_{i} x_i}}$,
\end{itemize}
on associe un $S$-point de $\mr{Gr}_{I,\on{\lesssim}\underline{\omega}}^{(I_{1},...,I_{k})}$ (par la construction \ref{rem-grassm}, après avoir étendu de manière arbitraire la trivialisation de $\mc G_{k}$ de $\Gamma_{\sum n_{i} x_i}$
à $\Gamma_{\sum \infty x_i}$). 
Par conséquent,  aux données  de a) et b) 
on associe un $S$-point de $\mr{Gr}_{I,\on{\lesssim}\underline{\omega}}^{(I_{1},...,I_{k})}/G_{\sum n_{i}x_{i}}$, sur lequel  le 
$G_{\sum n_{i}x_{i}}$-torseur tautologique est    $\restr{\mc G_{k}}{\Gamma_{\sum n_{i} x_i}}$. 
Dans la suite ce  $S$-point sera souvent noté 
$$(\mc G_{0} \xrightarrow{\phi_{1}}  \mc G_{1} \xrightarrow{\phi_{2}}
\cdots\xrightarrow{\phi_{k}}  \mc G_{k})\in \mr{Gr}_{I,\on{\lesssim}\underline{\omega}}^{(I_{1},...,I_{k})}/G_{\sum n_{i}x_{i}}.$$ 
\end{rem}

On va définir maintenant le   sous-schéma fermé  
$\mr{Gr}_{I,\underline{\omega}}^{(I_{1},...,I_{k})}\subset 
\mr{Gr}_{I,\on{\lesssim}\underline{\omega}}^{(I_{1},...,I_{k})}$. 
On commence par rappeler que si   $I$ est un singleton, $\mr{Gr}_{I}^{(I)}$ est un schéma sur $X$ dont la fibre en un  point géométrique $x$ de $X$ est la  grassmannienne affine habituelle $\mr{Gr}_{x}$, c'est-à-dire  le quotient  fpqc   $G(F_{x})/G(\mc O_{x})$. Il est bien connu que les  $G(\mc O_{x})$-orbites dans   $\mr{Gr}_{x}$ sont des  sous-schémas  localement fermés 
indexés par les copoids dominants $\omega$ de $G$. Plus précisément  on note  $\mr{Gr}_{x,\omega}^{0}$ la    $G(\mc O_{x})$-orbite  de $t^{\omega}$ où $t$ est une uniformisante de $\mc O_{x}$. 
On note $\mr{Gr}_{x,\omega} $ son  adhérence de Zariski, 
qui est la réunion des orbites $\mr{Gr}^{0}_{x,\lambda}$ telles que  $\lambda$ soit un poids de la représentation   irréductible de $\wh G$ de plus haut poids $\omega$. 

Habituellement $\mr{Gr}_{x,\omega} $ est noté $\overline{\mr{Gr}_{x,\omega} }$ et $\mr{Gr}^{0}_{x,\lambda}$  est noté $\mr{Gr} _{x,\lambda}$. 
Les notations utilisées ici sont mieux adaptées pour nous 
(car nous utilisons presque exclusivement les adhérences) 
et elles sont compatibles avec    \cite{var}. 

\begin{defi}
Quand  $I$ est un  singleton on note  $\mr{Gr}_{I,\omega}^{(I)} $ le  {\it   sous-schéma fermé réduit} de $\mr{Gr}_{I}^{(I)} $ dont la fibre en un  point $x\in X$ est $ \mr{Gr}_{x,\omega} $. 
En général soit $U\subset X^{I}$ le complémentaire de toutes les diagonales, c'est-à-dire  
$U=\{(x_{i})_{i\in I}, \forall i\neq j, x_{i}\neq x_{j}\}$. 
On a  un isomorphisme 
$$ \theta: \restr{\mr{Gr}_{I,\on{\lesssim}\underline{\omega}}^{(I_{1},...,I_{k})}}{U}
\isom \restr{\prod_{i\in I}\mr{Gr}_{\{i\},\on{\lesssim}\omega_{i}}^{(\{i\})}}{U}
$$ grâce à la remarque \ref{grassmann-disj-prod}. On définit alors 
$\mr{Gr}_{I,\underline{\omega}}^{(I_{1},...,I_{k})}$ comme le   sous-schéma fermé  réduit de $\mr{Gr}_{I,\on{\lesssim}\underline{\omega}}^{(I_{1},...,I_{k})}$ égal à l'adhérence de Zariski de 
$\theta^{-1}\Big(\restr{\prod_{i\in I}\mr{Gr}_{\{i\},\omega_{i}}^{(\{i\})}}{U}\Big)$.
\end{defi}

Il est bien connu qu'il existe un isomorphisme 
\begin{gather}\label{defi-beta-avant-Grass-infini}\beta_{I,\infty }^{(I_{1},...,I_{k})}:  \Hecke_{I}^{(I_{1},...,I_{k})}\isom \Big(
\mr{Gr}_{I}^{(I_{1},...,I_{k})}\times_{X^{I}} \Bun_{G,\sum \infty x_{i}}\Big)/G_{\sum \infty x_{i}}\end{gather}
où $G_{\sum \infty x_{i}}$ agit diagonalement. En fait on utilisera plutôt  
 l'isomorphisme \eqref{defi-beta-avant-Grass-lissite-lesssim} ci-dessous  où les modifications sont bornées  et $G_{\sum \infty x_{i}}$ est remplacé par $G_{\sum n_{i}x_{i}}$. 
L'avantage est que l'on reste dans le cadre des champs d'Artin et surtout que  l'isomorphisme \eqref{defi-beta-avant-Grass-lissite-lesssim}  va permettre dans le  \lemref{Hecke-Grass-lissite} de construire un  morphisme lisse, qui serait de dimension infinie si on avait utilisé 
$G_{\sum \infty x_{i}}$.

Soit $\underline n=(n_{i})_{i\in I}\in \N^{I}$ 
comme dans la proposition \ref{action-Gnx-Grass}.  
On a alors  un  isomorphisme 
\begin{gather}\label{defi-beta-avant-Grass-lissite-lesssim}\beta_{I,\on{\lesssim}\underline{\omega},\underline n}^{(I_{1},...,I_{k})}:  \Hecke_{I,\on{\lesssim}\underline{\omega}}^{(I_{1},...,I_{k})}\isom \Big(
\mr{Gr}_{I,\on{\lesssim}\underline{\omega}}^{(I_{1},...,I_{k})}\times_{X^{I}} \Bun_{G,\sum n_{i}x_{i}}\Big)/G_{\sum n_{i}x_{i}}\end{gather} (où $G_{\sum n_{i}x_{i}}$ agit diagonalement). 
Il est défini de la manière suivante. D'abord le  $G_{\sum n_{i}x_{i}}$-torseur  $\mc A$ sur  $\Hecke_{I,\on{\lesssim}\underline{\omega}}^{(I_{1},...,I_{k})}$ associé à ce morphisme est $\restr{\mc G_{k}}{\Gamma_{\sum n_{i}x_{i}}}$, pour 
$$\big((x_{i})_{i\in I}, \mc G_{0} \xrightarrow{\phi_{1}}  
\mc G_{1}\xrightarrow{\phi_{2}}
\cdots \xrightarrow{ \phi_{k}}   \mc G_{k}\big)\in \Hecke_{I,\on{\lesssim}\underline{\omega}}^{(I_{1},...,I_{k})}(S). $$ 
Autrement dit  l'espace total de $\mc A$ est constitué par les 
 trivialisations    $\kappa$ de
 $\restr{\mc G_{k}}{\Gamma_{\sum n_{i}x_{i}}}$. 
Ensuite  le   morphisme $\mc A\to \mr{Gr}_{I,\on{\lesssim}\underline{\omega}}^{(I_{1},...,I_{k})}\times_{X^{I}} \Bun_{G,\sum n_{i}x_{i}}$ envoie  
$$\Big(\big((x_{i})_{i\in I}, \mc G_{0} \xrightarrow{\phi_{1}}  
\mc G_{1}\xrightarrow{\phi_{2}}
\cdots \xrightarrow{ \phi_{k}}   \mc G_{k}\big), \kappa: \restr{\mc G_{k}}{\Gamma_{\sum n_{i}x_{i}}}\isom \restr{G}{\Gamma_{\sum n_{i}x_{i}}}\Big)$$ sur le  produit  
\begin{itemize}
\item du point de $\mr{Gr}_{I,\on{\lesssim}\underline{\omega}}^{(I_{1},...,I_{k})}$ associé par la 
\remref{rem-apres-action-Gnx-Grass} à
\begin{gather*}\Big((x_{i})_{i\in I}, \restr{\mc G_{0}}{\Gamma_{\sum \infty x_i}} \xrightarrow{\phi_{1}}  
\restr{\mc G_{1}}{\Gamma_{\sum \infty x_i}}\xrightarrow{\phi_{2}}
\cdots \xrightarrow{ \phi_{k}}   \restr{\mc G_{k}}{\Gamma_{\sum \infty x_i}} \Big) \text{ \ et  \ }\kappa, 
\end{gather*} 
  \item de $\big(\mc G_{k}, \kappa\big)$ dans $\Bun_{G,\sum n_{i}x_{i}}$. 
 \end{itemize} 

\begin{defi} On définit $ \Hecke_{I,\underline{\omega}}^{(I_{1},...,I_{k})}$ comme  l'image inverse  de 
$$  \Big(
\mr{Gr}_{I,\underline{\omega}}^{(I_{1},...,I_{k})}\times_{X^{I}} \Bun_{G,\sum n_{i}x_{i}}\Big)/G_{\sum n_{i}x_{i}}$$
par l'isomorphisme $\beta_{I,\on{\lesssim}\underline{\omega},\underline n}^{(I_{1},...,I_{k})}$. \end{defi}

Par définition on a donc un isomorphisme 
\begin{gather}\label{defi-beta-avant-Grass-lissite}
\beta_{I,\underline{\omega},\underline n}^{(I_{1},...,I_{k})}: 
 \Hecke_{I,\underline{\omega}}^{(I_{1},...,I_{k})}
 \isom 
 \Big(\mr{Gr}_{I,\underline{\omega}}^{(I_{1},...,I_{k})}\times_{X^{I}}
  \Bun_{G,\sum n_{i}x_{i}}\Big)/G_{\sum n_{i}x_{i}}\end{gather} 
  (où $G_{\sum n_{i}x_{i}}$ agit diagonalement). 
L'isomorphisme inverse de $\beta_{I,\underline{\omega},\underline n}^{(I_{1},...,I_{k})}$ envoie 
$(y,(\mc G, \kappa))\in \mr{Gr}_{I,\underline{\omega}}^{(I_{1},...,I_{k})}\times_{X^{I}}
  \Bun_{G,\sum n_{i}x_{i}}$ sur $\big((x_{i})_{i\in I}, \mc G_{0} \xrightarrow{\phi_{1}}  
\mc G_{1}\xrightarrow{\phi_{2}}
\cdots \xrightarrow{ \phi_{k}}   \mc G_{k}\big)$ où $ \mc G_{k}= \mc G$ et 
$(\mc G_{0} \xrightarrow{\phi_{1}}  
\mc G_{1}\xrightarrow{\phi_{2}}
\cdots \xrightarrow{ \phi_{k}}   \mc G_{k})$ est la modification de $\mc G$ associée au point $y\in \mr{Gr}_{I,\underline{\omega}}^{(I_{1},...,I_{k})}$ grâce à la structure de niveau $\kappa$. On vérifie que ce morphisme se factorise par le quotient par l'action diagonale de $G_{\sum n_{i}x_{i}}$.

Comme $\Bun_{G,\sum n_{i}x_{i}}$ est lisse sur $X^{I}$, cela entraîne immédiatement le lemme suivant. 
 
\begin{lem}\label{Hecke-Grass-lissite} 
Soit $\underline n=(n_{i})_{i\in I}\in \N^{I}$ comme dans la proposition précédente. 
 Alors le  morphisme \begin{gather}\label{morph-delta-defi-Grass}\delta_{(I),\underline{\omega},\underline n}^{(I_{1},...,I_{k})}: \Hecke_{I,\underline{\omega}}^{(I_{1},...,I_{k})} \to \mr{Gr}_{I,\underline{\omega}}^{(I_{1},...,I_{k})}/G_{\sum n_{i}x_{i}}\end{gather}
égal à la première composante de  $\beta_{I,\underline{\omega},\underline n}^{(I_{1},...,I_{k})}$ est lisse. \cqfd
 \end{lem}

 Le lemme suivant  apparaît dans  \cite{mv} et  résulte aussi des lemmes 3.1 et A.12 de \cite{var}. 

\begin{lem} \label{lemme-oubli}
Le morphisme $\pi^{(I_{1},...,I_{k})}_{(I'_{1},...,I'_{k'})}:\mr{Gr}_{I,\underline{\omega}}^{(I_{1},...,I_{k})}\to \mr{Gr}_{I,\underline{\omega}}^{(I'_{1},...,I'_{k'})}$  est  projectif, surjectif et petit. \cqfd
\end{lem}

On a un morphisme évident 
\begin{gather}\label{def-kappa-eps-J}\kappa_{I,(\omega_{i})_{i\in I}}^{(I_{1},...,I_{k})}: 
\mr{Gr}_{I,(\omega_{i})_{i\in I}}^{(I_{1},...,I_{k})}
\to 
\prod_{j=1}^{k}
\Big(\mr{Gr}_{I_{j},(\omega_{i})_{i\in I_{j}}}^{(I_{j})}/G_{\sum _{i\in I_{j}}n_{i}x_{i}}\Big)\end{gather}
qui envoie 
\begin{gather}\nonumber
\big( (x_i)_{i\in I}, \mc G_{0} \xrightarrow{\phi_{1}}  \mc G_{1} \xrightarrow{\phi_{2}}
\cdots\xrightarrow{\phi_{k}}  \mc G_{k}\isor{\theta} \restr{G}{\Gamma_{\sum \infty x_i}}\big)
\end{gather} 
sur le  produit des points $\big(\mc G_{j-1} \xrightarrow{\phi_{j}}  \mc G_{j}\big) \in 
\mr{Gr}_{I_{j},(\omega_{i})_{i\in I_{j}}}^{(I_{j})}/G_{\sum _{i\in I_{j}}n_{i}x_{i}}$
 (avec les  notations de la construction \ref{rem-grassm} et de la \remref{rem-apres-action-Gnx-Grass}, c'est-à-dire  
 que le  $G_{\sum _{i\in I_{j}}n_{i}x_{i}}$-torseur tautologique est $\restr{\mc G_{j}}
{\Gamma_{\sum _{i\in I_{j}}n_{i}x_{i}}}$). 
Si les  entiers  $m_{i}$ sont assez grands en fonction des $n_{i}$ et des $\omega_{i}$, il se factorise à travers un  morphisme 
\begin{gather}\label{kappa-tilde-I1-Ik}\wt \kappa_{I,(\omega_{i})_{i\in I}}^{(I_{1},...,I_{k})}: 
\mr{Gr}_{I,(\omega_{i})_{i\in I}}^{(I_{1},...,I_{k})}/G_{\sum _{i\in I}m_{i}x_{i}}\to 
\prod_{j=1}^{k}
\Big(\mr{Gr}_{I_{j},(\omega_{i})_{i\in I_{j}}}^{(I_{j})}/G_{\sum _{i\in I_{j}}n_{i}x_{i}}\Big). \end{gather}

\begin{lem}\label{kapp-eps-J-lemme-smooth}
Les morphismes $\kappa_{I,(\omega_{i})_{i\in I}}^{(I_{1},...,I_{k})}$ et 
$\wt \kappa_{I,(\omega_{i})_{i\in I}}^{(I_{1},...,I_{k})}$ sont lisses. 
\end{lem}
\dem Il suffit de prouver la lissité  de $\kappa_{I,(\omega_{i})_{i\in I}}^{(I_{1},...,I_{k})}$. 
On la démontre par récurrence sur $k$: 
$\kappa$ se factorise à travers le  morphisme évident 
$$\mr{Gr}_{I,(\omega_{i})_{i\in I}}^{(I_{1},...,I_{k})}\to 
\prod_{j=1}^{k-1}
\Big(\mr{Gr}_{I_{j},(\omega_{i})_{i\in I_{j}}}^{(I_{j})}/G_{\sum _{i\in I_{j}}n_{i}x_{i}}\Big)
\times \mr{Gr}_{I_{k},(\omega_{i})_{i\in I_{k}}}^{(I_{k})}
$$ qui est essentiellement une version tordue en famille, paramétrée par $\mr{Gr}_{I_{k},(\omega_{i})_{i\in I_{k}}}^{(I_{k})}$, de morphismes similaires à   
$\kappa_{I\setminus I_{k},(\omega_{i})_{i\in I\setminus I_{k}}}^{(I_{1},...,I_{k-1})}$. 
\cqfd

Nous rappelons maintenant l'équivalence de Satake géométrique. Les notations suivantes sont commodes. Si $W$ est une  représentation
$E$-linéaire  de dimension finie   de $(\wh G)^{I}$ on la décompose en somme d'irréductibles: 
\begin{gather}
\label{defi-W-Elin}
W=\bigoplus _{\underline{\omega}\in (X_{*}^{+}(T))^{I}} \Big(\bigotimes_{i\in I}V_{\omega_{i}} \Big)\otimes_{E} \mf W_{\underline{\omega}}
 \end{gather} où les  $\mf W_{\underline{\omega}}$ sont des $E$-espaces vectoriels de  dimension  finie, et sont presque tous nuls. On définit alors 
$\mr{Gr}_{I,W}^{(I_{1},...,I_{k})}$ comme la réunion  de 
tous les  sous-schémas $\mr{Gr}_{I,\underline{\omega}}^{(I_{1},...,I_{k})}\subset \mr{Gr}_{I}^{(I_{1},...,I_{k})}$ pour $\underline{\omega}=(\omega_{i})_{i\in I}$ tel que $\mf W_{\underline{\omega}}$ est non nul. 

Pour tout point géométrique $x\in X$, l'équivalence de Satake géométrique
 de  Lusztig, Drinfeld, Ginzburg, et Mirkovic--Vilonen \cite{lusztig-satake,ginzburg,hitchin,mv}  est une équivalence de la   catégorie des représentations de dimension finie  de $\wh G$ vers la  catégorie des  faisceaux pervers $G(\mc O_{x})$-équivariants 
 sur la  grassmannienne affine $\mr{Gr}_{x}=G(F_{x})/G(\mc O_{x})$. De plus cette  équivalence est un foncteur tensoriel, lorsque le but est muni du produit de convolution, ou de fusion (avec la modification de la contrainte de commutativité qui sera rappelée dans la \remref{contrainte-commutativite}).   
   
D'après \cite{mv,hitchin,ga-de-jong,ga-iwahori,var} nous avons  le théorème suivant, où nous nous limitons à ce qui est nécessaire pour notre  article: un foncteur de la  catégorie de  représentations de dimension finie de $(\wh G)^{I}$ 
 vers la  catégorie des 
 faisceaux pervers $G_{\sum \infty x_{i}}$-équivariants sur 
$\mr{Gr}_{I}^{(I_{1},...,I_{k})}$, vérifiant certaines propriétés. 
Dans  \cite{ga-de-jong} Gaitsgory énonce un résultat plus fort, qui décrit cette dernière catégorie uniquement à l'aide de  $\wh G$. 
 Dans le théorème suivant comme dans le reste de l'article nous dirons  
 ``à décalage près'' au lieu de ``à décalage et torsion à la Tate près''. 
 On va énoncer le théorème suivant avec des coefficients dans 
 $E$ ou $\mc O_{E}$  mais le cas de $\mc O_{E}$ ne servira que dans le chapitre \ref{mod-ell} et peut donc être sauté   en première lecture. 
 
\begin{thm}\label{thm-geom-satake} (Un sens de l'équivalence de  Satake géométrique, \cite{mv,hitchin,ga-de-jong}) Soit $\mc A$ égal à $E$ ou $\mc O_{E}$. 
On a un foncteur canonique 
$$W\mapsto \mc S_{I,W,\mc A}^{(I_{1},...,I_{k})}$$
de la catégorie des  représentations $\mc A$-linéaires de type fini de $(\wh G)^{I}$ vers la  catégorie des    $\mc A$-faisceaux pervers  
$G_{\sum \infty x_{i}}$-équivariants sur $\mr{Gr}_{I }^{(I_{1},...,I_{k})}$. 

De plus 
si $\mc A=E$, 
$\mc S_{I,W,\mc A}^{(I_{1},...,I_{k})}$ est supporté par 
$\mr{Gr}_{I,W}^{(I_{1},...,I_{k})}$ (défini juste après \eqref{defi-W-Elin}) et on peut donc le considérer comme un faisceau pervers 
(à un décalage près)   sur  $\mr{Gr}_{I,W}^{(I_{1},...,I_{k})}/G_{\sum n_{i}x_{i}}$  (où les entiers $n_{i}$ sont assez grands  et le décalage est déterminé par la condition que  l'image inverse sur  $\mr{Gr}_{I,W}^{(I_{1},...,I_{k})}$ est perverse relativement à  $X^{I}$). 
Cela reste vrai si $\mc A=\mc O_{E}$, avec $\mr{Gr}_{I,W}^{(I_{1},...,I_{k})}$ un sous-schéma fermé de type fini assez grand de $\mr{Gr}_{I}^{(I_{1},...,I_{k})}$, dépendant de $W$ d'une fa\c con qui ne sera pas précisée ici. 

Les $\mc S_{I,W,\mc A}^{(I_{1},...,I_{k})}$ sont universellement localement acycliques par rapport au morphisme vers $X^{I}$. De plus ils vérifient les  propriétés suivantes:  
\begin{itemize}
\item[] a) 
 Lorsque  $\mc A=E$ et que 
$W=\boxtimes_{i\in I} V_{\omega_{i}}$ est irréductible, le faisceau pervers (à décalage près) 
 $\mc S_{I,W,\mc A}^{(I_{1},...,I_{k})}$ sur 
$ \mr{Gr}_{I,\underline{\omega}}^{(I_{1},...,I_{k})}/G_{\sum _{i\in I}n_{i}x_{i}} $ est le faisceau d'intersection    de 
$\mr{Gr}_{I,\underline{\omega}}^{(I_{1},...,I_{k})}$, avec la  normalisation perverse relative à  $X^{I}$ et la structure  $G_{\sum _{i\in I}n_{i}x_{i}}$-équivariante naturelle. 
\item[] b)   On a un  isomorphisme canonique 
$$\Big(\pi^{(I_{1},...,I_{k})}_{(I'_{1},...,I'_{k'})}\Big)_{!}
\Big(\mc S_{I,W,\mc A}^{(I_{1},...,I_{k})}\Big)\simeq \mc S_{I,W,\mc A}^{(I'_{1},...,I'_{k'})},$$
\item[] c) Si $W=\boxtimes_{j\in \{1,...,k\}} W_{j}$ où $W_{j}$ est une  représentation  de $(\wh G)^{I_{j}}$ sur un $\mc A$-module libre de  type fini, 
on a un  isomorphisme canonique 
$$\mc S_{I,W,\mc A}^{(I_{1},...,I_{k})}\simeq \Big(
\wt \kappa_{I,W}^{(I_{1},...,I_{k})}\Big)^{*}
\Big(\boxtimes _{j\in \{1,...,k\}} \mc S_{I_{j},W_{j},\mc A}^{(I_{j})}\Big) $$
où on applique la contrainte de commutativité modifiée de \cite{mv} qui sera rappelée dans la \remref{contrainte-commutativite} ci-dessous, 
\item[] d)  Soient  $I,J$ des ensembles finis et  $\zeta: I\to J$ une application. 
On note  $$\Delta_{\zeta} : X^{J}\to X^{I},  \ \ (x_{j})_{j\in J}\mapsto (x_{\zeta(i)})_{i\in I}$$ le   morphisme diagonal associé à 
  $\zeta$. 
   Soit  $W$ une   représentation  $\mc A$-linéaire de type fini  de  $\wh G^{I}$. On note  $W^{\zeta}$ la   représentation de $\wh G^{J}$ qui est la composée de la   représentation $W$ avec le    morphisme diagonal $$\wh G^{J}\to \wh G^{I},  \ \ (g_{j})_{j\in J}\mapsto (g_{\zeta(i)})_{i\in I}  .$$ Soit $(J_{1},...,J_{k})$ une  partition de $J$. 
   On en déduit une  partition $(I_{1},...,I_{k})$ de $I$  en posant 
     $I_{j}=\zeta^{-1}(J_{j})$ pour tout $j$. On a alors un  isomorphisme canonique 
 \begin{gather}\label{coalescence-Gr-section1-thm} \Delta_{\zeta}^{*}\Big( \mc S_{I,W,\mc A}^{(I_{1},...,I_{k})} \Big)\simeq 
 \mc S_{J,W^{\zeta},\mc A}^{(J_{1},...,J_{k})}\end{gather} 
 où $\Delta_{\zeta}$ désigne (le  quotient par $G_{\sum n_{i}x_{i}}$ de) l'inclusion 
 $$\mr{Gr}_{J}^{(J_{1},...,J_{k})} =\mr{Gr}_{I}^{(I_{1},...,I_{k})}\times _{X^{I}}X^{J}\hookrightarrow \mr{Gr}_{I}^{(I_{1},...,I_{k})}.$$
De plus  \eqref{coalescence-Gr-section1-thm} est fonctoriel en  $W$ et compatible avec la composition pour $\zeta$. 
 \end{itemize}  
\end{thm}

\begin{rem}  Le cas intéressant dans la condition d) est celui où $\zeta$ est surjective. Il exprime   la compatibilité avec le produit de   fusion dans la grassmannienne affine de Beilinson-Drinfeld, c'est-à-dire le fait que  l'équivalence de Satake géométrique est un  foncteur tensoriel.  Les propriétés a) et b) auraient pu être énoncées avec des partitions plus générales, mais au prix de notations beaucoup plus lourdes. 
 \end{rem}
 
  \begin{rem}\label{contrainte-commutativite}
 La contrainte de commutativité est définie par la convention  {\it modifiée},  introduite dans la discussion qui précède la   proposition 6.3 de \cite{mv}. C'est seulement avec cette   convention modifiée que la catégorie tensorielle des  faisceaux pervers  $G(\mc O)$-équivariants   sur la   grassmannienne affine   (munie du   produit de  fusion) est équivalente à  la catégorie tensorielle des représentations de  type fini de $\wh G$, et c'est elle qui intervient dans 
   l'énoncé de la correspondance de   Langlands  géométrique. Pour plus de détails, voir \cite{mv} et \cite{hitchin}. Voici un bref rappel. 
   Il se trouve que 
   \begin{itemize}
   \item pour toute composante  connexe  de la grassmannienne affine, les strates sont toutes de dimension paire ou toutes de dimension impaire (on parle alors de composante paire ou impaire), 
   \item d'après le lemme 3.9 de  \cite{mv}, si un faisceau $\mc S_{I_{j},W_{j},\mc A}^{(I_{j})}$ est supporté sur une composante  paire ({\it resp.} impaire),   sa cohomologie totale est concentrée en degrés cohomologiques pairs ({\it resp.} impairs). 
   \end{itemize}
  La contrainte de commutativité modifée consiste à   ajouter aux signes habituels donnés par les règles de Koszul  un signe moins lorsque l'on permute deux faisceaux 
$\mc S_{I_{j},W_{j},\mc A}^{(I_{j})}$ et $\mc S_{I_{j'},W_{j'},\mc A}^{(I_{j'})}$ supportés 
sur deux composantes connexes impaires de la grassmannienne affine.  
Autrement dit c'est la contrainte de commutativité naturelle que l'on aurait si on 
normalisait les faisceaux $\mc S_{I_{j},W_{j},\mc A}^{(I_{j})}$ pour que leur cohomologie totale soit en degré pair. Le foncteur fibre donné par la cohomologie totale, de la catégorie des faisceaux pervers $G(\mc O)$-équivariants sur la grassmannienne affine, munie du produit de fusion (avec cette contrainte de commutativité modifée) vers la catégorie des $\mc A$-modules  de type fini est donc tensoriel, d'où une équivalence   
  avec  la catégorie tensorielle $\on{Rep}(\wh G)$  des représentations de  $\wh G$  sur des $\mc A$-modules de type fini (où $\wh G$ est  défini comme le groupe d'automorphismes du foncteur fibre).  
   \end{rem}
     
 \dem
 Le théorème résulte essentiellement du 
  théorème 2.6 de  \cite{ga-de-jong} (qui montre un résultat beaucoup plus fort). 
  L'argument que nous allons donner est inspiré  de la preuve du  théorème 2.6 de  \cite{ga-de-jong} (donnée dans l'appendice B de \cite{ga-de-jong}). 
  Plus précisément notre argument est identique à celui de \cite{ga-de-jong} quand $I$ est un singleton, et plus simple en général (ce qui est normal puisque nous montrons un résultat plus faible). 
    Mirkovic et Vilonen \cite{mv} associent à toute représentation $\mc A$-linéaire de type fini $W$ de $\wh G$  un faisceau pervers 
  $G_{\infty x}$-équivariant sur $\mr{Gr}_{\{0\}}^{(\{0\})}$. 
  On en déduit par convolution ou fusion, pour toute famille $(W_{i})_{i\in I}$ 
  de représentations $\mc A$-linéaires de type fini de $\wh G$  un  $\mc A$-faisceau  pervers  
$G_{\sum \infty x_{i}}$-équivariant  
$\mc S_{I,\boxtimes W_{i},\mc A}^{(I_{1},...,I_{k})}$
sur $\mr{Gr}_{I }^{(I_{1},...,I_{k})}$ et de plus ces faisceaux vérifient les propriétés a), b), c), d) dans le cas particulier des représentations $W$ de $(\wh G)^{I}$ 
de la forme $\boxtimes_{i\in I} W_{i}$. Soit $\mc R$ l'algèbre des fonctions régulières sur $\wh G$ (à coefficients dans $\mc A$), considérée comme ind-objet de $\on{Rep}(\wh G)$ par l'action de $\wh G$ sur lui-même par multiplication à gauche (on notera qu'on peut, si on veut,  écrire 
$\mc R$ comme une limite inductive de représentations de $\wh G$ sur des 
$\mc A$-modules libres de type fini). Alors 
$\boxtimes_{i\in I} \mc R$  est l'algèbre des fonctions régulières sur $\wh G^{I}$, 
considérée comme ind-objet de $\on{Rep}(\wh G^{I})$ par l'action de $\wh G^{I}$ sur lui-même par multiplication à gauche. 
On définit alors  
$$\mc S_{I,W,\mc A}^{(I_{1},...,I_{k})}=
\Big( \mc S_{I,\boxtimes _{i\in I}\mc R,\mc A}^{(I_{1},...,I_{k})} \otimes W\Big)^{\wh G^{I}}
.$$
Dans la formule précédente l'action de $\wh G^{I}$ sur $\boxtimes_{i\in I} \mc R$  
par multiplication à droite munit $ \mc S_{I,\boxtimes _{i\in I}\mc R,\mc A}^{(I_{1},...,I_{k})}$ d'une action de $\wh G^{I}$, et on prend les invariants par l'action diagonale de  $\wh G^{I}$ sur $\mc S_{I,\boxtimes_{i\in I}  \mc R,\mc A}^{(I_{1},...,I_{k})} \otimes W$. 
Les propriétés a), b), c), d) se déduisent aisément de cette formule. 
Pour d) on utilise le fait que pour toute application $I\to J$, $\boxtimes_{i\in I} \mc R$  
considéré comme ind-objet de $\on{Rep}(\wh G^{I})$ (par l'action de $\wh G^{I}$ sur lui-même par multiplication à gauche) est égale à l'induite de $\boxtimes_{i\in J} \mc R$  
 par le morphisme $\wh G^{J}\to \wh G^{I}$ (si $I\to J$ n'est pas surjective, ce morphisme n'est pas injectif et pour induire on commence par prendre les invariants par son noyau). De plus cette égalité est compatible avec les actions à droite de $\wh G^{J}$ et $\wh G^{I}$. 
 \cqfd
 
   \begin{rem}\label{rem-Satake-Gad-texte}
    Dans le théorème précédent $Z_{\sum_{i\in I} \infty x_{i}}\subset 
G_{\sum_{i\in I} \infty x_{i}}$ agit trivialement sur 
 $ \mr{Gr}_{I}^{(I_{1},...,I_{k})}$ et donc   sur tous les faisceaux $\mc S_{I,W,\mc A }^{(I_{1},...,I_{k})}$ (en effet il suffit de le démontrer pour $I$ singleton et alors 
 cela résulte de \cite{mv}, parce que les  projectifs $P_{Z}(\nu, \mc A)$ permettent de reconstruire toute la catégorie des faisceaux pervers $G(\mc O)$-équivariants sur la grassmannienne affine et que $Z(\mc O)$ agit trivialement sur eux d'après la formule (9.9) de \cite{mv}). On note  $G^{\mr{ad}}=G/Z$. On peut donc considérer $\mc S_{I,W,\mc A}^{(I_{1},...,I_{k})}$  comme un 
 faisceau pervers  (à un décalage près)
$G^{\mr{ad}}_{\sum \infty x_{i}}$-équivariant sur $\mr{Gr}_{I }^{(I_{1},...,I_{k})}$
 ou si on préfère comme  
 un faisceau pervers 
(à un décalage près)   sur  $\mr{Gr}_{I,W}^{(I_{1},...,I_{k})}/G^{\mr{ad}}_{\sum n_{i}x_{i}}$  (avec les entiers $n_{i}$   assez grands).   
    \end{rem}

 Pour toute place $v\in |X|$  et pour toute   représentation irréductible $V$ de  $\wh G$,  on note    $h_{V,v}\in C_{c}(G(\mc O_v)\backslash G(F_{v})/G(\mc O_v), \mc O_{E})$ la fonction sphérique associée à  $V$ (ou si on préfère au  caractère $\chi_{V}$) par l'isomorphisme de  Satake classique. La compatibilité entre l'isomorphisme de  Satake classique et l'équivalence de Satake géométrique s'exprime par le fait que, en notant $ \omega$ le plus haut poids de $V$ et 
 $\rho$ la demi-somme des coracines positives de $\wh G$, 
 $(-1)^{\s{2\rho, \omega}}h_{V,v}$ est égale à la  trace de  $\Frob_{\mr{Gr}_{v}/k(v)} $ sur le faisceau pervers $\restr{\mc S_{\{0\},V,E}^{(\{0\})}}{\mr{Gr}_{v}}$ 
 (qui est le faisceau d'intersection  de la strate fermée $\mr{Gr}_{v,\omega}$, où $\omega$ est le plus haut poids de $V$). Les fonctions  $h_{V,v}$  forment une base de    
    $C_{c}(G(\mc O_v)\backslash G(F_{v})/G(\mc O_v), \mc O_{E})$ sur $\mc O_{E}$ (cela est même vrai sur  $\Z[q^{1/2},q^{-1/2}]$, cf la  proposition 3.6 de \cite{gross}).

 \begin{rem}
 Le signe $(-1)^{\s{2\rho, \omega}}$   dans la définition de $h_{V,v}$ ci-dessus 
 est justifié par le fait que la strate $\mr{Gr}_{v,\omega}$ est de dimension $\s{2\rho, \omega}$. Autrement dit $h_{V,v}$ serait la trace de  $\Frob_{\mr{Gr}_{v}/k(v)} $ sur  $\restr{\mc S_{\{0\},V,E}^{(\{0\})}}{\mr{Gr}_{v}}$ 
si on l'avait normalisé pour que   sa cohomologie totale soit supportée en degrés cohomologiques pairs (mais sans changer la torsion à la Tate).
Ce choix  est cohérent avec  le fait que la contrainte de commutativité modifiée de \cite{mv} est celle que l'on obtiendrait naturellement avec de telles  normalisations, ainsi qu'on l'a rappelé dans la \remref{contrainte-commutativite}. La cohérence  de ce choix est la raison pour laquelle   il n'y aura pas de signe   dans l'égalité de la \propref{prop-coal-frob-cas-part}. 
 \end{rem}

\section{Champs de   $G$-chtoucas}\label{para-1-chtoucas}
 
 Ce chapitre est entièrement extrait de  \cite{var}, à l'exception des  propositions
   \ref{lissite-Cht-Grass} et \ref{lissite-Cht-Grass-gen}.  
On garde les notations du chapitre précédent.

Voici la  définition des champs classifiants de  $G$-chtoucas. 

\begin{defi} \label{defi-chtoucas}
Soit   $I$  un ensemble fini, $k\in \N$ et  $(I_{1}, ..., I_{k})$ une partition de $I$. Soit  $N$ un  sous-schéma fini de $X$. On note
  $$\Cht_{N,I}^{(I_{1},...,I_{k})} \text{({\it resp.} } 
  \Cht_{N,I,\underline{\omega}}^{(I_{1},...,I_{k})}, \  \Cht_{N,I,\underline{\omega}}^{(I_{1},...,I_{k}),\leq\mu})$$ l'ind-champ ({\it resp.} le champ) sur $(X\sm N)^{I}$ 
   qui classifie  les mêmes données   i)--iii) (de la  \defiref{defi-Hecke}) que  $$\Hecke_{N,I}^{(I_{1},...,I_{k})}  \text{({\it resp.} } 
\Hecke_{N,I,\underline{\omega}}^{(I_{1},...,I_{k})}, \ \Hecke_{N,I,\underline{\omega}}^{(I_{1},...,I_{k}),\leq\mu})$$
 plus un  isomorphisme $ \sigma: \ta \mc{G}_{0}\isom\mc{G}_{k}$, préservant les structures de niveau, c'est-à-dire   vérifiant 
 $\psi_{k}\circ \restr{\sigma}{N\times S}=\ta \psi_{0}$.  On note 
\begin{gather}\label{defi-gamma}
\gamma_{N,(I)}^{(I_{1},...,I_{k})}: \Cht_{N,I}^{(I_{1},...,I_{k})} \to \Hecke_{N,I}^{(I_{1},...,I_{k})}
\end{gather}
le morphisme tautologique qui consiste à oublier  $\sigma$. 
\end{defi}

Autrement dit  $\Cht_{N,I,\underline{\omega}} ^{(I_{1},...,I_{k})}$ est le produit fibré sur $\Bun_{G,N}\times \Bun_{G,N}$  de la  correspondance de Hecke  $\Hecke_{N,I,\underline{\omega}}^{(I_{1},...,I_{k})}$ avec le  graphe du morphisme  de Frobenius de $\Bun_{G,N}$. 

Pour récapituler, 
   $\Cht_{N,I} ^{(I_{1},...,I_{k})}$ est tel que pour tout schéma  $S$ sur  $\Fq$, $\Cht_{N,I}^{(I_{1},...,I_{k})}(S)$ classifie la donnée de    
\begin{gather}\label{donnee-chtouca}\big( (x_i)_{i\in I}, (\mc G_{0}, \psi_{0}) \xrightarrow{\phi_{1}}  (\mc G_{1}, \psi_{1}) \xrightarrow{\phi_{2}}
\cdots\xrightarrow{\phi_{k-1}}  (\mc G_{k-1}, \psi_{k-1}) \xrightarrow{ \phi_{k}}    (\ta{\mc G_{0}}, \ta \psi_{0})
\big)
\end{gather}
tels que   
 
 i)   $x_i\in (X\sm N)(S)$ pour $i\in I$, 

ii) $(\mc G_{0}, \psi_{0}),..., (\mc G_{k-1}, \psi_{k-1})\in \Bun_{G,N}(S)$, 
 
iii) en notant   $(\mc G_{k}, \psi_{k})=(\ta{\mc G_{0}}, \ta \psi_{0})$, 
pour  tout $j\in\{1,...,k\}$
 $$\phi_{j}:\restr{\mc G_{j-1}}{(X\times S)\sm(\bigcup_{i\in I_{j}}\Gamma_{x_i})}\isom \restr{\mc G_{j}}{(X\times S)\sm(\bigcup_{i\in I_{j}}\Gamma_{x_i})}$$ est un  isomorphisme tel que 
$\psi_{j}\circ \restr{\phi_{j}}{N\times S}=\psi_{j-1}$. 

De plus 
  $\Cht_{N,I,\underline{\omega}}^{(I_{1},...,I_{k})}$ est le    fermé  de $\Cht_{N,I}^{(I_{1},...,I_{k})}$ 
 défini par la condition que  
  $$\big( (x_i)_{i\in I}, (\mc G_{0}, \psi_{0}) \xrightarrow{\phi_{1}}  (\mc G_{1}, \psi_{1}) \xrightarrow{\phi_{2}}
\cdots\xrightarrow{\phi_{k-1}}  (\mc G_{k-1}, \psi_{k-1}) \xrightarrow{ \phi_{k}}    (\mc G_{k}, \psi_{k})
\big)
$$
appartient à  $\mr{Hecke}_{N,I,\underline{\omega}}^{(I_{1},...,I_{k})}$. 
  Enfin 
 $\Cht_{N,I,\underline{\omega}}^{(I_{1},...,I_{k}),\leq\mu}$ est l'ouvert de  $\Cht_{N,I,\underline{\omega}}^{(I_{1},...,I_{k})}$ défini par la  condition  $\mc G_{0}\in \Bun_{G}^{\leq\mu}$.

\begin{rem} \label{rem-chtoucas}
a) L'ind-champ   $\Cht_{N,I}^{(I_{1},...,I_{k})}$ est muni d'une action naturelle    du groupe $G(\mc{O}_N)$, qui fait agir $g\in G(\mc{O}_N)$ en rempla\c cant   $\psi_{j}$ par  $g\circ\psi_{j}$ 
pour tout   $j\in \{0,...,k\}$. Cette action préserve   $\Cht_{N,I,\underline{\omega}}^{(I_{1},...,I_{k})}$ et $\Cht_{N,I,\underline{\omega}}^{(I_{1},...,I_{k}),\leq\mu}$. 

b) Lorsque   $G=GL_r$,  $I=\{1,2\}$, $\omega_1=(1,0,\ldots,0)$ et $\omega _2=(0,\ldots,0,-1)$, les champs  
 $\Cht_{N,I,\underline{\omega}}^{(\{1\}, \{2\})}$, {\it resp.}  
  $\Cht_{N,I,\underline{\omega}}^{(\{2\}, \{1\})}$ sont les  champs de chtoucas  à gauche, {\it resp.} à  droite introduits par  Drinfeld \cite{Dr1}. Les troncatures ne sont pas les mêmes que dans  \cite{laurent-inventiones} mais les systèmes inductifs sont comparables. Dans ce cas  $x_{1}$ et $x_{2}$ sont appelés le  zéro et le  pôle du chtouca. 
\end{rem}

\begin{notation}
En général  les  $x_{i}$ seront appelés les  pattes du chtouca.  On notera    $$\mf p_{N,I} ^{(I_{1},...,I_{k})}: \Cht_{N,I} ^{(I_{1},...,I_{k})}\to (X\sm N)^{I}$$ le  morphisme correspondant (appelé morphisme caractéristique  par certains auteurs). On notera   encore $\mf p_{N,I} ^{(I_{1},...,I_{k})}$ sa restriction à 
$\Cht_{N,I, \underline{\omega}} ^{(I_{1},...,I_{k}) }$ et on notera 
 $ \mf p_{N,I} ^{(I_{1},...,I_{k}),\leq\mu}$ sa  restriction à l'ouvert  $\Cht_{N,I,\underline{\omega}} ^{(I_{1},...,I_{k}),\leq\mu}$. \end{notation}

\begin{defi} \label{defi-adm}
Un    $I$-uplet $\underline{\omega}=(\omega_i)_{i\in I}\in X_*(T)^I$ sera dit    {\em admissible}
si   $[\sum_{i\in I}\omega_i]=0$ dans   $\pi_{1}(G)$. 
\end{defi}

\begin{rem} \label{rem-adm}
Les   $I$-uplets de  copoids dominants de  $G$ sont en   bijection canonique avec les   $I$-uplets de  poids dominants de $\widehat{G}$, et donc avec les  $I$-uplets de   représentations irréductibles de   $\widehat{G}$.
Par cette  bijection, les    $I$-uplets admissibles correspondent aux  $I$-uplets de  représentations irréductibles de $\widehat{G}$ tels que  $Z_{\wh G}$ agisse trivialement sur leur produit tensoriel.
\end{rem}

La   proposition suivante, qui généralise les   propositions 2.3 et  3.2 de \cite{Dr1}, est contenue dans la   proposition 2.16 de \cite{var}. 

\begin{prop} (variante de la   proposition 2.16 de \cite{var}) \label{prop-chtoucas}
a) $\Cht_{N,I,\underline{\omega}}^{(I_{1},...,I_{k})}$ est un champ de Deligne-Mumford  sur  $(X\sm N)^I$,
et il  est localement de type fini. De plus les   composantes connexes de $\Cht_{N,I,\underline{\omega}}^{(I_{1},...,I_{k}),\leq\mu}$ sont  des  quotients  de schémas 
quasi-projectifs sur   $(X\sm N)^I$  par des groupes finis (si $N$ est non vide, sinon on doit se restreindre à $U^{I}$ avec $U \subsetneq X$ arbitraire), et sont des  schémas 
quasi-projectifs dès que   $\dim(\mc O_{N})$ est assez grand  en fonction de   $\mu$ et  
$\underline{\omega}$. 

b) $\Cht_{N,I,\underline{\omega}}^{(I_{1},...,I_{k})}$ ({\it resp.} $\Cht_{N,I,\underline{\omega}}^{(I_{1},...,I_{k}),\leq\mu}$) est un revêtement fini  étale  et galoisien de $\restr{\Cht_{I,\underline{\omega}}^{(I_{1},...,I_{k})}}{(X\sm N)^{I}}$ ({\it resp.} $\restr{\Cht_{I,\underline{\omega}}^{(I_{1},...,I_{k}),\leq\mu}}{(X\sm N)^{I}}$) de  groupe de Galois $G(\mc O_{N})$.

c) Si $J$ est une partie  de  $I$ telle que  $\omega_{i}=0$ pour $i\in I\sm J$, 
alors en notant  $J_{i}=I_{i}\cap J$ pour $i\in \{1,...,k\}$, 
le champ  $\Cht_{N,I,\underline{\omega}}^{(I_{1},...,I_{k})}$ est canoniquement isomorphe au produit   $(X\sm N)^{I\sm J}\times \Cht_{N,J,(\omega_{j})_{j\in J}}^{(J_{1},...,J_{k})}$,

d) Si  $\omega_i=0$ pour tout  $i\in I$, alors    $\Cht_{N,I,\underline{\omega}}^{(I_{1},...,I_{k})}$ est canoniquement isomorphe au produit de  $(X\sm N)^I$ avec le  champ discret  $\Bun_{G,N}(\Fq)$.

e) $\Cht_{N,I,\underline{\omega}}^{(I_{1},...,I_{k})}$ est non vide si et seulement si  $\underline{\omega}$ est admissible.
\end{prop}
\noindent{\bf Démonstration.} Tout est dans la   proposition 2.16 de  \cite{var}, à part    c) qui est évident. Grâce à    c), il suffit de montrer d)  dans le cas particulier  où  $I=\emptyset$, où il se réduit  à l'affirmation que 
 $\Cht_{N,\emptyset,0}^{(\emptyset)}$ est canoniquement isomorphe au  champ discret  $\Bun_{G,N}(\Fq)$. 
 \cqfd

\begin{construction}\label{construction-oubli-cht}
Comme dans   \eqref{oubli-Hecke}, soit 
$(I'_{1},...,I'_{k'})$ une  partition de  $I$ obtenue à partir de la  partition  $(I_{1},...,I_{k})$  en réunissant des parties dont les indices sont adjacents,  c'est-à-dire que l'on choisit des entiers  $0=j_{0}<j_{1}< \cdots < j_{k'}=k$
et que l'on pose   $I'_{j'}=\bigcup_{j_{j'-1}< j\leq j_{j'}}I_{j}$. 
La même   construction que dans   \eqref{oubli-Hecke} fournit un morphisme d'oubli 
 $$\pi^{(I_{1},...,I_{k})}_{(I'_{1},...,I'_{k'})}: 
\Cht_{N,I,\underline{\omega}}^{(I_{1},...,I_{k})} \to 
\Cht_{N,I,\underline{\omega}}^{(I'_{1},...,I'_{k'})} $$ 
qui conserve les   $(x_{i})_{i\in I}$,  les 
$(\mc G_{j_{j'}}, \psi_{j_{j'}})$ pour $j'\in \{0,...,k'\}$,  les  composées des morphismes $\phi_{j}$ qui les relient, et $ \sigma: \ta\mc{G}_{0}\isom\mc{G}_{k}$. 
\end{construction}

Soit  $\underline n=(n_{i})_{i\in I}\in \N^{I}$ comme dans la   \propref{action-Gnx-Grass}. 
On  pose   \begin{gather}\label{morph-epsilon-Cht-Gr}\epsilon_{(I),\underline{\omega},\underline n}^{(I_{1},...,I_{k})}
=
\delta_{(I),\underline{\omega},\underline n}^{(I_{1},...,I_{k})}\circ 
\gamma_{(I)}^{(I_{1},...,I_{k})}: \Cht_{I,\underline{\omega}}^{(I_{1},...,I_{k})}  \to \mr{Gr}_{I,\underline{\omega}}^{(I_{1},...,I_{k})}/G_{\sum n_{i}x_{i}}
\end{gather} 
(où l'on rappelle que  
$\delta_{(I),\underline{\omega},\underline n}^{(I_{1},...,I_{k})}$ et 
$\gamma_{(I)}^{(I_{1},...,I_{k})}$
ont été définis dans   \eqref{morph-delta-defi-Grass} et  \eqref{defi-gamma}). 
Autrement dit   
$\epsilon_{(I),\underline{\omega},\underline n}^{(I_{1},...,I_{k})}$ associe  à    un chtouca $$\big( (x_i)_{i\in I}, \mc G_{0} \xrightarrow{\phi_{1}}  \mc G_{1}\xrightarrow{\phi_{2}}
\cdots\xrightarrow{\phi_{k-1}}  \mc G_{k-1} \xrightarrow{ \phi_{k}}    \ta{\mc G_{0}}\big)$$ 
le   $G_{\sum n_{i}x_{i}}$-torseur égal à la   restriction de  $\mc G_{k}= \ta{\mc G_{0}}$ à  $\Gamma_{\sum n_{i}x_{i}}$ et le point de $\restr{\mc G_{k}}{\Gamma_{\sum n_{i}x_{i}}} \times_{G_{\sum n_{i}x_{i}}}\mr{Gr}_{I,\underline{\omega}}^{(I_{1},...,I_{k})}$  qui  est déterminé par  $(\mc G_{0} \xrightarrow{\phi_{1}}  \mc G_{1} \xrightarrow{\phi_{2}}
\cdots\xrightarrow{\phi_{k}}  \mc G_{k})$, grâce à la   construction 
\ref{rem-grassm},   à la \propref{action-Gnx-Grass} et à la \remref{rem-apres-action-Gnx-Grass}. 

\begin{prop}\label{lissite-Cht-Grass}
Le   morphisme $$\epsilon_{(I),\underline{\omega},\underline n}^{(I_{1},...,I_{k})}: \Cht_{I,\underline{\omega}}^{(I_{1},...,I_{k})}  \to \mr{Gr}_{I,\underline{\omega}}^{(I_{1},...,I_{k})}/G_{\sum n_{i}x_{i}}
$$ est lisse de  dimension $\dim G_{\sum n_{i}x_{i}}=(\sum n_{i})\dim G$. 
\end{prop}
\noindent {\bf Démonstration.} 
Comme l'énoncé  est local pour la  topologie lisse sur le but, il suffit de montrer,   pour tout schéma $S$ sur $X^{I}$ et tout point  $z$ de 
$\mr{Gr}_{I,\underline{\omega}}^{(I_{1},...,I_{k})}(S)$,  
la lissité sur  $S$ de 
\begin{gather}\label{Cht-prod-S}\Cht_{I,\underline{\omega}}^{(I_{1},...,I_{k})}\times _{ \mr{Gr}_{I,\underline{\omega}}^{(I_{1},...,I_{k})}/G_{\sum n_{i}x_{i}}}S.\end{gather}
Or ce champ est l'égalisateur  
\begin{itemize}
\item du morphisme d'oubli $$b_{1}: \Bun_{G,\sum n_{i}x_{i}}\times_{X^{I}}S \to \Bun_{G}\times S $$ qui est lisse de dimension $\dim G_{\sum n_{i}x_{i}}$,  
\item et du  morphisme $$b_{2 }: \big( \Frob_{\Bun_{G}}\times \Id_{S}\big)\circ  \big( a_{z},\Id_{S}): \Bun_{G,\sum n_{i}x_{i}}\times_{X^{I}}S \to \Bun_{G}\times S,$$ où 
$a_{z}$ est  la ``modification par $z$'', c'est-à-dire que 
$$a_{z}=p_{0}\circ \big(\beta_{I,\underline{\omega},\underline n}^{(I_{1},...,I_{k})}\big)^{-1}(z,\bullet): \Bun_{G,\sum n_{i}x_{i}}\times_{X^{I}}S\to \Bun_{G}$$
(où l'on rappelle que  $p_{0}$ et   $\beta_{I,\underline{\omega},\underline n}^{(I_{1},...,I_{k})}$ ont été définis dans la   notation \ref{N:proj} et dans  \eqref{defi-beta-avant-Grass-lissite}), ou en d'autres termes, en notant  
 $$z=\big((x_{i})_{i\in I}, \mc G_{0} \xrightarrow{\phi_{1}}  
\mc G_{1}\xrightarrow{\phi_{2}}
\cdots\xrightarrow{\phi_{k-1}}  \mc G_{k-1} \xrightarrow{ \phi_{k}}   \mc G_{k}\isor{\theta} \restr{G}{\Gamma_{\sum \infty x_i}}\big),$$ 
pour tout  $S$-schéma $S'$ 
et $(\mc F,\psi)\in \Bun_{G,\sum n_{i}x_{i}}(S')$, si $\wt\psi$ est une trivialisation de $\restr{\mc F}{\Gamma_{\sum \infty x_i}}$ qui raffine $\psi$, 
on  étend  les $\mc G_{i}$ (ou plutôt leurs images inverses par  
$ S'\to   S$) en des  $G$-torseurs $\check{\mc G}_{i}$ sur $X\times S'$  en recollant  $\mc G_{k}$ avec  $\mc F$ à l'aide  de $\wt\psi^{-1}\circ \theta$ et alors $a_{z}(\mc F,\psi)$ est égal à $\check{\mc G}_{0}$. 
 \end{itemize}
 Donc \eqref{Cht-prod-S} est l'égalisateur de morphismes $b_{1}$ et $b_{2}$
 de  $\Bun_{G,\sum n_{i}x_{i}}\times_{X^{I}}S$ vers $ \Bun_{G}\times S $, et
 \begin{itemize}
 \item ces deux champs sont lisses au-dessus de $S$, 
 \item $b_{1}$ et $b_{2}$ commutent avec la projection vers $S$, 
 \item $b_{1}$ est lisse et $b_{2}$ a une dérivée nulle dans les fibres au-dessus de $S$ (parce que   $\Frob_{\Bun_{G}}$ a une dérivée nulle). 
 \end{itemize}
  On en déduit que cet égalisateur est lisse sur $S$. 
\cqfd

 Plus tard on aura  besoin d'une variante de la  proposition préccédente,  où l'on casse  
$ \mr{Gr}_{I,\underline{\omega}}^{(I_{1},...,I_{k})}/G_{\sum n_{i}x_{i}}$ en morceaux. 

On note  $$\gamma_{(I_{1},...,I_{k})}^{(I_{1},...,I_{k})}: 
 \Cht_{N,I,\underline{\omega}} ^{(I_{1},...,I_{k})}
\to \prod_{j=1}^{k}
\Hecke_{I_{j},(\omega_{i})_{i\in I_{j}}}^{(I_{j})}$$
le   morphisme qui envoie un  chtouca 
$$\big( (x_i)_{i\in I}, \mc G_{0} \xrightarrow{\phi_{1}}  \mc G_{1}\xrightarrow{\phi_{2}}
\cdots\xrightarrow{\phi_{k-1}}  \mc G_{k-1} \xrightarrow{ \phi_{k}}    \ta{\mc G_{0}}\big)$$  sur la  famille 
$$\Big( (x_i)_{i\in I_{j}}, \mc G_{j-1} \xrightarrow{\phi_{j}}  \mc G_{j} \Big)_{j\in \{1,...,k\}}$$
en notant  $\mc G_{k}=\ta \mc G_{0}$. 

Par  composition avec le  produit des  morphismes 
$$\delta_{(I_{j}),(\omega_{i})_{i\in I_{j}},(n_{i})_{i\in I_{j}}}^{(I_{j})}: 
\Hecke_{I_{j},(\omega_{i})_{i\in I_{j}}}^{(I_{j})} 
\to \mr{Gr}_{I_{j},(\omega_{i})_{i\in I_{j}}}^{(I_{j})}/G_{\sum _{i\in I_{j}}n_{i}x_{i}}$$
on obtient 
\begin{gather}\label{epsilon-def-1k-1k}\epsilon_{(I_{1},...,I_{k}),\underline{\omega},\underline{n}}^{(I_{1},...,I_{k})}: 
\Cht_{N,I,\underline{\omega}} ^{(I_{1},...,I_{k})}
\to \prod_{j=1}^{k}
\Big(\mr{Gr}_{I_{j},(\omega_{i})_{i\in I_{j}}}^{(I_{j})}/G_{\sum _{i\in I_{j}}n_{i}x_{i}}\Big).\end{gather}

 La proposition \ref{lissite-Cht-Grass} admet alors la variante suivante.

\begin{prop}\label{lissite-Cht-Grass-gen}
Le morphisme $\epsilon_{(I_{1},...,I_{k}),\underline{\omega},\underline{n}}^{(I_{1},...,I_{k})}$ est lisse. 
\end{prop}
\noindent{\bf Démonstration. } On a établi dans le  \lemref{kapp-eps-J-lemme-smooth} la lissité du  morphisme 
$$\wt\kappa_{I,(\omega_{i})_{i\in I}}^{(I_{1},...,I_{k})}: 
\mr{Gr}_{I,(\omega_{i})_{i\in I}}^{(I_{1},...,I_{k})}/G_{\sum _{i\in I}m_{i}x_{i}}\to 
\prod_{j=1}^{k}
\Big(\mr{Gr}_{I_{j},(\omega_{i})_{i\in I_{j}}}^{(I_{j})}/G_{\sum _{i\in I_{j}}n_{i}x_{i}}\Big)$$ qui était défini  lorsque les  entiers  $m_{i}$ etaient  assez grands.

Comme  $\epsilon_{(I_{1},...,I_{k}),\underline{\omega},\underline{n}}^{(I_{1},...,I_{k})}=\wt\kappa_{I,(\omega_{i})_{i\in I}}^{(I_{1},...,I_{k})} \circ  \epsilon_{(I),\underline{\omega},\underline m}^{(I_{1},...,I_{k})}$, et  comme 
$ \epsilon_{(I),\underline{\omega},\underline m}^{(I_{1},...,I_{k})}$ est lisse par la \propref{lissite-Cht-Grass}, 
le  morphisme  $\epsilon_{(I_{1},...,I_{k}),\underline{\omega},\underline{n}}^{(I_{1},...,I_{k})}$ est lisse. 
\cqfd

\begin{rem} \label{rem-Jkl} Les  morphismes $\epsilon_{(I),\underline{\omega},\underline n}^{(I_{1},...,I_{k})}$ et 
$\epsilon_{(I_{1},...,I_{k}),\underline{\omega},\underline{n}}^{(I_{1},...,I_{k})}$ des propositions \ref{lissite-Cht-Grass} et \ref{lissite-Cht-Grass-gen} sont des cas particuliers d'un morphisme plus général, qui ne sera jamais utilisé. 
Soit $(I_{1},...,I_{k})$  une partition  de $I$ et pour tout $j\in \{1,...,k\}$ soit  $m_{j}\in \N^{*}$ et  
$(J_{j,1}, ..., J_{j,m_{j}})$ une   partition de  $I_{j}$. On note   $(J_{j,l})_{j,l}$ la   partition  de  $I$ égale à  $(J_{1,1}, ..., J_{1,m_{1}}, J_{2,1}, ..., J_{2,m_{2}}, ..., J_{k,1}, ..., J_{k,m_{k}})$. 
Lorsque les entiers $n_{i}$ sont assez grands on a    un   morphisme naturel 
$$\epsilon_{(I_{1},...,I_{k}),\underline{\omega},\underline{n}}^{(J_{j,l})_{j,l}}
: \Cht_{I,\underline{\omega}}^{(J_{j,l})_{j,l}}
\to 
\prod_{j=1}^{k}
\Big(\mr{Gr}_{I_{j},(\omega_{i})_{i\in I_{j}}}^{(J_{j,1}, ..., J_{j,m_{j}})}/G_{\sum _{i\in I_{j}}n_{i}x_{i}}\Big)
$$
et il est lisse. 
\end{rem}

Les propositions précédentes généralisent  la preuve donnée par   Drinfeld de la lissité du  champ des chtoucas dans le cas   b) de la   \remref{rem-chtoucas}. Par ailleurs elles impliquent que   localement pour la  topologie lisse, les champs  de $G$-chtoucas sont isomorphes à des  strates fermées des  grassmanniennes affines de Beilinson-Drinfeld 
(ou à des produits de telles strates fermées).

 La 
proposition  suivante  montre que cela est vrai localement pour la topologie étale
et en fournit une preuve constructive.  C'est une variante du 
théorème 2.20   de  \cite{var}. 
La preuve que nous allons donner est proche de la preuve figurant dans le paragraphe  4 de \cite{var}, mais mais elle n'utilise pas  le théorème de  Drinfeld-Simpson (le théorème  2 de  \cite{DS}, mentionné au début du chapitre précédent). La proposition suivante  peut être sautée par le lecteur car elle ne sera pas utilisée  mais elle est intéressante par elle-même.    
 
\begin{prop} (variante du théorème 2.20   de  \cite{var}) \label{thm-locisom} 
Localement  pour la topologie étale, $\Cht_{N,I,\underline{\omega}}^{(I_{1},...,I_{k})}$ est isomorphe à 
$\prod_{j=1}^{k}\mr{Gr}_{I_{j},(\omega_i)_{i\in I_{j}}}^{(I_{j})}$, de fa\c con compatible avec le  morphisme $\epsilon_{(I_{1},...,I_{k}),\underline{\omega},\underline{n}}^{(I_{1},...,I_{k})}$ de  \eqref{epsilon-def-1k-1k}.

Plus précisément, on suppose que les  entiers $(n_{i})_{i\in I}$ sont assez grands. Pour tout $j\in\{1,...,k\}$, soit 
\begin{gather}\label{alpha-x_i}
\alpha_{j}\times (x_{i})_{i\in I_{j}} : 
U_{j}\to\Bun_{G,N}\times (X\sm N)^{I_{j}}\end{gather}
 un  morphisme étale tel que, si $\mc G$ désigne le $G$-torseur universel sur $X\times \Bun_{G, N}$, on a une trivialisation 
$$\theta_{j}: \restr{(\Id_{X} \times \alpha_{j} )^{*}(\mc G)}{\Gamma_{\sum_{i\in I_{j}}n_{i}x_{i}}}\isom \restr{G}{\Gamma_{\sum_{i\in I_{j}}n_{i}x_{i}}}$$
sur le  sous-schéma fermé  
$\Gamma_{\sum_{i\in I_{j}}n_{i}x_{i}}\subset X \times U_{j} $. 

Soit $$\mc U=\Cht_{N,I,\underline{\omega}}^{(I_{1},...,I_{k})}\times_{\prod_{j=1}^{k}\big(\Bun_{G,N}\times (X\sm N)^{I_{j}}\big)}\prod_{j=1}^{k} U_{j}$$
obtenu par image inverse  du produit des  morphismes étales \eqref{alpha-x_i}
par le  morphisme 
$$\Cht_{N,I,\underline{\omega}}^{(I_{1},...,I_{k})}\to \prod_{j=1}^{k}\Big(\Bun_{G,N}\times (X\sm N)^{I_{j}}\Big)$$ qui envoie  
\begin{gather}\label{élément-Cht-proof-locisom}\big( (x_i)_{i\in I}, (\mc G_{0}, \psi_{0}) \xrightarrow{\phi_{1}}  (\mc G_{1}, \psi_{1}) \xrightarrow{\phi_{2}}
\cdots\xrightarrow{\phi_{k-1}}  (\mc G_{k-1}, \psi_{k-1}) \xrightarrow{ \phi_{k}}     (\mc G_{k}, \psi_{k})\simeq (\ta{\mc G_{0}}, \ta \psi_{0})
\big)
\end{gather} sur la  famille 
\begin{gather}
\label{fam-G-psi}
\Big( (\mc G_{j}, \psi_{j}), (x_{i})_{i\in I_{j}}  \Big)_{j\in \{1,...,k\}}.\end{gather}
Evidemment la première projection  $\mc U\to \Cht_{N,I,\underline{\omega}}^{(I_{1},...,I_{k})}$ est étale. 
Soit $$
\beta = \prod_{j=1}^{k} \beta_{j}:\mc U\to \prod_{j=1}^{k}\mr{Gr}_{I_{j},(\omega_i)_{i\in I_{j}}}^{(I_{j})}$$ le  morphisme tel que pour tout $j$  
et pour tout point $u$ dans $\mc U$ (d'images \eqref{élément-Cht-proof-locisom} dans $\Cht_{N,I,\underline{\omega}}^{(I_{1},...,I_{k})}$ et $(u_{j})_{j\in  \{1,...,k\}}$ dans  $\prod_{j=1}^{k} U_{j}$), 
$\beta_{j}(u)$ est le point de  $\mr{Gr}_{I_{j},(\omega_i)_{i\in I_{j}}}^{(I_{j})}$ associé à $ (\mc G_{j-1} \xrightarrow{\phi_{j}} \mc G_{j})$ et à la trivialisation  $\theta_{j}$ (au point  $u_{j}\in U_{j}$) grâce à la  \remref{rem-apres-action-Gnx-Grass}. Alors $\beta$ est étale. 
\end{prop}
\begin{rem} Ce théorème admet la généralisation suivante.  
Soit $(J_{j,l})_{j,l}$ comme dans la  \remref{rem-Jkl}. Alors,  
localement pour la topologie étale, $\Cht_{N,I,\underline{\omega}}^{(J_{j,l})_{j,l}}$ est isomorphe à 
$\prod_{j=1}^{k}\mr{Gr}_{I_{j},(\omega_i)_{i\in I_{j}}}^{(J_{j,1}, ..., J_{j,m_{j}})}$, de fa\c con  compatible avec le  morphisme $\epsilon_{(I_{1},...,I_{k}),\underline{\omega},\underline{n}}^{(J_{j,l})_{j,l}}$ de cette remarque. 
\end{rem}
\noindent {\bf Démonstration} 
(variante du paragraphe 4 de  \cite{var}).  On fixe  $\mu$ et on va démontrer la  proposition pour l'ouvert
  $\Cht_{N,I,\underline{\omega}}^{(I_{1},...,I_{k}), \leq \mu}$.
  Pour tout niveau  $N'\supset N$, 
 $\Cht_{N',I,\underline{\omega}}^{(I_{1},...,I_{k}), \leq \mu}$ est un revêtement étale galoisien  de  $\restr{\Cht_{N,I,\underline{\omega}}^{(I_{1},...,I_{k}), \leq \mu}}{(X\sm N')^{I}}$. Il suffit donc de montrer l'énoncé pour $N$ assez grand. 
  On suppose  $N$ assez grand  pour que  $\Bun_{G,N}^{\leq \mu+\kappa}$ soit un  schéma, où $\kappa$ est tel que pour tout point 
  \begin{gather}\label{élément-Hecke-proof-locisom}\big( (x_i)_{i\in I}, (\mc G_{0}, \psi_{0}) \xrightarrow{\phi_{1}}  (\mc G_{1}, \psi_{1}) \xrightarrow{\phi_{2}}
\cdots\xrightarrow{\phi_{k-1}}  (\mc G_{k-1}, \psi_{k-1}) \xrightarrow{ \phi_{k}}     (\mc G_{k}, \psi_{k})
\big)
\end{gather} 
  dans 
  $\Hecke_{N,I,\underline{\omega}}^{(I_{1},...,I_{k}), \leq \mu}$ tous les  $G$-torseurs $\mc G_{j}$ appartiennent à $\Bun_{G,N}^{\leq \mu+\kappa}$. La raison pour laquelle on augmente ainsi le niveau  est que le  \lemref{variante-lem-var-etale} ci-dessous est énoncé avec des schémas.

  Soit $$\mc W=\Hecke_{N,I,\underline{\omega}}^{(I_{1},...,I_{k}), \leq \mu}\times_{\prod_{j=1}^{k}\big(\Bun_{G,N}\times (X\sm N)^{I_{j}}\big)}\prod_{j=1}^{k} U_{j}$$
obtenu par image inverse du  produit des morphismes étales  \eqref{alpha-x_i}
par le  morphisme 
$$\Hecke_{N,I,\underline{\omega}}^{(I_{1},...,I_{k}), \leq \mu}\to \prod_{j=1}^{k}\Big(\Bun_{G,N}\times (X\sm N)^{I_{j}}\Big)$$ qui envoie \eqref{élément-Hecke-proof-locisom} 
sur la  famille \eqref{fam-G-psi}. 
 Evidemment la première  projection $\mc W\to \Hecke_{N,I,\underline{\omega}}^{(I_{1},...,I_{k}), \leq \mu}$ est étale. 
Soit $$
\beta = \prod_{j=1}^{k} \beta_{j}:\mc W\to \prod_{j=1}^{k}\mr{Gr}_{I_{j},(\omega_i)_{i\in I_{j}}}^{(I_{j})}$$   tel que pour tout $j$  
et pour tout point $w$ dans $\mc W$ (d'images \eqref{élément-Hecke-proof-locisom} dans  $\Hecke_{N,I,\underline{\omega}}^{(I_{1},...,I_{k}), \leq \mu}$ et $(w_{j})_{j\in  \{1,...,k\}}$ dans $\prod_{j=1}^{k} U_{j}$), 
$\beta_{j}(w)$ soit le  point de  $\mr{Gr}_{I_{j},(\omega_i)_{i\in I_{j}}}^{(I_{j})}$ associé à $ (\mc G_{j-1} \xrightarrow{\phi_{j}} \mc G_{j})$ et à la  trivialisation $\theta_{j}$ (au point $w_{j}\in U_{j}$).
On note $$p_{0}, p_{k}: \Hecke_{N,I,\underline{\omega}}^{(I_{1},...,I_{k}), \leq \mu}\to \Bun_{G,N}^{\leq \mu+\kappa}$$ le  morphisme qui envoie $w$ comme ci-dessus vers  $(\mc G_{0}, \psi_{0})$, {\it resp.} $(\mc G_{k}, \psi_{k})$. 
 Grâce à  l'isomorphisme 
\eqref{defi-beta-avant-Grass-lissite} on montre  (par récurrence sur $k$) que 
 $$(p_{k},\beta):\mc W\to \Bun_{G,N}^{\leq \mu+\kappa}\times 
\prod_{j=1}^{k}\mr{Gr}_{I_{j},(\omega_i)_{i\in I_{j}}}^{(I_{j})}$$ est étale. 
  On a 
  $$\mc U=\{w\in \mc W, 
  p_{k}(w)=\Frob_{\Bun_{G,N}^{\leq \mu+\kappa}}( p_{0}(w))\}.$$
    On applique le lemme suivant à 
  $$W=\mc W, \ Z=\Bun_{G,N}^{\leq \mu+\kappa}, \  
  T=\prod_{j=1}^{k}\mr{Gr}_{I_{j},(\omega_i)_{i\in I_{j}}}^{(I_{j})}, \  
  h=(p_{k},\beta)\text{ \  et  \ } f=p_{0}.$$
   Cela termine la preuve de la  proposition en montrant que $\beta$ est étale. \cqfd

\begin{lem}\label{variante-lem-var-etale} (petite généralisation  du lemme  4.3 de  \cite{var})
Soient  $W,Z,T$ des schémas de type fini sur  $\Fq$.  Soit  $h=(h_{1},h_{2}): W\to Z \times T$   un morphisme étale  et  $f: W\to Z$   un morphisme arbitraire. On pose 
$$\mc U=\{w\in W, h_{1}(w)=\Frob_{Z} (f(w))\}.$$
On suppose que  $Z$ est lisse sur  $\Fq$. Alors $h_{2}:\mc U\to T$ est étale. 
\end{lem}
\noindent {\bf Démonstration } (d'après la preuve du lemme  4.3 de  \cite{var}). 
Comme l'énoncé  est local pour la  topologie de Zariski de  $\mc U$
et donc pour celle de  $Z$, on peut supposer qu'il existe un  morphisme étale $\phi:Z\to \mb A^{m}$. 
Comme  $\mc U$ est ouvert et fermé dans 
$$\mc U'=\{w\in W, \phi\circ h_{1}(w)=\Frob_{\mb A^{m}} (\phi\circ  f(w))\},$$ quitte à  remplacer  $Z,h_{1},f$ par $\mb A^{m},\phi\circ h_{1}, \phi\circ f $, on est  ramené au cas où $Z=\mb A^{m}$. 

Quitte à remplacer $T$ et  $W$ par des ouverts de Zariski, on suppose 
 que $T$ est affine et que  $W$ est un ouvert de 
 $\on{Spec}\big(\mc  O(T) [x_{1},...,x_{m},y_{1},...,y_{n}]/(k_{1},...,k_{n})\big)$
sur lequel  $\det(\frac{\partial k_{i}}{\partial y_{j}})$ est  inversible, et que le morphisme étale $h:W\to Z\times T$ est donné par le morphisme évident 
$$ \mc  O(T) [x_{1},...,x_{m}]\to \mc  O(T) [x_{1},...,x_{m},y_{1},...,y_{n}]/(k_{1},...,k_{n})
.$$ On a noté  $\mc  O(T)$   l'anneau des  fonctions régulières sur  $T$, de sorte que 
$Z\times T=\on{Spec}\big(\mc  O(T) [x_{1},...,x_{m}]\big)$. Alors  
$\mc U\subset W$ est défini  par les  équations $g_{1},...,g_{m}$, où $g_{i}=x_{i}-f_{i}^{q}$ et les fonctions   $f_{i}$  sur $W$ sont telles que $f=(f_{1}, ..., f_{m})$ de $W$ dans $Z=\mb A^{m}$. Par conséquent   $\mc U$ est un ouvert de $\on{Spec}\big(\mc  O(T) [x_{1},...,x_{m},y_{1},...,y_{n}]/(k_{1},...,k_{n},g_{1},...,g_{m})\big)$ sur lequel la matrice 
$$\begin{pmatrix}
(\frac{\partial k_{i}}{\partial y_{j}}) & (\frac{\partial g_{k}}{\partial y_{j}}) \\
(\frac{\partial k_{i}}{\partial x_{l}}) & (\frac{\partial g_{k}}{\partial x_{l}})
\end{pmatrix}=\begin{pmatrix}
(\frac{\partial k_{i}}{\partial y_{j}}) & 0\\
(\frac{\partial k_{i}}{\partial x_{l}}) & \Id_{m}
\end{pmatrix}$$ est invertible, et  $h_{2}:\mc U\to T$ est  étale. 
\cqfd

On a une action évidente de $Z(F)\backslash Z(\mb A)$
sur $\Cht_{N,I,\underline{\omega}}^{(I_{1},...,I_{k})}$. En fait 
$Z(F)\backslash Z(\mb A)$ s'envoie naturellement dans $\Bun_{Z,N}(\Fq)$ et 
$\Bun_{Z,N}(\Fq)$ agit sur $\Cht_{N,I,\underline{\omega}}^{(I_{1},...,I_{k})}$
  par torsion d'un $G$-torseur par un $Z$-torseur. 
  Ces actions préservent les ouverts $\Cht_{N,I,\underline{\omega}}^{(I_{1},...,I_{k}), \leq \mu}$ 
  (car la troncature de Harder-Narasimhan se fait sur le $G^{\mr{ad}}$-torseur associé à $\mc G_{0}$)
  et elles commutent aux morphismes de Frobenius partiels.

Dans toute la suite on fixe un réseau  
 $\Xi$  dans $Z(F)\backslash Z(\mb A)$. Les quotients 
$\Cht_{N,I,\underline{\omega}}^{(I_{1},...,I_{k}), \leq \mu}/\Xi$ sont des champs de Deligne-Mumford de type fini.

Soit  $E$  une extension finie de  $\Ql$ contenant une racine carrée de  $q$. 

\begin{defi}\label{defi-F-NIomega}
On note $\mc F_{N,I,\underline{\omega},\Xi ,E}^{(I_{1},...,I_{k})}$  le faisceau d'intersection de  $\Cht_{N,I,\underline{\omega}}^{(I_{1},...,I_{k})}/\Xi$  normalisé relativement à   $(X\sm N)^{I}$, et à coefficients dans   $E$. 
\end{defi}

Dans les deux corollaires suivants nous allons donner des définitions équivalentes de   $\mc F_{N,I,\underline{\omega}, \Xi,E}^{(I_{1},...,I_{k})}$, qui sont mieux adaptées pour étudier la   coalescence des pattes et  l'action des  morphismes de  Frobenius partiels.

Par l'équivalence de Satake géométrique \cite{mv,ga-de-jong}, que nous avons rappelée dans le  \thmref{thm-geom-satake}, nous possédons le  $E$-faisceau pervers (à un décalage près) $\mc S_{I,\boxtimes_{i\in I} V_{\omega_{i}}}^{(I_{1},...,I_{k}),E}$ sur 
$\mr{Gr}_{I,\underline{\omega}}^{(I_{1},...,I_{k})}/G_{\sum _{i\in I}n_{i}x_{i}}$ dont l'image inverse à 
$\mr{Gr}_{I,\underline{\omega}}^{(I_{1},...,I_{k})}$ est le faisceau d'intersection   (avec la  normalisation perverse  relative à  $X^{I}$). D'après   la \remref{rem-Satake-Gad-texte}, $\mc S_{I,\boxtimes_{i\in I} V_{\omega_{i}}}^{(I_{1},...,I_{k}),E}$ est en fait un faisceau pervers (à un décalage près) sur 
$\mr{Gr}_{I,\underline{\omega}}^{(I_{1},...,I_{k})}/G^{\mr{ad}}_{\sum _{i\in I}n_{i}x_{i}}$, de fa\c con canonique. Ci-dessous on considère donc $\mc S_{I,\boxtimes_{i\in I} V_{\omega_{i}}}^{(I_{1},...,I_{k}),E}$ comme un faisceau pervers  (à un décalage près) sur 
$\mr{Gr}_{I,\underline{\omega}}^{(I_{1},...,I_{k})}/G^{\mr{ad}}_{\sum _{i\in I}n_{i}x_{i}}$.

Le morphisme $\epsilon_{N,(I),\underline{\omega},\underline n}^{(I_{1},...,I_{k})}$  
de  la proposition \ref{lissite-Cht-Grass}
ne se factorise pas par le quotient par $\Xi$ (comme me l'a fait remarquer un rapporteur anonyme), mais c'est  le cas de sa composée avec le morphisme d'oubli $ \mr{Gr}_{I,\underline{\omega}}^{(I_{1},...,I_{k})}/
  G_{\sum n_{i} x_i}\to  \mr{Gr}_{I,\underline{\omega}}^{(I_{1},...,I_{k})}/
  G^{\mr{ad}}_{\sum n_{i} x_i}$. 
  En effet l'action de $\Xi$ consiste à tordre par des $Z$-torseurs, et $G^{\mr{ad}}=G/Z$. 
  
  Autrement dit on possède un morphisme 
  \begin{gather}\label{lisse-chtouca-grass-texte-ad}
  \epsilon_{(I),\underline{\omega},\underline n}^{(I_{1},...,I_{k}), \Xi}:   \Cht_{N,I,\underline{\omega}}^{(I_{1},...,I_{k})}/\Xi\to \mr{Gr}_{I,\underline{\omega}}^{(I_{1},...,I_{k})}/
  G^{\mr{ad}}_{\sum n_{i} x_i}.  \end{gather}  
En le composant 
  avec le   morphisme  étale d'oubli du niveau 
$$\Cht_{N,I,\underline{\omega}}^{(I_{1},...,I_{k})} /\Xi \to \Cht_{I,\underline{\omega}}^{(I_{1},...,I_{k})}  /\Xi $$ 
on obtient donc     $$\epsilon_{N,(I),\underline{\omega},\underline n}^{(I_{1},...,I_{k}), \Xi}:\Cht_{N,I,\underline{\omega}}^{(I_{1},...,I_{k})} /\Xi   \to 
\mr{Gr}_{I,\underline{\omega}}^{(I_{1},...,I_{k})}/G^{\mr{ad}}_{\sum n_{i}x_{i}}
$$  qui est lisse  de  dimension $\dim G^{\mr{ad}}_{\sum n_{i}x_{i}}=(\sum n_{i})\dim G^{\mr{ad}}$.

 \begin{cor}\label{cor-IC-produit} (variante du  corollaire 2.21 c) de  \cite{var})
On a un  isomorphisme canonique 
\begin{gather} \label{isom-gamma}
\lambda_{N,I}^{(I_{1},...,I_{k})}: 
\mc F_{N,I,\underline{\omega}, \Xi,E}^{(I_{1},...,I_{k})}
\simeq 
\big(\epsilon_{N,(I),\underline{\omega},\underline n}^{(I_{1},...,I_{k}), \Xi}\big)^{*}\Big(\mc S_{I,\boxtimes_{i\in I} V_{\omega_{i}},E}^{(I_{1},...,I_{k})}\Big) \end{gather} 
  \end{cor}
 \cqfd

On utilisera aussi la variante suivante (qui résulte de  la \propref{lissite-Cht-Grass-gen}, ou   du \corref{cor-IC-produit} et de c) du \thmref{thm-geom-satake}, et où $\epsilon_{N,(I_{1},...,I_{k}),\underline{\omega},\underline n}^{(I_{1},...,I_{k}), \Xi}$ est défini de manière analogue à $\epsilon_{N,(I),\underline{\omega},\underline n}^{(I_{1},...,I_{k}), \Xi}$ ci-dessus). 

\begin{cor}\label{cor-IC-produit-gen} On a un  isomorphisme canonique  
\begin{gather} \label{isom-gamma2}
\lambda_{N,(I_{1},...,I_{k})}^{(I_{1},...,I_{k})}: 
\mc F_{N,I,\underline{\omega}, \Xi,E}^{(I_{1},...,I_{k})}
\simeq 
\big(\epsilon_{N,(I_{1},...,I_{k}),\underline{\omega},\underline n}^{(I_{1},...,I_{k}), \Xi}\big)^{*}\Big(\boxtimes _{j\in \{1,...,k\}}\mc S_{I_{j},\boxtimes_{i\in I_{j}} V_{\omega_{i}},E}^{(I_{j})}\Big) \end{gather} 
  \end{cor}
 \cqfd

\begin{rem} Les deux corollaires précédents sont des cas particuliers de l'énoncé suivant, qui ne sera jamais utilisé.   
Soit $(J_{j,l})_{j,l}$ comme  dans la remarque \ref{rem-Jkl}. 
On a alors un  isomorphisme canonique  \begin{gather} \label{isom-gamma2-bis}
\lambda_{N,(I_{1},...,I_{k})}^{(J_{j,l})_{j,l}}: 
\mc F_{N,I,\underline{\omega}, \Xi,E}^{(J_{j,l})_{j,l}}
\simeq 
\big(\epsilon_{N,(I_{1},...,I_{k}),\underline{\omega},\underline n}^{(J_{j,l})_{j,l}, \Xi}\big)^{*}\Big(\boxtimes _{j\in \{1,...,k\}}\mc S_{I_{j},\boxtimes_{i\in I_{j}} V_{\omega_{i}},E}^{(J_{j,l})_{j,l}}\Big). \end{gather} 
\end{rem}

\begin{cor} \label{cor-strates} (variante du corollaire 2.21 de \cite{var})
 Dans les  notations de la   \constructionref{construction-oubli-cht}, le   morphisme d'oubli 
 $$\pi^{(I_{1},...,I_{k})}_{(I'_{1},...,I'_{k'})}: 
\Cht_{N,I,\underline{\omega}}^{(I_{1},...,I_{k})} /\Xi  \to 
\Cht_{N,I,\underline{\omega}}^{(I'_{1},...,I'_{k'})} /\Xi $$
est projectif, surjectif et petit. 
Par conséquent  
$$R\big(\pi^{(I_{1},...,I_{k})}_{(I'_{1},...,I'_{k'})}\big)_{!}\big(
\mc F_{N,I,\underline{\omega}, \Xi,E}^{(I_{1},...,I_{k})}
\big)=\mc F_{N,I,\underline{\omega}, \Xi,E}^{(I'_{1},...,I'_{k'})}.$$
 \end{cor} 
\noindent{\bf Démonstration. } Cela résulte de la  
\propref{lissite-Cht-Grass} car  $$\Cht_{N,I,\underline{\omega}}^{(I_{1},...,I_{k})}=\Cht_{N,I,\underline{\omega}}^{(I'_{1},...,I'_{k'})}
 \times_{\mr{Gr}_{I,\underline{\omega}}^{(I'_{1},...,I'_{k'})}/G_{\sum n_{i}x_{i}}}\mr{Gr}_{I,\underline{\omega}}^{(I_{1},...,I_{k})}/G_{\sum n_{i}x_{i}}
$$ et la même propriété est connue pour les  strates fermées des grassmanniennes affines (voir \cite{mv}). \cqfd

 \begin{rem} (variante du   corollaire 2.21 de \cite{var})
 L'ouvert   $ \Cht_{N,I,\underline{\omega}}^{0}$ de $\Cht_{N,I,\underline{\omega}}^{(I)}$ 
formé des    
$$\big((x_i)_{i\in I}, (\mc G_{0}, \psi_{0}), (\mc G_{1}, \psi_{1}), 
 \phi_{1}, \sigma \big)$$  tels que les  $x_{i}$ soient deux à deux disjoints et la position relative de  $\mc G_{0}$ par rapport à    $\mc G_{1}\simeq  \ta \mc G_{0}$ en $x_{i}$ soit exactement  $\omega_{i}$,  est lisse et dense dans   $\Cht_{N,I,\underline{\omega}}^{(I)}$. Il est non vide si et seulement si  $\underline{\omega}$ est admissible. De plus  pour toute  partition 
 $(I_{1},...,I_{k})$  le morphisme d'oubli 
 $$\pi^{(I_{1},...,I_{k})}_{(I)} : \Cht_{N,I,\underline{\omega}}^{(I_{1},...,I_{k})} \to 
\Cht_{N,I,\underline{\omega}}^{(I)}$$ est un  isomorphisme au-dessus de cet ouvert,  et l'image inverse  de cet ouvert est  dense dans 
 $\Cht_{N,I,\underline{\omega}}^{(I_{1},...,I_{k})} $. On a les mêmes résultats après quotient par $\Xi$. 
\end{rem}

\begin{construction} \label{constr-corresp-hecke}{\bf : Correspondances de Hecke. }
On rappelle que  $K_{N}=\on{Ker}(G(\mb O)\to G(\mc O_{N}))$ est le  sous-groupe ouvert compact de $G(\mb A)$ associé à $N$ et  on note $K_{N,v}$ sa composante  en la place $v$, de sorte que   $K_{N}=\prod_{v}K_{N,v}$. On a   $K_{N,v}=G(\mc O_{v})$ pour tout  $v\in |X|\sm |N|$. 
Soit  $\mf T\subset |X|$   un ensemble fini de  places
et  $g=(g_{v})_{v\in |X|}\in G(\mb A)$    tel que  
$g_{v}\in K_{N,v}$ pour tout  $v\not \in \mf T$. Alors 
$g$ induit une  correspondance   représentable  finie étale  $\Gamma_{N}(g)$ entre l'ouvert $$
\restr{\Cht_{N,I,\underline{\omega}} ^{(I_{1},...,I_{k})}}{(X\sm (|N|\cup \mf T))^{I}}:=
(\mf p_{N,I}^{(I_{1},...,I_{k})})^{-1}((X\sm (|N|\cup \mf T))^{I}) \cap
\Cht_{N,I,\underline{\omega}} ^{(I_{1},...,I_{k})}$$ et lui-même. 
Cette correspondance dépend seulement de  la  double classe $K_{N}gK_{N}$
et  si  $N'$ est un  niveau  tel que  $K_{N'}\subset K_{N}\cap  gK_{N}g^{-1}\cap g^{-1}K_{N}g$ et  $|N'|=|N| \cup \mf T$, elle est donnée  par le diagramme 
\begin{gather}\label{corresp-hecke}
\begin{CD}
\Cht_{N',I,\underline{\omega}} ^{(I_{1},...,I_{k})}/\big((K_{N}\cap  gK_{N}g^{-1})/K_{N'} \big)  @>g>\sim>
\Cht_{N',I,\underline{\omega}} ^{(I_{1},...,I_{k})}  /\big((K_{N}\cap  g^{-1}K_{N}g)/K_{N'} \big)\\
 @VV\on{pr}_{1}V @VV\on{pr}_{2}V \\
\restr{\Cht_{N,I,\underline{\omega}} ^{(I_{1},...,I_{k})} }{(X\sm (|N|\cup \mf T))^{I}}
& &
\restr{\Cht_{N,I,\underline{\omega}} ^{(I_{1},...,I_{k})} }{(X\sm (|N|\cup \mf T))^{I}} \end{CD}
 \end{gather}
 où la flèche du haut est l'action à droite de $g$, qui envoie la structure de niveau $\psi$ sur $g^{-1}\circ \psi$.

De plus  il existe  $\kappa$ (dépendant de  $g$) tel que, pour tout  $\mu$, la  correspondance 
$\Gamma_{N}(g)$ et  sa transposée envoient  
$\restr{\Cht_{N,I,\underline{\omega}} ^{(I_{1},...,I_{k}),\leq\mu}}{(X\sm (|N|\cup \mf T))^{I}}$ dans   $\restr{\Cht_{N,I,\underline{\omega}} ^{(I_{1},...,I_{k}),\leq\mu+\kappa}}{(X\sm (|N|\cup \mf T))^{I}}$. 

Dans le cas particulier où $|N| \cap \mf T=\emptyset$, la  correspondance $\Gamma_{N}(g)$ admet une description plus simple. En effet 
 $K_{N}g K_{N}$ est déterminé par la donnée d'une famille  $(\lambda_{\mf t})_{\mf t\in \mf T}$ de copoids dominants de   $G$. Alors  $\Gamma_{N}(g)$  est le champ dont les  $S$-points classifient la donnée de  $(x_i)_{i\in I}:S\to (X\sm (|N|\cup \mf T))^{I}$ 
et d'un diagramme commutatif  
$$ \xymatrix{
 (\mc G'_{0}, \psi'_{0})   \ar[r]^-{\phi'_{1}}  & 
 (\mc G'_{1}, \psi'_{1})   \ar[r]^-{\phi'_{2}}   & 
  \cdots \ar[r]^-{\phi'_{k-1}} & 
 (\mc G'_{k-1}, \psi'_{k-1}) \ar[r]^-{ \phi'_{k}}   & 
 (\ta{\mc G_{0}'}, \ta \psi_{0}')  \\
 (\mc G_{0}, \psi_{0})   \ar[r]^-{\phi_{1}}   \ar[u]_-{\kappa_{0}} & 
 (\mc G_{1}, \psi_{1}) \ar[r]^-{\phi_{2}} \ar[u]_-{\kappa_{1}}& 
 \cdots \ar[r]^-{\phi_{k-1}}  & 
 (\mc G_{k-1}, \psi_{k-1}) \ar[r]^-{ \phi_{k}}   \ar[u]_-{\kappa_{k-1}}  & 
 (\ta{\mc G_{0}}, \ta \psi_{0}) \ar[u]_-{\ta \kappa_{0}} 
  } $$
tel que la ligne inférieure  \begin{gather}\label{GammaN-pr1}\big( (x_i)_{i\in I}, (\mc G_{0}, \psi_{0}) \xrightarrow{\phi_{1}}  (\mc G_{1}, \psi_{1}) \xrightarrow{\phi_{2}}
\cdots\xrightarrow{\phi_{k-1}}  (\mc G_{k-1}, \psi_{k-1}) \xrightarrow{ \phi_{k}}    (\ta{\mc G_{0}}, \ta \psi_{0})
\big)\end{gather} et la ligne supérieure 
\begin{gather}\label{GammaN-pr2}
 \big( (x_i)_{i\in I}, (\mc G'_{0}, \psi'_{0}) \xrightarrow{\phi'_{1}}  (\mc G'_{1}, \psi'_{1}) \xrightarrow{\phi'_{2}}
\cdots\xrightarrow{\phi'_{k-1}}  (\mc G'_{k-1}, \psi'_{k-1}) \xrightarrow{ \phi'_{k}}    (\ta{\mc G'_{0}}, \ta \psi'_{0})\big)
\end{gather}
appartiennent à    $\Cht_{N,I,\underline{\omega}} ^{(I_{1},...,I_{k})} (S)$
 et que 
 \begin{itemize}
 \item $\kappa_{0}:\restr{\mc G_{0}}{(X\sm \mf T)\times S}\isom \restr{\mc G'_{0}}{(X\sm \mf T)\times S}$ soit  un isomorphisme 
 entrela\c cant  $\psi_{0}$ et $\psi'_{0}$, 
 \item en tout point   $\mf t\in \mf T$,  la position relative 
 de  $\mc G_{0}$ par rapport à   $\mc G'_{0}$ soit {\it  égale} au copoids dominant  $\lambda_{\mf t}$, plus précisément en  tout point géométrique $s$ de $S$, les  restrictions de $(\mc G_{0}\xrightarrow{\phi_{k}\cdots \phi_{1}} \ta\mc G_{0})$ et $(\mc G'_{0}\xrightarrow{\phi'_{k}\cdots \phi'_{1}} \ta\mc G'_{0})$ 
 au voisinage formel  de $\mf t$  dans $X$  
peuvent être trivialisées de fa\c con unique modulo l'action de   $G(\mc O_{\mf t})$ et alors 
$\kappa_{0}:\restr{\mc G_{0}}{(X\sm \mf t)\times s}\isom \restr{\mc G_{0}'}{(X\sm \mf t)\times s}$ définit  un élément de 
$G(\mc O_{\mf t})\backslash G(F_{\mf t})/G(\mc O_{\mf t})$ qui corresponde  à  $\lambda_{\mf t}$ (c'est-à-dire que si on trivialise $(\mc G_{0},\phi_{k}\cdots \phi_{1})$ au voisinage formel de $t$ la restriction de $(\mc G'_{0},\phi'_{k}\cdots \phi'_{1})$ à ce voisinage formel
détermine un point de $G(F_{\mf t})/G(\mc O_{\mf t})$ appartenant à la 
$G(\mc O_{\mf t})$-orbite de $\lambda_{\mf t}$). 
 \end{itemize}
De plus  
  $\on{pr}_{1}$ et $\on{pr}_{2}$
  sont données par  \eqref{GammaN-pr1} et  \eqref{GammaN-pr2}. 
\end{construction}

\section{Morphismes de  Frobenius partiels}

Les morphismes de  Frobenius partiels ont été introduits par Drinfeld \cite{drinfeld78,Dr1}. La généralisation à notre situation  est évidente. 

Le morphisme  de  Frobenius partiel 
$$\on {Fr}_{I_{1}, N,I} ^{(I_{1},...,I_{k})}: \Cht_{N,I} ^{(I_{1},...,I_{k})}\to \Cht_{N,I} ^{(I_{2},...,I_{k},I_{1})},$$ 
est défini par 
\begin{gather}\on {Fr}_{I_{1}, N,I} ^{(I_{1},...,I_{k})}\big( (x_i)_{i\in I}, (\mc G_{0}, \psi_{0}) \xrightarrow{\phi_{1}}  (\mc G_{1}, \psi_{1}) \xrightarrow{\phi_{2}}
\cdots\xrightarrow{\phi_{k-1}}  (\mc G_{k-1}, \psi_{k-1}) \xrightarrow{ \phi_{k}}    (\ta{\mc G_{0}}, \ta \psi_{0})
\big)\nonumber \\
= \nonumber
\big( (x'_i)_{i\in I}, (\mc G_{1}, \psi_{1}) \xrightarrow{\phi_{2}}  (\mc G_{2}, \psi_{2}) \xrightarrow{\phi_{3}}
\cdots \xrightarrow{ \phi_{k}}    (\ta{\mc G_{0}}, \ta \psi_{0}) \xrightarrow{\ta  \phi_{1}   } (\ta{\mc G_{1}}, \ta \psi_{1})
\big). 
\end{gather}
Il est  au-dessus du morphisme $\Frob_{I_{1}}:(X\sm N)^{I}\to (X\sm N)^{I}$ qui envoie  $(x_i)_{i\in I}$ vers  $(x'_i)_{i\in I}$ avec  
\begin{gather}\label{defi-Frob-I1}x_{i}'=\Frob (x_{i})\text{ \  si \  }i\in I_{1}\text{ \  et   \  }x_{i}'=x_{i}\text{ \  sinon.} \end{gather}

\begin{lem}\label{frob-troncature}
Les morphismes de  Frobenius partiels ne préservent pas les  troncatures, mais agissent sur  les systèmes inductifs. Plus précisément on a  
\begin{gather}\label{incl-Frob-p1} \big(\on  {Fr}_{I_{1}, N,I} ^{(I_{1},...,I_{k})}\big)^{-1}\big( \Cht_{N,I,\underline{\omega}} ^{(I_{2},...,I_{k},I_{1}),\leq\mu}
\big) \subset \Cht_{N,I,\underline{\omega}} ^{(I_{1},...,I_{k}),\leq\mu+\sum_{i\in I_{1}}\omega_{i}}\end{gather}
et 
\begin{gather}\label{incl-Frob-p2}\on  {Fr}_{I_{1}, N,I} ^{(I_{1},...,I_{k})}
\big(\Cht_{N,I,\underline{\omega}} ^{(I_{1},...,I_{k}),\leq\mu}\big) \subset \Cht_{N,I,\underline{\omega}} ^{(I_{2},...,I_{k},I_{1}),\leq\mu-w_{0}(\sum_{i\in I_{1}}\omega_{i} )}
\end{gather}
où $w_{0}$ est l'élément le plus long dans le groupe de   Weyl. \cqfd
\end{lem}

\begin{rem} \label{composee-frob-partiels} La composée 
 \begin{gather}
 \Cht_{N,I,\underline{\omega}} ^{(I_{1},...,I_{k})} \xrightarrow{ \on  {Fr}_{I_{1}, N,I} ^{(I_{1},...,I_{k})}}  \Cht_{N,I,\underline{\omega}} ^{(I_{2},...,I_{k},I_{1})}  \xrightarrow{ \on  {Fr}_{I_{2}, N,I} ^{(I_{2},...,I_{1})}} \cdots\xrightarrow{ \on  {Fr}_{I_{k}, N,I} ^{(I_{k},...,I_{k-1})}} \Cht_{N,I,\underline{\omega}} ^{(I_{1},...,I_{k})}
 \end{gather}
  est le  morphisme de Frobenius (total) de   $\Cht_{N,I,\underline{\omega}} ^{(I_{1},...,I_{k})}$ sur $\Fq$. 
  \end{rem} 

 \begin{prop} \label{action-frob-partiels-IC} On a un  isomorphisme canonique  
 $$ \on F_{I_{1}, N,I} ^{(I_{1},...,I_{k})}: 
 \big(\on  {Fr}_{I_{1}, N,I} ^{(I_{1},...,I_{k})}\big)^{*}\Big( 
 \mc F_{N,I,\underline{\omega}, \Xi,E}^{(I_{2},...,I_{k},I_{1})}
\Big) \simeq 
\mc F_{N,I,\underline{\omega}, \Xi,E}^{(I_{1},...,I_{k})}.
 $$
 \end{prop}
 \noindent 
 Par convention on {\it incorpore} dans l'isomorphisme $\on F_{I_{1}, N,I} ^{(I_{1},...,I_{k})}$ 
 le scalaire $q^{-d/2}$, où $d$ est la dimension de la strate 
 $\on{Gr}_{I_{1},\omega_{1}}^{(I_{1})}$ (cela se comprend bien par la \propref{prop-IC-produit-Fr-partiels} ci-dessous qui pourrait être considérée comme une définition alternative de $\on F_{I_{1}, N,I} ^{(I_{1},...,I_{k})}$). 
 Pour les signes  on applique la contrainte de commutativité modifiée de 
 \cite{mv}, rappelée dans la  \remref{contrainte-commutativite}. 

\noindent{\bf Démonstration. } La remarque  \ref{composee-frob-partiels}  implique que les   morphismes de  Frobenius partiels sont  des homéomorphismes locaux complètement radiciels, et on rappelle que  $\mc F_{N,I,\underline{\omega}, \Xi,E}^{(I_{1},...,I_{k})}$ est (à un décalage près) le faisceau d'intersection de $\Cht_{N,I,\underline{\omega}} ^{(I_{1},...,I_{k})}$, et il en va de même pour $(I_{2},...,I_{k},I_{1})$.   
 \cqfd
 
 Dans le \corref{cor-IC-produit-gen} nous avons construit  un isomorphisme entre $\mc F_{N,I,\underline{\omega}, \Xi,E}^{(I_{1},...,I_{k})}$ et l'image inverse d'un produit extérieur de  faisceaux de Mirkovic-Vilonen sur 
 $\mr{Gr}_{I_{j},(\omega_{i})_{i\in I_{j}}}^{(I_{j})}/G^{\mr{ad}}_{\sum _{i\in I_{j}}n_{i}x_{i}}$. Comme ces images inverses se comportent plus  canoniquement par rapport à   la coalescence des pattes (grâce au  produit de  fusion de \cite{mv}), il est nécessaire de comprendre  $ \on F_{I_{1}, N,I} ^{(I_{1},...,I_{k})}$ en ces termes (car cela permettra de montrer la compatibilité entre les morphismes de Frobenius partiels et les isomorphismes de coalescence des pattes). 
  Dans la  proposition suivante nous l'exprimerons 
  à l'aide du  morphisme de  Frobenius de  $\mr{Gr}_{I_{1},(\omega_{i})_{i\in I_{1}}}^{(I_{1})}/G^{\mr{ad}}_{\sum _{i\in I_{1}}n_{i}x_{i}}$, que nous noterons simplement  
  $\Frob_{1}$, et de l'isomorphisme canonique  
 \begin{gather}\label{Frob-Gr-1} \on F_{1}:  
 \Frob_{1}^{*}\big( \mc S_{I_{1},\boxtimes_{i\in I_{1}} V_{\omega_{i}},E}^{(I_{1})}
\big) \simeq \mc S_{I_{1},\boxtimes_{i\in I_{1}} V_{\omega_{i}},E}^{(I_{1})}
\end{gather} 
 où $\mc S_{I_{1},\boxtimes_{i\in I_{1}} V_{\omega_{i}},E}^{(I_{1})}$  est  le  faisceau pervers (à un décalage près) sur 
 $\mr{Gr}_{I_{1},(\omega_{i})_{i\in I_{1}}}^{(I_{1})}/G^{\mr{ad}}_{\sum _{i\in I_{1}}n_{i}x_{i}}$ dont l'image inverse à $ \mr{Gr}_{I_{1},(\omega_{i})_{i\in I_{1}}}^{(I_{1})} $
 est le faisceau d'intersection.  
 
 \begin{prop}\label{prop-IC-produit-Fr-partiels}
Le diagramme 
\begin{gather*}
\begin{CD}
 \Cht_{N,I,\underline{\omega}} ^{(I_{1},...,I_{k})}/\Xi @>{\on  {Fr}_{I_{1}, N,I} ^{(I_{1},...,I_{k})}}>> \Cht_{N,I,\underline{\omega}} ^{(I_{2},...,I_{1})} /\Xi \\
 @VV{\epsilon_{N,(I_{1},...,I_{k}),\underline{\omega},\underline{n}}^{(I_{1},...,I_{k})}}V @VV{\epsilon_{N,(I_{2},...,I_{1}),\underline{\omega},\underline{n}}^{(I_{2},...,I_{1})}}V \\
\prod_{j=1}^{k}
\mr{Gr}_{I_{j},(\omega_{i})_{i\in I_{j}}}^{(I_{j})}/G^{\mr{ad}}_{\sum _{i\in I_{j}}n_{i}x_{i}}
 @>{\Frob_{1}
\times \Id}>>  \prod_{j=1}^{k} \mr{Gr}_{I_{j},(\omega_{i})_{i\in I_{j}}}^{(I_{j})}/G^{\mr{ad}}_{\sum _{i\in I_{j}}n_{i}x'_{i}}
 \end{CD}
 \end{gather*}   est commutatif et rend compatibles  les 
 isomorphismes de  comparaison 
 entre  $
\mc F_{N,I,\underline{\omega}, \Xi,E}^{(I_{1},...,I_{k})}$, $\mc F_{N,I,\underline{\omega}, \Xi,E}^{(I_{2},...,I_{k},I_{1})}$  et les faisceaux de  Mirkovic-Vilonen $\boxtimes _{j\in \{1,...,k\}}\mc S_{I_{j},\boxtimes_{i\in I_{j}} V_{\omega_{i}},E}^{(I_{j})}$,  donnés par  
 \begin{itemize}
 \item  $\on F_{I_{1}, N,I} ^{(I_{1},...,I_{k})}$ (défini dans la  \propref{action-frob-partiels-IC}),  
 \item  $F_{1}\times \Id$ (défini dans  \eqref{Frob-Gr-1}), \item et les  isomorphismes 
 $\lambda_{N,(I_{1},...,I_{k})}^{(I_{1},...,I_{k})}$, $\lambda_{N,(I_{2},...,I_{1})}^{(I_{2},...,I_{1})}$ (définis dans le    \corref{cor-IC-produit-gen}). \end{itemize}
  \end{prop}
\noindent{\bf Démonstration.} 
Il suffit de montrer la compatibilité sur un ouvert dense. 
Or sur le lieu de lissité les   faisceaux d'intersection sont triviaux à un décalage  et une torsion à la Tate près, 
 et les isomorphismes de comparaison sont l'identité. Le scalaire 
 $q^{-d/2}$ évoqué après l'énoncé de la \propref{action-frob-partiels-IC} apparaît dans  
 $F_{1}$. Dans la ligne du bas du diagramme ci-dessus on ne change pas l'ordre des pattes, donc il ne se pose pas de question de signes. Mais si on changeait l'ordre des pattes on appliquerait la contrainte de commutativité modifiée, comme on l'a dit après l'énoncé de la \propref{action-frob-partiels-IC}. 
  \cqfd 
 
  \begin{rem}
  Le diagramme de la  \propref{prop-IC-produit-Fr-partiels} est cartésien à des homéomorphismes locaux totalement radiciels près (qui sont localement le produit d'un isomorphisme avec le morphisme de Frobenius de 
$G^{\mr{ad}}_{\sum _{i\in I_{1}}n_{i}x_{i}}$). 
\end{rem}

 \section{Faisceaux de cohomologie et coalescence des pattes}
 \label{para-def-coh}

Ce chapitre ne contient rien de neuf et à part la    construction de l'action des morphismes de Frobenius partiels  il est entièrement extrait de~\cite{var} et~\cite{brav-var}. 
 
 \subsection{Rappels sur les correspondances cohomologiques. }
 
  A la suite de  \cite{sga5,varshavsky-fujiwara,brav-var} on appelle correspondance de $X_{1}$ vers $X_{2}$ un  morphisme $a=(a_{1},a_{2}):A\to X_{1}\times X_{2}$ de champs algébriques localement de type fini sur $\Fq$ tel que $a_{2}$ soit  représentable et de type fini. Alors, pour $\mc F_{1}\in D_{c}^{b}(X_{1}, E)$ et $\mc F_{2}\in D_{c}^{b}(X_{2}, E)$, un morphisme 
  $$u:a_{2,!}(a_{1}^{*}(\mc F_{1}))\to \mc F_{2}, \text{ \ ou, de fa\c con équivalente \ } 
 u:a_{1}^{*}(\mc F_{1})\to a_{2}^{!}(\mc F_{2})$$ 
  est appelé une correspondance cohomologique   de 
  $\mc F_{1}$ vers $\mc F_{2}$ (ou de 
    $(X_{1},\mc F_{1})$ vers $(X_{2},\mc F_{2})$) supportée par $A$.
 
 Si $a_{1}$ est propre elle   induit un  morphisme sur les  cohomologies à  support compact. En effet soit $Y$  un  champ et 
  $p_{i}:X_{i}\to Y$ (pour  $i=1,2$) des morphismes tels que $p_{1}\circ a_{1}=p_{2}\circ a_{2}$. 
  On va supposer que $X_{1},X_{2}$ et $Y$ sont des  champs de Deligne-Mumford parce que ce sera le cas lorsqu'on l'appliquera 
  et
  que cela permet de rester dans la catégorie dérivée bornée.  
   Alors 
la  correspondance cohomologique  $u$ induit un  morphisme 
  $H(u): p_{1,!}(\mc F_{1}) \to p_{2,!}(\mc F_{2})$ dans  $D_{c}^{b}(Y, E)$ donné  par 
  \begin{gather*}p_{1,!}(\mc F_{1}) \xrightarrow{\on{adj}} p_{1,!} a_{1,*} a_{1}^{*}(\mc F_{1})
  = p_{1,!} a_{1,!} a_{1}^{*}(\mc F_{1}) \\ 
  \xrightarrow{u} p_{1,!} a_{1,!} a_{2}^{!}(\mc F_{2})  =
 p_{2,!} a_{2,!} a_{2}^{!}(\mc F_{2}) \xrightarrow{\on{adj}} p_{2,!}(\mc F_{2})\end{gather*}
  où on a utilisé, dans la deuxième étape, l'hypothèse que  $a_{1}$ est propre. 
  
   Si $b=(b_{1},b_{2}):B\to X_{2}\times X_{3}$ 
   est une autre correspondance, 
  $\mc F_{3}\in D_{c}^{b}(X_{3}, E)$
   et $v:b_{2,!}(b_{1}^{*}(\mc F_{2}))\to \mc F_{3}$ est une correspondance cohomologique   de 
     $\mc F_{2}$ vers $\mc F_{3}$ supportée par $B$, la  composée de $v$ et $u$ est obtenue de la fa\c con suivante. 
     On pose $C=A\times_{X_{2}}B$. 
     On a un diagramme  commutatif
          $$   \xymatrix{
    &  &   C\ar[dl]^{\wt b_{2}} \ar[dr]_{\wt a_{2}} 
       & &
     \\
    &  A\ar[dl]^{a_{1}} \ar[dr]_{a_{2}}&    & B\ar[dl]^{b_{2}} \ar[dr]_{b_{3}} &
     \\
  X_{1}     &     & X_{2} 
  &&  X_{3}
    }$$     
    où le carré est cartésien.      D'où une  correspondance 
     $c= 
     (a_{1}\wt b_{2},b_{3}\wt a_{2}):C\to X_{1}\times X_{3}$. On définit la correspondance cohomologique   $v\circ u$ de $\mc F_{1}$ vers $\mc F_{3}$ supportée par  $C$ comme la  composée
     $$(b_{3}\wt a_{2})_{!}(a_{1}\wt b_{2})^{*}(\mc F_{1})=  b_{3,!}  \wt a_{2,!} \wt b_{2} ^{*} a_{1}^{*} (\mc F_{1})
     \simeq b_{3,!}b_{2}^{*}a_{2,!}a_{1}^{*}(\mc F_{1})\xrightarrow{u} 
     b_{3,!}b_{2}^{*}(\mc F_{2})\xrightarrow{v}  \mc F_{3}, $$
     où l'isomorphisme vient du changement de base propre $  \wt a_{2,!}\wt b_{2} ^{*}  \simeq b_{2}^{*}a_{2,!}$. 
   
   Si $a_{1}$ et $b_{2}$ sont propres et $p_{i}:X_{i}\to Y$ (pour $i=1,2,3$) sont des morphismes  tels que $p_{1}\circ a_{1}=p_{2}\circ a_{2}=p_{3}\circ a_{3}$ 
   (et $X_{1},X_{2},X_{3}$ et $Y$ sont des champs de Deligne-Mumford) 
   on a  
   $   H(v\circ u)=H(v)\circ H(u): p_{1,!}(\mc F_{1}) \to p_{3,!}(\mc F_{3})$. 
   
 \subsection{Définition des faisceaux de cohomologie}
 
 \begin{defi}\label{defi-coho}
 Comme  $\Cht_{N,I,\underline{\omega}}^{(I_{1},...,I_{k}),\leq\mu}/\Xi$
 est un   champ de Deligne-Mumford  de type fini, 
   d'après  \cite{laumon-moret-bailly,laszlo-olsson} on peut définir 
      \begin{gather} \label{def-H-Xi-E}
 \mc H _{ N,I,\underline{\omega}}^{\leq\mu,E}=  R
 \big(\mf p_{N,I}^{(I_{1},...,I_{k}),\leq\mu}\big)_{!}\big(\restr{\mc F_{N,I,\underline{\omega}, \Xi,E}^{(I_{1},...,I_{k})}}{\Cht_{N,I,\underline{\omega}}^{(I_{1},...,I_{k}),\leq\mu}/\Xi} \big) \in D^{b}_{c}((X\sm N)^{I}, E)
  \end{gather} où l'on rappelle que  $\mc F_{N,I,\underline{\omega}, \Xi,E}^{(I_{1},...,I_{k})}$ est le faisceau d'intersection  de $\Cht_{N,I,\underline{\omega}}^{(I_{1},...,I_{k})}/\Xi$ avec 
 le degré et le  poids normalisés relativement à  $ (X\sm N)^{I}$. 
Comme la   notation  $ \mc H _{ N,I,\underline{\omega}}^{\leq\mu,E}$ l'indique,   le membre de droite de \eqref{def-H-Xi-E}
  ne dépend pas du choix de la  partition
 $(I_{1},...,I_{k})$, en vertu du 
 \corref{cor-strates}. 
 Bien sûr $\mc H_{N,I,\underline{\omega}}^{\leq\mu,E}$ dépend de $\Xi$ mais on omet $\Xi$ de la notation pour raccourcir un peu.

 Soit  $i\in \Z$. 
La cohomologie de degré $i$ (pour la   $t$-structure ordinaire)
  \begin{gather} 
 \mc H _{ N,I,\underline{\omega}}^{i,\leq\mu,E}=H^{i}\big(  \mc H _{ N,I,\underline{\omega}}^{\leq\mu,E}\big) =  R^{i}
 \big(\mf p_{N,I}^{(I_{1},...,I_{k}),\leq\mu}\big)_{!} \big(\restr{\mc F_{N,I,\underline{\omega}, \Xi,E}^{(I_{1},...,I_{k})}}{\Cht_{N,I,\underline{\omega}}^{(I_{1},...,I_{k}),\leq\mu}/\Xi} \big)
  \end{gather}
 est  un  $E$-faisceau  constructible sur  $(X\sm N)^{I}$.  
 On a préféré commencer par cette définition qui est sans doute naturelle pour le lecteur familier des variétés de Shimura, mais on insiste sur le fait que 
 la bonne définition, équivalente à celle ci-dessus mais   fonctorielle en $W$ représentation de $(\wh G)^{I}$, sera donnée dans \eqref{def-bonne-HNIW} de la \defiref{prop-def-equiv-HcNIW}.  
    \end{defi}
 
 \begin{rem}
 Lorsque $I$ est vide et $\underline{\omega}$ est donc nul, 
$ \mc H _{ N,\emptyset,0}^{ \leq\mu,E} $ est le faisceau constant sur $\on{Spec} \Fq$ supporté en degré $0$ et égal à $C_{c}(\Bun_{G,N}^{\leq\mu}(\Fq)/\Xi ,E)$. 
 \end{rem}
 
 \begin{rem}
  Pour tout niveau   $N'\supset N$ de degré assez grand  en  fonction de $\mu$, le  morphisme 
  $$\Cht_{N',I,\underline{\omega}}^{(I_{1},...,I_{k}),\leq\mu}/\Xi
  \xrightarrow{\mf p_{N',I}^{(I_{1},...,I_{k}),\leq\mu}} (X\sm N')^{I}$$
 est  représentable  quasi-projectif de type fini, et 
  $\restr{\Cht_{N,I,\underline{\omega}}^{(I_{1},...,I_{k}),\leq\mu}/\Xi}{(X\sm N')^{I}}$
est le   quotient de  $\Cht_{N',I,\underline{\omega}}^{(I_{1},...,I_{k}),\leq\mu}/\Xi$
 par le groupe fini $\on{Ker}(G(\mc O_{N'})\to G(\mc O_{N}))$. 
 Par conséquent on pourrait définir les  faisceaux   $\mc H _{ N,I,\underline{\omega}}^{i,\leq\mu,E}$  
  en employant seulement la cohomologie  étale  des schémas (puis en prenant les  coinvariants par le groupe fini susmentionné). 
\end{rem}

 \subsection{Action des   morphismes de Frobenius partiels. }  
 
 Le  lemme   \ref{frob-troncature} et la \propref{action-frob-partiels-IC} fournissent, pour toute partie  $J$  de $I$
le morphisme suivant dans    $D^{b}_{c}((X\sm N)^{I}, E)
$ : 
\begin{gather}\label{defi-F-J}     F_{J}: 
 \Frob_{J}^{*}(\mc H _{ N,I,\underline{\omega}}^{\leq\mu,E})\to 
\mc H _{ N,I,\underline{\omega}}^{\leq\mu+\sum_{i\in J}\omega_{i},E}, 
\end{gather}
où   $\Frob_{J}:(X\sm N)^{I}\to (X\sm N)^{I}$ a été défini dans  \eqref{defi-Frob-I1}. 
Pour toute partition $(I_{1},...,I_{k})$ telle que $I_{1}=J$, ce morphisme est associé à la correspondance cohomologique 
de $(\Cht_{N,I,\underline{\omega}}^{(I_{2},...,I_{k},I_{1}),\leq\mu}/\Xi , \mc F_{N,I,\underline{\omega}, \Xi,E}^{(I_{2},...,I_{k},I_{1})})$
vers  $(\Cht_{N,I,\underline{\omega}}^{(I_{1},...,I_{k}),\leq\mu'}/\Xi, 
\mc F_{N,I,\underline{\omega}, \Xi,E}^{(I_{1},...,I_{k})})$ 
donnée par 
\begin{itemize}
\item le diagramme 
$$\xymatrix{ \Cht_{N,I,\underline{\omega}}^{(I_{2},...,I_{k},I_{1}),\leq\mu}/\Xi  \ar[d]_{\mf p_{N,I}^{(I_{2},...,I_{k},I_{1}),\leq\mu}} &  a_{1}^{-1}(
\Cht_{N,I,\underline{\omega}}^{(I_{2},...,I_{k},I_{1}),\leq\mu}/\Xi )
\ar[l]_-{a_{1}}
\ar@{^{(}->}[r] \ar[d]_{\mf p_{N,I}^{(I_{1},...,I_{k}),\leq\mu}}  & \Cht_{N,I,\underline{\omega}}^{(I_{1},...,I_{k}),\leq\mu'}/\Xi \ar[d]^{\mf p_{N,I}^{(I_{1},...,I_{k}),\leq\mu'}}
\\(X\sm N)^{I} &  
(X\sm N)^{I}\ar[l]_{\Frob_{J}} \ar[r]^{\Id} &  (X\sm N)^{I}}$$
où $a_{1}=\on  {Fr}_{I_{1}, N,I} ^{(I_{1},...,I_{k})}$, $\mu'=\mu+\sum _{i\in I_{1}}\omega_{i}$ et l'inclusion ouverte vient de  \eqref{incl-Frob-p1}, 
\item le morphisme $\on F_{I_{1}, N,I} ^{(I_{1},...,I_{k})}: a_{1}^{*}(\mc F_{N,I,\underline{\omega}, \Xi,E}^{(I_{2},...,I_{k},I_{1})})\to \mc F_{N,I,\underline{\omega}, \Xi,E}^{(I_{1},...,I_{k})}$
   (défini dans la 
\propref{action-frob-partiels-IC}). 
\end{itemize}
On peut remarquer qu'ici, bien que le carré de gauche ne soit pas cartésien, on est essentiellement dans une situation de changement de base propre car   $a_{1}=\on  {Fr}_{I_{1}, N,I} ^{(I_{1},...,I_{k})}$ et $\Frob_{J}$ sont des homéomorphismes  locaux totalement radiciels, et donc le morphisme de
$a_{1}^{-1}(
\Cht_{N,I,\underline{\omega}}^{(I_{2},...,I_{k},I_{1}),\leq\mu}/\Xi )$ 
dans le produit fibré associé au carré de gauche est un homéomorphisme  local totalement radiciel : 
 c'est pourquoi nous avions dit dans  l'introduction que  $F_{J}$ était obtenu  par  changement de base propre.
  Grâce au  \corref{cor-strates} le  morphisme $F_{J}$ ne dépend pas du choix de la   partition $(I_{1},...,I_{k})$ de $I$ telle que $I_{1}=J$. 
Pour toute partie   $J$ de $I$, $$F_{J}\Frob_{J}^{*}(F_{I\sm J}) : \Frob^{*}\mc H _{ N,I,\underline{\omega}}^{\leq\mu,E}\to 
\mc H _{ N,I,\underline{\omega}}^{\leq\mu+\sum_{i\in I}\omega_{i},E}, 
$$ 
est la  composée de l'action naturelle du morphisme de  Frobenius  (total) et d'une augmentation de la   troncature. 
  
 \subsection{Action des algèbres de   Hecke}\label{action-hecke-etale}

 Soit  $f\in C_{c}(K_{N}\backslash G(\mb A)/K_{N},E)$ et  $\mf T$   un ensemble fini de places de  $X$ contenant $|N|$ et tel que 
  \begin{gather}\label{condition-T}f=\bigotimes _{v\not\in \mf T}{\bf 1}_{G(\mc O_{v})} \otimes f'.\end{gather}  
 Alors la construction \ref{constr-corresp-hecke} fournit, pour $\kappa$ assez grand  en fonction de  $f$,  un morphisme
  $$T(f)\in \Hom_{D^{b}_{c}((X\sm \mf T)^{I}, E)} \big(\restr{\mc H _{ N,I,\underline{\omega}}^{\leq\mu,E}}{(X\sm \mf T)^{I}},  \restr{\mc H _{ N,I,\underline{\omega}}^{\leq\mu+\kappa,E}}{(X\sm \mf T)^{I}}\big) $$
  tel que si $f$ est la  fonction caractéristique de  $K_{N}gK_{N}$ 
  (multipliée  par la mesure de Haar   sur $G(\mathbb A)$ telle que $K_{N}$ ait pour  volume $1$), $T(f)$ est l'action de la  correspondance finie 
  $\Gamma_{N}(g)$, c'est-à-dire   $\on{pr}_{1,!}g^{*}\on{pr}_{2}^{*}$ 
  dans les  notations du  diagramme \eqref{corresp-hecke}.  
 Ces actions sont compatibles avec la   composition: si $f_{1}, f_{2}\in C_{c}(K_{N}\backslash G(\mb A)/K_{N},E)$ satisfont la   condition \eqref{condition-T} avec   $\mf T_{1}$ et $\mf T_{2}$, alors pour $\kappa$ assez grand  en fonction de  $f_{1}$ et  $f_{2}$, on a  
  \begin{gather}\label{Tf1f2}T(f_{1}f_{2})=T(f_{1})T(f_{2})\end{gather} dans   $$\Hom_{D^{b}_{c}((X\sm (\mf T_{1}\cup \mf T_{2}))^{I}, E)}\Big(\restr{\mc H _{ N,I,\underline{\omega}}^{\leq\mu,E}}{(X\sm (\mf T_{1}\cup \mf T_{2}))^{I}}, \restr{\mc H _{ N,I,\underline{\omega}}^{\leq\mu+\kappa,E}}{(X\sm (\mf T_{1}\cup \mf T_{2}))^{I}}\Big). $$
  On peut définir aussi les morphismes $T(f)$ à l'aide de correspondances cohomologiques (finies) et déduire 
  \eqref{Tf1f2} des règles de composition des corrrespondances. On note que ces correspondances cohomologiques sont définies sur $\mc O_{E}$ si $f$ l'est.

Par ailleurs l'action des algèbres de   Hecke  commute avec l'action des morphismes de Frobenius partiels. 
     
  Dans le  \corref{cor-hecke-etendus-compo} on verra que 
  ces morphismes $T(f)$ peuvent être étendus à $(X\sm N)^{I}$.

  \subsection{Définition fonctorielle en $W$}

 Soit  $W$  une   représentation $E$-linéaire de dimension finie de  $(\wh G)^{I}$. 
Le  \corref{cor-IC-produit} et l'équivalence  de   Satake géométrique 
vont fournir une  définition fonctorielle en $W$ de  $\mc H _{ N,I,W}^{\leq\mu,E}$  
 comme la  cohomologie d'un  certain faisceau pervers (à décalage près). Soit  $(I_{1},...,I_{k})$ une  partition de $I$. 
 
 \begin{notation} En tant que  représentation de   $(\wh G)^{I}$,   $W$ admet une   unique décomposition de la  forme 
\begin{gather}\label{decomp-W-Womega}
W=\bigoplus _{\underline{\omega}\in (X_{*}^{+}(T))^{I}} \Big(\bigotimes_{i\in I}V_{\omega_{i}} \Big)\otimes_{E} \mf W_{\underline{\omega}}
\end{gather} où les  $\mf W_{\underline{\omega}}$ sont des  $E$-espaces vectoriels de  dimension finie,  presque tous nuls. 
 On note 
$\Cht_{N,I,W}^{(I_{1},...,I_{k}),\leq\mu}$  la réunion des 
$\Cht_{N,I,\underline{\omega}}^{(I_{1},...,I_{k}),\leq\mu}$ 
pour   $\underline{\omega}$
tel que   $\mf W_{\underline{\omega}}$ est non nul. On rappelle que dans le chapitre \ref{rappels-Hecke-Gr-satake} on a défini de la même fa\c con  $\mr{Gr}_{I,W}^{(I_{1},...,I_{k})}/G_{\sum n_{i}x_{i}}$ (où les  entiers $n_{i}$ sont assez grands  en fonction de  $W$). 
\end{notation}
Par  l'équivalence  de   Satake géométrique  (sous la forme où on l'a rappelée dans  le \thmref{thm-geom-satake}, voir  \cite{mv, hitchin, ga-de-jong} pour plus de détails) on a le  faisceau pervers (à un décalage près)  $\mc S_{I,W,E}^{(I_{1},...,I_{k})}$ sur  $\mr{Gr}_{I,W}^{(I_{1},...,I_{k})}/G_{\sum n_{i}x_{i}}$
 (on rappelle qu'il est normalisé de telle sorte 	que son  image inverse  sur  $\mr{Gr}_{I,W}^{(I_{1},...,I_{k})}$ soit perverse relativement à  $X^{I}$). D'après 
   la \remref{rem-Satake-Gad-texte}, $\mc S_{I,W,E}^{(I_{1},...,I_{k})}$   est en fait un faisceau pervers (à un décalage près) sur 
$\mr{Gr}_{I,W}^{(I_{1},...,I_{k})}/G^{\mr{ad}}_{\sum _{i\in I}n_{i}x_{i}}$.

On rappelle, d'après la 
 proposition \ref{lissite-Cht-Grass}
et la notation qui précède le \corref{cor-IC-produit}, que 
   $$\epsilon_{N,(I),W,\underline n}^{(I_{1},...,I_{k}), \Xi} :\Cht_{N,I,W}^{(I_{1},...,I_{k})} /\Xi  \to 
\mr{Gr}_{I,W}^{(I_{1},...,I_{k})}/G^{\mr{ad}}_{\sum n_{i}x_{i}}
$$  est le morphisme lisse qui à un chtouca
 \eqref{donnee-chtouca} associe 
\begin{gather}\label{point-Gr-Gnixi}
(\mc G_{0} \xrightarrow{\phi_{1}}  \mc G_{1} \xrightarrow{\phi_{2}}
\cdots\xrightarrow{\phi_{k}}  
\mc G_{k})\in \mr{Gr}_{I,W}^{(I_{1},...,I_{k})}/G^{\mr{ad}}_{\sum n_{i}x_{i}}\end{gather} 
 (avec les notations de  la remarque \ref{rem-apres-action-Gnx-Grass}, en particulier le  $G^{\mr{ad}}_{\sum n_{i}x_{i}}$-torseur  tautologique sur le point \eqref{point-Gr-Gnixi} est déduit du $G_{\sum n_{i}x_{i}}$-torseur      $\restr{\mc G_{k}}{\Gamma_{\sum n_{i}x_{i}}}$).  

\begin{defi}\label{defi-F-E-Cht}
On définit  $\mc F_{N,I,W, \Xi,E}^{(I_{1},...,I_{k})}$ comme  le  faisceau pervers (à un décalage près) sur $\Cht_{N,I,W}^{(I_{1},...,I_{k}),\leq\mu}/\Xi$ égal à 
\begin{gather}\label{def-IC-V}\mc F_{N,I,W, \Xi,E}^{(I_{1},...,I_{k})}=
\big(\epsilon_{N,(I),W,\underline{n}}^{(I_{1},...,I_{k}), \Xi}\big)^{*}
\big(\mc S_{I,W,E}^{(I_{1},...,I_{k})}\big).\end{gather}
\end{defi}

\begin{rem} \label{rem-defi-eq-F-Vomega} Ceci généralise la \defiref{defi-F-NIomega}: lorsque 
$W=\boxtimes _{i}V_{\omega_{i}}$ est irréductible, 
$\mc F_{N,I,W, \Xi,E}^{(I_{1},...,I_{k})}$ est égal à 
$\mc F_{N,I,\underline{\omega}, \Xi,E}^{(I_{1},...,I_{k})}$, 
grâce au  \corref{cor-IC-produit}. 
Cela est vrai canoniquement (et non pas au produit tensoriel près par une droite vectorielle)  si on définit 
$V_{\omega_{i}}$ comme  la cohomologie totale du faisceau d'intersection de la strate fermée  associée à $\omega_{i}$ dans la  grassmannienne affine (car  le foncteur fibre donnant l'équivalence de Satake géométrique est la  cohomologie totale, comme on l'a rappelé dans le \thmref{thm-geom-satake}).
\end{rem}

La    définition   suivante de $\mc H _{ N,I,W}^{\leq\mu,E}$ est fonctorielle en $W$ et de plus   elle permettra de construire  l'isomorphisme de coalescence  de la  \propref{chgt-base-propre} d'une fa\c con {\it canonique}. 

\begin{defi} \label{prop-def-equiv-HcNIW}
On pose 
 \begin{gather} \label{def-bonne-HNIW}
 \mc H _{ N,I,W}^{\leq\mu,E}=  R
 \big(\mf p_{N,I}^{(I_{1},...,I_{k}),\leq\mu}\big)_{!}\Big(\restr{\mc F_{N,I,W,\Xi,E}^{(I_{1},...,I_{k})}} {\Cht_{N,I,W}^{(I_{1},...,I_{k}),\leq\mu}/\Xi}\Big) . 
  \end{gather}
  \end{defi}
 
 \begin{rem} En utilisant la décomposition  \eqref{decomp-W-Womega}
 avec la définition canonique des $V_{\omega_{i}}$ expliquée dans  la \remref{rem-defi-eq-F-Vomega}, on peut réécrire la définition précédente sous la forme 
 \begin{gather}\mc H _{ N,I,W}^{\leq\mu,E}=\bigoplus _{\underline{\omega}\in (X_{*}^{+}(T))^{I}} \mc H _{ N,I,\underline{\omega}}^{\leq\mu,E}\otimes_{E} \mf W_{\underline{\omega}}.\end{gather} 
   Cette formule permet de voir la fonctorialité en $W$, mais elle est beaucoup moins commode que   \eqref{def-bonne-HNIW}, notamment pour étudier la coalescence.  
\end{rem}

\begin{notation} \label{notation-Hu}
D'après   la définition précédente,  $W\mapsto \mc H _{ N,I,W}^{\leq\mu,E}$ est un foncteur $E$-linéaire. Si $u: W_{1}\to W_{2}$ est un  morphisme de représentations de  $(\wh G)^{I}$ on note  $$\mc H(u)\in \Hom_{D^{b}_{c}((X\sm N)^{I}, E)} \big(\mc H _{ N,I,W_{1}}^{\leq\mu,E},  \mc H _{ N,I,W_{2}}^{\leq\mu,E}\big)
$$
le morphisme associé à $u$. 
\end{notation}
 
 Une autre  définition  équivalente de $\mc F_{N,I,W, \Xi,E}^{(I_{1},...,I_{k})}$ est donnée par la proposition suivante. Elle permettra de réexprimer la \propref{prop-IC-produit-Fr-partiels} (qui explicite  les   actions des morphismes de  Frobenius partiels) de fa\c con fonctorielle en $W$. 
 \begin{prop}\label{prop-F-eq-def-Frob}
 Soit $W=\boxtimes_{j\in \{1,...,k\}} W_{j}$ où $W_{j}$ est une  représentation de $(\wh G)^{I_{j}}$. On a alors un   isomorphisme canonique  
\begin{gather}\label{def-IC-V-var}\mc F_{N,I,W, \Xi,E}^{(I_{1},...,I_{k})}=
\big(\epsilon_{N,(I_{1},...,I_{k}),W,\underline{n}}^{(I_{1},...,I_{k}), \Xi}\big)^{*}
\Big(\boxtimes_{j\in \{1,...,k\}} \mc S_{I_{j},W_{j},E}^{(I_{j})}\Big) 
.\end{gather}
Cet  isomorphisme est compatible avec l'action du morphisme de Frobenius   partiel associé à $I_{1}$. 
\end{prop}
\dem Cela résulte du 
lemme \ref{kapp-eps-J-lemme-smooth},    du
\corref{cor-IC-produit-gen} et de la \propref{prop-IC-produit-Fr-partiels}.  
\cqfd

\begin{rem} La  \defiref{defi-F-E-Cht}  et la  proposition    \ref{prop-F-eq-def-Frob}  admettent la généralisation suivante, qui ne sera jamais utilisée.  Soit $(J_{j,l})_{j,l}$ comme  dans la remarque \ref{rem-Jkl}. On a alors un   isomorphisme canonique 
$$\mc F_{N,I,W, \Xi,E}^{(J_{j,l})_{j,l}}=
\big(\epsilon_{N,(I_{1},...,I_{k}),W,\underline{n}}^{(J_{j,l})_{j,l}, \Xi}\big)^{*}
\Big(\boxtimes_{j\in \{1,...,k\}} \mc S_{I_{j},W_{j},E}^{(J_{j,1}, ..., J_{j,m_{j}})}\Big) 
.$$ 
\end{rem}

  La proposition suivante 
  jouera un rôle fondamental pour étudier la coalescence. 
   
 \begin{prop}\label{chgt-base-propre} (remarque 2.3.2 de \cite{brav-var}).  
 Soit  $I,J$ des ensembles finis et  $\zeta: I\to J$ une application. On note   $$\Delta_{\zeta} : X^{J}\to X^{I}, (x_{j})_{j\in J}\mapsto (x_{\zeta(i)})_{i\in I}$$ 
 le  morphisme diagonal associé à $\zeta$. 
 Soit  $W$ une    représentation  $E$-linéaire de  dimension finie  de $\wh G^{I}$. On note  $W^{\zeta}$ la   représentation de $\wh G^{J}$ qui est la composée de la   représentation $W$ avec le    morphisme diagonal $$\wh G^{J}\to \wh G^{I}, (g_{j})_{j\in J}\mapsto (g_{\zeta(i)})_{i\in I}. $$
  Alors l'isomorphisme canonique  rappelé dans  d) du \thmref{thm-geom-satake} fournit un   isomorphisme canonique  
 \begin{gather}\label{isom-chgt-base}
\chi_{\zeta} :  \Delta_{\zeta}^{*}(\mc H _{ N,I,W}^{\leq\mu,E})\isom \mc H _{ N,J,W^{\zeta}}^{\leq\mu,E} 
\text{ \ dans  \ }   D^{b}_{c}((X\sm N)^{J}, E)\end{gather}
que nous appellerons   isomorphisme de coalescence. Il est fonctoriel en $W$ et compatible avec la composition de $\zeta$. 
  \end{prop}
  \begin{rem}
  Cet énoncé est vrai sans décalage cohomologique ni torsion à la Tate parce que les  faisceaux d'intersection sur les champs  de chtoucas (ou sur les grassmanniennes affines de Beilinson-Drinfeld ) ont été normalisés relativement aux puissances de $X$.   Bien sûr on utilise  la contrainte de commutativité {\it modifiée}  de \cite{mv} (un exemple simple pour comprendre ce qui se passe est  le cas particulier où  $I=J$ et $\zeta $ est une permutation de $I$, de sorte que 
  $\Delta_{\zeta}:X^{I}\to X^{I}$ est le changement de l'ordre des pattes). 
     \end{rem}
    
 \noindent{\bf Démonstration.}  Il suffit de montrer le résultat quand  $\zeta$ est injective ou surjective.  Si $\zeta$ est  injective, cela résulte de  la propriété  c) de la \propref{prop-chtoucas}. 
On suppose maintenant que $\zeta$ est  surjective.  Soit $(J_{1},...,J_{k})$ une partition de $J$. On considère la partition $(I_{1},...,I_{k})$ de $I$ donnée par  $I_{j}=\zeta^{-1}(J_{j})$.   On a  un  isomorphisme  $$\kappa: \Cht_{N,I}^{(I_{1},...,I_{k})}\times_{(X\sm N)^{I}}(X\sm N)^{J}\isom \Cht_{N,J}^{(J_{1},...,J_{k})}
$$ défini tautologiquement  par 
\begin{gather*}\kappa\big(\big( ((x_{\zeta(i)})_{i\in I}), (\mc G_{0}, \psi_{0}) \xrightarrow{\phi_{1}}  
\cdots\xrightarrow{\phi_{k-1}}  (\mc G_{k-1}, \psi_{k-1}) \xrightarrow{ \phi_{k}}    (\ta{\mc G_{0}}, \ta \psi_{0})
\big),(x_j)_{j\in J} \big)\\
=\big( (x_j)_{j\in J}, (\mc G_{0}, \psi_{0}) \xrightarrow{\phi_{1}}  \cdots\xrightarrow{\phi_{k-1}}  (\mc G_{k-1}, \psi_{k-1}) \xrightarrow{ \phi_{k}}    (\ta{\mc G_{0}}, \ta \psi_{0})
\big).
\end{gather*}
On note   $$\on{pr}_{1}: \Cht_{N,I}^{(I_{1},...,I_{k})}
 \times_{(X\sm N)^{I}}(X\sm N)^{J}\to 
 \Cht_{N,I}^{(I_{1},...,I_{k})}$$
la première  projection.  
 Grâce à la \defiref{defi-F-E-Cht}  et à la propriété d) du \thmref{thm-geom-satake},  on obtient  un 
   isomorphisme    $\on{pr}_{1}^{*}(\mc F_{N,I,W,\Xi,E}^{(I_{1},...,I_{k})})\simeq \kappa^{*}(\mc F_{N,J,W^{\zeta},\Xi,E}^{(J_{1},...,J_{k})})$. Ces morphismes passent au quotient par $\Xi$. 
Grâce à la \defiref{prop-def-equiv-HcNIW} le   théorème  de changement de base propre fournit alors l'isomorphisme   \eqref{isom-chgt-base}. 
 \cqfd
 
  \begin{prop}\label{prop-coalescence-Hecke-Frob}
Les  actions des morphismes de  Frobenius partiels et des opérateurs de Hecke sont  compatibles  avec les morphismes $\mc H(u)$ de la notation \ref{notation-Hu} et avec  les isomorphismes de  coalescence de la proposition précédente. 
 \end{prop}
\noindent La      compatibilité entre les morphismes de Frobenius partiels et l'isomorphisme de coalescence \eqref{isom-chgt-base}   signifie 
    que pour tout $j\in J$, $\Delta_{\zeta}^{*}(F_{\zeta^{-1}(\{j\})})$ et 
    $F_{\{j\}}$ se correspondent par l'isomorphisme $\chi_{\zeta}$ 
     de \eqref{isom-chgt-base}. 
     
     \dem  Pour les opérateurs de Hecke cela est facile et pour les morphismes de  Frobenius partiels cela résulte de la  \propref{prop-IC-produit-Fr-partiels} réexprimée de fa\c con fonctorielle en $W$  grâce à la \propref{prop-F-eq-def-Frob}. 
\cqfd

     \section{Morphismes de création et d'annihilation}\label{para-creation-annih}
 
    Dans ce  chapitre on va  utiliser les   isomorphismes de coalescence  pour définir des morphismes de création et d'annihilation. 
    Ils joueront un rôle essentiel dans la  définition des opérateurs d'excursion.

Soit   $I$ et  $J$ des ensembles finis. 
 Dans la suite des réunions comme   $I\cup J$ et $I\cup \{0\}$ désigneront   {\it toujours} des réunions disjointes. 
 On introduit les applications évidentes 
  \begin{gather}\zeta_{J} :J \to \{0\}, \  \ 
 \zeta_{J}^{I}=(\Id_{I},\zeta_{J}): I\cup J\to I\cup \{0\} \ \text{ et }\zeta_{\emptyset}^{I}=(\Id_{I},\zeta_{\emptyset}): I \to I \cup \{0\}.\end{gather}

 Soit  $W$ et  $U$ des   représentations $E$-linéaires de dimension finie  de 
 $(\wh G)^{I}$ et   $(\wh G)^{J}$ respectivement. 
On rappelle que  $U^{\zeta_{J}}$ désigne  la représentation de  $\wh G$ obtenue en restreignant   $U$ à la diagonale. 
Soit  $x: \mbf 1\to U^{\zeta_{J}}$ et  $\xi: U^{\zeta_{J}}\to \mbf 1$ des morphismes de représentations de $\wh G$. En d'autres termes on se donne 
   $x\in U$ et $\xi\in U^{*}$     invariants par l'action diagonale de  $\wh G$.  

Alors  $W\boxtimes U$ est une   représentation de  $(\wh G)^{I\cup J}$,  et 
 $W\boxtimes U^{\zeta_{J}}$ et  $W\boxtimes \mbf 1$ sont des  représentations de  $(\wh G)^{I\cup \{0\}}$, reliées entre elles par les morphismes $\Id_{W}\boxtimes x$ et  $\Id_{W}\boxtimes \xi$. 
 
  On note  $\Delta:X\to X^{J}$ le  morphisme diagonal.

  On rappelle que  les   isomorphismes de coalescence    
 $$\mc H _{ N,I,W}^{\leq\mu,E} \boxtimes E_{X\sm N}\isor{\chi_{\zeta_{\emptyset}^{I}}}
  \mc H _{ N,I\cup\{0\},W\boxtimes \mbf 1}^{\leq\mu,E}$$ 
 et 
  $$   \restr{ \mc H _{ N,I\cup J,W\boxtimes U}^{\leq\mu,E}}{(X\sm N)^{I}\times \Delta(X\sm N)}
   \isor{ \chi_{\zeta_{J}^{I}}}
    \mc H _{ N,I\cup\{0\},W\boxtimes U^{\zeta_{J}}}^{\leq\mu,E}
 $$
   ont été  construits dans la  \propref{chgt-base-propre}. 
   
   Voici maintenant la définition des  morphismes de création et d'annihilation.  Ce sont des morphismes  dans 
  la catégorie  $D^{b}_{c}((X\sm N)^{I}\times (X\sm N), E)$ 
  mais on  gardera la même notation  pour les actions induites sur les  faisceaux de cohomologie   comme  $\mc H _{ N,I,W}^{0,\leq\mu,E}$ (au demeurant,  seules ces actions induites serviront dans la démonstration du théorème principal).

 \begin{defi}\label{morph-decoalescence}
Le morphisme  de création     $
\mc C_{x}^{\sharp}$ est défini  comme la  composée 
  \begin{gather*}
  \mc H _{ N,I,W}^{\leq\mu,E} \boxtimes E_{(X\sm N)} 
  \isor{\chi_{\zeta_{\emptyset}^{I}}}
 \mc H _{ N,I\cup\{0\},W\boxtimes \mbf 1}^{\leq\mu,E}  \\
 \xrightarrow{\mc H(\Id_{W}\boxtimes x) }
 \mc H _{ N,I\cup\{0\},W\boxtimes U^{\zeta_{J}}}^{\leq\mu,E} 
 \isor{ \chi_{\zeta_{J}^{I}}^{-1} }
 \restr{ \mc H _{ N,I\cup J,W\boxtimes U}^{\leq\mu,E}}{(X\sm N)^{I}\times \Delta(X\sm N)} . 
   \end{gather*}
    \end{defi}
    
     \begin{defi}\label{morph-coalescence}
Le  morphisme d'annihilation  $
\mc C_{\xi}^{\flat}$  est défini  comme la composée 
  \begin{gather*}
  \restr{ \mc H _{ N,I\cup J,W\boxtimes U}^{\leq\mu,E}}{(X\sm N)^{I}\times \Delta(X\sm N)} 
  \isor{ \chi_{\zeta_{J}^{I}}}
 \mc H _{ N,I\cup\{0\},W\boxtimes U^{\zeta_{J}}}^{\leq\mu,E} \\
 \xrightarrow{ \mc H(\Id_{W}\boxtimes \xi) }
  \mc H _{ N,I\cup\{0\},W\boxtimes \mbf 1}^{\leq\mu,E} 
 \isor{\chi_{\zeta_{\emptyset}^{I}}^{-1} }
   \mc H _{ N,I,W}^{\leq\mu,E} \boxtimes E_{(X\sm N)} . 
    \end{gather*}
       \end{defi}

Il résulte de la \propref{prop-coalescence-Hecke-Frob} que les  morphismes de création et d'annihilation commutent avec l'action des   opérateurs de Hecke et qu'ils sont compatibles avec l'action des  morphismes de  Frobenius partiels au sens suivant:  
\begin{itemize}
\item si  $K$ est une partie de  $I$, ils  commutent avec   $F_{K}$, 
\item ils entrelacent  $F_{J}$ et l'action naturelle du morphisme de  Frobenius partiel sur $E_{X\sm N}$ (qui correspond par  $\chi_{\zeta_{\emptyset}^{I}}$ à l'action de  $F_{\{0\}}$ sur  $ \mc H _{ N,I\cup\{0\},W\boxtimes \mbf 1}^{\leq\mu,E}$). 
\end{itemize} 

  Le lemme suivant est facile. Il dit que 
  \begin{itemize}
  \item 
 les  morphismes de création et d'annihilation commutent entre eux quand ils concernent  des sous-ensembles disjoints de pattes,    
 \item  la composée de deux morphismes de création
  associés à deux sous-ensembles disjoints de pattes placées au même  point de $X\sm N$  est le  morphisme de création
 associé à leur réunion, 
  \item  il en va de même pour les morphismes d'annihilation.
 \end{itemize}

  \begin{lem} \label{lem-lemme-4.3} Soit $I,J_{1},J_{2}$ des ensembles finis et  
     $W,U_{1},U_{2}$ des  représentations de
$(\wh G)^{I}$, $(\wh G)^{J_{1}}$ et $(\wh G)^{J_{2}}$. 
Pour $i=1,2$ on se donne des morphismes $x_{i}: \mbf 1\to U_{i}^{\zeta_{J_{i}}}$ et  $\xi_{i}: U_{i}^{\zeta_{J_{i}}}\to \mbf 1$ de représentations de  $\wh G$. 
 Alors   
 
    a) les composées 
     \begin{gather*}
    \restr{ \mc H _{ N,I\cup J_{1},W\boxtimes U_{1}}^{\leq\mu,E}}{(X\sm N)^{I}\times \Delta(X\sm N)}  \boxtimes E_{(X\sm N)} 
  \xrightarrow{\mc C_{ \xi_{1}}^{\flat }}
  \mc H _{ N,I,W}^{\leq\mu,E} \boxtimes E_{(X\sm N)^{2}}  \\
   \xrightarrow{\mc C_{  
x_{2}}^{\sharp }}
  \restr{ \mc H _{ N,I\cup J_{2},W\boxtimes U_{2}}^{\leq\mu,E}}{(X\sm N)^{I}\times \Delta(X\sm N)}    \boxtimes  
  E_{(X\sm N)} \end{gather*}
et 
 \begin{gather*}
    \restr{ \mc H _{ N,I\cup J_{1},W\boxtimes U_{1}}^{\leq\mu,E}}{(X\sm N)^{I}\times \Delta(X\sm N)}  
    \boxtimes E_{(X\sm N)}
    \xrightarrow{\mc C_{ 
x_{2}}^{\sharp } }   \\
\restr{ \mc H _{ N,I\cup J_{1}\cup J_{2},W\boxtimes U_{1}\boxtimes U_{2}}^{\leq\mu,E}}{(X\sm N)^{I}\times \Delta(X\sm N)^{2}} 
 \xrightarrow{\mc C_{ 
\xi_{1}}^{\flat } }
  \restr{ \mc H _{ N,I\cup J_{2},W\boxtimes U_{2}}^{\leq\mu,E}}{(X\sm N)^{I}\times \Delta(X\sm N)}  
   \boxtimes E_{(X\sm N)}   \end{gather*}
sont égales (dans ces formules les deux copies de $X\sm N$ ne sont pas bien distinguées mais la règle est simplement que   la création par $x_{2}$ agit sur l'une et l'annihilation par $\xi_{1}$ agit sur  l'autre).  

  b) la  composée 
  \begin{gather*}
   \mc H _{ N,I,W}^{\leq\mu,E} \boxtimes E_{(X\sm N)^{2}} 
    \xrightarrow{\mc C_{  
x_{1}}^{\sharp }}
 \restr{ \mc H _{ N,I\cup J_{1},W\boxtimes U_{1}}^{\leq\mu,E}}{(X\sm N)^{I}\times \Delta(X\sm N)}  \boxtimes E_{(X\sm N)} 
 \\
  \xrightarrow{\mc C_{  
x_{2}}^{\sharp }}
\restr{ \mc H _{ N,I\cup J_{1}\cup J_{2},W\boxtimes U_{1}\boxtimes U_{2}}^{\leq\mu,E}}{(X\sm N)^{I}\times (\Delta(X\sm N))^{2}} 
\end{gather*}
est  égale à la composée  
 \begin{gather*}
 \mc H _{ N,I,W}^{\leq\mu,E} \boxtimes E_{(X\sm N)^{2}} 
    \xrightarrow{\mc C_{  
x_{2}}^{\sharp }}
 \restr{ \mc H _{ N,I\cup J_{2},W\boxtimes U_{2}}^{\leq\mu,E}}{(X\sm N)^{I}\times \Delta(X\sm N)}  \boxtimes E_{(X\sm N)} \\
  \xrightarrow{\mc C_{  
x_{1}}^{\sharp } }
\restr{ \mc H _{ N,I\cup J_{1}\cup J_{2},W\boxtimes U_{1}\boxtimes U_{2}}^{\leq\mu,E}}{(X\sm N)^{I}\times (\Delta(X\sm N))^{2}} 
\end{gather*}
obtenue en changeant l'ordre des deux morphismes, 
et la  restriction de ces deux  composées à  $(X\sm N)^{I}\times \Delta(X\sm N)$ est égale à  
$\mc C_{  
x_{1}\boxtimes x_{2}}^{\sharp }$, 

  c) de   même, la  composée 
  \begin{gather*}
  \restr{ \mc H _{ N,I\cup J_{1}\cup J_{2},W\boxtimes U_{1}\boxtimes U_{2}}^{\leq\mu,E}}{(X\sm N)^{I}\times \Delta(X\sm N)^{2}} 
\xrightarrow{\mc C_{  
\xi_{2}}^{\flat } } \\
\restr{ \mc H _{ N,I\cup J_{1},W\boxtimes U_{1}}^{\leq\mu,E}}{(X\sm N)^{I}\times \Delta(X\sm N)}  \boxtimes E_{(X\sm N)}
  \xrightarrow{\mc C_{  
\xi_{1}}^{\flat }}
  \mc H _{ N,I,W}^{\leq\mu,E} \boxtimes E_{(X\sm N)^{2}} 
    \end{gather*}
est  égale à la  composée  obtenue en changeant l'ordre des deux morphismes, 
et la  restriction de ces deux  composées à   $(X\sm N)^{I}\times \Delta(X\sm N)$ est égale à 
$\mc C_{  
\xi_{1}\boxtimes \xi_{2}}^{\flat }$, 
\cqfd
  \end{lem}

  \begin{rem}\label{creation-annihilation-dualite}
Soit  $\theta$ une  involution de Chevalley  de  $\wh G$, telle que  
  $\theta(t)=t^{-1}$ pour 
  $t\in \wh T$. 
 Le  \thmref{thm-geom-satake} (qui rappelle 
 l'équivalence  de   Satake géométrique   \cite{mv, hitchin, ga-de-jong}) associe à toute représentation de dimension finie $W$ de $(\wh G)^{I}$
 un faisceau pervers  $G_{\sum n_{i}x_{i}}$-équivariant  $\mc S_{I,W,E}^{(I_{1},...,I_{k})}$ sur  $\mr{Gr}_{I,W}^{(I_{1},...,I_{k})}$
(avec la  normalisation perverse  relative à $X^{I}$). 
 On normalise la dualité de Verdier     relativement à  $X^{I}$. Alors 
  on a  un   isomorphisme canonique  
  $\mb D(\mc S_{I,W,E}^{(I_{1},...,I_{k})})\simeq 
  \mc S_{I,W^{*,\theta},E}^{(I_{1},...,I_{k})}
 $, fonctoriel en $W$, 
  où $W^{*,\theta}$ désigne la  représentation  contragrédiente de $(\wh G)^{I}$, composée avec   $\theta^{I}:(\wh G)^{I}\to (\wh G)^{I}$.   
Par conséquent  on a    $$\Cht_{N,I,W}^{(I_{1},...,I_{k}),\leq\mu}=\Cht_{N,I,W^{*,\theta}}^{(I_{1},...,I_{k}),\leq\mu}\text{ \ et  \ }
 \mb D(\mc F_{N,I,W,\Xi,E}^{(I_{1},...,I_{k})})\simeq\mc F_{N,I,W^{*,\theta},\Xi,E}^{(I_{1},...,I_{k})}, $$
 où le dernier  isomorphisme est canonique, fonctoriel en $W$ et compatible avec la  coalescence des pattes. 
On en déduit un  morphisme de faisceaux sur  $(X\sm N)^{I}$,  
\begin{gather}\label{bil-faisceaux} \mf B_{N,I,W}^{\Xi,E}: \mc H _{ N,I,W^{*,\theta}}^{0,\leq\mu,E}\otimes 
\mc H _{ N,I,W}^{0,\leq\mu,E}\to E_{(X\sm N)^{I}}.\end{gather}
Les  morphismes de création et d'annihilation sont transposés l'un de l'autre par rapport à  \eqref{bil-faisceaux}. 
Plus précisément, si $J$ est un autre ensemble fini, $U$ une représentation de dimension finie de $(\wh G)^{J}$ et $x\in U^{\zeta_{J}}$ alors pour toutes sections locales  $h$ et $h'$    de 
 $$  \mc H _{ N,I,W}^{0,\leq\mu,E} \boxtimes E_{(X\sm N)} \text{ \  et   \  }  \restr{ \mc H _{ N,I\cup J,W^{*,\theta}\boxtimes U^{*,\theta}}^{0,\leq\mu,E}}{(X\sm N)^{I}\times \Delta(X\sm N)}$$ sur   $(X\sm N)^{I}\times (X\sm N)$, 
on a  
 \begin{gather}\label{eg-transpose-t} 
 \Big(\restr{\mf B_{N,I\cup J,W\boxtimes U}^{\Xi,E}}{(X\sm N)^{I}\times \Delta(X\sm N)}\Big) \Big( h' \otimes  \mc C_{  x}^{\sharp }(h)\Big)=  
\mf B_{N,I,W}^{\Xi,E} \Big( \mc C_{  x}^{\flat }(h')\otimes  h\Big)\end{gather}
où, dans  le membre de droite, $x$ est considéré comme une forme  linéaire   sur  $U^{*,\theta}$ (invariante sous  l'action diagonale de  $\wh G$). 
\end{rem}

 \section{Morphismes de Frobenius  partiels, coalescence des pattes et opérateurs de Hecke}\label{frob-coalescence-cas part}

   Ce   chapitre est consacré à la démonstration de la  \propref{prop-coal-frob-cas-part}, qui affirme que les  
   opérateurs de Hecke en les  places $v$ de  $X\sm N$, qui sont des morphismes dans $D^{b}_{c}((X\sm (N\cup v))^{I}, E)$,  sont les restrictions de morphismes $S_{V,v}$  dans  \linebreak 
   $D^{b}_{c}((X\sm N)^{I}, E)$,  que nous allons définir comme la  composée    
 \begin{itemize}
 \item d'un  morphisme de création,
  \item de l'action d'un  morphisme de Frobenius partiel,
  \item d'un morphisme d'annihilation. 
  \end{itemize}
  La  \propref{prop-coal-frob-cas-part} servira dans le \lemref{S-non-ram}  
  pour montrer que les   opérateurs de Hecke en les places non ramifiées sont des cas particuliers d'opérateurs d'excursion  (ce qui justifiera que la  décomposition \eqref{intro1-dec-canonique} est compatible avec  l'isomorphisme de Satake en les  places non ramifiées). 
  Elle jouera aussi un rôle technique mais crucial dans le chapitre
   \ref{para-sous-faisceaux-Frob-partiels}, 
 grâce aux relations  d'Eichler-Shimura qui en seront déduites   dans le prochain chapitre. 

 \subsection{Enoncé des résultats}
 
  Soit  $I$   un ensemble fini et  $W$ une   représentation de dimension finie de 
$(\wh G)^{I}$.  Soit  $V$ une  représentation de dimension finie de  $\wh G$. 
 On considère  $W\boxtimes V\boxtimes V^{*}$  comme une    représentation de  $(\wh G)^{I\cup \{1,2\}}$.

 Soit   $v$ une   place de  $X\sm N$. 
     On considère  $v=\on{Spec}(k(v))$ comme un  sous-schéma de  $X$ et  on note  $E_{v}$ le faisceau constant sur  $v$. Soit $V$ une représentation irréductible  de  $\wh G$. La  fonction $h_{V,v}\in C_{c}(G(\mc O_v)\backslash G(F_{v})/G(\mc O_v), \mc O_{E})$ 
     désigne   la fonction sphérique associée à  $V$  par l'isomorphisme de  Satake, comme à la fin du premier chapitre.

On note $\Delta: X\to X\times X$ le   morphisme diagonal. On remarque que le point fermé   $\Delta(v)$ de  $X\times X$ est préservé  par 
 $\Frob_{X}^{\deg(v)}\times \mbf \Id$. 
 
 On a  des morphismes naturels de  $\wh G$-représentations   
 $$\mbf 1\xrightarrow{\delta_{V}} V\otimes V^{*}
 \text{ \ \ et  \ \ } V\otimes V^{*}\xrightarrow{\on{ev_{V}}}\mbf 1.$$ 
 
 \begin{rem} Les notations $\delta_{V}$ et $\on{ev_{V}}$, qui ont déjà été employées dans  l'introduction, viennent de   \cite{deligne-tens-fest}, à ceci près que nous nous autorisons à modifier arbitrairement l'ordre de   $V$ et $V^{*}$. Nous renvoyons  à la  \remref{contrainte-commutativite} pour   la contrainte de commutativité  {\it modifiée} en vigueur  pour les  faisceaux $\mc S_{I,W,E}^{(I_{1},...,I_{k})} $ et donc $\mc F_{N,I,W, \Xi,E}^{(I_{1},...,I_{k})}$. 
 \end{rem}     
 
  On a une équivalence de catégories 
 \begin{gather}\label{equiv-cat-Frob-v}
 D^{b}_{c}((X\sm N)^{I}\times v, E)^{F_{v}} \simeq D^{b}_{c}((X\sm N)^{I}, E)
 \end{gather}
 où la catégorie de gauche est formée des objets $\mc F$  munis d'un isomorphisme 
 $\theta:(\Id_{(X\sm N)^{I}} \times \Frob_{v})^{*}(\mc F)\simeq \mc F$ dont l'itérée $\deg(v)$ fois est
 l'identité. L'équivalence \eqref{equiv-cat-Frob-v} est une descente relativement à l'action de $\Z/\deg(v)\Z$ (donc elle est tautologique si $\deg(v)=1$). 
 Plus précisément elle
  associe (de la droite vers la gauche) à 
 un faisceau sur $(X\sm N)^{I}$ son image inverse par la première projection 
 $p_{1}: (X\sm N)^{I}\times v\to (X\sm N)^{I}$ et comme $p_{1}\circ (\Id_{(X\sm N)^{I}} \times \Frob_{v})=p_{1}$ cette image inverse est munie d'une donnée de descente et appartient donc à  $D^{b}_{c}((X\sm N)^{I}\times v, E)^{F_{v}}$. 
 Inversement, si $(\mc F,\theta) \in D^{b}_{c}((X\sm N)^{I}\times v, E)^{F_{v}}$, son image dans  $D^{b}_{c}((X\sm N)^{I}, E)$ est $\big(p_{1,*}(\mc F)\big)^{F_{v}}$ (c'est-à-dire que l'on prend les invariants par l'action de $\Z/\deg(v) \Z$ sur $p_{1,*}(\mc F)$ donnée par $\theta$).

      Pour  $\kappa$ assez grand  (en fonction de $\deg(v),V$), on définit 
        $S_{V,v}$ comme la  composée  
  \begin{gather}\label{def-SVv-text1}
 \mc H _{ N,I,W}^{\leq\mu,E} \boxtimes E_{v} \\ \label{def-SVv-text2}
  \xrightarrow{ \restr{\mc C_{  
    \delta_{V}}^{\sharp}}{(X\sm N)^{I}\times v}}
 \restr{ \mc H _{ N,I\cup\{1,2\},W\boxtimes V\boxtimes V^{*}}^{\leq\mu,E}}{(X\sm N)^{I}\times \Delta(v)} \\ \label{def-SVv-text3}
 \xrightarrow{ \restr{(F_{\{1\}})^{\deg(v)} }{(X\sm N)^{I}\times \Delta(v)}}
 \restr{ \mc H _{ N,I\cup\{1,2\},W\boxtimes V\boxtimes V^{*}}^{\leq\mu+\kappa,E}}{(X\sm N)^{I}\times \Delta(v)} \\ \label{def-SVv-text4}
  \xrightarrow{ \restr{\mc C_{  
    \on{ev_{V}}}^{\flat }}{(X\sm N)^{I}\times v}}
     \mc H _{ N,I,W}^{\leq\mu+\kappa,E}  \boxtimes E_{v} . 
 \end{gather}
  Le  morphisme  $S_{V,v}$ commute avec l'action naturelle du  morphisme de Frobenius partiel  sur  $E_{v}$ dans  \eqref{def-SVv-text1} et  \eqref{def-SVv-text4} parce que 
    \begin{itemize}
    \item  les morphismes de création et d'annihilation entrelacent  cette  action avec l'action de  $F_{\{1,2\}}$ sur \eqref{def-SVv-text2} et \eqref{def-SVv-text3}, 
        \item $F_{\{1\}}$, et donc   $F_{\{1\}}^{\deg(v)}$,  commutent  avec   $F_{\{1,2\}}$. 
        \end{itemize}
        Grâce à l'équivalence \eqref{equiv-cat-Frob-v},  $S_{V,v}$ se descend en un morphisme 
        $$S_{V,v}: \mc H _{ N,I,W}^{\leq\mu,E}\to 
         \mc H _{ N,I,W}^{\leq\mu+\kappa,E} $$  dans  $D^{b}_{c}((X\sm N)^{I},E)$.      
   
   A l'aide de la \propref{prop-coalescence-Hecke-Frob} et 
du \lemref{lem-lemme-4.3} on montre facilement que 
si $V$ et $V'$ sont des représentations de $\wh G$ et $v$ et $v'$ sont des places de $X\sm N$, 
$S_{V,v}$ et $ S_{V',v'}$ commutent et que si $v=v'$ leur composée est égale à $S_{V\otimes V',v}$.  
   
    \begin{prop}\label{prop-coal-frob-cas-part}
 Pour  $\kappa$ assez grand (en fonction de $\deg(v),V$), on a égalité dans 
 \begin{gather*}
 \Hom_{D^{b}_{c}((X\sm (N\cup v))^{I}, E)}\Big(\restr{ \mc H _{ N,I,W}^{\leq\mu,E}}{(X\sm (N\cup v))^{I}},
 \restr{ \mc H _{ N,I,W}^{\leq\mu+\kappa,E}}{(X\sm (N\cup v))^{I}} \Big)\end{gather*}  entre $T(h_{V,v})$ et la restriction de $S_{V,v}$ de $(X\sm N)^{I}$ à $(X\sm (N\cup v))^{I}$. 
 \end{prop}
 La proposition signifie que le morphisme  $S_{V,v}$ défini ci-dessus  prolonge à $(X\sm N)^{I}$ l'opérateur de Hecke $T(h_{V,v})$ qui était défini seulement 
 sur  
 $(X\sm (N\cup v))^{I}$ par la construction du paragraphe \ref{action-hecke-etale}.

  \begin{rem} La preuve du  résultat principal de cet article utilisera  seulement l'égalité entre les morphismes  induits de faisceaux sur $(X\sm (N\cup v))^{I}$ de $\restr{ \mc H _{ N,I,W}^{0,\leq\mu,E}}{(X\sm (N\cup v))^{I}}$ vers  $\restr{ \mc H _{ N,I,W}^{0,\leq\mu+\kappa,E}}{(X\sm (N\cup v))^{I}} $. 
 \end{rem}
La preuve de la \propref{prop-coal-frob-cas-part} lorsque  $V$ est  minuscule
 et  $\deg(v)=1$ a été esquissée dans  l'introduction. 
 La preuve du cas général sera donnée dans les paragraphes suivants. 
 
  \begin{rem}  Pour les variétés de Shimura sur les corps de nombres une  extension de  $T(h_{V,v})$ en   $v$ a été définie dans de nombreux cas, de fa\c con modulaire, par adhérence de Zariski, ou à l'aide de cycles proches  (voir  \cite{deligne-bki-mod,faltings-chai, genestier-tilouine}). 
    \end{rem}

   Sans aucun effort supplémentaire on peut étendre {\it tous} les   opérateurs de Hecke en des morphismes  dans  $D^{b}_{c}((X\sm N)^{I},E)$. 
      
   \begin{cor}    \label{cor-hecke-etendus-compo}
 Pour toute fonction $f\in C_{c}(K_{N}\backslash G(\mb A)/K_{N},E)$ 
 le  morphisme $T(f)$ défini dans le paragraphe  \ref{action-hecke-etale}
s'étend naturellement, pour $\kappa$ assez grand  en fonction de  $f$, en un morphisme  
 $$T(f)\in \Hom_{D^{b}_{c}((X\sm N)^{I}, E)} \big(\mc H _{ N,I,W}^{\leq\mu,E},  \mc H _{ N,I,W}^{\leq\mu+\kappa,E}\big) $$
dans  $D^{b}_{c}((X\sm N)^{I},E)$    
  de sorte que 
  \begin{itemize}
  \item $f\mapsto T(f)$ est $E$-linéaire, 
 \item
pour tout  $v\in |X|\sm N$, on a   $T(h_{V,v})=S_{V,v}$,  
   \item 
 ces extensions sont  compatibles avec la  composition: pour tout   $f_{1}, f_{2}\in C_{c}(K_{N}\backslash G(\mb A)/K_{N},E)$ et   $\kappa$ assez grand  en  fonction de  $f_{1}$ et  $f_{2}$, on a  
 $$T(f_{1}f_{2})=T(f_{1})T(f_{2})\text{  dans }\Hom_{D^{b}_{c}((X\sm N)^{I}, E)}\Big(\mc H _{ N,I,W}^{\leq\mu,E}, \mc H _{ N,I,W}^{\leq\mu+\kappa,E}\Big). $$
      \end{itemize}
  De plus ces    extensions commutent  avec l'action des morphismes de Frobenius partiels, et avec les  morphismes de création et d'annihilation. 
  \end{cor}
 \noindent{\bf Démonstration. } Comme $C_{c}(K_{N}\backslash G(\mb A)/K_{N},E)$ est le produit tensoriel restreint des $C_{c}(K_{N,v}\backslash G(F_{v})/K_{N,v},E)$, il suffit, pour toute place $v\in |X|$ et pour tout $f\in
C_{c}(K_{N,v}\backslash G(F_{v})/K_{N,v},E)$, 
d'étendre $$T(f)\in \Hom_{D^{b}_{c}((X\sm (N\cup v))^{I}, E)}\Big(\restr{\mc H _{ N,I,W}^{\leq\mu,E}}{(X\sm (N\cup v))^{I}}, \restr{\mc H _{ N,I,W}^{\leq\mu+\kappa,E}}{(X\sm (N\cup v))^{I}}\Big) $$ à 
$\Hom_{D^{b}_{c}((X\sm N)^{I}, E)}\Big(\mc H _{ N,I,W}^{\leq\mu,E}, \mc H _{ N,I,W}^{\leq\mu+\kappa,E}\Big)$. Il n'y a rien à faire si $v\in N$. Si $v\not\in N$ les fonctions $h_{V,v}$ forment une base de 
$C_{c}(K_{N,v}\backslash G(F_{v})/K_{N,v},E)$ et 
on  pose $T(h_{V,v})=S_{V,v}$. 
Il reste à montrer qu'étant donné  $v_{1},v_{2}\in |X|$ et $f_{i}\in C_{c}(K_{N,v_{i}}\backslash G(F_{v_{i}})/K_{N,v_{i}},E)$ pour $i=1,2$, les opérateurs $T(f_{1})$ et $T(f_{2})$ définis ainsi commutent si $v_{1} \neq v_{2}$ et que leur produit est $T(f_{1}f_{2})$ si $v_{1}=v_{2}$. 
Si $v_{1}$ et $v_{2}$ appartiennent à $|N|$ il n'y a rien de plus que les propriétés habituelles des opérateurs de Hecke. Si l'un appartient à $|N|$ et l'autre à $|X\sm N|$, on le sait déjà car les morphismes $S_{V,v}$ commutent aux morphismes $T(f)$ construits dans le paragraphe  \ref{action-hecke-etale}. Si $v_{1}$ et $v_{2}$ appartiennent à   $|X\sm N|$ on l'a déjà vu juste avant la \propref{prop-coal-frob-cas-part}.  \cqfd

   \subsection{La formule des traces de   Grothendieck-Lefschetz  par  les correspondances cohomologiques}
 \label{para-rappels}
 
 Ce paragraphe rappelle des résultats de   \cite{varshavsky-fujiwara,  brav-var}. 
 Soit  $Y$ un schéma de type fini sur un corps fini $k$ et $\Frob_{Y/k}$ le morphisme de Frobenius $k$-linéaire    (plus tard nous appliquerons ceci avec  $k=k(v)$). Soit $\mc F\in D^{b}_{c}(Y, E)$. On a un   isomorphisme canonique   $F_{\mc F}:\Frob_{Y/k}^{*}(\mc F)\rightarrow \mc F$.  On note  $E_{\on{Spec} k}$ le faisceau constant sur $\on{Spec} k$.
  On note  $\Delta : Y\to Y \times Y$ le   morphisme diagonal et   $\pi:Y\to \on{Spec} k$ le  morphisme évident. 
 On note  $K_{Y}=\pi^{!}(E_{\on{Spec} k})$ le faisceau dualisant de  $Y$ et  
 $\mb D\mc F=\on{RHom} (\mc F, K_{Y})\in D^{b}_{c}(Y, E)$ le dual de Verdier   de  $\mc F$. 
  On va étudier la composée des trois   correspondances cohomologiques suivantes. 
 
 \begin{itemize}
 \item On note $\mc C^{\Delta}_{\sharp}$ la  correspondance cohomologique  de 
 $(\on{Spec} k,E_{\on{Spec} k})$ vers $(Y\times Y,\mb D\mc F \boxtimes \mc F)$ supportée par la  correspondance  diagonale $$\on{Spec} k\xleftarrow{\pi} Y \xrightarrow{\Delta} 
 Y\times Y$$
  et donnée par le   morphisme tautologique 
 $$\alpha^{\sharp}:\pi^{*}(E_{\on{Spec} k})=E_{Y} \to \on{RHom}(\mc F,\mc F)=\Delta^{!}( \mb D\mc F\boxtimes \mc F).$$
 \item On note $\mc C_{\Id\times \Frob}$ la  correspondance cohomologique  de 
 $(Y\times Y,\mb D\mc F \boxtimes \mc F)$ vers lui-même  supportée par la  correspondance    ``image inverse''    
 $$ Y\times Y \xleftarrow{\Id_{Y} \times \Frob_{Y}} Y\times Y
 \xrightarrow{\Id_{Y}\times \Id_{Y}}  Y\times Y$$
 et donnée par le morphisme 
 $$ \Id_{\mb D\mc F}\boxtimes F_{\mc F}: \mb D\mc F \boxtimes \Frob_{Y}^{*}(\mc F)\to 
\mb D \mc F\boxtimes \mc F.$$
 \item  On note $\mc C^{\Delta}_{\flat}$ la correspondance cohomologique  de 
 $(Y\times Y, \mb D \mc F \boxtimes\mc F)$ vers $(\on{Spec} k,E_{\on{Spec} k})$ supportée par la  correspondance diagonale $$Y\times Y \xleftarrow{\Delta} Y \xrightarrow{\pi} \on{Spec} k
 $$
 et donnée par le  morphisme tautologique 
 $$\Delta^{*}( \mb D\mc F\boxtimes \mc F)=\mb D\mc F\otimes \mc F \to K_{Y}  =\pi^{!}(E_{\on{Spec} k}).$$
 \end{itemize}
 
 Sur l'ensemble  $Y(k)$ on possède la 
 fonction $f$ associé à $\mc F$ par le dictionnaire  ``faisceaux-fonctions'', 
 c'est-à-dire que pour 
 $y\in Y(k)$, $$f(y)=\on{Tr} (\Frob_{y}, \mc F_{y}):=
 \sum_{i\in \Z}(-1)^{i} \on{Tr} (\Frob_{y}, H^{i}(\mc F)_{y}).$$ 
 
 Dans la  proposition suivante on  considère  $Y(k)$ comme un schéma discret (égal à une réunion finie de copies de  $\on{Spec} k$ indexée par l'ensemble   $Y(k)$). 
On note  $p: Y(k)\to  \on{Spec} k$  le  morphisme évident. 

On définit  $\mc C_{f}$ comme la  correspondance cohomologique  supportée sur  $$ \on{Spec} k \xleftarrow{p}  Y(k) \xrightarrow{p}  \on{Spec} k
 $$ et donnée par le  morphisme 
 $$p^{*}(E_{\on{Spec} k})=E_{Y(k)} \to E_{Y(k)} =p^{*}(E_{\on{Spec} k})$$ de  multiplication par $f$.

La proposition suivante est  l'énoncé local sous-jacent à la formule des traces de   Grothendieck-Lefschetz, tel qu'il est formulé dans \cite{brav-var} et  \cite{varshavsky-fujiwara}. 
 
 \begin{prop}\label{trace-frob}
 (claim 1.5.5 de \cite{brav-var}, voir  aussi  \cite{varshavsky-fujiwara} 1.2.2, 1.5.7, 2.1.3, 2.2.4) 
 La composée  
 $\mc C^{\Delta}_{\flat}\circ 
 \mc C_{\Id \times \Frob} \circ 
 \mc C^{\Delta}_{\sharp}$ est égale à la 
 correspondance cohomologique $\mc C_{f}$ de $(\on{Spec} k,E_{\on{Spec} k})$ vers lui-même.  
 \end{prop}

   \subsection{Reformulation de la  \propref{prop-coal-frob-cas-part}}
 Pour raccourcir les formules on pose  
 $$\Check{X}=X\sm (N \cup v)$$
dans le reste de ce  chapitre (un autre avantage de cette  notation est qu'ici  il est clair que  $v$ est considéré comme un sous-schéma de $X$ alors que plus tard la lettre $v$ désignera souvent un  morphisme $S\to v\to X$, dont le graphe est strictement inclus dans   $v\times S\subset X\times S$ si  $\deg(v)>1$). 
 
 On peut supposer que  $V$ et $W$ sont  irréductibles. 
 Il suffit  de démontrer  le  lemme suivant.  
 \begin{lem}\label{lem-eq-coal-frob-cas-part}
 La restriction de 
  $S_{V,v}$ (définie comme la  composée \eqref{def-SVv-text1}-\eqref{def-SVv-text4}) à  $\Check{X}^{I}\times v$ est égale à 
 $T(h_{V,v})\boxtimes \Id_{E_{v}}$  dans  
 \begin{gather*}
 \Hom_{D^{b}_{c}(\Check{X}^{I}\times v, E)}
 \Big(\restr{ \mc H _{ N,I,W}^{\leq\mu,E}}{\Check{X}^{I}}
 \boxtimes E_{v},
 \restr{ \mc H _{ N,I,W}^{\leq\mu+\kappa,E}}{\Check{X}^{I}} \boxtimes E_{v}\Big). \end{gather*} 
\end{lem}

 On introduit quatre  correspondances cohomologiques: 
 \begin{itemize}
 \item [] a) $\mc C_{h_{V,v}}$ qui réalise l'opérateur de Hecke $T(h_{V,v})$, 
 \item    [] b) $\mc C_{\flat}$ qui réalise la restriction à   $\Check{X}^{I} \times v $ du  morphisme d'annihilation \eqref{def-SVv-text3}$\to$\eqref{def-SVv-text4}, 
 \item  [] c) $\mc C_{F}$ qui réalise la   restriction à  $\Check{X}^{I} \times v $ de l'action du morphisme de Frobenius partiel 
\eqref{def-SVv-text2}$\to$\eqref{def-SVv-text3}. 
\item  [] d) $\mc C_{\sharp}$ qui réalise la   restriction à  $\Check{X}^{I} \times v $ du  morphisme de création
\eqref{def-SVv-text1}$\to$\eqref{def-SVv-text2}. 
 \end{itemize}
 Le lemme \ref{lem-eq-coal-frob-cas-part} sera impliqué par le lemme \ref{lem-equiv-ST}, qui affirme que la  correspondance cohomologique $\mc C_{\flat}\circ \mc C_{F}\circ \mc C_{\sharp}$ (qui réalise donc la   restriction de $S_{V,v}$  à  $\Check{X}^{I} \times v $) est égale au  produit de $\mc C_{h_{V,v}}$ avec la  correspondance cohomologique identité  de $(v,E_{v})$ vers lui-même. 
 L'énoncé sans quotienter par $\Xi$ implique l'énoncé en quotientant par $\Xi$, donc on montrera le premier. Pour cette raison, dans la suite de la démonstration on quotientera par des restrictions à la Weil de $G$ et non de $G^{\mr{ad}}$. 
 
 \begin{rem}
 Dans l'esquisse de preuve donnée dans l'introduction (lorsque  $V$ est minuscule et $\deg(v)=1$), $\mc C_{\flat}$ était supportée par $\mc Y_{1}$ et  $\mc C_{F}\circ \mc C_{\sharp}$ par $\mc Y_{2}$ (et on pouvait calculer la composée sur des ouverts lisses où   ces  correspondances cohomologiques étaient déterminées par leurs supports). 
 \end{rem}

a) {\it Construction de $\mc C_{h_{V,v}}$ (réalisant  l'opérateur de Hecke $T(h_{V,v})$).} On pose 
\begin{gather} \mc Z^{(I)}=\restr{\Cht_{N,I,W }^{(I)}}
 {\Check{X}^{I}}\end{gather}
 et on notera   
  $\Gamma^{(I)}$ la correspondance de Hecke  de 
$ \mc Z^{(I)}$ vers lui-même. Le rapport avec la construction \ref{constr-corresp-hecke} est que $ \Gamma^{(I)}$ est la réunion des $\Gamma_{N}(g)$ pour $g$ parcourant $G(\mc O_{v})\backslash \mr{Orb}_{V}$ où  $\mr{Orb}_{V}\subset G(F_{v})/G(\mc O_{v})$ désigne la réunion finie des 
 $G(\mc O_{v})$-orbites associées aux   copoids dominants  de $G$ qui sont des poids de $V$. Alors 
 les  $S$-points du  champ  $\Gamma^{(I)}$
 classifient les données  de  $(x_i)_{i\in I}:S\to \Check{X}^{I}$ et d'un  diagramme 
\begin{gather} \label{texte-diag-Gamma}
 \xymatrix{
 (\mc G', \psi') \ar[r]^-{\phi'} & 
 (\ta{\mc G'}, \ta \psi')   \\
  (\mc G, \psi)  \ar[r]^-{\phi}  \ar[u]_-{\kappa}  &
 (\ta{\mc G}, \ta \psi) \ar[u]_-{\ta \kappa}
 } \end{gather}
 tel que   
 \begin{itemize}
 \item les lignes inférieures et supérieures appartiennent à   $\mc Z^{(I)}$, 
\item 
au-dessus de tout point géométrique $s$ de $S$, les  restrictions de $(\mc G\xrightarrow{\phi} \ta\mc G)$ et $(\mc G'\xrightarrow{\phi'} \ta\mc G')$ 
au voisinage fomel de $v$  dans  $X$ peuvent être trivialisées de fa\c con unique modulo l'action de    $G(\mc O_{v})$ et alors 
$\kappa:\restr{\mc G}{(X\sm v)\times s}\isom \restr{\mc G'}{(X\sm v)\times s}$ définit un élément de 
$G(\mc O_{v})\backslash G(F_{v})/G(\mc O_{v})$ qui appartient à   $G(\mc O_{v})\backslash \mr{Orb}_{V}$. 
\end{itemize}
Les projections 
 $\mr{pr}_{1}$,  {\it resp.} $\mr{pr}_{2}:  \Gamma^{(I)}\to \mc Z^{(I)}$ sont les morphismes  étales qui associent  au diagramme \eqref{texte-diag-Gamma} sa ligne 
 inférieure, {\it resp.} supérieure. 

 Pour toute fonction $f\in C_{c}(G(\mc O_v)\backslash G(F_{v})/G(\mc O_v), E)$  à support inclus dans $G(\mc O_{v})\backslash \mr{Orb}_{V}$, et en particulier pour   $f=h_{V,v}$, 
 l'opérateur de Hecke $T(f)$ est réalisé par la  correspondance cohomologique 
 $
 \mc C_{f} $ supportée par $\mc Z^{(I)} \xleftarrow{\mr{pr}_{2}} \Gamma^{(I)}\xrightarrow{\mr{pr}_{1}} \mc Z^{(I)}$
  et donnée  par la multiplication par $f$. Plus précisément on note  
  $\mc F^{(I)}$ 
    la  restriction 
 de $\mc F_{N,I,W,\Xi,E}^{(I)}$
 à l'ouvert $ \mc Z^{(I)}\subset \Cht_{N,I,W }^{(I)}$. 
 A un décalage près c'est le faisceau d'intersection de $ \mc Z^{(I)}$ et on a  un  isomorphisme canonique  $\mr{pr}_{1}^{*}(\mc F^{(I)})\simeq \mr{pr}_{2}^{*}(\mc F^{(I)})$ (que l'on comprend mieux en disant que les deux membres sont images inverses de $\mc S_{I,W,E}^{(I)}$). 
 Alors $
 \mc C_{f} $ est la  correspondance cohomologique  de 
 $(\mc Z^{(I)},\mc F^{(I)})$ vers lui-même définie comme 
le  morphisme de multiplication par $f$ de 
 $\mr{pr}_{1}^{*}(\mc F^{(I)})$ vers $\mr{pr}_{2}^{*}(\mc F^{(I)})=\mr{pr}_{2}^{!}(\mc F^{(I)})$. 
 
 b)  {\it Construction de $\mc C_{\flat}$ (réalisant la restriction à   $\Check{X}^{I} \times v $ du  morphisme d'annihilation $\eqref{def-SVv-text3} \to\eqref{def-SVv-text4}$).} On pose  
 \begin{gather*}\mc Z^{(\{1\},\{2\}, I)}=\restr{\Cht_{N,I \cup \{1,2\},W \boxtimes V\boxtimes V^{*}}^{(\{1\},\{2\},I)}}
 {\Check{X}^{I}\times \Delta(v)}.\end{gather*}
  C'est un  champ sur $\Check{X}^{I}\times v$ dont les points consistent en  la donnée de  $(x_i)_{i\in I}$ et d'un   diagramme  
 \begin{gather} \label{texte-diag-W}
 \xymatrix{
 & (\mc G_{1}, \psi_{1})\ar[d]^-{\phi_{2}}& 
  \\
 (\mc G_{0}, \psi_{0}) \ar[ru]^-{\phi_{1}} &
 (\mc G_{2}, \psi_{2}) \ar[r]^-{\phi_{3}}    &
 (\ta{\mc G_{0}}, \ta \psi_{0})
 } \end{gather}
 où $\phi_{1}$,  $\phi_{2}$ et $\phi_{3}$ sont  respectivement les  modifications associées à la patte  $1$, la patte  $2$ et aux pattes indexées  par  $I$.

 On note  $\mc Y_{\flat} \overset{i_{\flat}}{\hookrightarrow} 
 \mc Z^{(\{1\},\{2\}, I)}$ le  sous-champ fermé défini par la condition que dans  le diagramme  \eqref{texte-diag-W},  $\phi_{2}\phi_{1}: 
 \restr{\mc G_{0}}{(X-v)\times S}\to \restr{\mc G_{2}}{(X-v)\times S}$ s'étend en un   isomorphisme sur  $X\times S$. 
On a un   morphisme $$p_{\flat}: \mc Y_{\flat}\to \mc Z^{(I)}\times v 
 $$
dont la première composante   envoie    \begin{gather}\label{diag-Y-flat}
 \xymatrix{
 & (\mc G_{1}, \psi_{1})  \ar[d]^-{\phi_{2}}& 
  \\
 (\mc G_{0}, \psi_{0}) \ar[ru]^-{\phi_{1}} \ar[r]^-{\sim}&
 (\mc G_{2}, \psi_{2}) \ar[r]^-{\phi_{3}}    &
 (\ta{\mc G_{0}}, \ta \psi_{0}) 
 } \end{gather}
sur la ligne inférieure, c'est-à-dire   
 \begin{gather}\label{texte-ligne-bas}\big( (x_i)_{i\in I}, (\mc G_{0}, \psi_{0}) \xrightarrow{\phi_{3}(\phi_{2}\phi_{1})}    (\ta{\mc G_{0}}, \ta \psi_{0})
\big). \end{gather}
On pose  
 $$\mc F^{(\{1\},\{2\}, I)}=\restr{\mc F_{N,I \cup \{1,2\},W \boxtimes V\boxtimes V^{*} ,\Xi,E}^{(\{1\},\{2\},I)}}{\mc Z^{(\{1\},\{2\}, I)}}.$$   A un décalage près, $\mc F^{(\{1\},\{2\}, I)}$  est le faisceau d'intersection de $\mc Z^{(\{1\},\{2\}, I)}$.

On définit maintenant   une correspondance cohomologique  
$$\mc C_{\flat} \text{   \ 
de  \ }(\mc Z^{(\{1\},\{2\}, I)},\mc F^{(\{1\},\{2\}, I)})\text{ \  vers \  }(\mc Z^{(I)} \times v,\mc F^{(I)} \boxtimes E_{v})$$
supportée par la  correspondance 
$$\mc Z^{(\{1\},\{2\}, I)}
\overset{i_{\flat}}{\hookleftarrow}  
\mc Y_{\flat} 
\xrightarrow{p_{\flat}}
\mc Z^{(I)} \times v, 
  $$ et qui réalisera 
 la   restriction à  $\Check{X}^{I} \times v $ du morphisme d'annihilation     \eqref{def-SVv-text3}$\to$\eqref{def-SVv-text4}. 

On note $\mr{Gr}^{(\{1\},\{2\})}_{\Delta(v)}$ le schéma sur $v$ égal à 
la fibre de 
$\mr{Gr}_{\{1,2\},V\boxtimes V^{*}}^{(\{1\},\{2\})}$ sur $\Delta(v)$ et on définit de la même fa\c con  $\mr{Gr}^{(\{1,2\})}_{\Delta(v)}$. 
On rappelle que les  points de $\mr{Gr}^{(\{1\},\{2\})}_{\Delta(v)}$ consistent en  des diagrammes 
$(\mc G_{0} \xrightarrow{\phi_{1}}  \mc G_{1} \xrightarrow{\phi_{2}}
  \mc G_{2}\isor{\theta} \restr{G}{\Gamma_{\infty v}}\big)$
de $G$-torseurs sur le disque formel  $\Gamma_{\infty v}$ en $v$. 
Dans le  \thmref{thm-geom-satake} nous avons introduit le   faisceau pervers (à un décalage près) 
$\mc S_{\{1,2\},V\boxtimes V^{*},E}^{(\{1\},\{2\})}$ sur $\mr{Gr}_{\{1,2\},V\boxtimes V^{*}}^{(\{1\},\{2\})}$. Sa restriction à la fibre $\mr{Gr}^{(\{1\},\{2\})}_{\Delta(v)}$ en  $\Delta(v)$ 
n'est autre que  $IC_{\mr{Gr}^{(\{1\},\{2\})}_{\Delta(v)}}$ et nous préférons utiliser cette dernière notation car elle est plus simple. 
On note $\mf Y_{\flat}\overset{\mf i_{\flat}}{\hookrightarrow} \mr{Gr}^{(\{1\},\{2\})}_{\Delta(v)} $ le  sous-schéma fermé  défini  par la condition que   $\phi_{2}\phi_{1}$ est un isomorphisme. 
On a  
un carré cartésien 
 \begin{gather}\label{cart-square-Y-bemol-gr}\xymatrix{
 \mr{Gr}^{(\{1\},\{2\})}_{\Delta(v)} \ar[d]^-{\pi^{(\{1\},\{2\})}_{\{1,2\}}} &  \mf Y_{\flat} \ar[d]^-{\mf p_{\flat}}\ar@{_{(}->}[l]^-{\mf i_{\flat}}
  \\
 \mr{Gr}^{(\{1,2\})}_{\Delta(v)} &v \ar[l]^-{\nu_{\flat}}
 }\end{gather}
où 
$\pi^{(\{1\},\{2\})}_{\{1,2\}}$ est le  morphisme d'oubli et 
$\nu_{\flat}$ est  l'inclusion du point en  l'origine. 
Par le produit de fusion $\nu_{\flat,*}(E_{v})=\nu_{\flat,!}(E_{v})$ est un facteur direct de 
$(\pi^{(\{1\},\{2\})}_{\{1,2\}})_{!}\big(IC_{\mr{Gr}^{(\{1\},\{2\})}_{\Delta(v)}}\big)$. 
En effet, par le b) du \thmref{thm-geom-satake}, cette dernière expression est la restriction du faisceau pervers (à un décalage près) 
$\mc S_{\{1,2\},V\boxtimes V^{*},E}^{(\{1,2\})}$ sur $\mr{Gr}_{\{1,2\},V\boxtimes V^{*}}^{(\{1,2\})}$ à la fibre $\mr{Gr}^{(\{1,2\})}_{\Delta(v)}$ en $\Delta(v)$. Par le produit de fusion (c'est-à-dire le d) du \thmref{thm-geom-satake}) c'est aussi la restriction de $\mc S_{\{0\},V\otimes V^{*},E}^{(\{0\})}$ à la fibre de 
$\mr{Gr}_{\{0\},V\otimes V^{*}}^{(\{0\})}$ en  $v$. 
On a  donc un morphisme canonique 
$$ \varepsilon: (\pi^{(\{1\},\{2\})}_{\{1,2\}})_{!}\big(IC_{\mr{Gr}^{(\{1\},\{2\})}_{\Delta(v)}}\big)\to \nu_{\flat,*}(E_{v})$$ associé à $\on{ev}_{V}:V\otimes V^{*}\to \mbf 1$. 
Par adjonction il permet de définir  une correspondance cohomologique  
$\mc C_{\flat}^{\mr{Gr}}$ de $(\mr{Gr}^{(\{1\},\{2\})}_{\Delta(v)},IC_{\mr{Gr}^{(\{1\},\{2\})}_{\Delta(v)}})$ vers $(v,E_{v})$ supportée par 
$$\mr{Gr}^{(\{1\},\{2\})}_{\Delta(v)} \overset{\mf i_{\flat}}{\hookleftarrow} \mf Y_{\flat} \xrightarrow{\mf p_{\flat}} v.$$
En effet en appliquant le changement de base propre au  carré cartésien \eqref{cart-square-Y-bemol-gr} on obtient un  morphisme 
$$\mf p_{\flat,!}\mf i_{\flat}^{*}(IC_{\mr{Gr}^{(\{1\},\{2\})}_{\Delta(v)}})\simeq 
\nu_{\flat}^{*}(\pi^{(\{1\},\{2\})}_{\{1,2\}})_{!}(IC_{\mr{Gr}^{(\{1\},\{2\})}_{\Delta(v)}})\xrightarrow{\varepsilon} E_{v}.$$
Toute la construction ci-dessus était équivariante sous l'action  de $G_{nv}$ (la  restriction à la  Weil de $G$ de $nv$ vers le point) pourvu que  $n$ soit  assez grand  en fonction de $V$. 
On considère donc $\mc C_{\flat}^{\mr{Gr}}$ comme une correspondance cohomologique  
 de $(\mr{Gr}^{(\{1\},\{2\})}_{\Delta(v)}/G_{nv},IC_{\mr{Gr}^{(\{1\},\{2\})}_{\Delta(v)}})$ vers $(v/G_{nv},E_{v})$. 

On considère maintenant le  diagramme commutatif 
\begin{gather}\label{comm-diag-corresp-Y1}\xymatrix{
 \mc Z^{(\{1\},\{2\}, I)}  \ar[d]^-{}& 
\mc Y_{\flat}  \ar[l]_-{i_{\flat}} \ar[r]^-{p_{\flat}}  \ar[d]^-{} &\mc Z^{(I)} \times v \ar[d]^-{}
  \\
 \mr{Gr}^{(\{1\},\{2\})}_{\Delta(v)}/G_{nv} \times \mf K^{(I)}&  \mf Y_{\flat}/G_{nv}
\times \mf K^{(I)}\ar[l]_-{\mf i_{\flat}\times \Id}  \ar[r]^-{\mf p_{\flat}\times \Id}   & v/G_{nv}\times \mf K^{(I)}
 }\end{gather}
où 
$\mf K^{(I)}$ est un raccourci  pour $\restr{\Big(\mr{Gr}_{I,W}^{(I)}/G_{\sum n_{i}x_{i}}\Big)}{\Check{X}^{I}}$, et les entiers  $n_{i}$ sont assez grands. 
Les morphismes verticaux sont lisses par la  \propref{lissite-Cht-Grass}, car, d'après la remarque \ref{grassmann-disj-prod},  
$$\restr{\mr{Gr}_{\{1,2\}\cup I,W\boxtimes V\boxtimes V^{*}}^{(\{1\},\{2\},I)}}{\Check{X}^{I}\times \Delta(v)}= \mr{Gr}^{(\{1\},\{2\})}_{\Delta(v)}\times \restr{\mr{Gr}_{I,W}^{(I)}}{\Check{X}^{I}}.$$
 De plus les deux carrés sont cartésiens. 
On rappelle que la flèche verticale la plus à gauche  envoie  \eqref{texte-diag-W} sur le  produit de 
$$(\restr{\mc G_{0}}{\Gamma_{\infty v}} \xrightarrow{\phi_{1}} \restr{\mc G_{1}}{\Gamma_{\infty v}}  \xrightarrow{\phi_{2}}
  \restr{\mc G_{2}}{\Gamma_{\infty v}}) \text{ \ \ dans   \ \ } 
   \mr{Gr}^{(\{1\},\{2\})}_{\Delta(v)}/G_{nv}$$ (comme dans  la  \remref{rem-apres-action-Gnx-Grass}, c'est-à-dire que le 
    $G_{nv}$-torseur tautologique  est   $ \restr{\mc G_{2}}{\Gamma_{n v}}$)  et de 
$(\restr{\mc G_{2}}{\Gamma_{\sum n_{i}x_{i}}}
\xrightarrow{\phi_{3}}
\restr{\ta \mc G_{0}}{\Gamma_{\sum n_{i}x_{i}}})\in \mf K^{(I)}$ (de nouveau comme dans  la   \remref{rem-apres-action-Gnx-Grass} et nous ne le répèterons plus).  
Alors la correspondance cohomologique 
$\mc C_{\flat}$
est définie comme l'image inverse par le diagramme \eqref{comm-diag-corresp-Y1} du produit de 
$\mc C_{\flat}^{\mr{Gr}}$ et de la  correspondance cohomologique identité de $(\mf K^{(I)}, \mc S_{I,W,E}^{(I)})$ vers lui-même. 
En effet par la \defiref{defi-F-E-Cht}, 
$\mc F^{(\{1\},\{2\}, I)}$ est canoniquement isomorphe à l'image inverse  par la flèche verticale la plus à gauche  de $ IC_{\mr{Gr}^{(\{1\},\{2\})}_{\Delta(v)}} \boxtimes \mc S_{I,W,E}^{(I)}$,  et $ \mc F^{(I)}\boxtimes E_{v}$ est canoniquement isomorphe à l'image inverse  par la flèche verticale la plus à droite  de  $E_{v}\boxtimes \mc S_{I,W,E}^{(I)}$.

c)  {\it Construction de $\mc C_{F}$ (réalisant  la   restriction à  $\Check{X}^{I} \times v $ de l'action du morphisme de Frobenius partiel 
$\eqref{def-SVv-text2}\to\eqref{def-SVv-text3})$.} 
Quand on restreint  
$\Cht_{N,I \cup \{1,2\},W \boxtimes V\boxtimes V^{*}}^{(\{1\},\{2\},I)}$ à 
$ \Check{X}^{I}\times \Delta(v)$ seul importe l'ordre des  pattes $1$ et $2$ dans la   partition $(\{1\},\{2\},I)$. Quand on le restreint à  $\Check{X}^{I}\times (\Frob_{X}\times \Id_{X})^{i}\Delta(v)$, avec $i\in\{1,...,\deg(v)-1\}$ l'ordre de  $1$, $2$ et $I$ n'importe plus. 
On peut donc itérer  $\deg(v)$ fois le morphisme de Frobenius partiel en  la patte $1$ de la fa\c con suivante  
(dans  $D_{c}^{b}(\Check{X}^{I}\times v,E)$): 
\begin{gather}\label{iteration F1-degv}\restr{\Cht_{N,I \cup \{1,2\},W \boxtimes V\boxtimes V^{*}}^{(\{1\},\{2\},I)}}
 {\Check{X}^{I}\times \Delta(v)}\xrightarrow{\on  {Fr}_{\{1\}, N,I \cup \{1,2\}} ^{(\{1\},\{2\},I)}} \\ \nonumber 
 \restr{\Cht_{N,I \cup \{1,2\},W \boxtimes V\boxtimes V^{*}}^{(\{1\},\{2\},I)}}
 {\Check{X}^{I}\times (\Frob\times \Id)\Delta(v)}\xrightarrow{\on  {Fr}_{\{1\}, N,I \cup \{1,2\}} ^{(\{1\},\{2\},I)}} \cdots \\ \nonumber 
  \xrightarrow{\on  {Fr}_{\{1\}, N,I \cup \{1,2\}} ^{(\{1\},\{2\},I)}}
  \restr{\Cht_{N,I \cup \{1,2\},W \boxtimes V\boxtimes V^{*}}^{(\{1\},\{2\},I)}}
 {\Check{X}^{I}\times (\Frob\times \Id)^{\deg(v)-1}\Delta(v)} \\ \nonumber 
 \xrightarrow{\on  {Fr}_{\{1\}, N,I \cup \{1,2\}} ^{(\{1\},\{2\},I)}} 
 \restr{\Cht_{N,I \cup \{1,2\},W \boxtimes V\boxtimes V^{*}}^{(I,\{2\},\{1\})}}
 {\Check{X}^{I}\times \Delta(v)}
 \end{gather}
En copiant les  notations de  b) on pose 
\begin{gather*}\mc Z^{(I,\{2\},\{1\})}
=\restr{\Cht_{N,I \cup \{1,2\},W \boxtimes V\boxtimes V^{*}}^{(I,\{2\},\{1\})}}
 {\Check{X}^{I}\times \Delta(v)}
 \\ \text{ et   \ \ }\mc F^{(I,\{2\},\{1\})}
 =\restr{\mc F_{N,I \cup \{1,2\},W \boxtimes V\boxtimes V^{*} ,\Xi,E}^{(I,\{2\},\{1\})}}{\mc Z^{(I,\{2\},\{1\})}}. 
 \end{gather*}  
On note $\on {Fr}_{\{1\}}^{\deg(v)}: \mc Z^{(\{1\},\{2\},I)}\to \mc Z^{(I,\{2\},\{1\})}$ l'itérée \eqref{iteration F1-degv}. 
En itérant l'isomorphisme $\on F_{\{1\}, N,I \cup \{1,2\}} ^{(\{1\},\{2\},I)}$ (de la \propref{action-frob-partiels-IC}) on obtient  un isomorphisme 
$$\on {F}_{\{1\}}^{\deg(v)}:
 \big(\on {Fr}_{\{1\}}^{\deg(v)}\big)^{*}\Big( \mc F^{(I,\{2\},\{1\})}
\Big) \isom \mc F^{(\{1\},\{2\},I)}.$$

Alors $\mc C_{F}$ est définie comme la  correspondance cohomologique 
``image inverse'' de 
$(\mc Z^{(I,\{2\},\{1\})},\mc F^{(I,\{2\},\{1\})})$ vers 
$(\mc Z^{(\{1\},\{2\},I)},\mc F^{(\{1\},\{2\},I)})$
supportée par 
$$\mc Z^{(I,\{2\},\{1\})}\xleftarrow{\on {Fr}_{\{1\}}^{\deg(v)}} \mc Z^{(\{1\},\{2\},I)}=\mc Z^{(\{1\},\{2\},I)}$$ et donnée par l'isomorphisme $\on {F}_{\{1\}}^{\deg(v)}$. 

 d)  {\it Construction de $\mc C_{\sharp}$ (réalisant la   restriction à  $\Check{X}^{I} \times v $ du  morphisme de création
$\eqref{def-SVv-text1}\to\eqref{def-SVv-text2}$).} 
 A la permutation près des  pattes $1$ et $2$, $\mc C_{\sharp}$ n'est rien d'autre que le dual de Verdier  de $\mc C_{\flat}$, mais nous préférons en donner rapidement la construction. 
  
 On note 
$\mc Y_{\sharp} \overset{ i _{\sharp}}{\hookrightarrow} \mc Z^{(I,\{2\},\{1\})}$ le  sous-champ fermé défini par la condition que dans  le diagramme  
 \begin{gather}\label{diag-élément-Y-sharp}
 \xymatrix{
 (\mc G_{0}, \psi_{0}) \ar[r]^-{\phi_{1}} & 
 (\mc G_{1}, \psi_{1}) \ar[d]^-{\phi_{2}} & 
  (\ta{\mc G_{0}}, \ta \psi_{0}) \\
     &
 (\mc G_{2}, \psi_{2})\ar[ru]^-{\phi_{3}} &
 } \end{gather}
   $ \phi_{3}\phi_{2}: 
 \restr{\mc G_{1}}{(X-v)\times S}\to \restr{\ta{\mc G_{0}}}{(X-v)\times S}$ s'étend en un isomorphisme sur $X\times S$. 
 On a un   morphisme $$p_{\sharp}: \mc Y_{\sharp}\to \mc Z^{(I)}\times v
 $$ 
 dont la première composante  envoie  $$\xymatrix{
 (\mc G_{0}, \psi_{0}) \ar[r]^-{\phi_{1}} & 
 (\mc G_{1}, \psi_{1}) \ar[d]^-{\phi_{2}} \ar[r]^-{\sim}& 
 (\ta{\mc G_{0}}, \ta \psi_{0}) \\
     &
 (\mc G_{2}, \psi_{2})\ar[ru]^-{\phi_{3}} &
 }$$  vers la ligne supérieure, c'est-à-dire  
 \begin{gather}\label{texte-ligne-haut}\big( (x_i)_{i\in I}, (\mc G_{0}, \psi_{0}) \xrightarrow{(\phi_{3}\phi_{2})\phi_{1}}    (\ta{\mc G_{0}}, \ta \psi_{0})
\big). \end{gather}
On définit maintenant   une correspondance cohomologique  
$$\mc C_{\sharp} \text{   \ 
de  \ }(\mc Z^{(I)} \times v,\mc F^{(I)} \boxtimes E_{v})\text{ \  vers \  }(\mc Z^{(I,\{2\},\{1\})},\mc F^{(I,\{2\},\{1\})})$$
supportée par la  correspondance 
$$\mc Z^{(I)} \times v\xleftarrow{p_{\sharp}}\mc Y_{\sharp} \overset{i_{\sharp}}{\hookrightarrow} \mc Z^{(I,\{2\},\{1\})} ,  $$ et qui réalisera 
 la  restriction à   $\Check{X}^{I} \times v $ du   morphisme de création
 \eqref{def-SVv-text1}$\to$\eqref{def-SVv-text2}.

Comme dans  b), on associe à $\delta_{V}:\mbf 1\to V\otimes V^{*}$ une   
correspondance cohomologique 
$\mc C_{\sharp}^{\mr{Gr}}$ de $(v,E_{v})$ vers $(\mr{Gr}^{(\{2\},\{1\})}_{\Delta(v)},IC_{\mr{Gr}^{(\{2\},\{1\})}_{\Delta(v)}})$   supportée par 
$$v \xleftarrow{\mf p_{\sharp}} \mf Y_{\sharp}\overset{\mf i_{\sharp}}{\hookrightarrow}\mr{Gr}^{(\{2\},\{1\})}_{\Delta(v)}  $$
 où $\mf Y_{\sharp}$ est défini par le  carré cartésien 
  \begin{gather}\label{cart-square-Y-diese-gr}\xymatrix{
 \mr{Gr}^{(\{2\},\{1\})}_{\Delta(v)} \ar[d]^-{\pi^{(\{2\},\{1\})}_{\{1,2\}}} &  \mf Y_{\sharp} \ar[d]^-{\mf p_{\sharp}}\ar@{_{(}->}[l]^-{\mf i_{\sharp}}
  \\
 \mr{Gr}^{(\{1,2\})}_{\Delta(v)} &v \ar[l]^-{\nu_{\sharp}}
 }\end{gather}
Grâce au produit de fusion, on associe à $\delta_{V}$ un morphisme 
$$\nu_{\sharp,!}E_{v}\to (\pi^{(\{2\},\{1\})}_{\{1,2\}})_{*}IC_{ \mr{Gr}^{(\{2\},\{1\})}_{\Delta(v)}},$$ et donc  par adjonction et changement de base propre par le  carré cartésien  \eqref{cart-square-Y-diese-gr}, un morphisme 
$$E_{v}\to \nu_{\sharp}^{!}(\pi^{(\{2\},\{1\})}_{\{1,2\}})_{*}IC_{ \mr{Gr}^{(\{2\},\{1\})}_{\Delta(v)}}\simeq \mf p_{\sharp,*}\mf i_{\sharp}^{!}IC_{ \mr{Gr}^{(\{2\},\{1\})}_{\Delta(v)}}.$$
Cette correspondance cohomologique est équivariante par l'action   de $G_{nv}$. 
On considère maintenant   le diagramme commutatif 
\begin{gather}\label{comm-diag-corresp-Y2}\xymatrix{
\mc Z^{(I)} \times v \ar[d]^-{}& 
\mc Y_{\sharp}  \ar[l]_-{p_{\sharp}} \ar[r]^-{i_{\sharp}}  \ar[d]^-{} 
&\mc Z^{(I,\{2\},\{1\})}\ar[d]^-{}
  \\
 v/G_{nv}  \times \mf K^{(I)}&  \mf Y_{\sharp}/G_{nv} \times \mf K^{(I)}
\ar[l]_-{\mf p_{\sharp}\times \Id}  \ar[r]^-{\mf i_{\sharp}\times \Id}   & \mr{Gr}^{(\{2\},\{1\})}_{\Delta(v)}/G_{nv} \times \mf K^{(I)}
 }\end{gather}
 Les morphismes verticaux sont lisses d'après la  \propref{lissite-Cht-Grass}
 et les deux  carrés sont  cartésiens. 
On rappelle que la flèche verticale la plus à droite de 
\eqref{comm-diag-corresp-Y2} envoie  \eqref{diag-élément-Y-sharp} sur le  produit de 
$$(\restr{\mc G_{1}}{\Gamma_{\infty v}}
\xrightarrow{\phi_{2}}
\restr{\mc G_{2}}{\Gamma_{\infty v}}
\xrightarrow{ \phi_{3}}
\restr{\ta{\mc G_{0}}}{\Gamma_{\infty v}})
\text{ \ \  dans  \ \ }  
\mr{Gr}^{(\{2\},\{1\})}_{\Delta(v)}/G_{nv}
$$   et de 
$(\restr{\mc G_{0}}{\Gamma_{\sum n_{i}x_{i}}}
\xrightarrow{\phi_{1}}
\restr{\mc G_{1}}{\Gamma_{\sum n_{i}x_{i}}})\in \mf K^{(I)}$. 
Alors la  correspondance cohomologique 
$\mc C_{\sharp}$
est  définie comme l'image inverse  par le diagramme \eqref{comm-diag-corresp-Y2} du produit de  
$\mc C_{\sharp}^{\mr{Gr}}$ et de la  correspondance cohomologique identité de $(\mf K^{(I)}, \mc S_{I,W,E}^{(I)})$. 

\begin{rem} \label{rem-justif-norm-crea-annihil} {\it Justification des normalisations choisies pour $\mc C_{\flat}$  et $ \mc C_{\sharp}$. }
Il suffit de trouver une situation où apparaissent un opérateur de création et un opérateur d'annihilation du même type (avec des pattes crées ou annihilées apparaissant dans le même ordre) et où on sait calculer leur composée directement. 
   Voici une telle situation, où on part d'un champ de Hecke avec une patte indexée par $1$, puis on crée les pattes $2$ et $3$ par $\mathbf 1\to V^{*}\otimes V$, puis on annihile les pattes $1$ et $2$ par $V\otimes V^{*}\to \mathbf 1$. Plus précisément   on introduit le champ  de Hecke 
   $\Hecke_{\{1,2,3\}, V\boxtimes V^{*} \boxtimes V }^{\{1\},\{2\},\{3\} }$ muni 
   du faisceau image inverse de $\mc S  _{\{1,2,3\}, V\boxtimes V^{*} \boxtimes V,E}^{\{1\},\{2\},\{3\} }$ que l'on note encore par la même lettre. On introduit de même le champ  de Hecke 
   $\Hecke_{\{1 \}, V  }^{\{1\} }$ muni de $\mc S_{\{1 \}, V  }^{\{1\} }$
   (et aussi le même en rempla\c cant le singleton $\{1\}$ par $\{3\}$). On considère  alors 
        la composée 
   \begin{itemize}
   \item de la correspondance $\mc D_{\sharp}$ analogue à $\mc C_{\sharp}$ 
   $$\text{de   \ } (\Hecke_{\{1 \}, V  }^{\{1\} }, \mc S_{\{1 \}, V  }^{\{1\} })\text{ \  vers  \ } (\Hecke_{\{1,2,3\}, V\boxtimes V^{*} \boxtimes V }^{\{1\},\{2\},\{3\} },  \mc S  _{\{1,2,3\}, V\boxtimes V^{*} \boxtimes V,E}^{\{1\},\{2\},\{3\} })$$ créant  les pattes $2$ et $3$ par $\delta_{V}: \mathbf 1\to V^{*}\otimes V$, c'est-à-dire donnée par la correspondance
   $$\Hecke_{\{1 \}, V  }^{\{1\} } \leftarrow \mc X_{\sharp} \hookrightarrow 
   \Hecke_{\{1,2,3\}, V\boxtimes V^{*} \boxtimes V }^{\{1\},\{2\},\{3\}} $$
   où $\mc X_{\sharp}$ est formé des points 
   $$(\mc G_{0}\xrightarrow{\phi_{1}} \mc G_{1}\xrightarrow{\phi_{2}} \mc G_{2}\xrightarrow{\phi_{3}} \mc G_{3})\in 
   \Hecke_{\{1,2,3\}, V\boxtimes V^{*} \boxtimes V }^{\{1\},\{2\},\{3\}} $$
   tels que $\phi_{3}\phi_{2}$ soit un isomorphisme, 
   \item de la correspondance $\mc D_{\flat}$ analogue à $\mc C_{\flat}$ 
   $$\text{de   \ } (\Hecke_{\{1,2,3\}, V\boxtimes V^{*} \boxtimes V }^{\{1\},\{2\},\{3\} },  \mc S  _{\{1,2,3\}, V\boxtimes V^{*} \boxtimes V,E}^{\{1\},\{2\},\{3\} }) 
   \text{ \  vers  \ } 
   (\Hecke_{\{3 \}, V  }^{\{3\} }, \mc S_{\{3 \}, V  }^{\{3\} })
   $$ annihilant  les pattes $1$ et $2$ par $ \on{ev}_{V}: V\otimes V^{*}\to \mathbf 1$, c'est-à-dire donnée par la correspondance
   $$ \Hecke_{\{1,2,3\}, V\boxtimes V^{*} \boxtimes V }^{\{1\},\{2\},\{3\}}
   \hookleftarrow \mc X_{\flat} \rightarrow \Hecke_{\{3 \}, V  }^{\{3\} } 
   $$
   où $\mc X_{\flat}$ est formé des points 
   $$(\mc G_{0}\xrightarrow{\phi_{1}} \mc G_{1}\xrightarrow{\phi_{2}} \mc G_{2}\xrightarrow{\phi_{3}} \mc G_{3})\in 
   \Hecke_{\{1,2,3\}, V\boxtimes V^{*} \boxtimes V }^{\{1\},\{2\},\{3\}} $$
   tels que $\phi_{2}\phi_{1}$ soit un isomorphisme.    \end{itemize}
   On voit que sur les lieux lisses (correspondant aux $G(\mc O)$-orbites ouvertes dans les strates fermées de grassmanniennes affines associées à $V$ et $V^{*}$) la composition $\mc D_{\flat}\circ \mc D_{\sharp}$ est transverse et   égale à $\Id$. 
   D'un autre côté le lemme ``de Zorro'' dit que pour tout espace vectoriel V 
  \begin{gather}\label{zorro} \text{la composée  \ } V\xrightarrow{\Id_{V}\otimes \on{\delta}_{V}} 
    V\otimes V^{*} \otimes V\xrightarrow{ \on{ev}_{V} \otimes \Id_{V}} V \text{ \ \  est l'identité}\end{gather} (c'est un  lemme élémentaire dans la catégorie des espaces vectoriels et un des axiomes des catégories tannakiennes \cite{deligne-tens-fest}). La coïncidence des deux fa\c cons de calculer $\mc D_{\flat}\circ \mc D_{\sharp}$ (dont la première est limitée au lieu lisse) valide la normalisation choisie pour les correspondances $\mc C_{\sharp}$  et $\mc C_{\flat}$. 
\end{rem}
\medskip

On a donc terminé la  construction de $\mc C_{h_{V,v}}$, $\mc C_{\flat}$, $ \mc C_{F}$ et $ \mc C_{\sharp}$. 
 On voit facilement que   $\mc C_{\flat}\circ \mc C_{F}\circ \mc C_{\sharp}$ est supportée  par 
 le produit fibré   $\mc Y_{\flat}\times_{\mc Z^{(I,\{2\},\{1\})}}\mc Y_{\sharp}$, où les  morphismes sont  \begin{gather}\label{YbemolZ}\mc Y_{\flat}\xrightarrow{i_{\flat}} \mc Z^{(\{1\},\{2\},I)}\xrightarrow{\on {Fr}_{\{1\}}^{\deg(v)}}  \mc Z^{(I,\{2\},\{1\})}\end{gather}
  et $i_{\sharp}$. 
 \begin{lem}
 On a un   isomorphisme canonique  
\begin{gather}\label{Gamma-produit-fibre}
\Gamma^{(I)}\times v\simeq \mc Y_{\flat} \times_{\mc Z^{(I,\{2\},\{1\})}}\mc Y_{\sharp}. 
\end{gather}
\end{lem}
\dem Soit $S$ un  schéma sur $\Check{X}^{I}\times v$. On note $v$ le  $S$-point de $X$ donné  par $S\to v \to X$ et  on note 
$\tu{i} v$ son  image par  $\Frob_{X}^{i}$, pour $i\in \N$. 
Les $S$-points  $v$, $\ta v$, ..., $\tu{\deg(v)-1} v$ sont deux à deux disjoints    mais  $\tu{\deg(v)} v=v$. 
Soit 
 \begin{gather}\label{S-point-Y-flat}
 \xymatrix{
 & (\mc G_{1}, \psi_{1})  \ar[d]^-{\phi_{2}}& 
  \\
 (\mc G_{0}, \psi_{0}) \ar[ru]^-{\phi_{1}} \ar[r]^-{\sim}&
 (\mc G_{2}, \psi_{2}) \ar[r]^-{\phi_{3}}    &
 (\ta{\mc G_{0}}, \ta \psi_{0}) 
 } \end{gather}
un  $S$-point de $\mc Y_{\flat}$. On veut comprendre son image
 \begin{gather}\label{image-ZI12}
 \xymatrix{
 (\wt{\mc G}_{0}, \wt \psi_{0}) \ar[r]^-{\wt \phi_{1}} & 
 (\wt{\mc G}_{1}, \wt \psi_{1}) \ar[d]^-{\wt \phi_{2}} & 
  (\ta{\wt{\mc G}_{0}}, \ta {\wt \psi}_{0}) \\
     &
 (\wt{\mc G}_{2}, \wt \psi_{2})\ar[ru]^-{\wt \phi_{3}} &
 } \end{gather}
dans  $\mc Z^{(I,\{2\},\{1\})}$ par \eqref{YbemolZ}. 
 On remarque que  $\cup _{i=0}^{\deg(v)-1} \Gamma_{\infty \tu{i} v}\subset X\times S$ est la  complétion de $X\times S$ le long de  $v\times_{\on{Spec} \Fq} S$.  
Comme $\phi_{3}$ est une modification en les pattes indexées par $I$ et 
$\phi_{2}\phi_{1}$ est un isomorphisme sur $X \times S$, 
$\phi_{3}\phi_{2}\phi_{1}: \restr{\mc G_{0}}{\Gamma_{\infty \tu{i} v}}\to \restr{\ta \mc G_{0}}{\Gamma_{\infty \tu{i} v}}$ est   un isomorphisme  pour tout $i\in \N$. 

Par conséquent, pour tout $i\in \N$ on peut considérer $(\tu{i} \mc G_{0}\xrightarrow{\tu{i} \phi_{1}} \tu{i} \mc G_{1})$ comme une  modification de $\mc G_{0}$ (et aussi de $\ta \mc G_{0}$) en 
 $\tu{i} v$, et dans l'énoncé suivant on notera cette  modification  $A_{i}$. On remarque que  $A_{0}$ et $A_{\deg(v)}$ sont deux modifications en le même point $\tu{\deg(v)} v=v$  mais ne sont pas égales en général si on est parti d'un 
 $S$-point arbitraire  \eqref{S-point-Y-flat} de $\mc Y_{\flat}$.  
Avec ces notations, 
\begin{itemize}
\item $\wt{\mc G}_{0}$ est obtenu à partir de $\mc G_{0}$ par les modifications 
$A_{0}$, ..., $A_{\deg(v)-1}$, 
\item $\wt{\mc G}_{1}$ est obtenu  à partir de $\ta \mc G_{0}$ par les modifications 
$A_{0}$, ..., $A_{\deg(v)-1}$,  
\item $\wt{\mc G}_{2}$ est obtenu  à partir de  $\ta \mc G_{0}$ par les modifications 
$A_{1}$, ..., $A_{\deg(v)-1}$, 
\item $\ta \wt{\mc G}_{0}$ est obtenu  à partir de $\ta \mc G_{0}$ par les modifications 
$A_{1}$, ..., $A_{\deg(v)}$, 
\end{itemize}
et bien sûr les  structures de niveau en $N$ sont conservées.  
Par conséquent  
 \eqref{image-ZI12} appartient à  
$i_{\sharp}(\mc Y_{\sharp})\subset \mc Z^{(I,\{2\},\{1\})}$ 
si et seulement si  $A_{0}=A_{\deg(v)}$. 
 Si on note  $\kappa: \mc G_{0}\to 
\wt{\mc G}_{0}$ la  modification décrite  ci-dessus (et donnée par $A_{0}$, ..., $A_{\deg(v)-1}$), c'est exactement la condition pour que  le diagramme 
  \begin{gather*}
 \begin{CD} 
  (\wt{\mc G}_{0}, \wt \psi_{0}) 
  @> \wt \phi_{3}\wt\phi_{2}\wt\phi_{1}>> 
  (\ta{\wt{\mc G}_{0}}, \ta \wt \psi_{0})  \\
 @A\kappa AA 
 @A\ta \kappa AA \\
(\mc G_{0}, \psi_{0})
 @>\phi_{3}\phi_{2}\phi_{1}>>  
 (\ta{\mc G_{0}}, \ta \psi_{0}) 
 \end{CD}
 \end{gather*} soit commutatif  et alors il  appartient à  $ \Gamma^{(I)}\times v$
 (car on l'identifie avec   le diagramme \eqref{texte-diag-Gamma}). 
On a donc prouvé \eqref{Gamma-produit-fibre}. \cqfd

La correspondance cohomologique $\mc C_{\flat}\circ \mc C_{F}\circ \mc C_{\sharp}$  de  $(\mc Z^{(I)}\times v, \mc F^{(I)}\boxtimes E_{v})$ vers lui-même  est  donc 
supportée par $\Gamma^{(I)}\times v$. 
Comme les  projections $\pr_{1}$ et  $\pr_{2}$ de $\Gamma^{(I)}$ vers  $\mc Z^{(I)}$ sont étales et que $\mc F^{(I)}$ 
est  isomorphe à $IC_{\mc Z^{(I)}}$ à un décalage près, 
 $\mc C_{\flat}\circ \mc C_{F}\circ \mc C_{\sharp}$ est donné par la  multiplication par une  fonction  localement constante sur $\Gamma^{(I)}$
 (de fa\c con encore plus canonique $\pr_{1}^{*}(\mc F^{(I)})$ et $\pr_{1}^{*}(\mc F^{(I)})$ sont tous les deux  image inverse de $\mc S_{I,W,E}^{(I)}$). 
Le lemme suivant, qui  affirme qu'elle est égale à la   fonction   déduite de  $h_{V,v}$, implique donc le  \lemref{lem-eq-coal-frob-cas-part}.

\begin{lem}\label{lem-equiv-ST}
La correspondance cohomologique $\mc C_{\flat}\circ \mc C_{F}\circ \mc C_{\sharp}$  est égale au produit de $\mc C_{h_{V,v}}$ avec la correspondance cohomologique identité de $(v,E_{v})$ vers lui-même. 
\end{lem}
  
  On verra  que  le lemme précédent 
   résulte des deux lemmes suivants.

   \begin{lem}\label{lem-equiv-ST1}
Il existe une fonction $k_{V,v}\in C_{c}(G(\mc O_v)\backslash G(F_{v})/G(\mc O_v),E)$ (à support dans  $G(\mc O_v)\backslash \mr{Orb}_{V}$) qui dépend seulement de  $\on{Spec}(\mc O_{v})$, de la  restriction de $G$ à $\on{Spec}(\mc O_{v})$ et de $V$  (et en particulier est indépendante de $N$, $I$ et $W$) telle  que la  correspondance cohomologique $\mc C_{\flat}\circ \mc C_{F}\circ \mc C_{\sharp}$ soit  égale au produit de $\mc C_{k_{V,v}}$ avec  la correspondance cohomologique identité de $(v,E_{v})$ vers lui-même. 
\end{lem}
On peut paraphraser le 
 lemme précédent en disant que  l'opérateur $S_{V,v}$ est ``de nature locale en $v$''. C'est un énoncé extrêmement intuitif, dans la mesure où 
 le morphisme $S_{V,v}$ peut être défini au niveau des chtoucas locaux (une notion analogue à celle des groupes p-divisibles). 
  Dans la preuve on utilisera des  chtoucas restreints  (analogues aux   Barsotti-Tate tronqués) au lieu des  chtoucas locaux, pour des raisons techniques. En fait il s'agit du ``cas  non ramifié''  de \cite{genestier-lafforgue} où l'on montre que pour toute place $v$ (y compris dans $N$) les opérateurs d'excursions  $S_{I,W,x,\xi,(\gamma_{i})_{i\in I}}$ sont de  ``nature locale '' en $v$ dès lors que les  $\gamma_{i}$ appartiennent au groupe de Weil local en  $v$ (la preuve du  \lemref{S-non-ram} ci-dessous montre que    $S_{V,v}$ est un  opérateur d'excursion de ce type). 
  
  \begin{lem} \label{lem-equiv-ST2}
     Le lemme \ref{lem-equiv-ST} est vrai lorsque  $I=\emptyset $, $W=\mbf 1$ et $\deg(v)=1$. 
     \end{lem}

 Les lemmes  \ref{lem-equiv-ST1} et \ref{lem-equiv-ST2}
     impliquent  le  \lemref{lem-equiv-ST} car si on se donne  $\on{Spec}(\mc O_{v})$  on peut toujours trouver une courbe  $X$ sur $k(v)$ contenant un  point $v\in X(k(v))$ (de sorte que  
     $\deg(v)=1$), et étendre  $G$ de $\on{Spec}(\mc O_{v})$ à $X$ (ici c'est  automatique puisque $G$ est déployé), 
  et alors, par un argument classique de séries de Poincaré (rappelé dans \cite{genestier-lafforgue}),        une fonction 
     $k_{V,v}\in C_{c}(G(\mc O_v)\backslash G(F_{v})/G(\mc O_v), E)$ est caractérisée par son action sur 
     $C_{c}^{\mr{cusp}}(\Bun_{G,N}(\Fq)/\Xi ,E)$ lorsque   le  niveau  $N\subset X\sm v$  et $\Xi$ sont arbitraires. 	
     
    Il reste donc à montrer les lemmes \ref{lem-equiv-ST1} et \ref{lem-equiv-ST2}. 
     
     \subsection{Démonstration du \lemref{lem-equiv-ST1}}
  Pour établir la      nature locale  en $v$ du morphisme $S_{V,v}$, nous allons construire des champs classifiants de chtoucas restreints (c'est-à-dire  en gros des  chtoucas ``à coefficients'' sur un sous-schéma fini $nv\subset X$ pour $n$ assez grand). 
    Nous utilisons les chtoucas restreints plutôt que les  chtoucas locaux (c'est-à-dire  des chtoucas ``à coefficients'' sur $\infty v\subset X$) parce que  leurs champs classifiants sont des champs  d'Artin, et les  morphismes des champs  classifiant les   chtoucas globaux vers les  champs classifiant les  chtoucas restreints sont lisses. 
     
     On construira alors des analogues des  correspondances cohomologiques 
     $\mc C_{\flat}$,  $\mc C_{F}$, $\mc C_{\sharp}$ pour ces champs  de chtoucas restreints, 
    en utilisant la deuxième colonne du  diagramme  suivant (dont les trois premières lignes n'ont pas encore été définies).  
     Par rapport aux autres  diagrammes de ce  paragraphe, le diagramme  suivant ainsi que le diagramme \eqref{gros-diagramme-corresp} quelques pages plus loin sont tournés de  $\pi/2$ pour rentrer dans la  page.  
      Nous verrons plus tard que toutes les flèches horizontales (qui étaient donc verticales dans les diagrammes précédents) sont lisses et que  tous les carrés sont  cartésiens à part le carré  entre $a_{3}$ et $a_{4}$ qui est seulement commutatif.  
      Ainsi les  flèches horizontales  (lisses) situées entre le colonne de gauche et celle du milieu 
      entrelaceront les 
correspondances cohomologiques 
     $\mc C_{\flat}$,  $\mc C_{F}$, $\mc C_{\sharp}$ avec  leurs analogues pour les chtoucas restreints. 
    Voici   ce  diagramme, où tous les champs sont sur $\Check{X}^{I}\times v$: 
                    \begin{gather}\label{gros-diagramme}
 \xymatrix{
 \mc Z^{(I)} \times v    \ar[r]^-{a_{1}}& 
 \mr{Cht}^{\mr{res}(m,n)}_{\mbf 1} \times \mf K^{(I)}\ar[r]^-{b_{1}\times \Id}& v/G_{mv}\times \mf K^{(I)}
  \\
    \mc Y_{\flat}  \ar[d]^-{i_{\flat}}
   \ar[u]_-{p_{\flat}} \ar[r]^-{a_{2}}&  
   \mc Y_{\flat}^{\mr{res}(m,n)} \times \mf K^{(I)}
   \ar[u]^-{p_{\flat}^{\mr{res}(m,n)}\times \Id}  \ar[d]_-{i_{\flat}^{\mr{res}(m,n)}\times \Id} \ar[r]^-{b_{2}\times \Id}
   & \mf Y_{\flat}/G_{mv}
\times \mf K^{(I)}\ar[u]^-{\mf p_{\flat}\times \Id} \ar[d]_-{\mf i_{\flat}\times \Id} 
    \\
    \mc Z^{(\{1\},\{2\}, I)} \ar[d]^-{\on {Fr}_{\{1\}}^{\deg(v)}}
    \ar[r]^-{a_{3}} &  
    \mr{Cht}^{(\{1\},\{2\}),\mr{res}(m,n)}_{V\boxtimes V^{*}}\times \mf K^{(I)}
      \ar[d]_-{\on {Fr}_{\{1\}/k(v)}^{\mr{res}(m,n)}\times \Id} \ar[r]^-{b_{3}\times \Id}&  
    \mr{Gr}^{(\{1\},\{2\})}_{\Delta(v)}/G_{mv} \times \mf K^{(I)} 
      \\ 
   \mc Z^{(I,\{2\},\{1\})}   \ar[r]^-{a_{4}}& \mr{Gr}^{(\{2\},\{1\})}_{\Delta(v)}/G_{nv} \times \mf K^{(I)}
& 
   \\ 
  \mc Y_{\sharp} \ar[u]_-{i_{\sharp}} \ar[d]^-{p_{\sharp}} \ar[r]^-{a_{5}}
& \mf Y_{\sharp}/G_{nv} \times \mf K^{(I)} \ar[u]^-{\mf i_{\sharp}\times \Id} \ar[d]_-{\mf p_{\sharp}\times \Id}& 
  \\
  \mc Z^{(I)} \times v    \ar[r]^-{a_{6}}& v/G_{nv}  \times \mf K^{(I)} & 
   } \end{gather}
       On aurait pu écrire un  diagramme plus compliqué où les lignes $4,5,6$ auraient été similaires aux lignes $1,2,3$ du diagramme \eqref{gros-diagramme} (avec des  entiers $m,n$ plus petits que ceux des lignes $1,2,3$). En fait  nous pouvons  concevoir les 
     champs d'arrivée de  $a_{4}$, $a_{5}$ et $a_{6}$ dans le  diagramme  \eqref{gros-diagramme} comme des ``chtoucas restreints de niveau  $0$'' et nous considérons seulement ce cas particulier parce qu'il nous suffit.

 On commence par  définir le champ $ \mr{Cht}^{(\{1\},\{2\}),\mr{res}(m,n)}_{V\boxtimes V^{*}}$ sur $v$. C'est un champ classifiant de ``chtoucas restreints en  $v$'' c'est-à-dire en gros de   ``chtoucas à coefficients sur le  sous-schéma  $nv\subset X$''. Cependant  pour les définir on aura  besoin de choisir un entier subsidiaire 
  $m$ suffisamment grand par rapport à $n$, et de faire intervenir  des 
 $G$-torseurs sur $mv$.  
 
 \begin{rem}
 Les champs ${}^{r}\overline{\mc C}_{\emptyset}^{N}$ de \cite{laurent-inventiones} II.1.a
constituent  une excellente  définition des ``chtoucas à coefficients sur un sous-schéma  fini  $N$'' lorsque 
 $G=GL_{r}$ et $V=\mr{St}$, car 
 ils sont  vraiment   analogues aux  
 Barsotti-Tate tronqués 
 (pour les experts: une  action stricte au sens  de Faltings \cite{faltings-strict} est cachée dans leur  définition). 
 Nous ne savons pas comment  généraliser la construction de ces champs lorsque  $G$ et $V$ sont arbitraires.  \end{rem}

 Bien que nous ne les utilisions pas dans cet article, il est éclairant de commencer par définir les champs classifiants de   chtoucas locaux (qui ne sont pas des champs d'Artin). 
 Un chtouca local sur un schéma $S$ sur $v$ est la donnée   
   \begin{itemize}
   \item    d'un  morphisme \begin{gather}\label{kappa-locaux}\kappa:S\to  \mr{Gr}^{(\{1\},\{2\})}_{\Delta(v)}/G_{\infty v} \text {  donné par   } (\wh{\mc G}_{0} \xrightarrow{\wh \phi_{1}} \wh{\mc G}_{1}\xrightarrow{\wh \phi_{2}}
 \wh{\mc G}_{2})\end{gather}
       (où les  $\wh{\mc G}_{i}$ sont des $G$-torseurs sur $\Gamma_{\infty v}\subset X\times S$ et le  $   G_{\infty v}$-torseur tautologique sur le  quotient $\mr{Gr}^{(\{1\},\{2\})}_{\Delta(v)}/G_{\infty v}$ est $\wh{\mc G}_{2}$) 
        \item d'un isomorphisme \begin{gather}\label{theta-locaux}\theta: \wh{\mc G}_{2}\isom 
         (\Id_{\infty v}\times \Frob_{S/v})^{*}(\wh{\mc G}_{0}).\end{gather} 
        \end{itemize}
       On renvoie à \cite{fontaine-egale-car} pour une étude des  chtoucas locaux. 
      On note ici  $(\Id_{\infty v}\times \Frob_{S/v})^{*}(\wh{\mc G}_{0})$ au lieu  de $\ta \wh{\mc G}_{0}$ parce que le 
        morphisme  de  Frobenius est relatif à $k(v)$ et non à  $\Fq$. D'une fa\c con plus concise un  chtouca local peut être écrit sous la forme 
    \begin{gather}   \label{chtouca-local-concis}      (\wh{\mc G}_{0} \xrightarrow{\wh \phi_{1}} \wh{\mc G}_{1}\xrightarrow{\wh \phi_{2}}
(\Id_{\infty v}\times \Frob_{S/v})^{*}(\wh{\mc G}_{0})). \end{gather}

      Voici maintenant la définition des  chtoucas restreints. 
      Soit $n\geq 0$. Plus tard on demandera que $n$ soit assez grand pour que la condition \eqref{condition-n} soit satisfaite. 
      On a un   $G_{nv}$-torseur  naturel sur 
         $ \mr{Gr}^{(\{1\},\{2\})}_{\Delta(v)}$ dont la fibre en 
         \begin{gather}\label{élément-GR12Delta}(\wh{\mc G}_{0} \xrightarrow{\wh \phi_{1}} \wh{\mc G}_{1}\xrightarrow{\wh \phi_{2}}
 \wh{\mc G}_{2}\simeq \restr{G}{\Gamma_{\infty v}})\end{gather} (où 
 $\wh{\mc G}_{0}, \wh{\mc G}_{1}, \wh{\mc G}_{2}$ sont des $G$-torseurs sur $
\Gamma_{\infty v}$) est  $\restr{\wh{\mc G}_{0} }{\Gamma_{n v}}$. 
On note $\mr{Gr}^{(\{1\},\{2\})}_{\Delta(v), \mr{triv}_{n}}$ l'espace total de ce  $G_{nv}$-torseur, dont les  points consistent en  la donnée de \eqref{élément-GR12Delta} et d'une trivialisation de $\restr{\wh{\mc G}_{0} }{\Gamma_{n v}}$. 
 On choisit  $m\geq n$ assez grand (en fonction de $V$ et $n$) pour que l'action à droite naturelle  de 
         $G_{\infty v}$  sur $\mr{Gr}^{(\{1\},\{2\})}_{\Delta(v), \mr{triv}_{n}}$  (par changement  de la  trivialisation de $\wh{\mc G}_{2}$) se factorise à travers  $G_{mv}$.   
         On imposera plus tard la condition supplémentaire \eqref{condition-m} sur  $m$,  qui sera  satisfaite s'il est assez grand (en fonction de $n$). 
          Le quotient $ \mr{Gr}^{(\{1\},\{2\})}_{\Delta(v)}/G_{mv}$ est muni        \begin{itemize}
       \item d'un  $G_{mv}$-torseur  $\mc G_{2,mv}$ tautologiquement associé au quotient par $G_{mv}$
       \item d'un   $G_{nv}$-torseur   $\mc G_{0,nv}$
          dont l'espace total  est $\mr{Gr}^{(\{1\},\{2\})}_{\Delta(v), \mr{triv}_{n}}/G_{mv}$  (la   notation 
       $\restr{\wh{\mc G}_{0} }{\Gamma_{n v}}$ aurait constitué un abus car  $\wh{\mc G}_{0}$ existe sur $ \mr{Gr}^{(\{1\},\{2\})}_{\Delta(v)}$ et aussi sur $ \mr{Gr}^{(\{1\},\{2\})}_{\Delta(v)}/G_{\infty v}$ mais pas  sur $ \mr{Gr}^{(\{1\},\{2\})}_{\Delta(v)}/G_{mv}$).        \end{itemize}
       Alors, pour tout schéma $S$ sur $v$, 
     les  $S$-points de 
         $ \mr{Cht}^{(\{1\},\{2\}),\mr{res}(m,n)}_{V\boxtimes V^{*}}$ classifient la donnée         \begin{itemize}
         \item d'un  morphisme 
         \begin{gather}\label{kappa-def-Cht-res}\kappa:S\to   \mr{Gr}^{(\{1\},\{2\})}_{\Delta(v)}/G_{mv}\end{gather} de champs sur $v$, 
                \item d'une identification 
                \begin{gather}\label{theta-def-Cht-res}\theta :  \kappa^{*}(\restr{\mc G_{2,mv}}{nv})\isom  (\Id_{nv}\times \Frob_{S/v})^{*}\kappa^{*}(\mc G_{0,nv})\end{gather} de $G_{nv}$-torseurs sur $S$. 
                         \end{itemize}
         
        Le morphisme $$b_{3}: \mr{Cht}^{(\{1\},\{2\}),\mr{res}(m,n)}_{V\boxtimes V^{*}}\to  \mr{Gr}^{(\{1\},\{2\})}_{\Delta(v)}/G_{mv}$$ est simplement l'oubli de  $\theta$, donc il est lisse (c'est même un $G_{nv}$-bitorseur).

  Pour préparer la  construction de $a_{3}$ on commence par déclarer 
  que  le  chtouca local associé à un $S$-point \eqref{texte-diag-W} de $ \mc Z^{(\{1\},\{2\}, I)}$ est 
   \begin{gather}\label{cht-loc-associe}
   \restr{\mc G_{0}}{\Gamma_{\infty v}}
   \xrightarrow{\phi_{1}}
   \restr{\mc G_{1}}{\Gamma_{\infty v}}
\xrightarrow{
   \tu{\deg(v)-1}(\phi_{3}\phi_{2}\phi_{1})\cdots 
   \ta (\phi_{3}\phi_{2}\phi_{1}) 
   \phi_{3}\phi_{2}} 
   (\Id_{\infty v}\times \Frob_{S/v})^{*}(\restr{\mc G_{0}}{\Gamma_{\infty v}})
   \end{gather}
   où l'on a utilisé l'isomorphisme 
   $$  \restr{\tu{\deg(v)}\mc G_{0}}{\Gamma_{\infty v}} \simeq (\Id_{\infty v}\times \Frob_{S/v})^{*}(\restr{\mc G_{0}}{\Gamma_{\infty v}})  $$ et le fait que $$ \tu{\deg(v)-1}(\phi_{3}\phi_{2}\phi_{1})\cdots 
   \ta (\phi_{3}\phi_{2}\phi_{1}) 
   \phi_{3}:  \restr{\mc G_{2}}{\Gamma_{\infty v}}\to  \restr{\tu{\deg(v)}\mc G_{0}}{\Gamma_{\infty v}}
$$ est un isomorphisme. Pour construire $a_{3}$ on considère le chtouca restreint associé à ce chtouca local.

        Plus précisément la première composante  de $a_{3}$ est le  morphisme 
         $$   \mc Z^{(\{1\},\{2\}, I)}\to \mr{Cht}^{(\{1\},\{2\}),\mr{res}(m,n)}_{V\boxtimes V^{*}}$$
    qui     envoie   \eqref{texte-diag-W} vers le chtouca restreint associé au chtouca local \eqref{cht-loc-associe}, c'est-à-dire 
    \begin{itemize}
    \item le  morphisme  $\kappa:S\to \mr{Gr}^{(\{1\},\{2\})}_{\Delta(v)}/G_{mv}$ associé à $$(\restr{\mc G_{0}}{\Gamma_{\infty v}}
   \xrightarrow{\phi_{1}}
   \restr{\mc G_{1}}{\Gamma_{\infty v}} \xrightarrow{\phi_{2}}
   \restr{\mc G_{2}}{\Gamma_{\infty v}}), $$
   si bien que  $\kappa^{*}(\mc G_{0,nv})=\restr{\mc G_{0}}{\Gamma_{n v}}$ et $\kappa^{*}(\mc G_{2,mv})=\restr{\mc G_{2}}{\Gamma_{m v}}$
    \item et 
    $$\theta=\restr{\tu{\deg(v)-1}(\phi_{3}\phi_{2}\phi_{1})\cdots 
   \ta (\phi_{3}\phi_{2}\phi_{1}) \phi_{3}}{\Gamma_{nv}}: 
               \restr{\mc G_{2}}{\Gamma_{nv}}
     \isom 
      \restr{\tu{\deg(v)} \mc G_{0}}{\Gamma_{nv}}.$$ 
    \end{itemize}
        La deuxième composante  de $a_{3}$ envoie  
         \eqref{texte-diag-W} vers le  $S$-point de $\mf K^{(I)}$ associé à $(\mc G_{2}\xrightarrow{\phi_{3}} \ta \mc G_{0})$.

   On va montrer maintenant  la lissité  de $a_{3}$. On note  $\Bun_{G,[nv]}$ le champ sur $\Fq$ dont les  $S$-points classifient un   $G$-torseur $\mc G$ sur $X\times S$ et une trivialisation $\psi: \restr{\mc G}{nv\times S}\to \restr{G}{nv\times S}$ (autrement dit  $[nv]$ est considéré comme un diviseur sur  $X$). On a  un isomorphisme 
     $$\Bun_{G,[nv]}\times v=\Bun_{G,nv+n\ta v+\cdots + n\tu{\deg(v)-1}v}$$ de champs sur $v$  où, dans  le membre de droite, $v$ désigne le  morphisme $v\to X$.  
     L'action de $\Frob_{\Bun_{G,[nv]}}\times \Id_{v}$ sur le membre de gauche permute circulairement les   structures de niveau  en $v, \ta v, ..., \tu{\deg(v)-1}v$
     dans le membre de droite.

         Comme la  lissité de  $a_{3}$  est 
         une propriété locale pour la  topologie lisse sur le but, il suffit de montrer la  lissité sur $S$ de 
       \begin{gather}  \label{Z12-S}  \mc Z^{(\{1\},\{2\}, I)}\times_{    \mr{Cht}^{(\{1\},\{2\}),\mr{res}(m,n)}_{V\boxtimes V^{*}}\times \mf K^{(I)}} S \end{gather}
         où $S$ est un schéma sur $\Fq$ (ou en fait  sur $\Check{X}^{I}\times v$) et le  morphisme $S\to   \mr{Cht}^{(\{1\},\{2\}),\mr{res}(m,n)}_{V\boxtimes V^{*}}\times \mf K^{(I)}$ est donné par 
         \begin{itemize}
         \item [] a) un morphisme $\kappa$ comme dans  \eqref{kappa-def-Cht-res} provenant d'un morphisme $y:S\to \mr{Gr}^{(\{1\},\{2\})}_{\Delta(v)}$ 
                  \item [] b)  une identification $\theta$ comme dans  \eqref{theta-def-Cht-res} provenant d'une trivialisation $\lambda : y^{*}(\mc G_{0,nv})\isom \restr{G}{\Gamma_{nv}}$. 
                  \item [] c)  un morphisme $z:S\to \restr{\mr{Gr}_{I,W}^{(I)}}{\Check{X}^{I}}$. 
         \end{itemize}
         Comparée à  un $S$-point $(\kappa, \theta)$ de $\mr{Cht}^{(\{1\},\{2\}),\mr{res}(m,n)}_{V\boxtimes V^{*}}$ la donnée de a) et b) consiste en une  trivialisation $\lambda$
         de           $\mc G_{0,nv}$, puis un raffinement de  $nv$ à $mv$ de la  trivialisation $\theta\circ  (\Id_{nv}\times \Frob_{S/v})^{*}( \lambda)$ de $\restr{\mc G_{2,mv}}{nv}$. De plus 
         $\restr{\mr{Gr}_{I,W}^{(I)}}{\Check{X}^{I}}\to \mf K^{(I)}$ est évidemment lisse  puisque c'est un   $G_{\sum_{i\in I}n_{i}x_{i}}$-torseur. 
       On a donc justifié le fait que les  $S$ comme ci-dessus recouvrent  $ \mr{Cht}^{(\{1\},\{2\}),\mr{res}(m,n)}_{V\boxtimes V^{*}}\times \mf K^{(I)}$ pour la  topologie lisse. 
                  
               Il reste à montrer que  \eqref{Z12-S} est  lisse sur $S$. 
             Cela résulte du fait que    \eqref{Z12-S}  est l'égalisateur entre les deux morphismes suivants de champs
         sur $S$: 
      le   morphisme lisse d'oubli 
         \begin{gather*}\Bun_{G,mv+n\ta v+\cdots + n\tu{\deg(v)-1}v+\sum_{i\in I}n_{i}x_{i}}\times _{\big(\Check{X}^{I}\times v\big)}S \\ \to \Bun_{G,nv+n\ta v+\cdots + n\tu{\deg(v)-1}v}\times _{v}S =\Bun_{G,[nv]}\times_{\on{Spec} \Fq} S\end{gather*} 
         et la  composée 
          \begin{gather*} \Bun_{G,mv+n\ta v+\cdots + n\tu{\deg(v)-1}v+\sum_{i\in I}n_{i}x_{i}}\times _{\big(\Check{X}^{I}\times v\big)}S
           \\ \xrightarrow{a_{z}}  \Bun_{G,mv+n\ta v+\cdots + n\tu{\deg(v)-1}v}\times_{v}S
          \\ \xrightarrow{a_{(y,\lambda)}} \Bun_{G,nv+n\ta v+\cdots + n\tu{\deg(v)-1}v}\times _{v}S 
          =\Bun_{G,[nv]}\times_{\on{Spec} \Fq} S 
        \\  \xrightarrow{\Frob_{\Bun_{G,[nv]}}\times_{\on{Spec} \Fq} \Id_{S}}
          \Bun_{G,[nv]}\times_{\on{Spec} \Fq} S, 
          \end{gather*}  
               où 
       \begin{itemize}
         \item  $a_{z}$ 
         envoie   $(\mc G_{3},\psi_{3})$  vers  $(\mc G_{2},\psi_{2})$, où $\mc G_{2}$ est obtenu par la  modification de $\mc G_{3}$ en les  points $(x_{i})_{ i\in I}$ donnée par  $z$ et $\restr{\psi_{3}}{\Gamma_{\sum_{i\in I}n_{i}x_{i}}}$ (comme dans la preuve de la \propref{lissite-Cht-Grass}), c'est-à-dire  
         $$(\mc G_{2}\to \mc G_{3})=\big( \beta^{(I)}_{I,W,\underline{n}}\big)^{-1}
         \big(z, (\mc G_{3},\restr{\psi_{3}}{\Gamma_{\sum_{i\in I}n_{i}x_{i}}})\big)$$ où $ \beta^{(I)}_{I,W,\underline{n}}$
         a été défini dans \eqref{defi-beta-avant-Grass-lissite} et 
         $$\psi_{2}=\restr{\psi_{3}}{\Gamma_{mv+n\ta v+\cdots + n\tu{\deg(v)-1}v}}$$ (grâce au fait que les  $x_{i}$ restent disjoints  du  sous-schéma $v\subset X$), 
                  \item $a_{(y,\lambda)}$ envoie  $(\mc G_{2},\psi_{2})$ vers $(\mc G_{0},\psi_{0})$, où  $\mc G_{0}$  est obtenu par la modification de  $\mc G_{2}$ en $v$ donnée  par $y$ et $\restr{\psi_{2}}{\Gamma_{mv}}$ 
                  (de nouveau comme dans la preuve de la \propref{lissite-Cht-Grass}), c'est-à-dire  
                   $$(\mc G_{0}\to \mc G_{1}\to \mc G_{2})=\big( \beta^{(\{1\},\{2\})}_{\{1,2\},V\boxtimes V^{*},m}\big)^{-1}
         \big(y, (\mc G_{2},\restr{\psi_{2}}{\Gamma_{mv}})\big)$$         et       
              $$  \restr{\psi_{0}}{\Gamma_{nv}}=  \lambda \text{      et }  \restr{\psi_{0}}{\Gamma_{n\tu{i}v}}=   \restr{\psi_{2}}{\Gamma_{n\tu{i}v}} \text{        pour tout } i\in \{1,...,\deg(v)-1\}.$$ 
                       \end{itemize}
                       Donc  \eqref{Z12-S}  est lisse sur $S$ car c'est l'égalisateur 
                       entre un morphisme lisse
              entre deux champs lisses sur $S$ et un morphisme dont la dérivée relative à $S$ est nulle (à cause de $\Frob_{\Bun_{G,[nv]}}$). 
              Ceci termine la preuve de la lissité de $a_{3}$.

        Ensuite  $ \mc Y_{\flat}^{\mr{res}(m,n)}   
$ est le sous-champ fermé de $  \mr{Cht}^{(\{1\},\{2\}),\mr{res}(m,n)}_{V\boxtimes V^{*}}
$ défini par la condition que  $\phi_{2}\phi_{1}$ est un isomorphisme. Plus précisément il est tel que le  diagramme suivant (dont le  produit par $\mf K^{(I)}$ apparaît dans  \eqref{gros-diagramme})  soit cartésien: 
               $$    \xymatrix{
   \mc Y_{\flat}^{\mr{res}(m,n)}   
        \ar[d]_-{i_{\flat}^{\mr{res}(m,n)}} 
     \ar[r]^-{b_{2}}
   & \mf Y_{\flat}/G_{mv}
 \ar[d]_-{\mf i_{\flat}} 
    \\
      \mr{Cht}^{(\{1\},\{2\}),\mr{res}(m,n)}_{V\boxtimes V^{*}}
          \ar[r]^-{b_{3}}&  
    \mr{Gr}^{(\{1\},\{2\})}_{\Delta(v)}/G_{mv} }$$
Le champ $\mr{Cht}^{\mr{res}(m,n)}_{\mbf 1}$ classifie la donnée d'un  $G_{mv}$-torseur $\mc G$ et d'un isomorphisme $\theta : 
         (\Id_{nv}\times \Frob_{S/v})^{*}(\restr{\mc G}{nv})\isom \restr{\mc G}{nv}$ de $G_{nv}$-torseurs. Autrement dit  c'est le champ classifiant du groupe   $K_{m,n}$  défini  comme  l'image inverse  du groupe fini  $G(\mc O_{nv})$ par le  morphisme d'oubli      $G_{mv}\to G_{nv}$. L'inclusion naturelle de $K_{m,n}$   dans  $G_{mv}$ 
         fournit le morphisme lisse $b_{1}: v/K_{m,n}\to v/G_{mv}$. 
     Les lignes $4,5,6$ du diagramme \eqref{gros-diagramme} sont simplement obtenues en tournant le  diagramme \eqref{comm-diag-corresp-Y2}  de $\pi/2$. 
    
       Pour terminer la  construction du diagramme \eqref{gros-diagramme}, il ne reste plus qu'à définir         $$  \on {Fr}_{\{1\}/k(v)}^{\mr{res}(m,n)}:  \mr{Cht}^{(\{1\},\{2\}),\mr{res}(m,n)}_{V\boxtimes V^{*}}\to  \mr{Gr}^{(\{2\},\{1\})}_{\Delta(v)}/G_{nv}.  $$
      Pour comprendre l'idée, on commence par définir le  morphisme de Frobenius partiel analogue pour les   chtoucas locaux: il  envoie 
        \eqref{chtouca-local-concis}   
      vers   $$ ( \wh{\mc G}_{1}\xrightarrow{\wh \phi_{2}}
(\Id_{\infty v}\times \Frob_{S/v})^{*}(\wh{\mc G}_{0})
\xrightarrow{(\Id_{\infty v}\times \Frob_{S/v})^{*}(\wh \phi_{1})}
(\Id_{\infty v}\times \Frob_{S/v})^{*}( \wh{\mc G}_{1})). 
$$
Avec les  notations de \eqref{kappa-locaux}
et \eqref{theta-locaux} on peut le reformuler ainsi: on considère
$$( \wh{\mc G}_{1}\xrightarrow{\wh \phi_{2}}\wh{\mc G}_{2})\text{ et }
((\Id_{\infty v}\times \Frob_{S/v})^{*}(\wh{\mc G}_{0})
\xrightarrow{(\Id_{\infty v}\times \Frob_{S/v})^{*}(\wh \phi_{1})}
(\Id_{\infty v}\times \Frob_{S/v})^{*}( \wh{\mc G}_{1}))$$
qui sont tous les deux déterminés par $\kappa$, et on les recolle avec l'aide  de $\theta$. 

 Pour faire la même chose avec les  chtoucas restreints, on a  besoin de nouvelles notations.  On note $\mr{Gr}^{(\{1\})}_{v}$, {\it resp. } $\mr{Gr}^{(\{2\})}_{v}$ le  schéma sur $v$ égal à 
la fibre de 
$\mr{Gr}_{\{1\},V}^{(\{1\})}$, {\it resp. } $\mr{Gr}_{\{2\},V^{*}}^{(\{2\})}$,  sur $v$. On note $ \mr{Gr}^{(\{1\})}_{v, \mr{triv}_{n}}$ 
l'espace total du $G_{nv}$-torseur  $\mc G_{0,nv}$ sur 
         $\mr{Gr}_{\{1\},V}^{(\{1\})}$ dont la fibre en 
        $(\wh{\mc G}_{0} \xrightarrow{\wh \phi_{1}} \wh{\mc G}_{1}\simeq \restr{G}{\Gamma_{\infty v}})$  est  $\restr{\wh{\mc G}_{0} }{\Gamma_{n v}}$. 
        La condition sur $n$ mentionnée précédemment (et satisfaite s'il est assez grand) est que  
        \begin{gather}\label{condition-n}
   \text{       l'action de }G_{\infty v }\text{ sur }\mr{Gr}^{(\{2\})}_{v}\text{  se factorise par }G_{nv}.\end{gather} 
Alors $\mr{Gr}^{(\{2\},\{1\})}_{\Delta(v)}$ peut être identifié  au quotient de 
$\mr{Gr}^{(\{2\})}_{v}  
         \times_{v}   \mr{Gr}^{(\{1\})}_{v, \mr{triv}_{n}}$  par l'action diagonale de $G_{nv}$. On fixe  $r\geq n$ assez grand  en fonction de $n$ pour que 
         l'action naturelle de $G_{\infty v}$ sur $ \mr{Gr}^{(\{1\})}_{v, \mr{triv}_{n}}$ se factorise à travers $G_{r v}$.
Par conséquent 
         \begin{gather}\label{Gr12-recollement}\mr{Gr}^{(\{2\},\{1\})}_{\Delta(v)}/G_{rv}=\Big(\mr{Gr}^{(\{2\})}_{v}/G_{nv}\Big)
         \times _{\big(v/G_{nv}\big)} \mr{Gr}^{(\{1\})}_{v}/G_{rv}\end{gather} 
         où 
       le morphisme $\alpha_{2}:  \mr{Gr}^{(\{2\})}_{v}/G_{nv}\to v/G_{nv}$ est tautologique  et 
        le  morphisme $\alpha_{1}:  \mr{Gr}^{(\{1\})}_{v}/G_{rv}\to v/G_{nv}$ est donné  par 
         le $G_{nv}$-torseur $\mc G_{0,nv}$ sur $ \mr{Gr}^{(\{1\})}_{v}/G_{rv}$ 
         dont l'espace total  est $ \mr{Gr}^{(\{1\})}_{v, \mr{triv}_{n}}/G_{rv}$. 
     On note $\gamma: \mr{Gr}^{(\{1\},\{2\})}_{\Delta(v), \mr{triv}_{n}}/G_{\infty v}\to \mr{Gr}^{(\{1\})}_{v, \mr{triv}_{n}}/G_{rv}$  le morphisme naturel 
               qui provient du  morphisme évident 
         $$ \mr{Gr}^{(\{1\},\{2\})}_{\Delta(v)}/G_{\infty v}\to \mr{Gr}^{(\{1\})}_{v}/G_{rv}, 
     \ (\wh{\mc G}_{0} \xrightarrow{\wh \phi_{1}} \wh{\mc G}_{1}\xrightarrow{\wh \phi_{2}} \wh{\mc G}_{2})\mapsto (\wh{\mc G}_{0} \xrightarrow{\wh \phi_{1}} \wh{\mc G}_{1}). $$
    La condition sur   $m$ mentionnée précédemment (et satisfaite s'il est assez grand  en fonction de $r$, donc de $n$)  est que                    \begin{gather}\label{condition-m} \gamma \text{ se factorise à travers }\mr{Gr}^{(\{1\},\{2\})}_{\Delta(v), \mr{triv}_{n}}/G_{mv}. \end{gather}            
                  On obtient ainsi un morphisme $$\beta_{1}:  \mr{Gr}^{(\{1\},\{2\})}_{\Delta(v)}/G_{mv} \to \mr{Gr}^{(\{1\})}_{v}/G_{rv}  $$ et  une  identification 
     $\beta_{1}^{*}(\mc G_{0,nv})\isom   \mc G_{0,nv}$. 
     On introduit  aussi le  morphisme 
     $$\beta_{2}:  \mr{Gr}^{(\{1\},\{2\})}_{\Delta(v)}/G_{mv}\to \mr{Gr}^{(\{2\})}_{v}/G_{nv}$$ qui provient du  morphisme évident 
     $$\mr{Gr}^{(\{1\},\{2\})}_{\Delta(v)}\to \mr{Gr}^{(\{2\})}_{v}, 
     \ (\wh{\mc G}_{0} \xrightarrow{\wh \phi_{1}} \wh{\mc G}_{1}\xrightarrow{\wh \phi_{2}} \wh{\mc G}_{2}\simeq \restr{G}{\Gamma_{\infty v}})\mapsto (\wh{\mc G}_{1} \xrightarrow{\wh \phi_{2}} \wh{\mc G}_{2}\simeq \restr{G}{\Gamma_{\infty v}}).$$
     Alors $ \on {Fr}_{\{1\}/k(v)}^{\mr{res}(m,n)}$ est défini comme suit : 
     il envoie  un $S$-point de 
  $     \mr{Cht}^{(\{1\},\{2\}),\mr{res}(m,n)}_{V\boxtimes V^{*}}$ donné  par 
\begin{itemize}
\item  un $S$-point  $\kappa$ de $ \mr{Gr}^{(\{1\},\{2\})}_{\Delta(v)}/G_{mv}$ comme dans  \eqref{kappa-def-Cht-res}  
\item un isomorphisme $\theta$ 
  comme dans  \eqref{theta-def-Cht-res}  
  \end{itemize}
  vers le $S$-point de 
  $  \mr{Gr}^{(\{2\},\{1\})}_{\Delta(v)}/G_{nv}$ déduit du 
  $S$-point de 
  $  \mr{Gr}^{(\{2\},\{1\})}_{\Delta(v)}/G_{rv}$
  associé, avec l'aide  de 
  \eqref{Gr12-recollement},  
  \begin{itemize}
  \item au $S$-point de $\mr{Gr}^{(\{2\})}_{v}/G_{nv}$ égal à $\beta_{2}\circ  \kappa$, 
  \item  au $S$-point de $\mr{Gr}^{(\{1\})}_{v}/G_{rv}$ égal à $\beta_{1}\circ  \kappa \circ  \Frob_{S/k(v)}$
  \item à l'identification donnée  par $\theta$ de leurs  images respectives par $\alpha_{2}$ et $\alpha_{1}$ 
    comme  $S$-points de $v/G_{nv}$ (c'est-à-dire comme  $G_{nv}$-torseurs sur $S$).  
  \end{itemize}
          
                 Grâce à la lissité de  $b_{1},b_{2},b_{3}$ et au fait que le carré  entre $b_{1}$ et $b_{2}$ et le   carré  entre $b_{2}$ et $b_{3}$
                    sont  cartésiens, on définit une   correspondance cohomologique  
          $\mc C_{\flat}^{\mr{res}(m,n)}$ 
          $$\text{de \ }  ( \mr{Cht}^{(\{1\},\{2\}),\mr{res}(m,n)}_{V\boxtimes V^{*}}, 
   \mc F^{\mr{res}(m,n)}_{V\boxtimes V^{*}})     \text{ \ vers  \ }
                     (\mr{Cht}^{\mr{res}(m,n)}_{\mbf 1} ,
  \mc F^{\mr{res}(m,n)}_{\mbf 1})$$  supportée par 
  $$\mr{Cht}^{(\{1\},\{2\}),\mr{res}(m,n)}_{V\boxtimes V^{*}} 
  \xleftarrow{i_{\flat}^{\mr{res}(m,n)}}  \mc Y_{\flat}^{\mr{res}(m,n)}   \xrightarrow{p_{\flat}^{\mr{res}(m,n)}} 
  \mr{Cht}^{\mr{res}(m,n)}_{\mbf 1}$$
      comme l'image inverse  de $\mc C_{\flat}^{\mr{Gr}}$, où l'on pose 
  $$  \mc F^{\mr{res}(m,n)}_{V\boxtimes V^{*}}=b_{3}^{*}(IC_{\mr{Gr}^{(\{1\},\{2\})}_{\Delta(v)}})\text{ \  et \ }  
  \mc F^{\mr{res}(m,n)}_{\mbf 1}=b_{1}^{*}(E_{v})=E_{v}.$$ 
          
         On définit la correspondance cohomologique  $\mc C_{F}^{\mr{res}(m,n)}$ 
          $$\text{de \ }   (\mr{Gr}^{(\{2\},\{1\})}_{\Delta(v)}/G_{nv} ,IC_{\mr{Gr}^{(\{2\},\{1\})}_{\Delta(v)}})  \text{ \ vers \ }
                    ( \mr{Cht}^{(\{1\},\{2\}),\mr{res}(m,n)}_{V\boxtimes V^{*}}, 
   \mc F^{\mr{res}(m,n)}_{V\boxtimes V^{*}})   $$
         supportée par 
$$ \mr{Gr}^{(\{2\},\{1\})}_{\Delta(v)}/G_{nv} \xleftarrow{\on {Fr}_{\{1\}/k(v)}^{\mr{res}(m,n)}}
         \mr{Cht}^{(\{1\},\{2\}),\mr{res}(m,n)}_{V\boxtimes V^{*}}\xrightarrow{\Id}
         \mr{Cht}^{(\{1\},\{2\}),\mr{res}(m,n)}_{V\boxtimes V^{*}}$$
         comme le  morphisme 
         $\big(\on {Fr}_{\{1\}/k(v)}^{\mr{res}(m,n)}\big)^{*}\big(IC_{\mr{Gr}^{(\{2\},\{1\})}_{\Delta(v)}}\big)\to  \mc F^{\mr{res}(m,n)}_{V\boxtimes V^{*}}$
         (qui provient de l'énoncé analogue à  la \propref{prop-IC-produit-Fr-partiels} pour les chtoucas restreints).  
         
         Aux  morphismes lisses $a_{1},a_{3},a_{4},a_{6}$ on associe les   correspondances cohomologiques  ``images inverses''
               $\mc C^*_{a_{1}},\mc C^*_{a_{3}},\mc C^*_{a_{4}},\mc C^*_{a_{6}}$. 
         Par exemple $\mc C^*_{a_{3}}$ est la correspondance cohomologique  
         $$\text{de \ } (\mr{Cht}^{(\{1\},\{2\}),\mr{res}(m,n)}_{V\boxtimes V^{*}}\times \mf K^{(I)} , \mc F^{\mr{res}(m,n)}_{V\boxtimes V^{*}}\boxtimes  \mc S_{I,W,E}^{(I)}) \text{ \ vers \ } (\mc Z^{(\{1\},\{2\}, I)}, \mc F^{(\{1\},\{2\}, I)})$$ 
         supportée par 
       $$   \mr{Cht}^{(\{1\},\{2\}),\mr{res}(m,n)}_{V\boxtimes V^{*}}\times \mf K^{(I)} \xleftarrow{a_{3}} \mc Z^{(\{1\},\{2\}, I)}
         \xrightarrow{\Id} \mc Z^{(\{1\},\{2\}, I)}$$
         et donnée par l'isomorphisme 
         $$
  a_{3}^{*}( \mc F^{\mr{res}(m,n)}_{V\boxtimes V^{*}}\boxtimes  \mc S_{I,W,E}^{(I)})       \isom  \mc F^{(\{1\},\{2\}, I)}   $$
         provenant du fait que les deux  sont canoniquement isomorphes à 
           $(b_{3}a_{3})^{*}(IC_{\mr{Gr}^{(\{1\},\{2\})}_{\Delta(v)}} \boxtimes \mc S_{I,W,E}^{(I)})$.

          Finalement  on a  un diagramme commutatif de  correspondances cohomologiques: 
                         \begin{gather}\label{gros-diagramme-corresp}
 \xymatrix{
(\mc Z^{(I)} \times v,\mc F^{(I)} \boxtimes E_{v})   & 
 (\mr{Cht}^{\mr{res}(m,n)}_{\mbf 1} \times \mf K^{(I)},
  \mc F^{\mr{res}(m,n)}_{\mbf 1}\boxtimes  \mc S_{I,W,E}^{(I)})
  \ar[l]_-{\mc C^*_{a_{1}}}
    \\
  (\mc Z^{(\{1\},\{2\}, I)},\mc F^{(\{1\},\{2\}, I)})   \ar[u]^-{\mc C_{\flat}} 
   &  
   ( \mr{Cht}^{(\{1\},\{2\}),\mr{res}(m,n)}_{V\boxtimes V^{*}}\times \mf K^{(I)}, 
   \mc F^{\mr{res}(m,n)}_{V\boxtimes V^{*}}\boxtimes  \mc S_{I,W,E}^{(I)})
  \ar[u]^-{\mc C_{\flat}^{\mr{res}(m,n)}\times \Id} 
   \ar[l]_-{\mc C^*_{a_{3}}} 
    \\ 
 (\mc Z^{(I,\{2\},\{1\})},\mc F^{(I,\{2\},\{1\})})   
 \ar[u]^-{\mc C_{F}} 
 & (\mr{Gr}^{(\{2\},\{1\})}_{\Delta(v)}/G_{nv} \times \mf K^{(I)},IC_{\mr{Gr}^{(\{2\},\{1\})}_{\Delta(v)}}\boxtimes  \mc S_{I,W,E}^{(I)})
  \ar[l]_-{\mc C^*_{a_{4}}} 
  \ar[u]^-{\mc C_{F}^{\mr{res}(m,n)}\times \Id}
& 
    \\
 (\mc Z^{(I)} \times v,\mc F^{(I)} \boxtimes E_{v})   \ar[u]^-{\mc C_{\sharp}} & (v/G_{nv}  \times \mf K^{(I)} , E_{v}\boxtimes  \mc S_{I,W,E}^{(I)})  \ar[l]_-{\mc C^*_{a_{6}}} \ar[u]^-{\mc C_{\sharp}^{\mr{Gr}}\times \Id}
   } \end{gather}
       On rappelle que   $\mr{Orb}_{V}\subset G(F_{v})/G(\mc O_{v})$ est la réunion des  $G(\mc O_{v})$-orbites associées aux copoids dominants qui sont des   poids de $V$. 
         La correspondance cohomologique  $\mc C_{\flat}^{\mr{res}(m,n)}\circ 
         \mc C_{F}^{\mr{res}(m,n)}\circ 
         \mc C_{\sharp}^{\mr{Gr}}$ est supportée par 
         $$v/G_{nv} \leftarrow \mc Y_{\flat}^{\mr{res}(m,n)}
\times_{ \mr{Gr}^{(\{2\},\{1\})}_{\Delta(v)}/G_{nv}
}\mf Y_{\sharp}/G_{nv}=\mr{Orb}_{V}/K_{m,n}
\to \mr{Cht}^{\mr{res}(m,n)}_{\mbf 1}=v/K_{m,n} 
$$
 (où $K_{m,n}$ agit  sur $ \mr{Orb}_{V}$ à travers $G(\mc O_{nv})$). Elle est donc donnée par la    multiplication par une fonction $k_{V,v}\in C_{c}(G(\mc O_v)\backslash G(F_{v})/G(\mc O_v),E)$  à support dans  $
   \mr{Orb}_{V}/G(\mc O_v)$.  Cela termine la preuve  du 
   \lemref{lem-equiv-ST1}. \cqfd
                            
      \subsection{Démonstration du \lemref{lem-equiv-ST2}} 
      Nous sommes  ici dans le cas où 
      $I=\emptyset $, $W=\mbf 1$ et $\deg(v)=1$. 
      Nous commen\c cons par montrer le lemme à un signe et une puissance de $q^{1/2}$ près, puis nous vérifierons que la normalisation est correcte en nous  limitant aux lieux lisses, par un argument similiaire à la preuve du cas minuscule donnée dans l'introduction. 
      
     Adaptons rapidement les notations précédentes, avec en plus une  troncature  de Harder-Naramsimhan pour garantir que nous aurons affaire à des  schémas et non à des champs. Soit $\mu$ un copoids dominant arbitraire  de $G$ et $N\subset X\sm v$ un  niveau  assez grand  pour que  $\Bun_{G,N}^{\leq \mu}$  soit un  schéma. 
    En fait pour des raisons techniques nous demandons que  $\Bun_{G,N}^{\leq \mu+\kappa}$ soit un schéma (où $\kappa$, qui dépend seulement de  $V$, sera défini dans quelques lignes).  
       
       Pour tout singleton $J$  (qui sera $\{1\}$ pour la première copie de $\mr{H}_{N}$ et $\{2\}$ pour la  deuxième copie dans le  produit $\mr{H}_{N}\times \mr{H}_{N}$ qui apparaîtra dans le prochain diagramme) on définit  $\mr{H}_{N}$ comme  l'ouvert de $\restr{\Hecke_{N,J,V}^{(J)}}{v}$ formé des $(\mc G,\psi) \xrightarrow{\phi}  (\mc G',\psi')$ tels que $(\mc G,\psi)$ et $(\mc G',\psi')$ appartiennent à  $\Bun_{G,N}^{\leq \mu}$. 
C'est  un  schéma. Voici la  définition de $\kappa$ : il  est tel que 
pour tout point $(\mc G,\psi) \xrightarrow{\phi}  (\mc G',\psi')$ 
      de  $\restr{\Hecke_{N,J,V}^{(J)}}{v}$ vérifiant   
      $(\mc G',\psi')\in \Bun_{G,N}^{\leq \mu}$,  
      $(\mc G,\psi)$  appartient à  $\Bun_{G,N}^{\leq \mu+\kappa}$.  
      
On note 
$\mc Z$ et $\Gamma$ les   schémas discrets sur $\Fq$ 
qui sont des réunions de copies de $v=\on{Spec}(\Fq)$ indexées par les ensembles $\Bun_{G,N}^{\leq \mu}(\Fq)$ et 
  $\mr{H}_{N} (\Fq) $ (plus loin on écrira parfois 
  $\mc Z=\Bun_{G,N}^{\leq \mu}(\Fq)$ et $\Gamma=\mr{H}_{N} (\Fq) $).  On remarque que  $\Gamma$ est  la   correspondance de Hecke  entre $\mc Z$ et lui-même, qui est formée des  $((\mc E,\psi)\xrightarrow{\phi}  (\mc E',\psi'))$ avec  $(\mc E,\psi)$ et 
  $(\mc E',\psi')$ dans  $\mc Z$ et 
  $ \phi$ une  modification en  $v$ telle que l'élément associé 
  de $G(\mc O_{v})\backslash G(F_{v})/G(\mc O_{v})$  appartienne à   $G(\mc O_{v})\backslash \mr{Orb}_{V}$. On note $p:\mc Z\to v$ le  morphisme évident. 
  
  On définit    
   $   \mc Z^{(\{1\},\{2\})}$ comme  l'ouvert de $\restr{\Cht_{N, \{1,2\}, V\boxtimes V^{*}}^{(\{1\},\{2\})}}
 {\Delta(v)}$
formé des    diagrammes  
 \begin{gather}\label{diag-Cht-Ivide}(\mc G_{0}, \psi_{0})\xrightarrow{\phi_{1}} 
 (\mc G_{1}, \psi_{1})\xrightarrow{\phi_{2}} 
 (\ta{\mc G_{0}}, \ta \psi_{0})\end{gather} 
 tels que $(\mc G_{0}, \psi_{0})$ et $(\mc G_{1}, \psi_{1})$ appartiennent à  $\Bun_{G,N}^{\leq \mu}$. 
C'est   un schéma sur $v=\on{Spec}\Fq$. 
On note    $\mc Y_{\flat} \overset{i_{\flat}}{\hookrightarrow} 
 \mc Z^{(\{1\},\{2\})}$  le  sous-schéma fermé  défini par la condition que    $\phi_{2}\phi_{1}$ s'étend en  un  isomorphisme sur  $X\times S$. 
 On définit de fa\c con analogue  $\mc Z^{(\{2\},\{1\})}$ et $\mc Y_{\sharp} $. Le  morphisme de Frobenius partiel $\on {F}_{\{1\}}: \mc Z^{(\{1\},\{2\})}\to \mc Z^{(\{2\},\{1\})}$ envoie  \eqref{diag-Cht-Ivide} vers 
 $$(\mc G_{1}, \psi_{1})\xrightarrow{\phi_{2}} 
(\ta{\mc G_{0}}, \ta \psi_{0})\xrightarrow{\ta \phi_{1}} 
 (\ta{\mc G_{1}}, \ta \psi_{1}).$$ On rappelle que les  
  faisceaux d'intersection  $\mc F^{(\{1\},\{2\})}$ et $\mc F^{(\{2\},\{1\})}$ sont dotés des    correspondances cohomologiques 
 \begin{itemize}
  \item $\mc C_{\sharp}$  de  $(v,E_{v})$ vers     $(\mc Z^{(\{2\},\{1\})},\mc F^{(\{2\},\{1\})})$   supportée sur 
 $$v \xleftarrow{p_{\sharp}} \mc Y_{\sharp} \xrightarrow{i_{\sharp}} 
   \mc Z^{(\{2\},\{1\})} .$$  
  \item $\mc C_{F}$ de $(\mc Z^{(\{2\},\{1\})},\mc F^{(\{2\},\{1\})})$ vers $(\mc Z^{(\{1\},\{2\})},\mc F^{(\{1\},\{2\})})$ supportée sur 
 $$ \mc Z^{(\{2\},\{1\})} \xleftarrow{\on {F}_{\{1\}}}  \mc Z^{(\{1\},\{2\})} \xrightarrow{\Id} \mc Z^{(\{1\},\{2\})} ,$$ 
  \item $\mc C_{\flat}$ de $(\mc Z^{(\{1\},\{2\})},\mc F^{(\{1\},\{2\})})$ vers     $(v,E_{v})$  supportée par 
 $$\mc Z^{(\{1\},\{2\})} \xleftarrow{i_{\flat}} \mc Y_{\flat}  \xrightarrow{p_{\flat}} v , $$
\end{itemize}
 Leur composée $\mc C_{\flat}\circ\mc C_{F} \circ \mc C_{\sharp}$
      est supportée par $ \mc Y_{\sharp}\times_{\mc Z^{(\{2\},\{1\})}}   \mc Y_{\flat} =\Gamma$ et nous voulons montrer que c'est la  multiplication par $h_{V,v}$
      (à une puissance de $q^{1/2}$ et un signe près pour le moment). 
     Ces trois  correspondances cohomologiques  sont récapitulées par la première ligne du  diagramme suivant. D'autre part la   composée des trois  correspondances cohomologiques  de la  \propref{trace-frob}, appliquée au   schéma $\mr{H}_{N}$ sur $\Fq$ et  à son faisceau d'intersection,   est récapitulée par la dernière ligne de ce diagramme et d'après la  \propref{trace-frob} nous savons qu'elle est supportée par  $\Gamma=\mr{H}_{N} (\Fq) $ et donnée par la  multiplication par $h_{V,v}$ (à une puissance de $q^{1/2}$ et un signe près pour le moment). Nous allons utiliser ce  diagramme pour montrer que les  deux composées sont les mêmes (ou pour être plus précis  qu'elles sont données par la  multiplication par la même  fonction sur $\mr{H}_{N} (\Fq) $):  
             \begin{gather}
       \xymatrix{
\mc Z \ar@{=}[d]^-{} 
& \mc Y_{\sharp}  \ar@{=}[d]^-{}  \ar[l]_-{p_{\sharp}} 
\ar@{^{(}->}[r]^-{i_{\sharp}} 
 &  \mc Z^{(\{2\},\{1\})} \ar@{_{(}->}[d]^-{\alpha^{2,1}} 
 && \mc Z^{(\{1\},\{2\})} \ar[ll]_-{\on {F}_{\{1\}}}  \ar@{_{(}->}[d]^-{\alpha^{1,2}}
 &\mc Y_{\flat} \ar@{_{(}->}[l]_-{i_{\flat}} \ar[r]^-{p_{\flat}} \ar@{_{(}->}[d]^-{j_{\flat}}
 &\mc Z \ar[d]^-{p}
  \\
\mc Z \ar[d]^-{p} 
&  \mc Y_{\sharp}  \ar@{_{(}->}[d]^-{j_{\sharp}}  \ar[l]_-{p_{\sharp}} 
\ar@{^{(}->}[r]^-{i_{\sharp}^{\mc T}} 
& \mc T^{2,1}\ar@{_{(}->}[d]^-{\beta^{2,1}}
&& \mc T^{1,2}  \ar[ll]_-{\on {F}_{\{1\}}^{\mc T}}  \ar@{_{(}->}[d]^-{\beta^{1,2}} 
& \mr{H}_{N} \ar@{_{(}->}[l]_-{\Delta^{\mc T}} \ar[r]^-{p_{N}} \ar@{=}[d]^-{}  
&v  \ar@{=}[d]^-{}
\\
v 
& \mr{H}_{N} \ar[l]_-{p_{N}} \ar@{^{(}->}[r]^-{\Delta} 
&  \mr{H}_{N} \times  \mr{H}_{N} 
&& \mr{H}_{N} \times  \mr{H}_{N}  \ar[ll]_-{\Frob \times \Id}  
& \mr{H}_{N}  \ar@{_{(}->}[l]_-{\Delta} \ar[r]^-{p_{N}} 
&v 
 }\end{gather}
 On remarque  que $\mc Z^{(\{1\},\{2\})}$ s'identifie au   sous-schéma fermé de $\mr{H}_{N} \times  \mr{H}_{N} $ formé des  points 
\begin{gather}\label{point-HN}
\big( ((\mc G_{1},\psi_{1})\xrightarrow{\phi_{1}} (\mc G_{1}',\psi_{1}')), 
((\mc G_{2},\psi_{2})\xrightarrow{\phi_{2}} (\mc G_{2}',\psi_{2}')) \big)\end{gather}
tels que  \begin{gather}\label{cond-1-Z12}(\mc G_{1}',\psi_{1}')=(\mc G_{2}',\psi_{2}')\\ \label{cond-2-Z12} \text{ et \ \ }  
(\mc G_{2},\psi_{2})=(\ta \mc G_{1},\ta \psi_{1}) \end{gather} 
puisqu'ils fournissent alors le diagramme habituel 
$$(\mc G_{1},\psi_{1})\xrightarrow{\phi_{1}} (\mc G_{1}',\psi_{1}')\simeq (\mc G_{2}',\psi_{2}')\xrightarrow{\phi_{2}^{-1}}(\mc G_{2},\psi_{2})\simeq (\ta \mc G_{1},\ta \psi_{1}) .$$
On définit 
$\mc T^{1,2}$ comme le  sous-schéma fermé  de $\mr{H}_{N} \times  \mr{H}_{N} $ formé des  points \eqref{point-HN} tels que la condition \eqref{cond-1-Z12} soit satisfaite. On peut remarquer si on veut que  $\mc T^{1,2}$ est égal à un ouvert de la fibre de $\Hecke_{\{1,2\},V\boxtimes V^{*}}^{(\{1\},\{2\})}$ en $\Delta(v)$ mais cela  ne servira   pas.

De même $\mc Z^{(\{2\},\{1\})}$ s'identifie au   sous-schéma fermé de $\mr{H}_{N} \times  \mr{H}_{N} $ formé des  points 
\eqref{point-HN} tels que \begin{gather}\label{cond-1-Z21}(\mc G_{1}',\psi_{1}')=(\ta\mc G_{2}',\ta\psi_{2}')\\ \label{cond-2-Z21} \text{ et \ \ }  
(\mc G_{2},\psi_{2})=( \mc G_{1}, \psi_{1})\end{gather} 
puisqu'ils fournissent alors le diagramme habituel 
$$ (\mc G_{2}',\psi_{2}')\xrightarrow{\phi_{2}^{-1}}(\mc G_{2},\psi_{2}) \simeq (\mc G_{1},\psi_{1})\xrightarrow{\phi_{1}} (\mc G_{1}',\psi_{1}')\simeq (\ta\mc G_{2}',\ta\psi_{2}') .$$
On définit 
$\mc T^{2,1}$ comme le  sous-schéma fermé  de $\mr{H}_{N} \times  \mr{H}_{N} $ formé des  points \eqref{point-HN} tels que la  condition \eqref{cond-1-Z21} soit satisfaite. On définit $\on {F}_{\{1\}}^{\mc T}$ comme la restriction à $\mc T^{1,2}$ de $\mr{H}_{N} \times  \mr{H}_{N}  \xrightarrow{\Frob \times \Id}\mr{H}_{N} \times  \mr{H}_{N} $. 

Dans  le diagramme ci-dessus le  carré entre $j_{\sharp}$ et $\beta^{2,1}$ et le  carré entre $\alpha^{1,2}$ et $j_{\flat}$ sont cartésiens. 

 On définit $\mc F_{N}$ comme  le faisceau d'intersection de $\mr{H}_{N}$. 
 On note $\mr{Gr}_{V,v}$  le  schéma sur $v$ égal à 
la fibre de 
$\mr{Gr}_{J,V}^{(J)}$  sur $v$, où $J$ est un singleton (le  même  schéma avait été noté  $\mr{Gr}^{(\{1\})}_{v}$ dans le paragraphe précédent). 
 
On a  un morphisme évident $\gamma: \mr{H}_{N}\to \mr{Gr}_{V,v}/G_{nv}$ et  
 $\mc F_{N}$ est canoniquement isomorphe à $\gamma^{*}(IC_{\mr{Gr}_{V,v}})[d](d/2)$ où $d=\dim \Bun_{G,N}$. Bien sûr $\mathbb D{\mc F}_{N}\simeq \mc F_{N}$ mais il est préférable de garder la  notation $\mathbb D{\mc F}_{N}$ puisqu'elle figure dans l'énoncé de la  \propref{trace-frob}. 

On définit       $\mc F_{\mc T}^{2,1}=(\beta^{2,1})^{*}(\mc{F}_{N}\boxtimes \mathbb D \mc{F}_{N})$ et $\mc F_{\mc T}^{1,2}=(\beta^{1,2})^{*}(\mc{F}_{N}\boxtimes \mathbb D \mc{F}_{N})$. 
     D'où  une   correspondance cohomologique  ``image inverse'' évidente 
      $$\mc C_{\beta^{2,1}}^{*}
      \text{ de } (\mr{H}_{N} \times  \mr{H}_{N} ,\mc{F}_{N}\boxtimes \mathbb D \mc{F}_{N})
      \text{ vers }  (\mc T^{2,1}, \mc F_{\mc T}^{2,1}) $$
      supportée par  
      $$\mr{H}_{N} \times  \mr{H}_{N} \xleftarrow{\beta^{2,1}} \mc T^{2,1}
      \xrightarrow{\Id} \mc T^{2,1}.$$
     De même on une    correspondance cohomologique   ``image inverse'' 
 $\mc C_{\beta^{1,2}}^{*}$. 
      
      On a un   diagramme commutatif 
        $$   \xymatrix{
        \mc Z^{(\{2\},\{1\})}\ar[dr]_{\gamma^{2,1}}\ar[rr]_{\alpha^{2,1}}&    & \mc T^{2,1}\ar[dl]^{\delta^{2,1}}      \\
      & \mr{Gr}_{V,v}/G_{nv}\times \mr{Gr}_{V^{*},v}/G_{nv}
        &
    }$$     
     et  $\gamma^{2,1}$ et  $\delta^{2,1}$ sont lisses,  de dimension $2n\dim(G)$ et $2n\dim(G)+d$ et 
     $\alpha^{2,1}$ est, fibre par fibre un   plongement lisse de codimension $d$. 
     Plus précisément,  localement pour la topologie étale au voisinage de
      $\mc Z^{(\{2\},\{1\})}$ on peut construire un morphisme  $\kappa: \mc T^{2,1}\to \mathbb A^{d}$
      tel que 
           \begin{itemize}
      \item  le morphisme 
      \begin{gather}\label{lisse-Gr-Gr-Ad}
   (\delta^{2,1} ,\kappa):    \mc T^{2,1}\to \mr{Gr}_{V,v}/G_{nv}\times \mr{Gr}_{V^{*},v}/G_{nv}\times  \mathbb A^{d} 
      \end{gather}
      soit lisse de dimension $2n\dim(G)$, 
      \item $\alpha^{2,1}$ identifie $\mc Z^{(\{2\},\{1\})}$ avec l'image inverse de  $0\in  \mathbb A^{d}$. 
          \end{itemize}     
          En effet on peut construire $\kappa$ comme la différence entre des relevés de $( \mc G_{1}, \psi_{1})$ et $(\mc G_{2},\psi_{2})$ par un morphisme étale d'un ouvert de 
          $ \mathbb A^{d}$ dans $ \Bun_{G,N}^{\leq \mu}$ 
          (de sorte que l'annulation de $\kappa$ force l'équation 
          \eqref{cond-2-Z21}).

           Comme $$\mc F^{(\{2\},\{1\})}=
     ( \gamma^{2,1})^{*} (IC_{
      \mr{Gr}_{V,v} \times \mr{Gr}_{V^{*},v} 
            })\text{ \  et \ }  
     \mc F_{\mc T}^{2,1}= ( \delta^{2,1})^{*} (IC_{
      \mr{Gr}_{V,v} \times \mr{Gr}_{V^{*},v} 
      })[2d](d),$$ on obtient une   correspondance cohomologique
     ``image directe''
      canonique  
       $$\mc C_{!,\alpha^{2,1}}
      \text{ de } ( \mc Z^{(\{2\},\{1\})} , \mc F^{(\{2\},\{1\})})
      \text{ vers }  (\mc T^{2,1}, \mc F_{\mc T}^{2,1}) $$
      supportée par 
      $$ \mc Z^{(\{2\},\{1\})} \xleftarrow{\\Id}  \mc Z^{(\{2\},\{1\})}
      \xrightarrow{\alpha^{2,1}} \mc T^{2,1}.$$
      En effet 
      l'isomorphisme canonique $\iota^{!}\Big( E_{\mathbb A^{d}}[2d](d)\Big)\simeq E_{v}$ (où $\iota$ est l'inclusion de $0$ dans $\mathbb A^{d}$) fournit, par recollement pour la topologie étale, un isomorphisme canonique de faisceaux pervers (à décalage près) 
      \begin{gather}\label{alpha21-F-F-isom}(\alpha^{2,1})^{!}\Big(\mc F_{\mc T}^{2,1} \Big) \simeq \mc F^{(\{2\},\{1\})}.\end{gather}       
      
      De même on définit un isomorphisme 
       \begin{gather}\label{alpha12-F-F-isom}(\alpha^{1,2})^{!}\Big(\mc F_{\mc T}^{1,2} \Big) \simeq \mc F^{(\{1\},\{2\})} \end{gather}       
     et on en déduit une   correspondance cohomologique  ``image directe''
       $\mc C_{!,\alpha^{1,2}}$.
       
       La définition de $\mc C_{p,!}$ et $\mc C_{p}^{*}$ dans le diagramme ci-dessous est évidente puisque        $\mc Z$ est une réunion de copies de $v$.        On définit 
              les correspondances cohomologiques $\mc C_{\sharp}^{\mc T}$, $\mc C_{F}^{\mc T}$ et 
       $\mc C_{\flat}^{\mc T}$ du diagramme ci-dessous de la fa\c con suivante : 
       \begin{itemize}
       \item on définit    $\mc C_{\sharp}^{\mc T}$  comme 
       \begin{gather*}
    (i_{\sharp}^{\mc T})_{!}   p_{\sharp}^{*}E_{\mc Z}=
     (i_{\sharp}^{\mc T})_{!} j_{\sharp}^{*}  p_{N}^{*}E_{v}
     =(\beta^{2,1})^{*}\Delta_{!}p_{N}^{*}E_{v}\to (\beta^{2,1})^{*}(\mc{F}_{N}\boxtimes \mathbb D \mc{F}_{N})=\mc F_{\mc T}^{2,1}
       \end{gather*}
       où la deuxième égalité  est le changement de base propre 
    dans le   carré cartésien entre $j_{\sharp}$ et $\beta^{2,1}$, et le morphisme en troisième position  provient du  morphisme $\Delta_{!}p_{N}^{*}E_{v}\to \mc{F}_{N}\boxtimes \mathbb D \mc{F}_{N}$ associé à 
        $\mc C_{\sharp}^{\Delta}$, si bien que  $\mc C_{\sharp}^{\mc T}\circ \mc C_{p}^{*}=\mc C_{\beta^{2,1}}^{*}\circ \mc C_{\sharp}^{\Delta}$,   
        \item $\mc C_{F}^{\mc T}$ est associée à la correspondance image inverse 
         $$ \mc T^{2,1} \xleftarrow{\on {F}_{\{1\}}^{\mc T}}  \mc T^{1,2}\xrightarrow{\Id}  \mc T^{1,2} ,$$ et au morphisme 
         $ (\on {F}_{\{1\}}^{\mc T})^{*}(\mc F_{\mc T}^{1,2})\to \mc F_{\mc T}^{2,1}$
         qui vient par image inverse par $\beta_{1,2}$ et $\beta_{2,1}$
        du morphisme $(\Frob \times \Id)^{*}(\mc{F}_{N}\boxtimes \mathbb D \mc{F}_{N})\to \mc{F}_{N}\boxtimes \mathbb D \mc{F}_{N}$
        qui donne $\mc C_{\Frob \times \Id}$, si bien que 
        $ \mc C_{\beta^{1,2}}^{*}\circ  \mc C_{\Frob \times \Id} 
        =  \mc C_{F}^{\mc T}  \circ \mc C_{\beta^{2,1}}^{*}  $,       
       \item  on définit  $\mc C_{\flat}^{\mc T}$ comme la   correspondance cohomologique donnée par le même  morphisme 
       $$(\Delta^{\mc T})^{*}(\beta^{1,2})^{*}(\mc{F}_{N}\boxtimes \mathbb D \mc{F}_{N})=\Delta^{*}(\mc{F}_{N}\boxtimes \mathbb D \mc{F}_{N})\to p_{N}^{!}E_{v}$$ que  $\mc C_{\flat}^{\Delta}$, si bien que  $\mc C_{\flat}^{\mc T}\circ \mc C_{\beta^{1,2}}^{*}=\mc C_{\flat}^{\Delta}$. 
       \end{itemize}
       Alors le diagramme de correspondances cohomologiques suivant 
              est commutatif (à une puissance de $q^{1/2}$ près pour le carré du milieu en haut et à un signe près pour les carrés de gauche et de droite en haut):  
         \begin{gather*}
       \xymatrix{
(\mc Z,E_{\mc Z}) \ar[d]^-{\Id} \ar[r]^-{\mc C_{\sharp}}
 &  (\mc Z^{(\{2\},\{1\})},\mc F^{(\{2\},\{1\})}) \ar[d]^-{\mc C_{!,\alpha^{2,1}}} \ar[r]^-{\mc C_{F} }
 & (\mc Z^{(\{1\},\{2\})},\mc F^{(\{1\},\{2\})})   \ar[d]^-{\mc C_{!,\alpha^{1,2}}}  \ar[r]^-{\mc C_{\flat}}
  &(\mc Z,E_{\mc Z}) \ar[d]^-{\mc C_{p,!}}
    \\
(\mc Z,E_{\mc Z}) \ar[r]^-{\mc C_{\sharp}^{\mc T}}
& (\mc T^{2,1}, \mc F_{\mc T}^{2,1}) \ar[r]^-{\mc C_{F}^{\mc T}}
& (\mc T^{1,2}, \mc F_{\mc T}^{1,2}) \ar[r]^-{\mc C_{\flat}^{\mc T}}
&(v,E_{v})  
\\
(v,E_{v}) \ar[u]_-{\mc C_{p}^{*}} \ar[r]^-{\mc C_{\sharp}^{\Delta}}
&  (\mr{H}_{N} \times  \mr{H}_{N} ,\mc{F}_{N}\boxtimes \mathbb D \mc{F}_{N}) \ar[r]^-{\mc C_{\Frob \times \Id}} 
 \ar[u]_-{\mc C_{\beta^{2,1}}^{*}}
& (\mr{H}_{N} \times  \mr{H}_{N} ,\mc{F}_{N}\boxtimes \mathbb D \mc{F}_{N}) \ar[u]_-{\mc C_{\beta^{1,2}}^{*}}  \ar[r]^-{\mc C_{\flat}^{\Delta}}
&(v,E_{v}) \ar[u]_-{\Id}
 }\end{gather*}
Pour montrer la commutativité (à une puissance de $q^{1/2}$ près) du  carré entre 
$\mc C_{F}$ et $\mc C_{F}^{\mc T}$
on utilise les  morphismes vers $\mr{Gr}^{(\{1\})}_{v}/G_{nv}\times \mr{Gr}^{(\{2\})}_{v}/G_{nv}$, et des variantes évidentes de la  \propref{prop-IC-produit-Fr-partiels}.   Plus précisément, localement pour la topologie étale au voisinage de l'image de $\mc T^{2,1}$ on peut construire un diagramme commutatif 
   
         \begin{gather}\nonumber 
       \xymatrix{
      \mc T^{2,1}\ar[d] 
&&   \mc T^{1,2}  \ar[ll]_-{\on {F}_{\{1\}}^{\mc T}}  \ar[d]  
  \\
     \mr{Gr}_{V,v}/G_{nv}\times \mr{Gr}_{V^{*},v}/G_{nv}\times  \mathbb A^{d} 
 &&  \mr{Gr}_{V,v}/G_{nv}\times \mr{Gr}_{V^{*},v}/G_{nv}\times  \mathbb A^{d} 
     \ar[ll]_-{\Frob \times \Id\times \Id }  
 }\end{gather}
tel que les flèches verticales soient lisses et que 
$ \mc Z^{(\{2\},\{1\})}  \xleftarrow{\on {F}_{\{1\}}}   \mc Z^{(\{1\},\{2\})}  $
s'obtienne comme image inverse de $0\in \mathbb A^{d}$.  

 Il reste à montrer que les  carrés en haut à gauche et en haut à droite sont commutatifs (à un signe près). La  construction de 
$\mc C_{!,\alpha^{2,1}}$ ci-dessus fournit un   isomorphisme canonique  
$(\alpha^{2,1})^{!}(\mc F_{\mc T}^{2,1})\simeq \mc F^{(\{2\},\{1\})}$, d'où  un isomorphisme 
$(i_{\sharp}^{\mc T})^{!}(\mc F_{\mc T}^{2,1})\simeq i_{\sharp}^{!}(\mc F^{(\{2\},\{1\})})$ sur $\mc Y^{\sharp}$ 
 et par conséquent  il existe une  unique correspondance cohomologique 
$\mc C_{\sharp}'$ telle  que $\mc C_{!,\alpha^{2,1}}\circ \mc C_{\sharp}'=\mc C_{\sharp}^{\mc T}$.

 D'autre part on définit  $\mc C_{\flat}'$  à partir de  $\mc C_{\flat}^{\mc T}$ par changement de base propre 
 dans le   carré cartésien entre $\alpha^{1,2}$ et $j_{\flat}$, 
de sorte que  $\mc  C_{p,!}\circ \mc C_{\flat}'=\mc C_{\flat}^{\mc T}\circ
\mc C_{!,\alpha^{1,2}}$. 
Plus précisément $\mc C_{\flat}^{\mc T}$ est donnée par un morphisme 
$\mc F^{1,2}_{\mc T}\to \Delta^{\mc T}_{*}p_{N}^{!}E_{v}$, d'où l'on déduit la correspondance $\mc C_{\flat}'$ par la composition 
$$\mc F^{(\{1\},\{2\})} \isor{\eqref{alpha12-F-F-isom}}(\alpha^{1,2})^{!}\mc F^{1,2}_{\mc T}
\to (\alpha^{1,2})^{!}\Delta^{\mc T}_{*}p_{N}^{!}E_{v}
\simeq i_{\flat,*}j_{\flat}^{!} p_{N}^{!}E_{v}=i_{\flat,*}p_{\flat}^{!}E_{\mc Z} $$
où l'avant dernière étape est le changement de base propre. 

Nous sommes donc réduits à montrer que  
$\mc C_{\sharp}'=\mc C_{\sharp}$ et 
$\mc C_{\flat}'=\mc C_{\flat}$ (au signe près). Ces égalités sont vraies sur les lieux lisses  de $\mc Z^{(\{2\},\{1\})}$ et $\mc Z^{(\{1\},\{2\})}$ car ce sont les intersections avec le lieu lisse de  $H_{N}\times H_{N}$,  sur lequel tous les sous-schémas sont lisses  et les  deux  carrés cartésiens sont des intersections transverses. Maintenant on applique à  $\mc C_{\sharp}'$ et $\mc C_{\sharp}$ le lemme suivant  (l'argument pour $\mc C_{\flat}'$ et $\mc C_{\flat}$ ne sera pas répété car il est similaire et en résulte même par dualité de Verdier et   permutation de $1$ et $2$). 
\begin{lem} \label{lem-Z-Ysharp-correp}
Si dans  le diagramme ci-dessus on oublie les  troncatures par $\Bun_{G,N}^{\leq \mu}$, deux correspondances cohomologiques de $(\mc Z,E_{\mc Z})$ vers 
$(\mc Z^{(\{2\},\{1\})},\mc F^{(\{2\},\{1\})})$ supportées par  $\mc Y_{\sharp}$ 
sont égales si elles coïncident sur le lieu lisse. 
\end{lem}
Ce lemme suffit pour montrer que  $\mc C_{\sharp}'=\mc C_{\sharp}$. 
En effet on peut plonger le diagramme  ci-dessus dans le même  diagramme mais avec les  troncatures données par $\Bun_{G,N}^{\leq \mu+\kappa}$ (c'est la raison pour laquelle on avait demandé que  $\Bun_{G,N}^{\leq \mu+\kappa}$ soit un  schéma). Alors $\mc C_{\sharp}'$ et $\mc C_{\sharp}$ s'étendent en des  correspondances dans le nouveau diagramme et d'autre part l'image inverse toute entière de 
$\mc Z=\Bun_{G,N}^{\leq \mu}(\Fq)$ par $p_{\sharp}$ est incluse dans le nouveau diagramme (c'est-à-dire qu'elle n'est pas coupée par la troncature  par $\mu+\kappa$ grâce à l'hypothèse  sur $\kappa$). 

\noindent {\bf Démonstration du lemme \ref{lem-Z-Ysharp-correp}.} On a  
$   \mc Z^{(\{2\},\{1\})}=\restr{\Cht_{N, \{2,1\}, V\boxtimes V^{*}}^{(\{2\},\{1\})}}
 {\Delta(v)}$ et on pose  $$ \mc Z^{(\{2,1\})}=\restr{\Cht_{N, \{2,1\}, V\boxtimes V^{*}}^{(\{2,1\})}}
 {\Delta(v)}=\restr{\Cht_{N, \{0\}, V\otimes V^{*}}^{(\{0\})}}
 {v}.$$ On a  
 $$
(\pi^{(\{2\},\{1\})}_{(\{2,1\})})_{!}(\mc F^{(\{2\},\{1\})})=
\restr{
\mc F_{N,\{2,1\},V\boxtimes V^{*}, \Xi,E}^{\{2,1\}}
}{\Delta(v)}=\restr{\mc F_{N,\{0\},V\otimes V^{*}, \Xi,E}^{\{0\}}}{v}. 
$$
On a  une inclusion évidente $\mc Z\overset{\iota}{\hookrightarrow}  \mc Z^{(\{2,1\})}$ et un carré cartésien 
 \begin{gather}\label{cart-square-Y-sharp-gr}\xymatrix{
 \mc Y_{\sharp} \ar[d]^-{p_{\sharp}} \ar@{^{(}->}[r]^-{\mf i_{\sharp}} &   \mc Z^{(\{2\},\{1\})}  \ar[d]^-{\pi^{(\{2\},\{1\})}_{(\{2,1\})}}
  \\
\mc Z \ar[r]^-{\iota}
 &\mc Z^{(\{2,1\})}  }\end{gather}
 Une correspondance cohomologique 
 $$ \text{ de }(\mc Z,E_{\mc Z})\text{ vers 
}(\mc Z^{(\{2\},\{1\})},\mc F^{(\{2\},\{1\})})
\text{ supportée par }\mc Y_{\sharp}$$
 est la même chose qu'une    correspondance cohomologique  $$ \text{ de }(\mc Z,E_{\mc Z})\text{ vers 
}(\mc Z^{(\{2,1\})},\restr{\mc F_{N,\{0\},V\otimes V^{*}, \Xi,E}^{\{0\}}}{v})\text{  supportée par 
}\mc Z \xleftarrow{\Id} \mc Z \xrightarrow{\iota} \mc Z^{(\{2,1\})},$$ car par  changement de base propre dans  \eqref{cart-square-Y-sharp-gr} on a  
$$ (p_{\sharp})_{*}\mf i_{\sharp}^{!}(\mc F^{(\{2\},\{1\})})=
\iota^{!}(\pi^{(\{2\},\{1\})}_{(\{2,1\})})_{*}(\mc F^{(\{2\},\{1\})})=
\iota^{!}\Big(\restr{\mc F_{N,\{0\},V\otimes V^{*}, \Xi,E}^{\{0\}}}{v}\Big).$$
Or $\iota^{!}\Big(\restr{\mc F_{N,\{0\},V\otimes V^{*}, \Xi,E}^{\{0\}}}{v}
\Big)$ est supportée en  degré $\geq 0$ (puisque $\restr{\mc F_{N,\{0\},V\otimes V^{*}, \Xi,E}^{\{0\}}}{v}$ est pervers) et seul  $\mbf 1 \subset V\otimes V^{*}$ contribue à sa composante de degré $0$. 
Par conséquent  une telle correspondance cohomologique  est simplement donnée par la  multiplication par 
une  fonction sur $\mc Z$, et elle  est clairement déterminée par sa  restriction sur le  lieu lisse de $\mc Z^{(\{2\},\{1\})}$ (ou $\mc Y^{\sharp}$). \cqfd

La commutativité du  diagramme ci-dessus montre que les trois  composées
            $\mc C_{\flat}\circ \mc C_{F} \circ\mc C_{\sharp}$, 
            $\mc C_{\flat}^{\mc T}\circ \mc C_{F}^{\mc T} \circ\mc C_{\sharp}^{\mc T}$ et 
            $\mc C_{\flat}^{\Delta}\circ \mc C_{\Frob \times \Id}\circ
            \mc C_{\sharp}^{\Delta}$, sont égales (à une puissance de $q^{1/2}$ et un signe près pour le moment). 
           Nous savons qu'elles   sont supportées 
            par $\Gamma=H_{N}(\Fq)$, et d'après la 
            \propref{trace-frob} (appliqués à $\mr{H}_{N}$ et $\mc F_{N}$),  $\mc C_{\flat}^{\Delta}\circ \mc C_{\Frob \times \Id}\circ
            \mc C_{\sharp}^{\Delta}$  est  donné  par  la  multiplication par $h_{V,v}$ (à une puissance de $q^{1/2}$ et un signe près).
            Par conséquent $\mc C_{\flat}\circ \mc C_{F} \circ\mc C_{\sharp}$ est donné  par la multiplication par $h_{V,v}$ (à une puissance de $q^{1/2}$ et un signe près pour le moment). On pourrait expliciter tous les signes et les puissances de $q^{1/2}$ dans l'argument précédent mais pour terminer la preuve du  \lemref{lem-equiv-ST2} le plus simple est de remarquer que la restriction de $h_{V,v}$ à l'orbite correspondant au plus haut poids de $V$ dans 
            $G(\mc O_{v})\backslash G(F_{v})/G(\mc O_{v})$ est non nulle et correspond dans la composition $\mc C_{\flat}^{\Delta}\circ \mc C_{\Frob \times \Id}\circ
            \mc C_{\sharp}^{\Delta}$ à une situation où tous les lieux sont lisses et les intersections transverses. Cette situation est en tous points identique au cas minuscule traité dans l'introduction, dont il suffit de répéter les arguments. \cqfd
             
             La preuve de la \propref{prop-coal-frob-cas-part} était très longue. Les trois remarques ci-dessous proposent des idées de preuves différentes. Enfin dans la \remref{rem-varsha} on évoquera une preuve beaucoup plus simple (mais reposant sur un formalisme général qui sort du cadre de cet article) que Yakov Varshavsky a indiquée récemment à l'auteur \cite{varshav-com-perso}.

             \begin{rem} \label{rem-proof-ST}
 On suppose que  $\deg(v)=1$ car autrement les notations seraient trop compliquées. 
 Les   correspondances cohomologiques  $\mc C_{1}=\mc C_{\flat}$ et $\mc C_{2}=\mc C_{F}\circ \mc C_{\sharp}$ sont supportées par $\mc Y_{1}=\mc Y_{\flat}$ et $\mc Y_{2}=\mc Y_{\sharp}\times _{\mc Z^{(I,\{2\},\{1\})}}\mc Z^{(\{1\},\{2\},I)}$. 
 Par suite  $\mc C_{1}\circ \mc C_{2}=\mc C_{\flat}\circ \mc C_{F}\circ \mc C_{\sharp}$ est supportée par 
 $\mc Y_{1}  \times_{\mc Z^{(\{1\},\{2\},I)}}\mc Y_{2}=\Gamma^{(I)}$.  
 On peut montrer que cette intersection est transverse dans le  sens 
 où on est dans  une situation produit localement pour la topologie étale (voir la v5 de cet article sur arXiv pour plus de détails).  
  Cela devrait impliquer que les termes locaux  pour la  composée des  correspondances cohomologiques  $\mc C_{1}$ et $\mc C_{2}$ supportées par 
    $\mc Y_{1}$ et $\mc Y_{2}$ peuvent être calculés directement sur les fibres des  faisceaux  en les points de l'intersection  \eqref{Gamma-produit-fibre}, et on obtient facilement   $h_{V,v}$. Mais nous avons manqué de références. 
   \end{rem}
   
    \begin{rem} Une autre idée pour montrer la \propref{prop-coal-frob-cas-part} en évitant  ces difficultés techniques 
     serait d'utiliser les résolutions de   Bott-Samelson pour avoir des champs lisses et des intersections de sous-champs lisses qui soient  transverses au sens usuel.  Un inconvénient de cette approche est qu'elle requiert des arguments supplémentaires dans le cas non déployé considéré dans le chapitre \ref{para-non-deploye}. 
   \end{rem}
         
     \begin{rem}
     Dans la situation où $I$ est vide et $\deg(v)=1$, la composition $\mc C_{\flat}\circ \mc C_{F}\circ \mc C_{\sharp}$ 
apparaît, avec les notations de la \remref{rem-justif-norm-crea-annihil},  comme une sorte de restriction au graphe de 
Frobenius de $\Hecke_{\{1 \}, V  }^{\{1\} }$ vers 
$\Hecke_{\{3 \}, V  }^{\{3\} }$ de la composition 
$\mc D_{\flat}\circ \mc D_{\sharp}$. Mais nous ne sommes pas parvenus à mettre en forme cette idée.  \end{rem}

La conclusion de ces trois remarques est que l'auteur n'a  pas trouvé de  preuve plus simple de la     \propref{prop-coal-frob-cas-part}. 
Quelle que soit la méthode employée, il faut nécessairement montrer un résultat de transversalité. Enoncer et exploiter un tel résultat  est   rendu plus difficile par les singularités et  le cas où $\deg(v)>1$.     

\begin{rem} \label{rem-varsha}
Cependant Yakov Varshavsky~\cite{varshav-com-perso} a récemment indiqué à l'auteur que l'égalité 
schématique \eqref{Gamma-produit-fibre} suffit  en fait, grâce à des arguments généraux,   à impliquer un certain énoncé de transversalité et  à en déduire  que la composée des correspondances cohomologiques est égale à ce que l'on attend.  Cela repose sur un formalisme général qui dépasse le cadre ce cet article. 
\end{rem}
       
  \section{Relations d'Eichler-Shimura}
  \label{para-Relations d'Eichler-Shimura}
  
   Soit   $I$ un ensemble fini, 
 $W$ une représentation de dimension finie de $(\wh G)^{I}$ et  $V$ une  représentation  irréductible de $\wh G$.
  On considère  $W \boxtimes V$ comme une  représentation de 
$(\wh G)^{I\cup \{0\}}$. 
  Soit  $v$ une  place de  $X\sm N$. Comme précédemment   $E_{v}$ désigne  le faisceau constant en   $v$.

   Le  morphisme  $$F_{\{0\}}^{\deg(v)}: \restr{\mc H _{ N, I\cup\{0\}, W\boxtimes V}^{\leq\mu,E}}{(X\sm N)^{I}\times v}\to \restr{\mc H _{ N, I\cup\{0\}, W\boxtimes V}^{\leq\mu+\kappa,E}}{(X\sm N)^{I}\times v}$$
est bien défini  (pour $\kappa$ assez grand ) parce que   $\Frob_{v}^{\deg(v)}$ est l'identité en 
$v$. La  relation  d'Eichler-Shimura  exprime    qu'il est annulé par un  polynôme de degré  $\dim(V)$ dont les  coefficients sont 
les restrictions de $(X\sm N)^{I\cup \{0\}}$ 
à $(X\sm  N )^{I }\times v$   
des opérateurs de Hecke {\it étendus} $T(h_{\Lambda^{i}V,v})$ de la \propref{prop-coal-frob-cas-part}. Il est vraiment nécessaire de considérer les opérateurs étendus car les opérateurs de Hecke $T(h_{\Lambda^{i}V,v})$  du paragraphe  \ref{action-hecke-etale} 
n'étaient définis que sur  $(X\sm (N\cup v))^{I\cup \{0\}}$. On rappelle que dans la \propref{prop-coal-frob-cas-part} l'opérateur  de Hecke  étendu  $T(h_{\Lambda^{i}V,v})$ a été défini comme  le morphisme $S_{\Lambda^{i}V,v}$. 
 Pour cette raison on va énoncer la  proposition suivante avec  $S_{\Lambda^{i}V,v}$. 
  Un autre avantage est que cela la rend logiquement indépendante de la 
     \propref{prop-coal-frob-cas-part}. 
      En effet, grâce à la  définition des  morphismes  $S_{\Lambda^{i}V,v}$ par  \eqref{def-SVv-text1}-\eqref{def-SVv-text4}, la preuve va consister en un simple calcul d'algèbre tensorielle   (inspiré par une preuve du théorème  d'Hamilton-Cayley).

 \begin{prop}\label{prop-eichler-shimura}
Pour  $\kappa$ assez grand   (en fonction de  $\deg(v)$ et $V$), on a  
  \begin{gather}\label{egalite-prop-eichler-shimura} \sum_{i=0}^{\dim V} (-1)^{i} (F_{\{0\}}^{\deg(v)})^{i}\circ \restr{S_{\Lambda^{\dim V-i}V,v}}{(X\sm N)^{I}\times v}=0   \ \ \text{dans  }  
 \\ \nonumber  \Hom_{D^{b}_{c}((X\sm N)^{I}\times v, E)}\Big(
 \restr{ \mc H _{ N,I\cup \{0\},W\boxtimes V}^{\leq\mu,E}}{(X\sm N)^{I}\times v},
  \restr{ \mc H _{ N,I\cup \{0\},W\boxtimes V}^{\leq\mu+  \kappa,E}}{(X\sm N)^{I}\times v} 
  \Big). \end{gather} 
 \end{prop}
 
  Avant d'attaquer la démonstration on commence par rappeler une preuve ``tensorielle'' du théorème de Hamilton-Cayley, d'après le paragraphe 6.5 de~\cite{cvitanovic} et \cite{peterson}. En fait cette  preuve apparaît déjà  dans la démonstration du théorème  13.4.12
  de \cite{jannsen-murre} (dans un esprit assez proche du nôtre). 
  
  Soit $k$ un corps  de caractéristique $0$ (par exemple $E$). Pour tout ensemble fini $J$, on possède l'antisymétriseur 
$$\mc A_{J}=\frac{1}{(\sharp J)!}\sum_{\sigma\in \mf S (J)}s(\sigma)\sigma$$ dans l'algèbre du groupe des  permutations de $J$ (où $ s(\sigma)$ est la signature de $\sigma$). 
Pour tout  $k$-espace vectoriel $V$   de dimension finie,  $\mc A_{J}$ agit sur $V^{\otimes J}$ par un idempotent dont l'image est $\Lambda^{\sharp J} V$. 
  
  Soit $J$ un ensemble fini et $V$ un $k$-espace vectoriel     de dimension finie. 
  Soit $T\in \on{End}(V)$   et $n\in \N$. Pour tout $U\in 
  \on{End}(  V^{\otimes  \{0\}\cup J})$ on note $\mf C_{J}(T,U)\in \on{End}(V)$ 
   l'opérateur composé
  \begin{gather*}V\xrightarrow{\Id_{V}\otimes (\delta_{V})^{\otimes J}}
 V\otimes  (V  \otimes  V^{*})^{\otimes  J}
  =V^{\otimes \{0 \}\cup J } \otimes (V^{*})^{\otimes  J}
  \\
  \xrightarrow{U \otimes \Id_{(V^{*})^{\otimes   J}}}
 V^{\otimes \{0 \}\cup J } \otimes (V^{*})^{\otimes  J}
=  V\otimes  V^{\otimes   J} \otimes (V^{*})^{\otimes  J} 
   \\ \xrightarrow{ \Id_{V}\otimes T^{\otimes J}\otimes \Id_{(V^{*})^{\otimes  J}}}
    V\otimes  V^{\otimes   J} \otimes (V^{*})^{\otimes  J} 
    =V\otimes  (V  \otimes  V^{*})^{\otimes  J}
  \xrightarrow{\Id_{V}\otimes (  \on{ev}_{V})^{\otimes J}
}V  
  \end{gather*}
       où $  \{0 \}\cup J$ est bien sûr une réunion disjointe. 
       On va appliquer cette construction à $J=\{1,...,n\}$, de sorte que 
       $ \{0 \}\cup J=\{0,...,n\}$. 
On va montrer l'égalité   \begin{gather}\label{enonce-CnTA}\mf C_{\{ 1,...,n\}}(T,\mc A_{\{0,1,...,n\}})=\frac{1}{n+1}\sum_{i=0}^{n}(-1)^{i}\on{Tr}(\Lambda^{n-i}T) T^{i} \text{ \ dans \ } \rm{End}(V).\end{gather} 
    Cette égalité  fournit une preuve du théorème de Hamilton-Cayley car  si $n=\dim V$,  on a $\mc A_{\{0,1,...,n\}}=0$, donc 
  $\mf C_{\{ 1,...,n\}}(T,\mc A_{\{0,1,...,n\}})=0 $ et \eqref{enonce-CnTA} est la relation de Hamilton-Cayley 
  $\sum_{i=0}^{n}(-1)^{i}\on{Tr}(\Lambda^{n-i}T) T^{i}=0$. 
 
   Voici maintenant la démonstration de \eqref{enonce-CnTA}. 
     On développe $\mc A_{\{0,1,...,n\}}$ en une somme sur les permutations 
  $\sigma$ de $\{0,...,n\}$ et on distingue suivant la longueur $\ell(\sigma,0)$ du cycle contenant $0$. Pour $i\in \{0,...,n\}$ le nombre de permutations telles que $\ell(\sigma,0)=i+1$ vaut toujours $n!$, car il y a 
  $n(n-1)\cdots(n-i+1)$ possibilités pour le cycle, et $(n-i)!$ possibilités pour la permutation restante. Il suffit donc de  démontrer que pour tout $i\in \{0,...,n\}$, 
  $$\mf C_{\{ 1,...,n\}}\Big(T,\frac{1}{n!}\sum_{\sigma\in \mf S(\{0,...,n\}), \ell(\sigma,0)=i+1}s(\sigma)\sigma \Big)
  =(-1)^{i}\on{Tr}(\Lambda^{n-i}T) T^{i}.$$
  Comme $\mf C_{\{ 1,...,n\}}(T,U)$ est inchangé si on conjugue $U$ par une  
  permutation de $\{0,...,n\}$ fixant $0$,  on peut remplacer la moyenne sur les permutations $\sigma$ telles que $\ell(\sigma,0)=i+1$ par une moyenne sur les permutations $\sigma$ telles que $0\xrightarrow{\sigma} 1 \xrightarrow{\sigma} \cdots \xrightarrow{\sigma} i \xrightarrow{\sigma} 0$ et 
  il suffit donc de démontrer que     
  \begin{gather}\label{egalite-fin0}\mf C_{\{ 1,...,n\}}\Big(T,\frac{1}{(n-i)!}\sum_{\sigma\in \mf S(\{0,...,n\}), 0\xrightarrow{\sigma} 1 \xrightarrow{\sigma} \cdots \xrightarrow{\sigma} i \xrightarrow{\sigma} 0}s(\sigma)\sigma \Big)
  =(-1)^{i}\on{Tr}(\Lambda^{n-i}T) T^{i}.\end{gather}  
 Dans  \eqref{egalite-fin0} $\sigma$ est simplement la  juxtaposition de la permutation circulaire $\kappa_{i}: 0\to  1 \to  \cdots \to  i \to  0$ et d'une  permutation $\tau$
    de  $\{i+1,...,n\}$ et   la moyenne porte sur  $\tau$. On a  $s(\sigma)=(-1)^{i}s(\tau)$.  Donc le membre de gauche de \eqref{egalite-fin0} est égal au produit 
    \begin{itemize}
    \item de  $(-1)^{i}$
    \item du scalaire égal à  la composée 
         \begin{gather*}k\xrightarrow{  (\delta_{V})^{\otimes \{i+1,...,n\}}}
    (V  \otimes  V^{*})^{\otimes  \{i+1,...,n\}}
  =V^{\otimes    \{i+1,...,n\} } \otimes (V^{*})^{\otimes  \{i+1,...,n\}}
  \\
  \xrightarrow{\mc A_{\{i+1,...,n\}} \otimes \Id_{(V^{*})^{\otimes   \{i+1,...,n\}}}}
 V^{ \otimes  \{i+1,...,n\} } \otimes (V^{*})^{\otimes  \{i+1,...,n\}} \\
  \xrightarrow{T^{\otimes \{i+1,...,n\}} \otimes \Id_{(V^{*})^{\otimes   \{i+1,...,n\}}}}
 V^{  \otimes \{i+1,...,n\} } \otimes (V^{*})^{\otimes  \{i+1,...,n\}} \\
 =  (V  \otimes  V^{*})^{\otimes  \{i+1,...,n\}}
  \xrightarrow{  (  \on{ev}_{V})^{\otimes \{i+1,...,n\}}
}k  
  \end{gather*}
qui est évidemment $\on{Tr}(\Lambda^{n-i}T)$, 
    \item et de $\mf C_{\{1,...,i\}}(T, \kappa_{i} ):V\to V$, où $\kappa_{i}:V^{\otimes \{0,...,i\}}\to V^{\otimes \{0,...,i\}}$ agit par permutation circulaire des facteurs. 
    \end{itemize}
          On en déduit \eqref{enonce-CnTA},  car  on va montrer que 
  $\mf C_{\{1,...,i\}}(T,\kappa_{i} )$ est égal à $T^{i}:V\to V$. En effet, en changeant l'ordre des facteurs dans 
  $V\otimes V^{\otimes \{1,...,i\}} \otimes (V^{*})^{\otimes \{1,...,i\}}$
  en $V\otimes (V^{*}\otimes V)^{\otimes i}$ on obtient que 
   $\mf C_{\{1,...,i\}}(T, \sigma )$ est égal à la composée 
   \begin{gather}\label{composee-V1ii}
 V\xrightarrow{\Id_{V}\otimes \delta_V^{\otimes i}} V\otimes (V^*\otimes V)^{\otimes i}\xrightarrow{T \otimes (\Id_{V^*}\otimes T )\otimes \cdots \otimes
 (\Id_{V^*}\otimes T )
 }
 V\otimes (V^*\otimes V)^{\otimes i}\\
 \nonumber 
 =(V\otimes V^*)^{\otimes i}
 \otimes  V\xrightarrow{ \on{ev}_V^{\otimes i} \otimes\Id_{V}}  V. 
 \end{gather}
  On  montre par récurrence sur $i$ que \eqref{composee-V1ii} est égal à $T^{i}$ en utilisant le fait que   la composée
\begin{gather}\label{delta-ev-VVV}
 V\xrightarrow{\Id_{V}\otimes \delta_V} V\otimes (V^*\otimes V)=
 (V\otimes   V^{*})
 \otimes  V\xrightarrow{ \on{ev}_V \otimes\Id_{V}}  V\end{gather}   est égale à   $\Id_{V}$ d'après le lemme de Zorro \eqref{zorro}. Ceci termine l'explication de la 
 preuve ``tensorielle'' du théorème de Hamilton-Cayley.

\noindent {\bf Démonstration de la \propref{prop-eichler-shimura}.} On rappelle que pour tout ensemble fini $J$, on possède l'antisymétriseur 
$$\mc A_{J}=\frac{1}{(\sharp J)!}\sum_{\sigma\in \mf S (J)}s(\sigma)\sigma. $$      
   
   Pour tout entier  $n$ et tout  endomorphisme $U$ de $V^{\otimes \{0,1,...,n\}}=V^{\otimes (n+1)}$, on note  $\mf C_{n}(U)$ 
   (pour simplifier par rapport à la notation $\mf C_{\{0,1,...,n\}}(U)$ utilisée précédemment) 
   la  composée 
         \begin{gather}\label{morph-compose-Eichler-Shimura}
    \restr{ \mc H _{ N,I\cup \{0\},W\boxtimes V}^{\leq\mu,E}}{(X\sm N)^{I}\times v} 
   \xrightarrow{ \restr{\mc C_{  
       \delta_{V^{\otimes n}}}^{\sharp }}{(X\sm N)^{I}\times v}}
       \\ \nonumber
  \restr{ \mc H _{ N,I\cup \{0\}\cup \{1,...,n\}\cup \{n+1, ...,2n\},W\boxtimes V\boxtimes V^{\boxtimes n} \boxtimes (V^{*})^{\boxtimes n}}^{\leq\mu,E}}{(X\sm N)^{I}\times \Delta(v)} \\ \nonumber 
   \xrightarrow{ \mc H(\Id_{W}\boxtimes U\boxtimes 
   \Id_{(V^{*})^{\boxtimes n}})} \restr{ \mc H _{ N,I\cup \{0\}\cup \{1,...,n\}\cup \{n+1, ...,2n\},W\boxtimes V\boxtimes V^{\boxtimes n} \boxtimes (V^{*})^{\boxtimes n}}^{\leq\mu,E}}{(X\sm N)^{I}\times \Delta(v)} \\ \nonumber
 \xrightarrow{ \prod_{j\in \{1,...,n\}}(F_{\{j\}})^{\deg(v)} } 
  \restr{ \mc H _{ N,I\cup \{0\}\cup \{1,...,n\}\cup \{n+1, ...,2n\},W\boxtimes V\boxtimes V^{\boxtimes n} \boxtimes (V^{*})^{\boxtimes n}}^{\leq\mu+n\deg(v)\kappa,E}}{(X\sm N)^{I}\times \Delta(v)} \\ \nonumber
  \xrightarrow{ \restr{\mc C_{  
      \on{ev}_{V^{\otimes n}}}^{\flat }
     }{(X\sm N)^{I}\times v}}
     \restr{ \mc H _{ N,I\cup \{0\},W\boxtimes V}^{\leq\mu+n\deg(v)\kappa,E}}{(X\sm N)^{I}\times v} 
     \end{gather}
    Autrement dit, on part de   $I\cup \{0\}$ (avec  la patte $0$ fixée en  $v$)  et de  la  représentation $W\boxtimes V$, puis 
    \begin{itemize}
    \item on crée des paires  de pattes $(1,n+1), (2,n+2), ..., (n,2n)$ en  $v$ à l'aide de 
    $\delta_{V}:\mbf 1\to V\otimes V^{*}$, 
    \item on applique $U$ aux  pattes dans $\{0,...,n\}$ et à la  représentation $V^{\otimes \{0,1,...,n\}}$, 
    \item on applique le  morphisme de Frobenius partiel $(F_{\{j\}})^{\deg(v)}$ aux  pattes  $j\in \{1,...,n\}$
    \item on annihile les paires  de pattes $(1,n+1), (2,n+2), ..., (n,2n)$ en  $v$ à l'aide de 
    $\on{ev}_{V}: V\otimes V^{*}\to\mbf 1$. 
    \end{itemize}
    
    On prend  $n=\dim V$. 
    Comme  $\Lambda^{n+1}V=0$, $\mc A_{\{0,1,...,n\}}$ agit  par zéro sur $V^{\otimes \{0,...,n\}}$ et donc  $\mf C_{n}(\mc A_{\{0,1,...,n\}})$  est nul. Mais d'un autre côté, en développant la moyenne dans 
    $\mc A_{\{0,1,...,n\}}$ on va montrer que  $\mf C_{n}(\mc A_{\{0,1,...,n\}})$ est égal au membre de gauche de \eqref{egalite-prop-eichler-shimura} (divisé par $n+1$). 
    
   On  développe donc le morphisme composé \eqref{morph-compose-Eichler-Shimura} en écrivant   $\mc A_{\{0,1,...,n\}}$ comme une   moyenne sur les  permutations 
  $\sigma$ de  $\{0,...,n\}$ et on distingue suivant la longueur   $\ell(\sigma,0)$ du cycle contenant $0$. Pour tout $i\in \{0,...,n\}$ le nombre  de permutations telles  que $\ell(\sigma,0)=i+1$ est toujours $n!$, parce qu'il y a   $n(n-1)\cdots(n-i+1)$ possibilités pour le  cycle contenant $0$, et  $(n-i)!$ possibilités pour le reste de la   permutation. Il suffit donc de montrer que  pour tout   $i\in \{0,...,n\}$, 
  $$\mf C_{n}\Big(\frac{1}{n!}\sum_{\sigma\in \mf S(\{0,...,n\}), \ell(\sigma,0)=i+1}s(\sigma)\sigma \Big)
  =(-1)^{i} (F_{\{0\}}^{\deg(v)})^{i}\circ S_{\Lambda^{\dim V-i}V,v}
  .$$
  Comme  $\mf C_{n}(U)$ ne  change pas quand on conjugue  
   $U$ par une  
  permutation de  $\{0,...,n\}$ fixant $0$ (car cela revient seulement à changer les noms des indices dans $\{1,...,n\}$ et $\{n+1,...,2n\}$ qui apparaissent comme intermédiaires dans la construction), 
    on peut remplacer la moyenne sur les   permutations $\sigma$ telles  que  $\ell(\sigma,0)=i+1$ par une moyenne sur les  permutations $\sigma$ telles que le cycle contenant $0$ soit  $0\xrightarrow{\sigma} 1 \xrightarrow{\sigma} \cdots \xrightarrow{\sigma} i \xrightarrow{\sigma} 0$. Autrement dit il suffit de montrer que 
    \begin{gather}\label{egalite-fin}\mf C_{n}\Big(\frac{1}{(n-i)!}\sum_{\sigma\in \mf S(\{0,...,n\}), 0\xrightarrow{\sigma} 1 \xrightarrow{\sigma} \cdots \xrightarrow{\sigma} i \xrightarrow{\sigma} 0}s(\sigma)\sigma \Big)
  =(-1)^{i} (F_{\{0\}}^{\deg(v)})^{i}\circ S_{\Lambda^{\dim V-i}V,v}
.\end{gather}  
  Dans  \eqref{egalite-fin} $\sigma$ est simplement la  juxtaposition de la permutation circulaire $0\to  1 \to  \cdots \to  i \to  0$ et d'une  permutation $\tau$
    de  $\{i+1,...,n\}$ et   la moyenne porte sur  $\tau$. On a  $s(\sigma)=(-1)^{i}s(\tau)$ et cela explique le signe  $  (-1)^{i}$ dans le membre de droite de \eqref{egalite-fin}. 
 Les pattes $\{0\}\cup \{1,...,i\}\cup \{n+1,...,n+i\}$ d'un côté et $\{i+1,...,n\}\cup \{n+i+1,...,2n\}$ de l'autre côté jouent des rôles   complètement  indépendants (les pattes dans $I$ ne jouent quant à elles aucun rôle). 
  
    Ce qui arrive aux pattes dans $\{i+1,...,n\}\cup \{n+i+1,...,2n\}$ est la chose suivante:     \begin{itemize}
    \item on crée les paires  de pattes $(i+1,n+i+1), (i+2,n+i+2), ..., (n,2n)$ en $v$ à l'aide de 
    $\delta_{V}:\mbf 1\to V\otimes V^{*}$, 
    \item on antisymétrise  les pattes dans $\{i+1,...,n\}$ (grâce à  la moyenne  sur $\tau$),  
    \item on applique le  morphisme de Frobenius partiel en les  pattes de $\{i+1,...,n\}$
    \item on annihile les paires  de pattes $(i+1,n+i+1), (i+2,n+i+2), ..., (n,2n)$ en  $v$ à l'aide de     $\on{ev}_{V}: V\otimes V^{*}\to\mbf 1$. 
    \end{itemize}  
   Comme $\mc A_{\{i+1,...,n\}}$ agit  sur $V^{\otimes \{i+1,...,n\}}$ par un idempotent dont  l'image est $\Lambda^{\dim V-i}V$,  cette  composée est exactement le  morphisme 
   $ S_{\Lambda^{\dim V-i}V,v}$. Pour justifier  ceci formellement 
   on procède  comme dans la preuve de \eqref{SIW-p0-intro}, avec  $u$ égal à l'inclusion de 
  $\Lambda^{j}V\otimes (\Lambda^{j}V)^{*}$ dans $V^{\otimes j} \otimes (V^{*})^{\otimes j}$, qui est telle que 
       $$(\mc A_{j}\otimes 
 \Id_{(V^{*})^{\otimes j}})\circ \delta_{V^{\otimes j}}
 =u\circ \delta_{\Lambda^{j}V} \text{ \  et  \ } 
  \on{ev}_{V^{\otimes j}}\circ u= \on{ev}_{\Lambda^{j}V} . $$

   Ce qui arrive aux pattes dans $\{0\}\cup \{1,...,i\}\cup \{n+1,...,n+i\}$ est la chose suivante:     
 \begin{itemize}
    \item on crée les paires  de pattes $(1,n+1), (2,n+2), ..., (i,n+i)$ en  $v$ à l'aide de 
    $\delta_{V}:\mbf 1\to V\otimes V^{*}$, 
    \item on applique la permutation circulaire $0\to  1 \to  \cdots \to  i \to  0$  aux  pattes de $\{0,...,i\}$, 
    \item on applique le  morphisme de Frobenius partiel en  les pattes de $\{1,...,i\}$
    \item on annihile les paires  de pattes $(1,n+1), (2,n+2), ..., (i,n+i)$ en  $v$ à l'aide de 
    $\on{ev}_{V}: V\otimes V^{*}\to\mbf 1$. 
    \end{itemize}
Si on  numérote l'ensemble  $\{0\}\cup \{1,...,i\}\cup \{n+1,...,n+i\}$ de $0$ à $2i$ dans l'ordre mélangé   $(0,n+1,1,n+2,...,n+i,i)$, la  composée précédente devient égale à la suivante: partant de  $I\cup \{0\}$, 
\begin{itemize}
    \item on crée les paires  de pattes  $(1,2), (3,4), ..., (2i-1,2i)$ en  $v$ à l'aide de 
      $\delta_{V}:\mbf 1\to V^{*}\otimes V$, et on arrive dans 
    $$ \restr{ \mc H _{ N,I\cup \{0\}\cup \{1,...,2i\},W\boxtimes V\boxtimes V^{*}\cdots \boxtimes V\boxtimes V^{*}\boxtimes V}^{\leq\mu,E}}{(X\sm N)^{I}\times \Delta(v)}$$
       \item on applique le  morphisme de Frobenius partiel en  les pattes de $\{0,2,4,...,2i-2\}$
    \item on annihile les paires  de pattes $(0,1), (2,3), ..., (2i-2,2i-1)$ en $v$ à l'aide de     $\on{ev}_{V}: V\otimes V^{*}\to\mbf 1$. 
    \end{itemize}
Par récurrence  sur $i$ on montre facilement que cette dernière composée est 
égale à $(F_{\{0\}}^{\deg(v)})^{i}$, à l'aide de l'égalité entre \eqref{delta-ev-VVV} et   $\Id_{V}$.
Signalons que ce dernier argument  (c'est-à-dire l'utilisation  de la permutation circulaire pour obtenir  $(F_{\{0\}}^{\deg(v)})^{i}$) 
est très voisin du    théorème 1 et   du lemme  2 dans \cite{ngo-modif-sym}. 
\cqfd

\begin{rem} Dans \cite{xiao-zhu}, Liang Xiao et  Xinwen Zhu ont  
défini (dans un cadre un peu différent) un anneau de correspondances cohomologiques entre $\restr{ \Cht_{ N,I\cup \{0\},W\boxtimes V}^{(I,\{0\})}}{(X\sm N)^{I}\times v}$ et lui-même, dans lequel la relation d'Eichler-Shimura résulte formellement de Hamilton-Cayley. 
\end{rem}

  \section{Sous-faisceaux  constructibles stables par les  morphismes de  Frobenius partiels}
  \label{para-sous-faisceaux-Frob-partiels}

 On appelle  point géométrique  $\ov x$ d'un  schéma  $Y$ la donnée d'un corps  algébriquement clos  $k(\ov x)$ et d'un   morphisme $\on{Spec}(k(\ov x))\to Y$ 
 (dans \cite{grothendieck-sga4-2-VIII} la définition des points géométriques fait intervenir des extensions séparablement closes, mais nous ne suivons pas cette convention ici). 
  Autrement dit c'est la donnée d'un  point  $x$ de  $Y$ et d'une extension  algébriquement close     $k(\ov x)$ de  $k(x)$.  Soit  $A=\mc O_{E},E$  ou  $\Qlbar$. Si  $\mc F$ est un  $A$-faisceau constructible  sur un voisinage  de $x$ dans  $Y$,  on note  $\mc F_{\ov x}$ ou  $\restr{\mc F}{\ov x}$ la fibre de  $\mc F$ en   $\ov x$ (qui est un    $A$-module de type fini). On note  $Y_{(\ov x)}$ le localisé strict  (ou   hensélisé strict) de  $Y$ en   $\ov x$. 
  Autrement dit  $Y_{(\ov x)}$ est le spectre de l'anneau   $\varinjlim \Gamma(U,\mc O_{U})$, où la limite inductive est prise sur les  voisinages étales  $\ov x$-pointés de  $x$. C'est un anneau   local hensélien dont le corps résiduel est la  clôture séparable de  $k(x)$ dans  $k(\ov x)$. 
 Si  $\ov x$ et $\ov y$ sont deux points géométriques de $Y$, on appelle 
  flèche de spécialisation $\on{\mf{sp}}:\ov x\to \ov y$   un  morphisme  
 $Y_{(\ov x)}\to Y_{(\ov y)}$, ou de fa\c con équivalente un    morphisme  
 $\ov x\to Y_{(\ov y)}$ (une telle flèche existe si et seulement si $y$ est dans l'adhérence de Zariski de $x$). D'après le paragraphe 7 de \cite{grothendieck-sga4-2-VIII} une   flèche de spécialisation $\on{\mf{sp}}:\ov x\to \ov y$ induit pour tout $A$-faisceau constructible  $\mc F$ sur un ouvert de $Y$ contenant $y$ un homomorphisme de spécialisation $\on{\mf{sp}}^{*}: \mc F_{\ov y}\to \mc F_{\ov x}$ (qui découle simplement de la définition de la fibre d'un faisceau en un point géométrique). 

  Comme on l'avait déjà fait dans l'introduction, 
  on fixe une clôture algébrique $\ov F$ de  $F$ et on  note 
      $\ov \eta=\on{Spec}(\ov F)$ le point géométrique correspondant au-dessus du   point  générique
 $\eta$    de $X$.

  Soit 
  $I$ un ensemble fini. On note    $\Delta: X\to X^{I}$  le morphisme diagonal.  On note 
  $\eta^{I}=\on{Spec}(F^{I})$   le  point  générique de  $X^{I}$. 
  On fixe  un  point géométrique   $\ov{\eta^{I}}$ au-dessus de 
     $\eta^{I}$, muni d'une  flèche de spécialisation 
     $\on{\mf{sp}}: \ov{\eta^{I}}\to \Delta(\ov \eta)$. 
     
 \begin{rem}    \label{role-sp} 
      Le rôle de $\on{\mf{sp}}$ est le suivant. 
     Les faisceaux $\mc H _{ N, I, W}^{0,\leq\mu,E}$  sont lisses sur des ouverts de $X^{I}$ ne contenant pas forcément $\Delta(\eta)$ et c'est  donc  leur fibre  en $\ov{\eta^{I}}$ que l'on  étudie. D'un autre côté on construira certains sous-faisceaux qui se prolongent en des faisceaux lisses sur un ouvert de $X^{I}$ contenant $\Delta(\eta)$ et pour ces sous-faisceaux il est très important de considérer la fibre en $\Delta(\ov\eta)$ du prolongement car elle est plus canonique (par exemple elle est compatible avec l'action du groupe $\mf S(I)$ 
     des permutations de $I$ et avec la coalescence des pattes). Le rôle de      $\on{\mf{sp}}$ est d'identifier {\it canoniquement}  la fibre en $\ov{\eta^{I}}$  de ces sous-faisceaux avec la fibre 
      en $\Delta(\ov\eta)$ de leur prolongement.  
            L'idée de considérer la fibre  en un point géométrique 
de la diagonale apparaît dans le paragraphe 3.1 de \cite{ngo-modif-sym}.  \end{rem}

 Dans la suite on emploiera des lettres gothiques   $\mf E$, $\mf G$, $\mf M$ pour les $\mc O_{E}$-modules  ou  les $\mc O_{E}$-faisceaux et des lettres calligraphiées  $\mc E$, $\mc G$, $\mc M$ pour les $E$-espaces vectoriels ou les  $E$-faisceaux.

   On rappelle que pour toute partie $J\subset I$ on note 
  $\Frob_{J}:X^{I}\to X^{I}$ le morphisme qui à 
$(x_i)_{i\in I}$ associe $(x'_i)_{i\in I}$ avec 
\begin{gather}\label{defi-Frob-J}x_{i}'=\Frob (x_{i})\text{ \  si \  }i\in J\text{ \  et  \  }x_{i}'=x_{i}\text{ \  sinon.} \end{gather}

  \begin{lem} (Drinfeld, théorème 2.1 de \cite{drinfeld78} et proposition 6.1 de \cite{drinfeld-compact})   \label{lem-Frob-partiels-drinfeld}
  Soit  $U$ un ouvert dense de  $X$. On a une   équivalence entre 
    \begin{itemize}
  \item la catégorie $\mc C(U,I,\mc O_E)$ des  $\mc O_{E}$-faisceaux lisses (constructibles)   $\mf E$ sur  $U^{I}$, munis d'une action des morphismes de  Frobenius partiels, c'est-à-dire  d'isomorphismes 
  $F_{\{i\}}: \Frob_{\{i\}}^{*}(\mf E)\simeq \mf E$ commutant entre eux et dont la composée est l'isomorphisme  naturel 
  $\Frob_{U^{I}}^{*}(\mf E)\simeq \mf E$, 
  \item la  catégorie  des   représentations continues de   
  $\pi_{1}(U,\ov{\eta})^{I}$ sur des $\mc O_E$-modules de type fini,   
  \end{itemize}
  qui est caractérisée de manière unique par les deux  faits suivants
  \begin{itemize}
  \item la composée avec le foncteur de restriction des représentations de  $\pi_{1}(U,\ov{\eta})^{I}$  à $\pi_{1}(U,\ov{\eta})$  plongé diagonalement 
    est    le foncteur $\mf E\to \mf E_{\Delta(\ov{\eta})}$,   
  \item  si $(\mf F_{i})_{i\in I}$ est une famille de $\mc O_{E}$-faisceaux lisses (constructibles)     sur  $U$, 
  l'image par ce foncteur  de $\boxtimes_{i\in I}\mf F_{i}$ (muni de l'action naturelle des morphismes de Frobenius partiels) est 
  $(\boxtimes_{i\in I}\mf F_{i})_{\Delta(\ov{\eta})}  =\otimes _{i\in I} (\mf F_{i})_{ \ov{\eta} }$, muni de l'action de $\pi_{1}(U,\ov{\eta})^{I}$  venant du fait que chaque $(\mf F_{i})_{ \ov{\eta} }$ est muni d'une action de $\pi_{1}(U,\ov{\eta})$. 
  \end{itemize}
       \end{lem}
     
   \begin{rem} Cette équivalence de catégories est bien caractérisée de manière unique par les deux conditions à la fin du lemme car l'énoncé du lemme implique que tout objet  de la catégorie $\mc C(U,I,\mc O_E)$ est un quotient d'un objet de la forme   $\boxtimes_{i\in I}\mf F_{i}$ comme dans le lemme. En effet  c'est vrai dans la catégorie des   représentations continues de   
  $\pi_{1}(U,\ov{\eta})^{I}$ sur des $\mc O_E$-modules de type fini
  (cela est très facile à montrer pour des $\mc O_E$-modules de type fini de torsion et ce cas suffit pour que l'équivalence de catégories soit caractérisée de fa\c con unique par les deux conditions à la fin du lemme).    
   
   On peut décrire  le foncteur inverse de fa\c con explicite: 
   à une représentation  continue  de  
  $\pi_{1}(U,\ov{\eta})^{I}$ sur  un $\mc O_E$-module  de type fini il associe 
  un $\mc O_{E}$-faisceau  lisse  sur     $U^{I}$   grâce au   morphisme évident $\pi_{1}(U^{I},\Delta(\ov{\eta}))\to \pi_{1}(U,\ov{\eta})^{I}$ et l'action des      morphismes de  Frobenius partiels sur ce dernier est définie à l'aide de la construction suivante: 
  si on se donne des revêtements étales galoisiens $U_{i}$ de $U$ et  un ensemble fini $A$   muni d'une action de $\prod _{i\in I}\on{Gal}(U_{i}/U)$, alors 
  $(\prod_{i\in I}U_{i})\times_{\prod _{i\in I}\on{Gal}(U_{i}/U)} A$  est muni des morphismes de Frobenius partiels $\Frob_{U_{i}}\times \prod_{j\neq i }\Id_{U_{j}}\times \Id_{A}$. 
  
  Dans le lemme ci-dessus   on a présenté l'équivalence de catégories dans le sens où on l'utilise dans cet article mais, comme on vient de le voir,  c'est le foncteur inverse qui est construit simplement. Dans les lemmes ci-dessous on donnera des équivalences de catégorie dans l'autre sens. 
\end{rem}

     \begin{rem}
Dans le cadre du lemme précédent, 
    l'homomorphisme de spécialisation   $\mf{sp}^{*}:\restr{\mf E}{\Delta(\ov{\eta})}\to \restr{\mf E}{\ov{\eta^{I}}}$ fournit  un isomorphisme  entre les foncteurs  $\mf E\mapsto \restr{\mf E}{\Delta(\ov{\eta})}$ et  $\mf E\mapsto \restr{\mf E}{\ov{\eta^{I}}}$  de la  catégorie $\mc C(U,I,\mc O_E)$ vers la catégorie des $\mc O_{E}$-modules.
    \end{rem}
    
   \noindent{\bf Démonstration.}  Dans le cas   où $\sharp I=2$ l'énoncé figure dans le   théorème 2.1 de  
 \cite{drinfeld78}. La preuve est donnée dans  la proposition 6.1 de \cite{drinfeld-compact} (qui utilise la proposition 1.1 de \cite{Dr1}). Elle est également reprise dans  le  théorème 4 du paragraphe IV.2 de   
\cite{laurent-asterisque}. 
Le cas général est expliqué dans le théorème  8.1.4 de  \cite{eike-lau}, qui n'est pas publié. 

Pour la commodité du lecteur nous rappelons maintenant la preuve. 
La rédaction de la preuve donnée ci-dessous a été faite avec l'aide de Vladimir Drinfeld, Alain Genestier et  Gérard Laumon. 
On commence par le lemme suivant. 

\begin{lem} (Variante de la proposition 1.1 de \cite{Dr1}) Soit $k$ un corps séparablement clos contenant $\Fq$. 

1) Le foncteur qui à $V$ associe $W=V\otimes _{\Fq} k$ muni de $\varphi=\Id_{V}\otimes \Frob_{k}$ fournit une équivalence
\begin{itemize}
\item de la catégorie des $\Fq$-espaces vectoriels $V$ de dimension finie, 
\item vers la catégorie des couples $(W,\varphi)$ où $W$ est un $k$-espace vectoriel de dimension finie et $\varphi:(\Frob_{\on{Spec}k})^{*}( W)\to W$ est un isomorphisme
(si  $k$ est algébriquement clos,  $\varphi$ est simplement un isomorphisme $\Frob_{k}$-linéaire de $W$ dans $W$). 
 \end{itemize}

2) Soit $A $ une $\Fq$-algèbre. Alors le foncteur qui à $\mathcal M$ associe
$\mathcal M\otimes _{\Fq} k= \mathcal M\otimes _{A}(A\otimes k)$  induit une équivalence
 \begin{itemize}
\item de la catégorie des $A$-modules $\mathcal M$
\item vers la catégorie des  $A\otimes k$-modules  $\mc N$ munis d'un isomorphisme
$\varphi: (\Id_{\on{Spec}A}\times\Frob_{\on{Spec}k})^{*}(\mc N)\isom \mc N$ tel que tout élément de $\mathcal N$ soit inclus dans un sous-$k$-espace vectoriel de $\mathcal N$ de dimension finie stable par $\varphi$. 
\end{itemize}
\end{lem}
\noindent{\bf Démonstration (d'après \cite{Dr1}).} 
  Le 1) résulte de la surjectivité de l'isogénie de Lang pour $GL_{n}$ où $n=\dim V$. On a $V=W^{\varphi}$. Le 2) avec $A=\Fq$ découle de 1) en l'appliquant aux sous-$k$-espaces vectoriels de $\mathcal N$ de dimension finie stables par $\varphi$  (en effet la somme de deux sous-espaces de ce type est de ce type et 
  $\mathcal N$ est une réunion croissante de sous-espaces de ce type). 
  Le 2) avec $A$ arbitraire résulte du  2) avec $A=\Fq$  car l'action de $A$ se transmet par l'équivalence de catégories.    \cqfd
  
  \begin{rem}
  Voici un exemple de $(\mc N, \varphi)$  qui n'est pas dans l'image essentielle du foncteur de 2), bien que  $\mc N$ soit un module projectif de type fini sur $A\otimes k$,   parce que la condition à la fin de 2) n'est pas respectée: on prend $A=\Fq[t,t^{-1}]$
  (l'algèbre des fonctions sur $\mb{G}_{m}$),  $\mathcal N=A\otimes k=k[t,t^{-1}]$, 
  et $\varphi$ égal à la composée de $\Id_{A}\otimes \Frob_{k}$ et de la multiplication par $t$. 
  \end{rem}
  
    \begin{rem} Dans la situation de 2), si $\on{Spec} (A)$ est un ouvert affine $Y$ d'une variété projective $\ov Y$  sur $\Fq$ et $(\mc N,\varphi) $ vient d'un $\mc O$-module cohérent $\ov {\mc N}$
    sur $\ov Y \times \on{Spec} k$ muni d'un isomorphisme 
    $\ov{\varphi}:(\Id_{\ov Y}\times \Frob_{\on{Spec} k})^{*}(\ov {\mc N})\simeq \ov {\mc N}$, 
   , la condition à la fin de 2) est satisfaite (car, en notant $L$ un fibré ample sur $\ov Y$, les $H^{0}(\ov Y \times \on{Spec} k, \ov {\mc N}\otimes L^{n})$ sont de dimension finie et stables par l'action de Frobenius).  C'était la situation considérée par Drinfeld dans la proposition 1.1 de \cite{Dr1}
   (dont l'énoncé est que  la catégorie des 
   $(\ov {\mc N}, \ov{\varphi})$ comme ci-dessus est équivalente à celle des $\mc O$-modules cohérents  
    sur $\ov Y$).  
    \end{rem}
   
    \begin{lem} (généralisation de la   proposition 6.1 de \cite{drinfeld-compact}  et du  lemme 8.1.2 de \cite{eike-lau}) \label{Y-Z-T-genestier-normal}
 Soit $k$ un corps séparablement clos contenant $\Fq$. 
   Soit $Y$ un schéma  sur $\Fq$.   
   Alors le foncteur $Z\mapsto Z\times \on{Spec}k$ est une équivalence
   \begin{itemize}
   \item de la catégorie des revêtements étales  $Z\to Y$
   \item vers la catégorie des revêtements  étales 
   $T\to Y\times \on{Spec}k$
   munis d'un automorphisme   $\alpha:(\Id_{Y}\times \Frob_{\on{Spec}(k)})^{*}(T)\to T$ 

   \end{itemize}
    \end{lem}
    On appelle revêtement étale un morphisme   fini  étale. 
    \begin{rem}
    La conclusion est fausse si l'on omet l'hypothèse que les morphismes sont  étales. Voici un contre-exemple. On prend $Y=\mb{G}_{m}=\on{Spec}\Fq[t,t^{-1}]$, de sorte que $Y\times \on{Spec}k=\on{Spec}k[t,t^{-1}]$. 
    On prend alors $Z=\on{Spec}(k[t,t^{-1}])[\epsilon]/\epsilon^{2}$, et $\alpha$ tel que $\alpha^{*}(t)=t$ et $\alpha^{*}(\epsilon)=t\epsilon$ (et bien sûr $\alpha^{*}$ agit sur $k$ par $x\mapsto x^{q}$). 
        \end{rem}
   \dem Comme la catégorie des revêtements étales de $Y$ est équivalente à celle de son réduit, et de même pour $Y\times \on{Spec}k$,  il suffit de montrer le lemme pour $Y$ réduit. 
   En considérant les composantes irréductibles de $Y$ et leurs intersections, on voit qu'il suffit de montrer le lemme pour $Y$ irréductible.  En recollant il suffit de montrer le lemme pour $Y$ affine. On suppose donc $Y=\on{Spec}(A)$ avec $A$ une $\Fq$-algèbre intègre. 
   
   Il est clair que le foncteur est pleinement fidèle. Il reste à montrer qu'il est essentiellement surjectif. On se donne donc $T$ et $\alpha$ comme ci-dessus. 
   Les équations donnant $T$ et le morphisme $\alpha$ sont à coefficients dans une sous-algèbre de $A$ de type fini. On suppose donc que $A$  est une $\Fq$-algèbre intègre
    de type fini. 
   
 En langage géométrique la preuve consiste à compactifier $Y$, donc $Y\times \on{Spec} k$, puis à considérer sa normalisation à l'infini dans le corps des fractions de $T$. On va le faire de fa\c con complètement élémentaire.   
 
 On note $\mc B$ la $k$-algèbre des fonctions sur $T$. Grâce au 2) du lemme précédent il suffit de montrer que $\mc B$ est une réunion de sous-$k$-espaces vectoriels de dimension finie stables par $\alpha^{*}$. 
   Si $S$ est une partie finie de $A\sm \{0\}$, comme $\mc B\subset \mc B[S^{-1}]$ il suffit de montrer que $\mc B[S^{-1}]$ est une réunion de sous-$k$-espaces vectoriels de dimension finie stables par $\alpha^{*}$. Donc on peut remplacer $A$ par $A[S^{-1}]$. En choisissant $S$ convenablement on peut  supposer que $Y$ est lisse, ce qu'on fait désormais. Il suffit de montrer le résultat pour chaque composante connexe de $Z$. Or $Z$ est lisse puisque étale sur $Y\times \on{Spec}k$, et ses composantes connexes sont donc irréductibles. 
   
   Pour terminer la preuve il suffit donc de  montrer que si  $\mc B$ est une $k$-algèbre intègre contenant $A\otimes k$  et finie comme $A\otimes k$-module, telle que $\on{Frac}(\mc B)$ est une extension séparable de $\on{Frac}(A\otimes k)$, 
   et si $\varphi:\mc B\to \mc B$ est un morphisme relevant $\Id_{A}\otimes \Frob_{k}$ et induisant un isomorphisme $\mc B\otimes_{k,\Frob_{k}}k\to \mc B$, 
   alors $\mc B$ est une réunion de sous-$k$-espaces vectoriels de dimension finie stables par $\varphi$. 
     Le reste de la preuve consiste à montrer ce dernier énoncé. 
 
  Par le théorème de normalisation de Noether il existe $d\in \N$ et $x_{1},...,x_{d}$ algébriquement indépendants dans $A$ tels que $A$ soit finie sur la sous-algèbre $\Fq[x_{1},...,x_{d}]$, et que $\on{Frac}(A)$ soit une extension séparable de $\Fq(x_{1},...,x_{d})$ (voir par exemple le  corollaire 16.18 de \cite{eisenbud} ou le 
   théorème  4.2.2 de \cite{huneke}). 
   On est donc ramené à montrer l'énoncé avec $A=\Fq[x_{1},...,x_{d}]$. 
   Tout élément de $\mc B$ est entier sur $A\otimes k$ et 
   on possède $\on{Tr}:\mc B\to A\otimes k$ et $\beta_{1},..., \beta_{s}\in \mc B$ tels que 
   $\iota:b\mapsto (\on{Tr}(\beta_{i}b))_{i\in\{1,...,s\}}$ soit une injection de $(A\otimes k)$-modules de $\mc B$ dans $(A\otimes k)^{s}$.  On note $A_{n}$ le sous-$\Fq$-espace vectoriel de $A$ formé des polynômes de degré total $\leq n$.
      On note $\mc B_{n}$ le sous-$k$-espace vectoriel de $\mc B$ formé des éléments annulés par un polynôme de la forme 
      \begin{gather}\label{pol-annul}
      x^{k}+a_{1}x^{k-1}+ ...+ a_{k} \text{ avec } k\in \N^{*} \text{ et }a_{i}\in A_{in}\otimes k \text{ pour tout }i=1,...,k. \end{gather}
  Les formules classiques donnant à partir de deux polynômes $P$ et $Q$ de racines $X_{i}$ et $Y_{j}$ les polynômes de racines $X_{i}+Y_{j}$ et $X_{i}Y_{j}$ 
   montrent que $\mc B_{n}$ est un sous-$k$-espace vectoriel de $\mc B$ et que $\mc B_{m}\mc B_{n}\subset \mc B_{m+n}$.   De plus $\mc B=\bigcup_{n\in \N} \mc B_{n}$. On a $\on{Tr}(\mc B_{n})\subset A_{n}\otimes k$ car $\on{Tr}$ préserve la condition d'annulation par un polynôme de la forme 
   \eqref{pol-annul} et $A_{n}\otimes k$ est exactement l'ensemble des éléments de $A\otimes k$ vérifiant cette condition. Donc 
   si $m$ est un entier tel que $\beta_{1},..., \beta_{s}$ appartiennent à $\mc B_{m}$, 
   $\iota(\mc B_{n})\subset (A_{m+n}\otimes k)^{s}$. Donc $\mc B_{n}$ est un k-espace vectoriel  de dimension finie. Par construction il est stable par $\varphi$. Donc  on a fini. 
   \cqfd
 
 Le lemme précédent implique immédiatement le lemme suivant. 
 
  \begin{lem} \label{Y-Z-T}
  Soit  $Y$ un schéma  sur $\Fq$.  
  On note $F$ un  corps contenant $\Fq$   et $ F^{\mr{sep}}$ une clôture séparable  de $F$. 
   Le foncteur  
  \begin{itemize}
  \item de la catégorie des revêtements étales  $Z$ de $Y$, munis d'une action de 
 $\on{Gal}(F^{\mr{sep}}/F)$,  
 \item vers  la catégorie des revêtements étales
  $T\to   Y\times \on{Spec}(F)$  munis d'un automorphisme  $\alpha:(\Id_{Y}\times \Frob_{\on{Spec}(F)})^{*}(T)\to T$, 
 \end{itemize}
 qui à $Z$ muni d'une action de $\on{Gal}(F^{\mr{sep}}/F)$ associe 
 $(T,\alpha)$ avec 
 \begin{itemize}
 \item 
 $T$ égal au 
 quotient de $Z\times \on{Spec}(F^{\mr{sep}})$ par l'action diagonale de 
 $\on{Gal}(F^{\mr{sep}}/F)$, 
  \item  $\alpha$ donné par l'action de $\Id_{Z}\times \Frob_{\on{Spec}(F^{\mr{sep}})}$ sur ce quotient,  
  \end{itemize}
 est une équivalence de catégories.  
   \end{lem}
 \dem  
 On va construire un foncteur quasi-inverse qui à $(T,\alpha)$ associe $Z$ muni d'une action de $\on{Gal}(F^{\mr{sep}}/F)$. On voit que   $$\widetilde T=T\times _{ Y\times \on{Spec}(F)}(Y \times \on{Spec}(F^{\mr{sep}}))$$ est un revêtement étale de $Y \times \on{Spec}(F^{\mr{sep}})$ 
 muni d'un isomorphisme $$\widetilde \alpha:(\Id_{Y}\times \Frob_{\on{Spec}(F^{\mr{sep}})})^{*}(\widetilde T)\to \widetilde T.$$ Par le lemme précédent il existe $Z$ revêtement étale de $Y$ tel que $\widetilde T=Z\times \on{Spec} F^{\mr{sep}}$ et 
 $\widetilde \alpha=\Id_{Z}\times \Frob_{\on{Spec}  F^{\mr{sep}}}$. Comme $Z$ est égal à l'espace des coinvariants de  $\widetilde \alpha$ agissant sur $\widetilde T$ (c'est-à-dire au coégalisateur de $\Id$ et $\widetilde \alpha$), $\on{Gal}(F^{\mr{sep}}/F)$ agit naturellement sur $Z$. 
   \cqfd

    \begin{lem} (Drinfeld, théorème 2.1 de \cite{drinfeld78} et proposition 6.1 de \cite{drinfeld-compact}, et Eike Lau,  lemme 8.1.2 de \cite{eike-lau})   \label{lem-Frob-partiels-drinfeld-ensembles}
  Soient  $X_{1}$, ..., $X_{k}$ des schémas connexes sur $\Fq$.
Pour tout $i$, soit $\ov z_{i}$ un point géométrique de $X_{i}$. 
   On a une   équivalence entre 
    \begin{itemize}
     \item la  catégorie  des  actions  continues de   
  $\prod _{i=1}^{k} \pi_{1}(X_{i},\ov z_{i})$ sur des ensembles finis, 
  \item la catégorie des revêtements étales $T$ de $X_{1}\times ... \times X_{k}$, munis d'une action des morphismes de  Frobenius partiels, c'est-à-dire  de morphismes $F_{\{i\}}$ au-dessus de $\Frob_{X_{i}}\times \prod_{j\neq i}\Id_{X_{j}}$, commutant entre eux et dont la composée est  $\Frob_{T}$, 
    \end{itemize}
  où le foncteur  est le suivant: 
  si on se donne des revêtements étales galoisiens $U_{i}$ de $X_{i}$ et  un ensemble fini $A$   muni d'une action de $\prod _{i\in I}\on{Gal}(U_{i}/X_{i})$, alors 
  $(\prod_{i\in I}U_{i})\times_{\prod _{i\in I}\on{Gal}(U_{i}/X_{i})} A$,  muni des morphismes de Frobenius partiels $\Frob_{U_{i}}\times \prod_{j\neq i }\Id_{U_{j}}\times \Id_{A}$, est l'image par le foncteur de  l'action de $\prod _{i=1}^{k} \pi_{1}(X_{i},\ov z_{i})$ sur l'ensemble fini $\Big(\prod _{i=1}^{k} (U_{i})_{\ov z_{i}}\Big)\times_{\prod _{i\in I}\on{Gal}(U_{i}/X_{i})} A$. 
        \end{lem}
\dem Le cas où $k=2$ permet de montrer le cas général, par  récurrence sur $k$. On suppose donc que $k=2$.  Il est clair que le foncteur décrit dans l'énoncé est pleinement fidèle. Il reste donc à  montrer qu'il est essentiellement surjectif. 
On peut définir de fa\c con équivalente la première catégorie comme la catégorie des revêtements étales de $X_{2}$ munis d'une action continue de 
$\pi_{1}(X_{1},\ov z_{1})$. En recollant on peut supposer $X_{2}$ réduit, irréductible et affine. En échangeant les rôles de $X_{1}$ et $X_{2}$ on peut faire de même pour $X_{1}$. Désormais on suppose donc que 
 $X_{1}$ et $X_{2}$  sont  réduits, irréductibles et affines.

On se donne donc un revêtement étale 
$T$ de $X_{1}\times  X_{2}$, muni   de morphismes $F_{\{1\}}$ et  $F_{\{2\}}$ au-dessus de $\Frob_{X_{1}}\times  \Id_{X_{2}}$ et $\Id_{X_{1}}\times  \Frob_{X_{2}}$, commutant entre eux et dont la composée est  $\Frob_{T}$ et on veut montrer que $(T, F_{\{1\}}, F_{\{2\}})$  provient d'une action de $\pi_{1}(X_{1},\ov z_{1})\times 
\pi_{1}(X_{2},\ov z_{2})$ sur un ensemble fini. Comme les équations de $T$ et les morphismes $F_{\{1\}}$ et  $F_{\{2\}}$  s'expriment à l'aide d'un nombre fini d'éléments des algèbres de fonctions sur $X_{1}$ et $X_{2}$, on peut supposer $X_{1}$ et  $X_{2}$ de type fini, ce qu'on fait désormais. 

   En appliquant le  \lemref{Y-Z-T-genestier-normal}  à $X_{1}\times \ov z_{2}$ et        $\ov z_{1}\times X_{2}$ on voit que    la fibre $T_{\ov z_{1}\times \ov z_{2}}$ a un sens  (bien que $\ov z_{1}\times \ov z_{2}$ ne soit pas un point géométrique de $X_{1}\times X_{2}$) et qu'elle est munie d'actions de $\pi_{1}(X_{1},\ov z_{1})$ et $
\pi_{1}(X_{2},\ov z_{2})$. On veut montrer qu'elles commutent et que $T$ provient de cette action  de $\pi_{1}(X_{1},\ov z_{1})\times 
\pi_{1}(X_{2},\ov z_{2})$ sur $T_{\ov z_{1}\times \ov z_{2}}$. 
Si $\ov x$ et $\ov y$ sont deux points géométriques d'un schéma connexe $X$ on rappelle
qu'un chemin de $\ov x$ vers $\ov y$ est un isomorphisme entre les foncteurs fibres en 
$\ov x$ et $\ov y$ de la catégorie des revêtements étales de $X$ vers la catégorie des ensembles finis. Soient, pour $i=1,2$,  $\ov x_{i}$ et $\ov y_{i}$ des points géométriques de $X_{i}$ et $\gamma_{i}$ un chemin de $\ov x_{i}$ vers $\ov y_{i}$. On alors un diagramme
  $$ \xymatrixcolsep{3pc}
   \xymatrix{  
      T_{\ov x_{1}\times \ov x_{2}} \ar[r]^-{\gamma_{1}} \ar[d]^{\gamma_{2}} & 
       T_{\ov y_{1}\times \ov x_{2}} \ar[d]^{\gamma_{2}}
       \\   T_{\ov x_{1}\times \ov y_{2}} \ar[r]^-{\gamma_{1}} & 
        T_{\ov y_{1}\times \ov y_{2}}
       } $$   
       où les flèches sont obtenues en appliquant le  \lemref{Y-Z-T-genestier-normal}   
       à $X_{1}\times \ov x_{2}$, $X_{1}\times \ov y_{2}$, 
       $\ov x_{1}\times X_{2}$, $\ov y_{1}\times X_{2}$. Il reste à montrer que ce diagramme est commutatif. 
       On note $X_{2}^{\nu}$   le normalisé de $X_{2}$. Comme  $X_{2}$ est affine, intègre et  de type fini, le morphisme de $X_{2}^{\nu}$     vers $X_{2}$ est fini (corollaire 13.13 de \cite{eisenbud}) et évidemment surjectif.  
  Il suffit de montrer l'énoncé en supposant de plus que $\gamma_{2}$ se relève en un chemin sur $X_{2}^{\nu}$. En effet parmi les chemins de 
  $\ov x_{2}$ vers $\ov y_{2}$ (avec la topologie profinie)  les chemins pouvant être coupés en chemins qui se relèvent à 
       $X_{2}^{\nu}$ sont denses. 
       Cela est équivalent au fait que le morphisme $X_{2}^{\nu}\to X_{2} $ est de descente pour la catégorie fibrée des schémas étales sur des $X_{2}$-schémas, d'après le  corollaire 3.3 de l'exposé IX de \cite{grothendieck-sga1} (en fait la descente est même effective d'après le théorème 4.7 de loc. cit. mais on n'en a pas besoin).

      Par conséquent il suffit de montrer la commutativité du diagramme ci-dessus lorsque $X_{2}$ est intègre,  normal et de type fini. 
       Dans ce cas, en notant  $\ov {\eta_{X_{2}}}$ 
      un point géométrique  au-dessus du point générique $\eta_{X_{2}}$, on sait que 
   $\pi_{1}(X_{2}, \ov {\eta_{X_{2}}})$ est un quotient du groupe de Galois du corps de fractions $\pi_{1}(\eta_{X_{2}}, \ov {\eta_{X_{2}}})$. Donc le \lemref{Y-Z-T} implique que   $(T, F_{\{1\}}, F_{\{2\}})$ provient  d'une action de
$\pi_{1}(X_{1},\ov z_{1})\times 
\pi_{1}(X_{2},\ov z_{2})$    sur un ensemble fini 
(car un revêtement étale de $X_{1}
\times X_{2}$ est déterminé par sa restriction à $X_{1}\times \eta_{X_{2}}$). On voit alors que 
le diagramme ci-dessus est commutatif.   
\cqfd

 \noindent{\bf Fin de la démonstration du \lemref{lem-Frob-partiels-drinfeld}}
 On fixe une bijection $I=\{1,...,k\}$ et on applique le lemme précédent avec $X_{i}=U$ pour tout $i$. 
  \cqfd
 
   \begin{lem} \label{lem-U^I} (Drinfeld,   proposition 6.1 de \cite{drinfeld-compact}  et Eike Lau,  lemme 9.2.1 de \cite{eike-lau}) 
    Soit  $\Omega$ un ouvert dense de $X^{I}$ et 
  $ \mf E$ un $\mc O_{E}$-faisceau lisse 
  (constructible)  sur $\Omega$ muni 
   d'une action des morphismes de Frobenius partiels, c'est-à-dire d'isomorphismes 
  $F_{\{i\}}: \restr{\Frob_{\{i\}}^{*}(\mf E)}{\eta^{I}}\simeq \restr{\mf E}{\eta^{I}}$ commutant entre eux et dont la composée est l'isomorphisme  naturel 
  $\Frob_{\eta^{I}}^{*}(\mf E)\simeq \mf E$. Alors il existe un ouvert dense $U$ de $X$ tel que $\mf E$ s'étende en un $\mc O_{E}$-faisceau lisse  sur $U^{I}$. 
   \end{lem}
  \noindent{\bf Démonstration.}
     Soit $\Omegasour$ le plus grand ouvert de $X^{I}$ tel que $\mf E$ s'étende en un $\mc O_{E}$-faisceau lisse sur $\Omegasour$. 
    Alors  le complémentaire de $\Omegasour$ est un fermé strict de $X^{I}$, invariant par l'action des morphismes de Frobenius partiels $\Frob_{\{i\}}$ de $X^{I}$. D'après le 
    lemme 9.2.1 de \cite{eike-lau} ce fermé strict est inclus dans une réunion finie de diviseurs verticaux (c'est-à-dire  des images inverses d'un point fermé par l'une des projections  $X^{I}\to X$). Pour la commodité du lecteur on rappelle la preuve de ce lemme de Eike Lau. 
    Soit $Z\subset X^{I}$ un fermé irréductible de dimension $k<\sharp I$. 
    On considère les degrés de ses projections vers $X^{J}$ où $J$ parcourt l'ensemble des  parties de $I$ de cardinal $k$. Les morphismes de Frobenius partiels multiplient  ces degrés par des puissances de $q$, d'où l'on déduit  que si deux au moins sont non nuls alors l'orbite de $Z$ suivant les morphismes de Frobenius partiels est infinie. Mais si un seul de ces degrés est non nul, correspondant à une partie $J\subset I$,  
    alors $Z$ est vertical relativement à toutes les projections 
    vers les coordonnées dans $I\sm J$. 
    Il  existe donc un ouvert dense $U\subset X$ tel que $\Omegasour\supset U^{I}$.  
      \cqfd

\begin{cor}\label{equiv-cat-U-I-Omega}
L'équivalence de catégories du \lemref{lem-Frob-partiels-drinfeld} fournit un foncteur (dépendant seulement du choix de $\ov\eta$) 
\begin{itemize}
\item de la catégorie des  $\mc O_{E}$-faisceaux constructibles $\mf E$ sur $\eta^{I}$, munis  d'une action des morphismes de Frobenius partiels, et admettant un prolongement à un certain ouvert dense $\Omega$ de $X^{I}$
\item vers la catégorie des représentation continues de $\pi_{1}(\eta,\ov{\eta})^{I}$ dans un $\mc O_{E}$-module de type fini, se factorisant par $\pi_{1}(U, \ov{\eta})^{I}$ pour un certain ouvert dense $U$ de $X$, \end{itemize}
   et ce foncteur est une  équivalence de catégories.     
     \end{cor}
\noindent{\bf Démonstration.} Cela résulte des 
 lemmes \ref{lem-Frob-partiels-drinfeld} et \ref{lem-U^I}. Plus précisément le \lemref{lem-U^I} implique que   $\mf E$ se prolonge en un faisceau lisse 
 $\wt {\mf E}$ sur un ouvert de la forme $U^{I}$ et le  \lemref{lem-Frob-partiels-drinfeld} fournit une action de $\pi_{1}(U,\ov{\eta})^{I}$  sur 
 $\restr{\wt {\mf E}}{\Delta(\ov\eta)}$.   \cqfd

 \begin{rem} \label{rem-action-Gal^I-eta^I}
 Dans les notations de la démonstration précédente, 
 le choix de $\mf{sp}$ fournit un isomorphisme $\restr{\wt {\mf E}}{\Delta(\ov\eta)}
 \isor{ \mf{sp} ^{*}}
 \restr{\wt {\mf E}}{\ov{\eta^{I}}}= \restr{ \mf E}{\ov{\eta^{I}}}$, donc 
 $\restr{ \mf E}{\ov{\eta^{I}}}$ est muni d'une action de $\pi_{1}(\eta,\ov{\eta})^{I}$ dépendant du choix de $\mf{sp}$. 
\end{rem}
 
 Le lemme suivant sera utilisé ultérieurement (dans la preuve du \lemref{S-non-ram-prelim2}).

    \begin{lem}\label{S-non-ram-prelim}
      Soit  $v\in |X|$. Soit $I$ un ensemble fini et $\iota\in I$ un élément. 
      On fixe  un plongement $\ov F\subset \ov {F_{v}}$. Soit  $d\in \N$,  et 
    $\gamma\in \on{Gal}(\ov {F_{v}}/F_{v})\subset  \on{Gal}(\ov F/F)$ tel que $\deg(\gamma)=d$.  
    On définit $(\gamma_{i})_{i\in I}\in\on{Gal}(\ov F/F)^{I}$ en posant $\gamma_{\iota}=\gamma$ et $\gamma_{i}=1$ pour $i\neq \iota$. 
    On fixe un  point géométrique  $\ov v$ au-dessus de  $v$ et une  flèche de spécialisation 
    $\on{\mf{sp}}_{v}:\ov \eta\to \ov v$ associés au choix du plongement $\ov F\subset \ov {F_{v}}$. 
    Plus précisément $ \ov v$ est le spectre du corps résiduel de l'extension maximale non ramifiée de $F_{v}$ dans $\ov F_{v}$ et $\on{\mf{sp}}_{v}$ vient de l'inclusion de l'hensélisé strict de $F$ en $\ov v$ dans $\ov F\subset \ov {F_{v}}$. 
    On note encore   $\on{\mf{sp}}_{v}$ la  flèche de spécialisation de  
    $\Delta(\ov \eta)$ vers   $\Delta(\ov v)$ égale à l'image par $\Delta$ de cette dernière.  
    Soit $\mf E$ un faisceau comme dans le \lemref{lem-Frob-partiels-drinfeld}, c'est-à-dire lisse sur un ouvert de la forme $U^{I}$ et muni d'actions des morphismes de Frobenius partiels. On note $j^{I}:U^{I}\to X^{I}$ l'inclusion.   Alors on a la commutativité du  diagramme  
       $$ \xymatrixcolsep{1pc} \xymatrix{  
        \restr{(j^{I})_{*}\mf E}{\Delta(\ov v)} \ar[r]^-{\mf{sp}_{v}^{*}} \ar[d]^{F_{\{\iota\}}^{\deg(v)d}} & 
       \restr{\mf E}{\Delta(\ov\eta)} \ar[d]^{(\gamma_{i})_{i\in I}}
       \\  \restr{(j^{I})_{*}\mf E}{\Delta(\ov v)} \ar[r]^-{\mf{sp}_{v}^{*}} & 
        \restr{\mf E}{\Delta(\ov\eta)}  
       } $$   
       où la flèche verticale de droite est l'action de $\pi_{1}(U,\ov\eta)^{I}$ sur 
       $  \restr{\mf E}{\Delta(\ov\eta)}$ donnée par le \lemref{lem-Frob-partiels-drinfeld}. 
                \end{lem}    
  
   \begin{rem}  Dans l'énoncé précédent on ne suppose pas que $v$ appartient à $U$. \end{rem}
  
  \dem 
   Il suffit de le montrer avec  $\mf E$   de la forme 
   $\boxtimes_{i\in I}\mf E_{i}$ (comme dans l'énoncé du \lemref{lem-Frob-partiels-drinfeld}). Alors, en notant $j:U\to X$ l'inclusion, on a  
   $\restr{(j^{I})_{*}\mf E}{\Delta(\ov v)}=\bigotimes _{i\in I} \big(\restr{j_{*}\mf E_{i}}{\ov v}\big). $ On est donc ramené à montrer le lemme dans le cas où $I$ est un singleton, et alors cela résulte de la définition de $\deg: \on{Gal}(\ov {F_{v}}/F_{v})\to 
   \on{Gal}(\ov {k(v)}/k(v))$ par restriction de l'action sur l'extension non ramifiée maximale  et sur son corps résiduel. 
   \cqfd

  \begin{rem}  
      Le foncteur $\mc F\mapsto \mc F _{\ov{\eta^{I}}}$ induit une équivalence entre  la catégorie des $\mc O_{E}$-faisceaux 
      ({\it resp.} $E$-faisceaux) 
      constructibles sur $\eta^{I}$ et la catégorie des représentations continues de $\pi_{1}(\eta^{I},\ov{\eta^{I}})$ sur des $\mc O_{E}$-modules de type fini  
   ({\it resp.} $E$-espaces vectoriels de dimension finie). 
  L'image inverse par $\Frob_{\{i\}}$ est une auto-équivalence de la  catégorie des $\mc O_{E}$-faisceaux  ({\it resp.} $E$-faisceaux) constructibles sur $\eta^{I}$. 
  Si $\mc F$ est muni d'une action des morphismes de Frobenius partiels, et si $\mc M\subset \mc F _{\ov{\eta^{I}}}$  est une sous-$\pi_{1}(\eta^{I},\ov{\eta^{I}})$-représentation, et $\mc G$ le sous-faisceau 
   de  $\mc F$ sur $\eta^{I}$ tel que 
  $\mc M=
  \mc G_{\ov{\eta^{I}}}$, on dira que $\mc M$  est stable par  les morphismes de Frobenius partiels si c'est le cas de $\mc G$. 
  \end{rem}

  \begin{rem}\label{monoide}
  On va définir un groupe    $\on{FWeil}(\eta^{I},\ov{\eta^{I}})$ qui  
       \begin{itemize}
     \item est une  extension de $\Z^{I}$ par   $\on{Ker}(\pi_{1}(\eta^{I},\ov{\eta^{I}})\to \wh \Z)$, 
     \item lorsque  $I$ est un   singleton, 
 s'identifie au groupe de   Weil   usuel 
   $$\on{Weil}(\eta,\ov{\eta})=\pi_{1}(\eta,\ov{\eta})\times_{\wh \Z}\Z,   $$   
    \item lorsque   
     $\sharp I=2$,    est le groupe noté  $\on{FGal}(\ov K/K)$ dans le paragraphe  4.1 de \cite{drinfeld78}, $\wt W_{\Lambda}$ dans le  théorème 7 de 
     \cite{corvallis-kazhdan} et  $\Z W_{F^{2}}^{\eta}$ dans le lemme  VI.13 de \cite{laurent-inventiones}.  
     \end{itemize}
    On définit 
          \begin{gather*}\on{FWeil}(\eta^{I},\ov{\eta^{I}})=
       \big\{\varepsilon \in \on{Aut}_{\ov\Fq}(\ov{F^{I}}), \exists (n_{i})_{i\in I}\in \Z^{I}, \restr{\varepsilon}{(F^{I})^{\mr{perf}}}=\prod_{i\in I}(\Frob_{\{i\}})^{n_{i}}\big\} 
  .\end{gather*}
 Comme $
  \ov{\eta} \times_{\on{Spec}\Fqbar} \cdots\times_{\on{Spec}\Fqbar} \ov{\eta}$, muni du morphisme diagonal $$\Delta(\ov \eta)\to
  \ov{\eta} \times_{\on{Spec}\Fqbar} \cdots\times_{\on{Spec}\Fqbar} \ov{\eta},$$ est (à la perfectisation près) une limite projective de voisinages étales 
  $\Delta(\ov \eta)$-ponctués de $\Delta( \eta)$ dans $X^{I}$, 
  $\mf{sp}$ fournit 
  une inclusion  $$\ov{F}\otimes_{\ov\Fq} \cdots 
 \otimes_{\ov\Fq} \ov{F} \subset \ov{F^{I}} 
  .$$ 
  En effet $\ov{F}\otimes_{\ov\Fq} \cdots 
 \otimes_{\ov\Fq} \ov{F}$ est un anneau intègre: si $Y=\mr{Spec} A$ est une courbe
 affine irréductible sur $\ov\Fq$ munie d'un morphisme quasi-fini vers $X_{\ov\Fq}$, alors $A\otimes_{\ov\Fq} \cdots 
 \otimes_{\ov\Fq} A$ est intègre car c'est l'anneau des fonctions sur la variété irréductible $Y\times _{\ov\Fq} \cdots 
 \times_{\ov\Fq} Y$, et $\ov{F}\otimes_{\ov\Fq} \cdots 
 \otimes_{\ov\Fq} \ov{F}$ est la limite inductive de tels anneaux. 
  
  Par restriction des automorphismes, on  obtient donc  un 
   morphisme  surjectif 
    \begin{gather} \label{mor-surj-Fweil-Weil^I}
    \on{FWeil}(\eta^{I},\ov{\eta^{I}})\to \big(\on{Weil}(\eta,\ov{\eta})\big)^{I}\end{gather} (dépendant du choix de  $\mf{sp}$). 
   Il rend plus  explicite l'équivalence de catégories du \corref{equiv-cat-U-I-Omega}: 
   bien que ce morphisme ne soit pas injectif lorsque $\sharp I >1$, les deux groupes ont le même complété profini $\pi_{1}(\eta,\ov\eta)^{I}$    (à strictement parler cette  assertion est équivalente au corollaire~\ref{equiv-cat-U-I-Omega} pour les faisceaux de torsion  ou si on néglige les  conditions avec  $U$ et $\Omega$). 
    \end{rem}

Soit  $W$ une  représentation de dimension finie de 
$(\wh G)^{I}$.  
On a une action naturelle de $\on{FWeil}(\eta^{I},\ov{\eta^{I}})$ sur  $ \varinjlim _{\mu}\restr{\mc H _{ N, I, W}^{0,\leq\mu,E}}{\ov{\eta^{I}}}$ qui réunit  l'action de $\pi_1(\eta^{I},\ov{\eta^{I}})$ et celle des 
 morphismes de  Frobenius partiels. 
Cependant on ne peut pas appliquer le lemme de  Drinfeld à cet espace vectoriel parce qu'il est de  dimension infinie et d'autre part 
on ne peut pas appliquer le lemme de  Drinfeld à  $\mc H _{ N, I, W}^{0,\leq\mu,E}$ parce que l'action des morphismes de  Frobenius partiels augmente  $\mu$.

Il y a une équivalence entre 
\begin{itemize}
\item les sous-$\mc O_{E}$-faisceaux constructibles de 
$\varinjlim _{\mu}\restr{\mc H _{ N, I, W}^{0,\leq\mu,E}}{ \eta^{I}} $ stables par l'action des morphismes de Frobenius partiels, 
\item les sous-$\mc O_{E}$-modules de type fini 
de $ \varinjlim _{\mu}\restr{\mc H _{ N, I, W}^{0,\leq\mu,E}}{\ov{\eta^{I}}}$ stabilisés par l'action de $\on{FWeil}(\eta^{I}, \ov{\eta^{I}})$
 \end{itemize}
 et à de tels objets on peut appliquer le lemme de Drinfeld, c'est-à-dire que 
 l'action de $\on{FWeil}(\eta^{I},\ov{\eta^{I}})$ se factorise à travers  
     $\pi_{1}(\eta,\ov\eta)^{I}$.

    Pour pouvoir appliquer le lemme de Drinfeld,  on va    définir un sous-espace  $\Big( \varinjlim _{\mu}\restr{\mc H _{ N, I, W}^{0,\leq\mu,E}}{\ov{\eta^{I}}}\Big)^{\mr{Hf}}$ 
  et montrer qu'il est  une limite inductive   de sous-$\mc O_{E}$-modules de type fini stabilisés  par $\on{FWeil}(\eta^{I},\ov{\eta^{I}})$. Par le lemme de  Drinfeld,  l'action de $\on{FWeil}(\eta^{I},\ov{\eta^{I}})$ sur $\Big( \varinjlim _{\mu}\restr{\mc H _{ N, I, W}^{0,\leq\mu,E}}{\ov{\eta^{I}}}\Big)^{\mr{Hf}}$ se factorisera à travers  
     $\pi_{1}(\eta,\ov\eta)^{I}$.

\begin{defi}\label{hecke-fini}
Soit  $\ov x$ un  point géométrique  de $(X\sm N)^{I}$. 
On dit qu'un élément de   $\varinjlim _{\mu}\restr{\mc H _{ N, I, W}^{0,\leq\mu,E}}{\ov{x}}$
est  Hecke-fini s'il satisfait l'une des  conditions  équivalentes suivantes 
\begin{itemize}
\item  il appartient à  un  sous-$\mc O_{E}$-module $\mf M$ de type fini de  $\varinjlim _{\mu}\restr{\mc H _{ N, I, W}^{0,\leq\mu,E}}{\ov{x}}$ qui est 
stable par  $T(f)$ pour tout   $f\in C_{c}(K_{N}\backslash G(\mb A)/K_{N},\mc O_{E})$, 
\item il vérifie la même condition, et en plus $\mf M$ est   stable par  $\pi_{1}(x,\ov{x})$. 
\end{itemize}
On note  $\Big( \varinjlim _{\mu}\restr{\mc H _{ N, I, W}^{0,\leq\mu,E}}{\ov{x}}\Big)^{\mr{Hf}}$ l'ensemble de tous les  éléments Hecke-finis. 
C'est un  sous-$E$-espace vectoriel  de  $ \varinjlim _{\mu}\restr{\mc H _{ N, I, W}^{0,\leq\mu,E}}{\ov{x}}$ qui est  stable 
 par $\pi_{1}(x,\ov{x})$. 
\end{defi}
 
 \begin{rem} En pratique on appliquera cette  définition en des points  géométriques $\ov x$ dont  l'image sur chaque copie de $X\sm N$ est générique, et dans ce cas 
 la définition ne nécessite pas de connaître l'extension des  opérateurs de Hecke $T(f)$ à $(X\sm N)^{I}$ tout entier construite dans le \corref{cor-hecke-etendus-compo}. 
 \end{rem}
 
On remarque que  $\Big( \varinjlim _{\mu}\restr{\mc H _{ N, I, W}^{0,\leq\mu,E}}{\ov{x}}\Big)^{\mr{Hf}}$  est stable par  $C_{c}(K_{N}\backslash G(\mb A)/K_{N},E)$. 

Comme les  morphismes de création et d'annihilation commutent avec les  opérateurs de Hecke, ils  preservent les sous-espaces   Hecke-finis. 

Si  $\ov x$ et  $\ov y$ sont deux point géométriques   de $(X\sm N)^{I}$ et   $\on{\mf{sp}}:\ov x\to \ov y$ est une flèche de spécialisation,   l'homomorphisme de spécialisation 
$$\on{\mf{sp}}^{*}: \varinjlim _{\mu}\restr{\mc H _{ N, I, W}^{0,\leq\mu,E}}{\ov y}\to \varinjlim _{\mu}\restr{\mc H _{ N, I, W}^{0,\leq\mu,E}}{\ov x}$$ 
envoie  
$\Big( \varinjlim _{\mu}\restr{\mc H _{ N, I, W}^{0,\leq\mu,E}}{\ov y}\Big)^{\mr{Hf}}$ dans  $\Big(  \varinjlim _{\mu}\restr{\mc H _{ N, I, W}^{0,\leq\mu,E}}{\ov x}\Big)^{\mr{Hf}}$ (parce que  les opérateurs de Hecke sont des morphismes de faisceaux donc commutent avec $\on{\mf{sp}}^{*}$). 

Les sous-espaces  Hecke-finis sont  également stables sous l'action  des  morphismes de  Frobenius partiels: on a  
$$F_{\{i\}} \Big(\Big(  \varinjlim _{\mu}\restr{\mc H _{ N, I, W}^{0,\leq\mu,E}}{\Frob_{\{i\}}(\ov x)}\Big)^{\mr{Hf}}\Big)\subset \Big(  \varinjlim _{\mu}\restr{\mc H _{ N, I, W}^{0,\leq\mu,E}}{\ov x}\Big)^{\mr{Hf}}.$$ 
Ainsi, dans le cas particulier    où  $\ov x=\ov{\eta^{I}}$, $$\Big( \varinjlim _{\mu}\restr{\mc H _{ N, I, W}^{0,\leq\mu,E}}{\ov{\eta^{I}}}\Big)^{\mr{Hf}}\subset \varinjlim _{\mu}\restr{\mc H _{ N, I, W}^{0,\leq\mu,E}}{\ov{\eta^{I}}}$$ 
  est une  sous-représentation de $\on{FWeil}(\eta^{I},\ov{\eta^{I}})$.   

On commence par étudier le cas où $I=\emptyset$ et $W=\mbf  1$. D'après la \propref{prop-chtoucas}~d) (extraite de la  
proposition 2.16 de \cite{var})  $\Cht_{N,\emptyset,{\mbf  1} }^{(\{0\})}/\Xi$ est égal au champ constant 
 $\Bun_{G,N}(\Fq)$ sur  $\eta^{\emptyset}=\on{Spec}(\Fq)$.  Donc  
$$\varinjlim _{\mu}\restr{\mc H _{ N, \emptyset, \mbf  1}^{0,\leq\mu,E}}{\Fqbar}=C_{c}(\Bun_{G,N}(\Fq)/\Xi,E),  $$
où 
la  restriction à $\Fqbar$ dans le membre de gauche a un sens car  $\mc H _{ N, \emptyset, \mbf  1}^{0,\leq\mu,E}$ est un faisceau (au demeurant trivial) sur $\eta^{\emptyset}=\Fq$.

 Le cas où  $I$ est un  singleton et $W=\mbf  1$ est essentiellement le même  car  on a  un 
 isomorphisme 
  \begin{gather}\label{egalite-cht-singl-1}\Cht_{N,\{0\},{\mbf  1} }^{(\{0\})}/\Xi=\Big(\Cht_{N,\emptyset,{\mbf  1} }^{(\{0\})}/\Xi\Big)\times _{\Fq}(X\sm N)\end{gather} 
  (qui est à l'origine de l'isomorphisme  de  coalescence associé à $\zeta_{\emptyset}: \emptyset\to \{0\}$). Par conséquent 
$$\varinjlim _{\mu}\restr{\mc H _{ N, \{0\}, \mbf  1}^{0,\leq\mu,E}}{\ov\eta}=C_{c}(\Bun_{G,N}(\Fq)/\Xi,E) $$
et le  faisceau   $\mc H _{ N, \{0\}, \mbf  1}^{0,\leq\mu,E}$ est en fait constant sur  $X\sm N$.

 \begin{rem}\label{quotient-adelique-deploye}
  On a une inclusion évidente   
      \begin{gather}\label{inclusion-adlique}G(F)\backslash G(\mb A)/K_{N}\Xi\subset \Bun_{G,N}(\Fq)/\Xi \end{gather}
   dont l'image est formée exactement par les $G$-torseurs localement triviaux pour la topologie de Zariski. 
   On rappellera dans le chapitre \ref{para-non-deploye} que si $G$ est un groupe réductif général (non nécessairement déployé),    
      $\Bun_{G,N}(\Fq)/\Xi$ est une réunion finie, indexée par $\ker^{1}(F,G)$,  de tels quotients adéliques  pour certaines formes intérieures de $G$ 
  (ce qui donne un sens à la  notion de fonction cuspidale). 
   Quand $G$ est déployé  il s'avère que  \eqref{inclusion-adlique} est une bijection    car d'après   Kottwitz~\cite{kottwitz1,kottwitz2} (et le théorème 2.6.1 de  Nguyen Quoc Thang~\cite{thang}   pour l'adaptation en caractéristique $p$), $\ker^{1}(F,G)$ est le dual de $\ker^{1}(F,Z_{\wh G}(\Qlbar))$, qui est nul par le théorème de Tchebotarev dès lors que 
   $\on{Gal}(\ov F/F)$ agit trivialement sur $Z_{\wh G}$. 
 Nous utilisons ce fait de fa\c con extrêmement superficielle : cela nous a permis d'écrire  $G(F)\backslash G(\mb A)/K_{N}\Xi$ à la place  de $\Bun_{G,N}(\Fq)/\Xi$ dans le \thmref{intro-thm-ppal}, et de rendre ainsi l'introduction plus accessible. 
       \end{rem}

   \begin{notation} Une fonction $f\in C_{c} (G(F)\backslash G(\mb A)/K_N \Xi,E)$ 
   est dite cuspidale si pour tout 
    parabolique $P\subsetneq G$, de  Levi $M$ et de radical unipotent $U$,    le terme constant $f_{P}: g\mapsto \int_{U(F)\backslash U(\mb A)}f(ug)$ est  nul comme fonction sur $U(\mb A)M(F)\backslash G(\mb A)/K_{N}\Xi$.  
 On rappelle que $C_{c}^{\rm{cusp}}(G(F)\backslash G(\mb A)/K_N \Xi,E)$ est un $E$-espace vectoriel  de dimension finie (voir les références après 
 \eqref{def-cusp}) et  stable par   $T(f)$ pour $f\in C_{c}(K_{N}\backslash G(\mb A)/K_{N},E)$. 
Grâce à l'égalité $\Bun_{G,N}(\Fq)/\Xi=G(F)\backslash G(\mb A)/K_{N}\Xi$
on note 
       $C_{c}^{\rm{cusp}}(\Bun_{G,N}(\Fq)/\Xi,E)$ comme 
      $C_{c}^{\rm{cusp}}(G(F)\backslash G(\mb A)/K_N \Xi,E)$
      (plus généralement pour un groupe non déployé, on fera de même en utilisant la somme finie indexée par des formes intérieures de $G$).        \end{notation}

\begin{prop}\label{prop-cusp-hecke-finies}
L'espace des  fonctions cuspidales  $C_{c}^{\rm{cusp}}(\Bun_{G,N}(\Fq)/\Xi,E)$ est exactement le sous-espace 
$\Big(\varinjlim _{\mu}\restr{\mc H _{ N, \emptyset, \mbf  1}^{0,\leq\mu,E}}{\Fqbar}\Big)^{\mr{Hf}}$ 
des fonctions  Hecke-finies  de 
$C_{c}(\Bun_{G,N}(\Fq)/\Xi,E)$. 
\end{prop}
\noindent{\bf Démonstration.} 
On montre d'abord que toute fonction cuspidale est Hecke-finie. 
On pose 
$$C_{c}^{\rm{cusp}}(\Bun_{G,N}(\Fq)/\Xi,\mc O_{E})=C_{c}^{\rm{cusp}}(\Bun_{G,N}(\Fq)/\Xi,E)
\cap C_{c}(\Bun_{G,N}(\Fq)/\Xi,\mc O_{E}).$$
Alors $C_{c}^{\rm{cusp}}(\Bun_{G,N}(\Fq)/\Xi,\mc O_{E})$ est un sous-$\mc O_{E}$-module de type fini. Il est stable par tous les opérateurs $T(f)$ pour $f\in C_{c}(K_{N}\backslash G(\mb A)/K_{N},\mc O_{E})$ car c'est le cas de $C_{c}^{\rm{cusp}}(\Bun_{G,N}(\Fq)/\Xi,E)
$ et de $C_{c}(\Bun_{G,N}(\Fq)/\Xi,\mc O_{E})$. 

Il reste à montrer que toute fonction Hecke-finie est cuspidale. Cela résulte du lemme suivant. 
\cqfd

  \begin{lem}
        Soit $v\in |X|\sm |N|$ et $H$ un sous-espace de dimension finie de 
 $C_{c}(\Bun_{G,N}(\Fq)/\Xi,E)$, stable par tous les opérateurs de Hecke en $v$. Alors $H$ est inclus dans  $C_{c}^{\mr{cusp}}(\Bun_{G,N}(\Fq)/\Xi,E)$.      \end{lem}     
 \noindent
 Dans ce lemme et le suivant, l'hypothèse de Hecke-finitude (en $v$) est au sens $E$-linéaire, ce qui est évidemment plus faible qu'au sens $\mc O_{E}$-linéaire, mais suffit pour entraîner la cuspidalité. 
 
        \noindent{\bf Démonstration.} 
        D'après la \remref{quotient-adelique-deploye} il suffit de montrer le résultat pour le quotient adélique $G(F)\backslash G(\mb A)/K_{N}\Xi$. 
        Le lemme résulte alors du lemme suivant, qui servira de nouveau lorsqu'on étudiera le cas où $G$ n'est pas déployé.
        En fait il est inutilement compliqué d'invoquer la \remref{quotient-adelique-deploye} car il suffirait de dire que 
        $\Bun_{G,N}(\Fq)$ est une réunion finie de quotients adéliques pour des formes intérieures de $G$ et d'appliquer le lemme suivant à chacun d'entre eux. 
         \cqfd

\begin{lem}\label{lem-TVv-cusp}
Soit $G$ un  groupe réductif sur $F$ et soit  $K_{N}=\prod_{v\in |X|} K_{N,v}$ un  sous-groupe ouvert compact de $G(\mathbb A)$. 
On note  $Z_{G}$ le centre  de $G$. Soit $\Xi$ un réseau dans $Z_{G}(\mb A)/Z_{G}(F)$. 
Soit $v$ une  place où $G$ est déployé et $K_{N,v}=G(\mc O_{v})$. Soit $H\subset C_{c}(G(F) \backslash G(\mb A)/K_{N}\Xi,E)$  un sous-espace vectoriel  de dimension finie. On suppose que  $H$ est stable par tous les opérateurs de Hecke en  $v$. Alors $H$ est inclus dans $C_{c}^{\mr{cusp}}(G(F) \backslash G(\mb A)/K_{N}\Xi,E)$. 
\end{lem}
       \noindent{\bf Démonstration.} 
      On raisonne par l'absurde et on suppose que  $H$ n'est pas inclus dans $C_{c}^{\mr{cusp}}(G(F) \backslash G(\mb A)/K_{N}\Xi,E)$. Il existe alors un parabolique $P\subsetneq G$, de radical unipotent $U$, tel que le terme constant $f_{P}: g\mapsto \int_{U(F)\backslash U(\mb A)}f(ug)$ ne s'annule pas identiquement pour  $f\in H$. Soit $M$ un sous-groupe de Levi de $P$.  On rappelle que pour tout $f\in H$, $f_{P}$ est une fonction localement constante sur $U(\mb A)M(F)\backslash G(\mb A)/K_{N}\Xi$.
      Soit $S$ le tore déployé maximal de la partie connexe du centre $Z_{M}$ de $M$. 
      D'après (i) de la proposition 20.6 de \cite{borel-alg-groups}, $M$ est égal au centralisateur de $S$ dans $G$ et d'après la preuve de cette assertion   
      il existe $\mb G_{m}\subset S$ agissant sur $\on{Lie}(U)$ par des caractères strictement positifs. On fixe un élément de $F_{v}^{*}$ de norme $>1$  et on note $a\in S(F_{v})$ son image par l'inclusion ci-dessus. Il en résulte que pour tout caractère $\chi:S\to \mb G_{m}$ apparaissant dans l'action de $S$ sur $\on{Lie}(U)$ la norme  de $\chi(a)\in F_{v}^{*}$ est $>1$. Par conséquent la conjugaison par $a^{-1}$ contracte $U(F_{v})$ 
      (au sens où pour tout compact $\mathfrak W\subset U(F_{v})$ et tout ouvert 
      $\mathfrak V\subset U(F_{v})$ contenant $1$, on a $a^{-n}\mathfrak W a^{n}\subset \mathfrak V$ pour $n$ assez grand). 
      On en déduit  que pour tout $g\in G(\mb A)$ il existe $n_{0}\in \Z
$ tel que pour tout $f\in H$ et pour tout $n\geq n_{0}$, $f_{P}(a^{n}g)=f(a^{n}g)$. 
En effet on fixe un 
domaine fondamental compact $\mathfrak U=\prod_{w}\mathfrak U _{w}$ pour l'action de $U(F)$ sur $U(\mb A)$, tel que,  pour tout $w\neq v$,  
$\mathfrak U _{w}$ est un sous-groupe ouvert compact de 
$U(F_{w})$ inclus dans 
$g K_{N,w} g^{-1} $. Alors  pour $n$ assez grand  on a $a^{-n}    \mathfrak U_{v} a^{n}\subset 
g K_{N,v} g^{-1} $, et donc $a^{-n}    \mathfrak U a^{n}\subset 
g K_{N} g^{-1} $. En notant $du$ la mesure sur  $\mathfrak U$ provenant de la mesure de Haar sur $U(F)\backslash U(\mb A)$, on a donc 
$$f_{P}(a^{n}g)=\int_{u\in \mathfrak U} f(ua^{n}g)du=\int_{u\in \mathfrak U} f(a^{n}g (g^{-1}a^{-n}    u a^{n} g ))du=f(a^{n}g)$$ où la dernière égalité résulte  du fait que  $f$ est invariante à droite 
  par $K_{N}$. 
  
   Comme tout $f\in H$ est à support compact (c'est-à-dire à support fini) sur $G(F) \backslash G(\mb A)/K_{N}\Xi$ et que $H$ est de dimension finie, 
   il existe $n_{0}\in \Z
$ tel que pour tout $f\in H$ et pour tout $n\geq n_{0}$, 
   $f(a^{n}g)=0$.  
   En effet quand $n$ tend vers l'infini $a^{n}g$ 
   sort de toute partie finie   de 
   $G(F) \backslash G(\mb A)/K_{N}\Xi$. 
   En effet, quitte à augmenter $N$ il suffit de le montrer pour $g=1$. On fixe un plongement $\iota: G^{\on{ad}}\hookrightarrow SL_{r}$. 
   Comme  $\mb G_{m}\subset S$  agit non trivialement par conjugaison il n'est pas inclus dans  $Z_{G}$, donc, à un quotient près par un $\mu_{m}$,  il     se plonge dans $SL_{r}$. Par conséquent   le point de $\Bun_{SL_{r}}(\Fq)$ associé à $\iota(a^{n})$ sort de toute partie finie quand $n$ tend vers l'infini
   (car le polygône de Harder-Narasimhan
   associé  sort de tout ensemble fini de polygônes de Harder-Narasimhan).

   On déduit des deux faits précédents que, pour  tout $g\in G(\mb A)$,  il existe $n_{0}\in \Z
$ tel que pour tout $f\in H$ et pour tout $n\geq n_{0}$, $f_{P}(a^{n}g)=0$.

       On a une action naturelle de $M(\mb A)$ sur $U(\mb A)\backslash G(\mb A)/K_{N} $ par translation à gauche. Donc il existe un point  $g_{1}\in U(\mb A) \backslash G(\mb A)/K_{N}$ de stabilisateur $K'$ (sous-groupe ouvert compact de $M(\mb A)$)  
       tel que la restriction $f_{P}'$  de $f_{P}$ à la $M(\mb A)$-orbite passant par ce point, c'est-à-dire à  l'image de l'inclusion 
     $$ M(F)\backslash M(\mb A)/K'\Xi\subset U(\mb A)M(F)\backslash G(\mb A)/K_{N}\Xi, \ \ \ m\mapsto mg_{1} ,$$
soit non identiquement nulle pour au moins un $f\in H$. 
On a \begin{gather}\nonumber U(\mb A) \backslash G(\mb A)/K_{N} =\prod_{w} U(F_{w})\backslash G(F_{w})/K_{N,w} \end{gather}  et  le facteur en $w=v$ est égal à $M(F_{v})/M(\mc O_{v})$, donc quitte à multiplier $g_{1}$ à gauche par un élément de $M(F_{v})$
on peut supposer que  $K'$ contient $M(\mc O_{v})$, ce que l'on fait désormais.  
On fixe un tel $g_{1}$ et on le relève en $g_{1}\in G(\mb A)$.

   On sait que l'algèbre de Hecke sphérique pour $M$ en $v$ est un module de type fini
   sur l'algèbre de Hecke sphérique pour $G$ en $v$
   (par exemple parce que
      l'anneau des représentations   de   $\wh M$ est un module de type fini sur  l'anneau des représentations de $\wh G$).  
    Donc   la Hecke-finitude pour $G$ en $v$ implique  la Hecke-finitude pour $M$ en $v$  et donc l'espace des       $f_{P}'$  pour $f\in H$ peut être inclus dans  un sous-espace de dimension finie de 
      $C(M(F)\backslash M(\mb A)/K'\Xi,E)$ stable par 
      $C_{c}(M(\mc O_{v})\backslash M(F_{v})/M(\mc O_{v}),E)$. Ces opérateurs de  Hecke   contiennent comme cas particuliers les translations par  
       $Z_{M}(F_{v})/Z_{M}(\mc O_{v})$, et donc en particulier par $S(F_{v})/S(\mc O_{v})$. 
         
         Soit $m_{1}\in M(\mb A)$ tel qu'il existe $f\in H$ tel que  $f'_{P}(m_{1})$ ne soit pas nul. Pour tout $f\in H$, on note $f_{P}''$ la fonction $n\mapsto f'_{P}(a^{n}m_{1})=f_{P}(a^{n}m_{1}g_{1})$ qui appartient donc à l'espace $C^{-}(\Z, E)$ des fonctions $k:\Z\to E$ telles que $k(n)=0$ pour $n$ assez grand. 
         Sur $C^{-}(\Z, E)$  on a l'opérateur de translation $T$ tel que $T(f): n\mapsto f(n+1)$. 
         L'image de $H$ par $f\mapsto f_{P}''$ est donc un sous-$E$-espace non nul de 
         $C^{-}(\Z, E)$ et on sait qu'il peut être inclus dans un sous-$E$-espace de dimension finie stable par $T$. On aboutit à une contradiction car pour tout élément non nul $k\in 
C^{-}(\Z, E)$ on montre (en considérant le maximum de son support) que  les $T^{n}(k)$ pour $n\in \Z$ sont linéairement indépendants et engendrent donc un espace vectoriel de dimension infinie. 
             \cqfd
       
      \begin{rem}
    En fait  $f_{P}$ et $f_{P}'$ sont supportés sur les composantes indexées par un translaté d'un  cône dans le  dual du réseau des  caractères de $M/Z_{G}$, qui s'identifie (après tensorisation par $\mb Q$) au réseau des cocaractères de $S/Z_{G}$. Ce cône est engendré par les  opposées  des projections  des coracines associées à $U$ sur le réseau des cocaractères de $S/Z_{G}$ (tensorisé par $\mb Q$). 
    Lorsque $G$ est déployé cela est montré dans le   lemme 5.3.1 de \cite{wang-bilinear}.      
    \end{rem}

La  proposition suivante jouera un   rôle crucial  dans  la  construction des opérateurs d'excursion  qui sera l'objet du  prochain chapitre.

On a vu que le $E$-espace vectoriel    $\Big(\varinjlim _{\mu}\restr{\mc H _{ N, I, W}^{0,\leq\mu,E}}{\ov{\eta^{I}}}\Big)^{\mr{Hf}}$ est stable sous l'action  de $\pi_{1}(\eta^{I},\ov{\eta^{I}})$ et des  morphismes de  Frobenius partiels. Autrement dit  il est stable sous l'action  de $\on{FWeil}(\eta^{I},\ov{\eta^{I}})$. 

\begin{prop}\label{prop-action-Hf}
L'action de $\on{FWeil}(\eta^{I},\ov{\eta^{I}})$ sur $\Big(\varinjlim _{\mu}\restr{\mc H _{ N, I, W}^{0,\leq\mu,E}}{\ov{\eta^{I}}}\Big)^{\mr{Hf}}$ se factorise (de fa\c con unique) à travers  $(\pi_{1}(\eta,\ov{\eta}))^{I}$.    
 
De plus  tout sous-$\mc O_{E}$-module $\mf M$ de type fini  de $\Big(\varinjlim _{\mu}\restr{\mc H _{ N, I, W}^{0,\leq\mu,E}}{\ov{\eta^{I}}}\Big)^{\mr{Hf}}$ est inclus dans un  sous-$\mc O_{E}$-module $\wt{\mf M}$ de type fini  stable par   $(\pi_{1}(\eta,\ov{\eta}))^{I}$ 
 (ou de fa\c con équivalente  stable par  $\pi_{1}(\eta^{I},\ov{\eta^{I}})$ et par les   morphismes de  Frobenius partiels) et 
 l'action de $(\pi_{1}(\eta,\ov{\eta}))^{I}$  sur $\wt{\mf M}$ se factorise à travers $(\pi_{1}(U, \ov\eta))^{I}$ pour un ouvert dense   $U\subset X$ (dépendant {\it a priori} de  $\wt{\mf M}$, donc de $\mf M$). 
  \end{prop}
 
 \begin{rem} On peut exprimer cette  action de $(\pi_{1}(\eta,\ov{\eta}))^{I}$  sur $\wt{\mf M}$ de fa\c con plus géométrique: il existe 
 \begin{itemize}
  \item un object $(\mf E,(F_{\{i\}})_{i\in I})$ de la  catégorie 
  $\mc C(U,I,\mc O_{E})$ définie dans le  \lemref{lem-Frob-partiels-drinfeld}, 
  \item un  copoids dominant $\mu_{1}$,  un ouvert dense 
  $\Omega_{1}\subset U^{I}$  sur lequel   $\mc H _{ N, I, W}^{0,\leq\mu_{1},E}$ est lisse, et une  injection  $\iota_1: \restr{\mf E}{\Omega_{1}}\hookrightarrow \restr{\mc H _{ N, I, W}^{0,\leq\mu_{1},E}}{\Omega_{1}}$ de $\mc O_{E}$-faisceaux lisses sur  $\Omega_{1}$,
     \end{itemize}
 tels que   
 \begin{itemize}
 \item $\restr{\iota_1}{\eta^{I}}$ 
 est compatible avec l'action des  morphismes de  Frobenius partiels  sur  $\varinjlim_{\mu} \restr{\mc H _{ N, I, W}^{0,\leq\mu,E}}{\eta^{I}}$
 \item $\restr{\iota_1}{\ov{\eta^{I}}}: \restr{\mf E}{\ov{\eta^{I}}} \to \varinjlim_{\mu} \restr{\mc H _{ N, I, W}^{0,\leq\mu,E}}{\ov{\eta^{I}}}$ est une  injection de $\mc O_{E}$-modules, dont l'image est $\wt{\mf M}$, 
 \item l'action de  $(\pi_{1}(U, \ov\eta))^{I}$ sur $\restr{\mf E}{\ov{\eta^{I}}} $ (et donc  sur $\wt{\mf M}$ par l'intermédiaire de  $\restr{\iota_1}{\ov{\eta^{I}}}$) 
 se déduit de  l'équivalence de categories du  \corref{equiv-cat-U-I-Omega}. 
    \end{itemize}    
    On ne sait pas si l'on peut espérer que   ci-dessus  on puisse avoir 
     $\Omega_{1}=U^{I}$.  Heureusement on n'en a pas besoin. 
         \end{rem}

\begin{rem} Dans \cite{these-cong}, Cong Xue a montré que  $\Big(\varinjlim _{\mu}\restr{\mc H _{ N, I, W}^{0,\leq\mu,E}}{\ov{\eta^{I}}}\Big)^{\mr{Hf}}$ est  de dimension finie, mais la preuve est difficile (et écrite seulement pour $G$ déployé) et  on ne l'utilise pas dans cet article.  \end{rem}

  \noindent{\bf Démonstration de  la \propref{prop-action-Hf}. }
 Il suffit de montrer l'énoncé pour   $W$   irréductible. 
 On écrit 
      $W=\boxtimes_{i\in I}W_{i}$. 
      
    On fixe  $\kappa$ assez grand   (en   fonction de $W$) pour que l'action des morphismes de  Frobenius partiels  soit donnée par des morphismes  
   \begin{gather}
   \label{isom-frob-partiels}
   F_{\{i\}}: \Frob_{\{i\}}^{*}(\mc H _{ N, I, W}^{0,\leq\mu,E})\to 
     \mc H _{ N, I, W}^{0,\leq\mu+\kappa,E} \end{gather}
de $E$-faisceaux constructibles sur  $(X\sm N)^{I}$, pour tout $\mu$.

   Soit  $\mf M$ un  sous-$\mc O_{E}$-module de type fini de 
  $\varinjlim _{\mu}\restr{\mc H _{ N, I, W}^{0,\leq\mu,E}}{\ov{\eta^{I}}}$,  stable par  
$\pi_{1}(\eta^{I},\ov{\eta^{I}})$ et   $C_{c}(K_{N}\backslash G(\mb A)/K_{N},\mc O_{E})$. 
Comme   $\mf M$ est de type fini,  il existe $\breve \mu$ tel que 
   $\mf M$ est inclus dans l'image de $\restr{\mc H _{ N, I, W}^{0,\leq\breve\mu,E}}{\ov{\eta^{I}}}$ dans 
      $\varinjlim_{\mu }\restr{\mc H _{ N, I, W}^{0,\leq\mu,E}}{\ov{\eta^{I}}}$. On choisit   $\mu_{0}\geq \breve \mu$ assez grand  pour que cette  image soit égale à l'image de  $\restr{\mc H _{ N, I, W}^{0,\leq\breve\mu,E}}{\ov{\eta^{I}}}$ dans  $\restr{\mc H _{ N, I, W}^{0,\leq\mu_{0},E}}{\ov{\eta^{I}}}$. Alors 
         $\mf M$ est canoniquement un  sous-$\mc O_{E}$-module de       $\restr{\mc H _{ N, I, W}^{0,\leq\mu_{0},E}}{\ov{\eta^{I}}}$. 
               On note  $\mf G$
 le sous-$\mc O_{E}$-faisceau (sur $\eta^{I}$) de 
  $ \restr{\mc H _{ N, I, W}^{0,\leq\mu_{0},E}}{\eta^{I}}$ tel que $\restr{\mf G}{\ov {\eta^{I}}}=\mf M$.        Alors, pour tout 
   $(n_{i})_{i\in I}\in \N^{I}$, 
   $$\prod_{i\in I} F_{\{i\}}^{n_{i}}\big((\prod_{i\in I} \Frob_{\{i\}}^{n_{i}})^{*}(\mf G)\big)$$ est un sous-$\mc O_{E}$-faisceau  de 
  $ \restr{\mc H _{ N, I, W}^{0,\leq\mu_{0}+(\sum_{i\in I}n_{i})\kappa,E}}{\eta^{I}}$.

    Soit  $\wt {\mf M}$ la somme sur $(n_{i})_{i\in I}\in \N^{I}$  des sous-$\mc O_{E}$-représentations de $\pi_{1}(\eta^{I},\ov{\eta^{I}})$ dans le $E$-espace vectoriel $\varinjlim_{\mu} \restr{\mc H _{ N, I, W}^{0,\leq\mu,E}}{\ov{\eta^{I}}}$ données par les 
   $$\restr{\Big(\prod_{i\in I} F_{\{i\}}^{n_{i}}\big((\prod_{i\in I} \Frob_{\{i\}}^{n_{i}})^{*}(\mf G)\big)\Big)}{\ov{\eta^{I}}}.$$  
   
   \begin{lem}\label{lem-frob-tf}
Le $\mc O_{E}$-module 
   $\wt {\mf M}$ est  de type fini   et est réalisé comme un sous-$\mc O_{E}$-module de $
      \restr{\mc H _{ N, I, W}^{0,\leq\wt \mu_{0},E}}{\ov{\eta^{I}}}$ pour $\wt \mu_{0}$ assez grand. De plus 
              le sous-$\mc O_{E}$-faisceau
            $\wt{\mf G}$ de $ \restr{\mc H _{ N, I, W}^{0,\leq\wt \mu_{0},E}}{\eta^{I}}$             tel que $\restr{\wt{\mf G}}{\ov {\eta^{I}}}=\wt{\mf M}$
  est  stable par  l'action des morphismes de Frobenius partiels  dans 
  $\varinjlim_{\mu} 
      \restr{\mc H _{ N, I, W}^{0,\leq\mu,E}}{\eta^{I}}$, et 
      les 
      $F_{\{i\}}: \Frob_{\{i\}}^{*}(\wt {\mf G})\to \wt {\mf G}$ 
      sont des isomorphismes.          \end{lem}
    \noindent{\bf Démonstration du \lemref{lem-frob-tf}. } Il suffit  de montrer que $\wt {\mf M}$ est de type fini comme  $\mc O_{E}$-module. En effet 
    la stabilité de   $\wt{\mf G}$  par les actions des morphismes de Frobenius partiels vient de la construction même de $\wt{\mf M}$ et les 
    $F_{\{i\}}: \Frob_{\{i\}}^{*}(\wt {\mf G})\to \wt {\mf G}$ 
      sont alors  des isomorphismes car leur produit dans n'importe quel ordre 
     est l'action du  morphisme de Frobenius global. 
     
     Il reste donc à montrer que $\wt {\mf M}$ est de type fini. 
  Soit $\mu_{0}$ comme ci-dessus, c'est-à-dire tel que 
  $\mf G$ se réalise comme  un sous-$\mc O_{E}$-faisceau  de 
  $ \restr{\mc H _{ N, I, W}^{0,\leq\mu_{0},E}}{\eta^{I}}$. Soit $\Omega_{0}$ un ouvert non vide  de $(X\sm N)^{I}$ sur lequel 
  $\mc H _{ N, I, W}^{0,\leq\mu_{0},E}$ est lisse.   On prolonge $\mf G$ de manière unique en un  sous-$\mc O_{E}$-faisceau de $\restr{\mc H _{ N, I, W}^{0,\leq\mu_{0},E}}{\Omega_{0}}$, que l'on note encore $\mf G$, et qui est donc un  $\mc O_{E}$-faisceau lisse constructible sur $\Omega_{0}$ tel que $ \restr{\mf G}{\ov{\eta^{I}}}=\mf M$.

   Pour toute famille  $(v_{i})_{i\in I}$ de   points fermés de $X\sm N$, on note   ${\times_{i\in I} v_{i}}$ leur produit, qui est une réunion finie  de points fermés de $(X\sm N)^{I}$.   Par récurrence  sur  $\sharp I$ on montre que tout ouvert dense de $X^{I}$ contient un tel  produit ${\times_{i\in I} v_{i}}$. 
On fixe 
 $(v_{i})_{i\in I}$ tel que ${\times_{i\in I} v_{i}}$ soit inclus dans  $\Omega_0$. 

 Soit $i\in I$.   D'après la 
               relation d'Eichler-Shimura en  $v_{i}$ (\propref{prop-eichler-shimura}), on a 
               $$\sum_{\alpha=0}^{\dim W_{i}} (-1)^{\alpha} (F_{\{i\}}^{\deg(v_{i})})^{\alpha}\circ S_{\Lambda^{\dim W_{i}-\alpha}W_{i},v_i}=0\text{ \ \  dans \ \ }
  \on{End}\big(\varinjlim_{\mu} \restr{\mc H _{ N, I, W}^{0,\leq\mu,E}}{\times_{i\in I} v_{i}}\big).$$             
               Par conséquent 
 \begin{gather}\label{incl-v-I}(F_{\{i\}}^{\deg(v_{i})})^{\dim W_{i}}(\restr{\mf G}{\times_{i\in I} v_{i}})
   \subset 
   \sum_{\alpha=0}^{\dim W_{i}-1}  (F_{\{i\}}^{\deg(v_{i})})^{\alpha}(S_{\Lambda^{\dim W_{i}-\alpha}W_{i},v_i}(\restr{\mf G}{\times_{i\in I} v_{i}}))
  \end{gather} dans $
  \varinjlim_{\mu} \restr{\mc H _{ N, I, W}^{0,\leq\mu,E}}{\times_{i\in I} v_{i}}$. 
On va voir que grâce à la lissité de $\mf G$ sur $\Omega_{0}$ cette inclusion se propage  à $\eta^{I}$.  

    On choisit un point géométrique $\ov{\times_{i\in I} v_{i}}$ au-dessus d'un des points fermés de ${\times_{i\in I} v_{i}}$ et une flèche de spécialisation 
 $\on{\mf{sp}}_{v,I}$ de $\ov{\eta^{I}}$ 
 vers $\ov{\times_{i\in I} v_{i}}$.   Pour tout $n\in \N$, la lissité de $\mf G$ sur $\Omega_{0}$ implique la lissité de $(\Frob_{\{i\}}^{n})^{*}\mf G$ sur 
 $(\Frob_{\{i\}}^{n})^{-1}(\Omega_{0})$ qui contient aussi le produit ${\times_{i\in I} v_{i}}$, et donc 
   \begin{gather}
   \label{sp-v-I-surjectif}\on{\mf{sp}}_{v,I}^{*}:\restr{(\Frob_{\{i\}}^{n})^{*}\mf G}{\ov{\times_{i\in I} v_{i}}}\to \restr{(\Frob_{\{i\}}^{n})^{*}\mf G}{\ov{\eta^{I}}}\text{ \ \ \   est bijectif.}\end{gather} 
   
 En prenant la fibre en $\ov{\times_{i\in I} v_{i}}$,  il  résulte de \eqref{incl-v-I} que 
  \begin{gather}\label{incl-v-I-geom}
  (F_{\{i\}}^{\deg(v_{i})})^{\dim W_{i}}(\restr{(\Frob_{\{i\}}^{\deg(v_{i})\dim W_{i}})^{*}\mf G}{\ov{\times_{i\in I} v_{i}}}) \\ \nonumber
   \subset 
   \sum_{\alpha=0}^{\dim W_{i}-1}  (F_{\{i\}}^{\deg(v_{i})})^{\alpha}\Big(\restr{(\Frob_{\{i\}}^{\deg(v_{i})\alpha})^{*}\big(S_{\Lambda^{\dim W_{i}-\alpha}W_{i},v_i}(\mf G)\big)}{\ov{\times_{i\in I} v_{i}}}\Big)
  \end{gather} dans $
  \varinjlim_{\mu} \restr{\mc H _{ N, I, W}^{0,\leq\mu,E}}{\ov{\times_{i\in I} v_{i}}}$. 
    Comme  $F_{\{i\}}$ et $S_{\Lambda^{\dim W_{i}-\alpha}}$ sont des morphismes de faisceaux, ils commutent avec $\on{\mf{sp}}_{v,I}^{*}$. En appliquant \eqref{sp-v-I-surjectif} à $n=\deg(v_{i})\dim(W_{i})$  on obtient 
     \begin{gather*}(F_{\{i\}}^{\deg(v_{i})})^{\dim W_{i}}\big(\restr{(\Frob_{\{i\}}^{\deg(v_{i})\dim W_{i}})^{*}\mf G}{\ov{\eta^{I}}}\big)
 \\  \subset 
   \sum_{\alpha=0}^{\dim W_{i}-1}  (F_{\{i\}}^{\deg(v_{i})})^{\alpha}\big(\restr{(\Frob_{\{i\}}^{\deg(v_{i})\alpha})^{*}(S_{\Lambda^{\dim W_{i}-\alpha}W_{i},v_i}(\mf G))}{\ov{\eta^{I}}}\big) 
   \end{gather*} dans 
   $
  \varinjlim_{\mu} \restr{\mc H _{ N, I, W}^{0,\leq\mu,E}}{\ov{\eta^{I}}}$. 
       Par l'hypothèse et la proposition \ref{prop-coal-frob-cas-part}, 
       $\mf M$  et donc $\restr{\mf G}{\eta^{I}}$ sont  stables par les   morphismes  $S_{\Lambda^{\dim W_{i}-\alpha}W_{i},v_i}$,  
puisque 
 $$h_{\Lambda^{\dim W_{i}-\alpha}W_{i},v_i}\in C_{c}(G(\mc O_{v_{i}})\backslash G(F_{v_{i}})/G(\mc O_{v_{i}}),\mc O_{E})\subset C_{c}(K_{N}\backslash G(\mb A)/K_{N},\mc O_{E}). $$
                     On en déduit  immédiatement 
        $$(F_{\{i\}}^{\deg(v_{i})})^{\dim W_{i}}\big(\restr{(\Frob_{\{i\}}^{\deg(v_{i})\dim W_{i}})^{*}\mf G}{\ov{\eta^{I}}}\big)
   \subset 
   \sum_{\alpha=0}^{\dim W_{i}-1}  (F_{\{i\}}^{\deg(v_{i})})^{\alpha}\big(\restr{(\Frob_{\{i\}}^{\deg(v_{i})\alpha})^{*}\mf G}{\ov{\eta^{I}}}\big) 
   $$ dans 
   $
  \varinjlim_{\mu} \restr{\mc H _{ N, I, W}^{0,\leq\mu,E}}{\ov{\eta^{I}}}$. 
Il en résulte que   $$\wt{\mf M}=\sum_{(n_{i})_{i\in I}\in \prod _{i\in I}\{0,...,\deg(v_{i})\dim(W_{i})-1\}}\prod_{i\in I} F_{\{i\}}^{n_{i}}\Big(\restr{\prod_{i\in I} (\Frob_{\{i\}}^{n_{i}})^{*}(\mf G)}{\ov{\eta^{I}}}\Big)$$  et donc $\wt{\mf M}$ est 
  un sous-$\mc O_{E}$-module de type fini  de 
  $\varinjlim _{\mu}\restr{\mc H _{ N, I, W}^{0,\leq\mu,E}}{\ov{\eta^{I}}}$, inclus dans l'image de 
  $ \restr{\mc H _{ N, I, W}^{0,\leq\mu_{0}+(\sum_{i\in I}(\deg(v_i)\dim W_{i}-1))\kappa,E}}{\ov{\eta^{I}}}$ (et stable par $C_{c}(K_{N}\backslash G(\mb A)/K_{N},\mc O_{E})$).  
                  \cqfd

        \noindent{\bf Fin de la démonstration de la \propref{prop-action-Hf}. }    
        D'après le \lemref{lem-frob-tf} il existe $\wt \mu_{0}$, un ouvert 
    $\wt \Omega_{0}\subset X^{I}$ tel que 
        $ \restr{\mc H _{ N, I, W}^{0,\leq\wt \mu_{0},E}}{\wt \Omega_{0}}$ soit
         lisse  sur $\wt \Omega_{0}$, et 
        $\wt{\mf G}$ s'étend en un sous-$\mc O_{E}$-faisceau lisse sur $\wt \Omega_{0}$
        tel que $\restr{\wt{\mf G}}{\ov{\eta^{I}}}=\wt {\mf M} $. 
        De plus  $\wt{\mf G}$ est muni d'une action des morphismes de Frobenius partiels, plus précisément pour tout $i$ on a un isomorphisme 
     $$F_{\{i\}}:     \restr{\Frob_{\{i\}}^{*}(\wt{\mf G})}{\wt \Omega_{0}\cap \Frob_{\{i\}}^{-1}(\wt \Omega_{0})}
  \to   \restr{\wt{\mf G}}{\wt \Omega_{0}\cap \Frob_{\{i\}}^{-1}(\wt \Omega_{0})}$$
  compatible avec l'action des morphismes de Frobenius partiels sur  $ \varinjlim_{\mu} \restr{\mc H _{ N, I, W}^{0,\leq\mu,E}}{\eta^{I}}$. 
 Le \lemref{lem-U^I}
       implique  qu'il existe un ouvert dense $U\subset X\sm N$ tel que 
               $(\wt{\mf G},(F_{\{i\}})_{i\in I})$ s'étende  à $U^{I}$ et fournisse un objet de   la catégorie 
  $\mc C(U,I,\mc O_{E})$.  Le \corref{equiv-cat-U-I-Omega} fournit alors 
   l'action de  $\pi_{1}(U,\ov{\eta})^{I}$ sur 
$\wt{\mf M}=\restr{\wt{\mf G}}{\ov{\eta^{I}}}$.        \cqfd

\begin{prop}\label{surjectivite-Hecke-fini}
L'image de l'homomorphisme de spécialisation 
 \begin{gather}\label{sp*-sans-Hf}\on{\mf{sp}}^{*}: \varinjlim _{\mu}\restr{\mc H _{ N, I, W}^{0,\leq\mu,E}}{\Delta(\ov{\eta})} \to \varinjlim _{\mu}\restr{\mc H _{ N, I, W}^{0,\leq\mu,E}}{\ov{\eta^{I}}}\end{gather} contient 
 $\Big( \varinjlim _{\mu}\restr{\mc H _{ N, I, W}^{0,\leq\mu,E}}{\ov{\eta^{I}}}\Big)^{\mr{Hf}}$. 
  \end{prop}
 \noindent{\bf Démonstration de  la \propref{surjectivite-Hecke-fini}.}
 D'après la  \propref{prop-action-Hf}, 
  $\Big(\varinjlim _{\mu}\restr{\mc H _{ N, I, W}^{0,\leq\mu,E}}{\ov{\eta^{I}}}\Big)^{\mr{Hf}}$   est la réunion  de sous-$\mc O_{E}$-modules $\mf M=\restr{\mf G}{\ov{\eta^{I}}}$ où 
   $\mf G$ est un sous-$\mc O_{E}$-faisceau constructible   de 
  $\varinjlim _{\mu}\restr{\mc H _{N, I, W}^{0,\leq\mu,E}}{\eta^{I}}$ 
     stable sous l'action des morphismes de  Frobenius partiels
     (pour éviter toute confusion on signale que les $\mf M$ comme ci-dessus étaient notés $\wt{\mf M}$ dans l'énoncé de la \propref{prop-action-Hf}). 
     Il suffit donc de montrer qu'un tel $\mf M$ est inclus dans l'image de  \eqref{sp*-sans-Hf}. 
           Soit $\mu_0$ assez grand  pour que   $\mf G$  soit un  sous-$\mc O_{E}$-faisceau    de    
  $\restr{\mc H _{ N, I, W}^{0,\leq\mu_0,E}}{\ov{\eta^{I}}}$.       Soit  $\Omega_0$ un ouvert dense  de  $X^{I}$ tel que  $\restr{\mc H _{ N, I, W}^{0,\leq\mu_0,E}}{\Omega_0}$ soit lisse. Alors   $\mf G$ se prolonge en un  sous-$\mc O_{E}$-faisceau lisse  de 
  $\restr{\mc H _{ N, I, W}^{0,\leq\mu_0,E}}{\Omega_0}$. 
      D'après le  lemme 9.2.1 de \cite{eike-lau} (dont l'argument a été rappelé dans la preuve du \lemref{lem-U^I}), l'ensemble des  
   $\big(\prod_{i\in I}\Frob_{\{i\}}^{n_{i}}\big)(\Delta(\eta))$ pour $(n_{i})_{i\in I}\in \N^{I}$ est  Zariski dense dans $X^{I}$. Il existe  donc 
  $(n_{i})_{i\in I}\in \N^{I}$ tel que 
  $\big(\prod_{i\in I}\Frob_{\{i\}}^{n_{i}}\big)(\Delta(\eta))$ appartienne à   $\Omega_0$. 
  
 Alors 
  $\restr{\mf G}{\big(\prod_{i\in I}\Frob_{\{i\}}^{n_{i}}\big)(\ov{\eta^{I}})}$ est inclus dans l'image de    \begin{gather}\label{texte-sp-Frob}\on{\wt{\mf{sp}}}^{*}: \varinjlim _{\mu}\restr{\mc H _{ N, I, W}^{0,\leq\mu,E}}{\big(\prod_{i\in I}\Frob_{\{i\}}^{n_{i}}\big)(\Delta(\ov{\eta}))}\to 
  \varinjlim _{\mu}\restr{\mc H _{ N, I, W}^{0,\leq\mu,E}}{\big(\prod_{i\in I}\Frob_{\{i\}}^{n_{i}}\big)(\ov{\eta^{I}})}\end{gather}
 pour toute flèche de spécialisation $$\wt{\mf{sp}}: \big(\prod_{i\in I}\Frob_{\{i\}}^{n_{i}}\big)(\ov{\eta^{I}})\to \big(\prod_{i\in I}\Frob_{\{i\}}^{n_{i}}\big)(\Delta(\ov{\eta})).$$  On prend pour 
 $\wt{\mf{sp}}$ l'image de  $\mf{sp}$ par $\prod_{i\in I}\Frob_{\{i\}}^{n_{i}}$. 
 D'où le diagramme commutatif 
  \begin{gather*}  \xymatrix{
  \varinjlim _{\mu}\restr{\mc H _{ N, I, W}^{0,\leq\mu,E}}{ (\prod \Frob  )(\Delta(\ov{\eta}))} 
   \ar[d]^-{\prod_{i\in I}F_{\{i\}}^{n_{i}} } \ar[r]^-{  \eqref{texte-sp-Frob}}
 &   \varinjlim _{\mu}\restr{\mc H _{ N, I, W}^{0,\leq\mu,E}}{(\prod \Frob )(\ov{\eta^{I}})} \ar[d]^-{\prod_{i\in I}F_{\{i\}}^{n_{i}} }  & \ar@{_{(}->}[l]^{}
 \restr{\mf G}{(\prod \Frob )(\ov{\eta^{I}})} \ar[d]^-{\wr }
 \\  \varinjlim _{\mu}\restr{\mc H _{ N, I, W}^{0,\leq\mu,E}}{\Delta(\ov{\eta})} 
  \ar[r]^-{ \eqref{sp*-sans-Hf}}  
& \varinjlim _{\mu}\restr{\mc H _{ N, I, W}^{0,\leq\mu,E}}{\ov{\eta^{I}}} 
& \ar@{_{(}->}[l]^{} \restr{\mf G}{\ov{\eta^{I}}}}
 \end{gather*}
  où $(\prod \Frob)$ est un raccourci pour $\prod_{i\in I}\Frob_{\{i\}}^{n_{i}}$. 
 La bijectivité de la flèche verticale la plus à droite vient du fait que 
 $\mf G$ est stable par l'action  des  morphismes de  Frobenius partiels. 
 On conclut que 
  $\mf M=\restr{\mf G}{\ov{\eta^{I}}}$ est inclus dans l'image de \eqref{sp*-sans-Hf}. \cqfd
  
\begin{prop}\label{injectivite-sp}
L'homomorphisme de spécialisation 
 \begin{gather}\label{sp*-sans-Hf-intro2}\on{\mf{sp}}^{*}: 
 \varinjlim _{\mu}\restr{\mc H _{N, I, W}^{0,\leq\mu,E}}{\Delta(\ov{\eta})} \to
  \varinjlim _{\mu}\restr{\mc H _{N, I, W}^{0,\leq\mu,E}}{\ov{\eta^{I}}}\end{gather} est injectif.   \end{prop}
\noindent{\bf Démonstration. }
Soit $a$ dans le noyau de \eqref{sp*-sans-Hf-intro2}. 
On choisit $\mu_{0}$ et $\wt a\in \restr{\mc H _{N, I, W}^{0,\leq\mu_{0},E}}{\Delta(\ov{\eta})}$ tels que 
 $a$ soit  l'image de $\wt a$ dans 
 $\varinjlim _{\mu}\restr{\mc H _{N, I, W}^{0,\leq\mu,E}}{\Delta(\ov{\eta})} $. 
Soit $\Omega_{0}$ un ouvert dense de $X\sm N$ sur lequel 
$\Delta^{*}\big( \mc H _{N, I, W}^{0,\leq\mu_{0},E}\big)$ est lisse. Soit $v\in |\Omega_{0}|$. On pose $d=\deg(v)$ pour raccourcir les formules. Soit $\ov v$ un point géométrique au-dessus de $v$. 
Soit $\on{\mf{sp}}_{v}:\ov \eta\to \ov v$ une flèche de spécialisation. On note encore  $\on{\mf{sp}}_{v}:\Delta(\ov \eta)\to \Delta(\ov v)$ la flèche de spécialisation qui s'en déduit. 
Grâce à la lissité de $\Delta^{*}\big( \mc H _{N, I, W}^{0,\leq\mu_{0},E}\big)$  sur 
 $\Omega_{0}$ on possède un unique élément 
 $\wt b\in \restr{\mc H _{N, I, W}^{0,\leq\mu_{0},E}}{\Delta(\ov{v})}$ tel que $\wt a =\on{\mf{sp}}_{v}^{*}(\wt b)$. 
 On note $b$ l'image de  $\wt b$ dans $\varinjlim _{\mu}\restr{\mc H _{N, I, W}^{0,\leq\mu,E}}{\Delta(\ov{v})}$, de sorte que $a$ est l'image de $b$ par 
 $$\on{\mf{sp}}_{v}^{*}: \varinjlim _{\mu}\restr{\mc H _{N, I, W}^{0,\leq\mu,E}}{\Delta(\ov{v})}\to 
  \varinjlim _{\mu}\restr{\mc H _{N, I, W}^{0,\leq\mu,E}}{\Delta(\ov{\eta})} .$$
  
L'action des morphismes de Frobenius partiels fournit  pour tout $\mu$  et pour tout 
$(n_{i})_{i\in I}\in \N^{I}$  
un morphisme de faisceaux sur $(X\sm N)^{I}$ 
\begin{gather}\label{action-Frob-partiels}\prod_{i\in I}F_{\{i\}}^{dn_{i}}:(\prod_{i\in I}\Frob_{\{i\}}^{dn_{i}})^{*} ( \mc H _{N, I, W}^{0,\leq\mu,E})\to 
 \mc H _{N, I, W}^{0,\leq\mu+\kappa(\sum n_{i}),E}\end{gather} 
avec  $\kappa $ assez grand en fonction de $W$ et de $d$.   
Comme $ \prod_{i\in I}\Frob_{\{i\}}^{dn_{i}} $ agit trivialement sur $\Delta(v)$, 
 $\prod_{i\in I}F_{\{i\}}^{dn_{i}}$ agit sur 
$\varinjlim _{\mu}\restr{\mc H _{N, I, W}^{0,\leq\mu,E}}{\Delta(\ov{v})}$. 

Pour tout $(n_{i})_{i\in I}\in \N^{I}$ on note 
$$b_{(n_{i})_{i\in I}} =\prod_{i\in I}F_{\{i\}}^{dn_{i}}(b)\in 
\varinjlim _{\mu}\restr{\mc H _{N, I, W}^{0,\leq\mu,E}}{\Delta(\ov{v})}. $$
En particulier 
$b_{(0)_{i\in I}}=b$. 

On pose 
\begin{gather}\label{def-a-n-i}a_{(n_{i})_{i\in I}} =\on{\mf{sp}}_{v}^{*}(b_{(n_{i})_{i\in I}} )\in 
\varinjlim _{\mu}\restr{\mc H _{N, I, W}^{0,\leq\mu,E}}{\Delta(\ov{\eta})}, \end{gather}
de sorte que $a_{(0)_{i\in I}}=a$. 

La suite $a_{(n_{i})_{i\in I}}$ vérifie  les deux propriétés énoncées dans le lemme suivant. La première affirme que cette suite est ``multirécurrente'', c'est-à-dire récurrente en chaque variable $n_{i}$, et la seconde implique qu'elle est ``presque partout'' nulle. On déduira aisément de la conjonction des deux propriétés que cette suite est partout nulle, et donc en particulier que 
$a=a_{(0)_{i\in I}}$ est nul. 

Pour énoncer la seconde propriété on remarque que   $\on{\mf{sp}}^{*}(\on{\mf{sp}}_{v}^{*}(b))=\on{\mf{sp}}^{*}(a)=0$  
 dans $\varinjlim _{\mu}\restr{\mc H _{N, I, W}^{0,\leq\mu,E}}{\ov{\eta^{I}}}
$. 
Donc il existe $\mu_{1}\geq \mu_{0}$ tel que 
$\on{\mf{sp}}^{*}(\on{\mf{sp}}_{v}^{*}(\wt b))\in \restr{\mc H _{N, I, W}^{0,\leq\mu_{0},E}}{\ov{\eta^{I}}}$
ait une image   nulle dans $\restr{\mc H _{N, I, W}^{0,\leq\mu_{1},E}}{\ov{\eta^{I}}}
$. Autrement dit en notant 
$\wh b$ l'image de $\wt b$ dans 
$\restr{\mc H _{N, I, W}^{0,\leq\mu_{1},E}}{\Delta(\ov{v})}$, 
on a 
$\on{\mf{sp}}^{*}(\on{\mf{sp}}_{v}^{*}(\wh b))=0$ dans 
$\restr{\mc H _{N, I, W}^{0,\leq\mu_{1},E}}{\ov{\eta^{I}}}$. 
Soit $\Omega_{1}\subset (X\sm N)^{I}$ un ouvert dense sur lequel 
$\mc H _{N, I, W}^{0,\leq\mu_{1},E}$ est lisse. 

\begin{lem}\label{lem-a-b}
a) Pour tout $j\in I$ et pour tout $(n_{i})_{i\in I}\in \N^{I}$, 
\begin{gather} \label{relation-ani}
\sum_{\alpha=0}^{\dim W_{j}}(-1)^{\alpha} S_{\Lambda^{\dim W_{j}-\alpha}W_{j},v} (a_{(n_{i}+\alpha \delta_{i,j})_{i\in I}})=0
\end{gather}
dans $\varinjlim _{\mu}\restr{\mc H _{N, I, W}^{0,\leq\mu,E}}{\Delta(\ov{\eta})}$. 

b)  
Pour tout $(n_{i})_{i\in I}\in \N^{I}$ tel que $\prod_{i\in I}\Frob_{\{i\}}^{dn_{i}}
(\Delta(\ov\eta))\in \Omega_{1}$, on a $a_{(n_{i})_{i\in I}}=0$ dans 
$\varinjlim _{\mu}\restr{\mc H _{N, I, W}^{0,\leq\mu,E}}{\Delta(\ov{\eta})} $. 
\end{lem}

\noindent{\bf Démonstration de a). }
Les $b_{(n_{i})_{i\in I}} $ satisfont une  relation identique à 
\eqref{relation-ani} 
(dans $\varinjlim _{\mu}\restr{\mc H _{N, I, W}^{0,\leq\mu,E}}{\Delta(\ov{v})}$), à savoir la  relation d'Eichler-Shimura en la patte $j$ (\propref{prop-eichler-shimura}). 
Alors \eqref{relation-ani} s'obtient en appliquant 
$\on{\mf{sp}}_{v}^{*}$ à 
 cette relation (ce qui est légitime puisque les $S_{\Lambda^{\dim W_{j}-\alpha}W_{j},v}$ sont des morphismes de faisceaux).  \cqfd

\noindent{\bf Démonstration de b). } Soit $(n_{i})_{i\in I}$ satisfaisant l'hypothèse
 de b). 
Comme 
\eqref{action-Frob-partiels} 
est un morphisme de faisceaux sur $(X\sm N)^{I}$, 
on peut intervertir les homomorphismes de spécialisation et les morphismes de Frobenius partiels. 
Autrement dit 
on a un diagramme commutatif 
  \begin{gather*} \xymatrixcolsep{5pc} \xymatrix{
     \restr{ \mc H _{N, I, W}^{0,\leq\mu_{1},E}}{\Delta(\ov v)}
=
\restr{(\prod_{i\in I}\Frob_{\{i\}}^{dn_{i}})^{*} ( \mc H _{N, I, W}^{0,\leq\mu_{1},E})}{\Delta(\ov v)}\ar[d]^-{ \on{\mf{sp}}_{v,(n_{i})_{i\in I}}^{*} } \ar[r]^-{\prod_{i\in I}F_{\{i\}}^{dn_{i}}}
 & \varinjlim _{\mu}\restr{\mc H _{N, I, W}^{0,\leq\mu,E}}{\Delta(\ov{v})} \ar[d]^-{\on{\mf{sp}}_{v}^{*}} 
 \\
\restr{(\prod_{i\in I}\Frob_{\{i\}}^{dn_{i}})^{*} ( \mc H _{N, I, W}^{0,\leq\mu_{1},E})}{\Delta(\ov \eta)}
 \ar[r]^-{\prod_{i\in I}F_{\{i\}}^{dn_{i}}}
& \varinjlim _{\mu}\restr{\mc H _{N, I, W}^{0,\leq\mu,E}}{\Delta(\ov{\eta})}  }\end{gather*}
 où la notation $ \on{\mf{sp}}_{v,(n_{i})_{i\in I}}^{*}$ indique que l'homomorphisme de spécialisation associé à la flèche $\on{\mf{sp}}_{v}:\Delta(\ov \eta)\to \Delta(\ov v)$ {\it est appliqué au  faisceau} $(\prod_{i\in I}\Frob_{\{i\}}^{dn_{i}})^{*} ( \mc H _{N, I, W}^{0,\leq\mu_{1},E})$ (et non pas à 
 $\mc H _{N, I, W}^{0,\leq\mu_{1},E}$). 
 Le diagramme précédent donne lieu à 
   \begin{gather*} \xymatrixcolsep{5pc} \xymatrix{
\wh b \ar@{|->}[d]^-{\on{\mf{sp}}_{v,(n_{i})_{i\in I}}^{*}} \ar@{|->}[r]^-{\prod_{i\in I}F_{\{i\}}^{dn_{i}}}
 &b_{(n_{i})_{i\in I}} \ar@{|->}[d]^-{\on{\mf{sp}}_{v}^{*}} 
 \\
\on{\mf{sp}}_{v,(n_{i})_{i\in I}}^{*}(\wh b) \ar@{|->}[r]^-{\prod_{i\in I}F_{\{i\}}^{dn_{i}}}
&a_{(n_{i})_{i\in I}}  }\end{gather*}

Donc  
pour montrer $a_{(n_{i})_{i\in I}} =0$ (et terminer ainsi la preuve de b)) il suffit de montrer que     
\begin{gather}\label{fibre-Frob-Delta}
  \on{\mf{sp}}_{v,(n_{i})_{i\in I}}^{*}(\wh b) \in 
\restr{ \mc H _{N, I, W}^{0,\leq\mu_{1},E}}{(\prod_{i\in I}\Frob_{\{i\}}^{dn_{i}})(\Delta(\ov \eta))}=
\restr{(\prod_{i\in I}\Frob_{\{i\}}^{dn_{i}})^{*} ( \mc H _{N, I, W}^{0,\leq\mu_{1},E})}{\Delta(\ov \eta)}\end{gather}
est nul. 
 Or \eqref{fibre-Frob-Delta} peut aussi être considéré comme l'image de $\wh b$ par un   homomorphisme de spécialisation pour le faisceau
 $\mc H _{N, I, W}^{0,\leq\mu_{1 },E}$ mais associé  à une flèche de spécialisation $(\prod_{i\in I}\Frob_{\{i\}}^{dn_{i}})(\Delta(\ov \eta))\to \Delta(\ov v)$. On en déduit que  \eqref{fibre-Frob-Delta} est nul car 
 \begin{itemize}
 \item 
  $\prod_{i\in I}\Frob_{\{i\}}^{dn_{i}}
(\Delta(\ov\eta))$ appartient à $ \Omega_{1}$ par hypothèse 
\item pour tout point géométrique $\ov x$ de 
$\Omega_{1}$
et toute flèche de spécialisation $\on{\mf{sp}}_{\ov x}: \ov x\to \Delta(\ov v)$,     $\on{\mf{sp}}_{\ov x}^{*}(\wh  b )$ s'annule dans 
$\restr{\mc H _{N, I, W}^{0,\leq\mu_{1 },E}}{\ov x}$. 
\end{itemize}
 Cette dernière assertion résulte du fait que   $\mc H _{N, I, W}^{0,\leq\mu_{1},E}$ est lisse sur $\Omega _{1}$ 
et que  l'image de $\wh  b$ par tout homomorphisme de spécialisation vers 
$\restr{\mc H _{N, I, W}^{0,\leq\mu_{1},E}}{\ov{\eta^{I}}}
$ est nulle (puisque c'est le cas de $\on{\mf{sp}}^{*}(\on{\mf{sp}}_{v}^{*}(\wh  b))$ et que 
$\pi_{1}(\eta^{I},\ov{\eta^{I}})$ agit transitivement sur les flèches de spécialisation de $\ov{\eta^{I}}$ vers $\Delta(\ov v)$).   \cqfd

\noindent{\bf Fin de la démonstration de la \propref{injectivite-sp}. }
 Comme $\prod _{i\in I}\Frob_{\{i\}}$ est le Frobenius total, 
 $\prod_{i\in I}F_{\{i\}}^{dn }$ agit de fa\c con bijective sur 
 $\varinjlim _{\mu}\restr{\mc H _{N, I, W}^{0,\leq\mu,E}}{\Delta(\ov{\eta})} $ et envoie $a_{(n_{i} )_{i\in I}}$ sur $a_{(n_{i}+n)_{i\in I}}$. De ceci et du a) du \lemref{lem-a-b} on  
   déduit facilement que pour montrer que $a=a_{(0)_{i\in I}}$ est nul 
   (et même que toute la suite $a_{(n_{i})_{i\in I}}$ est nulle) il suffit de  trouver 
 $(n_{i})_{i\in I}\in \N^{I}$ tel que 
 \begin{gather*} 
 a_{(n_{i}+\alpha_{i})_{i\in I}}=0 \text{ pour tout } (\alpha_{i})_{i\in I}\in \prod _{i\in I} \{0,..., \dim W_{i}-1\}.
 \end{gather*}
Or cela est possible d'après le b) du \lemref{lem-a-b}, car on peut trouver  $(n_{i})_{i\in I}\in \N^{I}$ tel que 
\begin{gather*} 
 \prod_{i\in I}\Frob_{\{i\}}^{d(n_{i}+\alpha_{i})}
(\Delta(\ov\eta))\in \Omega_{1} \text{ pour tout } (\alpha_{i})_{i\in I}\in \prod _{i\in I} \{0,..., \dim W_{i}-1\}.  \end{gather*} 
En effet  la densité de l'ouvert $\Omega_{1}$ implique  la densité de l'ouvert 
$$\bigcap_{ (\alpha_{i})_{i\in I}\in \prod _{i\in I} \{0,..., \dim W_{i}-1\}}
\Big(  \prod_{i\in I}\Frob_{\{i\}}^{d \alpha_{i} }\Big) ^{-1} (\Omega_{1})$$
et les $ \prod_{i\in I}\Frob_{\{i\}}^{d n_{i} }
(\Delta(\ov\eta))$ sont Zariski-denses lorsque $(n_{i})_{i\in I}$ parcourt $ \N^{I}$.  
 Ceci termine la preuve de la \propref{injectivite-sp}. \cqfd
 
 Les propositions  \ref{surjectivite-Hecke-fini} et \ref{injectivite-sp}   entraînent le corollaire suivant. 
 
 \begin{cor}\label{bijectivite-Hecke-fini}
L'homomorphisme de spécialisation 
 \begin{gather}\label{sp*-sans-Hf2}\on{\mf{sp}}^{*}: \Big( \varinjlim _{\mu}\restr{\mc H _{N, I, W}^{0,\leq\mu,E}}{\Delta(\ov{\eta})} \Big)^{\mr{Hf}}\to 
 \Big( \varinjlim _{\mu}\restr{\mc H _{N, I, W}^{0,\leq\mu,E}}{\ov{\eta^{I}}}\Big)^{\mr{Hf}}\end{gather}
 est une bijection.     \end{cor}
\dem
L'injectivité résulte de la \propref{injectivite-sp}. Voici la preuve de la surjectivité. Soit $c\in  \Big( \varinjlim _{\mu}\restr{\mc H _{N, I, W}^{0,\leq\mu,E}}{\ov{\eta^{I}}}\Big)^{\mr{Hf}}$. D'après la proposition \ref{surjectivite-Hecke-fini} il existe $a\in  \varinjlim _{\mu}\restr{\mc H _{N, I, W}^{0,\leq\mu,E}}{\Delta(\ov{\eta})}$ tel que $\on{\mf{sp}}^{*}(a)=c$. L'injectivité de $\on{\mf{sp}}^{*}$ montrée dans 
la \propref{injectivite-sp} implique que $a$ est Hecke-fini. \cqfd

                \section{Opérateurs d'excursion}\label{construction-S-proprietes}

        Soit  $I$ un ensemble fini et $W$ une représentation de $(\wh G)^{I}$. 
  Soit   $x\in W$ et  $\xi\in W^{*}$  invariants par l'action diagonale de $\wh G$. 
            Soit  $(\gamma_i)_{i\in I}\in 
( \pi_{1}(\eta,\ov\eta))^{I}$. On va rappeler (d'une fa\c con un peu différente) la construction des  opérateurs d'excursion 
 $$S_{I,W,x,\xi,(\gamma_i)_{i\in I}}\in \on{End}_{C_{c}(K_{N}\backslash G(\mb A)/K_{N},E)}(C_{c}^{\rm{cusp}}(\Bun_{G,N}(\Fq)/\Xi,E)) $$
 qui a déjà été expliquée dans l'introduction.  
 
  On rappelle que l'on possède le  morphisme de création
     $$\mc C_{  x}^{\sharp }
    :C_{c}(\Bun_{G,N}^{\leq\mu,E}(\Fq)/\Xi,E) \boxtimes E_{X\sm N} \to 
    \restr{ \mc H _{ N, I, W}^{0,\leq\mu,E}}{ \Delta(X\sm N)}$$
 qui est un morphisme de   $E$-faisceaux  constructibles  
sur $X\sm N$. 
   En le restreignant à  $\ov\eta$  et en passant à la  limite inductive   sur  $\mu$  on obtient le  morphisme 
 \begin{gather}\label{creation-limite-inductive}\mc C_{  x}^{\sharp } : C_{c}(\Bun_{G,N}(\Fq)/\Xi,E)\to   \varinjlim _{\mu} \restr{ \mc H _{ N, I, W}^{0,\leq\mu,E}}{\Delta(\ov\eta)}. \end{gather} 
 On rappelle que ce morphisme est la  composée
  \begin{gather*}C_{c}(\Bun_{G,N}(\Fq)/\Xi,E)= \varinjlim _{\mu}\restr{\mc H _{ N, \{0\},\mbf  1}^{0,\leq\mu,E}}{\ov \eta} \\
 \xrightarrow{\mc H(x)} 
\varinjlim _{\mu}\restr{\mc H _{ N, \{0\},W^{\zeta_{I}}}^{0,\leq\mu,E}}{\ov \eta}
 \isor{\chi_{\zeta_{I}}^{-1}}
\varinjlim _{\mu}\restr{\mc H _{ N, I,W}^{0,\leq\mu,E}}{\Delta(\ov \eta)},\end{gather*}  
D'après la  \propref{prop-cusp-hecke-finies} il envoie   $C_{c}^{\mr{cusp}}(\Bun_{G,N}(\Fq)/\Xi,E)$ dans 
 $\Big( \varinjlim _{\mu} \restr{ \mc H _{ N, I, W}^{0,\leq\mu,E}}{\Delta(\ov\eta)}\Big)^{\mr{Hf}}$. 

 De même on rappelle que l'on possède le  morphisme d'annihilation   
    $$\mc C_{  \xi}^{\flat }
    : 
    \restr{ \mc H _{ N, I, W}^{0,\leq\mu,E}}{ \Delta(X\sm N)}
    \to 
    C_{c}(\Bun_{G,N}^{\leq\mu}(\Fq)/\Xi,E) \boxtimes E_{X\sm N}$$
  qui est un morphisme de   $E$-faisceaux  constructibles  
sur $X\sm N$. 
   En le restreignant à  $\ov\eta$  et en passant à la  limite inductive   sur  $\mu$  on obtient le  morphisme $$ \mc C_{   
\xi}^{\flat }
:      \varinjlim _{\mu} \restr{\mc H _{ N, I, W}^{0,\leq\mu,E}}{\Delta(\ov\eta)} \to
        C_{c}(\Bun_{G,N}(\Fq)/\Xi,E).$$
On rappelle que ce morphisme est la  composée
      \begin{gather*} \varinjlim _{\mu}\restr{\mc H _{ N,I,W}^{0,\leq\mu,E}}{\Delta(\ov \eta)}
     \isor{\chi_{\zeta_{I}}}
\varinjlim _{\mu}\restr{\mc H _{ N, \{0\},W^{\zeta_{I}}}^{0,\leq\mu,E}}{\ov \eta}
 \\
\xrightarrow{\mc H(\xi )} \varinjlim _{\mu}\restr{\mc H _{ N, \{0\},\mbf  1}^{0,\leq\mu,E}}{\ov \eta}
= C_{c}(\Bun_{G,N}(\Fq)/\Xi,E).\end{gather*}
Il envoie 
$ \Big( \varinjlim _{\mu} \restr{\mc H _{ N, I, W}^{0,\leq\mu,E}}{\Delta(\ov\eta)}\Big)^{\mr{Hf}} 
   $ sur $    C_{c}^{\mr{cusp}}((\Bun_{G,N}(\Fq)/\Xi,E)  $. 
Enfin, d'après le
  \corref{bijectivite-Hecke-fini}, on a un isomorphisme
\begin{gather}
\label{isom-Hf-diag-gen-rappel}
\on{\mf{sp}}^{*}: \Big( \varinjlim _{\mu}\restr{\mc H _{N, I, W}^{0,\leq\mu,E}}{\Delta(\ov{\eta})} \Big)^{\mr{Hf}}\isom  
 \Big( \varinjlim _{\mu}\restr{\mc H _{N, I, W}^{0,\leq\mu,E}}{\ov{\eta^{I}}}\Big)^{\mr{Hf}}.\end{gather}
 On rappelle que la \propref{prop-action-Hf} munit le membre de droite de \eqref{isom-Hf-diag-gen-rappel} d'une action de $ \pi_{1}(\eta,\ov\eta)^{I}$ 
 (dépendant du choix de $\mf{sp}$). 
  Dans le \lemref{lem-HIW-action-indep-sp} ci-dessous on justifiera le fait (déjà expliqué dans l'introduction) que l'action de $ \pi_{1}(\eta,\ov\eta)^{I}$
  sur le membre de gauche de \eqref{isom-Hf-diag-gen-rappel} qui s'en déduit  ne dépend du choix de $\ov{\eta^{I}}$ et $\on{\mf{sp}}$. 
  Donc la construction suivante ne dépend pas du choix de $\ov{\eta^{I}}$ et $\on{\mf{sp}}$.

   \begin{def-prop}\label{antecedent-def-prop}
    On définit l'opérateur d'excursion    $S_{I,W,x,\xi,(\gamma_i)_{i\in I}}$ 
   comme la composée 
        \begin{gather}\label{diag-S}
    C_{c}^{\mr{cusp}}(\Bun_{G,N}(\Fq)/\Xi,E)
     \xrightarrow{ \mc C_{  x}^{\sharp}} 
   \Big( \varinjlim _{\mu} \restr{\mc H _{ N, I, W}^{0,\leq\mu,E}}{\Delta(\ov\eta)}\Big)^{\mr{Hf}}\\ 
   \nonumber    \xrightarrow{\mf{sp}^{*}}
   \Big( \varinjlim _{\mu}  \restr{\mc H _{ N, I, W}^{0,\leq\mu,E}}{ \ov{  \eta^{I}}}\Big)^{\mr{Hf}} 
    \xrightarrow{(\gamma_i)_{i\in I}}  
    \Big( \varinjlim _{\mu}  \restr{\mc H _{ N, I, W}^{0,\leq\mu,E}}{ \ov{  \eta^{I}}}\Big)^{\mr{Hf}} \\ \nonumber 
   \xrightarrow{(\mf{sp}^{*})^{-1}}
    \Big( \varinjlim _{\mu} \restr{\mc H _{ N, I, W}^{0,\leq\mu,E}}{\Delta(\ov\eta)}\Big)^{\mr{Hf}} 
    \xrightarrow{ \mc C_{  
\xi}^{\flat }}
    C_{c}^{\mr{cusp}}((\Bun_{G,N}(\Fq)/\Xi,E). 
    \end{gather}
      
    Alors 
    $S_{I,W,x,\xi,(\gamma_i)_{i\in I}}$
    appartient à  
 \begin{gather}\label{End-cusp-Xi}\on{End}_{C_{c}(K_{N}\backslash G(\mb A)/K_{N},E)}\big(C_{c}^{\rm{cusp}}(\Bun_{G,N}(\Fq)/\Xi,E)\big).\end{gather}
 De plus  il existe un  ouvert dense $U\subset X\sm N$ (dépendant seulement de   $I, W$ et $N$) tel que 
 \begin{itemize}
 \item $S_{I,W,x,\xi,(\gamma_i)_{i\in I}}$ dépend seulement de  l'image de  $(\gamma_i)_{i\in I}$ dans 
   $\pi_{1}(U, \ov\eta)^{I}$ 
   \item    $ (\gamma_i)_{i\in I}\mapsto S_{I,W,x,\xi,(\gamma_i)_{i\in I}}$ est continu du groupe profini  $\pi_{1}(U, \ov\eta)^{I}$ vers la $E$-algèbre de  dimension finie  
    \eqref{End-cusp-Xi} munie de la topologie   $E$-adique. 
    \end{itemize}
    \end{def-prop}
        \noindent{\bf Démonstration.} 
        L'opérateur     $S_{I,W,x,\xi,(\gamma_i)_{i\in I}}$ commute avec les opérateurs de Hecke  $T(f)$ pour $f\in C_{c}(K_{N}\backslash G(\mb A)/K_{N},E)$ car c'est le cas de tous les  morphismes apparaissant dans  sa définition.   Enfin  $\on{\mf{sp}}^{*}(\mc C_{  x}^{\sharp }
     (C_{c}^{\rm{cusp}}(\Bun_{G,N}(\Fq)/\Xi,E)))$ est  un sous-$E$-espace vectoriel de dimension finie   de $\Big( \varinjlim _{\mu} \restr{ \mc H _{ N, I, W}^{0,\leq\mu,E}}{\ov{\eta^{I}}}\Big)^{\mr{Hf}}$, donc  par la  \propref{prop-action-Hf} il existe un  ouvert dense $U\subset X$ tel que la dernière assertion soit vérifiée.   \cqfd

        On va donner dans la remarque suivante une caractérisation de 
        $S_{I,W,x,\xi,(\gamma_i)_{i\in I}}$ 
        indépendante du fait que ${\mf{sp}}^{*}$ est un isomorphisme sur les parties Hecke-finies (par le  \corref{bijectivite-Hecke-fini}). Dans la version 4 de cet article sur arXiv
         (où l'on savait seulement que l'image  ${\mf{sp}}^{*}$ contenait la partie Hecke-finie) cela était vraiment utilisé dans la construction de  $S_{I,W,x,\xi,(\gamma_i)_{i\in I}}$, alors que maintenant ce n'est plus qu'une remarque. 
        
     \begin{rem}    Grâce à la  \remref{creation-annihilation-dualite}  on a  pour tout $\mu$ un   morphisme de faisceaux  
 $$\mf B_{N,I,W}^{\Xi,E}:\mc H _{ N, I, W}^{0,\leq\mu,E}\otimes \mc H _{ N, I, W^{*,\theta}}^{0,\leq\mu,E}\to E_{(X\sm N)^{I}}.$$ En particulier 
 il fournit des  formes   bilinéaires  
  \begin{gather}\label{forme-bil-Delta-1}\s{.,.}_{\Delta(\ov\eta)} :  \Big( \varinjlim _{\mu} \restr{\mc H _{ N, I, W^{*,\theta}}^{0,\leq\mu,E}}{\Delta(\ov\eta)} \Big) \otimes 
 \Big( \varinjlim _{\mu}  \restr{ \mc H _{ N, I, W}^{0,\leq\mu,E}}{\Delta(\ov\eta)}\Big)\to E\end{gather} et 
 $$\s{.,.}_{\ov{\eta^{I}}} :  \Big( \varinjlim _{\mu} \restr{\mc H _{ N, I, W^{*,\theta}}^{0,\leq\mu,E}}{\ov{\eta^{I}}}\Big) \otimes 
 \Big( \varinjlim _{\mu} \restr{ \mc H _{ N, I, W}^{0,\leq\mu,E}}{\ov{\eta^{I}}} \Big)\to E$$
  telles  que \begin{gather}
  \label{pairing-sp-sp}
  \s{\on{\mf{sp}}^{*}(a),\on{\mf{sp}}^{*}(b)}_{\ov{\eta^{I}}}=\s{a,b}_{\Delta(\ov\eta)} \text{ pour }a\otimes b\text{ dans  le membre de gauche de \eqref{forme-bil-Delta-1}}.\end{gather}

  De la même fa\c con que   \eqref{creation-limite-inductive}, on définit     $$\mc C_{  \xi}^{\sharp }:  C_{c}(\Bun_{G,N}(\Fq)/\Xi,E)\to  \varinjlim _{\mu} \restr{\mc H _{ N, I, W^{*,\theta}}^{0,\leq\mu,E}}{\Delta(\ov\eta)} .$$
   D'après la  \remref{creation-annihilation-dualite}, pour tout 
$$\delta\in  \varinjlim _{\mu} \restr{\mc H _{ N, I, W}^{0,\leq\mu,E}}{\Delta(\ov\eta)} \text{ \  et \ } 
       \check h\in C_{c}(\Bun_{G,N}(\Fq)/\Xi,E),$$         on a  
      \begin{gather}\label{adj-check-h-delta}
      \s{\mc C_{  \xi}^{\sharp }
     (\check h), \delta}_{\Delta(\ov\eta)}=\int_{\Bun_{G,N}(\Fq)/\Xi}  \check  h \mc C_{  
\xi}^{\flat }( \delta).  \end{gather}

 Soit    $x\in W$ et  $\xi \in W^{*}$  invariants par l'action diagonale  de $\wh G$,  $$(\gamma_{i})_{i\in I}\in 
 ( \pi_{1}(\eta,\ov\eta))^{I} ,   h\in C_{c}^{\rm{cusp}}(\Bun_{G,N}(\Fq)/\Xi,E) \text{ et }   \check  h\in C_{c}(\Bun_{G,N}(\Fq)/\Xi,E). $$
Alors 
 \begin{gather}\label{S-bil-carac}
     \int_{\Bun_{G,N}(\Fq)/\Xi} \check h  S_{I,W,x,\xi,(\gamma_i)_{i\in I}}(h)
  \\ \nonumber  
 \overset{\eqref{diag-S}}{=}   \int_{\Bun_{G,N}(\Fq)/\Xi} \check h \Big( \mc C_{  \xi}^{\flat }\big((\on{\mf{sp}}^{*})^{-1 }((\gamma_i)_{i\in I}\cdot (\on{\mf{sp}}^{*}(\mc C_{  x}^{\sharp }
     (h))))\big)\Big)
   \\ \nonumber  
      \overset{\eqref{adj-check-h-delta}}{=}      \s{\mc C_{  \xi}^{\sharp }
     (\check h) , (\on{\mf{sp}}^{*})^{-1 }((\gamma_i)_{i\in I}\cdot (\on{\mf{sp}}^{*}(\mc C_{  x}^{\sharp }
     (h))))}_{\Delta(\ov\eta)} \\ \label{S-bil-carac2}
    \overset{\eqref{pairing-sp-sp}}{=}      \s{\on{\mf{sp}}^{*}(\mc C_{   \xi}^{\sharp }
     (\check h)) , (\gamma_i)_{i\in I}\cdot (\on{\mf{sp}}^{*}(\mc C_{  x}^{\sharp }
     (h)))}_{\ov{\eta^{I}}}. 
               \end{gather}

     L'égalité entre \eqref{S-bil-carac} et \eqref{S-bil-carac2} implique que $S_{I,W,x,\xi,(\gamma_i)_{i\in I}}$ est caractérisé par la forme bilinéaire \eqref{S-bil-carac2}  pour 
     $h, \check  h \in C_{c}^{\rm{cusp}}(\Bun_{G,N}(\Fq)/\Xi,E) $. 
     On remarque de plus que si $(\gamma_i)_{i\in I}$ appartient à 
     $\big(\on{Weil}(\eta,\ov{\eta})\big)^{I}
$, et si on le remplace dans \eqref{S-bil-carac2}
par un relevé arbitraire dans 
  $  \on{FWeil}(\eta^{I},\ov{\eta^{I}})$ par      la surjection 
     \eqref{mor-surj-Fweil-Weil^I}, la caractérisation ci-dessus s'exprime sans  même utiliser l'énoncé de la \propref{prop-action-Hf}, ni le \lemref{lem-Frob-partiels-drinfeld}. 
               \end{rem}
     
   On peut résumer la construction de 
     $S_{I,W,x,\xi,(\gamma_i)_{i\in I}}$
   sous la forme 
          $$   \xymatrix{
    \star   \ar[r]^{\mc C^{\sharp}_{x}}    &    
    \star   \ar[r]^{\mf{sp}^{*}} & 
     \star  \ar[r]^{\gamma} & 
    \star  \ar[r]^{(\mf{sp}^{*})^{-1}} &
     \star   \ar[r]^{\mc C^{\flat}_{\xi}} & 
     \star
    }$$           
    (avec la   notation $\gamma=(\gamma_i)_{i\in I}$).

     Le reste de ce  chapitre est consacré aux   propriétés des opérateurs d'excursion, qui sont les mêmes que   \eqref{SIW-p0-intro}, \eqref{SIW-p1-intro}, \eqref{SIW-p2-intro} et \eqref{SIW-p3-intro} dans l'introduction.

    D'après le \corref{bijectivite-Hecke-fini}, la restriction de l'homomorphisme $\on{\mf{sp}}^{*}$   aux parties Hecke-finies est  un isomorphisme 
   \begin{gather}\label{isom-sp*-avec-Hf} 
  \Big( \varinjlim _{\mu}\restr{\mc H _{N, I, W}^{0,\leq\mu,E}}{\Delta(\ov{\eta})}\Big)^{\mr{Hf}} \isor{\on{\mf{sp}}^{*}}
 \Big(  \varinjlim _{\mu}\restr{\mc H _{N, I, W}^{0,\leq\mu,E}}{\ov{\eta^{I}}}\Big)^{\mr{Hf}}. \end{gather}

   Comme dans l'introduction on introduit la notation courte $H_{I,W}$
    pour ce  $E$-espace vectoriel    (on omet $N$ dans la notation $ H_{I,W}$ pour limiter la taille des diagrammes dans le chapitre suivant). 
   
   \begin{defi}\label{def-HIW-texte}
    On définit $H_{I,W}$ comme le {\it membre de gauche} de \eqref{isom-sp*-avec-Hf}. 
  \end{defi}
 
 \begin{lem}
 \label{lem-HIW-action-indep-sp}
 L'action de $\on{Gal}(\ov F/F)^{I}=\pi_{1}(\eta,\ov\eta)^{I}$ sur 
 $H_{I,W}$ fournie par la \propref{prop-action-Hf}  ne dépend pas du choix de $\ov{\eta^{I}}$ et $\on{\mf{sp}}$. \end{lem}
 \dem En effet on peut reformuler ce qui précède en disant que d'après la \propref{prop-action-Hf}   on peut trouver 
 \begin{itemize}
 \item une réunion croissante 
 (indexée par $\lambda\in \N$) 
 de sous-$\mc O_{E}$-faisceaux constructibles 
 $\mf F_{\lambda}\subset \varinjlim _{\mu}\restr{\mc H _{N, I, W}^{0,\leq\mu,E}}{\eta^{I}}$ stables par les morphismes de Frobenius partiels
 \item une suite décroissante d'ouverts  denses $U_{\lambda}\subset X\sm N$ tels que 
 $\mf F_{\lambda}$ se prolonge en un faisceau lisse sur $(U_{\lambda})^{I}$
 \end{itemize}
 de sorte que \begin{gather}\label{egalite-Flambda-H}
 \bigcup _{\lambda\in \N}  \restr{\mf F_{\lambda}}{\ov{\eta^{I}}} 
 =
  \Big(  \varinjlim _{\mu}\restr{\mc H _{N, I, W}^{0,\leq\mu,E}}{\ov{\eta^{I}}}\Big)^{\mr{Hf}}.\end{gather} 
   Alors le \corref{bijectivite-Hecke-fini}   
   implique   que  le morphisme naturel 
  \begin{equation}\label{mor-HIW-Flambda}
  H_{I,W}= \Big( \varinjlim _{\mu}\restr{\mc H _{N, I, W}^{0,\leq\mu,E}}{\Delta(\ov{\eta})}\Big)^{\mr{Hf}} \to \bigcup _{\lambda\in \N}  \restr{\mf F_{\lambda}}{\Delta(\ov{\eta})} \end{equation}
(qui vient de la lissité de $\mf F_{\lambda}$ sur $(U_{\lambda})^{I}\ni 
\Delta(\ov{\eta})$)
  est un isomorphisme. 
  Plus précisément \eqref{mor-HIW-Flambda} est défini comme la composée
  \begin{gather*} 
  \Big( \varinjlim _{\mu}\restr{\mc H _{N, I, W}^{0,\leq\mu,E}}{\Delta(\ov{\eta})}\Big)^{\mr{Hf}}
  \isor{\mf{sp}^{*}} 
   \Big( \varinjlim _{\mu}\restr{\mc H _{N, I, W}^{0,\leq\mu,E}}{\ov{\eta^{I}}}\Big)^{\mr{Hf}}\overset{\eqref{egalite-Flambda-H}}{=}
\bigcup _{\lambda\in \N}  \restr{\mf F_{\lambda}}{\ov{\eta^{I}}} \isor{(\mf{sp}^{*})^{-1}} 
 \bigcup _{\lambda\in \N}  \restr{\mf F_{\lambda}}{\Delta(\ov{\eta})}
   \end{gather*}
où la flèche de gauche est un isomorphisme par le \corref{bijectivite-Hecke-fini} 
et celle de droite en est un par la lissité de $\mf F_{\lambda}$ sur $(U_{\lambda})^{I}\ni 
\Delta(\ov{\eta})$. Comme ces deux flèches sont définies à l'aide du même choix de $\mf{sp}$, la composée n'en dépend pas.

  Or l'action de $\on{Gal}(\ov F/F)^{I}$   sur le membre de droite de \eqref{mor-HIW-Flambda}, obtenue en appliquant le lemme de Drinfeld (le \lemref{lem-Frob-partiels-drinfeld}) à $\mf F_{\lambda}$, ne dépend pas  du choix de $\ov{\eta^{I}}$ et $\on{\mf{sp}}^{*}$, et donc l'action  de $\on{Gal}(\ov F/F)^{I}$ sur le membre de gauche n'en dépend pas non plus. \cqfd
  
 \begin{rem}
 Dans cet article nous montrons seulement que  $H_{I,W}$ est une limite inductive de $E$-espaces vectoriels de dimension finie munis de représentations continues de $\on{Gal}(\ov F/F)^{I}$. En fait  Cong Xue a montré dans \cite{these-cong} que 
 $H_{I,W}$ est de dimension finie. 

 \end{rem}
 
   Pour toute application  $\zeta:I\to J$, l'isomorphisme de coalescence 
   \eqref{intro-isom-coalescence} respecte trivialement les  parties Hecke-finies et induit donc  un isomorphisme 
    \begin{gather} \label{isom-chi-zeta}
    H_{I,W}= \Big( \varinjlim _{\mu}\restr{\mc H _{N, I, W}^{0,\leq\mu,E}}{\Delta(\ov{\eta})}\Big)^{\mr{Hf}}
    \isor{\chi_{\zeta}} 
    \Big( \varinjlim _{\mu}\restr{\mc H _{N, J,W^{\zeta}}^{0,\leq\mu,E}}{\Delta(\ov{\eta})}\Big)^{\mr{Hf}}=H_{J,W^{\zeta}}
    \end{gather}
    où $\Delta$ désigne le morphisme diagonal $X\to X^{I}$ ou $X\to X^{J}$. 
    
    La définition ci-dessous reprend celle déjà donnée dans  \eqref{def-chi-zeta-intro} de l'introduction. 
    
    \begin{defi} On définit l'isomorphisme de coalescence
    $$ \chi_{\zeta}:    H_{I,W} \isom H_{J,W^{\zeta}}$$  par \eqref{isom-chi-zeta}. 
    \end{defi}
    
    L'isomorphisme \eqref{isom-chi-zeta} est $\on{Gal}(\ov F/F)^{J}$-équivariant,  où $\on{Gal}(\ov F/F)^{J}$ agit sur le membre de gauche par le   morphisme diagonal  
\begin{gather}\nonumber 
\on{Gal}(\ov F/F)^{J}\to \on{Gal}(\ov F/F)^{I},  \ (\gamma_{j})_{j\in J}\mapsto (\gamma_{\zeta(i)})_{i\in I}. 
\end{gather}
En effet, soit $\Delta_{\zeta}:X^{J}\to X^{I}$  le morphisme diagonal 
    \eqref{morph-giad-X-intro}. En appliquant le \corref{bijectivite-Hecke-fini} à $I$ et $J$ on voit que, pour toute flèche de spécialisation de $\ov{\eta^{I}}$ vers 
    $\Delta_{\zeta}(\ov{\eta^{J}})$, le morphisme image inverse 
    $ \Big( \varinjlim _{\mu}\restr{\mc H _{N, I, W}^{0,\leq\mu,E}}{\Delta_{\zeta}(\ov{\eta^{J}})}\Big)^{\mr{Hf}}
  \to    \Big( \varinjlim _{\mu}\restr{\mc H _{N, I, W}^{0,\leq\mu,E}}{\ov{\eta^{I}}}\Big)^{\mr{Hf}}$ est un isomorphisme. Donc 
    si la suite  $(\mf F_{\lambda})_{\lambda\in \N}$ est comme ci-dessus relativement à $I$ et $W$, alors la suite  $(\Delta_{\zeta}^{*}(\mf F_{\lambda}))_{\lambda\in \N}$ vérifie les mêmes propriétés relativement à $J$ et $W^{\zeta}$, donc 
       $$\chi_{\zeta}: H_{I,W}\overset{\eqref{mor-HIW-Flambda}}{\simeq}\bigcup _{\lambda\in \N}  \restr{\mf F_{\lambda}}{\Delta(\ov{\eta})} =\bigcup _{\lambda\in \N}  \restr{\Delta_{\zeta}^{*}(\mf F_{\lambda})}{\Delta(\ov{\eta})}
       \overset{\eqref{mor-HIW-Flambda}}{\simeq}
       H_{J,W^{\zeta}}
       $$ est $\on{Gal}(\ov F/F)^{J}$-équivariant
       (puisque les actions de $\on{Gal}(\ov F/F)^{I}$ et $\on{Gal}(\ov F/F)^{J}$ 
       sur les termes de gauche et de droite 
       ont été défines en appliquant le lemme de Drinfeld aux termes centraux, et que l'égalité centrale est clairement $\on{Gal}(\ov F/F)^{J}$-équivariante). 
    
    La proposition suivante est identique à la \propref{prop-a-b-c}. 
    
      \begin{prop}\label{prop-a-b-texte}   Les $H_{I,W}$ vérifient   les 
  propriétés   suivantes : 
        \begin{itemize}
    \item[] {\bf a)} pour tout ensemble fini   $I$,      $$W\mapsto  
    H_{I,W},  \ \ u\mapsto \mc H(u)$$  est un foncteur  $E$-linéaire  de la  catégorie des  représentations $E$-linéaires de dimension finie de  $(\wh G)^{I}$ vers la  catégorie des  limites inductives de représentations $E$-linéaires continues de dimension finie de     $\on{Gal}(\ov F/F)^{I}$,    
              \item[] {\bf b)} pour toute application   $\zeta: I\to J$, 
 on possède  un isomorphisme 
      \begin{gather}\nonumber
           \chi_{\zeta}: H_{I,W}\isom 
 H_{J,W^{\zeta}},\end{gather} 
 qui est 
 \begin{itemize}
 \item  fonctoriel en   $W$, où  $W$ est une  représentation de $(\wh G)^{I}$ et  $W^{\zeta}$ désigne la   représentation de $(\wh G)^{J}$ sur  $W$ obtenue en composant avec le  morphisme  diagonal $$  (\wh G)^{J}\to (\wh G)^{I}, (g_{j})_{j\in J}\mapsto (g_{\zeta(i)})_{i\in I} $$ 
 \item $\on{Gal}(\ov F/F)^{J}$-équivariant, où $\on{Gal}(\ov F/F)^{J}$ agit sur le membre de gauche par le   morphisme diagonal  
 \begin{gather}\nonumber 
\on{Gal}(\ov F/F)^{J}\to \on{Gal}(\ov F/F)^{I},  \ (\gamma_{j})_{j\in J}\mapsto (\gamma_{\zeta(i)})_{i\in I}, 
\end{gather}
 \item   et compatible avec la  composition, c'est-à-dire   que pour 
 $I\xrightarrow{\zeta} J\xrightarrow{\eta} K$ on a 
 $\chi_{\eta\circ \zeta}=\chi_{\eta}\circ\chi_{\zeta}$,
    \end{itemize}
    \item[] {\bf c)} pour $I=\emptyset$ et  $W=\mbf  1$, on a un isomorphisme     \begin{gather}\nonumber
       H_{\emptyset,\mbf  1}=C_{c}^{\mr{cusp}}(G(F)\backslash G(\mb A)/K_{N}\Xi,E). \end{gather}
    \end{itemize}
    
  Par ailleurs les $H_{I,W}$ sont des modules sur 
  $C_{c}(K_{N}\backslash G(\mb A)/K_{N},E)$, de fa\c con compatible  avec 
les propriétés a), b), c) ci-dessus. 
  \end{prop}

\noindent {\bf Démonstration de la \propref{prop-a-b-texte}.  }    Les propriétés a) et b) ont déjà été expliquées. 
   En appliquant   b) à l'application évidente  $\zeta_{\emptyset}: \emptyset \to \{0\}$, on obtient un  isomorphisme \begin{gather}
    \label{isom-chi-singleton-empty}\chi_{\zeta_{\emptyset}}: H_{\emptyset,\mbf 1}\isom 
 H_{\{0\},\mbf 1} \end{gather} 
que l'on connaissait  déjà comme conséquence de \eqref{egalite-cht-singl-1}.  
 La propriété c) résulte de la \propref{prop-cusp-hecke-finies}.  \cqfd

   \section{Propriétés des opérateurs d'excursion}

     Soit  $I$   un ensemble fini et  $W$ une représentation  $E$-linéaire de  
    $(\wh G)^{I}$. 
    On note $\zeta_{I}:I\to \{0\}$ l'application évidente, si bien que 
    $ W^{\zeta_{I}}$ est simplement  $W$ muni de l'action diagonale  de $\wh G$. Soit 
$x: \mbf 1\to W^{\zeta_{I}}$ et  $\xi :  W^{\zeta_{I}}\to  \mbf 1$
des morphismes de  représentations de  $\wh G$ (autrement dit  $x\in W$ et  $\xi\in W^{*}$  sont invariants sous  l'action diagonale de  $\wh G$). Soit   $(\gamma_{i})_{i\in I}\in \on{Gal}(\ov F/F)^{I}$. 
 
 Il résulte immédiatement de la 
définition-proposition \ref{antecedent-def-prop}
 que  l'opérateur  d'excursion 
  \begin{gather*}S_{I,W,x,\xi,(\gamma_{i})_{i\in I}}\in 
  \on{End}(H_{\{0\},\mbf 1})  \end{gather*}
 est égal à  la composée 
  \begin{gather}\label{excursion-def-texte}
  H_{\{0\},\mbf  1}\xrightarrow{\mc H(x)}
 H_{\{0\},W^{\zeta_{I}}}\isor{\chi_{\zeta_{I}}^{-1}} 
  H_{I,W}
  \xrightarrow{(\gamma_{i})_{i\in I}}
  H_{I,W} \isor{\chi_{\zeta_{I}}} H_{\{0\},W^{\zeta_{I}}}  
  \xrightarrow{\mc H(\xi)} 
  H_{\{0\},\mbf  1}.
    \end{gather} 
    C'était la définition qui avait été donnée dans l'introduction. 
    Le \lemref{lem-HIW-action-indep-sp} montre que 
 $S_{I,W,x,\xi,(\gamma_{i})_{i\in I}}$ ne dépend pas du choix de $\ov{\eta^{I}}$ et $\on{\mf{sp}}^{*}$. 

   Le lemme suivant (qui reprend le \lemref{lem-intro-ptes-SIW} de l'introduction) va résulter des propriétés a) et  b) de la \propref{prop-a-b-texte}. 
   \begin{lem} 
   \label{lem-ptes-SIW}
   Les  opérateurs d'excursion $S_{I,W,x,\xi,(\gamma_{i})_{i\in I}}$ vérifient 
    les propriétés suivantes : 
     \begin{gather}
     \label{SIW-p0}      S_{I,W,x,{}^{t} u(\xi'),(\gamma_i)_{i\in I}}=S_{I,W',u(x),\xi',(\gamma_i)_{i\in I}} 
 \end{gather}
où $u:W\to W'$ est un   morphisme $(\wh G)^{I}$-équivariant et  $x\in W$ et $\xi'\in (W')^{*}$ sont $\wh G$-invariants,  
  \begin{gather}   \label{SIW-p1}
     S_{J,W^{\zeta},x,\xi,(\gamma_j)_{j\in J}}=S_{I,W,x,\xi,(\gamma_{\zeta(i)})_{i\in I}},
     \\
      \label{SIW-p2}
 S_{I_{1}\cup I_{2},W_{1}\boxtimes W_{2},x_{1}\boxtimes x_{2},\xi_{1}\boxtimes \xi_{2},(\gamma^{1}_i)_{i\in I_{1}}\times (\gamma^{2}_i)_{i\in I_{2}}}= S_{I_{1},W_{1},x_{1},\xi_{1},(\gamma^{1}_i)_{i\in I_{1}}}\circ 
S_{I_{2},W_{2},x_{2},\xi_{2},(\gamma^{2}_i)_{i\in I_{2}}}, 
\\ 
\label{SIW-p3} 
S_{I,W,x,\xi,(\gamma_i(\gamma'_i )^{-1}\gamma''_i)_{i\in I}}=
    S_{I\cup I \cup I,W\boxtimes W^{*}\boxtimes W,\delta_{W} \boxtimes x,
    \xi \boxtimes \on{ev}_{W},
    (\gamma_i)_{i\in I} \times (\gamma'_i)_{i\in I} \times (\gamma''_i)_{i\in I}
    }
             \end{gather}
  où la plupart des  notations sont évidentes, 
    $I_{1}\cup I_{2}$ et  $I\cup I\cup I$ désignent des réunions disjointes, et  $\delta_{W}: \mbf 1\to W\otimes W^{*}$ et $\on{ev}_{W}:W^{*} \otimes W \to \mbf 1$ sont les  morphismes naturels.      
  \end{lem}
   
     \noindent {\bf Démonstration de  \eqref{SIW-p0}.}  
    On pose  $x'=u(x)\in W'$ et $\xi={}^{t}u(\xi')\in W^{*}$. 
Le diagramme 
$$   \xymatrix{
     &  H_{\{0\},(W')^{\zeta_{I}}}  
    \ar[r]^{\chi_{\zeta_{I}}^{-1}}&  H_{I,W'} \ar[r]^{(\gamma_{i})_{i\in I}} 
     &
    H_{I,W'} \ar[r]^{\chi_{\zeta_{I}}} &  
    H_{\{0\},(W')^{\zeta_{I}}}  \ar[dr]^{\mc H(\xi')} 
     \\
  H_{\{0\},\mbf  1}  \ar[r]_{\mc H(x)}    \ar[ur]^{\mc H(x')}  &    
    H_{\{0\},W^{\zeta_{I}}}   \ar[u]_{\mc H(u)}
    \ar[r]_{\chi_{\zeta_{I}}^{-1}} & 
    H_{I,W}  \ar[u]_{\mc H(u)} \ar[r]_{(\gamma_{i})_{i\in I}} & 
   H_{I,W}  \ar[u]^{\mc H(u)} \ar[r]_{\chi_{\zeta_{I}}} &
     H_{\{0\},W^{\zeta_{I}}} \ar[u]^{\mc H(u)} \ar[r]_{\mc H(\xi)} & 
     H_{\{0\},\mbf  1}
    }$$           
est commutatif. Or la ligne du bas est égale à  
   $S_{I,W,x,\xi,(\gamma_{i})_{i\in I}}
$ et celle du haut est égale à 
   $S_{I,W',x',\xi',(\gamma_{i})_{i\in I}}
$. \cqfd
       
   \noindent {\bf Démonstration  de  \eqref{SIW-p1}.}  
Le diagramme 
$$   \xymatrix{
     & &  H_{J,W^{\zeta}} \ar[rr]^{(\gamma_{j})_{j\in J}}  &&
    H_{J,W^{\zeta}} \ar[dr]^{\chi_{\zeta_{J}}}
     \\
  H_{\{0\},\mbf  1}  \ar[r]_{\mc H(x)}    &    
    H_{\{0\},W^{\zeta_{I}}}   \ar[ur]^{\chi_{\zeta_{J}}^{-1}} 
    \ar[r]_{\chi_{\zeta_{I}}^{-1}} & 
    H_{I,W}  \ar[u]_{\chi_{\zeta}} \ar[rr]_{(\gamma_{\zeta(i)})_{i\in I}} && 
   H_{I,W}  \ar[u]^{\chi_{\zeta}} \ar[r]_{\chi_{\zeta_{I}}} &
     H_{\{0\},W^{\zeta_{I}}}  \ar[r]_{\mc H(\xi)} & 
     H_{\{0\},\mbf  1}
    }$$           
est commutatif. Or la ligne du bas est égale à 
   $S_{I,W,x,\xi,(\gamma_{\zeta(i)})_{i\in I}}
$ et celle du haut est égale à 
   $    S_{J,W^{\zeta},x,\xi,(\gamma_j)_{j\in J}}$.  \cqfd
       
   \noindent {\bf Démonstration  de  \eqref{SIW-p2}.}  
   L'application évidente 
   $\{0\}\cup \{0\} \to \{0\}$ 
  donne un  isomorphisme $ H_{\{0\}\cup \{0\},\mbf  1}\simeq  H_{\{0\},\mbf  1}$. En notant   $\zeta_{1}:I_{1}\to \{0\}$ et  $\zeta_{2}:I_{2}\to \{0\}$ les applications évidentes,  le membre de gauche de  \eqref{SIW-p2} est égal à la composée 
       \begin{gather*}
 H_{\{0\}\cup \{0\},\mbf  1}
 \xrightarrow{\mc H(x_{1}\boxtimes x_{2})}
 H_{\{0\}\cup \{0\},W_{1}^{\zeta_{1}}\boxtimes W_{2}^{\zeta_{2}}}
 \isor{\chi_{\zeta_{1}\times \zeta_{2}}^{-1}} 
  H_{I_{1}\cup I_{2},W_{1}\boxtimes W_{2}}\\
  \xrightarrow{(\gamma^{1}_{i})_{i\in I_{1}}\times (\gamma^{2}_{i})_{i\in I_{2}}}
  H_{I_{1}\cup I_{2},W_{1}\boxtimes W_{2}} \isor{\chi_{\zeta_{1}\times \zeta_{2}}} 
  H_{\{0\}\cup \{0\},W_{1}^{\zeta_{1}}\boxtimes W_{2}^{\zeta_{2}}}  
  \xrightarrow{\mc H(\xi_{1}\boxtimes \xi_{2})} 
  H_{\{0\}\cup \{0\},\mbf  1}. 
   \end{gather*}
En regroupant  $x_{1},\chi_{\zeta_{1}}^{-1},(\gamma^{1}_{i})_{i\in I_{1}} ,\chi_{\zeta_{1}},\xi_{1}$ d'un côté  et 
$x_{2},\chi_{\zeta_{2}}^{-1},(\gamma^{2}_{i})_{i\in I_{2}} ,\chi_{\zeta_{2}},\xi_{2}$ de l'autre on trouve  le membre de droite. 
On peut le faire car  dans le diagramme  suivant  (où l'on note 
$\gamma^{1}=(\gamma^{1}_{i})_{i\in I_{1}}$ et $\gamma^{2}=(\gamma^{2}_{i})_{i\in I_{2}}$) 
tous les carrés et les triangles commutent. 
\cqfd

   {  \resizebox{14cm}{!}{ $$\!\!\!\!\!\!\!\!\!\!\!\!\!\!\!    \xymatrix{
     {  H_{\{0\}\cup \{0\},\mbf 1} }
           \ar[d]_{{\mc H(x_{1}\boxtimes 1)} }
           \ar[rd]^-{{\mc H(x_1\boxtimes x_2)}}
             \\
      {    H_{\{0\}\cup \{0\},W_{1}^{\zeta_{1}}\boxtimes \mbf 1}}
           \ar[d]_{{\chi_{\zeta_1\times \Id}^{-1}}}
           \ar[r]^<<<{{\mc H(\Id\boxtimes x_{2} )} }
       &  { H_{\{0\}\cup \{0\},W_{1}^{\zeta_{1}}\boxtimes W_{2}^{\zeta_{2}}}}
        \ar[d]_{{\chi_{\zeta_1\times \Id}^{-1}} }
        \ar[rd]^-{{\chi_{\zeta_1\times \zeta_{2}}^{-1}}}
            \\
    {     H_{I_{1}\cup \{0\},W_{1}\boxtimes \mbf 1} }
         \ar[d]_{{\gamma^{1}\times 1}}
         \ar[r]^-{{ \mc H(\Id\boxtimes x_{2} )} }
       &   { H_{I_{1}\cup \{0\},W_{1} \boxtimes W_{2}^{\zeta_{2}}}}
        \ar[d]_{{\gamma^{1}\times 1}}
        \ar[r]^-{{\chi_{\Id\times \zeta_2}^{-1}} }
       &{  H_{I_{1}\cup I_{2},W_{1} \boxtimes W_{2} }  }
         \ar[d]_{{\gamma^{1}\times 1} }
         \ar[rd]^-{ {\gamma^{1}\times \gamma^{2}}}
   \\
   { H_{I_{1}\cup \{0\},W_{1}\boxtimes \mbf 1} }
    \ar[d]_{{\chi_{\zeta_1\times \Id}}}
  \ar[r]^-{{\mc H(\Id\boxtimes x_{2} )} }&  
  { H_{I_{1}\cup \{0\},W_{1} \boxtimes W_{2}^{\zeta_{2}}} }
    \ar[d]_{{\chi_{\zeta_1\times \Id}}}
  \ar[r]^-{{\chi_{\Id\times \zeta_2}^{-1}} }
  & 
 {   H_{I_{1}\cup I_{2},W_{1} \boxtimes W_{2}} }
    \ar[d]_{{\chi_{\zeta_1\times \Id}}}
   \ar[r]^-{{ 1\times  \gamma^{2}}  }
    &{  H_{I_{1}\cup I_{2},W_{1} \boxtimes W_{2}}   }
     \ar[d]_{{\chi_{\zeta_1\times \Id}}}
   \ar[rd]^-{{\chi_{\zeta_1\times \zeta_{2}}}       }
 \\
  {  H_{\{0\}\cup \{0\},W_{1}^{\zeta_{1}}\boxtimes \mbf 1} }
  \ar[d]_{{\mc H(\xi_{1}\boxtimes 1)}  }
  \ar[r]^<<<{{\mc H(\Id\boxtimes x_{2} )} }
  &  { H_{\{0\}\cup \{0\},W_{1}^{\zeta_{1}} \boxtimes W_{2}^{\zeta_{2}}}}
     \ar[d]_{{\mc H(\xi_{1}\boxtimes \Id)} } 
  \ar[r]^-{{\chi_{\Id\times \zeta_2} ^{-1}}}
  & { H_{\{0\}\cup I_{2},W_{1}^{\zeta_{1}} \boxtimes W_{2}}}
    \ar[r]^-{{ 1\times  \gamma^{2}}  }
     \ar[d]_{{\mc H(\xi_{1}\boxtimes \Id)} } &
 {  H_{\{0\}\cup I_{2},W_{1}^{\zeta_{1}} \boxtimes W_{2}}}
      \ar[d]_{{\mc H(\xi_{1}\boxtimes \Id)}  }
       \ar[r]^-{{\chi_{\Id\times \zeta_{2}}}   }
  & {  H_{\{0\}\cup \{0\},W_{1}^{\zeta_{1}} \boxtimes W_{2}^{\zeta_{2}} }}
  \ar[d]_{{\mc H(\xi_{1}\boxtimes \Id)}  }
 \ar[rd]^-{{\mc H(\xi_{1} \boxtimes \xi_2)} }
 \\           
 { H_{\{0\}\cup \{0\},\mbf 1} }
  \ar[r]^-{{ \mc H(1\boxtimes x_{2} )}} &   
 { H_{\{0\}\cup \{0\},\mbf 1  \boxtimes W_{2}^{\zeta_{2}}} }
  \ar[r]^-{{ \chi_{\Id\times \zeta_2} ^{-1}}}&  
{  H_{\{0\}\cup I_{2},\mbf 1  \boxtimes W_{2}} }
   \ar[r]^-{{ 1\times  \gamma^{2}}   } &  
 { H_{\{0\}\cup I_{2},\mbf 1  \boxtimes W_{2}} }
 \ar[r]^-{{\chi_{\Id\times \zeta_{2}}}   }
&{  H_{\{0\}\cup \{0\},\mbf 1  \boxtimes W_{2}^{\zeta_{2}}} }
\ar[r]^-{{\mc H(1\boxtimes \xi_{2} )}} & 
{ H_{\{0\}\cup \{0\},\mbf 1}} }$$}
  }
          
  \vskip1mm        
  \noindent{\bf Démonstration  de   \eqref{SIW-p3}.}  
   Pour tout  $(g_i)_{i\in I}\in (\wh G)^{I}$, 
\begin{itemize}
\item   $ \xi \boxtimes \on{ev}_{W}$     est  invariant par $(1)_{i\in I}\times (g_i)_{i\in I} \times  (g_i)_{i\in I} $
\item $\delta_{W} \boxtimes x$   est  invariant par      $(g_i)_{i\in I} \times (g_i)_{i\in I} \times  (1)_{i\in I} $.
    \end{itemize}    
Donc pour tous   $(\alpha_{i})_{i\in I}$ et $(\beta_{i})_{i\in I}$ dans $\on{Gal}(\ov F/F)^{I}$, le membre de droite de \eqref{SIW-p3} est égal à \begin{gather}\label{SIW-p3-alpha-beta}
 S_{I\cup I \cup I,W\boxtimes W^{*}\boxtimes W,\delta_{W} \boxtimes x,
    \xi \boxtimes \on{ev}_{W},
    (\gamma_i\beta_{i})_{i\in I} \times (\alpha_{i}\gamma'_i\beta_{i})_{i\in I} \times (\alpha_{i}\gamma''_i)_{i\in I}
   . } \end{gather}
  Pour le montrer de fa\c con formelle on factorise  le membre de droite de \eqref{SIW-p3} à travers 
  $$H_{I, \mbf 1} \xrightarrow{ \mc H(\delta_{W})}
     H_{ I,(W\boxtimes W^{*})^{\zeta}} \text{ \  et  \   } H_{ I,(W^{*}\boxtimes W)^{\zeta}} \xrightarrow{ \mc H(\on{ev}_{W})}
    H_{I, \mbf 1},$$ où  $\zeta:I\cup I\to I$ est l'application évidente,   
et on utilise le fait que   $\on{Gal}(\ov F/F)^{I}$ agit trivialement sur  $ H_{I,\mbf 1}\simeq  H_{\emptyset,\mbf 1}$. On prend $\alpha_{i}=\gamma_i (\gamma'_i)^{-1}$ et 
$\beta_{i}=(\gamma'_i)^{-1}\gamma''_i $. Alors   \eqref{SIW-p3-alpha-beta} est égal à  
           \begin{gather}\label{eq-lem-b-2-ggg-texte}
                 S_{I\cup I \cup I,W\boxtimes W^{*}\boxtimes W,\delta_{W} \boxtimes x,\xi \boxtimes \on{ev}_{W},
                 (\gamma_i(\gamma'_i )^{-1}\gamma''_i)_{i\in I} \times 
                 (\gamma_i(\gamma'_i )^{-1}\gamma''_i)_{i\in I} \times
                  (\gamma_i(\gamma'_i )^{-1}\gamma''_i)_{i\in I}
                  }.                   \end{gather}
En appliquant   \eqref{SIW-p1} à  l'application évidente  $\zeta:I\cup I \cup I\to I$, on voit que  \eqref{eq-lem-b-2-ggg-texte} est égal à 
  \begin{gather}\label{eq-lem-b-2-ggg2-texte}S_{I,W\otimes W^{*}\otimes W,\delta_{W} \otimes x,\xi \otimes \on{ev}_{W},(\gamma_i(\gamma'_i )^{-1}\gamma''_i)_{i\in I}}.\end{gather}
Finalement on montre que 
 \eqref{eq-lem-b-2-ggg2-texte} est égal au  membre de gauche de \eqref{SIW-p3}
 en appliquant  \eqref{SIW-p0} 
 à   l'injection  $(\wh G)^{I}$-linéaire $$u:W=\mbf  1\otimes W
    \xrightarrow{  \delta_{W}    \otimes\Id_{W}} W\otimes W^{*}\otimes W, $$
        qui vérifie   $\delta_{W} \otimes x =u( x)$ et  
                  ${}^{t}u(\xi\otimes\on{ev}_{W}) =\xi$, puisque la composée                   $$ W
    \xrightarrow{  \delta_{W}    \otimes\Id_{W}} W\otimes W^{*}\otimes W
     \xrightarrow{ \Id_{W}\otimes\on{ev}_{W}} W$$
                  est égale à $\Id_{W}$ d'après le lemme de Zorro \eqref{zorro}.             \cqfd

                  Le lemme suivant affirme que les   opérateurs de Hecke en les places de $X\sm N$ sont des cas particuliers d'opérateurs d'excursion. Il sera utilisé pour montrer que la  décomposition \eqref{intro1-dec-canonique} (que nous construirons dans le prochain chapitre) est compatible avec  l'isomorphisme de Satake en les places non ramifiées. 
                            
      \begin{lem}\label{S-non-ram}
      Soit  $v\in |X|\sm N$. 
      On fixe  un plongement $\ov F\subset \ov {F_{v}}$. Soit  $d\in \N$,  et 
    $\gamma\in \on{Gal}(\ov {F_{v}}/F_{v})\subset  \on{Gal}(\ov F/F)$ tel que $\deg(\gamma)=d$.  
    Alors  $S_{\{1,2\},V \boxtimes V^{*},\delta_{V},\on{ev}_{V},(\gamma,1)}$   dépend seulement de  $d$, et si  $d=1$ il est égal à 
    $T(h_{V,v})$.      \end{lem}    
    \noindent{\bf Démonstration.}
 On fixe un  point géométrique  $\ov v$ au-dessus de  $v$ et une  flèche de spécialisation 
    $\on{\mf{sp}}_{v}:\ov \eta\to \ov v$, associés au plongement $\ov F\subset \ov {F_{v}}$ choisi dans l'énoncé. On note encore   $\on{\mf{sp}}_{v}$ la  flèche de spécialisation de  
    $\Delta(\ov \eta)$ vers   $\Delta(\ov v)$ égale à l'image par $\Delta$ de cette dernière. Pour que le diagramme suivant tienne dans la page on  pose   $I=\{1,2\}$ et $W=V\boxtimes V^{*}$. 
    Dans le diagramme  
       $$ \xymatrixcolsep{1pc} \xymatrix{  C_{c}^{\mr{cusp}}(\Bun_{G,N}(\Fq)/\Xi,E)\ar[d]_-{\restr{\mc C_{ \delta_{V}}^{\sharp }}{\ov v}} \ar[dr]^{ \mc C_{  \delta_{V}}^{\sharp }} & &
       \\ \Big( \varinjlim _{\mu} \restr{\mc H _{ N, I, W}^{0,\leq\mu,E}}{\Delta(\ov v)}\Big)^{\mr{Hf}} \ar[r]^-{\mf{sp}_{v}^{*}} \ar[d]^{F_{\{1\}}^{\deg(v)d}} & 
       \Big( \varinjlim _{\mu} \restr{\mc H _{ N, I, W}^{0,\leq\mu,E}}{\Delta(\ov\eta)}\Big)^{\mr{Hf}}\ar@{=}[r] & H_{I,W} \ar[d]^{(\gamma,1)}
       \\  \Big( \varinjlim _{\mu} \restr{\mc H _{ N, I, W}^{0,\leq\mu,E}}{\Delta(\ov v)}\Big)^{\mr{Hf}} \ar[r]^-{\mf{sp}_{v}^{*}}\ar[d]_-{\restr{\mc C_{  \on{ev}_{V}}^{\flat }}{\ov v}}   & 
        \Big( \varinjlim _{\mu} \restr{\mc H _{ N, I, W}^{0,\leq\mu,E}}{\Delta(\ov\eta)}\Big)^{\mr{Hf}} \ar@{=}[r]\ar[dl]^{ \mc C_{  \on{ev}_{V}}^{\flat }} & 
       H_{I,W}
       \\ C_{c}^{\mr{cusp}}(\Bun_{G,N}(\Fq)/\Xi,E)  & &
       } $$
       la commutativité des triangles est évidente, et celle du grand rectangle 
     résulte   du lemme suivant appliqué à $\iota=1$ dans  $I=\{1,2\}$.           
     On en déduit que  $S_{\{1,2\},V \boxtimes V^{*},\delta_{V},\on{ev}_{V},(\gamma,1)}$ 
     (qui par construction est la composée suivant le chemin le plus à droite) 
     est égal à la  composée donnée par la colonne de gauche. 
     Par conséquent  il dépend seulement de   $d$,  
     et par  la \propref{prop-coal-frob-cas-part}  il est égal à   
      $T(h_{V,v})$ si $d=1$. \cqfd
    
    \begin{rem} 
   Par un calcul formel d'algèbre tensorielle, on peut déterminer 
  les opérateurs   $S_{\{1,2\},V \boxtimes V^{*},\delta_{V},\on{ev}_{V},(\gamma,1)}$ pour des valeurs arbitraires de $\deg(\gamma)$. 
    Ils s'expriment comme des combinaisons des $T(h_{W,v})$ pour $W$ irréductible et ils n'apportent donc aucune information supplémentaire.  
      \end{rem}

     \begin{lem}\label{S-non-ram-prelim2}
      Soit  $v\in |X|\sm N$. Soit $I$ un ensemble fini et $\iota\in I$ un élément. 
      On fixe  un plongement $\ov F\subset \ov {F_{v}}$. Soit  $d\in \N$,  et 
    $\gamma\in \on{Gal}(\ov {F_{v}}/F_{v})\subset  \on{Gal}(\ov F/F)$ tel que $\deg(\gamma)=d$.  
    On définit $(\gamma_{i})_{i\in I}\in\on{Gal}(\ov F/F)^{I}$ en posant $\gamma_{\iota}=\gamma$ et $\gamma_{i}=1$ pour $i\neq \iota$. 
    On fixe un  point géométrique  $\ov v$ au-dessus de  $v$ et une  flèche de spécialisation 
    $\on{\mf{sp}}_{v}:\ov \eta\to \ov v$ associés au choix du plongement $\ov F\subset \ov {F_{v}}$. On note encore   $\on{\mf{sp}}_{v}$ la  flèche de spécialisation de  
    $\Delta(\ov \eta)$ vers   $\Delta(\ov v)$ égale à l'image par $\Delta$ de cette dernière. 
      Alors on a la commutativité du  diagramme  
       $$ \xymatrixcolsep{1pc} \xymatrix{  \Big( \varinjlim _{\mu} \restr{\mc H _{ N, I, W}^{0,\leq\mu,E}}{\Delta(\ov v)}\Big)^{\mr{Hf}} \ar[r]^-{\mf{sp}_{v}^{*}} \ar[d]^{F_{\{\iota\}}^{\deg(v)d}} & 
       \Big( \varinjlim _{\mu} \restr{\mc H _{ N, I, W}^{0,\leq\mu,E}}{\Delta(\ov\eta)}\Big)^{\mr{Hf}}\ar@{=}[r] & H_{I,W} \ar[d]^{(\gamma_{i})_{i\in I}}
       \\  \Big( \varinjlim _{\mu} \restr{\mc H _{ N, I, W}^{0,\leq\mu,E}}{\Delta(\ov v)}\Big)^{\mr{Hf}} \ar[r]^-{\mf{sp}_{v}^{*}} & 
        \Big( \varinjlim _{\mu} \restr{\mc H _{ N, I, W}^{0,\leq\mu,E}}{\Delta(\ov\eta)}\Big)^{\mr{Hf}} \ar@{=}[r] & 
       H_{I,W}
       } $$
                \end{lem}    

            \dem 
           On va utiliser les deux faits généraux suivants. Soit $\Omega$ un ouvert dense d'un schéma $Y$ de type fini sur $\Fq$. On note $i:\Omega\to Y$ l'inclusion et on désigne par $i_{*}$   l'image directe non dérivée. 
    \begin{itemize}
    \item [] 1) Si $\mathcal L$ est un système local $\ell$-adique sur $Y$, 
    le morphisme d'adjonction $\mc L\to i_{*}i^{*}(\mc L)=i_{*}(\restr{\mc L}{\Omega})$ est un isomorphisme de faisceaux sur $Y$.  
      \item []
  2)  Si   $ \mathcal L_1 \subset \mathcal L_2$ sont deux systèmes locaux $\ell$-adiques sur  $\Omega$ et  si $\mf{sp}:\ov x\to \ov y$ est 
  une flèche de spécialisation dans $Y$ d'un point géométrique $\ov x$ de $\Omega$ 
 vers  un point géométrique $\ov y$ de $Y$,  alors le morphisme 
 $(i_*(\mathcal L_1))_{\ov y}\to (i_*(\mathcal L_2))_{\ov y}$ est injectif et son image est formée exactement des éléments $a\in (i_*(\mathcal L_2))_{\ov y}$ tels que 
 $\mf{sp}^{*}(a)\in (\mathcal  L_2)_{\ov x}$  appartienne  à $(\mathcal  L_1)_{\ov x}$. 
     \end{itemize}
    Ces deux assertions se montrent pour des faisceaux de torsion puis, en passant à la limite, pour des faisceaux $\ell$-adiques. 
            
            On commence maintenant la démonstration. 
          On rappelle que l'on peut trouver  une réunion croissante 
 (indexée par $\lambda\in \N$) 
 de sous-$\mc O_{E}$-faisceaux constructibles 
 $\mf F_{\lambda}\subset \varinjlim _{\mu}\restr{\mc H _{N, I, W}^{0,\leq\mu,E}}{\eta^{I}}$ stables par les morphismes de Frobenius partiels
  et  une suite décroissante d'ouverts  denses $U_{\lambda}\subset X\sm N$ tels que 
 $\mf F_{\lambda}$ se prolonge en un faisceau lisse 
 (encore noté $\mf F_{\lambda}$) sur $(U_{\lambda})^{I}$, 
  de sorte que \eqref{egalite-Flambda-H} et 
  \eqref{mor-HIW-Flambda} soient vrais. 
              
       Soit $a\in    \Big( \varinjlim _{\mu} \restr{\mc H _{ N, I, W}^{0,\leq\mu,E}}{\Delta(\ov v)}\Big)^{\mr{Hf}} $. On fixe $\lambda$ tel que, via l'isomorphisme 
       \eqref{mor-HIW-Flambda},  
       $\restr{\mc F_{\lambda}}{\Delta(\ov\eta)} $ contienne 
     $  \mf{sp}^{*}_{v}(a)$. On fixe $\mu$  tel que 
      $\mf F_{\lambda}$   soit un sous-faisceau lisse
de $\mc H _{N, I, W}^{0,\leq\mu,E}$ sur un ouvert $\Omega$ où celui-ci est lisse. 
Quitte  à restreindre $\Omega$ on suppose qu'il est inclus dans  
 $(U_{\lambda})^I$.
     
      On note $j^{I}:U^I\to X^I$  et 
      $j_{\Omega}:\Omega\to X^I$  les inclusions ouvertes. On fixe une flèche de spécialisation 
      $\mf{sp}$ de $\ov{\eta^{I}}$ vers $\Delta(\ov \eta)$ et on note  $  {\mf{sp}}_{v}'$ 
      la    flèche de spécialisation
      de $\ov{\eta^{I}}$ vers $\Delta(\ov v)$
      égale à la composée  de $    \mf{sp}$  et de $  {\mf{sp}}_{v}$. 
On a le diagramme commutatif  
$$ \xymatrixcolsep{3pc} \xymatrix{  
  \restr{(j^{I})_{*}(\mc F_{\lambda})}{\Delta(\ov v)}
  \ar@{=}[d] \ar[dr]^-{ ({\mf{sp}}_{v}')^{*}} & &
       \\ \restr{(j_{\Omega})_{*}(\restr{\mc F_{\lambda}}{\Omega})}{\Delta(\ov v)} \ar@{^{(}->}[r]^-{  ({\mf{sp}}_{v}')^{*}} \ar@{^{(}->}[d] & 
      \restr{\mc F_{\lambda}}{\ov {\eta^{I}}} \ar@{^{(}->}[d]
        \\  \restr{(j_{\Omega})_{*}(\restr{\mc H _{N, I, W}^{0,\leq\mu,E}}{\Omega})}{\Delta(\ov v)} \ar@{^{(}->}[r]^-{ ({\mf{sp}}_{v}')^{*}}& 
      \restr{\mc H _{N, I, W}^{0,\leq\mu,E}}{\ov {\eta^{I}}} 
            } $$
La flèche verticale en haut à gauche est bien une égalité : en effet  en appliquant  1)   à $\mathcal L=\mf F_{\lambda}$ et à l'inclusion $i:\Omega\to 
(U_{\lambda})^I$  on obtient l'égalité (de faisceaux sur $(U_{\lambda})^I$)
$\mf F_{\lambda}=i_{*}(\restr{\mc F_{\lambda}}{\Omega})$ et on remarque que 
 $j_{\Omega}=j^{I}\circ i$. L'injectivité des $ ({\mf{sp}}_{v}')^{*}$ découle du fait que 
$\mc F_{\lambda}$ et $\mc H _{N, I, W}^{0,\leq\mu,E}$ sont lisses sur $\Omega$. On a un morphisme d'adjonction \begin{gather}\label{morph-adjonc}\mc H _{N, I, W}^{0,\leq\mu,E}\to (j_{\Omega})_{*}(j_{\Omega})^{*}(\mc H _{N, I, W}^{0,\leq\mu,E})=(j_{\Omega})_{*}(\restr{\mc H _{N, I, W}^{0,\leq\mu,E}}{\Omega}). \end{gather}
On note encore $a\in  \restr{(j_{\Omega})_{*}(\restr{\mc H _{N, I, W}^{0,\leq\mu,E}}{\Omega})}{\Delta(\ov v)} $ l'image de $a$ par  la fibre en $\Delta(\ov v)$ du morphisme \eqref{morph-adjonc}. 
Par hypothèse $({\mf{sp}}_{v}')^{*}(a)\in  \restr{\mc H _{N, I, W}^{0,\leq\mu,E}}{\ov {\eta^{I}}} $ est l'image d'un élément  $b\in  \restr{\mc F_{\lambda}}{\ov {\eta^{I}}}$.      
Comme on l'a rappelé en détails dans la preuve du \lemref{lem-HIW-action-indep-sp},  l'action de $(\gamma_{i})_{i\in I}$ sur $ \mf{sp}^{*}_{v}(a)\in H_{I,W}$ 
    (ou sur   $ (\mf{sp}_{v}')^{*}(a)\in  \Big( \varinjlim _{\mu} \restr{\mc H _{ N, I, W}^{0,\leq\mu,E}}{\ov{\eta_{I}}}\Big)^{\mr{Hf}}$   qui est  son image  par 
    $\mf{sp}^{*}$) 
    est donnée par l'application du lemme de Drinfeld à 
          $\mc F_{\lambda}$. Plus précisément c'est l'action de $(\gamma_{i})_{i\in I}$ sur  $b\in \restr{\mc F_{\lambda}}{\ov {\eta^{I}}}$. 
          En appliquant  2) à $Y=X^{I}$, 
          $\ov x=\ov {\eta^{I}}$, $\ov y=\Delta(\ov v)$, 
          $ \mathcal L_1=\restr{\mf F_{\lambda}}{\Omega}$ et  
$\mathcal L_2=\restr{\mc H _{N, I, W}^{0,\leq\mu,E}}{\Omega}$ on voit qu'il existe un unique $c\in  \restr{(j_{\Omega})_{*}(\restr{\mc F_{\lambda}}{\Omega})}{\Delta(\ov v)} $
dont l'image par $ ({\mf{sp}}_{v}')^{*}$ est $b$ et dont l'image dans  
$\restr{(j_{\Omega})_{*}(\mc H _{N, I, W}^{0,\leq\mu,E})}{\Delta(\ov v)}$ est $a$. 
De plus  $F_{\{\iota\}}^{\deg(v)d} (a) $ est l'image de     $F_{\{\iota\}}^{\deg(v)d} (c) $    (où $F_{\{\iota\}}$ est ici le morphisme de Frobenius partiel pour 
$\mc F_{\lambda}$). 
          En appliquant le  \lemref{S-non-ram-prelim}
       à $\mc E=\mc F_{\lambda}$  et à $c$ on obtient finalement 
       la commutativité du diagramme. 
                  \cqfd
          
  \begin{notation} \label{alg-B}
    On note  $\mc B$ la  sous-$E$-algèbre de 
   $$\on{End}_{C_{c}(K_{N}\backslash G(\mb A)/K_{N},E)}(C_{c}^{\rm{cusp}}(\Bun_{G,N}(\Fq)/\Xi,E))$$  engendrée par tous les  opérateurs d'excursion  
      $S_{I,W,x,\xi,(\gamma_i)_{i\in I}}$ (elle dépend bien sûr de $N$,  $\Xi$ et $E$). 
          D'après   \eqref{SIW-p2}, $\mc B$ est  commutative. Elle est évidemment 
             de dimension finie. D'après le  \lemref{S-non-ram} elle contient tous les  opérateurs de Hecke en les  places de $X\sm N$. 
 \end{notation}

   Les fonctions \begin{gather}\label{art-def-f}
   f: (g_{i})_{i\in I}\mapsto \s{\xi, (g_{i})_{i\in I}\cdot x}\end{gather} que l'on obtient en faisant varier  $W$, 
 $x$, et  $\xi $ sont exactement les  fonctions régulières sur le quotient grossier 
  de   $(\wh G_E)^{I}$ par les actions par translation à gauche et à droite de   $\wh G_E$ diagonal, que l'on notera $\wh G_E\backslash (\wh G_E)^{I}/\wh G_E$. 

Le lemme suivant  reprend le \lemref{f-W,x,xi} de l'introduction. 

 \begin{lem}\label{f-W,x,xi-texte}
L'opérateur $S_{I,W,x,\xi,(\gamma_{i})_{i\in I}}$ dépend seulement de   $I$, $f$, et $(\gamma_{i})_{i\in I}$, où $f$ est donnée par \eqref{art-def-f}. 
\end{lem}
\dem Soit $W,x,\xi$ comme précédemment et soit $f\in \mc O(\wh G_E\backslash (\wh G_E)^{I}/\wh G_E)$ donnée  par \eqref{art-def-f}. 
 On note $W_{f}$ le sous-$E$-espace vectoriel de dimension finie   de 
 $\mc O((\wh G_E)^{I}/\wh G_E)$ engendré les  translatées à gauche de $f$ par $(\wh G_E)^{I}$. On pose  $x_{f}=f\in W_{f}$ et on note  $\xi_{f}$ la forme linéaire  sur $W_{f}$ donnée par l'évaluation en $1\in (\wh G_E)^{I}/\wh G_E$. 
Alors  $W_{f}$ est un  sous-quotient de $W$: si $W_{x}$ est la 
sous-$(\wh G_E)^{I}$-représentation   de $W$ engendrée par $x$, 
 $W_{f}$ est le  quotient de $W_{x}$ par la plus grande sous-$(\wh G_E)^{I}$-représentation $E$-linéaire sur laquelle $\xi$ s'annule. 
 On a alors les   diagrammes $$W \overset{\alpha}{\hookleftarrow} W_{x} \overset{\beta}{\twoheadrightarrow} W_{f}, 
 \ \ x \overset{\alpha}{\longleftarrow\!\shortmid} x \overset{\beta}{\shortmid\!\longrightarrow} x_{f}, 
 \ \ \xi \overset{{}^{t}\alpha}{\shortmid\!\longrightarrow}  \restr{\xi}{W_{x}} \overset{{}^{t}\beta}{\longleftarrow\!\shortmid} \xi_{f}$$
 de  $(\wh G_E)^{I}$-représentations, de 
  vecteurs $\wh G_E$-invariants  et de  formes linéaires  $\wh G_E$-invariantes. 
  En appliquant  \eqref{SIW-p0} 
  à $u=\alpha$ et $u=\beta$, on obtient 
      \begin{gather}\nonumber 
  S_{I,W,x,\xi,(\gamma_{i})_{i\in I}}=
  S_{I,W_{x},x,\xi |_{W_{x}},(\gamma_{i})_{i\in I}}
  =S_{I,W_{f},x_{f},\xi_{f},(\gamma_{i})_{i\in I}}.\end{gather}
Cela montre que  $S_{I,W,x,\xi,(\gamma_{i})_{i\in I}}$ dépend seulement de  $I$, $f$, et $(\gamma_{i})_{i\in I}$. 
\cqfd

 \begin{notation} \label{S-If}
 Pour toute fonction  $f\in  \mc O(\wh G_E\backslash (\wh G_E)^{I}/\wh G_E)$ on pose 
 \begin{gather} S_{I,f,(\gamma_{i})_{i\in I}}=S_{I,W,x,\xi,(\gamma_{i})_{i\in I}} \in \mc B\end{gather}
  où $W,x,\xi$ sont tels que  $f$ satisfasse  
 \eqref{art-def-f}. 
  \end{notation}

 Grâce aux  notations \ref{alg-B} et \ref{S-If}, la  proposition  suivante (qui reprend la \propref{prop-SIf-i-ii-iii} de l'introduction) 
 reformule de fa\c con plus synthétique toutes les  propriétés  des opérateurs d'excursion, et servira de référence dans le chapitre suivant. 
 
  \begin{prop} \label{art-prop-SIf-i-ii-iii} Les opérateurs d'excursion   $S_{I,f,(\gamma_{i})_{i\in I}}$ vérifient les  propriétés suivantes:  
  \begin{itemize}
  \item [] (i) pour tout  $I$ et 
 $(\gamma_{i})_{i\in I}\in  \on{Gal}(\ov F/F)^{I}$, 
  $$f\mapsto 
  S_{I,f,(\gamma_{i})_{i\in I}}$$ est un  morphisme 
  d'algèbres commutatives  $\mc O(\wh G_E\backslash (\wh G_E)^{I}/\wh G_E)\to 
  \mc B$, 
  \item [] (ii) pour toute  application 
  $\zeta:I\to J$,  tout   $f\in \mc O(\wh G_E\backslash (\wh G_E)^{I}/\wh G_E)$ et  
  tout $(\gamma_{j})_{j\in J}\in  \on{Gal}(\ov F/F)^{J}$, on a  
  $$S_{J,f^{\zeta},(\gamma_{j})_{j\in J}}=S_{I,f,(\gamma_{\zeta(i)})_{i\in I}}$$
   où $f^{\zeta}\in \mc O(\wh G_E\backslash (\wh G_E)^{J}/\wh G_E)$ est définie par    $$f^{\zeta}((g_{j})_{j\in J})=f((g_{\zeta(i)})_{i\in I}),$$
   \item [] (iii) 
      pour tout   $f\in \mc O(\wh G_E\backslash (\wh G_E)^{I}/\wh G_E)$
  et  $(\gamma_{i})_{i\in I},(\gamma'_{i})_{i\in I},(\gamma''_{i})_{i\in I}$ dans $  \on{Gal}(\ov F/F)^{I}$ on a     $$S_{I\cup I\cup I,\wt f,(\gamma_{i})_{i\in I}\times (\gamma'_{i})_{i\in I}\times (\gamma''_{i})_{i\in I}}=
  S_{I,f,(\gamma_{i}(\gamma'_{i})^{-1}\gamma''_{i})_{i\in I}}$$
   où  
   $\wt f\in \mc O(\wh G_E\backslash (\wh G_E)^{I\cup I\cup I}/\wh G_E)$ est définie par  
   $$\wt f((g_{i})_{i\in I}\times (g'_{i})_{i\in I}\times (g''_{i})_{i\in I})=f((g_{i}(g'_{i})^{-1}g''_{i})_{i\in I}).$$

        \item [] (iv) pour tout  $I$ et $f\in \mc O(\wh G_E\backslash (\wh G_E)^{I}/\wh G_E)$, il existe un ouvert dense $U\subset X$ tel que 
     $
  S_{I,f,(\gamma_{i})_{i\in I}}$ dépend seulement de  l'image de 
 $(\gamma_{i})_{i\in I}$ dans   $\pi_{1}(U, \ov\eta)^{I}$, et 
 $(\gamma_{i})_{i\in I}\mapsto S_{I,f,(\gamma_{i})_{i\in I}}$ est continue du  groupe profini $\pi_{1}(U, \ov\eta)^{I}$ vers la $E$-algèbre de  dimension finie   
   $\mc B$ munie de la topologie   $E$-adique. 
     \item [] (v) 
           Soit  $v\in |X|\sm N$. Soit $V$ une   représentation $E$-linéaire  irréductible de $\wh G$. Soit $f\in \mc O(\wh G_E\backslash (\wh G_E)^{\{1,2\}}/\wh G_E)$ la  fonction 
           $(g,g')\mapsto \chi_{V}(gg'^{-1})$, où $\chi_{V}$ désigne  le  caractère de $V$. 
     On fixe  un plongement $\ov F\subset \ov {F_{v}}$, qui induit donc un plongement
      $\on{Gal}(\ov {F_{v}}/F_{v})\hookrightarrow \on{Gal}(\ov F/F)$. 
      Soit  $d\in \N$. Soit 
    $\gamma\in \on{Gal}(\ov {F_{v}}/F_{v})$ tel que $\deg(\gamma)=d$.  
    Alors  $S_{\{1,2\},f,(\gamma,1)}$   dépend seulement de  $d$, et est égal à 
    $T(h_{V,v})$ si  $d=1$.    
           \end{itemize}
   \end{prop}
    
      \noindent {\bf Démonstration de la \propref{art-prop-SIf-i-ii-iii}.}    
On déduit   (ii) de \eqref{SIW-p1}. Pour montrer   (i) on remarque d'abord que la linéarité de $S_{I,f,(\gamma_{i})_{i\in I}}$ en $f$ se déduit facilement de la linéarité de $S_{I,W,x,\xi,(\gamma_{i})_{i\in I}}$ en $x$ et $\xi$, et pour la multiplicativité 
 on utilise  
\eqref{SIW-p2} et  on applique  (ii) à l'application évidente $\zeta :I\cup I\to I$
en remarquant que 
$$\s{ \xi_{1} \boxtimes \xi_{2}, \big((g_{i})_{i\in I} 
  \boxtimes  (g_{i})_{i\in I}  \big)\cdot ( x_{1} \boxtimes x_{2})}
  =\s{ \xi_{1}  ,  (g_{i})_{i\in I} 
  \cdot   x_{1}  } \s{ \xi_{2}  ,  (g_{i})_{i\in I} 
   \cdot   x_{2}  } 
.$$
 L'assertion (iii) découle de \eqref{SIW-p3}, en remarquant que  pour $(g_{i})_{i\in I}$, $(g'_{i})_{i\in I}$  et $(g''_{i})_{i\in I}$ dans  $(\wh G_E)^{I}$ on a  
   $$\s{ \xi \boxtimes \on{ev}_{W}, \big((g_{i})_{i\in I} 
  \boxtimes  (g'_{i})_{i\in I} \boxtimes  (g''_{i})_{i\in I} \big)\cdot ( \delta_{W} \boxtimes x)}
  =\s{\xi, (g_{i}(g'_{i})^{-1}g''_{i})_{i\in I} \cdot x}.$$
     L'assertion  (iv) est incluse dans la  définition-proposition \ref{antecedent-def-prop}. Enfin l'assertion (v) est une reformulation du  \lemref{S-non-ram}. 
 \cqfd

\begin{rem} \label{rem-gamma-en-plus}
Pour tout $\gamma\in  \on{Gal}(\ov F/F)$, on a 
  $S_{I,f,(\gamma_{i})_{i\in I}}=S_{I,f,(\gamma_{i}\gamma)_{i\in I}}$.   On  le vérifie très facilement à l'aide de la définition des opérateurs d'excursion, ou bien  on applique (iii) avec $\gamma'_{i}=1$ et $\gamma''_{i}=\gamma$, puis on applique (ii) aux applications 
  $I\cup I \cup I\to I\cup \{0\} \cup \{0\}$ et $I\to I\cup \{0\} \cup \{0\}$. 
  De même on a  $S_{I,f,(\gamma_{i})_{i\in I}}=S_{I,f,(\gamma\gamma_{i})_{i\in I}}$. 
\end{rem}

La proposition suivante améliore (iv) de la proposition précédente. 
Cet énoncé plus fort (qui donnera une preuve plus limpide que dans les versions antérieures de la non-ramification sur $X\sm N$ des paramètres de Langlands que nous construirons) a été trouvé par    B\"ockle,    Harris,  Khare et  Thorne  et sert   dans  leur preprint \cite{boeckle-harris...}. 

           \begin{prop} \label{prop-harris} (trouvée par B\"ockle,    Harris,  Khare et  Thorne) 
          Pour tout  $I$ et $f\in \mc O(\wh G_E\backslash (\wh G_E)^{I}/\wh G_E)$,    $S_{I,f,(\gamma_{i})_{i\in I}}$ dépend seulement de  l'image de 
 $(\gamma_{i})_{i\in I}$ dans   $\pi_{1}(X\sm N, \ov\eta)^{I}$, et 
 $(\gamma_{i})_{i\in I}\mapsto S_{I,f,(\gamma_{i})_{i\in I}}$ est continue du  groupe profini $\pi_{1}(X\sm N, \ov\eta)^{I}$ vers la $E$-algèbre de  dimension finie   
   $\mc B$ munie de la topologie   $E$-adique. 
          \end{prop}
    \dem 
    La démonstration reprend d'abord celle du \lemref{lem-sigma-nonram}  de l'introduction. 
       Soit  $v$ une  place de $X\sm N$. On fixe  un plongement $\ov F\subset \ov {F_{v}}$, d'où une  inclusion 
            $ \on{Gal}(\ov {F_{v}}/F_{v})\subset  \on{Gal}(\ov F/F)$. Soit  $I_{v}=\on{Ker}( \on{Gal}(\ov {F_{v}}/F_{v})\to \on{Gal}(\ov {k(v)}/k(v)))$ le groupe d'inertie en   $v$. Alors 
pour tous $I, W,x,\xi$, l'image de  
la composée \begin{gather}\label{compo-crea-inv}  H_{\{0\},\mbf  1}\xrightarrow{\mc H(x)}
 H_{\{0\},W^{\zeta_{I}}}\isor{\chi_{\zeta_{I}}^{-1}} 
  H_{I,W}\end{gather}   est formée d'éléments invariants par 
    $(I_{v})^{I}$. En effet 
      les opérateurs 
    de création sont des morphismes de faisceaux  sur $\Delta(X\sm N)$ tout entier (et en particulier en $\Delta(v)$).  Donc on a un triangle commutatif 
    (que l'on avait déjà utilisé dans le \lemref{S-non-ram})
     $$ \xymatrixcolsep{1pc} \xymatrix{  C_{c}^{\mr{cusp}}(\Bun_{G,N}(\Fq)/\Xi,E)\ar[d]_-{\restr{\mc C_{ \delta_{V}}^{\sharp }}{\ov v}} \ar[dr]^{ \mc C_{  \delta_{V}}^{\sharp }} & &
       \\ \Big( \varinjlim _{\mu} \restr{\mc H _{ N, I, W}^{0,\leq\mu,E}}{\Delta(\ov v)}\Big)^{\mr{Hf}} \ar[r]^-{\mf{sp}_{v}^{*}}  & 
       \Big( \varinjlim _{\mu} \restr{\mc H _{ N, I, W}^{0,\leq\mu,E}}{\Delta(\ov\eta)}\Big)^{\mr{Hf}}\ar@{=}[r] & H_{I,W}  }
      $$
      où $\mf{sp}_{v}$ est la flèche de spécialisation de $\ov\eta$ vers $\ov v=\on{Spec}(\ov{k(v)})$ associée à l'inclusion $\ov F\subset \ov{F_{v}}$. 
   On applique alors le \lemref{S-non-ram-prelim2} avec $d=0$ et on en déduit  
    l'invariance par 
    $(I_{v})^{I}$ de l'image de \eqref{compo-crea-inv}. 
    
      Donc pour $(\gamma_{i})_{i\in I}\in \on{Gal}(\ov F/F)^{I}$ et 
    $(\delta_{i})_{i\in I}\in (I_{v})^{I}$ on a 
    \begin{gather}\label{rel-S-gammai-deltai-texte}S_{I,W,x,\xi,(\gamma_{i})_{i\in I}}=S_{I,W,x,\xi,(\gamma_{i}\delta_{i})_{i\in I}}.\end{gather} 
    Cela est vrai pour tout plongement $\ov F\subset \ov {F_{v}}$ 
    (en fait,    grâce à la \remref{rem-gamma-en-plus},  \eqref{rel-S-gammai-deltai-texte}  pour 
    un plongement implique \eqref{rel-S-gammai-deltai-texte} pour tous les plongements). Or pour tout ouvert $U\subset X\sm N$, $\pi_{1}(X\sm N, \ov\eta)$ est le quotient topologique    de 
    $\pi_{1}(U, \ov\eta)$ par le sous-groupe fermé engendré  
        par les $I_{v}$ pour $v\in (X\sm N)\sm U$ et leurs conjugués. \cqfd
    
             \section {Décomposition suivant les paramètres de Langlands}
        \label{para-dec-param-Langlands}

                 L'idée  se résume ainsi :  grâce à \eqref{isom-chi-singleton-empty} et au   c) de la \propref{prop-a-b-c}, on a 
    $$H_{\{0\},\mbf 1}=C_{c}^{\mr{cusp}}(G(F)\backslash G(\mb A)/K_{N}\Xi,E).$$
   Pour obtenir la décomposition \eqref{intro1-dec-canonique}  il est donc  équivalent de construire (quitte à augmenter $E$) 
   une  décomposition canonique 
     \begin{gather} \nonumber 
  H_{\{0\},\mbf 1}=\bigoplus_{\sigma}
 \mf H_{\sigma}.\end{gather}

     \begin{notation}    Pour tout entier   $n\in \N^{*}$ on note   $\mc O((\wh G_E)^{n}\modmod \wh G_E)$ la  $E$-algèbre des fonctions régulières  sur $(\wh G_E)^{n}$ qui sont invariantes par  conjugaison diagonale, c'est-à-dire   telles  que 
      $f(hg_{1}h^{-1},...,hg_{n}h^{-1})=
        f(g_{1},...,g_{n})$ pour tout  $h,g_{1},...,g_{n}\in \wh G_E$. Hilbert a montré  dans  \cite{hilbert} que cette  $E$-algèbre est de type fini (voir \cite{git} I.2). 
   Autrement dit  $(\wh G_E)^{n}\modmod \wh G_E=
    \mr{Spec}\big(\mc O\big((\wh G_E)^{n}\modmod \wh G_E\big)\big)$ est le  quotient grossier de $(\wh G_E)^{n}$ par la conjugaison diagonale par  $\wh G_E$. 
                 \end{notation}
                 
 \begin{rem}   En caractéristique  $0$ (ce qui est notre cas) ces quotients ont été étudiés notamment  dans \cite{procesi, richardson, vinberg}. 
 L'idée de considérer ces  quotients pour étudier les espaces de modules de représentations de groupes  (par exemple de groupes discrets  de type fini) à valeurs dans des  groupes réductifs complexes est ancienne, voir \cite{lub-magid} et les réferences incluses dans la dernière partie, en particulier   \cite{chandler} I 1.6E qui contient la  réference à Poincaré citée dans l'introduction. 
\end{rem}

On a  un isomorphisme 
\begin{gather}\label{isom-I-1..n}\beta : 
(\wh G_E)^{n}\modmod \wh G_E\isom \wh G_E\backslash (\wh G_E)^{\{0,...,n\}}/\wh G_E, \ (g_{1},...,g_{n})\mapsto (1,g_{1},...,g_{n}) \end{gather}
  dont  l'inverse est donné par 
         $(g_{0},...,g_{n})\mapsto (g_{0}^{-1}g_{1},...,g_{0}^{-1}g_{n})
         $. 
        On va reformuler les  propriétés des opérateurs d'excursion  récapitulées dans la  \propref{art-prop-SIf-i-ii-iii} à l'aide de ces  quotients $(\wh G_E)^{n}\modmod \wh G_E$ et de nouvelles notations. 
        Le passage de l'indexation des quotients $\wh G_E\backslash (\wh G_E)^{I}/\wh G_E$   par des ensembles finis abstraits à l'indexation des  quotients $(\wh G_E)^{n}\modmod \wh G_E$  par des entiers est motivé par les deux raisons suivantes: cela évite toute confusion entre les deux sortes de quotients, et   la nouvelle fa\c con d'indexer sera plus commode pour la  construction des paramètres de Langlands. 

L'identification \eqref{isom-I-1..n} induit l'isomorphisme d'algèbres 
\begin{gather}\label{isom-I-1..n-alg}
\mc O((\wh G_E)^{n}\modmod \wh G_E) \isom \mc O(\wh G_E\backslash (\wh G_E)^{\{0,...,n\}}/\wh G_E), \ f\mapsto f\circ \beta^{-1}. \end{gather}

 Pour tout  groupe  profini $\Gamma$, on note  
   $C(\Gamma,\mc B)$ l'algèbre des fonctions continues de   $\Gamma$  vers   $\mc B$ muni de la topologie $E$-adique.  
            
       \begin{def-prop}\label{cor-annul-Lambda-n} 
           Pour tout $n\in \N^{*}$  on définit 
                     $$ \Theta_{n}:    \mc O((\wh G_E)^{n}\modmod \wh G_E)\to C(\on{Gal}(\ov F/F)^{n},\mc B)$$   par 
        \begin{gather}\label{cond-cor-annul-Lambda-n}
     \Theta_{n}(f): (\gamma_{1},..., \gamma_{n})\mapsto 
     S_{\{0,...,n\},f\circ \beta^{-1} ,(1,\gamma_{1},...,\gamma_{n})}
     \in 
   \mc B. 
     \end{gather}
     
       Ces applications ont les  propriétés suivantes: 
        \begin{itemize}
        \item  [] a) pour tout  $n$, $\Theta_{n}$ est un  morphisme d'algèbres, 
                  \item  [] b) pour tout   $n$,  $\Theta_{n}$ prend ses valeurs dans  
   $C(\pi_{1}(X\sm N, \ov\eta)^{n},\mc B)$, et son image est formée de fonctions invariantes par conjugaison diagonale par $\pi_{1}(X\sm N, \ov\eta)$, 
                  \item [] c) la suite  $(\Theta_{n})_{n\in \N^{*}}$ est fonctorielle par rapport à  toutes les applications entre les ensembles  $\{1,...,n\}$, c'est-à-dire  que  pour $m,n\in \N^{*}$, 
               $\zeta: \{1,...,m\}\to \{1,...,n\}$ arbitraire,   
               $f\in \mc O((\wh G_E)^{m}\modmod \wh G_E)$ et 
               $(\gamma_{1},...,\gamma_{n })\in \on{Gal}(\ov F/F)^{n }$, 
               on a  
             $$\Theta_{n}( f^{\zeta}) ((\gamma_{j})_{j\in \{1,...,n\}})=
             \Theta_{m}(f)((\gamma_{\zeta(i)})_{i\in \{1,...,m\}})$$ 
            où  $f^{\zeta}\in \mc O((\wh G_E)^{n}\modmod \wh G_E)$  est définie par 
               $$f^{\zeta}((g_{j})_{j\in \{1,...,n\}})=f((g_{\zeta(i)})_{i\in \{1,...,m\}}), $$
 \item []    d) pour $n\geq 1$, 
  $f\in \mc O((\wh G_E)^{n}\modmod \wh G_E)$  
    et   $(\gamma_{1},...,\gamma_{n+1})\in \on{Gal}(\ov F/F)^{n+1}$ on a  
   $$\Theta_{n+1}( \wh f)(\gamma_{1},...,\gamma_{n+1})=
   \Theta_{n}( f)(\gamma_{1},...,\gamma_{n}\gamma_{n+1}) $$
 où $\wh f\in  \mc O((\wh G_E)^{n+1}\modmod \wh G_E)$ est définie par 
   $$\wh f(g_{1},...,g_{n+1})=f(g_{1},...,g_{n}g_{n+1}),$$
       \item [] e)    pour toute représentation $E$-linéaire  irréductible  $V$ de $\wh G$, 
       de caractère $\chi_{V}\in \mc O(\wh G_E \modmod \wh G_E)$, pour toute     place $v\in |X|\sm N$, et pour tout élément de Frobenius $\Frob_{v}\in \pi_{1}(X\sm N, \ov\eta)$,
          on a 
   $\Theta_{1}(\chi_{V})(\Frob_{v})=T(h_{V,v})$.   
   \end{itemize}      
      \end{def-prop}
       
         \begin{rem} 
           La propriété  e) ne dépend pas du choix de l'élément $\Frob_{v}$ car, par b),  pour $f\in \mc O(\wh G_E\modmod \wh G_E)$, $\Theta_{1}(f) $ est une fonction centrale sur  $\on{Gal}(\ov F/F)$.           \end{rem}
          
 \begin{rem}\label{rem-S-gamma0=1} La \remref{rem-gamma-en-plus} implique que  tous les  opérateurs d'excursion  apparaissent dans  \eqref{cond-cor-annul-Lambda-n} et donc en considérant la suite $(\Theta_{n})_{n\in \N^{*}}$  on ne perd pas  d'information.  
 \end{rem}

      \noindent {\bf Démonstration de la définition-proposition  \ref{cor-annul-Lambda-n}.}  
     La propriété  a) est simplement  (i) de la  \propref{art-prop-SIf-i-ii-iii} 
   et b) résulte de la \propref{prop-harris} et de la \remref{rem-gamma-en-plus}. Pour démontrer   c) on applique   (ii) de la \propref{art-prop-SIf-i-ii-iii}   à  $I=\{0,...,m\}$, $J=\{0,...,n\}$ et  $\zeta: I\to J$ égal à  l'extension de  $\zeta: \{1,...,m\}\to \{1,...,n\}$ par $\zeta(0)=0$. 
Pour justifier d), on applique la  propriété 
   (iii) de la \propref{art-prop-SIf-i-ii-iii}  à   $$I=\{0,...,n \}, 
   (\gamma_{i})_{i\in I}=(1,\gamma_{1},...,\gamma_{n}), (\gamma'_{i})_{i\in I}=(1)_{i\in I}, (\gamma''_{i})_{i\in I}=(1, ...,1,\gamma_{n+1})$$
    et on utilise   (ii) de la \propref{art-prop-SIf-i-ii-iii} pour supprimer tous les   $1$ sauf le premier dans  
   $(\gamma_{i})_{i\in I}\times (\gamma'_{i})_{i\in I}\times (\gamma''_{i})_{i\in I}$. 
Enfin 
e) est une simple reformulation  de (v) de la \propref{art-prop-SIf-i-ii-iii} dans le cas  où $d=1$ (on ne retient que ce cas car c'est le seul qui sera utilisé). 
         \cqfd

      On rappelle que   $\mc B$ est de dimension finie. 
      On remplace  $E$ par une  extension finie telle que 
    tous les  caractères 
       $\mc B \to \Qlbar$ prennent leurs valeurs  dans  $E$. 
      Pour tout caractère 
       $\nu: \mc B \to E$ on pose  
       $$\Theta_{n}^{\nu}=\nu\circ \Theta_{n}:  \mc O((\wh G_E)^{n}\modmod \wh G_E)\to C(\pi_{1}(X\sm N, \ov\eta)^{n},E).$$ 
      Evidemment la donnée des   $\Theta_{n}^{\nu}$ est équivalente à celle de 
              $$\Theta_{n}^{\rm{red}}:  \mc O((\wh G_E)^{n}\modmod \wh G_E)\to C(\pi_{1}(X\sm N, \ov\eta)^{n},\mc B^{\rm{red}}).$$
       
   \begin{rem}  \label{rem-B-reduite?}  
      On ne sait pas  si $\mc B$ est réduite. La sous-algèbre $\mc B^{1}$ de $\mc B$ engendrée par les  fonctions  $\Theta_{1}(f)$ pour $f\in \mc O(\wh G_E\modmod \wh G_E)$ est en fait  engendrée par les  opérateurs de Hecke en  les places non ramifiées (par la propriété e) et le  théorème  de  Tchebotarev). Par conséquent  $\mc B^{1}$ est définie sur $\ov \Q$, et après changement des scalaires à   $\C$ c'est une  $C^{*}$-algèbre commutative (parce que l'adjoint d'un opérateur de Hecke pour la structure hermitienne usuelle sur 
   $C_{c}^{\mr{cusp}}(G(F)\backslash G(\mb A)/K_{N}\Xi,\C)$  est un opérateur de Hecke).  Donc $\mc B^{1}$ est réduite (en termes plus élémentaires les éléments de $\mc B^{1}$ sont normaux, donc diagonalisables).  Par le même argument, la conjecture  \ref{conj-algebre-B-Qbar}  ci-dessous impliquerait que   $\mc B$ est réduite. 
 \end{rem}   
       
     On appliquera la  proposition suivante à chacun des  $\Theta_{n}^{\nu}$  avec  $$\Gamma= \pi_{1}(X\sm N, \ov\eta), \ \  \text{ et } \ \  
    H=H^{0}=\wh G_{E}.$$ 
 Dans la proposition suivante on ne suppose pas $H$ connexe car on l'utilisera  de nouveau dans le 
     chapitre  \ref{para-non-deploye},  avec 
    $H$   égal au  $L$-groupe d'un groupe réductif non nécessairement déployé.

       \begin{prop}\label{Xi-n}
     Soit $\Gamma$  un groupe profini et  $H$ un groupe réductif sur $E$ non nécessairement  connexe    tel que $H^{0}$ est déployé. 
    On se donne    pour tout   $n\in \N^{*}$    un   morphisme d'algèbres 
          $$\Xi_{n}:   \mc O((H)^{n}\modmod H^{0})\to C(\Gamma^{n},E)$$
 de telle sorte que    \begin{itemize}
              \item [] a)  la suite $(\Xi_{n})_{n\in \N^{*}}$ est fonctorielle relativement aux applications entre les ensembles $\{1,...,n\}$, c'est-à-dire que pour $m,n\in \N^{*}$, 
                $\zeta: \{1,...,m\}\to \{1,...,n\}$ arbitraire,    
                           $f\in \mc O((H)^{m}\modmod H^{0})$ et 
                             $(\gamma_{1},...,\gamma_{n })\in \Gamma^{n }$,  on  a 
             $$\Xi_{n}( f^{\zeta}) ((\gamma_{j})_{j\in \{1,...,n\}})=
             \Xi_{m}(f)((\gamma_{\zeta(i)})_{i\in \{1,...,m\}}),$$
             où $f^{\zeta}\in \mc O((H)^{n}\modmod H^{0})$  est définie par  
             $$f^{\zeta}((g_{j})_{j\in \{1,...,n\}})=f((g_{\zeta(i)})_{i\in \{1,...,m\}}), $$
          \item [] b) 
          pour tout $n\geq 1$, 
  $f\in \mc O((H)^{n}\modmod H^{0})$ et 
 $(\gamma_{1},...,\gamma_{n+1})\in \Gamma^{n+1}$ on a   
   $$\Xi_{n+1}( \wh f)(\gamma_{1},...,\gamma_{n+1})=
   \Xi_{n}( f)(\gamma_{1},...,\gamma_{n}\gamma_{n+1}) $$
    où  $\wh f\in  \mc O((H)^{n+1}\modmod H^{0})$ est définie par  
   $$\wh f(g_{1},...,g_{n+1})=f(g_{1},...,g_{n}g_{n+1}). $$
  \end{itemize}
 
 Alors il existe un  morphisme continu  $\sigma:\Gamma\to H(E')$, où  $E'$  est une  extension finie de $E$, tel que 
l'adhérence de Zariski de son image soit un sous-groupe réductif de     $H_{E'}$ et    tel que 
 pour tout   $n\in \N^{*}$, $f\in \mc O((H)^{n}\modmod H^{0})$ et $(\gamma_{1},...,\gamma_{n})\in \Gamma^{n}$, on ait   
 \begin{gather}\label{équation-Xi-sigma} f(\sigma(\gamma_{1}),...,\sigma(\gamma_{n}))=
 \big(\Xi_{n}(f)\big)(\gamma_{1},...,\gamma_{n}).\end{gather} De plus les   morphismes $\sigma$ ayant  ces  propriétés forment une   unique classe de conjugaison par   $H^{0}(\Qlbar)$.  
         \end{prop}

        \begin{rem} \label{rem-pseudo-car-Taylor}
Dans le cas  où $H=GL_{r}$ la  proposition précédente est bien connue grâce aux résultats de \cite{taylor} sur les pseudo-caractères.  
Plus précisément on note   $\on{St}$ la    représentation standard de $GL_r$. 
Alors il  résulte de   \cite{procesi} que dans la   proposition précédente, $\tau=\Theta_{1}(\chi_{\on{St}})$ détermine  $\Theta_{n}$ pour tout  $n$. 
Comme  la    représentation standard engendre 
(avec l'inverse du déterminant) 
toutes les  représentations (de dimension finie) de $GL_r$,  il suffit de connaître   
$\Theta_{n}(f)$ pour $f$ de la forme $(g_{1},...,g_{n})\mapsto
\on{Tr}(T\cdot (g_{1},...,g_{n}))$ avec  $T\in \on{End}(\on{St}^{\otimes n})^{GL_{r}}$. Or  $\on{End}(\on{St}^{\otimes n})^{GL_{r}}$
 est l'image de l'algèbre du groupe $\mf S_{n}$, et on peut donc supposer que   $T$ est l'action d'une permutation $\sigma$. On calcule 
 $$\on{Tr}(\sigma \cdot (g_{1},...,g_{n}))=\prod _{(i_{1},...,i_{k}) \text{ cycle de } \sigma}\chi_{\on{St}}(g_{i_{k}}\cdots g_{i_{1}})$$   et donc  $$\Theta_{n}(f)(\gamma_{1},...,\gamma_{n})=
 \prod _{(i_{1},...,i_{k}) \text{ cycle de } \sigma}\tau (\gamma_{i_{k}}\cdots \gamma_{i_{1}}). $$
  Le fait que   $\Lambda^{r+1}\on{St}=0$ implique   l'annulation de la  fonction
  $$(g_{1},...,g_{r+1})\mapsto \on{Tr}_{\on{St}^{\otimes (r+1)}}\Big((\sum_{\sigma\in \mf S_{r+1}}s(\sigma) \sigma)(g_{1}\otimes \cdots \otimes g_{r+1})\Big)$$  dans  $\mc O((\wh G_E)^{r+1}\modmod \wh G_E)$.  
 En décomposant   $\sigma$ en un produit de cycles, on obtient que $\tau =\Theta_{1}(\chi_{\on{St}})$, qui est une  fonction centrale sur  $\on{Gal}(\ov F/F)$,   vérifie la relation de {\it pseudo-caractère} suivante: 
\begin{gather}\label{rel-GLN}\sum_{\sigma\in \mf S_{r+1}}s(\sigma) \Big(   
\prod _{(i_{1},...,i_{k}) \text{ cycle de } \sigma}\tau (\gamma_{i_{k}}\cdots \gamma_{i_{1}})
\Big)=0\end{gather} comme  fonction de   $(\gamma_{1},...,\gamma_{r+1})\in \on{Gal}(\ov F/F)^{r+1}$. 
 La relation \eqref{rel-GLN} est due à Frobenius \cite{frobenius-pseudo-car} et elle a été  étudiée et généralisée dans de nombreux travaux, notamment  \cite{wiles, taylor,rouquier,bellaiche-chenevier,chenevier}. En utilisant des résultats de Procesi \cite{procesi,procesi87}, Taylor a montré dans   \cite{taylor} que tout pseudo-caractère à coefficients dans  un corps algébriquement clos de caractéristique  $0$ est le  caractère d'une  représentation semi-simple de dimension $r$ (on renvoie à la   proposition 2.3 de \cite{chenevier-hab} et à la  discussion qui la précède pour un énoncé plus schématique). 
\end{rem}

  \noindent{\bf Démonstration de la \propref{Xi-n}.} 
Pour tout $n$-uplet 
       $(\gamma_{1},...,\gamma_{n})\in \Gamma^{n}$, on note  
       $\xi_{n}(\gamma_{1},...,\gamma_{n})$ le   point de   $(H)^{n}\modmod H^{0}$ sur  $E$   correspondant   au caractère 
        $$\mc O((H)^{n}\modmod H^{0})\to E, \ \  f\mapsto \big( \Xi_{n}(f)\big)(\gamma_{1},...,\gamma_{n}).$$  
     Suivant   \cite{richardson} on dit qu'un   $n$-uplet $(g_{1},...,g_{n})$ de points de $H(\Qlbar)$ est semi-simple si l'adhérence de Zariski $\ov{<g_{1},...,g_{n}>}$   du sous-groupe  $<g_{1},...,g_{n}>$ de $H$ engendré par  $g_{1},...,g_{n}$ est réductive. 
   
\begin{lem}\label{lem-richardson} (Richardson)
Soit $(g_{1},...,g_{n})\in H(\Qlbar)^{n}$. Alors les  assertions suivantes sont équivalentes 
\begin{itemize}
\item [] (i) $(g_{1},...,g_{n})$ est un   $n$-uplet  semi-simple, 
\item [] (ii) la  classe de conjugaison de 
$(g_{1},...,g_{n})$ par $H^{0}_{\Qlbar}$ est fermée  dans  $(H_{\Qlbar})^{n}$. 
\end{itemize} 
De plus les  $H^{0}_{\Qlbar}$-orbites fermées  dans  $(H_{\Qlbar})^{n}$ correspondent aux points sur  $\Qlbar$ du quotient grossier 
$(H)^{n}\modmod H^{0}$, c'est-à-dire  aux  caractères  
$\mc O((H)^{n}\modmod H^{0})\to \Qlbar$. 
\end{lem}
 \noindent{\bf Démonstration.} La dernière   assertion résulte des propriétés générales des quotients grossiers (cf 1.3.2 de \cite{richardson}). Il reste à montrer  l'équivalence de $(i)$ et  $(ii)$. C'est le  théorème 3.6 de \cite{richardson}, 
 à ceci près que l'on quotiente ici  par conjugaison  par $H^{0}$ et non par $H$. 
 Cependant la preuve du  théorème 3.6 de \cite{richardson}  s'adapte sans difficultés:   \begin{itemize}
 \item  dans  $(ii) \Rightarrow (i)$ il n'y a rien à changer puisque la contradiction est amenée par un  cocaractère, qui prend donc  ses valeurs dans $H^{0}$, 
 \item dans  $(i) \Rightarrow (ii)$ on applique le théorème  de  Hilbert-Mumford 
 (théorème 2.1 de \cite{richardson}) à l'action de $H^{0}$ (au lieu  de $H$) sur 
  $(H)^{n}$.  \cqfd
 \end{itemize}

  \noindent{\bf Suite de la démonstration  de  la \propref{Xi-n}.} 
D'après  le  lemme précédent on a une  bijection entre 
\begin{itemize}
\item les  points sur $\Qlbar$ du quotient  grossier 
       $(H)^{n}\modmod H^{0}$
        \item  les classes  de conjugaison  par $H^{0}(\Qlbar)$ de  $n$-uplets semi-simples  de points de $H(\Qlbar)$.
        \end{itemize} 
     On note   $\xi_{n}^{ss}(\gamma_{1},...,\gamma_{n})$ la  classe  de conjugaison par $H^{0}(\Qlbar)$ de  $n$-uplets  semi-simples dans $H(\Qlbar)$ qui est associée à  $\xi_{n}(\gamma_{1},...,\gamma_{n})\in \big((H)^{n}\modmod H^{0}\big)(\Qlbar)$. Autrement dit $\xi_{n}^{ss}(\gamma_{1},...,\gamma_{n})$ est, au choix 
     \begin{itemize}
     \item  l'ensemble des  points sur  $\Qlbar$ de  l'unique orbite fermée au-dessus de $\xi_{n}(\gamma_{1},...,\gamma_{n})$
     \item l'ensemble des    $n$-uplets semi-simples $(g_{1},...,g_{n})$ de $H(\Qlbar)^{n}$ satisfaisant 
      \begin{gather}\label{cond-g1-gn}
   \forall f\in 
     \mc O((H)^{n}\modmod H^{0}),    \ \  f(g_{1},...,g_{n})=
     \Xi_{n}(f)(\gamma_{1},...,\gamma_{n}) . 
     \end{gather} 
     \end{itemize}
        Soit  $(g_{1},...,g_{n})\in \xi_{n}^{ss}(\gamma_{1},...,\gamma_{n})$. 
     On note  $C(g_{1},...,g_{n})$ le centralisateur  de $(g_{1},...,g_{n})$ dans 
     $H^{0}_{\Qlbar}$
     (c'est-à-dire   le centralisateur dans  $H^{0}_{\Qlbar}$ du  sous-groupe  réductif 
     $\ov{<g_{1},...,g_{n}>}$ de $H_{\Qlbar}$). On note  $D(g_{1},...,g_{n})$ le centralisateur    de $C(g_{1},...,g_{n})$ 
    dans  $H_{\Qlbar}$. 
    Comme le centralisateur d'un  sous-groupe  réductif dans un   groupe réductif  est réductif, $C(g_{1},...,g_{n})$ et  $D(g_{1},...,g_{n})$ 
    sont réductifs. 
   On note  $\mf N$ l'ensemble de  $(n,(\gamma_{1},...,\gamma_{n}))$ avec $n\in \N^{*}$ et  $(\gamma_{1},...,\gamma_{n})\in \Gamma^{n}$. 
   Pour  $(n,(\gamma_{1},...,\gamma_{n})) \in \mf N$, on choisit  
   $(g_{1},...,g_{n}) \in \xi_{n}^{ss}(\gamma_{1},...,\gamma_{n})$ (on rappelle que  $(g_{1},...,g_{n})$, et donc  $C(g_{1},...,g_{n})$, sont déterminés de fa\c con unique  à conjugaison près par  $H^{0}_{\Qlbar}$). 
   On note  
           \begin{itemize}
  \item $\mf N^{1}$ le sous-ensemble   de $\mf N$ formé des  $(n,(\gamma_{1},...,\gamma_{n}))$ tels que  $\dim (\ov{<g_{1},...,g_{n}>})$ soit maximal,
  
  \item $\mf N^{2}$ le sous-ensemble  de $\mf N^{1}$ formé des  $(n,(\gamma_{1},...,\gamma_{n}))$ tels que  $\dim (C(g_{1},...,g_{n}))$ soit   minimal,
  
  \item $\mf N^{3}$ le sous-ensemble  de   $\mf N^{2}$ formé des  $(n,(\gamma_{1},...,\gamma_{n}))$  tels que le nombre  de composantes connexes de  $C(g_{1},...,g_{n})$ soit  minimal. 
   \end{itemize}
 Il est clair que  pour $(n,(\gamma_{1},...,\gamma_{n}))\in \mf N^{3}$, 
   $C(g_{1},...,g_{n})$ est minimal (à conjugaison près) parmi tous les sous-groupes  de   $H^{0}_{\Qlbar}$ 
     construits de cette fa\c con lorsque   $(n,(\gamma_{1},...,\gamma_{n}))$ parcourt   $\mf N^{1}$.   
       
       On choisit $(n,(\gamma_{1},...,\gamma_{n}))\in \mf N^{3}$ et on fixe 
        $(g_{1},...,g_{n}) \in \xi_{n}^{ss}(\gamma_{1},...,\gamma_{n})$.
        
     \begin{lem}  \label{lem-gammaCD}
      Pour tout  $\gamma\in \Gamma$, il existe un  unique    $g\in G$
 tel que        $(g_{1},...,g_{n},g)$ appartienne  à   $\xi_{n+1}^{ss}(\gamma_{1},...,\gamma_{n},\gamma)$. De plus 
  $C(g_{1},...,g_{n},g)=C(g_{1},...,g_{n})$ et $g$   appartiennent à  $D(g_{1},...,g_{n})$. 
  \end{lem}
 \noindent{\bf Démonstration. } Soit 
 $(h_{1},...,h_{n},h)\in \xi^{ss}_{n+1}(\gamma_{1},...,\gamma_{n},\gamma)$.  On ne sait pas  {\it a priori} que   
 $(h_{1},...,h_{n})$ est  semi-simple mais, grâce à la  condition a) appliquée à l'inclusion $\zeta: \{1,...,n\}\to \{1,...,n+1\}$,   $(h_{1},...,h_{n})$ est au-dessus de  $\xi_{n}(\gamma_{1},...,\gamma_{n})$. Le théorème 5.2 de 
 \cite{richardson} implique donc  que   $\ov{<g_{1},...,g_{n}>}$ est conjugué à un sous-groupe de Levi  de  $\ov{<h_{1},...,h_{n}>}$ et donc que  
 \begin{gather}\label{ineg-hh-gg}\dim (\ov{<g_{1},...,g_{n}>})\leq \dim (\ov{<h_{1},...,h_{n}>}).\end{gather}
  Par  définition  de $\mf N^{1}$, et comme  $\dim (\ov{<h_{1},...,h_{n}>}) \leq \dim (\ov{<h_{1},...,h_{n},h>})$, l'égalité a lieu dans \eqref{ineg-hh-gg}
   et donc  $(h_{1},...,h_{n})$ est semi-simple et conjugué à  $(g_{1},...,g_{n})$
   (de plus on en déduit que $(n+1,(\gamma_{1},...,\gamma_{n},\gamma))\in \mf N^{1}$). Quitte à conjuguer  $(h_{1},...,h_{n},h)$, on  suppose que   $(h_{1},...,h_{n})=(g_{1},...,g_{n})$. On pose alors   $g=h$, si bien que  $(g_{1},...,g_{n},g)\in \xi_{n+1}^{ss}(\gamma_{1},...,\gamma_{n},\gamma)$. On a évidemment 
 $C(g_{1},...,g_{n},g)\subset C(g_{1},...,g_{n})$,  et comme   $(n,(\gamma_{1},...,\gamma_{n}))\in \mf N^{3}$  et  $(n+1,(\gamma_{1},...,\gamma_{n},\gamma))\in \mf N^{1}$, 
cette  inclusion est une égalité, 
 donc  $g$ appartient à   $D(g_{1},...,g_{n})$  et il est déterminé de manière unique  (en effet il était déterminé  à conjugaison près par $C(g_{1},...,g_{n})$  mais on vient  de montrer qu'il  appartient au centralisateur de $C(g_{1},...,g_{n})$). 
 \cqfd
 
   \noindent{\bf Fin de la preuve  de la \propref{Xi-n}.}  On note  $\sigma:\Gamma\to H(\Qlbar)$ l'application  $\gamma\mapsto g$  que l'on a construite dans le lemme précédent. Il reste à montrer que 
 \begin{itemize}
\item $ \sigma$ prend ses valeurs  dans $H(E')$ où $E'$ est une extension  finie de $E$, 
\item $\sigma$ est un   morphisme de groupes, 
\item $\sigma$ est continue. 
\end{itemize}
 
 Le premier point est clair: si  $E'$ est une  extension finie de $E$ 
 telle  que  $g_{1},...,g_{n}$ appartiennent  à $H(E')$ 
 alors pour tout  $\gamma$,  $g$ appartient à   $H(E')$ (en effet $g\in  D(g_{1},...,g_{n})$ est déterminé de manière unique par l'image de 
 $(g_{1},...,g_{n},g)$ dans $(H)^{n+1}\modmod H^{0}$, qui est définie sur $E\subset E'$). 
 
 Pour montrer le deuxième point soit $\gamma, \gamma'\in \Gamma$. 
 Une variante immédiate  du  \lemref{lem-gammaCD} montre qu'il existe  $g,g'\in H$ uniques tels  que $(g_{1},...,g_{n},g,g')$ appartienne à   $\xi_{n+2}^{ss}(\gamma_{1},...,\gamma_{n},\gamma, \gamma')$. 
 Par les mêmes arguments que  dans la preuve du  \lemref{lem-gammaCD}, les  $(n+1)$-uplets 
 $(g_{1},...,g_{n},g)$, $(g_{1},...,g_{n},g')$ et  $(g_{1},...,g_{n},gg')$ sont semi-simples. La condition a) implique que  les deux premiers appartiennent à   $\xi_{n+1}^{ss}(\gamma_{1},...,\gamma_{n},\gamma )$ et 
 $\xi_{n+1}^{ss}(\gamma_{1},...,\gamma_{n},  \gamma')$. La condition b) implique que le troisième appartient à   $\xi_{n+1}^{ss}(\gamma_{1},...,\gamma_{n},\gamma \gamma')$. 
 Par conséquent   $\sigma(\gamma)=g$,   
 $\sigma(\gamma')=g'$ et   $\sigma(\gamma \gamma')=gg'$. Finalement on a montré que $\sigma(\gamma \gamma')=
  \sigma(\gamma)\sigma( \gamma')$. 

 Pour le troisième  point on commence par rappeler que 
 $C(g_{1},...,g_{n})$ et  $D(g_{1},...,g_{n})$ sont des groupes réductifs définis  sur  $E'$. Comme  $\sigma$ prend ses valeurs  dans  $D(g_{1},...,g_{n})$, on doit  montrer que  pour toute  fonction  $f\in \mc O(D(g_{1},...,g_{n}))$, $f\circ \sigma$ appartient à  l'algèbre  $C(\Gamma,E')$ des fonctions continues de   $\Gamma$ dans  $E'$. 
  Or 
    \begin{align*}
 q: \mc O((H_{E'})^{n+1}\modmod H^{0}_{E'}) &\to  \mc O(D(g_{1},...,g_{n}))
  \\
  f& \mapsto  [g\mapsto f(g_{1},...,g_{n},g)]
 \end{align*}
    est un  morphisme surjectif  puisqu'il est la  composée de deux morphismes surjectifs:    \begin{itemize}
     \item 
    le  morphisme évident 
        \begin{align*}
\mc O((H_{E'})^{n+1}\modmod H^{0}_{E'}) &\to  \mc O(H_{E'} \modmod C(g_{1},...,g_{n})) =\mc O(H_{E'})^{C(g_{1},...,g_{n})}  \\
  f& \mapsto  [g\mapsto f(g_{1},...,g_{n},g)]
 \end{align*}  
 (où $C(g_{1}, ...,g_{n})$ agit  par conjugaison sur $H_{E'}$) 
 est surjectif car l'orbite de  $(g_{1},...,g_{n})$ par conjugaison par 
 $H^{0}_{E'}$ est une sous-variété affine fermée de  $(H_{E'})^{n}$ qui s'identifie  au  quotient  $H^{0}_{E'}/C(g_{1},...,g_{n})$ (on rappelle que le quotient d'un groupe réductif par un sous-groupe réductif est affine), 
 \item 
la  restriction 
     $\mc O(H_{E'})^{C(g_{1},...,g_{n})} \to \mc O(D(g_{1},...,g_{n}))$ est surjective parce que la   restriction $\mc O(H_{E'}) \to \mc O(D(g_{1},...,g_{n}))$  est  évidemment   surjective, $C(g_{1},...,g_{n})$ agit trivialement sur  $\mc O(D(g_{1},...,g_{n}))$  et toute représentation du groupe réductif   $C(g_{1},...,g_{n})$ est complètement réductible (voir la  discussion au sujet de  l'opérateur  de Reynolds dans le paragraphe  I.1 de \cite{git}).  
     \end{itemize}
 Le morphisme 
   \begin{align*}
\mc O((H_{E'})^{n+1}\modmod H^{0}_{E'}) &\to  C(\Gamma,E')  \\
  f& \mapsto  [\gamma\mapsto \Xi_{n+1}(f)(\gamma_{1},...,\gamma_{n},\gamma)]
 \end{align*}     
 se factorise  par  $q$ car, dans les  notations précédentes, 
 $\Xi_{n+1}(f)(\gamma_{1},...,\gamma_{n},\gamma)=f(g_{1},...,g_{n},g)$ et  
on a vu que  $g\in D(g_{1},...,g_{n})$ pour tout  $\gamma\in \Gamma$.     
  Comme $q$ est surjectif,   $f\mapsto f\circ \sigma$ est un morphisme  d'algèbres 
  $\mc O(D(g_{1},...,g_{n}))
\to C(\Gamma,E')$ bien défini. Le troisième point est démontré.  

On montre  maintenant que pour tout   
 $m\in \N^{*}$, $f\in \mc O((H)^{m}\modmod H^{0})$ et  $(\delta_{1},...,\delta_{m})\in \Gamma^{m}$, on a  
\begin{gather}\label{egalite-f-sigma-delta-m} f(\sigma(\delta_{1}),...,\sigma(\delta_{m}))=
 \big(\Xi_{m}(f)\big)(\delta_{1},...,\delta_{m}).\end{gather}
Par les mêmes  arguments que  dans la preuve du  \lemref{lem-gammaCD} il existe  $h_{1},...,h_{m}\in D(g_{1},...,g_{n})$ tels que 
 $(g_{1},...,g_{n},h_{1},...,h_{m})\in \xi_{n+m}^{ss}(\gamma_{1},...,\gamma_{n},\delta_{1},...,\delta_{m} )$ et de plus  on a  
 $h_{j}=\sigma(\delta_{j})$ pour tout  $j\in \{1,...,m\}$. 
En appliquant la  condition a) à l'injection  
$\{n+1,...,n+m\}\subset \{1,...,n+m\}$ on voit que 
 $(\sigma(\delta_{1}),...,\sigma(\delta_{m}))=(h_{1},...,h_{m})$ est au-dessus de  $\xi_{m}(\delta_{1},...,\delta_{m} )$ et   l'égalité  \eqref{egalite-f-sigma-delta-m} en résulte. 
 
 Le fait que pour $m,\delta_{1},...,\delta_{m}$ comme ci-dessus, le $(n+m)$-uplet 
 $(\sigma(\gamma_{1}),...,\sigma(\gamma_{n}),\sigma(\delta_{1}),...,\sigma(\delta_{m}))$ soit semi-simple
 implique que l'adhérence de l'image de $\sigma$ est un sous-groupe réductif de $H$. La construction précédente de $\sigma(\gamma)$ était entièrement nécessaire, et dépendait seulement du choix de $(g_{1},...,g_{n})$ dans la classe de conjugaison semi-simple associée  à $\xi^{ss}_{n}(\gamma_{1}, ..., \gamma_{n})$. Donc $\sigma$ est unique à conjugaison près par $H^{0}(\Qlbar)$.      \cqfd
       
   L'action de 
       $\mc B$ sur 
       $C_{c}^{\rm{cusp}}(\Bun_{G,N}(\Fq)/\Xi,\Qlbar)$ induit évidemment une action  de $\mc B^{\rm{red}} $, d'où une  décomposition        \begin{gather}\label{dec-param} 
       C_{c}^{\rm{cusp}}(\Bun_{G,N}(\Fq)/\Xi,\Qlbar)=\oplus_{\nu} \mf H_{\nu}
       \end{gather}
       où $\nu$ parcourt   les  caractères de  $\mc B^{\rm{red}} $ 
       (c'est-à-dire aussi   de  $\mc B$). 
      C'est ce que nous appelons    ``décomposition spectrale'' dans  le \thmref{intro-thm-ppal}. 
           Par la définition-proposition \ref{cor-annul-Lambda-n} et  la \propref{Xi-n} (appliquée à   $H=H^{0}=\wh G_E$) on associe à chaque  caractère  $\nu$ un  morphisme 
       $\sigma:\on{Gal}(\ov F/F)\to \wh G(\Qlbar)$
tel que 
\begin{itemize}
\item [] (C1) $\sigma$ prend ses valeurs  dans  $\wh G(E')$, où $E'$ est une extension finie  de $E$ (donc de $\Ql$), et il est continu, 
\item [] (C2) l'adhérence de Zariski de son  image est réductive, 
\item[] (C3) 
pour tout  $n\in \N^{*}$, on a   
\begin{gather} \label{Theta-nu}
 f(\sigma(\gamma_{1}),...,\sigma(\gamma_{n}))= \big(\Theta_{n}^{\nu}(f)\big)(\gamma_{1},...,\gamma_{n}), 
\end{gather}
\item [] (C4) $\sigma$ se factorise à travers $\pi_{1}(X\sm N, \ov\eta)$. 
\end{itemize}

De plus  $\nu$ et  la  classe de conjugaison de $\sigma$ par 
$\wh G(\Qlbar)$ se déterminent mutuellement de fa\c con unique grâce à la   condition (C3).  
Le  théorème suivant est le \thmref{intro-thm-ppal} de l'introduction
(avec un énoncé plus précis du lien avec  les opérateurs d'excursion, 
qui sont  maintenant construits). 
    
       \begin{thm} \label{dec-param-cor-thm}    On a une  décomposition de 
       $C_{c}(K_{N}\backslash G(\mb A)/K_{N},\Qlbar)$-modules  
        \begin{gather}\label{dec-param-cor} 
       C_{c}^{\rm{cusp}}(\Bun_{G,N}(\Fq)/\Xi,\Qlbar)=\oplus_{\sigma} \mf H_{\sigma}
       \end{gather} 
       où la somme  est indexée par  les classes de conjugaison  par $\wh G(\Qlbar)$  de morphismes
        $$\sigma:\pi_{1}(X\sm N, \ov\eta)\to \wh G(\Qlbar) $$
 vérifiant les   conditions (C1) et (C2). 
        Cette  décomposition est déterminée de manière unique par la décomposition   \eqref{dec-param} grâce à la  condition (C3).  Elle  est compatible avec  l'isomorphisme de Satake en les places non ramifiées: pour tout   $\sigma$ tel que $\mf H_{\sigma}\neq 0$, 
     pour toute représentation  irréductible  $V$ de $\wh G$ et toute  place  $v\in |X|\sm |N|$,  $T(h_{V,v})$ 
     agit  sur  $\mf H_{\sigma}$ par le scalaire  $\chi_{V}(\sigma(\Frob_{v}))$, où 
     $\Frob_{v}\in \pi_{1}(X\sm N, \ov\eta)$ est un élément de Frobenius en $v$. 
           \end{thm}
       \noindent{\bf Démonstration.}      Il reste seulement à montrer la compatibilité avec Satake. 
   Soit  $V$ une représentation  irréductible  de $\wh G$. 
   On sait que $T(h_{V,v})$ respecte $\mf H_{\sigma}$ et d'après 
le   e) de la  définition-proposition \ref{cor-annul-Lambda-n}  et   la  condition (C3), $T(h_{V,v})$ agit sur $\mf H_{\sigma}$ avec l'unique valeur propre généralisée $\chi_{V}(\sigma(\Frob_{v}))$. Mais, comme on l'a déjà dit dans la \remref{rem-B-reduite?},  
les opérateurs de Hecke en les places non ramifiées sont diagonalisables 
(car ils sont normaux pour la structure hermitienne standard sur $C_{c}^{\rm{cusp}}(\Bun_{G,N}(\Fq)/\Xi,\C)$). 
Donc $T(h_{V,v})$ agit sur $\mf H_{\sigma}$ par multiplication par le scalaire 
      $\chi_{V}(\sigma(\Frob_{v}))$ et on a montré que   $\sigma$ est compatible en  $v$ avec  l'isomorphisme de Satake.
    \cqfd

  \section{Cas des groupes non nécessairement déployés}
  \label{para-non-deploye}
    
    \subsection{Enoncé du théorème principal}
    \label{para-non-deploye-statement}
    
    Soit $G$ un groupe réductif lisse et géométriquement connexe  sur  $F$. 
On note $U$ l'ouvert maximal de $X$ tel que $G$ se prolonge en un   schéma  en groupes lisse et réductif  sur  $U$. D'après les paragraphes  5.1.9 et   4.6 de \cite{bruhat-tits}
 (qui m'ont été indiqués par Jochen Heinloth), 
 on peut choisir  
 un modèle entier  parahorique de $G$ en tous les  points de  $X\sm U$. 
 En recollant ces modèles entiers sur  $U$ et sur les voisinages formels des   points de $X\sm U$, on obtient  un schéma en groupes lisse sur  $X$ que l'on note encore  $G$. Ainsi  $G$ est   un schéma en groupes lisse sur  $X$,  réductif sur $U$, de type parahorique   en les  points de $X\sm U$, et dont toutes les fibres  sont géométriquement connexes. 
 
  D'après la  proposition 1 de \cite{heinloth-unif}, le champ  $\Bun_{G}$, dont les points sur un schéma $S$ classifient les   $G$-torseurs sur $X\times S$, est un champ algébrique lisse, localement  de type fini.  
 Une autre preuve de ce résultat apparaît  dans \cite{hartl}. On rappelle les énoncés dans le lemme suivant. 
 
 \begin{lem} \label{lem-behrend-hartl-heinloth}
 (\cite{behrend-thesis}, \cite{heinloth-unif} exemple (1) page 504,  et \cite{hartl}, proposition 2.2 et théorème 2.5)
 Soit $G$ comme ci-dessus. Il existe un fibré vectoriel $\mc V$ de rang $r$ sur $X$ muni d'une trivialisation de  $\det(\mc V)$ et une représentation fidèle $\rho:G\to SL(\mc V)$ telle que les quotients $SL(\mc V)/G$ (et donc $GL(\mc V)/G$) soient quasi-affines. 
 Alors le morphisme $\rho_{*}:\Bun_{G}\to Bun_{GL_{r}}^{0}$ est représentable, quasi-affine et de présentation finie (on a noté $\Bun_{GL_{r}}^{0}$ la composante connexe de $\Bun_{GL_{r}}$ classifiant les fibrés de degré $0$).

Pour tout copoids 
   dominant $\mu$  pour $GL_{r}$,   
$\Bun_{GL_{r}}^{0,\leq \mu}$ est de type fini (voir \cite{wang-bundle} pour une preuve et les références) et donc $\rho_{*}^{-1}(\Bun_{GL_{r}}^{0,\leq \mu})$ est un ouvert de type fini de $\Bun_{G}$. Pour deux choix différents de $\mc V$ et $\rho$ les systèmes inductifs d'ouverts de $\Bun_{G}$ sont comparables. 
  \end{lem}
  
  Seule la dernière phrase n'apparaît pas dans les références, mais elle résulte immédiatement du fait que ces ouverts sont de type fini.

 Pour définir des troncatures $ \Bun_{G}^{\leq \mu}$ il serait inutilement compliqué de chercher à comprendre les troncatures de Harder-Narasimhan dans un cadre non déployé et on prefère utiliser un plongement dans $SL_{r}$ comme dans le lemme précédent. De plus on va définir ces troncatures à l'aide  de $G^{\mr{ad}}$ au lieu de $G$ pour que les troncatures qui s'en déduisent dans les champs de chtoucas soient invariantes par l'action de $\Xi$. 
 
  On applique le lemme précédent  à $G^{\mr{ad}}$,  et on choisit donc un fibré vectoriel $\mc V$ de rang $r$ sur $X$ muni d'une trivialisation de  $\det(\mc V)$ et un plongement $\rho:G^{\mr{ad}}\to SL(\mc V)$.  Pour tout copoids 
   dominant $\mu$  pour $GL_{r}$  
     on définit l'ouvert de type fini $ \Bun_{G^{\mr{ad}}}^{\leq \mu}$  comme $(\rho_{*})^{-1}(\Bun_{GL_{r}}^{0,\leq \mu})$.     On définit  $ \Bun_{G }^{\leq \mu}$  comme  l'image inverse de  $ \Bun_{G^{\mr{ad}}}^{\leq \mu}$. C'est un ouvert de $\Bun_{G}$ qui n'est   de type fini  que si $G$ est  semi-simple. Des choix différents de $\mc V$ et $\rho$ conduiraient à des systèmes inductifs comparables. 
     
 De plus, en anticipant sur les notations du paragraphe \ref{reminder-satake-twisted}, 
 pour toute représentation $W$ de ${}^{L} G$ il existe $\kappa$ tel que pour tout $\mu$ et tout point 
 $(\mc G_{0}\to \mc G_{1})$  de 
 $\Hecke_{\{0\},W} ^{(\{0\})}$, 
  si 
   $\mc G_{0}$ appartient à   $ \Bun_{G}^{\leq \mu}$ alors $ \mc G_{1}$ appartient à $ \Bun_{G}^{\leq \mu+\kappa}$.  
     
  Soit $N\subset X$ un sous-schéma fini. On note   
  $K_{N}$  le noyau  de $G(\mathbb O)\to G(\mc O_N)$. 
  On note 
  $\Bun_{G,N}$ le champ  dont les  points sur un  schéma $S$ classifient la donnée d'un   $G$-torseur $\mc G$ 
  sur $X\times S$ et d'une  trivialisation de  $\restr{\mc G}{N\times S}$. 
   Comme les fibres    de $G$ sont géométriquement connexes, 
   la restriction  à la  Weil  $G_{N}$ de $G$ de  $N$ à $\on{Spec} \Fq$ est connexe. On peut donc    appliquer le   théorème de Lang à  $G_{N}$. Par conséquent  $\Bun_{G,N}(\Fq)$ est un   
  $G(\mc O_N)$-torseur sur  le  groupoïde $\Bun_{G}(\Fq)$. 
  
  On note $\wh N=|N|\cup (X\sm U)$, de sorte que $X\sm \wh N$ est l'ouvert non ramifié.

         Le  $L$-groupe ${}^{L }G$ est un  produit   semi-direct
   $\wh G\rtimes \on{Gal}(\wt F/F)$ où  $\wt F$ est l'extension finie galoisienne de $F$  telle  que  $\on{Gal}(\wt F/F)$ soit l'image  de $\on{Gal}(\ov F/F)$ dans le groupe des  automorphismes du diagramme de  Dynkin  de $G$. 
 Le  produit   semi-direct est pris  pour l'action de $\on{Gal}(\wt F/F)$ sur 
 $\wh G$ qui préserve un épinglage, voir \cite{borel-corvallis}.

  Pour tout $v\in |U|$, $G(\mc O_{v})$ est  hyperspécial, donc   $G(F_{v})$ est quasi-déployé et déployé sur une extension  non ramifiée de $F_{v}$. Donc   $\wt F/F$ est non ramifié sur  $U$. On note  $\wt{U}$ le revêtement galoisien  de $U$, dont le corps des  fonctions est  $\wt F$. 
   On considère  ${}^{L }G$ comme  un  groupe algébrique sur $\Ql$ et pour toute extension   $E\supset \Ql$ on note  ${}^{L }G_{E}$ le groupe obtenu par extension des scalaires à   $E$ et on note 
  ${}^{L }G(E)=\wh G(E)\rtimes \on{Gal}(\wt F/F)$  le  groupe de ses points à valeurs dans $E$.

          Pour toute place  $v\in |X|$, on note   $\wt F_{v}$ l'extension galoisienne finie  de  $F_{v}$ telle  que   $\on{Gal}(\wt F_{v}/F_{v})$ soit  l'image  de $\on{Gal}(\ov F_{v}/F_{v})$ dans le groupe des  automorphismes du   diagramme de  Dynkin de $G$. 
On note  ${}^{L}G_{v}=\wh G \rtimes \on{Gal}(\wt F_{v}/F_{v})$ le  $L$-groupe local. 
     On a  un plongement  
     $$\on{Gal}(\wt F_{v}/F_{v})\subset 
   \on{Gal}(\wt F/F)  \text{ \ \ et \ \ } {}^{L}G_{v} \subset {}^{L}G$$
  bien défini  à conjugaison près et  déterminé de fa\c con unique  par le choix      
     d'un plongement $\ov F\subset \ov F_{v}$. 
    Si $v\in |U|$, $\wt F_{v}/F_{v}$ est non ramifié et donc   $\on{Gal}(\wt F_{v}/F_{v})$
     est cyclique, avec $\Frob_{v}$  comme générateur canonique.

     \begin{rem}\label{rem-une-classe-ker1}   On a  une  inclusion évidente 
   \begin{gather}\label{incl-double-quot}
   G(F)\backslash G(\mb A)/K_{N}\subset 
 \Bun_{G,N}(\Fq)\end{gather}
dont  l'image est formée des  $G$-torseurs   localement triviaux pour la  topologie de Zariski. En général, comme on l'a déjà mentionné dans la \remref{quotient-adelique-deploye},     \begin{gather}\label{dec-alpha-general-ker1}  \Bun_{G,N}(\Fq)=
  \bigcup_{\alpha\in \ker^{1}(F,G)}G_{\alpha}(F)\backslash G_{\alpha}(\mb A)/K_{N}\end{gather} où la réunion est  disjointe, $\ker^{1}(F,G)$ est fini et $G_{\alpha}$ est la  forme intérieure pure  de $G$  obtenue par torsion  par $\alpha$. On rappelle que $$ \ker^{1}(F,G)=\mr{Ker}(
  H^{1}(F,G)\to \prod_{v} H^{1}(F_{v},G)). $$ Pour tout $\alpha\in \ker^{1}(F,G)$  
  on fixe un $G$-torseur sur $F$ ayant $\alpha$ comme classe d'isomorphisme, 
  on l'étend à $X\sm \wh  N$, et on le munit d'une trivialisation sur chaque corps local $F_{v}$ (venant en presque toute place $v$ de $X\sm \wh N$ d'une trivialisation sur $\mc O_{v}$). On note $G_{\alpha}$ le groupe d'automorphismes de ce $G$-torseur. On possède donc un isomorphisme $G_{\alpha}(\mathbb A)=G(\mathbb A)$, qui donne un sens au quotient par $K_{N}$ dans le membre de droite de 
  \eqref{dec-alpha-general-ker1}. L'égalité   \eqref{dec-alpha-general-ker1} résulte facilement du fait que, pour tout $v\in |X|$,  $H^{1}(\mc O_{v},G)=0$ (grâce au théorème de Lang, comme on l'a montré plus haut). Pour plus de détails on renvoie au lemme 1.1 de \cite{ngo-invent2006}.

   D'après    Kottwitz~\cite{kottwitz1,kottwitz2} et  Nguyen Quoc Thang~\cite{thang} théorème 2.6.1   pour l'extension en caractéristique $p$, $\ker^{1}(F,G)$ est le  dual de $\ker^{1}(F,Z_{\wh G}(\Qlbar))$.      
\end{rem}
    
      On fixe un réseau $\Xi\subset Z(\mb A)/Z(F)$. 
   On définit    $C_{c}^{\mr{cusp}}(\Bun_{G,N}(\Fq)/\Xi,E)$ comme le sous-espace  de   $C_{c}(\Bun_{G,N}(\Fq)/\Xi,E)$ formé des    éléments dont  l'image par la  correspondance 
   $$\Bun_{G,N}(\Fq)/\Xi\leftarrow (  \Bun_{P,N}(\Fq)\times_{P(\mc O_{N})}G(\mc O_{N}) ) /\Xi \rightarrow (  \Bun_{ M,N}(\Fq)\times_{P(\mc O_{N})}G(\mc O_{N}) )/\Xi$$ s'annule  dans $C(\Bun_{ M,N}(\Fq)/\Xi,E)$ pour tout sous-groupe parabolique  $P$ de $G$ de Levi $ M$ (on renvoie à \cite{these-cong} pour un cadre général, où de telles correpondances sont étudiées). Grâce à la remarque précédente on a la définition équivalente  \begin{gather}\label{def-cusp}C_{c}^{\mr{cusp}}(\Bun_{G,N}(\Fq)/\Xi,E)= \bigoplus_{\alpha\in \ker^{1}(F,G)}C_{c}^{\mr{cusp}}(G_{\alpha}(F)\backslash G_{\alpha}(\mb A)/K_{N}\Xi,E).\end{gather}
          C'est un $E$-espace vectoriel  muni d'une  action de       $C_{c}(K_{N}\backslash G(\mb A)/K_{N},E)$. Il est  de dimension finie par \cite{harder} et la proposition 5.2 de \cite{borel-jacquet}:  l'idée est que la cuspidalité implique l'annulation sur les $G$-torseurs avec niveau $N$ très instables (c'est-à-dire possédant, pour un certain sous-groupe parabolique propre une réduction au sous-groupe de  Levi associé admettant un unique  relèvement au sous-groupe parabolique), donc il existe $\mu$ tel que toute fonction cuspidale soit supportée sur $\Bun_{G,N}^{\leq \mu}(\Fq)/\Xi$ qui est fini.       
  
    Soit $v\in |U|$ (en fait on n'utilisera ce qui va suivre que  pour $v\in |X|\sm \wh N\subset  |U|$).   L'isomorphisme  de  Satake   
(voir \cite{satake,cartier-satake,borel-corvallis,blasius-rogawski-pspm}) 
est un  isomorphisme d'anneaux 
\begin{gather}\label{Isom-Satake-non-deploye}\mathscr S: \mc O(\wh G_E\rtimes \Frob_{v}\modmod \wh G_E )\to C_{c}(G(\mc O_{v})\backslash G(F_{v})/G(\mc O_{v}),E)\end{gather}
où le membre de gauche est l'anneau des  fonctions régulières  sur la classe à gauche (ou à droite) par 
$\wh G_E$ égale à  $\wh G_E\rtimes \Frob_{v}\subset {}^{L}G_{v,E}$,  qui sont 
invariantes par conjugaison par $\wh G_E$. 
Pour toute représentation $E$-linéaire 
$V$ de ${}^{L} G_{v}$, on note   $h_{V,v}\in C_{c}(G(\mc O_{v})\backslash G(F_{v})/G(\mc O_{v}),E)$ l'image par l'isomorphisme de Satake $\mathscr S$ de $\restr{\chi_{V}}{\wh G_E\rtimes \Frob_{v}}$. 
Ces fonctions $\restr{\chi_{V}}{\wh G_E\rtimes \Frob_{v}}$  engendrent  $\mc O(\wh G_E\rtimes \Frob_{v}\modmod \wh G_E )$ quand  $V$ varie. 
En effet les $\chi_{V}$ engendrent le terme de gauche dans 
$$\mc O({}^{L} G_{v,E}\modmod {}^{L} G_{v,E})\twoheadrightarrow 
\mc O(\wh G_E\rtimes \Frob_{v}\modmod {}^{L} G_{v,E})\isom 
\mc O(\wh G_E\rtimes \Frob_{v}\modmod \wh G_E ),$$
où le premier morphisme est simplement la restriction. Le deuxième est un isomorphisme car  les  classes de conjugaison
par  $\wh G_{E}$ et par  ${}^{L} G_{v,E}$ d'un élément de la forme  $g \rtimes \Frob_{v}$  sont égales (puisque 
$g \rtimes \Frob_{v}$ commute avec lui-même), donc  ${}^{L} G_{v,E}/ \wh G_E$
agit trivialement sur $\wh G_E\rtimes \Frob_{v}\modmod \wh G_E$. 
Par conséquent  
$[V]\mapsto h_{V,v}$ est un homomorphisme surjectif   d'anneaux  
 $\on{Rep}_{E}({}^{L}G_{v}) \to C_{c}(G(\mc O_{v})\backslash G(F_{v})/G(\mc O_{v}),E)$. La fonction  $h_{V,v}$ appartient à 
 $C_{c}(G(\mc O_v)\backslash G(F_{v})/G(\mc O_v), \mc O_{E})$ 
 lorsque $V$ est défini sur   $\mc O_{E}$.

 Dans le paragraphe \ref{excur-op-non-split} nous définirons une famille commutative d'opérateurs d'excursion 
    $$S_{I,W,x,\xi,(\gamma_i)_{i\in I}}\in \on{End}_{C_{c}(K_{N}\backslash G(\mb A)/K_{N},E)}\Big( C_{c}^{\mr{cusp}}(\Bun_{G,N}(\Fq)/\Xi,E)\Big)$$
       où $I$ est un ensemble fini,  $W$ est une   représentation $E$-linéaire de $({}^{L }G)^{I}$, $x\in W$ et  $\xi\in W^{*}$ sont invariants sous l'action diagonale de  $\wh G$, et $(\gamma_i)_{i\in I}\in (\on{Gal}(\ov F/F))^{I}$.  
      
       \begin{thm}\label{dec-param-cor-thm-non-deploye}     On a une   décomposition canonique de  $C_{c}(K_{N}\backslash G(\mb A)/K_{N},\Qlbar)$-modules 
        \begin{gather}\label{dec-param-cor-non-deploye} 
      C_{c}^{\mr{cusp}}(\Bun_{G,N}(\Fq)/\Xi,\Qlbar)=\oplus_{\sigma} \mf H_{\sigma}
       \end{gather} 
      indexée par  les paramètres de Langlands, c'est-à-dire les 
classes de conjugaison par $\wh G(\Qlbar)$  de   morphismes
        $$\sigma:\on{Gal}(\ov F/F) \to {}^{L} G(\Qlbar) $$
      tels que 
\begin{itemize}
\item [] (C'1) $\sigma$ prend ses valeurs dans  ${}^{L}G(E')$, où $E'$ est une extension finie  de $E$ (donc de $\Ql$), et il est continu  et non ramifié en dehors de  $\wh N$ (c'est-à-dire qu'il se factorise à travers $\pi_{1}(X\sm \wh N, \ov\eta)$), 
\item [] (C'2) l'adhérence de Zariski de son  image est réductive, 
\item []  (C'5)  on a  la commutativité du   diagramme 
 \begin{gather}\label{diag-sigma}
 \xymatrix{
\on{Gal}(\ov F/F) \ar[rr] ^{\sigma}
\ar[dr] 
&& {}^{L} G(\Qlbar) \ar[dl] 
 \\
& \on{Gal}(\wt F/F) }\end{gather}
\end{itemize}
Cette  décomposition découle de  la  décomposition spectrale de la famille commutative des opérateurs d'excursion, au sens suivant: 
$\mf H_{\sigma}$ est l'espace propre généralisé     correspondant au système de valeurs propres $\s{\xi, (\sigma(\gamma_{i}))_{i\in I}.x}$ pour 
$S_{I,W,x,\xi,(\gamma_i)_{i\in I}}$ quand  $I,W,x,\xi$ et $(\gamma_i)_{i\in I}$ varient. 
Il est  compatible avec  l'isomorphisme de Satake: pour tout  $\sigma$,  toute place $v\in |X|\sm \wh N$, et pour toute représentation irréductible $V$ de ${}^{L} G_{v}$, $h_{V,v}$ agit  sur $\mf H_{\sigma}
$ par multiplication par le scalaire  $\chi_{V}(\sigma(\Frob_{v}))$. 
Enfin la  décomposition \eqref{dec-param-cor-non-deploye} est compatible avec la limite sur $N$. 
           \end{thm}

    \begin{rem}
En passant à la limite    sur $N$, 
 on obtient une 
 décomposition canonique,  indexée par les paramètres de Langlands, de la 
   représentation de $G(\mb A)$ $$\varinjlim_{N } C_{c}^{\mr{cusp}}(\Bun_{G,N}(\Fq)/\Xi,\Qlbar)= \bigoplus_{\alpha\in \ker^{1}(F,G)}C_{c}^{\mr{cusp}}(G_{\alpha}(F)\backslash G_{\alpha}(\mb A)/\Xi,E). $$ 
  Autrement dit  pour toute représentation lisse   irréductible $\pi$ de $G(\mb A)$, le  $\Qlbar$-espace vectoriel de dimension finie    
$$\on{Hom}_{G(\mb A)}\Big(\pi, \varinjlim_{N } C_{c}^{\mr{cusp}}(\Bun_{G,N}(\Fq)/\Xi,\Qlbar)\Big)$$
admet une   décomposition canonique indexée par les paramètres de Langlands qui sont  compatibles  par   l'isomorphisme de Satake avec $\pi$  en toutes les places où $\pi$ est non ramifié. La conjecture \ref{conj-algebre-B-Qbar} ci-dessous affirme que 
(si $\pi$ est défini sur $\ov\Q$) 
cette décomposition est définie sur  $\ov\Q$, indexée par des paramètres de Langlands motiviques et ne dépend pas de   $\ell$  ni du plongement $\ov\Q \hookrightarrow \Qlbar$. 
\end{rem}

     \subsection{Compléments, remarques et conjectures}
     \label{complements-non-split}
     
    \subsubsection{Compatibilité avec la fonctorialité}      La  proposition suivante 
    (dont l'énoncé a émergé  après des discussions avec  Vladimir Drinfeld et Erez Lapid) affirme que la  décomposition \eqref{dec-param-cor-non-deploye} 
                est compatible avec les cas triviaux  de fonctorialité. 
         Bien sûr on espère aussi qu'elle est  compatible avec les cas non triviaux, 
         comme la    correspondance theta, 
         et d'ailleurs grâce au lien entre notre construction et  \cite{brav-var} (expliqué dans le chapitre  \ref{subsection-link-langl-geom}) cela devrait résulter 
         de la géométrisation du noyau theta  par Sergey Lysenko \cite{sergey-theta,sergey-theta-SO-Sp}.

Soit $G'$ un autre  groupe réductif sur $F$ et $\Upsilon:G\to G'$ un morphisme de groupes sur $F$ dont  l'image est un sous-groupe distingué de $G'$. La fonctorialité ``triviale'' pour les  $L$-groupes 
expliquée dans le paragraphe  2.5 de \cite{borel-corvallis} 
fournit un morphisme ${}^{L}\Upsilon: {}^{L}G'\to {}^{L}G$ (en supposant les  $L$-groupes construits à l'aide d'une même extension galoisienne finie  $\wt F$ déployant à la fois  $G$ et $G'$). 
      
      Soit $\Xi$ et $\Xi'$ des réseaux  dans $Z(F)\backslash Z(\mathbb A)$ et 
      $Z'(F)\backslash Z'(\mathbb A)$ (où $Z'$ désigne le centre  de $G'$) tels  que $\Upsilon(\Xi)\subset \Xi'$. Alors $\Upsilon$ induit 
      $$\beta_{\Upsilon}: \varprojlim_{N}  \Bun_{G,N}(\Fq)/\Xi\to \varprojlim_{N}  \Bun_{G',N}(\Fq)/\Xi'$$ et donc  
     \begin{gather}\label{restriction-G-G'}\beta_{\Upsilon}^{*}: \varinjlim_{N} C_{c}^{\mr{cusp}}(\Bun_{G',N}(\Fq)/\Xi',E)\to \varinjlim_{N}  C_{c}^{\mr{cusp}}(\Bun_{G,N}(\Fq)/\Xi,E)\end{gather}
     défini par 
     $\beta_{\Upsilon}^{*}(h)= h\circ \beta_{\Upsilon}$. 
       On considère la  décomposition \eqref{dec-param-cor-non-deploye} 
      (après être passé à la limite sur $N$) 
      et la  décomposition  analogue 
               \begin{gather}\label{dec-param-cor-non-deploye-G'} 
    \varinjlim_{N}  C_{c}^{\mr{cusp}}(\Bun_{G',N}(\Fq)/\Xi',E)=\oplus_{\sigma'} \mf H'_{\sigma'}. 
       \end{gather}

      \begin{prop}\label{functoriality-L-groups-H-sigma}
  Pour tout paramètre de Langlands $\sigma'$ pour $G'$,  on a 
  $\beta_{\Upsilon}^{*}(\mf H'_{\sigma'})\subset \mf H_{\sigma}$ 
avec  $\sigma={}^{L}\Upsilon\circ \sigma'$.
 \end{prop}\noindent 
La preuve sera donnée juste après la preuve du \thmref{dec-param-cor-thm-non-deploye}.

    \subsubsection{Paramètres d'Arthur}\label{subsubsection-Arthur}

 Bien sûr on aimerait montrer que les  paramètres de Langlands $\sigma$ qui apparaissent dans la décomposition \eqref{dec-param-cor-non-deploye}  proviennent de  paramètres d'Arthur elliptiques.  On rappelle qu'un paramètre d'Arthur  est une classe de conjugaison par $\wh G(\Qlbar)$ de 
   morphisme $$
  \psi : \on{Gal}(\ov F/F) \times SL_{2}(\Qlbar)\to 
 {}^{L} G(\Qlbar) \text{ (algébrique sur  $SL_{2}(\Qlbar)$),}$$ dont la restriction à  $\on{Gal}(\ov F/F)$ est non ramifiée sur un ouvert dense et prend ses valeurs  dans une extension finie de  $\Ql$,  et qui est continu, fait commuter le diagramme \eqref{diag-sigma} et est    tel que 
 
  (A1) la  restriction de $\psi$ à $\on{Gal}(\ov F/F)$ est pure de poids  $0$ (c'est-à-dire  que son  image par toute représentation de ${}^{L}G$ est pure de poids $0$). 
 
 De plus  $ \psi$ est dit elliptique  si 

 (A2) le centralisateur  de $\psi$ dans  $\wh  G(\Qlbar)$ est fini  modulo $(Z(\wh  G)(\Qlbar))^{\on{Gal}(\wt F/F)}$. 
  
 Le  paramètre de Langlands associé à $\psi$ est  $\sigma_{\psi}: \on{Gal}(\ov F/F) \to {}^{L} G(\Qlbar)$ défini par 
 $$\sigma_{\psi}(\gamma)=\psi\Big(\gamma, \begin{pmatrix} |\gamma|^{1/2} & 0 \\
 0 & |\gamma|^{-1/2}
 \end{pmatrix}\Big)$$  où $|\gamma|^{1/2}$ est bien défini grâce au choix d'une racine carrée  de $q$. 
 
Le paramètre d'Arthur, c'est-à-dire la 
 classe de conjugaison par  $
  \wh G(\Qlbar)$  de  $\psi$,  est déterminé  de manière unique par la  classe de conjugaison par $
  \wh G(\Qlbar)$ de   $\sigma_{\psi}$ (Malcev,  cf le  corollaire 4.2 de \cite{kostant-betti}). 
  
  \begin{rem}\label{rem-A1-A2}
  Soit $\psi$ un paramètre d'Arthur  tel que $\sigma_{\psi}$ apparaisse dans 
  \eqref{dec-param-cor-non-deploye}. Par la théorie du corps de classe on sait que  
   pour tout caractère $\nu$ de ${}^{L} G(\Qlbar)$,  $\nu\circ \sigma_{\psi}$ est pure de poids  $0$. Par conséquent la  condition (A2) implique la  condition (A1). 
   En effet si on fixe   $\iota:\Qlbar\isom\C$ et si la  restriction de  $\psi$ à  $\on{Gal}(\ov F/F)$ n'est pas  $\iota$-pure de poids  $0$, alors le   $\iota$-poids  fournit un cocaractère non trivial à valeurs dans le centralisateur  de $\psi$, dont la  composée avec  tout caractère $\nu$ de  ${}^{L} G(\Qlbar)$ est  trivial, et cela contredit la  condition (A2). 
    \end{rem}
  
  \begin{conj}\label{conj-arthur}
  Tout  paramètre de Langlands  $\sigma$ apparaissant dans la  décomposition \ref{dec-param-cor-non-deploye}  est de la forme   $\sigma_{\psi}$  avec  $\psi$ un  paramètre d'Arthur elliptique (bien défini à conjugaison près par  $\wh G(\Qlbar)$). 
    \end{conj}

\begin{rem}\label{rem-apres-conj-arthur}
Pour montrer la conjecture précédente on    peut espérer utiliser  les actions du $SL_{2} $ de Lefschetz 
sur la cohomologie de compactifications des champs de chtoucas
(il faudrait évidemment considérer la cohomologie en tous degrés et non pas seulement en degré moitié comme dans cet article). Cela était discuté dans la remarque 10.8 et  le chapitre 11 de la version 3 de cet article sur arXiv mais on signale que le lemme 11.5 était erroné car l'action des morphismes de Frobenius partiels manquait dans les hypothèses (et il n'est pas évident de définir une telle action sur la cohomologie d'une compactification). 
   \end{rem}

 \begin{rem} \label{remark-discrete-part} On aimerait en fait  obtenir une  décomposition comme \eqref{dec-param-cor-non-deploye}  pour toute la partie discrète (et non pas seulement la partie cuspidale). Cette décomposition devrait être indexée par les paramètres d'Arthur elliptiques. 
       \end{rem}                 
                
           \subsubsection{Lien avec les   formules de multiplicités d'Arthur}
             
           Il est tentant de conjecturer une formule de multiplicités pour les espaces $\mf H_{\sigma}$ apparaissant dans la décomposition \eqref{dec-param-cor-non-deploye}.                Comme la  décomposition \eqref{dec-param-cor-non-deploye}  n'est construite que pour la partie cuspidale, et pas encore pour la partie discrète (cf la remarque \ref{remark-discrete-part}), il est préférable pour le moment de ne discuter
     les   formules de multiplicités que pour les  paramètres d'Arthur elliptiques pour lesquels l'action de $SL_{2}$ est triviale, car dans ce cas  le paquet global d'Arthur est conjecturalement  inclus dans la partie cuspidale.          
              
      La formule de multiplicités que nous conjecturons 
      pour   $\mf H_{\sigma}$ est exactement le terme correspondant à $\sigma$ dans la formule de multiplicités d'Arthur usuelle. 
            La raison est la ``compensation'' suivante. 
On rappelle que  dans \eqref{incl-double-quot} la somme  est indexée par 
           $\ker^{1}(F,G)$   et comme les   formules de multiplicités d'Arthur 
          concernent le quotient adélique $G(F)\backslash G(\mb A)/K_{N}\Xi$, 
           on a  l'impression que l'on doit multiplier toutes les multiplicités    par $\sharp   \ker^{1}(F,G)$. Mais les formules de multiplicités d'Arthur sont écrites pour une relation  d'équivalence   $\mc R$ sur l'ensemble des  paramètres d'Arthur elliptiques, qui est plus faible que la    conjugaison par $\wh G(\Qlbar)$, et fait intervenir en plus la torsion  
           par un cocyle  localement trivial  à valeurs dans  $Z_{\wh G}(\Qlbar)$,  c'est-à-dire   un élément de 
  $\ker^{1}(F,Z_{\wh G}(\Qlbar))$. 
Or     $\ker^{1}(F,G)$ est le  dual de $\ker^{1}(F,Z_{\wh G}(\Qlbar))$ donc  ils ont le même cardinal. 

\begin{rem} Pour des groupes autres que $GL_{r}$, 
des  multiplicités $>1$  peuvent apparaître  dans chaque  $\mf H_{\sigma}$. Cela est indépendant du fait que différents $\mf H_{\sigma}$ peuvent correspondre au même système de valeurs propres de Hecke (comme on l'a expliqué dans le paragraphe \ref{intro-rem-suppl}). 
   \end{rem}

Lorsque   $G$ est un tore, les sous-espaces $\mf H_{\sigma}$ ont pour  dimension  $1$, et  la  décomposition \eqref{dec-param-cor-non-deploye}
 est simplement la  bijection entre caractères de  $\varprojlim_{N}\Bun_{G,N}(\Fq)$ et paramètres de Langlands modulo conjugaison par  $\wh G(\Qlbar)$. 
  
    \begin{rem}  Il est clair que  la  décomposition \eqref{dec-param-cor-non-deploye}
 implique une  décomposition plus grossière de 
  $C_{c}^{\mr{cusp}}(\Bun_{G,N}(\Fq)/\Xi,\Qlbar)$ indexée par les paramètres de Langlands  modulo 
 la  relation d'équivalence plus faible  $\mc R$ mentionnée ci-dessus. 
 Cependant on ne sait pas si cette relation d'équivalence plus faible  est compatible avec  l'inclusion 
        $C_{c}^{\rm{cusp}}(G(F)\backslash G(\mb A)/K_{N}\Xi,\Qlbar)\subset C_{c}^{\mr{cusp}}(\Bun_{G,N}(\Fq)/\Xi,\Qlbar)$ qui vient de  \eqref{incl-double-quot}.  
     On ne sait donc pas construire en général  une   décomposition canonique 
        de    $C_{c}^{\rm{cusp}}(G(F)\backslash G(\mb A)/K_{N}\Xi,\Qlbar)$ indexée par les paramètres de Langlands pour la relation d'équivalence plus faible  $\mc R$. 
       \end{rem}

\subsubsection{Indépendance par rapport à  $\ell$}
La philosophie des motifs indique qu'il doit 
exister une notion de ``paramètre de Langlands motivique'' définie sur $\ov\Q$. Si l'on admet les conjectures standard,  un paramètre de Langlands  motivique est  défini comme une classe d'isomorphisme de couples $(T,\kappa)$,   où 
\begin{itemize}
\item 
 $T$ est un  foncteur tensoriel  de la catégorie de  représentations $\ov\Q$-linéaires de dimension finie 
 de ${}^{L}G$ 
vers la  catégorie  $\ov\Q$-linéaire des motifs  purs sur  $F$,
\item $\kappa$ est un  isomorphisme entre la   restriction de $T$ aux représentations de $ {}^{L}G/\wh G=\on{Gal}(\wt F/F)$ et le foncteur naturel des représentations 
de  $\on{Gal}(\wt F/F)$ vers les motifs d'Artin   sur $F$. 
\end{itemize}
    
    Dans un article récent \cite{drinfeld-pro-completion}, Drinfeld a donné un sens inconditionnel à la notion de  paramètre de Langlands motivique. 
 
\begin{conj}\label{conj-algebre-B-Qbar} La décomposition \eqref{dec-param-cor-non-deploye} est définie sur $\ov\Q$ et indexée  par des paramètres de Langlands motiviques. De plus la décomposition    est indépendante de $\ell$ et du plongement $\ov\Q \hookrightarrow \Qlbar$. 
  \end{conj}

\begin{rem}\label{rem-lapid} Bien sûr dans le cas où  $G=GL_{r}$, il est évident que la  décomposition \eqref{dec-param-cor-non-deploye}  est définie sur  $\ov\Q$  
 et  indépendante de 
  $\ell$ et du plongement $\ov\Q \hookrightarrow \Qlbar$, puisque les espaces  $\mf H_{\sigma}$ sont  obtenus par diagonalisation des opérateurs de Hecke en les  places non ramifiées. En effet, grâce au   théorème de Tchebotarev,  une  représentation 
  linéaire   semi-simple $\sigma$ est déterminée par son caractère.  On pourrait aussi invoquer le théorème de multiplicité un fort 
   \cite{piat71,shalika,jacquet-shalika-euler}. 
   Grâce aux explications  de Erez Lapid, nous savons que  pour $G=SL_{r}$  la décomposition \eqref{dec-param-cor-non-deploye}  est aussi définie sur  $\ov\Q$  
 et indépendante de 
  $\ell$ et  du plongement $\ov\Q \hookrightarrow \Qlbar$. En effet la \propref{functoriality-L-groups-H-sigma} appliquée à  l'inclusion $\Upsilon:SL_{r}\hookrightarrow GL_{r}$ montre que la décomposition \eqref{dec-param-cor-non-deploye} pour $SL_{r}$ est déterminée par la  décomposition pour $GL_{r}$, puisque   $\beta_{\Upsilon}^{*}$ est surjective \cite{hiraga-saito,labesse-langlands}. D'après  \cite{larsen1,larsen2} et des variantes de  l'argument précédent on peut obtenir un résultat analogue pour d'autres groupes. 
   \end{rem}

\begin{rem}  Si on admet la conjecture de Tate en plus des conjectures standard, on pourrait espérer démontrer la conjecture précédente en rendant motiviques toutes les constructions de cet article. Cependant cela serait très difficile car on a utilisé la cohomologie à support compact du faisceau d'intersection de champs de Deligne-Mumford qui ne sont pas propres. La conjecture de Tate permettrait d'obtenir  
 à l'aide du  lemme de Drinfeld  (le \corref{equiv-cat-U-I-Omega}) 
des représentations $\ov\Q$-linéaires de 
$G_{\mr{mot}}^{W}(F)^{I}$, où $G_{\mr{mot}}(F)$ est le groupe de Galois motivique de $F$ (associé au choix d'un foncteur fibre sur $\ov\Q$ dont l'existence résulte de \cite{deligne-tens-fest}), et 
$G_{\mr{mot}}^{W}(F)=G_{\mr{mot}}(F)\times_{\wh\Z}\Z$ pourrait être appelé   
``groupe de Weil motivique''. On utiliserait alors la variante évidente de  la \propref{Xi-n} où l'on  remplace $\Gamma$ par $G_{\mr{mot}}^{W}(F)$. 
\end{rem} 

  \begin{rem}   La conjecture \ref{conj-algebre-B-Qbar} impliquerait l'existence d'une algèbre $\mc B_{\ov \Q}$ telle que $\mc B\otimes_{E} \Qlbar=\mc B_{\ov \Q}\otimes_{\ov \Q}\Qlbar$ pour tout plongement $\ov\Q\hookrightarrow \ov \Q_{\ell}$. Grâce à la remarque   \ref{creation-annihilation-dualite} et  à l'égalité entre \eqref{S-bil-carac} et \eqref{S-bil-carac2},  le transposé d'un opérateur d'excursion  est un opérateur d'excursion, par conséquent $\mc B_{\ov \Q}\otimes_{\ov \Q}\C$  serait alors une $C^{*}$-algèbre  commutative et 
l'algèbre $\mc B_{\ov \Q}$ serait donc  réduite. 
 \end{rem}

 \subsubsection{Digression sur le cas des  corps de nombres}
 Il paraît totalement hors d'atteinte d'appliquer les méthodes de cet article aux  corps de nombres. Cependant on peut se demander s'il est raisonnable d'espérer une décomposition analogue à la  décomposition  {\it canonique} 
  \eqref{dec-param-cor-non-deploye}, une formule de multiplicités d'Arthur pour chacun des espaces  $\mf H_{\sigma}$, et une conjecture analogue à la conjecture 
\ref{conj-algebre-B-Qbar}. 

Quand  $F$ est un corps de fonctions comme dans cet article, 
la limite   $ \varprojlim_{N } \Bun_{G,N}(\Fq)$ est égale à  $\big(G(\ov F)\backslash G(\mb A\otimes_{F}\ov F)   \big)^{\on{Gal}(\ov F/F)}$. Or cette dernière expression garde un sens 
pour les corps de nombres et pour éviter des problèmes topologiques on peut remarquer qu'elle est aussi égale à  
    $\big(G(\check F)\backslash G(\mb A\otimes_{F}\check F)\big)^{\on{Gal}(\check F/F)}$ où  $\check F$ est une  extension finie galoisienne de $F$ sur laquelle  $G$ est déployé. 
 
 On peut donc espérer que si   $F$ est un corps de nombres  et si 
$\Xi$ est un réseau  dans  $Z(F)\backslash Z(\mb  A)$, 
la partie discrète de  $L^{2}\Big(\big(G(\ov F)\backslash G(
\mb A\otimes_{F}\ov F)\big)^{\on{Gal}(\ov F/F)}/\Xi,\C\Big)$ admet une décomposition   {\it canonique}   indexée par les classes de conjugaison par  $\wh G(\C)$ de 
paramètres d'Arthur elliptiques. 
 Comme dans le cas des corps de  fonctions, on a 
$$\big(G(\ov F)\backslash G(\mb A\otimes_{F}\ov F)\big)^{\on{Gal}(\ov F/F)}
= \bigcup_{\alpha\in \ker^{1}(F,G)}G_{\alpha}(F)\backslash G_{\alpha}(\mb A)
$$
où $\ker^{1}(F,G)$ est fini et $G_{\alpha}$ est une forme intérieure de $G$. Le cas particulier qui ressemble le plus au cas des corps de fonctions est celui des formes automorphes cohomologiques. 
En effet la partie cohomologique  de $L^{2}_{\mr{disc}}\Big(\big(G(\ov F)\backslash G(
\mb A\otimes_{F}\ov F)\big)^{\on{Gal}(\ov F/F)}/\Xi,\C\Big)$ est définie sur  $\ov\Q$ et on peut se demander si elle admet une  décomposition {\it canonique}   sur  $\ov \Q$ indexée par les  classes d'équivalence de paramètres  d'Arthur  elliptiques 
(pour rendre cela précis il faut utiliser  \cite{buzzard-gee}).

    \subsection{Démonstration du \thmref{dec-param-cor-thm-non-deploye}}
 On va expliquer les ingrédients nouveaux nécessaires pour la preuve du   \thmref{dec-param-cor-thm-non-deploye}
 comparée à la preuve du  \thmref{dec-param-cor-thm} qui a fait l'objet des chapitres   \ref{rappels-Hecke-Gr-satake} à \ref{para-dec-param-Langlands}  de cet article.

    Si  ${}^{L}G=\wh G$,  les modifications sont très limitées. 
       Mais pour traiter le cas  où ${}^{L}G\neq \wh G$ on aura besoin d'une variante tordue  sur la courbe de l'équivalence de  Satake géométrique.   On suppose  $E$  assez grand  pour que toutes les  représentations irréductibles de ${}^{L }G$ soient définies sur $E$. 
    
       \subsubsection{Rappels sur les  grassmanniennes affines et l'équivalence de Satake  géométrique}\label{reminder-satake-twisted}
   De la même fa\c con que dans le  chapitre \ref{rappels-Hecke-Gr-satake}, on définit l'ind-schéma
    $\mr{Gr}_{I} ^{(I_{1},...,I_{k})}$ sur $U^{I}$ (on rappelle que $G$ est réductif sur $U$). On a $X\sm \wh N\subset U$  et en fait dans la suite on utilisera uniquement la restriction de $\mr{Gr}_{I} ^{(I_{1},...,I_{k})}$
     à  $(X\sm \wh N)^{I}$, que l'on notera de la même fa\c con. 
    Soit $\wh U$ un revêtement  de $U$ (et même  de $\wt U$) tel que l'image inverse  de $G$ à $\wh U$ soit déployée. Soit $W$ une représentation $E$-linéaire de 
    $({}^{L }G)^{I}$. 
Alors $\mr{Gr}_{I,W} ^{(I_{1},...,I_{k})}$ est défini  comme le   sous-schéma fermé  de  
$\mr{Gr}_{I} ^{(I_{1},...,I_{k})}$ formé des   points 
dont l'image inverse à  $(\wh U)^{I}$ appartient à 
$\mr{Gr}_{\wh U,I,\restr{W}{(\wh G)^I}} ^{(I_{1},...,I_{k})}$. 
La définition des grassmanniennes affines de Beilinson-Drinfeld  garde un sens pour les courbes ouvertes et la notation précédente signifie 
que l'on a appliqué cette définition à la courbe  $\wh U$. On a pu restreindre  $W$ à $(\wh G)^I$ parce que l'image inverse  de $G$ sur $\wh U$ est déployée. On rappelle que par définition ce  schéma est la réunion des  schémas $\mr{Gr}_{\wh U,I,V} ^{(I_{1},...,I_{k})}$ pour tous les constituants irréductibles   $V$ de $\restr{W}{(\wh G)^I}$. 
On définit alors $\Hecke_{I,W} ^{(I_{1},...,I_{k})}$ comme  l'image inverse  de 
$\mr{Gr}_{I,W} ^{(I_{1},...,I_{k})}/G_{\sum \infty x_{i}}$ par le morphisme naturel 
$\Hecke_{I} ^{(I_{1},...,I_{k})}\to \mr{Gr}_{I} ^{(I_{1},...,I_{k})}/G_{\sum \infty  x_{i}}$. 
Toutes les propriétés restent vraies, et en particulier on possède l'extension suivante du \thmref{thm-geom-satake}  au cas  non déployé  (la structure du $L$-groupe comme produit semi-direct pour l'action respectant un épinglage vient du fait qu'un épinglage de $\wh G$ apparaît naturellement dans la preuve de l'équivalence de Satake géométrique, voir \cite{hitchin}).   

Contrairement au cas déployé, les faisceaux pervers que nous allons considérer ne sont pas forcément des faisceaux intersection complète, mais sont issus du théorème suivant, où  $\mr{Gr}_{I,W} ^{(I_{1},...,I_{k})}$ joue simplement le rôle de support. 
C'est pourquoi nous avons supposé ci-dessus que $W$ est une représentation de 
$({}^{L}G)^{I}$ (et non pas seulement de $(\wh G)^{I}$) et n'avons pas cherché à  introduire de corps réflex.  

\begin{thm}\label{thm-geom-satake-non-deploye} (variante tordue de l'équivalence de Satake géométrique  \cite{zhu,richarz}) 
On a un  foncteur canonique 
$$W\mapsto \mc S_{I,W,E}^{(I_{1},...,I_{k})}$$
de la catégorie des  représentations $E$-linéaires de dimension finie de $({}^{L}G)^{I}$ vers la  catégorie des    $E$-faisceaux pervers  
$G^{\mr{ad}}_{\sum \infty x_{i}}$-équivariants sur $\mr{Gr}_{I }^{(I_{1},...,I_{k})}$. De plus $\mc S_{I,W,E}^{(I_{1},...,I_{k})}$ est supporté par 
$\mr{Gr}_{I,W}^{(I_{1},...,I_{k})}$ et on peut donc le considérer comme un faisceau pervers 
   (à un décalage près)   sur  $\mr{Gr}_{I,W}^{(I_{1},...,I_{k})}/G^{\mr{ad}}_{\sum n_{i}x_{i}}$  (où les  entiers $n_{i}$ sont assez grands  et le décalage est déterminé par la condition que  l'image inverse sur  $\mr{Gr}_{I,W}^{(I_{1},...,I_{k})}$ est perverse relativement à  $(X\sm  \wh N)^{I}$). Les  $E$-faisceaux pervers  (à un décalage près) $\mc S_{I,W,E}^{(I_{1},...,I_{k})}$ sont universellement localement acycliques relativement au morphisme vers $(X\sm  \wh N)^{I}$. Ils vérifient les mêmes propriétés b), c), d) que dans    le \thmref{thm-geom-satake}. Si $W$ se factorise à travers 
$({}^{L}G/\wh G)^{I}=\on{Gal}(\wt F/F)^{I}$, $\mc S_{I,W,E}^{(I_{1},...,I_{k})}$ est supporté sur la section  évidente   $(X\sm  \wh N)^{I}\to  \mr{Gr}_{I}^{(I_{1},...,I_{k})}$ 
(dont  l'image est le  sous-schéma où tous  les  $\phi_{j}$ sont des isomorphismes) 
et elle provient du  $E$-faisceau lisse sur $(X\sm  \wh N)^{I}$ correspondant à $W$
(et dont la  fibre en $\Delta(\ov\eta)$ est canoniquement égale à $W$). 
\end{thm}
 
 \begin{rem} Dans le théorème précédent on peut  remplacer $E$ par $\mc O_{E}$ (on a toujours $\mr{Gr}_{I,W}^{(I_{1},...,I_{k})}$ de type fini mais   on n'en donne plus de description explicite en fonction en $W$). Il en va de même dans la discussion qui suit. 
 \end{rem}
 
 En fait nous aurons besoin d'une petite généralisation du  théorème précédent, où certains $x_{i}$ restent des points variables   sur $X\sm  \wh N$ mais les autres sont fixés en des  places dans $|X|\sm \wh N$, et en contrepartie on demande seulement que  $W$ soit une  représentation du $L$-groupe {\it  local } (au lieu  de ${}^{L}G$) pour ces indices.

  Soit $ I$ un ensemble fini,  et  
  \begin{gather}\label{alpha-pattes-places}
  \alpha:I \to  (|X|\sm  \wh N)\cup \{\bigstar\}\end{gather}
une application  (l'idée  est que, pour $i\in I$,  le  point $x_{i}$ restera un point variable  sur  $X\sm  \wh N$ si $\alpha(i)=\bigstar$ et sera fixé en la place  $v$ si 
  $\alpha(i)=v$). On pose 
     $X(\bigstar) =X\sm  \wh N$, et  pour $v\in |X|\sm  \wh N$ on écrit    $X(v)=v$, considéré comme un   sous-schéma fermé   de $X$.
         Soit   $W$ une    représentation $E$-linéaire de dimension finie de $         \prod_{i\in I} {}^{L} G_{\alpha(i)}$, où 
     $ {}^{L}G_{\bigstar}={}^{L}G$. 
        Soit   $(I_{1},...,I_{k})$ une  partition de $I$ telle que $\alpha$ soit constant sur chaque $I_{j}$.      Alors on définit naturellement le  sous-schéma fermé  
     $ \mr{Gr}_{I,W}^{(I_{1},...,I_{k}),\alpha}$ de 
     $ \mr{Gr}_{I}^{(I_{1},...,I_{k}),\alpha}=\restr{ \mr{Gr}_{I}^{(I_{1},...,I_{k})}}{\prod_{i\in I}X(\alpha(i))}$  sur 
     $\prod_{i\in I}X(\alpha(i))$ 
         et le 
      faisceau pervers (à un décalage près)  
      \begin{gather}\label{F-alpha-Gr}\mc S_{I,W,E}^{(I_{1},...,I_{k}),\alpha} \text{ \ \   sur \ \  } \mr{Gr}_{I,W}^{(I_{1},...,I_{k}),\alpha}/G^{\mr{ad}}_{\sum n_{i}x_{i}},\end{gather} qui est fonctoriel en  $W$ et vérifie les mêmes  propriétés que dans le théorème ci-dessus, donc en particulier les propriétés analogues à b), c), d) du \thmref{thm-geom-satake}. On peut le vérifier   facilement en plongeant $W$ dans une représentation de $\prod_{i\in I} {}^{L} G_{\alpha(i)}$ qui admet un prolongement à $({}^{L} G)^{I}$.

 Le lien avec    l'isomorphisme de Satake classique \eqref{Isom-Satake-non-deploye}
 est le suivant. Soit $v$ une  place dans $|X|\sm  \wh N$. 
On prend  $I=\{1\}$ un  singleton et $\alpha$ défini par $\alpha(1)=v$. 
On écrit  $\mr{Gr}_{v}$ pour $\restr{\mr{Gr}_{\{1\}}^{(\{1\})}}{v}$
qui est la fibre de  la  grassmannienne affine sur $v$. 
  L'ensemble des  $\Fq$-points de cet ind-schéma s'identifie à 
   $ G(F_{v})/G(\mc O_v)$. 
   Alors, 
 pour toute  représentation  irréductible  $V$ de  ${}^{L} G_{v}$,  
 la  trace de $\Frob_{\mr{Gr}_{v}/k(v)}$ sur le  faisceau pervers $\mc S_{\{1\},V,E}^{(\{1\}),\alpha}$ est égale à $(-1)^{\s{2\rho, \omega}}h_{V,v}$, où   $h_{V,v}\in C_{c}(G(\mc O_v)\backslash G(F_{v})/G(\mc O_v), \mc O_{E})$ a été introduite juste après \eqref{Isom-Satake-non-deploye} et $\omega$ est le plus haut poids de n'importe quel constituant irréductible de $\restr{V}{\wh G}$ (de sorte que  
 $\s{2\rho, \omega}$ est la dimension de la $G(\mc O)$-orbite correspondante dans $\mr{Gr}_{v}$). 

  \subsubsection{Chtoucas}\label{chtoucas-non-split}
On construit maintenant le champ classifiant  des chtoucas dans le cas des groupes  non nécessairement déployés. Cela est facile et parallèle au cas déployé parce que les  pattes restent dans  $X\sm  \wh N\subset U$ sur lequel  $G$ est réductif. 
Une  construction sans cette hypothèse est proposée dans \cite{hartl}, mais n'est pas nécessaire ici. 
Soit   $I$ un ensemble fini, $k\in \N$ et  $(I_{1}, ..., I_{k})$ une  partition de $I$. On note
  $$\Cht_{N,I}^{(I_{1},...,I_{k})} \text{({\it resp.} } 
  \Cht_{N,I,W}^{(I_{1},...,I_{k})}, \  \Cht_{N,I,W}^{(I_{1},...,I_{k}),\leq\mu})$$ le champ sur $(X\sm  \wh N)^{I}$ dont les  points sur un  schéma $S$ classifient   un $S$-point 
  \begin{gather}\label{donnee-Hecke-non-split}\big( (x_i)_{i\in I}, (\mc G_{0}, \psi_{0}) \xrightarrow{\phi_{1}}  (\mc G_{1}, \psi_{1}) \xrightarrow{\phi_{2}}
\cdots\xrightarrow{\phi_{k-1}}  (\mc G_{k-1}, \psi_{k-1}) \xrightarrow{ \phi_{k}}    (\mc G_{k}, \psi_{k})
\big)
\end{gather}
  dans  $$\Hecke_{N,I}^{(I_{1},...,I_{k})}  \text{({\it resp.} } 
\Hecke_{N,I,W}^{(I_{1},...,I_{k})}, \ \Hecke_{N,I,W}^{(I_{1},...,I_{k}),\leq\mu})$$
 plus un  isomorphisme $ \sigma: \ta \mc{G}_{0}\isom\mc{G}_{k}$, préservant les   structures de niveau, c'est-à-dire   vérifiant 
 $\psi_{k}\circ \restr{\sigma}{N\times S}=\ta \psi_{0}$. 
 
 Comme les fibres   de $G$ sont géométriquement connexes, on peut appliquer le  théorème de Lang à la     restriction à la Weil de $G$ 
 de n'importe quel sous-schéma fini  de $X$ à  $\Fq$. Donc  pour tout niveau  $N'\supset N$,   $\Cht_{N',I,V}^{(I_{1},...,I_{k}),\leq\mu}$ est un revêtement fini étale du champ 
     $\restr{\Cht_{N,I,V}^{(I_{1},...,I_{k}),\leq\mu}}{(X\sm  \wh{N'})^{I}}$  de groupe de Galois  
  $K_{N}/K_{N'}$.  
 
 \begin{rem}
 On peut montrer que, comme les pattes varient dans $X\sm \wh N$,   $\Cht_{N,I,V}^{(I_{1},...,I_{k})}$ ne dépend en fait que de la donnée de $G$ groupe réductif sur $X\sm \wh  N$ et du sous-groupe compact ouvert $\Ker(\prod_{v\in \wh N}G(\mc O_{v})\to G(\mc O_{N}))$ de $\prod_{v\in \wh  N}G(F_{v})$. 
 \end{rem}
 
 La lemme suivant a bénéficié de discussions avec Alain Genestier, Urs Hartl et Jochen Heinloth. 
 On renvoie au théorème 3.14 et à la remarque 3.18 de \cite{hartl} pour plus de détails. 
 
 \begin{lem}
 Le quotient $\Cht_{N,I,V}^{(I_{1},...,I_{k}),\leq\mu}/\Xi$ est un champ de Deligne-Mumford de type fini. 
  \end{lem}
 \dem 
 On rappelle que $\mu$ est un copoids dominant de $G^{\mr{ad}}$. 
 Il suffit de montrer l'énoncé pour $(I_{1},...,I_{k})=(I)$, ce qu'on suppose désormais. Par ailleurs il suffit de montrer le résultat avec $N$ assez grand (en fonction de $\mu$). 
 
 a) Cas où $G$ est adjoint. On a rappelé au paragraphe \ref{para-non-deploye-statement}  que, d'après \cite{hartl}, 
 $\Bun_{G}^{\leq \mu}$ est un champ d'Artin de type fini, et pour $N$ assez grand   $\Bun_{G,N}^{\leq \mu}$, $\Bun_{G,N}^{\leq \mu+\kappa}$ et 
 $\Hecke_{N,I,V}^{(I),\leq\mu}$ sont des schémas de type fini, donc 
  $\Cht_{N,I,V}^{(I),\leq\mu}$ est un sous-schéma fermé de 
 $\Hecke_{N,I,V}^{(I),\leq\mu}$ (image inverse du graphe de Frobenius de 
 $\Bun_{G,N}^{\leq \mu+\kappa}$) et il est donc un schéma de  type fini.

 b) Cas où $G$ est un tore $T$. Alors il existe un revêtement fini $Y$ de $X\sm  \wh N$ tel que   $\Cht_{N,I,V}^{(I)}\times_{(X\sm  \wh N)^{I}}Y^{I}$ soit une réunion finie de $\Bun_{T,N}(\Fq)$-torseurs (cette réunion finie a pour cardinal 
 le nombre de caractères différents de $\wh T$ dans la restriction de $V$ à $\wh T\subset {}^{L}T$). Or le morphisme naturel    $\Xi\to \Bun_{T,N}(\Fq)$ a un noyau fini et une image  d'indice fini, donc
 (en supposant, quitte à diminuer $\Xi$,  que le noyau soit nul), 
 $\Cht_{N,I,V}^{(I)}/\Xi$ est un revêtement fini de $(X\sm  \wh N)^{I}$. 
 Plus précisément on montre qu'il existe $U\subset  \Cht_{N,I,V}^{(I)}$ ouvert fermé de type fini tel que la réunion de ses translatés par $\Xi$ recouvre 
 $\Cht_{N,I,V}^{(I)}$. 
 
 c) Cas général. On considère le tore $G^{\mr{ab}}=G/G^{\mr{der}}$ abélianisé de $G$. Alors le morphisme $(\beta^{\mr{ad}}, \beta^{\mr{ab}}):G\to G^{\mr{ad}}\times G^{\mr{ab}}$ est surjectif sur $\ov F$ et a pour noyau un groupe abélien fini $K$. 
 
 On choisit les modèles parahoriques lisses de $G$, $G^{\mr{ad}}$ et $G^{\mr{ab}}$ pour que le morphisme $G \to G^{\mr{ad}}\times G^{\mr{ab}}$ s'étende à    $X$ tout entier.

 On note encore $V$ la représentation de ${}^{L}G^{\mr{ad}}\times {}^{L}G^{\mr{ab}}$ (et donc de chacun des deux facteurs) qui se déduit de $V$ en composant par 
 $({}^{L}\beta^{\mr{ad}}, {}^{L}\beta^{\mr{ab}}): {}^{L}G^{\mr{ad}}\times {}^{L}G^{\mr{ad}}\to{}^{L}G$. D'où un morphisme 
 $$(\beta^{\mr{ad}}_{*}, \beta^{\mr{ab}}_{*}):\Cht_{G,N,I,V}^{(I)}\to 
 \Cht_{G^{\mr{ad}},N,I,V}^{(I)}\times \Cht_{G^{\mr{ab}},N,I,V}^{(I)}.$$

 On applique b) à $G^{\mr{ab}}$ et on note $U$ un ouvert fermé de type fini
 de  $\Cht_{G^{\mr{ab}}, N,I,V}^{(I)}$
  tel que la réunion de ses translatés par $\Xi$ recouvre 
 $\Cht_{G^{\mr{ab}},N,I,V}^{(I)}$. La réunion des translatés de 
 $(\beta^{\mr{ab}}_{*})^{-1}(U)$ 
 par $\Xi$ recouvre $\Cht_{G,N,I,V}^{(I)}$. Il suffit donc de montrer que 
  $(\beta^{\mr{ab}}_{*})^{-1}(U)$ est de type fini. 

 Soit, comme dans le \lemref{lem-behrend-hartl-heinloth}, 
un fibré vectoriel $\mc V$ de rang $r$ sur $X$ muni d'une trivialisation de  $\det(\mc V)$ et une représentation fidèle $\rho:G\to SL(\mc V)$ telle que les quotients $SL(\mc V)/G$ (et donc $GL(\mc V)/G$) soient quasi-affines. 
En prenant $n=\sharp K$ on voit que la représentation $\rho^{\otimes n}$  de $G$ sur $  {\mc V}^{\otimes n}$ se factorise par $(\beta^{\mr{ad}}, \beta^{\mr{ab}})$ au point générique de $X$. On en déduit facilement que $(\rho^{\otimes n})_{*}((\beta^{\mr{ab}}_{*})^{-1}(U))$ est inclus dans $\Bun_{GL_{r^{n}}}^{0,\leq \nu}$ pour $\nu$ assez grand. 
Par conséquent  $\rho _{*}((\beta^{\mr{ab}}_{*})^{-1}(U))$ est inclus 
 dans $\Bun_{GL_{r}}^{0,\leq \nu}$ pour $\nu$ assez grand. 
 Le même argument que dans a) montre alors que $(\beta^{\mr{ab}}_{*})^{-1}(U)$ est de type fini. 
 \cqfd
 
    Le champ discret $\Cht_{N,\emptyset,\mbf 1}^{(\emptyset),\leq\mu}$
 est constant sur $\Fq$ et égal à   $\Bun_{G,N}^{\leq \mu}(\Fq)$. L'argument pour le démontrer est le même que  dans la preuve de la proposition 2.16 c): 
 on applique le  lemme  3.3 b) de \cite{var}
      à $\Bun_{G,N'}^{\leq \mu}$, où $N'\supset N$ est assez grand  (en  fonction  de $\mu$) pour que  toute composante connexe  de $\Bun^{\leq \mu}_{G,N'}$ soit  un  schéma. 

 Si  les  entiers $n_{i}$ sont assez grands, le morphisme naturel 
 $$\epsilon_{N,(I),W,\underline{n}}^{(I_{1},...,I_{k}), \Xi}: 
\Cht_{N,I,W} ^{(I_{1},...,I_{k})}/\Xi
\to \mr{Gr}_{I,W}^{(I_{1},...,I_{k})}/G^{\mr{ad}}_{\sum _{i\in I}n_{i}x_{i}}$$ est lisse.

Comme dans la  \defiref{defi-F-E-Cht}
on définit  $\mc F_{N,I,W, \Xi,E}^{(I_{1},...,I_{k})}$ comme le  faisceau pervers (à un décalage près) sur $\Cht_{N,I,W}^{(I_{1},...,I_{k})}/\Xi$ égal à 
\begin{gather}\label{def-IC-V-nonsplit}\mc F_{N,I,W, \Xi,E}^{(I_{1},...,I_{k})}=
\big(\epsilon_{N,(I),W,\underline{n}}^{(I_{1},...,I_{k}), \Xi}\big)^{*}
\big(\mc S_{I,W,E}^{(I_{1},...,I_{k})}\big).\end{gather}
On définit 
 \begin{gather} \label{defi-H-alpha-Cht}
 \mc H _{ N,I,W}^{\leq\mu,E}=  R
 \big(\mf p_{N,I}^{(I_{1},...,I_{k}),\leq\mu}\big)_{!}\Big(\restr{\mc F_{N,I,W,\Xi,E}^{(I_{1},...,I_{k})}} {\Cht_{N,I,W}^{(I_{1},...,I_{k}),\leq\mu}/\Xi}\Big)  
  \end{gather}
   qui appartient à $D^{b}_{c}((X\sm  \wh N)^{I}, E)$ 
   (grâce au lemme précédent) et ne dépend pas du choix de la   partition $(I_{1},...,I_{k})$.  
  Il  dépend bien sûr de $\Xi$, que l'on omet pour raccourcir les formules. 
   En particulier on possède   le faisceau  constructible  
       $\mc H _{ N,I,W}^{0,\leq\mu,E}$ sur  $(X\sm  \wh N)^{I}$. 

\subsubsection{Une petite généralisation, et les opérateurs $S_{V,v}$} 
\label{mild-gen-SVv}
Pour construire les  opérateurs analogues aux $S_{V,v}$, on a besoin d'une petite  généralisation des champs de chtoucas,  associée à la donnée d'une  application 
 $ \alpha$ comme dans  \eqref{alpha-pattes-places}.   Soit   $(I_{1},...,I_{k})$ une partition de $I$ telle que $\alpha$ soit constant sur chaque  $I_{j}$.   

  Soit   $W$ une  représentation  $E$-linéaire de $\prod_{i\in I} {}^{L} G_{\alpha(i)}$. 
On définit  $\Hecke_{N,I,W}^{(I_{1},...,I_{k}),\alpha}$ de la même fa\c con que 
 $\on{Gr}_{I,W}^{(I_{1},...,I_{k}),\alpha}$ introduit ci-dessus. 
On note
  $\Cht_{N,I,W}^{(I_{1},...,I_{k}),\alpha}$ le champ sur $\prod_{i\in I}X(\alpha(i))$  classifiant un  point \eqref{donnee-Hecke-non-split} dans 
  $\Hecke_{N,I,W,}^{(I_{1},...,I_{k}),\alpha}$
 plus un  isomorphisme $ \sigma: \ta \mc{G}_{0}\isom\mc{G}_{k}$, préservant les    structures de niveau, c'est-à-dire   vérifiant 
 $\psi_{k}\circ \restr{\sigma}{N\times S}=\ta \psi_{0}$. Si les entiers $n_{i}$ sont assez grands on a un morphisme naturel $\Cht_{N,I,W}^{(I_{1},...,I_{k}),\alpha}/\Xi\to 
 \on{Gr}_{I,W}^{(I_{1},...,I_{k}),\alpha}/G^{\mr{ad}}_{\sum n_{i}x_{i}}$ qui est lisse. 
 
 Comme dans  \eqref{def-IC-V-nonsplit} on définit 
  $\mc F_{N,I,W, \Xi,E}^{(I_{1},...,I_{k}),\alpha}$ comme le  faisceau pervers (à un décalage près) sur $\Cht_{N,I,W}^{(I_{1},...,I_{k}),\alpha}/\Xi$ égal à 
 l'image inverse du  faisceau pervers $\mc S_{I,W,E}^{(I_{1},...,I_{k}),\alpha} $ introduit dans    \eqref{F-alpha-Gr}. Comme dans  \eqref{defi-H-alpha-Cht}  on pose  
 \begin{gather} \label{defi-HNIW-non-deploye}
 \mc H _{ N,I,\alpha,W}^{\leq\mu,E}=  R
 \big(\mf p_{N,I}^{(I_{1},...,I_{k}),\alpha,\leq\mu}\big)_{!}\Big(\restr{\mc F_{N,I,W,\Xi,E}^{(I_{1},...,I_{k}),\alpha}} {\Cht_{N,I,W}^{(I_{1},...,I_{k}),\alpha,\leq\mu}/\Xi}\Big)  
  \end{gather} dans $D^{b}_{c}(\prod_{i\in I}X(\alpha(i)), E)$. 
   Le fait que 
\eqref{defi-HNIW-non-deploye} 
est indépendant du choix de la partition $(I_{1},...,I_{k})$
vient de ce que 
$\mc S_{I,W,E}^{(I_{1},...,I_{k}),\alpha}$ 
défini dans 
\eqref{F-alpha-Gr}   vérifie l'analogue de  b)  du \thmref{thm-geom-satake}
              
   Soit  de plus $v\in |X|\sm \wh N$ et  $V$ une   représentation de 
     ${}^{L} G_{v}$. On prolonge  $\alpha$ en      $\alpha' : I\cup  \{1,2\}\to (|X|\sm \wh N)\cup \{\bigstar\}$ en posant  $\alpha'(1)=\alpha'(2)=v$. On définit    alors    $$S_{V,v}\in \Hom_{D^{b}_{c}(\prod_{i\in I}X(\alpha(i)), E)}\Big(
   \mc H _{ N,I,\alpha,W}^{\leq\mu,E},
 \mc H _{ N,I,\alpha,W}^{\leq\mu+\kappa,E}
  \Big)
$$ comme la  composée  
       \begin{gather}\label{compo-SVv-non-deploye}
 \mc H _{ N,I,\alpha,W}^{\leq\mu,E} \boxtimes E_{v} 
  \xrightarrow{ \mc C_{  \delta_{V},v}^{\sharp}}
 \restr{ \mc H _{ N,I\cup\{1,2\},\alpha',W\boxtimes V\boxtimes V^{*}}^{\leq\mu,E}}{\prod_{i\in I}X(\alpha(i))\times \Delta(v)} \\ \nonumber 
 \xrightarrow{ (F_{\{1\}})^{\deg(v)} }
 \restr{ \mc H _{ N,I\cup\{1,2\},\alpha',W\boxtimes V\boxtimes V^{*}}^{\leq\mu+\kappa,E}}{\prod_{i\in I}X(\alpha(i))\times \Delta(v)} 
  \xrightarrow{\mc C_{
    \on{ev}_{V},v}^{\flat}}
   \mc H _{ N,I,\alpha,W}^{\leq\mu+\kappa,E} \boxtimes E_{v} . 
 \end{gather}
 où les morphismes de création et d'annihilation $ \mc C_{  \delta_{V},v}^{\sharp}$ et $\mc C_{ \on{ev}_{V},v}^{\flat}$ sont évidents
 (la notation $(F_{\{1\}})^{\deg(v)}$ ci-dessus est un peu abusive dans la mesure où $F_{\{1\}}$ n'existe pas dans ce cadre). Plus précisément la  composée  \eqref{compo-SVv-non-deploye} commute avec l'action du  morphisme de  Frobenius partiel sur $E_{v}$, et par descente on obtient $S_{V,v}$.

La même preuve que celle  de  la \propref{prop-coal-frob-cas-part} donne l'égalité  \begin{gather}\label{T-S-non-deploye}
T(h_{V,v})= \restr{S_{V,v}}{\prod_{i\in I}\big(X(\alpha(i))\cap (X\sm v)\big) }. 
\end{gather}
   Comme dans 
 le \corref{cor-hecke-etendus-compo}, on utilise les  morphismes  $S_{V,v}$ pour prolonger l'action de  $C_{c}(K_{N}\backslash G(\mb A)/K_{N})$ sur   $\mc H _{ N,I,\alpha,W}^{\leq\mu,E} $ à    $\prod_{i\in I}X(\alpha(i))$ tout entier. 
 De la même fa\c con que dans le chapitre \ref{para-Relations d'Eichler-Shimura} on montre les relations d'Eichler-Shimura  
  (elles existent en toutes les places de $X\sm \wh N$ mais on n'en aura besoin qu'aux  places où $G$ est déployé, et dans ce cas il n'y a aucune  différence avec le  
  chapitre \ref{para-Relations d'Eichler-Shimura}).

    \subsubsection{Opérateurs d'excursion}\label{excur-op-non-split}

Pour tout point géométrique   $\ov x$  dans $(X\sm \wh N)^{I}$, 
on définit le $E$-espace vectoriel des  éléments  Hecke-finis   $$\Big(\varinjlim_{\mu }\restr{\mc H _{ N,I,W}^{0,\leq\mu,E}}{\ov x} \Big)^{\mr{Hf}}\subset \varinjlim_{\mu }\restr{\mc H _{ N,I,W}^{0,\leq\mu,E}}{\ov x} $$ de la même fa\c con que dans la   \defiref{hecke-fini}.

         Dans le cas   où $I$ est vide,
               \begin{gather*}\varinjlim_{\mu }\restr{\mc H _{ N,\{0\},\mbf 1}^{0,\leq\mu,E}}{\ov\eta}=\varinjlim_{\mu }\restr{\mc H _{ N,\emptyset,\mbf 1}^{0,\leq\mu,E}}{\Fqbar}=
      C_{c}(\Bun_{G,N}(\Fq)/\Xi,E)\end{gather*} et 
      la \propref{prop-cusp-hecke-finies}  implique que 
   \begin{gather}\label{Cccusp-Hf-non-deploye}\Big(\varinjlim_{\mu }\restr{\mc H _{ N,\{0\},\mbf 1}^{0,\leq\mu,E}}{\ov\eta}\Big)^{\mr{Hf}}=
       \Big(\varinjlim_{\mu }\restr{\mc H _{ N,\emptyset,\mbf 1}^{0,\leq\mu,E}}{\Fqbar}\Big)^{\mr{Hf}}=C_{c}^{\mr{cusp}}(\Bun_{G,N}(\Fq)/\Xi,E). \end{gather}
     On choisit une flèche de spécialisation  $\on{\mf{sp}}$ de $\ov{\eta^{I}}$ vers $\Delta(\ov \eta)$.   Les mêmes  arguments que dans  la \propref{prop-action-Hf} et le \corref{bijectivite-Hecke-fini} montrent que   
     $$\on{\mf{sp}}^{*}: \Big( \varinjlim _{\mu}\restr{\mc H _{N, I, W}^{0,\leq\mu,E}}{\Delta(\ov{\eta})} \Big)^{\mr{Hf}}\to 
 \Big( \varinjlim _{\mu}\restr{\mc H _{N, I, W}^{0,\leq\mu,E}}{\ov{\eta^{I}}}\Big)^{\mr{Hf}} $$
 est une bijection et que cet espace est muni 
      d'une action de 
       $\pi_{1}(\eta,\ov\eta)^{I}=\on{Gal}(\ov F/F)^{I}$. 
  Pour la partie de l'argument qui utilise les  relations  d'Eichler-Shimura,  on peut choisir  (si on le souhaite) les  places $v_{i}$ telles que $G$ y soit déployé. 
        
Pour alléger les  notations, nous ne définirons les  morphismes de création et d'annihilation que sur  les  faisceaux  $\mc H _{ N,I,W}^{0,\leq\mu,E}$
      (au lieu du cadre plus général des  faisceaux  $\mc H _{ N,I,\alpha,W}^{0,\leq\mu,E}$ qui contiendrait comme cas particulier les  morphismes $\mc C_{  \delta_{V},v}^{\sharp}$ et  $\mc C_{
    \on{ev}_{V},v}^{\flat}$ de \eqref{compo-SVv-non-deploye}). Soit $J$ un ensemble fini et  $U$ une  représentation de $({}^{L }G)^{J}$. On rappelle la notation $\zeta_J:J\to \{0\}$, de sorte que   $U^{\zeta_J}$ est  la   représentation diagonale de ${}^{L }G$. Le sous-espace  $(U^{\zeta_J})^{\wh G}$, 
    formé des  vecteurs  de  $U$ qui sont invariants sous l'action diagonale  de $\wh G$,       est une   représentation de         $\on{Gal}(\wt F/F)$ (par l'action diagonale). 
    On note 
  $\iota $ l'inclusion $(U^{\zeta_J})^{\wh G}\subset U^{\zeta_J}$. 
     On note 
      $\restr{ (U^{\zeta_J})^{\wh G}}{X\sm  \wh N}$ le  $E$-faisceau lisse sur $X\sm  \wh N$ qui est associé à la  représentation de 
      $\on{Gal}(\wt F/F)$ sur $ (U^{\zeta_J})^{\wh G}$, 
          et dont la fibre en  $\ov \eta$ est  canoniquement égale à $(U^{\zeta_J})^{\wh G}$. 
   Le  morphisme de création
 est  la  composée 
     \begin{gather}\label{def-creation-non-deploye}
  \mc H _{ N,I,W}^{0,\leq\mu,E} \boxtimes \restr{ (U^{\zeta_J})^{\wh G}}{X\sm  \wh N}
  \xrightarrow{\sim}
 \mc H _{ N,I\cup\{0\},W\boxtimes \mbf (U^{\zeta_J})^{\wh G}}^{0,\leq\mu,E}  \\ \nonumber 
 \xrightarrow{\mc H(\Id_{W}\boxtimes \,  \iota) }
 \mc H _{ N,I\cup\{0\},W\boxtimes U^{\zeta_{J}}}^{0,\leq\mu,E} 
 \xrightarrow{\sim}
 \restr{ \mc H _{ N,I\cup J,W\boxtimes U}^{0,\leq\mu,E}}{(X\sm  \wh N)^{I}\times \Delta(X\sm  \wh N)}.  
   \end{gather}
On remarque que le vecteur $x$, qui apparaissait  dans la  notation  $\mc C_{  x}^{\sharp }$ dans le cas  déployé, ne peut être introduit ici puisque le faisceau 
    $\restr{ (U^{\zeta_J})^{\wh G}}{X\sm  \wh N}$ n'est pas  nécessairement constant.  Il réapparaîtra  dans \eqref{creation-Ivide-x} quand on restreindra ce faisceau à  $\ov\eta$. 

         De la même fa\c con, en notant  $\nu : U^{\zeta_J}\to (U^{\zeta_J})_{\wh G} $  le quotient  maximal sur lequel l'action diagonale de $\wh G$ est triviale, on définit 
          le  morphisme d'annihilation comme  la  composée
           \begin{gather}\label{def-annihil-non-deploye}
            \restr{ \mc H _{ N,I\cup J,W\boxtimes U}^{0,\leq\mu,E}}{(X\sm  \wh N)^{I}\times \Delta(X\sm  \wh N)}
            \xrightarrow{\sim}
        \mc H _{ N,I\cup\{0\},W\boxtimes U^{\zeta_{J}}}^{0,\leq\mu,E} 
  \\ \nonumber 
       \xrightarrow{\mc H(\Id_{W}\boxtimes \nu) }
   \mc H _{ N,I\cup\{0\},W\boxtimes \mbf (U^{\zeta_J})_{\wh G}}^{0,\leq\mu,E} 
  \xrightarrow{\sim}
  \mc H _{ N,I,W}^{0,\leq\mu,E} \boxtimes \restr{ (U^{\zeta_J})_{\wh G}}{X\sm  \wh N}.  
   \end{gather}

        On prend maintenant  $I=\emptyset$  et  $W=\mbf 1$ et on restreint à  $\ov\eta$ les   morphismes \eqref{def-creation-non-deploye} et 
       \eqref{def-annihil-non-deploye}. En rempla\c cant les lettres $(J,U)$ par 
$(I,W)$ (qui ne sont plus utilisées), on obtient
\begin{itemize}
\item pour tout  $x\in W$   invariant par l'action diagonale de $\wh G$, le   morphisme de création 
\begin{gather}\label{creation-Ivide-x} \restr{\mc H _{ N,\emptyset,\mbf 1}^{0,\leq\mu,E}}{\Fqbar}\to 
\restr{\mc H _{ N,I,W}^{0,\leq\mu,E}}{\Delta(\ov\eta)}\end{gather}
\item et pour tout 
 $\xi \in W^{*}$
 qui est  invariant par l'action diagonale  de $\wh G$ le    morphisme 
 d'annihilation
$$ \restr{\mc H _{ N,I,W}^{0,\leq\mu,E}}{\Delta(\ov\eta)}\to 
\restr{\mc H _{ N,\emptyset,\mbf 1}^{0,\leq\mu,E}}{\Fqbar}. 
$$
\end{itemize}

Grâce à  \eqref{Cccusp-Hf-non-deploye}, 
on construit les opérateurs d'excursion exactement comme dans la  définition-proposition  \ref{antecedent-def-prop}. Soit  $I$ un ensemble fini, $W$ une   représentation $E$-linéaire de dimension finie  de $({}^{L}G)^{I}$, $x\in W $ et  $\xi\in W^{*}$ invariants sous l'action diagonale  de $\wh G$, 
    et $(\gamma_i)_{i\in I}\in \pi_{1}(\eta,\ov\eta)^{I}$.
  Alors   $$S_{I,W,x,\xi,(\gamma_i)_{i\in I}}
\in \mr{End}_{C_{c}(K_{N}\backslash G(\mb A)/K_{N},E)}
\Big(   C_{c}^{\mr{cusp}}(\Bun_{G,N}(\Fq)/\Xi,E)\Big)$$ est défini  comme la  composée \eqref{diag-S}.

 Ces opérateurs d'excursion engendrent une  sous-algèbre commutative
     $$\mc B\subset  \on{End}_{C_{c}(K_{N}\backslash G(\mb A)/K_{N},E)}\Big(C_{c}^{\mr{cusp}}(\Bun_{G,N}(\Fq)/\Xi,E)\Big).
$$
        Pour toute fonction $f\in \mc O(\wh G_E\backslash ({}^{L}G_E)^{I}/\wh G_E)$, on pose 
        $$S_{I,f,(\gamma_i)_{i\in I}}=S_{I,W,x,\xi,(\gamma_i)_{i\in I}}$$
            où $W,x,\xi$ sont tels que 
            $f((g_{i})_{i\in I})=\s{\xi, (g_{i})_{i\in I} . x}$. 
        Ces opérateurs  vérifient des propriétés similaires à  (i), (ii), (iii) et (iv) de la \propref{art-prop-SIf-i-ii-iii}, et à la \propref{prop-harris}. L'analogue  de (v) de la \propref{art-prop-SIf-i-ii-iii}  sera considéré dans la propriété e') ci-dessous. 
       
   On a une  nouvelle propriété 
   
   (vi') Si $f$ provient d'une  fonction sur 
            $  (\on{Gal}(\wt F/F))^{I}$ (grâce à l'égalité  ${}^{L}G/\wh G=\on{Gal}(\wt F/F)$), $S_{I,f,(\gamma_i)_{i\in I}}$ est le scalaire $f((\gamma_i)_{i\in I})$
            (où les $\gamma_i$ désignent ici par abus leurs images dans $\on{Gal}(\wt F/F)$). Plus généralement si $J$ est un sous-ensemble  de $I$ et $f$ provient d'une  fonction sur $$\wh G_E\backslash ({}^{L}G_E)^{J}/\wh G_E \times (\on{Gal}(\wt F/F))^{I\sm J},$$ alors $$S_{I,f,(\gamma_i)_{i\in I}}=S_{J,\check f,(\gamma_j)_{j\in J}}\text{ où }\check f((g_{j})_{j\in J})=f((g_{j})_{j\in J},(\gamma_i)_{i\in I\sm J}).$$

            \subsubsection{Construction des paramètres de Langlands}
            \label{subsubsection-construction-Langlands-non-split}
         Pour tout entier   $n$ on note   $ \mc O(({}^{L }G_E)^{n}\modmod \wh G_E)$ la  $E$-algèbre de type fini formée par les  fonctions régulières sur  $({}^{L }G_{E})^{n}$ qui sont invariantes par  conjugaison diagonale par  $\wh G_{E}$, c'est-à-dire  telles  que  $f(hg_{1}h^{-1},...,hg_{n}h^{-1})=
        f(g_{1},...,g_{n})$ pour tout  $h\in \wh G_{E}$ et $g_{1},...,g_{n}\in {}^{L }G_{E}$. 
        
   Comme dans la   définition-proposition \ref{cor-annul-Lambda-n} on obtient pour tout $n\in \N$ un   morphisme d'algèbres 
       $$ \Theta_{n}:   \mc O(({}^{L }G_E)^{n}\modmod \wh G_E)
       \to C( \pi_{1}(X\sm  \wh N, \ov\eta)^{n},\mc B)$$   
 tel que, si on pose   $I=\{0,...,n\}$, 
   $$\Theta_{n}(f)(\gamma_1,...,\gamma_{n})=S_{I,\wt f,(1,\gamma_1,...,\gamma_{n})}$$
   où $\wt f\in \mc O(\wh G_E\backslash ({}^{L}G_E)^{I}/\wh G_E)$ est définie par 
   $\wt f(g_{0},g_{1},...,g_{n})=f(g_{0}^{-1}g_{1},..., g_{0}^{-1}g_{n})$.

    \begin{rem} A l'aide de (vi')  on peut montrer que les  opérateurs $S_{I,\wt f,(1,\gamma_1,...,\gamma_{n})}$ ci-dessus engendrent la   $E$-algèbre $\mc B$ quand  $n$, $f$ et $(\gamma_1,...,\gamma_{n})$ varient. 
    En effet $\mc O({}^{L}G_E\backslash ({}^{L}G_E)^{I}/\wh G_E)$ et 
    $\mc O(({}^{L}G_E/\wh G_E)^{I})$ 
    engendrent 
   $ \mc O(\wh G_E\backslash ({}^{L}G_E)^{I}/\wh G_E)$, car l'action de 
   ${}^{L}G_E/\wh G_E$ par translation à gauche sur $({}^{L}G_E/\wh G_E)^{I}=\pi_{0}(\wh G_E\backslash ({}^{L}G_E)^{I}/\wh G_E)$ est libre. 
    Donc en considérant la suite  $\Theta_{n}$ on ne perd aucune  information.  
       \end{rem}
       
    La suite  $\Theta_{n}$   vérifie des  propriétés similaires à  a), b), c) et d) de la définition-proposition  \ref{cor-annul-Lambda-n} ainsi qu'une  variante  e')  de la propriété  e) 
   et une nouvelle  propriété  f'): 
     \begin{itemize}
     \item []       e')    pour toute    place $v\in |X|\sm  \wh N$, pour tout plongement  $\ov F\subset \ov {F_{v}}$ (induisant des  plongements  $\on{Gal}(\ov {F_{v}}/F_{v})\hookrightarrow \on{Gal}(\ov F/F)$ et ${}^{L }G_{v,E}\subset {}^{L }G_{E}$), 
et pour toute représentation irréductible $V$ de ${}^{L }G_{v,E}$, 
  de caractère  $\chi_{V}\in \mc O({}^{L }G_{v,E} \modmod {}^{L }G_{v,E})$,  
      pour toute fonction $f\in  \mc O({}^{L }G_{E} \modmod {}^{L }G_{v,E})$ prolongeant $\chi_{V}$,  et pour tout  $\gamma \in \on{Gal}(\ov {F_{v}}/F_{v})$ avec  $\deg(\gamma)=1 $, 
   $\Theta_{1}(f)(\gamma)=T(h_{V,v})$.  
      \item []   f')
      si  $f$ est une fonction sur le  groupe fini  $\on{Gal}(\wt F/F)^{n}$ et  que l'on considère $f$ comme  un élément de  $\mc O(({}^{L }G_E)^{n}\modmod \wh G_E)$ alors 
       $\Theta_{n}(f)$ est la   composée 
       $$
   \on{Gal}(\ov F/F)^{n}\to \on{Gal}(\wt F/F)^{n} \xrightarrow{f} E \subset \mc B.$$
        \end{itemize}
        La propriété   f') est une reformulation de la première partie  de (vi') ci-dessus.  La propriété   e') résulte  de \eqref{T-S-non-deploye}, grâce à une compatibilité évidente entre les 
        morphismes de création et d'annihilation  $ \mc C_{  \delta_{V},v}^{\sharp}$ et  $\mc C_{
    \on{ev}_{V},v}^{\flat}$
de    \eqref{compo-SVv-non-deploye}
et les  restrictions à $v$ des morphismes \eqref{def-creation-non-deploye} et \eqref{def-annihil-non-deploye} (associés à $J=\{1,2\}$ et à une représentation $U$ de $({}^{L }G_{E})^{2}$ dont la restriction à $({}^{L }G_{v,E})^{2}$ contient 
$V\boxtimes V^{*}$). 
  En effet, par un argument similaire  au  \lemref{S-non-ram}, $\Theta_{1}(f)(\gamma)$ dépend seulement de  la  restriction de $f$ à 
${}^{L }G_{v,E} \modmod {}^{L }G_{v,E}$
et toute fonction dans $\mc O({}^{L }G_{v,E} \modmod {}^{L }G_{v,E})$ est une combinaison de caractères de représentations  irréductibles de ${}^{L }G_{v,E}$.

    On a une  décomposition spectrale (en espaces propres généralisés)
       \begin{gather}\label{dec-param-non-deploye} 
      C_{c}^{\mr{cusp}}(\Bun_{G,N}(\Fq)/\Xi,\Qlbar)=\oplus_{\nu} \mf H_{\nu}
       \end{gather}
       où $\nu$ parcourt les  caractères de  $\mc B$. 
       
      On applique  la \propref{Xi-n} à  $H={}^{L}G_E$ (et $H^{0}=\wh G_E$). 
   Ainsi  pour tout caractère 
 $\nu$ de $\mc B$, on obtient un   morphisme 
       $\sigma:\on{Gal}(\ov F/F)\to {}^{L}G(\Qlbar )$ bien défini à conjugaison près par   $\wh G(\Qlbar)$, et 
tel que 
\begin{itemize}
\item [] (C'1) $\sigma$ prend ses valeurs  dans  ${}^{L}G(E')$, où $E'$ est une extension finie  de $E$ (et donc de $\Ql$), et il  est continu, 
\item [] (C'2) l'adhérence de Zariski de son image est réductive, 
\item[] (C'3) 
pour    $n\in \N^{*}$, $f\in  \mc O(({}^{L }G_E)^{n}\modmod \wh G_E)$ et  $(\gamma_{1},...,\gamma_{n})
\in \on{Gal}(\ov F/F)^{n}$
on a  \begin{gather} \label{Theta-nu-non-deploye}
 f(\sigma(\gamma_{1}),...,\sigma(\gamma_{n}))=\big(\nu \circ \Theta_{n}(f)\big)(\gamma_{1},...,\gamma_{n}), 
\end{gather}
\item [] (C'4)   $\sigma$  se factorise à travers 
$ \pi_{1}(X\sm  \wh N, \ov\eta)$. 
\end{itemize}

La condition  f') ci-dessus implique: 
  \begin{itemize}
\item []  (C'5)  le  diagramme \eqref{diag-sigma} est commutatif. 
\end{itemize}

De plus   $\nu$ et la classe de conjugaison par  $\wh G(\Qlbar)$ de $\sigma$ se  déterminent mutuellement.  La propriété e') ci-dessus montre la compatibilité à l'isomorphisme de Satake en toutes les places de $X\sm  \wh N$, par le même argument que dans la preuve du \thmref{dec-param-cor-thm}. 

Cela termine la preuve  du \thmref{dec-param-cor-thm-non-deploye}. \cqfd 

 \begin{rem}\label{anisotropic-local-systems}
 Si $G$ est anisotrope, la démonstration  du \thmref{dec-param-cor-thm-non-deploye} pourrait être simplifiée car il n'y a plus besoin de  troncatures par $\mu$,   la cohomologie Hecke-finie est égale à la cohomologie totale et le lemme de Drinfeld peut être appliqué directement. En effet tous les champs $\Cht_{N,I,W}^{(I_{1},...,I_{k})}/\Xi$ sont  de type fini.  
 L'argument le plus simple pour le démontrer semble être de choisir un revêtement galoisien $\wh U\to U$  tel que l'image inverse de $G$ soit déployée, noter $\wh X$ la courbe projective lisse contenant $\wh U$ comme ouvert dense, puis appliquer à l'image inverse d'un chtouca à $\wh X$ (qui aura donc des pattes supplémentaires en les points 
 de $\wh X\sm \wh U$, contrôlées par des copoids dominants ne dépendant que de $G, X$ et $\wh X$) 
 la théorie de la filtration de Harder-Narasimhan   dans le cas déployé \cite{behrend-thesis,var,schieder} et utiliser les mêmes arguments que dans la preuve du théorème 2.25 de \cite{var} pour montrer que   des pentes trop importantes entraîneraient l'existence d'un sous-groupe parabolique propre  dans $G$ (cet argument apparaît déjà dans  \cite{laumon-rapoport-stuhler,laurent-asterisque,eike-lau} pour les algèbres à  division). 
 Donc 
 la limite inductive  
 $\varinjlim_{\mu }\mc H _{ N,I,W}^{0,\leq\mu,E}$ est un $E$-faisceau constructible  sur $(X\sm  \wh N)^{I}$, muni d'une action des    morphismes de  Frobenius partiels 
 (qui respectent de plus une  $\mc O_{E}$-structure naturelle). Ce faisceau est donc  lisse  sur $\mf U^{I}$ pour un certain ouvert dense  $\mf U\subset X$ et on peut appliquer directement le lemme de  Drinfeld.  Sa fibre en  $\Delta(\ov\eta)$ vérifie  les  propriétés a), b), c) de la \propref{prop-a-b-c}.     \end{rem}

  \noindent{ \bf Démonstration de la \propref{functoriality-L-groups-H-sigma}.}   
      Soit $W$ une  représentation
 de $({}^{L}G)^I$, $x\in W $ et  $\xi\in W^{*}$ invariants sous l'action diagonale  de $\wh G$, 
    et $(\gamma_i)_{i\in I}\in \pi_{1}(\eta,\ov\eta)^{I}$. 
      On note $W_{{}^{L}\Upsilon}$ la  représentation
de $({}^{L}G')^I$ obtenue en composant  $W$ avec  ${}^{L}\Upsilon:{}^{L}G'\to {}^{L}G$. 
Alors le morphisme de  restriction  
$\beta_{\Upsilon}^{*}$ de \eqref{restriction-G-G'} 
entrelace les opérateurs $S_{W_{{}^{L}\Upsilon},x,\xi,(\gamma_i)}$
et $S_{W,x,\xi,(\gamma_i)}$. 
En effet, grâce à une ``fonctorialité triviale'' de l'équivalence  de Satake géométrique
vis à vis de $\Upsilon$ (qui est laissée au lecteur), on définit naturellement un morphisme de faisceaux 
$$\beta_{\Upsilon}^{*}: \varinjlim _{\mu} \mc H _{ N, I, W_{{}^{L}\Upsilon}}^{0,\leq\mu,E}
\to \varinjlim _{\mu} \mc H _{ N, I, W}^{0,\leq\mu ,E}$$
(où le premier faisceau est relatif à $G'$ et $ \Xi' $ et le second à $G$ et $ \Xi$). Ces morphismes 
 entrelacent les   opérateurs  de  création et d'annihilation, ainsi que les  actions de Galois sur les parties Hecke-finies des fibres en $\Delta(\ov\eta)$.  La proposition \ref{functoriality-L-groups-H-sigma}   en résulte. \cqfd

\section{Résultats à coefficients dans  $\ov {\mathbb F_{\ell}}$}\label{mod-ell}
\label{chapitre-coeff-finis}

Je remercie Jean-Pierre Serre pour un commentaire qui m'a 
amené à réfléchir au problème considéré dans ce chapitre. 
Dans le cas  où $G$  est déployé, le théorème ci-dessous fournit une   décomposition canonique
de 
 $C_{c}^{\mr{cusp}}(G(F)\backslash G(\mb A)/K_{N}\Xi,\overline {\mathbb F_\ell})$
 indexée par les  classes de conjugaison  de représentations  complètement réductibles 
  $\pi_{1}(X\sm N, \ov\eta)\to \widehat G(\overline {\mathbb F_\ell})$. La notion de complète réductibilité est expliquée par  Serre dans \cite{bki-serre}. 

Soit $G$ un groupe  réductif comme dans le chapitre précédent. 

\begin{prop}
Le $\mc O_{E}$-réseau 
$$ C_{c}^{\mr{cusp}}(\Bun_{G,N}(\Fq)/\Xi,\mc O_{E})\subset 
 C_{c}^{\mr{cusp}}(\Bun_{G,N}(\Fq)/\Xi,E)$$ est préservé par tous les opérateurs d'excursion  
$S_{I,W,x,\xi,(\gamma_i)_{i\in I}}$ lorsque  $I$ et $(\gamma_i)_{i\in I}$ sont  arbitraires, et $W,x,\xi$ sont définis sur $\mc O_{E}$. 
\end{prop}
\dem 
Cela  résulte du fait que  l'équivalence de Satake géométrique a été construite  à coefficients dans   $\mc O_{E}$ par Mirkovic et Vilonen \cite{mv} (voir aussi \cite{ga-de-jong}).  On commence par le cas  où $G$ est déployé. 
Comme on l'a rappelé dans le  \thmref{thm-geom-satake} et la  \remref{rem-Satake-Gad-texte}, 
on a un  foncteur canonique $W\mapsto \mc S_{I,W,\mc O_{E}}^{(I_{1},...,I_{k})}$
de la catégorie des  représentations de $(\wh G)^{I}$ dans 
des $\mc O_{E}$-modules de type fini 
vers la  catégorie des   faisceaux pervers (à un décalage près) $G^{\mr{ad}}_{\sum \infty x_{i}}$-équivariants 
à coefficients dans  $\mc O_{E}$   sur  $\mr{Gr}_{I}^{(I_{1},...,I_{k})}$.   
Mirkovic et Vilonen  utilisent la  t-structure  perverse standard dans la catégorie dérivée des  $\mc O_{E}$-faisceaux. 
Le fait qu'elle n'est pas préservée par la dualité de  Verdier ne   pose pas de problème car celle-ci n'est pas utilisée dans l'argument. 
En fait on n'utilisera ce foncteur $W\mapsto \mc S_{I,W,\mc O_{E}}^{(I_{1},...,I_{k})}$
 que lorsque $W$ est sans torsion. 
Avec les  notations de la \propref{defi-F-E-Cht}
on obtient le  $\mc O_{E}$-faisceau  pervers  $\mc F_{N,I,W, \Xi,\mc O_{E}}^{(I_{1},...,I_{k})}$  sur $\Cht_{N,I,W}^{(I_{1},...,I_{k}),\leq\mu}/\Xi$ défini  comme l'image inverse de $\mc S_{I,W,\mc O_{E}}^{(I_{1},...,I_{k})}$. 
 On  définit alors 
 $ \mc H _{ N,I,W}^{0,\leq\mu,\mc O_{E}}$ comme le  sous-$\mc O_{E}$-faisceau de $ \mc H _{ N,I,W}^{0,\leq\mu,E}$ égal à {\it l'image de } 
 $$R^{0}
 \big(\mf p_{N,I}^{(I_{1},...,I_{k}),\leq\mu}\big)_{!}\Big(\restr{\mc F_{N,I,W,\Xi,\mc O_{E}}^{(I_{1},...,I_{k})}} {\Cht_{N,I,W}^{(I_{1},...,I_{k}),\leq\mu}/\Xi}\Big)$$ dans 
 $$
\mc H _{ N,I,W}^{0,\leq\mu,E}=  R^{0}
 \big(\mf p_{N,I}^{(I_{1},...,I_{k}),\leq\mu}\big)_{!}\Big(\restr{\mc F_{N,I,W,\Xi,E}^{(I_{1},...,I_{k})}} {\Cht_{N,I,W}^{(I_{1},...,I_{k}),\leq\mu}/\Xi}\Big)$$
 (autrement dit on tue la torsion éventuelle). 
Ces sous-$\mc O_{E}$-faisceaux 
$ \mc H _{ N,I,W}^{0,\leq\mu,\mc O_{E}}$ sont préservés par les morphismes de création et d'annihilation et, avec augmentation de $\mu$,   par l'action des opérateurs de Hecke à coefficients dans $\mc O_{E}$ et des morphismes de  Frobenius partiels. 

Pour tout point géométrique   $\ov x$  de $(X\sm N)^{I}$ on définit 
$$\Big( \varinjlim _{\mu}\restr{\mc H _{ N, I, W}^{0,\leq\mu,\mc O_{E}}}{\ov{x}}\Big)^{\mr{Hf}}=\Big( \varinjlim _{\mu}\restr{\mc H _{ N, I, W}^{0,\leq\mu,E}}{\ov{x}}\Big)^{\mr{Hf}}\cap \Big( \varinjlim _{\mu}\restr{\mc H _{ N, I, W}^{0,\leq\mu,\mc O_{E}}}{\ov{x}}\Big).$$
On vérifie facilement qu'un élément de   $\varinjlim _{\mu}\restr{\mc H _{ N, I, W}^{0,\leq\mu,\mc O_{E}}}{\ov{x}}$
est   Hecke-fini si et seulement si il vérifie les  conditions équivalentes suivantes: 
\begin{itemize}
\item il  appartient à  un  sous-$\mc O_{E}$-module $\mf M$ de type fini de  $\varinjlim _{\mu}\restr{\mc H _{ N, I, W}^{0,\leq\mu,\mc O_{E}}}{\ov{x}}$ qui est 
stable par  $T(f)$ pour tout   $f\in C_{c}(K_{N}\backslash G(\mb A)/K_{N},\mc O_{E})$, 
\item il vérifie la même  condition, et en plus $\mf M$ est   stable par  $\pi_{1}(x,\ov{x})$. 
\end{itemize}
Le même  argument que dans  
la \propref{surjectivite-Hecke-fini}
 montre que  l'image de    l'homomorphisme de spécialisation  
  \begin{gather}\label{texte-sp}\on{\mf{sp}}^{*}: \varinjlim _{\mu}\restr{\mc H _{ N, I, W}^{0,\leq\mu,\mc O_{E}}}{\Delta(\ov{\eta})}\to 
  \varinjlim _{\mu}\restr{\mc H _{ N, I, W}^{0,\leq\mu,\mc O_{E}}}{\ov{\eta^{I}}}\end{gather}
 contient 
  $\Big( \varinjlim _{\mu}\restr{\mc H _{ N, I, W}^{0,\leq\mu,\mc O_{E}}}{\ov{\eta^{I}}}\Big)^{\mr{Hf}}$. 
  Grâce à la \propref{injectivite-sp}  
    on en déduit que 
   \begin{gather}\label{texte-sp-isom}\on{\mf{sp}}^{*}: \Big(\varinjlim _{\mu}\restr{\mc H _{ N, I, W}^{0,\leq\mu,\mc O_{E}}}{\Delta(\ov{\eta})}\Big)^{\rm{Hf}}\to \Big(
  \varinjlim _{\mu}\restr{\mc H _{ N, I, W}^{0,\leq\mu,\mc O_{E}}}{\ov{\eta^{I}}}\Big)^{\rm{Hf}}\end{gather}
  est un isomorphisme. 
  Par conséquent  on peut remplacer $E$ par $\mc O_{E}$ partout dans la définition \eqref{diag-S} des opérateurs d'excursion. 
 Lorsque  $G$ n'est pas  nécessairement déployé, les arguments sont les mêmes, avec les  adaptations du chapitre précédent. 
  \cqfd

Les  fonctions sur $({}^{L} G)^{n}\modmod\widehat G $
 de la forme  $(g_{1},...,g_{n})\mapsto \s{\xi, (1,g_{1},...,g_{n}).x}$ où $I=\{0,...,n\}$ et $W,x,\xi$ sont définis sur $\mc O_{E}$ sont exactement les  fonctions régulières  sur $({}^{L} G)^{n}\modmod\widehat G $ définies sur $\mc O_{E}$.  Le même énoncé est vrai  à coefficients dans  $\Zlbar$ et $\Flbar$. 

On note $\mc B_{ \Flbar}\in \on{End}( C_{c}^{\mr{cusp}}(\Bun_{G,N}(\Fq)/\Xi,\Flbar))$ la  sous-algèbre commutative engendrée par 
les  réductions modulo l'idéal maximal de $\mc O_{E}$ des opérateurs 
$S_{I,W,x,\xi,(\gamma_i)_{i\in I}}$ quand  $I$ et $(\gamma_i)_{i\in I}$ sont  arbitraires, et $W,x,\xi$ sont définis sur $\mc O_{E}$.

Les  propriétés de   la définition-proposition \ref{cor-annul-Lambda-n}  satisfaites  par les  opérateurs d'excursion 
  sont préservées  par la  réduction modulo l'idéal maximal de $\mc O_{E}$. 
  
Pour construire les   représentations galoisiennes $\sigma$  à coefficients dans  $\Flbar$ à partir des   caractères de $\mc B_{ \Flbar}$, on utilise, lorsque $G$ est déployé,  
le  théorème 3.1 de 
Bate, Martin et  R\" ohrle \cite{bmr}, qui complète, en  caractéristique   $\ell\neq 0$,  les résultats  de  Richardson \cite{richardson} (la combinaison des deux est récapitulée dans le  théorème 3.7 de \cite{bki-serre}). Pour plus de détails on renvoie à la preuve du théorème 4.5 de \cite{boeckle-harris...}. 

Quand  $G$ n'est pas  déployé, les propriétés satisfaites  par les  opérateurs d'excursion (à savoir  les analogues  de a), b), c), d) de la définition-proposition \ref{cor-annul-Lambda-n}  ainsi que  
e') et f') du paragraphe \ref{subsubsection-construction-Langlands-non-split})
 sont également préservées  par la  réduction modulo l'idéal maximal de $\mc O_{E}$. Pour construire les morphismes $\sigma$ 
on a besoin de l'extension de la  notion de complète réductibilité aux groupes réductifs  non nécessairement connexes, qui est donnée dans le    paragraphe  6 de \cite{bmr}. En particulier l'image d'un   morphisme 
 $$\sigma:\pi_{1}(X\sm \wh N, \ov\eta)\to {}^{L} G(\Flbar) $$
 tel que le diagramme 
  \begin{gather}\label{diag-sigma-mod-l}
 \xymatrix{
\on{Gal}(\ov F/F) \ar[rr] ^{\sigma}
\ar[dr] 
&& {}^{L} G(\Flbar) \ar[dl] 
 \\
& \on{Gal}(\wt F/F) }\end{gather}
est commutatif est complètement réductible au  sens de  \cite{bmr} si et seulement si   tout  sous-groupe parabolique   $P$ de $ {}^{L} G(\Flbar) $ contenant cette  image 
admet un sous-groupe de Levi la contenant. 
En effet tous les  sous-groupes paraboliques  $P$ de ${}^{L}G$ se surjectant sur  ${}^{L}G/\widehat G$ sont de  type de Richardson dans le sens de \cite{bmr}
(dans le  paragraphe  I.3 de  \cite{borel-corvallis} Borel étudie   
  ces  sous-groupes dans le cadre des  $L$-groupes complexes, mais la description qu'il en donne s'étend à notre situation : leurs classes de conjugaison par $\wh G$ sont en bijection avec les parties $\on{Gal}(\wt F/F) $-invariantes de l'ensemble des racines simples et par la preuve du  lemme 3.5 de \cite{borel-corvallis} ils sont de Richardson). 

 On démontre donc  le  théorème suivant. 
 
  \begin{thm}\label{dec-param-cor-thm-non-deploye-mod-l}     
   On a une   décomposition canonique de  $C_{c}(K_{N}\backslash G(\mb A)/K_{N},\Flbar)$-modules 
        \begin{gather}\label{dec-param-cor-non-deploye-mod-l} 
  C_{c}^{\mr{cusp}}(\Bun_{G,N}(\Fq)/\Xi,\Flbar)=\oplus_{\sigma} \mf H_{\sigma}^{\Flbar}
       \end{gather} 
       indexée par  les 
classes de conjugaison par $\wh G(\Flbar)$    de morphismes
        $$\sigma:\pi_{1}(X\sm  \wh N, \ov\eta)\to {}^{L} G(\Flbar) $$
      tels que  
    \begin{itemize}
  \item $\sigma$ prend ses valeurs  dans ${}^{L}G(\mathbb F')$, où $\mathbb F'$  est une extension finie  de $\Fl$, et il est continu, 
  \item le diagramme \eqref{diag-sigma-mod-l} est commutatif, 
\item $\sigma$ est complètement réductible dans le sens précédent.  
    \end{itemize}
 Cette  décomposition est obtenue par la  décomposition spectrale de la famille commutative des opérateurs d'excursion au sens suivant: 
$\mf H_{\sigma}^{\Flbar}$ est l'espace propre généralisé   correspondant au système de valeurs propres $\s{\xi, (\sigma(\gamma_{i}))_{i\in I}.x}$ pour les  opérateurs obtenus à partir de $S_{I,W,x,\xi,(\gamma_i)_{i\in I}}$  par  réduction  
modulo l'idéal maximal de $\mc O_{E}$,  lorsque  $I$ et $(\gamma_i)_{i\in I}$ sont arbitraires et $W,x,\xi$  sont définis sur $\mc O_{E}$. 

 Cette  décomposition  est  compatible avec  l'isomorphisme de Satake au sens suivant: pour tout  $\sigma$, toute place $v\in |X|\sm  \wh N$ et toute représentation $V$ de ${}^{L} G_{v}$ définie sur $\mc O_{E}$, 
 $\mf H_{\sigma}
$ est inclus dans l'espace propre généralisé de 
 $h_{V,v}$ pour la valeur propre  $\chi_{V}(\sigma(\Frob_{v}))$. 
              \end{thm}

\begin{rem} Contrairement à la situation du \thmref{dec-param-cor-thm},
il n'est sans doute pas vrai en général que $\mf H_{\sigma}
$ soit inclus dans l'espace propre  de 
 $h_{V,v}$,  car les opérateurs de Hecke non ramifiés modulo $\ell$ ne sont sans doute pas  toujours diagonalisables. 
\end{rem}

Le lien entre les décompositions du théorème précédent et du \thmref{dec-param-cor-thm} est le suivant. 
Soit    $\sigma:\pi_{1}(X\sm  \wh N, \ov\eta)\to {}^{L} G(\Zlbar) $
 tel que 
  $$\sigma^{\mr{rat}}: \pi_{1}(X\sm  \wh N, \ov\eta)\xrightarrow{\sigma} {}^{L} G(\Zlbar) \hookrightarrow {}^{L}  G(\Qlbar)$$ satisfasse les  conditions du 
 \thmref{dec-param-cor-thm}. 
 Soit  
$\ov \sigma^{\mr{ss}}$ la semi-simplification de 
  $$\ov \sigma : \pi_{1}(X\sm  \wh N, \ov\eta)\xrightarrow{\sigma} {}^{L} G(\Zlbar) \twoheadrightarrow {}^{L} G(\Flbar),$$ c'est-à-dire que 
  l'on prend un   sous-groupe parabolique le plus petit possible  de ${}^{L} G$ contenant l'image de $\ov \sigma $ et alors $\ov \sigma^{\mr{ss}}$ est la  projection  de $\ov\sigma$ sur le  Levi. 
  Le théorème \ref{dec-param-cor-thm-non-deploye-mod-l}  implique alors  que 
   l'image de 
$\mf H_{\sigma^{\mr{rat}}}\cap C_{c}^{\mr{cusp}}(\Bun_{G,N}(\Fq)/\Xi,\Zlbar)$ dans $C_{c}^{\mr{cusp}}(\Bun_{G,N}(\Fq)/\Xi,\Flbar)$ est incluse dans $ \mf H_{\ov \sigma}^{\Flbar}$.

\section{Indications sur le cas des groupes métaplectiques}
\label{para-meta}

Ce chapitre indique sommairement comment les résultats de cet article s'étendent au cas des groupes métaplectiques.   

Les constructions ci-dessous sont inspirées de \cite{sergey-theta} 
et utilisent  de fa\c con essentielle l'extension au cas métaplectique 
de l'équivalence de Satake géométrique due à Finkelberg 
et  Lysenko \cite{finkelberg-lysenko},  généralisée par Reich \cite{reich} et étendue par Lysenko  au cas réductif \cite{lysenko-red}, puis traitée dans le cas le plus général possible par Gaitsgory et Lysenko \cite{dennis-sergey}. La littérature sur les groupes métaplectiques est très vaste et on renvoie à \cite{gan-gao,weissman-split, weissman-L} pour des références récentes. Ces références considèrent des extensions métaplectiques très générales, associées à la donnée d'une extension de $G$ par $K_{2}$ sur le corps de fonctions, et d'un entier $N$ divisant $q-1$. 
On va se contenter d'un cadre plus modeste, où 
$N$ est comme ci-dessus et l'extension métaplectique provient des racines $N$-ièmes de fibrés en droites de degré pair avec de bonnes propriétés de factorisation (ces fibrés sont notés  $\mc A$, $\mc B$  et $\mc B_{I }^{(I)}$ ci-dessous). 
Ce cadre 
  contient  le cas de l'extension métaplectique usuelle de $Sp_{2n}$. 
  Pour traiter le cas le plus général, il faudrait 
  sans doute considérer des gerbes de factorisation (au sens de \cite{dennis-sergey}) définies sur $\Fq$ (la \remref{rem-generale-meta} ci-dessous décrit la gerbe de factorisation dans le cas traité ici).    

Soient $G$, $U$, $\wt U$, $\wt F, N$ et $\wh N$ comme dans le chapitre \ref{para-non-deploye}. Les champs $\mr{Gr}_{I}^{(I_{1},...,I_{k})}$ et 
$\Hecke_{I}^{(I_{1},...,I_{k})}$ sont définis au-dessus de $U^{I}$. 
Soit  $\mc A$ un fibré en droite sur $\Bun_{G}$ 
et $\mc B$ un fibré $G(\mc O)$-équivariant sur 
la grassmannienne affine  muni  d'une  donnée  de factorisation. 
Autrement dit on se donne, pour tout ensemble fini $I$ un fibré en droites $\mc B_{I }^{(I)}$ sur 
$\mr{Gr}_{I }^{(I)}/G_{\sum \infty x_{i}}$, 
qui soit 
\begin{itemize}
\item  compatible à la fusion, c'est-à-dire que pour toute application $\zeta:I\to J$, on a un isomorphisme canonique 
 \begin{gather}\nonumber
  \Delta_{\zeta}^{*}\Big( \mc B_{I}^{(I)} \Big)\simeq 
 \mc B_{J}^{(J)}\end{gather} 
 où $\Delta_{\zeta}$ désigne (le  quotient par $G_{\sum \infty x_{i}}$ de) l'inclusion 
 $$\mr{Gr}_{J}^{(J)} =\mr{Gr}_{I}^{(I)}\times _{X^{I}}X^{J}\hookrightarrow \mr{Gr}_{I}^{(I)} ,$$
\item compatible à la convolution c'est-à-dire que pour toute partition $(I_{1},...,I_{k})$ de $I$, on ait un isomorphisme canonique entre 
 $$\mc B_{I }^{(I_{1},...,I_{k})}:=\Big(\pi^{(I_{1},...,I_{k})}_{(I)}\Big)^{*}\Big(\mc B_{I }^{(I)} \Big), 
$$ et l'image inverse de $\prod _{j=1}^{k}\mc B_{I_{j} }^{(I_{j})}$ par le morphisme évident 
$$\mr{Gr}_{I }^{(I_{1},...,I_{k})}/G_{\sum \infty x_{i}}\to \prod _{j=1}^{k}\mr{Gr}_{I_{j} }^{(I_{j})}/G_{\sum_{i\in I_{j}} \infty x_{i}}.$$
\end{itemize}
et que $B_{I }^{(I)}$ soit égal à $\mc B$ lorsque $I$ est un singleton. 
On suppose de plus que pour tout ensemble fini $I$
et   toute partition $(I_{1},...,I_{k})$ de $I$, 
on a 
 un isomorphisme entre
 \begin{itemize}
\item 
 le fibré sur $
\Hecke_{I}^{(I_{1},...,I_{k})}$ dont la fibre en $\big((x_i)_{i\in I},(\mc G_{0}\xrightarrow{\phi_{1}}\mc G_{1} ... \xrightarrow{\phi_{k}}\mc G_{k})\big)$   est 
$\mc A_{\mc G_{0}}\otimes \mc A_{\mc G_{k}}  ^{-1}$
\item l'image inverse de 
$\mc B_{I }^{(I_{1},...,I_{k})}$ par le morphisme évident 
$\Hecke_{I}^{(I_{1},...,I_{k})}\to \mr{Gr}_{I }^{(I_{1},...,I_{k})}/G_{\sum \infty x_{i}}$. 
\end{itemize}
De plus on demande que ce dernier  isomorphisme soit compatible à la fusion (c'est-à-dire fonctoriel en $I$) et à la convolution.

Par exemple si $G$ est déployé de tore maximal $T$, on obtient un tel couple $(\mc A,\mc B)$ 
pour toute représentation $W$ de $G$ (éventuellemement virtuelle) telle que 
\begin{gather}\label{condition-reich}\sum_{\alpha} \alpha\otimes \alpha\in 2 X^{*} (T)\otimes X^{*} (T)\end{gather} (où la somme porte sur les poids de $W$, comptés avec multiplicités). 
En effet on pose 
$\mc A_{\mc G}=\det(R\Gamma(X,W_{\mc G}))$ pour $\mc G\in \Bun_{G}$, 
et $\mc B_{(x,\mc G\xrightarrow{\phi}\mc G')} $ égal au déterminant relatif des réseaux associés aux restrictions de  $W_{\mc G}$ et $W_{\mc G'}$
au voisinage formel de $x$ (grâce à la condition \eqref{condition-reich} ces fibrés sont pairs et donc les règles de Koszul ne font pas apparaître  de signe négatif, cf \cite{reich} Proposition II.6.8). 

 On se donne  un entier $N$ divisant $q-1$. On note $\mu_{N}$ le groupe fini des racines $N$-ièmes de l'unité dans $\Fq^{*}$. 
 On se donne un caractère primitif $\zeta: \mu_{N}\to E^{*}$. 
 On note 
 $\wt{\Bun}_{G}$ la $\mu_{N}$-gerbe sur $\Bun_{G}$  
 associée au fibré en droites $\mc A$. Autrement dit     les $S$-points de  $\wt{\Bun}_{G}$
classifient  la donnée d'un morphisme $u:S\to \Bun_{G}$, d'un fibré en droites $\mf A$ sur $S$ et d'un isomorphisme $\mf A^{N}\simeq u^{*}(\mc A)$. De même on introduit la  $\mu_{N}$-gerbe  $\wt{\mr{Gr}}_{I}^{(I_{1},...,I_{k})}/G_{\sum_{i\in I} \infty x_{i}}$  sur 
 $\mr{Gr}_{I}^{(I_{1},...,I_{k})}/G_{\sum_{i\in I} \infty x_{i}}$  
associée  au fibré en droite  
   $ \mc B_{I }^{(I_{1},...,I_{k})}$. 
 Son image inverse 
$\wt{\Hecke}_{I}^{(I_{1},...,I_{k})}$ sur 
 $\Hecke_{I}^{(I_{1},...,I_{k})}$ est également associée au fibré en droites 
 $\mc A_{\mc G_{0}}\otimes \mc A_{\mc G_{k}}  ^{-1} $.

 \begin{rem}
 \label{rem-generale-meta}
 La $\mu_{N}$-gerbe ci-dessus est une gerbe de factorisation au sens de \cite{dennis-sergey}. Grâce à l'hypothèse que $N$ divise $q-1$, elle est définie sur $\Fq$. 
  \end{rem}
 
Les groupes duaux des groupes métaplectiques ont été introduits par Savin après des travaux de Kazhdan et Patterson. On renvoie à 
\cite{finkelberg-lysenko, reich, lysenko-red, mcnamara, weissman-split, weissman-L,dennis-sergey}   pour la définition de ces groupes (dans le cadre très général  de \cite{dennis-sergey} la notion  de donnée de Langlands duale métaplectique est introduite, mais nous ne la considérons pas ici). 
 
On suppose donc  qu'il existe un groupe 
${}^{L }\wt G$ extension de $\on{Gal}(\wt F/F)$ par ${}^{L }\wt G^{0}$,   tel que 
l'analogue du théorème \ref{thm-geom-satake-non-deploye} soit vraie. Plus précisément on suppose que l'on a   un  foncteur canonique 
$$W\mapsto \mc S_{I,W,E}^{(I_{1},...,I_{k})}$$
de la catégorie des représentations $E$-linéaires de dimension finie  de $({}^{L}\wt  G)^{I}$ vers la  catégorie de    $E$-faisceaux pervers  (à un décalage près)   sur  $\wt{\mr{Gr}}_{I}^{(I_{1},...,I_{k})}/G_{\sum_{i\in I} \infty x_{i}}$   sur lesquels $\mu_{N}$ agit par $\zeta$ (et où  le décalage est déterminé par la condition que  l'image inverse sur  $\wt{\mr{Gr}}_{I}^{(I_{1},...,I_{k})}$ est perverse relativement à  $(X\sm \wh N)^{I}$). Les  $E$-faisceaux pervers  (à un décalage près) $\mc S_{I,W,E}^{(I_{1},...,I_{k})}$ doivent être  universellement localement acycliques relativement au morphisme vers $(X\sm \wh N)^{I}$. Ils doivent vérifier  les mêmes propriétés b), c), d) que dans    le théorème \ref{thm-geom-satake-non-deploye}.  

L'existence d'un tel groupe est justifiée dans \cite{finkelberg-lysenko, reich, lysenko-red, dennis-sergey} lorsque $G$ est déployé.  
 
  Soit   $I$ un ensemble fini, $(I_{1},...,I_{k})$ une partition de $I$ et $W$ une représentation irréductible de $({}^{L }G)^{I}$. 
     
   Le  morphisme évident 
   $\Cht_{N,I } ^{(I_{1},...,I_{k})}\to  \mr{Gr}_{I}^{(I_{1},...,I_{k})}/G_{\sum_{i\in I} \infty x_{i}}$ 
   se relève en un morphisme 
   \begin{gather}\label{mor-Cht-gerbe}
   \Cht_{N,I } ^{(I_{1},...,I_{k})}\to \wt{\mr{Gr}}_{I}^{(I_{1},...,I_{k})}/G_{\sum_{i\in I} \infty x_{i}}. \end{gather}
   En effet la restriction de la $\mu_{N}$-gerbe $\wt{\Hecke}_{I}^{(I_{1},...,I_{k})}$ à 
$\Cht_{N,I } ^{(I_{1},...,I_{k})}$  est  trivialisée par  $$\mc A_{\mc G_{0}}\otimes \mc A_{\ta{\mc G_{0}}}  ^{-1} =
(\mc A_{\mc G_{0}})^{1-q} =
\big((\mc A_{\mc G_{0}})^{-\frac{q-1}{N}}\big)^{N}
.$$ 
 
On note 
 $\Cht_{N,I,W } ^{(I_{1},...,I_{k})}$ le support de l'image inverse de 
$\mc S_{I,W,E}^{(I_{1},...,I_{k})}$ par \eqref{mor-Cht-gerbe}. 

Contrairement au cas non métaplectique, il n'est pas vrai en général que l'action de $G_{\sum_{i\in I} \infty x_{i}}$
 sur $\wt{\mr{Gr}}_{I}^{(I_{1},...,I_{k})}$ 
 se factorise par $G^{\mr{ad}}_{\sum_{i\in I} \infty x_{i}}$, c'est-à-dire soit triviale sur $Z_{\sum_{i\in I} \infty x_{i}}$. En revanche, en faisant comme \cite{sergey-tores} 5.2.3 (qui traite le cas des tores)  on montre qu'il existe une isogénie $Z^{\sharp}\to Z$ telle que 
 \begin{itemize}
 \item l'action de $Z^{\sharp}_{\sum_{i\in I} \infty x_{i}}$ sur $\wt{\mr{Gr}}_{I}^{(I_{1},...,I_{k})}$  soit trivialisée, 
 \item $\Bun_{Z^{\sharp},N}$ agisse sur $\wt \Bun_{G,N}$ de fa\c con compatible avec la trivialisation précédente. 
 \end{itemize}
 On suppose que  l'équivalence de Satake géométrique admise plus haut est telle que $Z^{\sharp}_{\sum_{i\in I} \infty x_{i}}$ agit naturellement sur les faisceaux
 $\mc S_{I,W,E}^{(I_{1},...,I_{k})}$. 
 
 On choisit alors un réseau $\Xi$ dans $\Bun_{Z^{\sharp},N}(\Fq)$
   tel qu'il s'injecte dans $\Bun_{Z,N}(\Fq)$, et on note encore $\Xi$ son image.

Alors l'image inverse de 
$\mc S_{I,W,E}^{(I_{1},...,I_{k})}$ par \eqref{mor-Cht-gerbe} est naturellement 
$\Xi$-équivariante et on note 
 $\mc F_{N,I,W,\Xi,E}^{(I_{1},...,I_{k})}$ le faisceau pervers (à un décalage près) 
 sur $\Cht_{N,I,W } ^{(I_{1},...,I_{k})}/\Xi$ qui en résulte.
  
 On pose alors 
    \begin{gather} 
 \wt{\mc H }_{ N,I,W}^{\leq\mu,E}=  R
 \big(\mf p_{N,I}^{(I_{1},...,I_{k}),\leq\mu}\big)_{!}\Big(\restr{\mc F_{N,I,W,\Xi,E}^{(I_{1},...,I_{k})}} {\Cht_{N,I,W}^{(I_{1},...,I_{k}),\leq\mu}/\Xi}\Big) . 
  \end{gather}
       qui ne dépend pas du choix de la partition $(I_{1},...,I_{k})$.  On a en particulier le faisceau constructible  
       $\wt{\mc H} _{ N,I,V}^{0,\leq\mu,\Qlbar}$ sur $(X\sm \wh N)^{I}$. 

   Dans le cas où $I$ est vide, $\wt{\mc H} _{ N,\emptyset,\mbf 1}^{0,\leq\mu,E} $  est un faisceau constructible (trivial) sur $\Fq$ et on a  
       \begin{gather}\label{defi-HIW-meta} \varinjlim_{\mu }\restr{\wt{\mc H }_{ N,\emptyset,\mbf 1}^{0,\leq\mu,E}}{\Fqbar}=
    C_{ \zeta}(\wt{\Bun_{G,N } (\Fq)}/\Xi,E)\end{gather}
       où 
       $\wt{\Bun_{G,N } (\Fq)}$ est le $\mu_{N}$-torseur sur $\Bun_{G,N } (\Fq)$ dont la fibre 
       en $\mc G$ est l'image 
       par le morphisme surjectif $\Fq^{*}\to   \mu_{N}, a \mapsto a^{-\frac{q-1}{N}}$ 
       du $\Fq^{*}$-torseur des trivialisations de la $\Fq$-droite 
       $\mc A_{\mc G}$, et où 
       le membre de droite de \eqref{defi-HIW-meta} désigne l'espace des fonctions à support compact qui sont $\zeta$-équivariantes.

        La notion de cuspidalité s'étend sans modification au cas métaplectique et permet   de définir un sous-espace de dimension finie $C_{ \zeta}^{\mr{cusp}}(\wt{\Bun_{G,N } (\Fq)}/\Xi,E)$. 
     Les mêmes arguments que dans  le  chapitre \ref{para-non-deploye} fournissent alors une décomposition de $C_{ \zeta}^{\mr{cusp}}(\wt{\Bun_{G,N } (\Fq)}/\Xi,E)$ suivant des  paramètres de Langlands, c'est-à-dire des classes de conjugaison par ${}^{L }\wt G^{0}(\Qlbar)$ de morphismes
        $\sigma:\pi_{1}(X\sm \wh N, \ov\eta)\to {}^{L }\wt G(\Qlbar) $ 
    définis sur une     extension finie de $\Ql$,  continus et tels que 
 l'adhérence de Zariski de l'image soit réductive. Cette décomposition est compatible avec l'isomorphisme de Satake (dont une définition possible consiste à prendre les traces de Frobenius dans l'équivalence de Satake géométrique  discutée ci-dessus) et réalise donc la correspondance de Langlands dans le sens
  ``automorphe vers Galois''. On renvoie à \cite{mcnamara, gan-gao,  weissman-L} pour l'énoncé de l'isomorphisme de Satake classique et le lien avec 
  \cite{finkelberg-lysenko, lysenko-red, dennis-sergey} lorsque $G$ est déployé. 
   Avec les notations de  \eqref{dec-alpha-general-ker1} on a 
    \begin{gather} 
    \wt{\Bun_{G,N } (\Fq)}=
  \bigcup_{\alpha\in \ker^{1}(F,G)}G_{\alpha}(F)\backslash \wt{G_{\alpha}(\mb A)}/K_{N}\end{gather} où   $\wt {G_{\alpha}(\mb A)}$ est une extension de $G_{\alpha}(\mb A)$ par $\mu_{N}$. 
  
  \begin{rem}
  On peut espérer que les résultats à coefficients dans $\Flbar$
  du chapitre \ref{chapitre-coeff-finis} 
  s'étendent au cas métaplectique mais on ne peut pas le montrer pour le moment parce que  \cite{finkelberg-lysenko} et  \cite{lysenko-red} ne traitent pas les coefficients dans $\mc O_{E}$.   \end{rem}

      \section{Lien avec le programme de Langlands géométrique}
     \label{subsection-link-langl-geom}
  Il est évident que la coalescence   et la  permutation des pattes 
sont reliées aux  structures de  factorisation introduites par  
  Beilinson et Drinfeld
 \cite{chiral} et en effet notre article utilise de fa\c con essentielle   le produit de fusion  sur la grassmannienne affine de Beilinson-Drinfeld dans l'équivalence de  Satake  géométrique     \cite{mv,ga-de-jong}.  
 Par ailleurs l'idée de décomposition spectrale est familière dans le programme de Langlands géométrique, cf \cite{beilinson-heisenberg} et surtout 
 le corollaire  4.5.5 de 
 \cite{dennis-laumon} qui affirme,  dans le cadre du programme de Langlands géométrique pour les $D$-modules (où la courbe $X$ est définie sur un corps algébriquement clos de caractéristique $0$) que la DG-catégorie des $D$-modules sur 
 $\Bun_{G}$ est ``au-dessus'' du champ des systèmes locaux pour 
 $\wh G$.  On notera curieusement que l'on ne sait pas formuler d'énoncé analogue avec les faisceaux $\ell$-adiques lorsque $X$ est sur $\Fq$
 (même si la conjecture d'annulation montrée par Gaitsgory \cite{ga-vanishing} va dans ce sens), et cependant notre article peut être considéré comme une version ``classique'' ou plutôt ``arithmétique'' d'un tel  énoncé.  
 
 En fait le lien est bien plus direct qu'une simple analogie: 
 nous allons voir   que les conjectures du programme de Langlands géométrique $\ell$-adique  permettent  de comprendre les opérateurs d'excursion et fournissent une explication très éclairante de notre approche grâce à une construction de Braverman et Varshavsky \cite{brav-var} qui généralise le fait qu'un faisceau sur $\Bun_{G}$ donne par les  traces de Frobenius une fonction sur $\Bun_{G}(\Fq)$. 
On prend ici  $G$ déployé et $N$ vide, c'est-à-dire  $K_{N}=G(\mathbb O)$ mais les considérations qui suivent resteraient valables pour tout niveau $N$ (et même pour les groupes réductifs non déployés). 

   Les conjectures du programme de Langlands   géométrique sont formulées à l'aide des foncteurs de  Hecke:  pour toute représentation $W$ de $(\wh G)^{I}$  le foncteur de  Hecke
     $$\phi_{I,W}:D^{b}_{c}(\Bun_{G},\Qlbar)\to  D^{b}_{c}(\Bun_{G}\times X^{I}, \Qlbar)$$      
    est donné  par 
     $$\phi_{I,W}(\mc F)=q_{1,!}\big(q_{0}^{*}(\mc F)\otimes \mc F_{I,W}\big)$$
 où $\Bun_{G}\xleftarrow{q_{0}}\Hecke_{I,W}^{(I)}\xrightarrow{q_{1}}\Bun_{G}\times X^{I}$ est la   correspondance de Hecke    et  
 \begin{itemize}
 \item quand  $W$ est irréductible, $\mc F_{I,W}$ est égal,  à un décalage près,  au faisceau d'intersection  de $\Hecke_{I,W}^{(I)}$   
 \item en général il est défini, de fa\c con fonctorielle en $W$,  comme  l'image inverse de 
 $\mc S_{I,W,\Qlbar}^{(I)}$ par le morphisme lisse naturel 
 $\Hecke_{I,W}^{(I)}\to \mr{Gr}_{I,W }^{(I)}/G_{\sum n_{i}x_{i}}$ (où les $n_{i}$ sont assez grands).  
 \end{itemize}
Soit $\mc E$  un $\wh G$-système local sur $X$. Alors  $\mc F\in D^{b}_{c}(\Bun_{G},\Qlbar)$ est dit propre pour  $\mc E$  si l'on possède,  pour tout ensemble fini $I$ et toute représentation $W$ de 
  $(\wh G)^{I}$,  un isomorphisme 
  $\phi_{I,W}(\mc F)\isom \mc F\boxtimes W_{\mc E}$,   fonctoriel en  $W$, et compatible aux produits extérieurs et à la  fusion (c'est-à-dire  à l'image inverse par le   morphisme diagonal $X^{J}\to X^{I}$ associé à n'importe quelle application $I\to J$). Les  conjectures du programme de Langlands   géométrique impliquent l'existence d'un objet $\mc F$ propre pour  $\mc E$   (et vérifiant une condition de normalisation de Whittaker qui l'empêche en particulier d'être nul). Dans le programme de Langlands   géométrique  $X$ et $\Bun_{G}$ sont habituellement définis sur un corps algébriquement clos, mais ici nous les prenons définis sur $\Fq$. 
 
Soit $\mc F$ propre pour  $\mc E$. On note  $f\in C(\Bun_{G}(\Fq),\Qlbar)$ la fonction associée à $\mc F$ par le dictionnaire  faisceaux-fonctions, c'est-à-dire que  pour $x\in \Bun_{G}(\Fq)$, $f(x)=\on{Tr}(\Frob_{x}, \restr{\mc F}{x})$. 
Soit $\Xi\subset Z(F)\backslash Z(\mathbb A)$ un réseau. On suppose que  $\mc F$ est $\Xi$-équivariant, si bien que  $f\in  C(\Bun_{G}(\Fq)/\Xi,\Qlbar)$ (quitte à diminuer $\Xi$ cela est impliqué par une condition   sur $\mc E$, en fait sur son image par le morphisme de $\wh G$ vers son abélianisé). 
Il est bien connu que  $f$ est un vecteur propre pour tous les  opérateurs de Hecke: pour toute place $v$ et toute représentation  irréductible $V$ de $\wh G$, 
$T(h_{V,v})(f)=\on{Tr}(\Frob_{v},\restr{V_{\mc E}}{v})f$, où 
$\Frob_{v}$ est un  élément de  Frobenius en $v$. 

La proposition suivante (reposant sur un résultat non encore rédigé) exprime la compatibilité entre  le programme de Langlands géométrique et la décomposition \eqref{intro1-dec-canonique}. 

\begin{prop}\label{prop-compat-langl-geom} (conditionnelle à un résultat non encore rédigé) 
Etant donné $\mc F$ propre pour $\mc E$ et $\Xi$-équivariant tel que   la fonction $f$ associée à $\mc F$ soit cuspidale, alors $f$ appartient à $\mc H_{\sigma}$  où 
 $\sigma:\pi_{1}(X,\ov\eta)\to 
\wh G(\Qlbar)$ est la   représentation galoisienne correspondant au système local $\mc E$. 
\end{prop}
\dem 
Dans  \cite{brav-var}, Braverman et Varshavsky utilisent un morphisme de trace très général, et le fait que  
$\Cht_{I,W}^{(I)}$ est l'intersection de la  correspondance de Hecke avec le  graphe de l'endomorphisme de  Frobenius  de $\Bun_{G}$, pour construire  un  morphisme de faisceaux sur $X^{I}$
\begin{gather}\label{trace-brav-var-piIW}\pi^{\mc F,\mc E}_{I,W}: \varinjlim_{\mu}\mc H_{N,I,W}^{0,\leq\mu,\Qlbar}\to W_{\mc E}.\end{gather}
Ces morphismes sont fonctoriels et $\Qlbar$-linéaires en  $W$, et compatibles avec la  coalescence des pattes et avec l'action des morphismes de  Frobenius partiels (ce dernier point n'a pas encore été rédigé). 
De plus   $\pi^{\mc F,\mc E}_{\emptyset,\mbf 1}:C_{c}(\Bun_{G}(\Fq)/\Xi,\Qlbar)\to \Qlbar$ n'est autre que  $h\mapsto \int_{\Bun_{G}(\Fq)/\Xi}     fh$. 
Alors on déduit des propriétés de ces  morphismes $\pi^{\mc F,\mc E}_{I,W}$ que pour tout  $I,W,x,\xi,(\gamma_{i})_{i\in I}$, 
on a 
$$S_{I,W,x,\xi,(\gamma_{i})_{i\in I}}(f)=
\s{\xi, (\sigma(\gamma_{i}))_{i\in I}.x}f.$$ 
Ceci termine la preuve de la \propref{prop-compat-langl-geom}. \cqfd

\section{Cas de $GL_{r}$}\label{GL-previous-works}
          
      La correspondance de Langlands pour $GL_{r}$ sur les corps de fonctions  est établie dans 
       \cite{laurent-inventiones} (après  \cite{drinfeld-proof-peterson} pour $GL_{2}$) à l'aide du principe de récurrence de Piatetski-Shapiro et Deligne. Dans ce chapitre,  qui n'apporte aucun résultat nouveau, on montre   comment le  \thmref{intro-thm-ppal}
       fournit une nouvelle preuve de l'ingrédient de récurrence  (\eqref{hyp-laurent} ci-dessous), et donc une nouvelle preuve de la correspondance de Langlands pour $GL_{r}$, si on ajoute aux arguments ci-dessous le paragraphe 6.1 et l'appendice B de \cite{laurent-inventiones}, qui expliquent la récurrence.

\begin{rem} Certaines parties de la preuve du \thmref{intro-thm-ppal}
se simplifient dans le cas de $GL_{r}$. D'abord il suffit d'établir la proposition 
\ref{prop-coal-frob-cas-part} dans le cas où $V=\Lambda^{k}\mr{St}$ pour $k=1,...,r$ (car ces représentations engendrent tout l'anneau des représentations de $GL_{r}$).  Dans ce cas $V$  est minuscule, mais comme $\deg(v)$ n'est pas nécessairement égal à $1$, la preuve est 
un petit peu plus compliquée que celle donnée dans l'introduction
mais quand même beaucoup plus simple que la preuve de la 
\propref{prop-coal-frob-cas-part}. 
D'autre part, comme on l'a dit dans la remarque \ref{rem-pseudo-car-Taylor}, dans le cas particulier où $G=GL_{r}$,  la \propref{Xi-n} (c'est-à-dire la 
\propref{intro-Xi-n} de l'introduction) avait déjà été montrée par Taylor \cite{taylor} grâce à  la théorie des pseudo-caractères.  
\end{rem}

  La suite de ce chapitre fournit les quelques arguments nécessaires
  pour donner une nouvelle preuve de l'ingrédient de récurrence de 
  \cite{laurent-inventiones} à l'aide du \thmref{intro-thm-ppal}, en faisant attention à ne pas utiliser les résultats connus comme conséquences de  \cite{laurent-inventiones} pour éviter toute circularité. C'est pour éviter une telle circularité
  que l'on  ajoute la dernière assertion du \lemref{lem-sigma-dans-coho} ci-dessous  (qui  résulte {\it a posteriori} de  \cite{laurent-inventiones} par  la 
\remref{rem-eviter-circ}) et qu'on montre le \lemref{lem-irred} ci-dessous 
(qui lui aussi  résulte {\it a posteriori} de  \cite{laurent-inventiones}). 

Bien que l'on soit intéressé ici par $G=GL_{r}$ on énonce le lemme suivant pour $G$ quelconque (déployé pour simplifier). 
        
             \begin{lem}\label{lem-sigma-dans-coho} 
 Soit $\sigma$ 
 apparaissant dans la décomposition  \eqref{dec-param-cor} du 
 \thmref{dec-param-cor-thm}. 
    Soit  $V$ une  représentation irréductible de $\wh G$ et 
     $V_{\sigma}=\oplus_{\tau} \tau \otimes \mf V_{\tau}$ la  décomposition de la   représentation semi-simple $V_{\sigma}$ indexée par les 
     classes d'isomorphisme  de  représentations  irréductibles $\tau $ de $\pi_{1}(\eta, \ov\eta)$. Alors si  $\mf V_{\tau}\neq 0$, $\tau\boxtimes \tau ^{*}$ apparaît comme un sous-quotient de la   représentation 
     \begin{gather}\label{rep-12-VV*}  H_{\{1,2\}, V\boxtimes V^{*}}=    \Big(\varinjlim _{\mu}\restr{\mc H _{N, \{1,2\}, V\boxtimes V^{*}}^{0,\leq\mu,E}}{\ov{\eta^{\{1,2\}}}}\Big)^{\mr{Hf}}\end{gather}  de 
  $(\pi_{1}(\eta,\ov\eta))^{2}$. 
   Plus précisément pour tout    $h\in \mf H_{\sigma}$ non nul, il existe un  isomorphisme entre 
  $\tau\boxtimes \tau ^{*}$ et  un  quotient de la sous-représentation de dimension finie   de 
     \eqref{rep-12-VV*} 
      engendrée par      $f=\on{\mf{sp}}^{*}(\mc C_{  \delta_{V}}^{\sharp }(h))$, qui envoie  
       $\Id_{\tau}\in  \tau\boxtimes \tau ^{*}$ sur   l'image de   $f$ dans ce  quotient.

        De plus $\tau$ est $\iota$-pure pour tout isomorphisme $\iota:\Qlbar \isom \C$.      \end{lem}

\begin{rem} \label{rem-eviter-circ} Bien sûr la dernière assertion n'est pas un résultat nouveau  car d'après le théorème VII.6 de   \cite{laurent-inventiones}   toute représentation  irréductible (définie sur une extension finie de $\Ql$ et continue) de $\pi_{1}(X\sm N, \ov\eta)$ est $\iota$-pure pour tout $\iota$. 
\end{rem}

   \noindent{\bf Démonstration. } Quitte à augmenter $E$ on suppose $\sigma$ et $\mf H_{\sigma}$ définis sur $E$. Soit $h\neq 0$ dans $\mf H_{\sigma}$. 
  On rappelle qu'on ne sait pas si $\mc B$ est réduite, mais comme $\mc B$ est commutative, le sous-$\mc B$-module $\mc B.h$ de $\mf H_{\sigma}$ engendré par $h$ admet un quotient de dimension $1$, sur lequel $\mc B$ agit forcément par le caractère $\nu$ associé à $\sigma$ par (C3). Il existe donc une forme linéaire $\lambda$ sur $\mc B.h$ telle que 
  \begin{gather}\label{lambda-restriction}
  \lambda(bh)=\nu(b) \text{ \  pour tout  \ } b\in \mc B. \end{gather} 
        Soit  $\check h\in   C_{c}^{\rm{cusp}}(G(F)\backslash G(\mb A)/K_{N}\Xi,E)$ tel que la forme linéaire sur 
     $C_{c}^{\rm{cusp}}(G(F)\backslash G(\mb A)/K_{N}\Xi,E)$ donnée par 
     \begin{gather}\label{int-check-h-h}k\mapsto \int_{G(F)\backslash G(\mb A)/K_{N}\Xi} \check h \,  k\end{gather}  
     prolonge $\lambda$. 
     
L'élément  $f=\on{\mf{sp}}^{*}(\mc C_{  \delta_{V}}^{\sharp }(h))$ 
  est aussi  l'élément de \eqref{rep-12-VV*} image de $h$ par la composée 
 $$C_{c}^{\rm{cusp}}(G(F)\backslash G(\mb A)/K_{N}\Xi,E)=
  H_{\{0\},\mbf  1}\xrightarrow{\mc H(\delta_{V})}
 H_{\{0\},V\otimes V^{*}}\isor{\chi_{\zeta_{\{1,2\}}}^{-1}} 
 H_{\{1,2\}, V\boxtimes V^{*}}$$
 On note  $\check f$ la forme linéaire   sur  \eqref{rep-12-VV*}
  égale à la composée de  
 $$H_{\{1,2\}, V\boxtimes V^{*}} \isor{\chi_{\zeta_{\{1,2\}}}} H_{\{0\},V\otimes V^{*}}  \xrightarrow{\mc H(\on{ev}_{V})} 
  H_{\{0\},\mbf  1} =C_{c}^{\rm{cusp}}(G(F)\backslash G(\mb A)/K_{N}\Xi,E)$$
     et de la forme linéaire
     $$C_{c}^{\rm{cusp}}(G(F)\backslash G(\mb A)/K_{N}\Xi,E)\to E, \ \ g\mapsto \int_{G(F)\backslash G(\mb A)/K_{N}\Xi} \check h g.$$
        Alors $f$  et $\check f$ 
      sont invariants par l'action diagonale de $\pi_{1}(\eta,\ov\eta)$ (cf \remref{rem-gamma-en-plus}). 
    Pour tout  $(\gamma, \gamma')\in (\pi_{1}(\eta,\ov\eta))^{2}$ on a  
        \begin{gather*}\s{\check f , (\gamma, \gamma') \cdot f}=
       \int_{G(F)\backslash G(\mb A)/K_{N}\Xi} \check h S_{\{1,2\}, V\boxtimes V^{*},\delta_{V},\on{ev}_{V},(\gamma,\gamma')}(h)\\
    =\nu(S_{\{1,2\}, V\boxtimes V^{*},\delta_{V},\on{ev}_{V},(\gamma,\gamma')})       =
        \chi_{V}(\sigma(\gamma\gamma'^{-1}))=  \chi_{V_{\sigma}}(\gamma\gamma'^{-1}),\end{gather*}  où 
   \begin{itemize}
   \item  la première égalité  vient de la définition \eqref{excursion-def-texte}  des opérateurs d'excursion, 
      \item la deuxième égalité résulte de   l'hypothèse que  \eqref{int-check-h-h} prolonge $\lambda$, et de \eqref{lambda-restriction}, 
      \item la troisième égalité vient du fait que $\nu$ correspond  à $\sigma$ par (C3). 
        \end{itemize}
  
    Le   quotient de la   représentation de $(\pi_{1}(\eta,\ov\eta))^{2}$ engendrée par  $f$ par la plus grande sous-représentation sur laquelle  $\check f$ s'annule  est alors isomorphe à la sous-représentation  engendrée par 
       $\chi_{V_{\sigma}}$  dans  
     $C( \pi_{1}(\eta,\ov\eta), E)$ muni  de l'action  par translations à gauche et à droite de $(\pi_{1}(\eta,\ov\eta))^{2}$. D'après \cite{Bki-A8} chapitre 20.5 théorème 1, cette dernière représentation est isomorphe à 
   $\oplus_{\tau, \mf V_{\tau}\neq 0} \tau\boxtimes \tau ^{*}$    et on peut normaliser cet  isomorphisme de telle sorte qu'il envoie   $\chi_{V_{\sigma}}$   sur   $\sum_{\tau, \mf V_{\tau}\neq 0} \Id_{\tau}$.    
    On a donc montré que  si   $\mf V_{\tau}\neq 0$, $\tau\boxtimes \tau ^{*}$ est un  quotient d'une  sous-représentation de  \eqref{rep-12-VV*} contenant $f$, de telle sorte que  l'image de 
     $f$ soit  $\Id_{\tau}$. 
      
     On en déduit maintenant que $\tau$ est $\iota$-pure.  Comme $\tau\boxtimes \tau^{*}$ est un sous-quotient de \eqref{rep-12-VV*}, 
          il résulte de Weil II \cite{weil2} que $\tau\boxtimes \tau ^{*}$ est 
  $\iota$-pure de poids $\leq 0$ comme représentation de $\pi_{1}((X\sm N)^{2}, \Delta(\ov\eta))$. 
 Donc pour presque toute place $v$ les valeurs propres de $\tau(\Frob_{v})$ ont  des $\iota$-poids égaux, dont la valeur est déterminée par le $\iota$-poids de $\det(\tau)$. 
     \cqfd

      On rappelle que dans \cite{laurent-inventiones} la correspondance de Langlands est obtenue par récurrence sur $r$,  à l'aide du   principe de récurrence de Deligne, qui combine 
                                    \begin{itemize}
                                  \item
    les équations fonctionnelles
de fonctions $L$ de Grothendieck \cite{sga5}, 
\item la formule du produit de Laumon \cite{laumon-produit}, 
\item les théorèmes de multiplicité $1$ \cite{piat71,shalika}  
et les  théorèmes réciproques de Hecke, Weil,  Piatetski-Shapiro et Cogdell \cite{inverse-thm}.  
\end{itemize}
    La récurrence est expliquée dans le paragraphe 6.1 et l'appendice B de 
    \cite{laurent-inventiones} et admet comme ingrédient
    \begin{gather}\label{hyp-laurent}
  \text{  l'hypothèse de la   proposition VI.11 (ii) de   \cite{laurent-inventiones}}
  \end{gather}
   à savoir qu'à toute représentation  automorphe cuspidale $\pi$ pour $GL_{r}$ de niveau $N$ on peut associer   $\sigma: \pi_{1}(X\sm N, \ov\eta)\to GL_{r}(\Qlbar) $ défini sur une extension finie de  $\Ql$, continu, 
    pur de poids $0$ et correspondant à $\pi$ au sens de Satake en toutes les places de $X\sm N$. Le théorème
    \ref{intro-thm-ppal} (plus précisément le \thmref{dec-param-cor-thm})  fournit une nouvelle démonstration de  \eqref{hyp-laurent}, grâce au  lemme suivant. 
    
    \begin{lem} \label{lem-irred} On prend $G=GL_{r}$. 
    On suppose la correspondance de Langlands connue pour  $GL_{r'}$ pour tout $r'<r$. 
         Alors tout $\sigma$ 
    apparaissant dans la décomposition \eqref{intro1-dec-canonique} est irréductible et pur  de poids  $0$. 
    \end{lem}
      
 \noindent{\bf Démonstration (extraite de \cite{laurent-inventiones}). }   
 Soit $(H_{\pi},\pi)$ une représentation automorphe cuspidale pour $GL_{r}$ 
 telle que $(H_{\pi})^{K_{N}}$ soit non nul et apparaisse dans $\mf H_{\sigma}$. 
      On note  
   \begin{gather}\label{dec-sigma}\sigma=\oplus_{\tau} \tau \otimes \mf V_{\tau}\end{gather} 
   la décomposition de la représentation semi-simple $\sigma$ suivant les classes d'équivalence de  représentations irréductibles $\tau $ de $\pi_{1}(X\sm N, \ov\eta)$. 
    On suppose par l'absurde que $\sigma$ n'est pas irréductible. 
   Toute représentation   $\tau$ telle que $\mf V_{\tau}\neq 0$ est alors de rang  $r_{\tau}<r$ et,  par 
   l'hypothèse de récurrence incluse dans l'énoncé, il existe une   représentation automorphe cuspidale $\pi_{\tau}$  pour $GL_{r_{\tau}}$  associée à $\tau$ par la correspondance de Langlands pour $GL_{r_{\tau}}$.  On choisit un ensemble fini de places $S$ en dehors duquel  les représentations $\pi,\sigma$, et $\tau $, $\pi_{\tau}$ telles que $\mf V_{\tau}\neq 0$ (qui sont en nombre fini)  sont non ramifiées et se correspondent par l'isomorphisme de Satake. 
 D'après la méthode de Rankin-Selberg pour $GL_{r}\times GL_{r_{\tau}}$ (due à Jacquet, Piatetski-Shapiro, Shalika \cite{jacquet-shalika-euler,rankin-selberg}), 
$L(\check{\pi} \times \pi_{\tau},Z)$ est un polynôme en $Z$  
donc a fortiori 
 $L_{S}(\check{\pi} \times \pi_{\tau},Z)$ est un polynôme en $Z$ 
 (puisque les facteurs locaux ont des pôles mais jamais de zéros). 
  Donc 
$$L_{S}(\check{\pi}  \times \pi,Z)=L_{S}(\check{\pi} \times \sigma,Z)=\prod_{\tau }L_{S}(\check{\pi} \times \tau,Z) ^{\dim \mf V_{\tau}}=\prod_{\tau }L_{S}(\check{\pi} \times \pi_{\tau},Z)^{\dim \mf V_{\tau}}$$  n'a pas de pôle. 
 Pourtant d'après le théorème B.10 de \cite{laurent-inventiones} (dû à Jacquet, Shahidi, Shalika), 
 $L_{S}(\check{\pi}  \times \pi,Z)$ a un pôle en $Z=q^{-1}$, ce qui amène  une contradiction. 
     Donc on sait maintenant que $\sigma$ est irréductible. Pour tout $\iota$ on sait d'après le \lemref{lem-sigma-dans-coho} que $\sigma$ est $\iota$-pur et  la connaissance de son déterminant (par la théorie du corps de classe) implique alors que le $\iota$-poids est nul. Donc $\sigma$ est pur de poids $0$. \cqfd
  
  \begin{rem} On mentionne pour le lecteur que les conséquences déjà très importantes de la correspondance de Langlands pour $GL_{r}$, expliquées dans le chapitre VII de \cite{laurent-inventiones},  ont été encore  étendues dans des travaux récents de Deligne et Drinfeld 
  \cite{deligne-finitude,esnault-deligne,drinfeld-del-conj}
  sur les questions de rationalité 
  des traces des Frobenius et sur l'indépendance de $\ell$ pour les faisceaux $\ell$-adiques sur les  variétés lisses sur $\Fq$. 
\end{rem}

\end{document}